\documentclass[11pt,twoside,]{amsart}
\usepackage[dvips]{graphicx}
\usepackage{tikz,tkz-berge,tkz-graph}
\usepackage{longtable}
\usepackage{amsmath, amsthm, amscd, amsfonts, amssymb, color, comment,verbatim}
\usepackage[bookmarksnumbered, plainpages]{hyperref}
 \usepackage[usenames,dvipsnames]{pstricks}
 \usepackage{epsfig}
 \usepackage{pst-grad} 
 \usepackage{pst-plot} 
\newtheorem{thm}{Theorem}[section]
\newtheorem{cor}[thm]{Corollary}
\newtheorem{lem}[thm]{Lemma}
\newtheorem{prop}[thm]{Proposition}
\newtheorem{defn}[thm]{Definition}
\newtheorem{conj}[thm]{Conjecture}
\newtheorem{rem}[thm]{Remark}

\numberwithin{equation}{section}

\begin{document}

\oddsidemargin 0mm
\evensidemargin 0mm

\thispagestyle{plain}

\vspace{5cc}
\begin{center}


{\large\bf Kaplansky's zero divisor and unit conjectures on elements\\ with supports of size $3$}
\rule{0mm}{6mm}\renewcommand{\thefootnote}{}
\footnotetext{
{\rule{2.4mm}{0mm}Keywords and Phrases.  Kaplansky's zero divisor conjecture, Kaplansky's unit conjecture, group ring, torsion-free group, zero divisor.}}

\vspace{1cc}
{\large\it Alireza Abdollahi and Zahra Taheri}


\vspace{1cc}
\parbox{27cc}{{\small

\textbf{Abstract.} Kaplansky's zero divisor conjecture (unit conjecture, respectively)  states that for a torsion-free group $G$ and a field $\mathbb{F}$, the group ring $\mathbb{F}[G]$ has no zero divisors (has no unit with support of size greater than $1$).  In this paper, we study  possible zero divisors and units in $\mathbb{F}[G]$ whose supports have size $3$.
For any field $\mathbb{F}$ and all torsion-free groups $G$, we prove that if $\alpha \beta=0$ for some non-zero $\alpha, \beta \in \mathbb{F}[G]$ such that $|supp(\alpha)|=3$, then $|supp(\beta)|\geq 10$. If $\mathbb{F}=\mathbb{F}_2$ is the field with 2 elements, the latter result can be improved so that $|supp(\beta)|\geq 20$. This improves a result in [J. Group Theory, 16 (2013), no. 5, 667-693]. 
Concerning the unit conjecture, we  prove that if $\alpha \beta=1$ for some $\alpha, \beta \in \mathbb{F}[G]$ such that $|supp(\alpha)|=3$, then $|supp(\beta)|\geq 9$. The latter improves a part of a result in  [Exp. Math., 24 (2015), 326-338] to arbitrary fields. }}

\end{center}

\vspace{1cc}


\section{\bf Introduction and Results}\label{Int}
Let $R$ be a ring. A non-zero element $\alpha$ of  $R$ is called a zero divisor if $\alpha\beta= 0$ or $\beta\alpha = 0$ for
some non-zero element $\beta \in R$. Let $G$ be a group. Denote by $R[G]$ the group ring of $G$ over $R$. If $R$ contains a zero divisor, then clearly so does $R[G]$. Also, if $G$ contains a non-identity torsion element $x$ of finite order $n$, then $R[G]$ contains zero divisors $\alpha=1-x$ and $\beta=1+x+\cdots+x^{n-1}$, since $\alpha \beta=0$.  Around $1950$, Irving Kaplansky conjectured that existence of a zero divisor in a group ring depends only on the existence of such elements in the ring or non-trivial torsions in the group by stating one of the most challenging problems in the field of group rings \cite{kap}.
\begin{conj}[Kaplansky's zero divisor conjecture] \label{conj-zero} Let $\mathbb{F}$ be a field and $G$ be a torsion-free group. Then $\mathbb{F}[G]$ does not contain a zero divisor.
\end{conj}

Another famous problem, namely the unit conjecture, also proposed by Kaplansky \cite{kap}, states that:
\begin{conj}[Kaplansky's unit conjecture] \label{conj-unit} Let $\mathbb{F}$ be a field and $G$ be a torsion-free group. Then $\mathbb{F}[G]$ has
no non-trivial units (i.e., non-zero scalar multiples of group elements).
\end{conj}
It can be shown that the zero divisor conjecture is true if the unit conjecture has an affirmative solution (see Lemma 13.1.2 in \cite{pass1}).

Over the years, some partial results have been obtained on Conjecture \ref{conj-zero} and it has been confirmed for special classes of groups which are torsion free. One of the first known special families which satisfy Conjecture \ref{conj-zero} are unique product groups \cite[Chapter 13]{pass1}, in particular ordered groups. Furthermore, by the fact that Conjecture \ref{conj-zero} is known to hold valid for amalgamated free products when the group ring of the subgroup over which the amalgam is formed satisfies the Ore condition \cite{lewin}, it is proved by Formanek \cite{form} that supersolvable groups are another families which satisfy Conjecture \ref{conj-zero}. Another result, concerning large major sorts of groups for which Conjecture \ref{conj-zero} holds in the affirmative, is obtained for elementary amenable groups \cite{krop}. The latter result covers the cases in which the group is polycyclic-by-finite, which was firstly studied in \cite{brown} and \cite{farkas}, and then extended in \cite{sni}. Some other affirmative results are obtained on congruence subgroups in \cite{laz} and \cite{farkas2}, and certain hyperbolic groups \cite{del}. 
Nevertheless, Conjecture \ref{conj-zero} has not been confirmed for any fixed field and it seems that confirming the conjecture even for the smallest finite field $\mathbb{F}_2$  with two elements is still out of reach.

The support of an element $\alpha=\sum_{x\in G}{a_x x}$ of $R[G]$, denoted by $supp(\alpha)$, is the set $\{x \in G \mid a_x\neq 0\}$. For any division ring $\mathbb{K}$ and all torsion-free group $G$, it is known that $\mathbb{K}[G]$ does not contain a zero divisor whose support is of size at most $2$ (see  \cite[Proposition 2.6]{Arxiv} and also \cite[Theorem 2.1]{pascal} when $\mathbb{K}$ is assumed to be a field), but it is not  known a similar result for group ring elements with the support of size $3$. By describing a combinatorial structure, named matched rectangles, Schweitzer \cite{pascal} showed that if $\alpha\beta=0$ for $\alpha,\beta\in \mathbb{F}_2[G]\setminus \{0\}$ when $\vert supp(\alpha)\vert=3$, then $\vert supp(\beta)\vert>6$. Also, with a computer-assisted approach, he showed that if $\vert supp(\alpha)\vert=3$, then $\vert supp(\beta)\vert>16$. 

Let $G$ be an arbitrary torsion-free group and let $\alpha \in \mathbb{F}[G]$ be a possible zero divisor such that  $\vert supp(\alpha)\vert=3$ and $\alpha \beta=0$ for some non-zero $\beta \in \mathbb{F}[G]$.  In this paper, we study the minimum possible size of the support of such an element  $\beta$. Let $\beta$ have minimum possible support size and $\mathbb{F}=\mathbb{F}_2$. In \cite[Definition 4.1]{pascal} a graph is associated to the non-degenerate $3\times |supp(\beta)|$ matched rectangle corresponding to $\alpha$ and $\beta$ and it is proved in \cite[Theorem 4.2]{pascal} that the graph is a simple cubic one without triangles. We call the graph Kaplansky graph of $(\alpha, \beta)$ over $\mathbb{F}_2$ and it is denoted  by $K_{\mathbb{F}_2}(\alpha,\beta)$. We extend such definition to the case that $\mathbb{F}$ is an arbitrary field and the corresponding Kaplansky graph is denoted  by $K_{\mathbb{F}}(\alpha,\beta)$. So, any Kaplansky graph is derived from a possible zero divisor with support of size $3$ in the group algebra of a torsion-free group over the field $\mathbb{F}$. In fact $K_{\mathbb{F}}(\alpha,\beta)$ is the induced subgraph  on the set $supp(\beta)$ of the Cayley graph $Cay(G,S)$ \footnote{By a Cayley graph $Cay(G,S)$ for a group $G$ and a subset $S$ of $G$ with $1\not\in S=S^{-1}$, is the graph whose vertex set is $G$ and two vertices $g_1,g_2$ are adjacent if $g_1g_2^{-1} \in S$.}, where $S=\{h^{-1} h' \;|\; h,h'\in supp(\alpha), h\not=h' \}$. Here we study  forbidden subgraphs of Kaplansky graphs. Our main results on Conjecture \ref{conj-zero} are the followings.

\begin{figure}[ht]
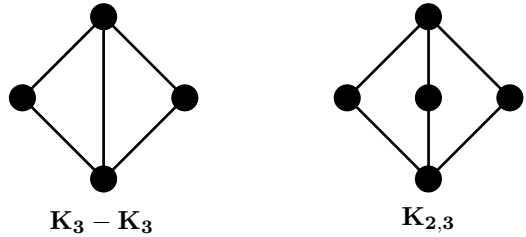

\centering
\psscalebox{0.9 0.9} 
{

}
\caption{Two forbidden subgraphs of the Kaplansky graph over $\mathbb{F}$}\label{forbiddens-F}
\end{figure}

\begin{thm}\label{thm-F}
None of the graphs in Figure \ref{forbiddens-F} can be isomorphic to a subgraph of any Kaplansky graph over any field $\mathbb{F}$.
\end{thm}

\begin{thm}
Let $\alpha$ and $\beta$ be non-zero elements of the group algebra of any torsion-free group over an arbitrary field. If $|supp(\alpha)|=3$ and $\alpha \beta=0$ then $|supp(\beta)|\geq 10$. 
\end{thm}

\begin{thm}\label{mainthm}
None of the graphs in Table \ref{tab-forbiddens} can be isomorphic to a subgraph of any Kaplansky graph over $\mathbb{F}_2$.
\end{thm}
In Appendix \ref{S2-app} some details of our computations needed in the proof of Theorem \ref{mainthm} are given for the reader's convenience. 

The following result improves a result in \cite{pascal}. 
\begin{thm}
Let $\alpha$ and $\beta$ be non-zero elements of the group algebra of any torsion-free group over the field with two elements. If $|supp(\alpha)|=3$ and $\alpha \beta=0$ then $|supp(\beta)|\geq 20$. 
\end{thm}

The best known result on Conjecture \ref{conj-unit}, which has the purely group-theoretic approach, is concerned with unique product groups \cite{pass1,pass2}. 
The latter result covers ordered groups, in particular torsion-free nilpotent groups. Nevertheless, it is still unknown whether or not Conjecture \ref{conj-unit} do hold true for supersolvable torsion-free groups.  Dykema et al. \cite{dyk} have shown that there exist no $\gamma , \delta \in \mathbb{F}_2[G]$ such that $\gamma \delta=1$, where $|supp(\gamma)|=3$ and $|supp(\delta)|\leq 11$. Concerning Conjecture \ref{conj-unit},  we prove the following result which improves a part of the result in \cite{dyk} to arbitrary fields. 

\begin{thm}
Let $\gamma$ and $\delta$ be elements of the group algebra of any torsion-free group over an arbitrary field. If $|supp(\gamma)|=3$ and $\gamma \delta =1$ then $|supp(\delta)|\geq 9$. 
\end{thm}

It is known that $\mathbb{F}[G]$ contains a zero divisor if and only if it contains a non-zero element whose square is zero (see \cite{pass1}).  Using the latter fact,  it is mentioned in \cite[p. 691]{pascal} that it is sufficient to check  Conjecture \ref{conj-zero} only for the case that $|supp(\alpha)|=|supp(\beta)|$, but in the construction that, given a zero divisor produces an element of square zero, it is not clear how the length changes.  We clarify the latter by the following.
\begin{prop}
If $\mathbb{F}[G]$ has no non-zero element $\alpha$ with $|supp(\alpha)|\leq k$ such that $\alpha^2=0$, then there exist no non-zero elements $\alpha_1,\alpha_2 \in \mathbb{F}[G]$ such that $\alpha_1\alpha_2=0$ and $|supp(\alpha_1)||supp(\alpha_2)|\leq k$.
\end{prop}

\begin{center}
\begin{longtable}{|c|c|c|c|}
\caption{{Forbidden subgraphs of  Kaplansky graphs over $\mathbb{F}_2$}}\label{tab-forbiddens}\\
\hline 
\endfirsthead
\multicolumn{4}{c}%
{{\tablename\ \thetable{} -- continued from previous page}} \\
\endhead
\hline \multicolumn{4}{|r|}{{Continued on next page}} \\\hline
\endfoot
\hline 
\endlastfoot
\psscalebox{0.62 0.62} 
{
\begin{pspicture}(0,-2.0985577)(3.594231,2.0985577)
\psdots[linecolor=black, dotsize=0.4](3.3971152,0.30144233)
\psdots[linecolor=black, dotsize=0.4](1.7971153,-1.2985576)
\psdots[linecolor=black, dotsize=0.4](0.19711533,0.30144233)
\psdots[linecolor=black, dotsize=0.4](1.7971153,1.9014423)
\psdots[linecolor=black, dotsize=0.4](1.7971153,0.30144233)
\psline[linecolor=black, linewidth=0.04](1.7971153,1.9014423)(0.19711533,0.30144233)(1.7971153,-1.2985576)(3.3971152,0.30144233)(1.7971153,1.9014423)(1.7971153,0.30144233)(1.7971153,-1.2985576)
\rput[bl](1.2,-2.0985577){$\mathbf{1)\ K_{2,3}}$}
\end{pspicture}
}&\psscalebox{0.62 0.62} 
{
\begin{pspicture}(0,-2.0985577)(3.594231,2.0985577)
\psdots[linecolor=black, dotsize=0.4](3.3971152,0.30144233)
\psdots[linecolor=black, dotsize=0.4](1.7971153,-1.2985576)
\psdots[linecolor=black, dotsize=0.4](0.19711533,0.30144233)
\psdots[linecolor=black, dotsize=0.4](1.7971153,1.9014423)
\psdots[linecolor=black, dotsize=0.4](1.7971153,0.30144233)
\psdots[linecolor=black, dotsize=0.4](2.5971153,-0.49855766)
\psline[linecolor=black, linewidth=0.04](1.7971153,1.9014423)(0.19711533,0.30144233)(1.7971153,-1.2985576)(2.5971153,-0.49855766)(3.3971152,0.30144233)(1.7971153,1.9014423)(1.7971153,0.30144233)(1.7971153,-1.2985576)
\rput[bl](0.7,-2.0985577){$\mathbf{2) \ C_4--C_5}$}
\end{pspicture}
}&\psscalebox{0.62 0.62} 
{
\begin{pspicture}(0,-2.0985577)(3.594231,2.0985577)
\psdots[linecolor=black, dotsize=0.4](3.3971152,0.30144233)
\psdots[linecolor=black, dotsize=0.4](1.7971153,-1.2985576)
\psdots[linecolor=black, dotsize=0.4](0.19711533,0.30144233)
\psdots[linecolor=black, dotsize=0.4](1.7971153,1.9014423)
\psdots[linecolor=black, dotsize=0.4](1.7971153,0.30144233)
\psdots[linecolor=black, dotsize=0.4](2.5971153,-0.49855766)
\psline[linecolor=black, linewidth=0.04](1.7971153,1.9014423)(0.19711533,0.30144233)(1.7971153,-1.2985576)(2.5971153,-0.49855766)(3.3971152,0.30144233)(1.7971153,1.9014423)(1.7971153,0.30144233)(1.7971153,-1.2985576)
\psdots[linecolor=black, dotsize=0.4](2.5971153,1.1014423)
\rput[bl](0.7,-2.0985577){$\mathbf{3) \ C_4--C_6}$}
\end{pspicture}
}&\psscalebox{0.62 0.62} 
{
\begin{pspicture}(0,-1.8985577)(3.594231,1.8985577)
\psdots[linecolor=black, dotsize=0.4](0.9971155,1.7014424)
\psdots[linecolor=black, dotsize=0.4](2.5971155,1.7014424)
\psdots[linecolor=black, dotsize=0.4](1.7971154,1.3014424)
\psdots[linecolor=black, dotsize=0.4](1.7971154,0.5014423)
\psdots[linecolor=black, dotsize=0.4](2.5971155,-0.29855767)
\psdots[linecolor=black, dotsize=0.4](1.7971154,-1.0985577)
\psdots[linecolor=black, dotsize=0.4](0.9971155,-0.29855767)
\psdots[linecolor=black, dotsize=0.4](3.3971155,0.90144235)
\psdots[linecolor=black, dotsize=0.4](0.19711548,0.90144235)
\psline[linecolor=black, linewidth=0.04](0.9971155,1.7014424)(0.19711548,0.90144235)(0.9971155,-0.29855767)(1.7971154,-1.0985577)(2.5971155,-0.29855767)(3.3971155,0.90144235)(2.5971155,1.7014424)(1.7971154,1.3014424)(0.9971155,1.7014424)
\psline[linecolor=black, linewidth=0.04](1.7971154,1.3014424)(1.7971154,0.5014423)(0.9971155,-0.29855767)
\psline[linecolor=black, linewidth=0.04](1.7971154,0.5014423)(2.5971155,-0.29855767)
\rput[bl](0.3,-1.8985577){$\mathbf{4) \ C_4-C_5(-C_5-)}$}
\end{pspicture}
}\\
\end{longtable}
\end{center}


\section{\bf Kaplansky graphs over $\mathbb{F}$ and some of their properties}\label{S1-2}
Throughout this paper let $G$ be a torsion-free group and  $\alpha=\alpha_1 h_1+\alpha_2h_2+\alpha_3h_3 \in \mathbb{F}[G]$ such that $\vert supp(\alpha)\vert=3$. Suppose further that  $\alpha \beta =0$ for some non-zero $\beta \in \mathbb{F}[G]$ and assume that  $n:=|supp(\beta)|$ is minimum with respect to the latter property and $\beta=\beta_1 g_1+\beta_2g_2+\cdots+\beta_ng_n$. So, $n\geq 3$ \cite{pascal}. 

\begin{lem}[See also Lemma 2.1 of \cite{Arxiv}]\label{supp}
$\langle h_i^{-1} supp(\alpha) \rangle=\langle  supp(\beta) g_j^{-1} \rangle$ for all $i\in \{1,2,3\}$ and $j\in\{1,2,\ldots,n\}$.
\end{lem}
\begin{proof}
Let $i\in \{1,2,3\}$, $j\in\{1,\dots,n\}$, $H=\langle h_i^{-1} supp(\alpha) \rangle$, $K=\langle  supp(\beta) g_j^{-1} \rangle$ and $\{t_1,t_2,\ldots,t_k\}$ be a set of right coset representatives of $H$ in $G$ such that if $k_{j'}\in supp(\beta) g_j^{-1}$, then $k_{j'}\in Ht_{i'}$ for some $i'\in \{1,2,\ldots,k\}$. Suppose that $k>1$. Since $\alpha \beta=0$ and $Ht_{l_1}\cap Ht_{l_2}=\varnothing$ for all distinct $l_1,l_2\in \{1,2,\ldots,k\}$,  $(h_i^{-1}\alpha )(a_1h_1'+a_2h_2'+\cdots+a_lh_{i_l}')t_l=0$ for some $l\in \{1,2,\ldots,k\}$, $h_1',h_2',\ldots,h_{i_l}' \in H$ and $\{a_1,a_2,\ldots,a_{i_l}\}\subseteq supp(\beta)$, where $i_l<n$, that is a contradiction with the minimality of $n$ because $\alpha(a_1h_1'+a_2h_2'+\cdots+a_lh_{i_l}')=0$. So $k=1$, $\alpha(\beta_1h_1't_1+\beta_2h_2't_1+\cdots+\beta_nh_n't_1)=0$ and $supp(\beta) g_j^{-1}=\{h_1't_1,h_2't_1,\ldots,h_n't_1\}$ where $h_1',h_2',\ldots,h_n' \in H$. Since $1\in supp(\beta) g_j^{-1}\subseteq Ht_1$, $1=ht_1$ for some $h\in H$ and so $t_1\in H$. Therefore, $supp(\beta) g_j^{-1} \subseteq H$ and so $K\leq H$.

Since $\alpha \beta =0$, $h_ig_j=h_{i'}g_{j'}$ for some $(i',j')\in A=\{1,2,3\}\times \{1,2,\ldots,n\}$ such that $i\neq i'$ and $j\neq j'$. So, $h_i^{-1}h_{i'}=g_jg_{j'}^{-1}=(g_{j'}g_j^{-1})^{-1}\in K$. Furthermore, $h_i^{-1}h_i=1=g_{j}g_j^{-1}\in K$. Therefore, $h_i^{-1}h_i,h_i^{-1}h_{i'}\in K$ where $supp(\alpha)=\{h_i,h_{i'},h_{i''}\}$. Similarly, $h_{i''}g_j=h_sg_t$ for some $(s,t)\in A$ such that $s\neq i''$ and $t\neq j$. Since $h_s\in supp(\alpha)\setminus \{h_{i''}\}$, $h_s=h_{i}$ or $h_s=h_{i'}$. If $h_s=h_{i}$, then $h_i^{-1}h_{i''}=g_tg_{j}^{-1}\in K$. If $h_s=h_{i'}$, then $h_{i'}^{-1}h_{i''}=g_tg_{j}^{-1}\in K$. So, $h_i^{-1}h_{i''}\in K$ because $h_{i'}^{-1}h_{i''},h_i^{-1}h_{i'}\in K$ and $h_i^{-1}h_{i''}=h_i^{-1}h_{i'}h_{i'}^{-1}h_{i''}$. Therefore, $h_i^{-1} supp(\alpha) \subseteq K$ which implies that $H=K$ because $K\leq H$.
\end{proof}

\begin{rem}\label{r-G}
{\rm
By Lemma \ref{supp}, it can be supposed that $G=\langle h_i^{-1} supp(\alpha) \rangle$ for all $i\in \{1,2,3\}$. Since $1\in h_i^{-1} supp(\alpha)$, without loss of generality we may assume that $supp(\alpha)=\{1,h_2,h_3\}$ and $G=\langle supp(\alpha) \rangle$. 
}
\end{rem}

The size of $S=\{h^{-1} h' \;|\; h,h'\in supp(\alpha), h\not=h' \}$ is at most $6$ and we prove that $S$ should have its largest possible size.
\begin{lem}\label{size of S}
$|S|=6$.
\end{lem}
\begin{proof}
 Suppose, for a contradiction, that $|S|<6$. Then $h_{i}^{-1}h_j=h_{i'}^{-1}h_{j'}$ for some $(i,j)\not= (i',j')$ and $(i,i')\not=(j,j')$.
 It follows that $(i',j')=(j,i)$ or $(j,k)$ or $(k,i)$, where $k \in\{1,2,3\}\setminus\{i,j\}$. If $(i',j')=(j,i)$, then $(h_i^{-1}h_j)^2=1$ and since the group is torsion-free, $h_i=h_j$, a contradiction. If $(i',j')=(j,k)$, then $h_i^{-1}h_k=(h_i^{-1}h_j)^2$ and if $(i',j')=(k,i)$ then 
 $h_i^{-1}h_k=h_i^{-1} h_j$. It follows that $H:=\langle h_i^{-1} supp(\alpha) \rangle=\langle h_i^{-1} h_j \rangle$ is the infinite cyclic group.
 It follows from  Lemma \ref{supp} that  $h_i^{-1} \alpha, \beta g_1^{-1} \in \mathbb{F}_2[H]$. 
 Now $(h_i^{-1} \alpha) (\beta g_1^{-1})=0$ contradicts the fact that the group algebra of the infinite cyclic group has no zero-divisor \cite{pass1}. This completes  the proof.   
\end{proof}

If $A=\{1,2,3\}\times \{1,2,\ldots,n\}$, for all $(i,j)\in A$ there must be an $(i',j')\in A$ such that $i\neq i'$, $j\neq j'$ and  $h_ig_j=h_{i'}g_{j'}$ because $\alpha\beta=(\alpha_1 h_1+\alpha_2h_2+\alpha_3h_3)(\beta_1 g_1+\beta_2g_2+\cdots+\beta_ng_n)=0$. 

The Kaplansky graph of  $\alpha$ and $\beta$ over $\mathbb{F}$ can be defined and it is denoted by $K_{\mathbb{F}}(\alpha,\beta)$. The vertex set of $K_{\mathbb{F}}(\alpha,\beta)$ is $supp(\beta)$ and two vertices $g_i$ and $g_j$ are adjacent, denoted by $g_i\sim g_j$, whenever $h_{i'}g_i=h_{j'} g_j$ for some distinct $i',j' \in \{1,2,3\}$. So, any Kaplansky graph over $\mathbb{F}$ is derived from a possible zero divisor with support of size $3$ in the group algebra of a torsion-free group over the field $\mathbb{F}$. In the following, we give some properties of Kaplansky graphs. 

\begin{lem}
$K_{\mathbb{F}}(\alpha,\beta)\cong K_{\mathbb{F}}(x \alpha, \beta y)$ for all $x,y \in G$. 
\end{lem}  
\begin{proof}
Note that the being minimum of the support size of $\beta y$ in respect to the condition $(x \alpha) (\beta y)=0$ is obvious as the support size of $\beta$ in $\alpha \beta=0$. Note that if $S$ is equal to the set $\{h^{-1} h' \;|\; h,h'\in supp(x\alpha), h\not=h' \}$, the map on  $G$ defined by $g\mapsto gy$ for all $g\in G$ is a graph isomorphism on $Cay(G,S)$. This completes the proof.
\end{proof}


\begin{figure}[ht]
\psscalebox{0.9 0.9} 
{
\begin{pspicture}(0,-2.105)(4.72,2.105)
\psdots[linecolor=black, dotsize=0.4](2.4,1.495)
\psdots[linecolor=black, dotsize=0.4](1.2,1.095)
\psdots[linecolor=black, dotsize=0.4](0.8,-0.105)
\psdots[linecolor=black, dotsize=0.4](1.2,-1.305)
\psdots[linecolor=black, dotsize=0.4](2.4,-1.705)
\psdots[linecolor=black, dotsize=0.4](3.6,-1.305)
\psdots[linecolor=black, dotsize=0.4](4.0,-0.105)
\psdots[linecolor=black, dotsize=0.4](3.6,1.095)
\psline[linecolor=black, linewidth=0.04](3.6,1.095)(2.4,1.495)(1.2,1.095)(0.8,-0.105)(1.2,-1.305)(2.4,-1.705)(3.6,-1.305)(4.0,-0.105)
\psline[linecolor=black, linewidth=0.04, linestyle=dotted, dotsep=0.10583334cm](3.6,1.095)(4.0,-0.105)
\end{pspicture}
}
\caption{A cycle of length $k$ in $K_{\mathbb{F}}(\alpha,\beta)$}\label{f-cycle}
\end{figure}
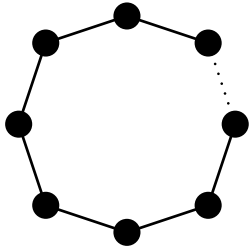

\begin{defn}\label{tuples}
 {\rm To any cycle $C$ of $K_{\mathbb{F}}(\alpha,\beta)$ of length $k$ as Figure \ref{f-cycle}, we assign a $2k$-tuple
$T_C=[a_1,b_1,a_2,b_2,\dots,a_k,b_k]$,  where $a_1, b_1,a_2,b_2,\dots,a_k,b_k \in supp(\alpha)$ satisfying the following relations:
\begin{equation}\label{e-cycle1}
R(T_C): \left\{
\begin{array}{l}
a_1g_1'=b_1g_2'\\
a_2g_2'=b_2g_3'\\
\vdots\\
a_kg_k'=b_kg_1'\\
\end{array} \right.
\end{equation}
where $g_1',g_2',g_3',\dots,g_k'\in supp(\beta)$ are vertices of $C$ such that $g_i'\sim g_{i+1}'$, for all $i\in \{1,2,\ldots, k-1\}$, and $g_1'\sim g_k'$.
Also, we can derive from the relations \ref{e-cycle1} that 
$r(T_C)=(a_1^{-1}b_1)(a_2^{-1}b_2)\cdots(a_k^{-1}b_k)$ is equal to $1$. It follows from Lemma \ref{size of S} that if $[a'_1,b'_1,\dots,a'_k,b'_k]$ is another $2k$-tuple assigning to $C$ as above, then $[a'_1,b'_1,\dots,a'_k,b'_k]$ is one of the following $2k$-tuples:

\begin{equation*}
\begin{matrix}
[a_1,b_1,a_2,b_2,\ldots, a_{k-1},b_{k-1}, a_k,b_k],\\
[a_k,b_k,a_1,b_1,\ldots, a_{k-2},b_{k-2}, a_{k-1},b_{k-1}],\\
\vdots\\
[a_2,b_2,a_3,b_3,\ldots, a_{k},b_{k},a_1,b_1],\\
[b_1,a_1,b_k,a_k, \ldots, b_3,a_3,b_2,a_2],\\
[b_2,a_2,b_1,a_1,\ldots, b_4,a_4,b_3,a_3],\\
\vdots\\
[b_k,a_k,b_{k-1},a_{k-1},\ldots, b_2,a_2,b_1,a_1].
\end{matrix}
\end{equation*}
The set of all such $2k$-tuples will be denoted by $\mathcal{T}(C)$. Also, a member of the set $\mathcal{R}(C)=\{R(T_C) | T_C\in \mathcal{T}(C)\}$ is called the corresponding relations of $C$.
}
\end{defn}
\begin{defn}{\rm
Let $C$ be a cycle of $K_{\mathbb{F}}(\alpha,\beta)$ of length $k$. Since $r(T_1)=1$ if and only if $r(T_2)=1$, for all $T_1,T_2\in \mathcal{T}(C)$, a member of $\{r(T_C) | T_C\in \mathcal{T}(C)\}$ is given as a representative and denoted by $r(C)$, and $r(C)=1$ is called the relation of $C$. 
}\end{defn}

\begin{defn}\label{def-equ}{\rm
Let $C$ and $C'$ be two cycles of length $k$ in $K_{\mathbb{F}}(\alpha,\beta)$. We say that these two cycles  are equivalent, if $\mathcal{T}(C)\cap \mathcal{T}(C')\not=\varnothing$.
}\end{defn}
\begin{rem}{\rm
If $C$ and $C'$ are two equivalent cycles of length $k$ in $K_{\mathbb{F}}(\alpha,\beta)$, then  $\mathcal{T}(C)=\mathcal{T}(C')$.
}\end{rem}


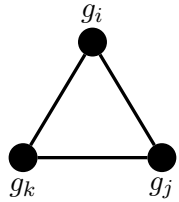
\begin{figure}[ht]
\centering
\begin{tikzpicture}[scale=0.7]
\node (2) [circle, minimum size=3pt, fill=black, line width=0.625pt, draw=black] at (125.0pt, -125.0pt)  {};
\node (3) [circle, minimum size=3pt, fill=black, line width=0.625pt, draw=black] at (50.0pt, -125.0pt)  {};
\node (1) [circle, minimum size=3pt, fill=black, line width=0.625pt, draw=black] at (87.5pt, -62.5pt)  {};
\draw [line width=1.25, color=black] (1) to  (2);
\draw [line width=1.25, color=black] (1) to  (3);
\draw [line width=1.25, color=black] (3) to  (2);
\node at (125.0pt, -141.375pt) {\textcolor{black}{$g_j$}};
\node at (50.0pt, -141.375pt) {\textcolor{black}{$g_k$}};
\node at (87.5pt, -47.125pt) {\textcolor{black}{$g_i$}};
\end{tikzpicture}
\caption{The triangle $K_3$ in a Kaplansky graph}\label{f-1}
\end{figure}

Now let $\mathbb{F}=\mathbb{F}_2$. Obviously, for each $(i,j)$ and $(i,j')$ in $A$ where $j\neq j'$ we have  $h_ig_j\not=h_ig_{j'}$. Also, for each $(i,j)$ and $(i',j)$ in $A$ where $i\neq i'$ we have  $h_ig_j\not=h_{i'}g_j$. Therefore, there is a matched rectangle $M$ corresponding to $(\alpha,\beta)$ (see \cite[Definition 4.1]{pascal}) that is non-degenerate and so the underlying graph $K(M)$ following \cite[Definition 4.1]{pascal} can be defined. We call $K(M)$, the Kaplansky graph of $(\alpha, \beta)$ over $\mathbb{F}_2$ and it is denoted by $K_{\mathbb{F}_2}(\alpha,\beta)$. The vertex set of the Kaplansky graph is $supp(\beta)$ and two vertices $g_i$ and $g_j$ are adjacent whenever $h_{i'}g_i=h_{j'} g_j$ for some distinct $i',j' \in \{1,2,3\}$.  

The following theorem is obtained in \cite{pascal}.
\begin{thm}[Theorem 4.2 of \cite{pascal}]\label{thm-graph}
Any Kaplansky graph over $\mathbb{F}_2$ is a connected simple cubic one  containing no subgraph isomorphic to a triangle. 
\end{thm} 
\begin{proof}
The proof is essentially the same as the proof of \cite[Theorem 4.2]{pascal}, but only note that the connectedness follows from the way we have chosen $\beta$ of minimum support size with respect to the property $\alpha \beta=0$.  
\end{proof}
\begin{rem}
{\rm By a triangle in Theorem \ref{thm-graph}, we mean a subgraph such as Figure \ref{f-1}, where $g_i,g_j,g_k\in supp(\beta)$, with the corresponding relations as
\begin{equation}\label{e-1}
\left\{
\begin{array}{l}
a_1g_i=b_1g_j\\
a_2g_j=b_2g_k\\
a_3g_k=b_3g_i,\\
\text{for some } a_s,b_t \in supp(\alpha) \text{ where } s,t\in \{1,2,3\} \\
\text{and } a_1\not=b_1\not=a_2\not=b_2\not=a_3\not=b_3.
\end{array} \right.
\end{equation}
}
\end{rem}

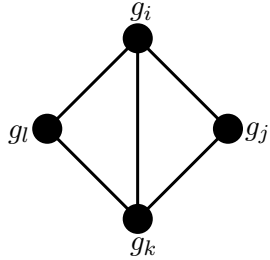
\begin{figure}[h]
\centering
\psscalebox{1.0 1.0} 
{
\begin{pspicture}(0,-1.905)(3.91,1.905)
\psdots[linecolor=black, dotsize=0.4](2.0,1.295)
\psdots[linecolor=black, dotsize=0.4](2.0,-1.105)
\psdots[linecolor=black, dotsize=0.4](0.8,0.095)
\psdots[linecolor=black, dotsize=0.4](3.2,0.095)
\psline[linecolor=black, linewidth=0.04](2.0,1.295)(0.8,0.095)(2.0,-1.105)(3.2,0.095)(2.0,1.295)(2.0,-1.105)(2.0,-1.105)
\rput[bl](1.9,1.53){$g_i$}
\rput[bl](3.43,-0.1){$g_j$}
\rput[bl](1.9,-1.6){$g_k$}
\rput[bl](0.28,-0.1){$g_l$}
\end{pspicture}
}
\caption{Two triangles with a common edge in the Kaplansky graph over $\mathbb{F}$}\label{f-K3-K3}
\end{figure}

\begin{thm}\label{K3-K3}
Kaplansky graphs over $\mathbb{F}$ contain no subgraphs isomorphic to the graph in Figure \ref{f-K3-K3} i.e. two triangles with one edge in common. 
\end{thm}
\begin{proof}
Similar to Theorem \ref{thm-graph}, it can be seen that a Kaplansky graph over $\mathbb{F}$ is also a connected simple graph containing no subgraph isomorphic to a triangle as Figure \ref{f-1} with the corresponding relations \ref{e-1}.  So, if $K_{\mathbb{F}}(\alpha,\beta)$ contains a subgraph isomorphic to a triangle with vertices $g_i,g_j$ and $g_k$, the corresponding relations of such triangle are as follows
\begin{equation}\label{e-K3 over F}
ag_i=bg_j=cg_k\text{, where } \{a,b,c\}=supp(\alpha)
\end{equation}
Suppose that $K_{\mathbb{F}}(\alpha,\beta)$ contains two triangle with one edge in common as Figure \ref{f-K3-K3}, where $g_i,g_j,g_k,g_l\in supp(\beta)$. With the discussion above and by the relation \ref{e-K3 over F},  $a_1g_i=b_1g_k=c_1g_j$  and $a_1g_i=b_1g_k=c_2g_l$ where $\{a_1,b_1,c_1\}=\{a_1,b_1,c_2\}= supp(\alpha)$. So,  $c_2g_l=a_1g_i=b_1g_k=c_1g_j$ where $\{a_1,b_1,c_1\}=\{a_1,b_1,c_2\}= supp(\alpha)$. Therefore, $c_1=c_2$ and so $g_j=g_l$, a contradiction. So, the graph $K_{\mathbb{F}}(\alpha,\beta)$ contains no two triangles with one edge in common. 
\end{proof}

\begin{figure}[ht]
\centering
\begin{tikzpicture}[scale=0.7]
\node (1) [circle, minimum size=3pt, fill=black, line width=0.625pt, draw=black] at (75.0pt, -25.0pt)  {};
\node (2) [circle, minimum size=3pt, fill=black, line width=0.625pt, draw=black] at (125.0pt, -75.0pt)  {};
\node (3) [circle, minimum size=3pt, fill=black, line width=0.625pt, draw=black] at (25.0pt, -75.0pt)  {};
\node (5) [circle, minimum size=3pt, fill=black, line width=0.625pt, draw=black] at (75.0pt, -125.0pt)  {};
\draw [line width=1.25, color=black] (1) to  (2);
\draw [line width=1.25, color=black] (1) to  (3);
\draw [line width=1.25, color=black] (2) to  (5);
\draw [line width=1.25, color=black] (3) to  (5);
\node at (75.0pt, -10.625pt) {\textcolor{black}{$g_i$}};
\node at (140.875pt, -75.0pt) {\textcolor{black}{$g_j$}};
\node at (10.125pt, -75.0pt) {\textcolor{black}{$g_l$}};
\node at (75.0pt, -141.375pt) {\textcolor{black}{$g_k$}};
\end{tikzpicture}
\caption{An square in the Kaplansky graph over $\mathbb{F}$}\label{f-C4}
\end{figure}
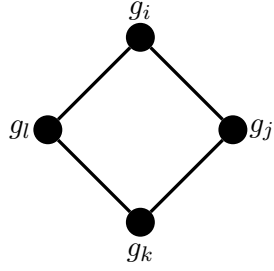

\begin{rem}\label{r-K3-K3}
{\rm It follows from Theorem \ref{K3-K3} that if $K_{\mathbb{F}}(\alpha,\beta)$ contains an square  as Figure \ref{f-C4}, where  $g_i,g_j,g_k,g_l\in supp(\beta)$, then $g_i \not\sim g_k$ and $g_j \not\sim g_l$. So, the corresponding relations of such square are as follows:
\begin{equation}\label{e-4}
\left\{
\begin{array}{l}
a_1g_i=b_1g_j\\
a_2g_j=b_2g_k\\
a_3g_k=b_3g_l\\
a_4g_l=b_4g_i,\\
\text{where } a_s,b_s \in supp(\alpha) \text{ for all } s\in \{1,2,3,4\}\\
\text{and } a_1\not= b_1\not= a_2 \not= b_2\not= a_3\not= b_3\not= a_4\not= b_4\not= a_1.
\end{array} \right.
\end{equation} }
\end{rem}

To finish this section, we consider Kaplansky graphs over $\mathbb{F}$ containing a subgraph isomorphic to an square. We have not been able as \cite[Theorem 4.2]{pascal} to prove  that squares are forbidden subgraphs for Kaplansky graphs even for the case that $\mathbb{F}=\mathbb{F}_2$, as triangles are so. However we show that by existence of squares, Kaplansky graphs over $\mathbb{F}$ gives us certain slightly significant relations on  elements of the support of a possible zero-divisor (see below, Theorem \ref{C4}).

\begin{rem}{\rm
The Baumslag-Solitar group $BS(m,n)$ is the group given by the presentation $\langle a,b\mid b a^m b^{-1}=a^n\rangle $. Such groups are solvable if and only if $\vert m\vert=1$ or $\vert  n\vert =1$ \cite{baumslag}. So, such latter groups and their quotients satisfy Conjecture \ref{conj-zero} \cite{pass1}.}
\end{rem}

\begin{table}[h]
\centering
\caption{The possible relations of squares in $K_{\mathbb{F}}(\alpha,\beta)$}\label{tab-C4}
\begin{tabular}{|c|l|l||c|l|l|}\hline
$n$&$R$&$ \ E \ $&$n$&$R$&$ \ E \ $\\\hline
$1$&$h_2^4=1$& $T $&$19$&$h_2 h_3^{-1} (h_2^{-1} h_3)^2=1$& $BS(2,1) $\\
$2$&$h_2^3 h_3=1$& $ A$&$20$&$h_2 h_3^{-2} h_2^{-1} h_3=1$& $BS(1,2)$\\
$3$&$h_2^3 h_3^{-1} h_2=1$& $ A$&$21$&$h_2 h_3^{-3} h_2=1$& $ *$\\
$4$&$h_2^2 h_3^2=1$& $BS(1,-1) $&$22$&$h_2 h_3^{-2} h_2 h_3=1$& $ *$\\
$5$&$h_2^2 h_3 h_2^{-1} h_3=1$& $ *$&$23$&$h_2 h_3^{-1} (h_3^{-1} h_2)^2=1$& $BS(-2,1) $\\
$6$&$h_2^2 h_3^{-1} h_2^{-1} h_3=1$& $BS(1,2) $&$24$&$(h_2 h_3^{-1} h_2)^2=1$& $ A$\\
$7$&$h_2^2 h_3^{-2} h_2=1$& $* $&$25$&$h_2 h_3^{-1} h_2 h_3^2=1$& $ *$\\
$8$&$h_2^2 h_3^{-1} h_2 h_3=1$& $ BS(1,-2)$&$26$&$h_2 h_3^{-1} h_2 h_3 h_2^{-1} h_3=1$& $ *$\\
$9$&$h_2 (h_2 h_3^{-1})^2 h_2=1$& $BS(1,-1) $&$27$&$(h_2 h_3^{-1})^2 h_2^{-1} h_3=1$& $BS(2,1) $\\
$10$&$(h_2 h_3)^2=1$& $ A$&$28$&$(h_2 h_3^{-1})^2 h_3^{-1} h_2=1$& $ BS(1,-2)$\\
$11$&$h_2 h_3 h_2 h_3^{-1} h_2=1$& $ BS(-2,1)$&$29$&$(h_2 h_3^{-1})^2 h_2 h_3=1$& $* $\\
$12$&$h_2 h_3^3=1$& $ A$&$30$&$(h_2 h_3^{-1})^3 h_2=1$& $A $\\
$13$&$h_2 h_3^2 h_2^{-1} h_3=1$& $BS(1,-2) $&$31$&$h_3^4=1$& $ T$\\
$14$&$h_2 h_3 h_2^{-2} h_3=1$& $ *$&$32$&$h_3^3 h_2^{-1} h_3=1$& $ A$\\
$15$&$h_2 h_3 h_2^{-1} h_3^{-1} h_2=1$& $ BS(2,1)$&$33$&$h_3 (h_3 h_2^{-1})^2 h_3=1$& $BS(1,-1) $\\
$16$&$h_2 h_3 h_2^{-1} h_3^2=1$& $BS(-2,1) $&$34$&$(h_3 h_2^{-1} h_3)^2=1$& $ A$\\
$17$&$h_2 (h_3 h_2^{-1})^2 h_3=1$& $* $&$35$&$(h_3 h_2^{-1})^3 h_3=1$& $ A$\\
$18$&$h_2 h_3^{-1} h_2^{-1} h_3^2=1$& $BS(2,1) $&$36$&$(h_2^{-1} h_3)^4=1$& $A $\\\hline
\end{tabular}
\end{table}

\begin{thm}\label{C4-1}
If a Kaplansky  graph over an arbitrary field has an square $C$, then there are $9$ non-equivalent cases for $C$ which $r(C)$ is one of the relations $5$, $7$, $14$, $17$, $21$, $22$, $25$, $26$ and $29$  in Table \ref{tab-C4}.
\end{thm}
\begin{proof} 
Let $C$ be an square as Figure \ref{f-C4} with the $8$-tuple $T_C=[a_1,b_1,a_2,b_2, a_3,b_3,a_4,b_4]$ and $R(T_C)$ as \ref{e-4}. Using GAP \cite{gap}, it can be seen that there are $258$ cases for $T_C$. By Definitions \ref{tuples} and \ref{def-equ}, and by using GAP \cite{gap}, it can be seen that these $258$ cases can be categorized into $36$ non-equivalent cases. The relations of such non-equivalent cases are listed in the column labelled by $R$ of Table \ref{tab-C4}. It can be shown that all the relations in this table, except the relations marked by ``$*$''s in the column labelled by $E$, lead to contradictions because the group $G$ generated by $h_2$ and $h_3$ with one of such relations has at least one of the following properties:
\begin{enumerate}
\item 
It is an abelian group,
\item
It is a quotient of $BS(1,k)$ or $BS(k,1)$ where $k\in \{-2,-1,1,2\}$,
\item
It has a non-trivial torsion element. 
\end{enumerate}
Each relation which leads to being $G$ an abelian group or having a non-trivial torsion element is marked by an $A$ or a $T$ in the column labelled by $E$, respectively. Also, if the group $G$ is a quotient of $BS(1,k)$ or $BS(k,1)$ where $k\in \{-2,-1,1,2\}$, we show this in the column $E$. Therefore, there are $9$ non-equivalent cases for existence of an square in a Kaplansky  graph over an arbitrary field.
\end{proof}

\begin{thm}\label{C4}
If a Kaplansky graph over an arbitrary field has an square $C$, then there exist non-trivial group elements $x$ and $y$ such that  $x^2=y^3$ and either $\{1,x,y\}$ or $\{1,y,y^{-1}x\}$ is the support of a zero divisor in $\mathbb{F}[G]$.
\end{thm}
\begin{proof}
By Theorem \ref{C4-1}, there are $9$ non-equivalent cases for $C$ which $r(C)$ is one of the relations $5$, $7$, $14$, $17$, $21$, $22$, $25$, $26$ and $29$  in Table \ref{tab-C4}. So, we have the followings:
\begin{enumerate}
\item[(5)]
$h_2^2 h_3 h_2^{-1} h_3=1$: Let $x=h_2^{-1}h_3$ and $y=h_2^{-1}$. So, $x^2=y^3$. Also, since $\alpha\beta=0$, we have $h_2^{-1}(\alpha_1\cdot 1 +\alpha_2h_2+\alpha_3h_3 )(\beta_1 g_1+\beta_2g_2+\cdots+\beta_ng_n)=0$. So,  $\alpha_2\cdot 1+\alpha_3x+\alpha_1y$ is a zero divisor with the support $\{1,x,y\}$.
\item[(7)]
$h_2^2 h_3^{-2} h_2=1$: Let $x=h_3$ and $y=h_2$. So $x^2=y^3$ and $\alpha_1\cdot 1 +\alpha_2y+\alpha_3x$ is a zero divisor with the support $\{1,x,y\}$.
\item[(14)]
$h_2 h_3 h_2^{-2} h_3=1$: Let $x=h_2h_3$ and $y=h_2$. So $x^2=y^3$ and $\alpha_1\cdot 1+\alpha_2y+\alpha_3y^{-1}x=\alpha_1\cdot 1 +\alpha_2h_2+\alpha_3h_3$ is a zero divisor with the support $\{1,y,y^{-1}x\}$.
\item[(17)] 
$h_2 (h_3 h_2^{-1})^2 h_3=1$: Let $x=h_2^{-1}$ and $y=h_3h_2^{-1}$. So $x^2=y^3$. Also, since $\alpha\beta=0$, we have $\alpha h_2^{-1}h_2\beta=0$. So, $\alpha_2\cdot 1 +\alpha_1x+\alpha_3y$ is a zero divisor with the support $\{1,x,y\}$.
\item[(21)]
$h_2 h_3^{-3} h_2=1$: By interchanging $h_2$ and $h_3$ in (7) and with the same discussion, the statement is true.
\item[(22)]
$h_2 h_3^{-2} h_2 h_3=1$: By interchanging $h_2$ and $h_3$ in (14) and with the same discussion, the statement is true.
\item[(25)]
$h_2 h_3^{-1} h_2 h_3^2=1$: By interchanging $h_2$ and $h_3$ in (5) and with the same discussion, the statement is true.
\item[(26)] 
$h_2 h_3^{-1} h_2 h_3 h_2^{-1} h_3=1$: Let $x=h_2h_3^{-1}h_2$ and $y=h_3^{-1}h_2$. So, $x^2=y^3$. Also, since $\alpha_1\cdot 1 +\alpha_2h_2+\alpha_3h_3=1+xy^{-1}+xy^{-2}=\alpha_2xy^{-1}+\alpha_3xy^{-2}+\alpha_1 x^2y^{-3}$, we have $x^{-1}(\alpha_2xy^{-1}+\alpha_3xy^{-2}+\alpha_1 x^2y^{-3})yy^{-1}\beta=0$. Therefore, $\alpha_2\cdot 1+\alpha_3y^{-1}+\alpha_1 xy^{-2}$ is also a zero divisor with the support of size $3$. Furthermore, $(\alpha_3\cdot 1+\alpha_2y+\alpha_1y^{-1}x)(y^{-3}\beta)=(y^{-1}(\alpha_2\cdot 1+\alpha_3y^{-1}+\alpha_1 xy^{-2})y^2)(y^{-3}\beta)=0$. Hence, $\alpha_3\cdot 1+\alpha_2y+\alpha_1y^{-1}x$ is a zero divisor with the support $\{1,y,y^{-1}x\}$.
\item[(29)]
$(h_2 h_3^{-1})^2 h_2 h_3=1$: By interchanging $h_2$ and $h_3$ in (17) and with the same discussion, the statement is true.
\end{enumerate}
\end{proof}

In the following, we discuss about the existence of two squares in $K_{\mathbb{F}}(\alpha,\beta)$. 
\begin{lem}\label{C4,C4}
Suppose that there exist two squares in $K_{\mathbb{F}}(\alpha,\beta)$. Then, such cycles are equivalent.
\end{lem}
\begin{proof}
By Theorem \ref{C4-1}, if there exist two squares in  $K_{\mathbb{F}}(\alpha,\beta)$, then these two cycles must be between $9$ non-equivalent cases with the relations $5$, $7$, $14$, $17$, $21$, $22$, $25$, $26$ or $29$  in Table \ref{tab-C4}. We may choose two relations similar to or different from each other.  When choosing two relations different from each other, there are $\binom{9}{2}= 36$ cases. Using GAP \cite{gap}, each group with two generators $h_2$ and $h_3$, and two relations which is between $36$ latter cases is finite and solvable, that is a contradiction with the assumptions. So by Theorem \ref{C4-1}, if there exist two squares in the graph $K_{\mathbb{F}}(\alpha,\beta)$, then such cycles must be equivalent.
\end{proof}

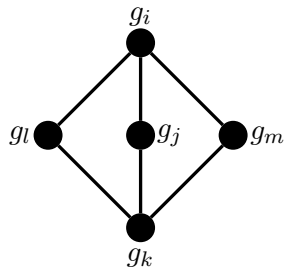
\begin{figure}[ht]
\centering
\begin{tikzpicture}[scale=0.7]
\node (1) [circle, minimum size=3pt, fill=black, line width=0.625pt, draw=black] at (75.0pt, -25.0pt)  {};
\node (2) [circle, minimum size=3pt, fill=black, line width=0.625pt, draw=black] at (125.0pt, -75.0pt)  {};
\node (3) [circle, minimum size=3pt, fill=black, line width=0.625pt, draw=black] at (25.0pt, -75.0pt)  {};
\node (4) [circle, minimum size=3pt, fill=black, line width=0.625pt, draw=black] at (75.0pt, -75.0pt)  {};
\node (5) [circle, minimum size=3pt, fill=black, line width=0.625pt, draw=black] at (75.0pt, -125.0pt)  {};
\draw [line width=1.25, color=black] (1) to  (2);
\draw [line width=1.25, color=black] (1) to  (4);
\draw [line width=1.25, color=black] (1) to  (3);
\draw [line width=1.25, color=black] (4) to  (5);
\draw [line width=1.25, color=black] (2) to  (5);
\draw [line width=1.25, color=black] (3) to  (5);
\node at (75.0pt, -10.625pt) {\textcolor{black}{$g_i$}};
\node at (143.875pt, -75.0pt) {\textcolor{black}{$g_m$}};
\node at (9.875pt, -75.0pt) {\textcolor{black}{$g_l$}};
\node at (90.875pt, -75.0pt) {\textcolor{black}{$g_j$}};
\node at (75.0pt, -141.375pt) {\textcolor{black}{$g_k$}};
\end{tikzpicture}
\caption{The complete bipartite graph $K_{2,3}$ in the Kaplansky graph over $\mathbb{F}$}\label{f-K2,3}
\end{figure}
\begin{thm}\label{K2,3}
The Kaplansky graph over $\mathbb{F}$ contains no subgraph isomorphic to the complete bipartite graph $K_{2,3}$.
\end{thm}
\begin{proof}
Suppose that $K_{\mathbb{F}}(\alpha,\beta)$ contains $K_{2,3}$ as a subgraph. So, it contains $2$ squares $C$ and $C'$ with two edges in common as Figure \ref{f-K2,3}, where $g_i,g_j,g_k,g_l,g_m\in supp(\beta)$. Let $T_C=[a_1,b_1,a_2,b_2,a_3,b_3,a_4,b_4]$ and $T_{C'}=[a_1,b_1,a_2,b_2,a_3',b_3',a_4',b_4']$ be the $8$-tuples of $C$ and $C'$, respectively, with the corresponding relations as follows:
\begin{equation}\label{e-7}
R(T_C): \left\{
\begin{array}{l}
a_1g_i=b_1g_j\\
a_2g_j=b_2g_k\\
a_3g_k=b_3g_l\\
a_4g_l=b_4g_i\\
\end{array} \right.
\qquad \qquad
R(T_{C'}): \left\{
\begin{array}{l}
a_1g_i=b_1g_j\\
a_2g_j=b_2g_k\\
a_3'g_k=b_3'g_m\\
a_4'g_m=b_4'g_i
\end{array} \right.
\end{equation}
where $a_s,b_s,a_t',b_t' \in supp(\alpha)$ for all $s\in \{1,2,3,4\}$ and $t\in \{3,4\}$.

By Remark \ref{r-K3-K3}, we have the following inequalities.
\begin{align}\label{e-7-2}
a_1\not= b_1\not= a_2 \not= b_2\not= a_3\not= b_3\not= a_4\not= b_4\not= a_1 \text{ and } b_2\not= a_3'\not= b_3'\not= a_4'\not= b_4'\not= a_1.
\end{align}

Since the graph with the set of vertices $\{g_i,g_m,g_k,g_l\}$ in $K_{2,3}$ is also an square, by Remark \ref{r-K3-K3} we have the following inequalities,
\begin{align}\label{e-5}
a_3 \neq a_3'.
\end{align}
\begin{align}\label{e-6}
b_4 \neq b_4'.
\end{align}

By Lemma \ref{C4,C4}, the cycles $C$ and $C '$ are equivalent.  So, $T_{C'}$ must be in $\mathcal{T}(C)$. In the following, we show that this leads to contradictions.
\begin{enumerate}
\item
$[a_1,b_1,a_2,b_2,a_3',b_3',a_4',b_4']=[a_1,b_1,a_2,b_2,a_3,b_3,a_4,b_4]$: So $a_3=a_3'$, that is a contradiction with the relation \ref{e-5}.
\item
$[a_1,b_1,a_2,b_2,a_3',b_3',a_4',b_4']=[a_4,b_4,a_1,b_1,a_2,b_2,a_3,b_3]$: Therefore, $T_C=[a_1,b_1,a_1,b_1,a_3,b_3,a_1,b_1]$ and  $T_{C'}=[a_1,b_1,a_1,b_1,a_1,b_1,a_3,b_3]$. By the relations \ref{e-5} and \ref{e-6} we have $a_3\neq a_1$ and $b_3\neq b_1$. Also, in such $8$-tuples we have  $a_3\neq b_1$, $b_3\neq a_1$ and $a_1\neq b_1$. Therefore, $a_3=b_3$ since $a_3,b_3 \in supp(\alpha)$, that is a contradiction.
\item
$[a_1,b_1,a_2,b_2,a_3',b_3',a_4',b_4']=[a_3,b_3,a_4,b_4,a_1,b_1,a_2,b_2]$: So $a_3=a_1=a_3'$, that is a contradiction with the relation \ref{e-5}.
\item
$[a_1,b_1,a_2,b_2,a_3',b_3',a_4',b_4']=[a_2,b_2,a_3,b_3,a_4,b_4,a_1,b_1]$: Therefore, $T_C=[a_1,b_1,a_1,b_1,a_1,b_1,a_4,b_4]$ and  $T_{C'}=[a_1,b_1,a_1,b_1,a_4,b_4,a_1,b_1]$. By the relations \ref{e-5} and \ref{e-6} we have $a_4\neq a_1$ and $b_4\neq b_1$. Also, in such $8$-tuples we have  $a_4\neq b_1$, $b_4\neq a_1$ and $a_1\neq b_1$. Therefore, $a_4=b_4$ since $a_4,b_4 \in supp(\alpha)$, that is a contradiction.
\item
$[a_1,b_1,a_2,b_2,a_3',b_3',a_4',b_4']=[b_1,a_1,b_4,a_4,b_3,a_3,b_2,a_2]$: So $a_1=b_1$, that is a contradiction with the relations \ref{e-7-2}.
\item
$[a_1,b_1,a_2,b_2,a_3',b_3',a_4',b_4']=[b_2,a_2,b_1,a_1,b_4,a_4,b_3,a_3]$: So $b_1=a_2$, that is a contradiction with the relations \ref{e-7-2}.
\item
$[a_1,b_1,a_2,b_2,a_3',b_3',a_4',b_4']=[b_3,a_3,b_2,a_2,b_1,a_1,b_4,a_4]$: So $a_2=b_2$, that is a contradiction with the relations \ref{e-7-2}.
\item
$[a_1,b_1,a_2,b_2,a_3',b_3',a_4',b_4']=[b_4,a_4,b_3,a_3,b_2,a_2,b_1,a_1]$: So $b_2=a_3$, that is a contradiction with the relations \ref{e-7-2}.
\end{enumerate}
So, the graph $K_{\mathbb{F}}(\alpha,\beta)$ contains no subgraph isomorphic to the graph $K_{2,3}$.
\end{proof}

\section{\bf Possible zero divisors with supports of size $3$ in $\mathbb{F}[G]$}\label{S1-1}
If $A=\{1,2,3\}\times \{1,2,\ldots,n\}$, then for all $(i,j)\in A$ there must be an $(i',j')\in A$ such that $i\neq i'$, $j\neq j'$ and  $h_ig_j=h_{i'}g_{j'}$ because $\alpha\beta=(\alpha_1 h_1+\alpha_2h_2+\alpha_3h_3)(\beta_1 g_1+\beta_2g_2+\cdots+\beta_ng_n)=0$. Also, $n\geq 3$ \cite{pascal}. Firstly in this section we show that $n>3$. Then, we examine some small positive integers greater that $3$ as the possible values of $n$ and show that $n$ must be at least $10$. 

\subsection{The support of $\beta$ cannot be of size $3$}
Let $|supp(\beta)|=3$. Since $\alpha \beta =0$, we must have  $(\alpha_1\beta_1h_1g_1+\alpha_1\beta_2h_1g_2+\alpha_1\beta_3h_1g_3)+
(\alpha_2\beta_1h_2g_1+\alpha_2\beta_2h_2g_2+\alpha_2\beta_3h_2g_3)+
(\alpha_3\beta_1h_3g_1+\alpha_3\beta_2h_3g_2+\alpha_3\beta_3h_3g_3)=0$. Therefore, $h_1g_1=h_ig_j$ for some $(i,j)\in A$ where $i\not=1$ and $j\not=1$. Also, $h_2g_1=h_{i'}g_{j'}$ for some $(i',j')\in A$ where $i'\not=2$ and $j'\not=1$. Furthermore, $h_3g_1=h_{i''}g_{j''}$ for some $(i'',j'')\in A$ where $i''\not=3$ and $j''\not=1$. Note that $j'\not=j$, $j''\not=j$ and $j''\not=j'$ because the Lemma \ref{size of S} states that the set $S=\{h^{-1} h' \;|\; h,h'\in supp(\alpha), h\not=h' \}$ has size $6$. Hence, $g_{j''}\in supp(\beta)$ and $g_{j''}\notin \{g_1,g_2,g_3\}$, a contradiction. So, $|supp(\beta)|$ must be at least $4$.

\subsection{The support of $\beta$ must be of size greater than or equal to $10$}
Abelian groups satisfy Conjecture \ref{conj-zero}. So, $G$ must be a nonabelian torsion-free group. The following theorem is obtained in \cite{ham}. 

\begin{thm}[Corollary 11 of \cite{ham}]\label{hamidoune}
If $C$ is a finite generating subset of a nonabelian torsion-free group $G$ such that $1\in C$ and $|C|\geq 4$, then $|BC|\geq |B|+|C|+1$ for all $B\subset G$ with $|B|\geq 3$.
\end{thm}
Without loss of generality we may assume that $G$ is generated by $supp(\alpha)\cup supp(\beta)$, since otherwise we replace $G$ by the subgroup generated by this set. Also by Lemma \ref{supp}, $\langle h_i^{-1} supp(\alpha) \rangle=\langle  supp(\beta) g_j^{-1} \rangle$ for all $i\in \{1,2,3\}$ and $j\in\{1,2,\ldots,n\}$. Therefore by Theorem \ref{hamidoune}, $|supp(\alpha) supp(\beta)|\geq |supp(\alpha)|+|supp(\beta)|+1$. 
Hence, $3n\geq |supp(\alpha) supp(\beta)|\geq 4+n$.
\begin{thm}[Proposition 4.12 of \cite{dyk}]\label{dyk-rem}
There exist no $\gamma , \delta \in \mathbb{F}_2[G]$ such that $\gamma \delta=1$, where $|supp(\gamma)|=3$ and $|supp(\delta)|\geq 13$ is an odd integer. 
\end{thm}
\begin{enumerate}
\item
Let $n=4$. Then with the discussion above $12\geq |supp(\alpha) supp(\beta)|\geq 8$. Since $12-8=4$, there is an $(i,j)\in A$ such that $h_ig_j\not=h_{i'}g_{j'}$ for all $(i',j')\in A$ where $i\neq i'$ and $j\neq j'$, a contradiction with $\alpha \beta =0$. So, $|supp(\beta)|$ must be at least $5$.
\item
Let $n=5$. Then with the discussion above $15\geq |supp(\alpha) supp(\beta)|\geq 9$. Since $15-9=6$, there is an $(i,j)\in A$ such that $h_ig_j\not=h_{i'}g_{j'}$ for all $(i',j')\in A$ where $i\neq i'$ and $j\neq j'$, a contradiction with $\alpha \beta =0$. So, $|supp(\beta)|$ must be at least $6$.
\item
Let $n=6$. Then with the discussion above $18\geq |supp(\alpha) supp(\beta)|\geq 10$. Since $18-10=8$, there is an $(i,j)\in A$ such that $h_ig_j\not=h_{i'}g_{j'}$ for all $(i',j')\in A$ where $i\neq i'$ and $j\neq j'$, a contradiction with $\alpha \beta =0$. So, $|supp(\beta)|$ must be at least $7$.
\item
Let $n=7$. Then with the discussion above $21\geq |supp(\alpha) supp(\beta)|\geq 11$. Since $21-11=10$, there is an $(i,j)\in A$ such that $h_ig_j\not=h_{i'}g_{j'}$ for all $(i',j')\in A$ where $i\neq i'$ and $j\neq j'$, a contradiction with $\alpha \beta =0$. So, $|supp(\beta)|$ must be at least $8$.
\item
Let $n=8$. Then with the discussion above $24\geq |supp(\alpha) supp(\beta)|\geq 12$. Let $|supp(\alpha) supp(\beta)|> 12$. Then $|supp(\alpha) supp(\beta)|\geq 13$. Since $24-13=11$, there is an $(i,j)\in A$ such that $h_ig_j\not=h_{i'}g_{j'}$ for all $(i',j')\in A$ where $i\neq i'$ and $j\neq j'$, a contradiction with $\alpha \beta =0$. So, $|supp(\alpha) supp(\beta)|=12$ and because $\alpha \beta =0$, there is a partition $\pi$ of $A$ with all sets containing two elements, such that if $(i, j)$ and $(i',j')$ belong to the same set of $\pi$, then $h_ig_j=h_{i'}g_{j'}$. Let $\alpha'=\sum_{a\in supp(\alpha)}{a}$ and $\beta'=\sum_{b\in supp(\beta)}{b}$. So, $\alpha', \beta' \in \mathbb{F}_2[G]$, $|supp(\alpha')|=3$ and $|supp(\beta')|=8$ and with the above discussion we have $\alpha'\beta'=0$, that is a contradiction (see below, Corollary \ref{maincor}). Therefore, $|supp(\alpha) supp(\beta)|\not=12$ and so $|supp(\beta)|$ must be at least $9$.
\item
Let $n=9$. Then with the discussion above $27\geq |supp(\alpha) supp(\beta)|\geq 13$. Let $|supp(\alpha) supp(\beta)|> 13$. Then $|supp(\alpha) supp(\beta)|\geq 14$. Since $27-14=13$, there is an $(i,j)\in A$ such that $h_ig_j\not=h_{i'}g_{j'}$ for all $(i',j')\in A$ where $i\neq i'$ and $j\neq j'$, a contradiction with $\alpha \beta =0$. So, $|supp(\alpha) supp(\beta)|=13$ and because $\alpha \beta =0$, there is a partition $\pi$ of $A$ with one set of size $3$ and all other sets containing two elements, such that if $(i, j)$ and $(i',j')$ belong to the same set of $\pi$, then $h_ig_j=h_{i'}g_{j'}$. With the discussion above, $\left(\sum_{a\in supp(\alpha)}{a}\right)\left(\sum_{b\in supp(\beta)}{b}\right)x^{-1}=1$ where $x=h_ig_j$ for some $(i,j)$ belongs to the set  of size $3$ in $\pi$. Hence, there are $\gamma , \delta \in \mathbb{F}_2[G]$ such that $\gamma \delta=1$, where $\gamma=\sum_{a\in supp(\alpha)}{a}$, $\delta=\sum_{b\in supp(\beta)}{bx^{-1}}$, $|supp(\gamma)|=3$ and $|supp(\delta)|=|supp(\beta)|=9$, that is a contradiction with Theorem \ref{dyk-rem}. Therefore, $|supp(\alpha) supp(\beta)|\not=13$ and so $|supp(\beta)|$ must be at least $10$.
\end{enumerate}
With the discussion above, we have the following theorem.
\begin{thm}\label{main-F}
Let $\alpha$ and $\beta$ be non-zero elements of the group algebra of any torsion-free group over an arbitrary field. If $|supp(\alpha)|=3$ and $\alpha \beta=0$ then $|supp(\beta)|\geq 10$. 
\end{thm}

\begin{prop}\label{supp2}
If $\mathbb{F}[G]$ has no non-zero element $\alpha$ with $|supp(\alpha)|\leq k$ such that $\alpha^2=0$, then there exist no non-zero elements $\alpha_1,\alpha_2 \in \mathbb{F}[G]$ such that $\alpha_1\alpha_2=0$ and $|supp(\alpha_1)||supp(\alpha_2)|\leq k$.
\end{prop}
\begin{proof}
Suppose, for a contradiction, that $\alpha_1,\alpha_2 \in \mathbb{F}[G]\setminus \{0\}$ such that $\alpha_1\alpha_2=0$ and $|supp(\alpha_1)||supp(\alpha_2)|\leq k$. We may assume that $1\in supp(\alpha_1) \cap supp(\alpha_2)$, since $(a^{-1} \alpha_1 ) (\alpha_2 b^{-1})=0$ for any $a\in supp(\alpha_1)$ and $b\in supp(\alpha_2)$. 

Suppose, for a contradiction, that $\alpha_2 x \alpha_1=0$ for all $x\in G$. Then it follows from 
\cite[Lemma 1.3, p. 3]{pass2} that $\theta(\alpha_2) \theta(\alpha_1)=0$, where $\theta$ is the projection $\theta:\mathbb{F}[G]\rightarrow \mathbb{F}[\Delta]$ given by 
$\beta= \sum_{x\in G} f_x x \mapsto \theta(\beta)=\sum_{x\in \Delta} f_x x$, where $\Delta$ is the subgroup of all
elements of $G$ having a finite number of conjugates in $G$ (see \cite[p. 3]{pass2}). Now it follows from \cite[Lemma 2.2, p. 5]{pass2} and \cite[Lemma 2.4, p. 6]{pass2} that $\theta(\alpha_1)=0$ or $\theta(\alpha_2)=0$, which are both contradiction since  $1\in supp(\alpha_1) \cap supp(\alpha_2)$.  
Therefore, there exists an element $x\in G$ such that $\beta=\alpha_2 x \alpha_1\not=0$. Now    
$$\beta^2=(\alpha_2 x \alpha_1)^2=\alpha_2 x \alpha_1 \alpha_2 x \alpha_1=0$$
and $$|supp(\beta)|\leq |supp(\alpha_2)||supp(x \alpha_1)|=|supp(\alpha_2)||supp(\alpha_1)|\leq k,$$
which is a contradiction. This completes the proof. 
\end{proof}

In the next three sections, we discuss about Kaplansky graphs over  $\mathbb{F}_2$ and give some forbidden subgraphs for such graphs.
\section{\bf Kaplansky graphs over $\mathbb{F}_2$ and some of their subgraphs containing an square}\label{S-C4}
Throughout the rest of this paper, except in Section \ref{S-unit}, let $\mathbb{F}$ be the finite field $\mathbb{F}_2$ and  $\alpha=h_1+h_2+h_3 \in \mathbb{F}_2[G]$ such that $\vert supp(\alpha)\vert=3$. Suppose further that  $\alpha \beta =0$ for some non-zero $\beta \in \mathbb{F}_2[G]$ and assume that  $n:=|supp(\beta)|$ is minimum with respect to the latter property.  Let   $\beta=g_1+g_2+\cdots+g_n$. If there is no ambiguity, we denote the Kaplansky graph of $(\alpha, \beta)$ over $\mathbb{F}_2$ by $K(\alpha,\beta)$ and simply call it Kaplansky graph.

\begin{lem}\label{Cay}
The Kaplansky graph $K(\alpha,\beta)$ is isomorphic to the induced subgraph on the set $supp(\beta)$ of the Cayley graph $Cay(G,S)$, where 
$S=\{h^{-1} h' \;|\; h,h'\in supp(\alpha), h\not=h' \}$.
\end{lem}
\begin{proof}
Let  $M$ be the matched rectangle corresponding to $(\alpha,\beta)$ (see \cite[Definition 4.1]{pascal}). The vertex set of $K(M)=K(\alpha,\beta)$ is the columns of $M$ which are labelled by the elements of $supp(\beta)$ and two distinct columns $c$ and $c'$ are adjacent whenever their labels $g$ and $g' \in supp(\beta)$ respectively, satisfying $hg=h'g'$ for some  $h,h'\in supp(\alpha)$; Or equivalently the columns $c$ and $c'$ are adjacent whenever  $gg'^{-1} \in S$. Hence the map $\psi$ from the columns of $M$ to $supp(\beta)$ which sends each column to its label is a graph isomorphism from $K(M)=K(\alpha,\beta)$ to the induced subgraph on the set $supp(\beta)$ of the Cayley graph $Cay(G,S)$.   
\end{proof}

\begin{rem}\label{r-F2}
{\rm It follows from Lemma \ref{Cay} and Theorem \ref{thm-graph} that the induced subgraph of $Cay(G,S)$ on $supp(\beta)$ is a cubic graph having no subgraph isomorphic to a triangle. Also $n=|supp(\beta)|$ which is the number of vertices of $K(\alpha,\beta)$ is always an even number, since the number of vertices of any cubic graph is even.}
\end{rem}

In the rest of this section, we consider Kaplansky graphs containing a subgraph isomorphic to an square. 

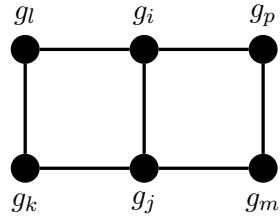
\begin{figure}[ht]
\centering
\begin{tikzpicture}[scale=0.9]
\node (2) [circle, minimum size=3pt, fill=black, line width=0.625pt, draw=black] at (100.0pt, -75.0pt)  {};
\node (4) [circle, minimum size=3pt, fill=black, line width=0.625pt, draw=black] at (100.0pt, -25.0pt)  {};
\node (1) [circle, minimum size=3pt, fill=black, line width=0.625pt, draw=black] at (50.0pt, -25.0pt)  {};
\node (3) [circle, minimum size=3pt, fill=black, line width=0.625pt, draw=black] at (50.0pt, -75.0pt)  {};
\node (5) [circle, minimum size=3pt, fill=black, line width=0.625pt, draw=black] at (150.0pt, -25.0pt)  {};
\node (6) [circle, minimum size=3pt, fill=black, line width=0.625pt, draw=black] at (150.0pt, -75.0pt)  {};
\draw [line width=1.25, color=black] (1) to  (3);
\draw [line width=1.25, color=black] (3) to  (2);
\draw [line width=1.25, color=black] (1) to  (4);
\draw [line width=1.25, color=black] (4) to  (5);
\draw [line width=1.25, color=black] (6) to  (5);
\draw [line width=1.25, color=black] (2) to  (6);
\draw [line width=1.25, color=black] (2) to  (4);
\node at (100.0pt, -89.375pt) {\textcolor{black}{$g_j$}};
\node at (100.0pt, -10.625pt) {\textcolor{black}{$g_i$}};
\node at (50.0pt, -10.625pt) {\textcolor{black}{$g_l$}};
\node at (50.0pt, -89.375pt) {\textcolor{black}{$g_k$}};
\node at (150.0pt, -10.625pt) {\textcolor{black}{$g_p$}};
\node at (150.0pt, -89.375pt) {\textcolor{black}{$g_m$}};
\end{tikzpicture}
\caption{Two squares with one common edge in $K(\alpha,\beta)$}\label{f-C4-C4}
\end{figure}

\begin{thm}\label{C4-C4}
Suppose that $K(\alpha,\beta)$ contains two squares with exactly one edge in common. Then exactly one of the relations $14$, $22$ or $26$ of Table \ref{tab-C4} will be satisfied in $G$.
\end{thm}
\begin{proof}
Suppose that the graph $K(\alpha,\beta)$ contains two squares $C$ and $C'$ with exactly one common edge as Figure \ref{f-C4-C4}, where $g_i,g_j,g_k,g_l,g_m,g_p\in supp(\beta)$. Let $T_C=[a_1,b_1,a_2,b_2,a_3,b_3,a_4,b_4]$ and $T_{C'}=[a_1,b_1,a_2',b_2',a_3',b_3',a_4',b_4']$ be the $8$-tuples of $C$ and $C'$, respectively, with the corresponding relations as follows:
\begin{equation}\label{e-8}
R(T_C): \left\{
\begin{array}{l}
a_1g_i=b_1g_j\\
a_2g_j=b_2g_k\\
a_3g_k=b_3g_l\\
a_4g_l=b_4g_i\\
\end{array} \right.
\qquad \qquad
R(T_{C'}): \left\{
\begin{array}{l}
a_1g_i=b_1g_j\\
a_2'g_j=b_2'g_m\\
a_3'g_m=b_3'g_p\\
a_4'g_p=b_4'g_i\\
\end{array} \right.
\end{equation}
where $a_z,b_z, a_t',b_t' \in supp(\alpha)$ for all $z\in \{1,2,3,4\}$ and $t\in \{2,3,4\}$.

By Remark \ref{r-K3-K3}, $a_1\neq b_1\neq \cdots \neq a_4\neq b_4\neq a_1$ and $b_1\neq  a_2'\neq \cdots \neq a_4' \neq b_4'\neq a_1$. We want to prove that $a_2 \neq a_2'$ and $b_4 \neq b_4'$. Suppose that  $a_2=a_2'$.  So we have $a_2g_j=a_2'g_j$. Since $a_2g_j=b_2g_k$ and $a_2'g_j=b_2'g_m$, we have $b_2g_k=a_2'g_j=b_2'g_m$. Since $(1+h_2+h_3)(g_1+g_2+\cdots+g_n)=0$ in $\mathbb{F}_2[G]$, there must exist $g_a \in \{g_1,g_2,\ldots,g_n\}\setminus \{g_j,g_k,g_m\}$ and $h_a \in supp(\alpha)$ such that $b_2g_k=a_2'g_j=b_2'g_m=h_ag_a$. Since the set $supp(\alpha)$ has size $3$ and $ \{a_2',b_2,b_2'\}=supp(\alpha)$, we have $h_a \in \{a_2',a_3,a_3'\}$. Without loss of generality we may assume that $h_a=a_2'$. So $g_j=g_a$, that is a contradiction. Hence, 
\begin{align}\label{e-9}
a_2 \neq a_2'.
\end{align}
Also with the same discussion such as above, we have
\begin{align}\label{e-10}
b_4 \neq b_4'.
\end{align}
By Lemma \ref{C4,C4}, the cycles $C$ and $C '$ are equivalent. So, $T_{C'}$ must be in $\mathcal{T}(C)$. In the following, we discuss about all the possible cases in details.
\begin{enumerate}
\item
$[a_1,b_1,a_2',b_2',a_3',b_3',a_4',b_4']=[a_1,b_1,a_2,b_2,a_3,b_3,a_4,b_4]$: So $a_2= a_2'$, that is a contradiction with the relation \ref{e-9}.
\item
$[a_1,b_1,a_2',b_2',a_3',b_3',a_4',b_4']=[a_4,b_4,a_1,b_1,a_2,b_2,a_3,b_3]$: So, $T_{C}=[a_1,b_1,a_2,b_2,a_3,b_3,a_1,b_1]$ and  $T_{C'}=[a_1,b_1,a_1,b_1,a_2,b_2,a_3,b_3]$. By the relations  \ref{e-8}, \ref{e-9} and \ref{e-10}, we have $a_2\neq a_1$, $a_2\neq b_1$, $b_3\neq a_1$, $b_3\neq b_1$ and $a_1\neq b_1$. Therefore, $a_2=b_3$ because $a_2,b_3 \in supp(\alpha)$. So,  $T_{C}=[a_1,b_1,a_2,b_2,a_3,a_2,a_1,b_1]$ and  $T_{C'}=[a_1,b_1,a_1,b_1,a_2,b_2,a_3,a_2]$. Since by the relation \ref{e-8}, $\{a_1,b_1,a_2\}=supp(\alpha)$, $b_2\neq a_2$, $a_3\neq a_2$ and $b_2\neq a_3$, there are just two cases for choosing $b_2$ and $a_3$. In the following, we show that both of them lead to contradictions.
\begin{enumerate}
\item[i)]
$b_2=a_1$ and $a_3=b_1$: So, $T_C=[a_1,b_1,a_2,a_1,b_1,a_2,a_1,b_1]$ and the relation of such cycle is $a_1^{-1}b_1a_2^{-1}a_1b_1^{-1}a_2a_1^{-1}b_1=1$. Since $\{a_1,b_1,a_2\}=supp(\alpha)$, just one of the following cases may be happened:
\begin{enumerate}
\item[a)]
$a_1=1$: $a_2^{-1}b_1a_2=b_1^{\ 2}$ and so $\langle h_2,h_3\rangle=BS(1,2)$, that is a contradiction.
\item[b)]
$b_1=1$: $a_2^{-1}a_1a_2=a_1^{\ 2}$ and so $\langle h_2,h_3\rangle=BS(1,2)$, that is a contradiction.
\item[c)]
$a_2=1$: $a_1^{-1}b_1a_1b_1^{-1}a_1^{-1}b_1=1$. If $x=a_1^{-1}b_1$ and $y=b_1^{-1}$, then $y^{-1}xy=x^2$ and $\langle h_2,h_3\rangle=\langle x,y\rangle=BS(1,2)$, that is a contradiction.
\end{enumerate}
\end{enumerate}
\begin{enumerate}
\item[ii)]
$b_2=b_1$ and $a_3=a_1$: So, $T_C=[a_1,b_1,a_2,b_1,a_1,a_2,a_1,b_1]$ and the relation of such cycle is $a_1^{-1}b_1a_2^{-1}b_1a_1^{-1}a_2a_1^{-1}b_1=1$. With the same discussion as in item (i), exactly one of the following cases may be happened:
\begin{enumerate}
\item[a)]
$a_1=1$: $a_2^{-1}b_1a_2=b_1^{-2}$ and so $\langle h_2,h_3\rangle=BS(1,-2)$, that is a contradiction.
\item[b)]
$b_1=1$: $a_2^{-1}a_1a_2=a_1^{-2}$ and so $\langle h_2,h_3\rangle=BS(1,-2)$, that is a contradiction.
\item[c)]
$a_2=1$: $a_1^{-1}b_1^{\ 2}a_1^{-2}b_1=1$. If $x=b_1^{-1}a_1$ and $y=b_1^{-1}$, then $y^{-1}xy=x^{-2}$ and $\langle h_2,h_3\rangle=\langle x,y\rangle=BS(1,-2)$, that is a contradiction.
\end{enumerate}
\end{enumerate}
Hence, $[a_1,b_1,a_2',b_2',a_3',b_3',a_4',b_4']\neq [a_4,b_4,a_1,b_1,a_2,b_2,a_3,b_3]$.
\item
$[a_1,b_1,a_2',b_2',a_3',b_3',a_4',b_4']=[a_3,b_3,a_4,b_4,a_1,b_1,a_2,b_2]$: So, $T_{C}=[a_1,b_1,a_2,b_2,a_1,b_1,a_4,b_4]$ and  $T_{C'}=[a_1,b_1,a_4,b_4,a_1,b_1,a_2,b_2]$. By the relation \ref{e-8}, we have $a_2 \neq b_1$. In the following, we show that $a_2\neq a_1$. \\
Suppose that $a_2=a_1$. So, $T_{C}=[a_1,b_1,a_1, b_2,a_1,b_1,a_4,b_4]$ and  $T_{C'}=[a_1,b_1,a_4,b_4,a_1,b_1,a_1,b_2]$. By the relations \ref{e-8}, \ref{e-9} and \ref{e-10}, we have $b_1\neq a_1$, $b_1\neq a_4$, $b_4\neq a_1$, $b_4\neq a_4$ and $a_1\neq a_4$. Therefore, $b_1=b_4$ since $b_1,b_4 \in supp(\alpha)$. Now by the relations \ref{e-8}, \ref{e-9} and \ref{e-10}, we have $b_2\neq a_1$, $b_2\neq b_1$, $a_4\neq a_1$, $a_4\neq b_1$ and $a_1\neq b_1$. Hence, $b_2=a_4$ because $b_2,a_4 \in supp(\alpha)$. So $T_C=[a_1,b_1,a_1,b_2,a_1,b_1,b_2,b_1]$ and the relation of $C$ is $a_1^{-1}b_1a_1^{-1}b_2a_1^{-1}b_1b_2^{-1}b_1=1$. Since $\{a_1,b_1,b_2\}=supp(\alpha)$, just one of the following cases may be happened:
\begin{enumerate}
\item[a)]
$a_1=1$: $b_2^{-1}b_1^{-2}b_2=b_1$ and so $\langle h_2,h_3\rangle=BS(-2,1)$, that is a contradiction.
\item[b)]
$b_1=1$: $b_2^{-1}a_1^{-2}b_2=a_1$ and so $\langle h_2,h_3\rangle=BS(-2,1)$, that is a contradiction.
\item[c)]
$b_2=1$: $a_1^{-1}b_1a_1^{-2}b_1^{\ 2}=1$. If $x=b_1a_1^{-1}$ and $y=a_1^{-1}$, then $y^{-1}x^{-2}y=x$ and $\langle h_2,h_3\rangle=\langle x,y\rangle=BS(-2,1)$, that is a contradiction.
\end{enumerate}
Hence, $a_2\neq a_1$. Also, by the relations \ref{e-8}, \ref{e-9} and \ref{e-10} and by $T_{C}$ and $T_{C'}$, we have $a_1\neq b_1$, $a_4\neq a_2$, $a_4\neq b_1$ and $b_1\neq a_2$. Therefore, $a_1=a_4$ because $a_1,a_4 \in supp(\alpha)$. Now by the relations \ref{e-8}, \ref{e-9} and \ref{e-10}, we have $a_2\neq a_1$, $a_2\neq b_2$, $b_4\neq a_1$, $b_4\neq b_2$ and $a_1\neq b_2$. So, $a_2=b_4$ because $a_2,b_4 \in supp(\alpha)$. Also, $b_2\neq a_1$, $b_2\neq a_2$, $b_1\neq a_1$, $b_1\neq a_2$ and $a_1\neq a_2$. Therefore, $b_1=b_2$ because $b_1,b_2 \in supp(\alpha)$. Hence, $T_C=[a_1,b_1,a_2,b_1,a_1,b_1,a_1,a_2]$ and the relation of such cycle is $a_1^{-1}b_1a_2^{-1}b_1a_1^{-1}b_1a_1^{-1}a_2=1$. Since $\{a_1,b_1,a_2\}=supp(\alpha)$,  exactly one of the following cases may be happened:
\begin{enumerate}
\item[a)]
$a_1=1$: $a_2^{-1}b_1^{-2}a_2=b_1$ and so $\langle h_2,h_3\rangle=BS(-2,1)$, that is a contradiction.
\item[b)]
$b_1=1$: $a_2^{-1}a_1^{-2}a_2=a_1$ and so $\langle h_2,h_3\rangle=BS(-2,1)$, that is a contradiction.
\item[c)]
$a_2=1$: $a_1^{-1}b_1^{\ 2}a_1^{-1}b_1a_1^{-1}=1$. If $x=a_1b_1^{-1}$ and $y=a_1^{-1}$, then $y^{-1}x^{-2}y=x$ and $\langle h_2,h_3\rangle=\langle x,y\rangle=BS(-2,1)$, that is a contradiction.
\end{enumerate}
Hence, $[a_1,b_1,a_2',b_2',a_3',b_3',a_4',b_4']\neq [a_3,b_3,a_4,b_4,a_1,b_1,a_2,b_2]$.
\item
$[a_1,b_1,a_2',b_2',a_3',b_3',a_4',b_4']=[a_2,b_2,a_3,b_3,a_4,b_4,a_1,b_1]$: So, $T_{C}=[a_1,b_1,a_1,b_1,a_3,b_3,a_4,b_4]$ and  $T_{C'}=[a_1,b_1,a_3,b_3,a_4,b_4,a_1,b_1]$. By the relations  \ref{e-8}, \ref{e-9} and \ref{e-10}, we have $a_3\neq a_1$, $a_3\neq b_1$, $b_4\neq a_1$, $b_4\neq b_1$ and $a_1\neq b_1$. Therefore, $a_3=b_4$ because $a_3,b_4 \in supp(\alpha)$. Now because $\{a_1,b_1,a_3\}=supp(\alpha)$, $b_3\neq a_3$, $a_4\neq a_3$ and $b_3\neq a_4$, there are two cases for choosing $b_3$ and $a_4$. In the following, we show that both of them lead to contradictions.
\begin{enumerate}
\item[i)]
$b_3=a_1$ and $a_4=b_1$: So, $T_C=[a_1,b_1,a_1,b_1,a_3,a_1,b_1,a_3]$ and the relation of such cycle is $a_1^{-1}b_1a_1^{-1}b_1a_3^{-1}a_1b_1^{-1}a_3=1$. Since $\{a_1,b_1,a_3\}=supp(\alpha)$, exactly one of the following cases may be happened:
\begin{enumerate}
\item[a)]
$a_1=1$: $a_3^{-1}b_1a_3=b_1^{\ 2}$ and so $\langle h_2,h_3\rangle=BS(1,2)$, that is a contradiction.
\item[b)]
$b_1=1$: $a_3^{-1}a_1a_3=a_1^{\ 2}$ and so $\langle h_2,h_3\rangle=BS(1,2)$, that is a contradiction.
\item[c)]
$a_3=1$: $a_1^{-1}b_1a_1b_1^{-1}a_1^{-1}b_1=1$. If $x=a_1^{-1}b_1$ and $y=b_1^{-1}$, then $y^{-1}xy=x^2$ and $\langle h_2,h_3\rangle=\langle x,y\rangle=BS(1,2)$, that is a contradiction.
\end{enumerate}
\end{enumerate}
\begin{enumerate}
\item[ii)]
$b_3=b_1$ and $a_4=a_1$: So, $T_C=[a_1,b_1,a_1,b_1,a_3,b_1,a_1,a_3]$ and relation of such cycle is $a_1^{-1}b_1a_1^{-1}b_1a_3^{-1}b_1a_1^{-1}a_3=1$. With the same discussion as in item (i), exactly one of the following cases may be happened:
\begin{enumerate}
\item[a)]
$a_1=1$: $a_3^{-1}b_1a_3=b_1^{-2}$ and so $\langle h_2,h_3\rangle=BS(1,-2)$, that is a contradiction.
\item[b)]
$b_1=1$: $a_3^{-1}a_1a_3=a_1^{-2}$ and so $\langle h_2,h_3\rangle=BS(1,-2)$, that is a contradiction.
\item[c)]
$a_3=1$: $a_1^{-1}b_1^{\ 2}a_1^{-2}b_1=1$. If $x=b_1^{-1}a_1$ and $y=b_1^{-1}$, then $y^{-1}xy=x^{-2}$ and $\langle h_2,h_3\rangle=\langle x,y\rangle=BS(1,-2)$, that is a contradiction.
\end{enumerate}
\end{enumerate}
Hence, $[a_1,b_1,a_2',b_2',a_3',b_3',a_4',b_4']\neq [a_2,b_2,a_3,b_3,a_4,b_4,a_1,b_1]$.
\item
$[a_1,b_1,a_2',b_2',a_3',b_3',a_4',b_4']=[b_1,a_1,b_4,a_4,b_3,a_3,b_2,a_2]$: So $a_1 =b_1$, that is a contradiction with the relation \ref{e-8}.
\item
$[a_1,b_1,a_2',b_2',a_3',b_3',a_4',b_4']=[b_4,a_4,b_3,a_3,b_2,a_2,b_1,a_1]$: So $a_1=b_4$, that is a contradiction with the relation \ref{e-8}.
\item
$[a_1,b_1,a_2',b_2',a_3',b_3',a_4',b_4']=[b_3,a_3,b_2,a_2,b_1,a_1,b_4,a_4]$: So, $T_{C}=[a_1,b_1,a_2,b_2,b_1,a_1,a_4,b_4]$ and  $T_{C'}=[a_1,b_1,b_2,a_2,b_1,a_1,b_4,a_4]$. Since by the relation \ref{e-8}, $a_1 \neq b_1$ and $a_2 \neq b_1$, there are two cases, namely $a_2=a_1$ and $a_2\neq a_1$.
\begin{enumerate}
\item[A)] $a_2=a_1$: So, $T_{C}=[a_1,b_1,a_1,b_2,b_1,a_1, a_4,b_4]$ and  $T_{C'}=[a_1,b_1,b_2,a_1,b_1,a_1,b_4,a_4]$. By the relation \ref{e-8}, we have $b_2\neq a_1$, $b_2\neq b_1$ and $a_1\neq b_1$. So, $\{a_1,b_1,b_2\}=supp(\alpha)$. Also by this relation, $a_4\neq a_1$, $b_4\neq a_1$ and $a_4\neq b_4$. Therefore, there are two cases for choosing $a_4$ and $b_4$ which we show that one of them leads to a contradiction with our assumptions.\\
Suppose that $a_4=b_2$ and $b_4=b_1$. So, $T_C=[a_1,b_1,a_1,b_2,b_1,a_1,b_2,b_1]$ and the relation of such cycle is $a_1^{-1}b_1a_1^{-1}b_2b_1^{-1}a_1b_2^{-1}b_1=1$. Since $\{a_1,b_1,b_2\}=supp(\alpha)$, just one of the following cases may be happened:
\begin{enumerate}
\item[a)]
$a_1=1$: $b_2^{-1}b_1^{\ 2}b_2=b_1$ and so $\langle h_2,h_3\rangle=BS(2,1)$, that is a contradiction.
\item[b)]
$b_1=1$: $b_2^{-1}a_1^{\ 2}b_2=a_1$ and so $\langle h_2,h_3\rangle=BS(2,1)$, that is a contradiction.
\item[c)]
$b_2=1$: $a_1^{-1}b_1a_1^{-1}b_1^{-1}a_1b_1=1$. If $x=b_1a_1^{-1}$ and $y=a_1^{-1}$, then $y^{-1}x^2y=x$ and $\langle h_2,h_3\rangle=\langle x,y\rangle=BS(2,1)$, that is a contradiction.
\end{enumerate}
Therefore, if $a_2=a_1$, then $a_4=b_1$ and $b_4=b_2$. So, $T_C=[a_1,b_1,a_1,b_2,b_1,a_1,b_1,b_2]$ and the relation of such cycle is $a_1^{-1}b_1a_1^{-1}b_2b_1^{-1}a_1b_1^{-1}b_2=1$. Since $\{a_1,b_1,b_2\}=supp(\alpha)$, exactly one of the following cases may be happened:
\begin{enumerate}
\item[a)]
$a_1=1$: $b_1b_2b_1^{-2}b_2=1$, where $\{b_1,b_2\}=\{h_2,h_3\}$.
\item[b)]
$b_1=1$: $a_1b_2a_1^{-2}b_2=1$, where $\{a_1,b_2\}=\{h_2,h_3\}$.
\item[c)]
$b_2=1$: $b_1 a_1^{-1} b_1 a_1 b_1^{-1} a_1=1$, where $\{a_1,b_1\}=\{h_2,h_3\}$.
\end{enumerate}
Hence, if $a_2=a_1$, then $T_{C}=[a_1,b_1,a_1,b_2,b_1, a_1,b_1,b_2]$ and  $T_{C'}=[a_1,b_1,b_2,a_1,b_1,a_1,b_2,b_1]$, where $\{a_1,b_1, b_2\}=supp(\alpha)$. Also with this assumption, exactly one of the relations $14$, $22$ or $26$ of Table \ref{tab-C4} will be satisfied in $G$.
\item[B)] $a_2\neq a_1$: By the relation \ref{e-8}, we have $a_2\neq a_1$, $a_2\neq b_1$ and $a_1\neq b_1$. So, $\{a_1,b_1,a_2\}=supp(\alpha)$. Also $b_2= a_1$ because $b_2\neq b_1$ and $b_2\neq a_2$, by the relation \ref{e-8}. So, $T_{C}=[a_1,b_1,a_2,a_1,b_1,a_1,a_4,b_4]$ and  $T_{C'}=[a_1,b_1,a_1,a_2,b_1,a_1,b_4,a_4]$. Also by the relation \ref{e-8}, $a_4\neq a_1$, $b_4\neq a_1$ and $a_4\neq b_4$. Therefore, there are two cases for choosing $a_4$ and $b_4$ which we show that one of them leads to a contradiction with our assumptions.\\
Suppose that $a_4=b_1$ and $b_4=a_2$. So,  $T_C=[a_1,b_1,a_2,a_1,b_1,a_1,b_1,a_2]$ and the relation of such cycle is $a_1^{-1}b_1a_2^{-1}a_1b_1^{-1}a_1b_1^{-1}a_2=1$. Since $\{a_1,b_1,a_2\}=supp(\alpha)$, exactly one of the following cases may be happened:
\begin{enumerate}
\item[a)]
$a_1=1$: $a_2^{-1}b_1^{\ 2}a_2=b_1$ and so $\langle h_2,h_3\rangle=BS(2,1)$, that is a contradiction.
\item[b)]
$b_1=1$: $a_2^{-1}a_1^{\ 2}a_2=a_1$ and so $\langle h_2,h_3\rangle=BS(2,1)$, that is a contradiction.
\item[c)]
$a_2=1$: $a_1^{-1}b_1a_1b_1^{-1}a_1b_1^{-1}=1$. If $x=a_1b_1^{-1}$ and $y=b_1^{-1}$, then $y^{-1}x^2y=x$ and $\langle h_2,h_3\rangle=\langle x,y\rangle=BS(2,1)$, that is a contradiction.
\end{enumerate}
Therefore, if $a_2\neq a_1$, then $a_4=a_2$ and $b_4=b_1$. So, $T_C=[a_1,b_1,a_2,a_1,b_1,a_1,a_2,b_1]$ and the relation of such cycle is $a_1^{-1}b_1a_2^{-1}a_1b_1^{-1}a_1a_2^{-1}b_1=1$. Since $\{a_1,b_1,a_2\}=supp(\alpha)$, exactly one of the following cases may be happened:
\begin{enumerate}
\item[a)]
$a_1=1$: $b_1a_2b_1^{-2}a_2=1$, where $\{b_1,a_2\}=\{h_2,h_3\}$.
\item[b)]
$b_1=1$: $a_1a_2a_1^{-2}a_2=1$, where $\{a_1,a_2\}=\{h_2,h_3\}$.
\item[c)]
$a_2=1$: $b_1 a_1^{-1} b_1 a_1 b_1^{-1} a_1=1$, where $\{a_1,b_1\}=\{h_2,h_3\}$.
\end{enumerate}
Hence, if $a_2\neq a_1$, then $T_{C}=[a_1,b_1,a_2,a_1,b_1, a_1,a_2,b_1]$ and  $T_{C'}=[a_1,b_1,a_1,a_2,b_1,a_1,b_1,a_2]$, where $\{a_1,b_1, a_2\}=supp(\alpha)$. Also with this assumption, exactly one of the relations $14$, $22$ or $26$ of Table \ref{tab-C4} will be satisfied in $G$.
\end{enumerate}
\item
$[a_1,b_1,a_2',b_2',a_3',b_3',a_4',b_4']=[b_2,a_2,b_1,a_1,b_4,a_4,b_3,a_3]$: So $b_1=a_2$, that is a contradiction with the relation \ref{e-8}.
\end{enumerate}
Hence by the discussion above, if $K(\alpha,\beta)$ contains two squares with exactly one edge in common, then exactly one of the relations $14$, $22$ or $26$ of Table \ref{tab-C4} will be satisfied in $G$.
\end{proof}

\begin{rem}\label{r-1}
{\rm By the proof of Theorem \ref{C4-C4}, if the Kaplansky graph contains two squares $C$ and $C'$ with exactly one common edge, then the relations and the $8$-tuples of such cycles must be of just one of the following forms, where the first two components of these $8$-tuples are corresponded to the common edge of $C$ and $C'$:
\begin{enumerate}
\item[1)]
$T_C=[h_2,1,h_2,h_3,1,h_2,1,h_3]$ and $T_{C'}=[h_2,1,h_3,h_2,1,h_2,h_3,1]$, and vice versa, with relation $14$ of Table \ref{tab-C4}.
\item[2)]
$T_{C}=[1,h_2,1,h_3,h_2,1,h_2,h_3]$ and $T_{C'}=[1,h_2,h_3,1,h_2,1,h_3,h_2]$, and vice versa, with relation $14$ of Table \ref{tab-C4}.
\item[3)]
$T_{C}=[h_3,1,h_3,h_2,1,h_3,1,h_2]$ and $T_{C'}=[h_3,1,h_2,h_3,1,h_3,h_2,1]$, and vice versa, with relation $22$ of Table \ref{tab-C4}.
\item[4)]
$T_{C}=[1,h_3,1,h_2,h_3,1,h_3,h_2]$ and $T_{C'}=[1,h_3,h_2,1,h_3,1,h_2,h_3]$, and vice versa, with relation $22$ of Table \ref{tab-C4}.
\item[5)]
$T_{C}=[h_2,h_3,h_2,1,h_3,h_2,h_3,1]$ and $T_{C'}=[h_2,h_3,1,h_2,h_3,h_2,1,h_3]$, and vice versa, with relation $26$ of Table \ref{tab-C4}.
\item[6)]
$T_{C}=[h_3,h_2,h_3,1,h_2,h_3,h_2,1]$ and $T_{C'}=[h_3,h_2,1,h_3,h_2,h_3,1,h_2]$, and vice versa, with relation $26$ of Table \ref{tab-C4}.
\end{enumerate} }
\end{rem}

\begin{figure}[ht]
\centering
\psscalebox{0.8 0.8} {\begin{tikzpicture}
\node (14) [circle, minimum size=3pt, fill=black, line width=0.625pt, draw=black] at (75.0pt, -62.5pt)  {};
\node (12) [circle, minimum size=3pt, fill=black, line width=0.625pt, draw=black] at (75.0pt, -100.0pt)  {};
\node (15) [circle, minimum size=3pt, fill=black, line width=0.625pt, draw=black] at (112.5pt, -62.5pt)  {};
\node (16) [circle, minimum size=3pt, fill=black, line width=0.625pt, draw=black] at (112.5pt, -100.0pt)  {};
\node (10) [circle, minimum size=3pt, fill=black, line width=0.625pt, draw=black] at (150.0pt, -100.0pt)  {};
\node (9) [circle, minimum size=3pt, fill=black, line width=0.625pt, draw=black] at (150.0pt, -62.5pt)  {};
\node (19) [circle, minimum size=3pt, fill=black, line width=0.625pt, draw=black] at (187.5pt, -100.0pt)  {};
\node (20) [circle, minimum size=3pt, fill=black, line width=0.625pt, draw=black] at (187.5pt, -62.5pt)  {};
\node (22) [circle, minimum size=3pt, fill=black, line width=0.625pt, draw=black] at (225.0pt, -62.5pt)  {};
\node (21) [circle, minimum size=3pt, fill=black, line width=0.625pt, draw=black] at (225.0pt, -100.0pt)  {};
\node (2) [circle, minimum size=3pt, fill=black, line width=0.625pt, draw=black] at (262.5pt, -100.0pt)  {};
\node (4) [circle, minimum size=3pt, fill=black, line width=0.625pt, draw=black] at (262.5pt, -62.5pt)  {};
\node (6) [circle, minimum size=3pt, fill=black, line width=0.625pt, draw=black] at (300.0pt, -100.0pt)  {};
\node (5) [circle, minimum size=3pt, fill=black, line width=0.625pt, draw=black] at (300.0pt, -62.5pt)  {};
\draw [line width=1.25, color=black] (14) to  (15);
\draw [line width=1.25, color=black] (16) to  (15);
\draw [line width=1.25, color=black] (12) to  (16);
\draw [line width=1.25, color=black] (12) to  (14);
\draw [line width=1.25, color=black] (2) to  (6);
\draw [line width=1.25, color=black] (5) to  (6);
\draw [line width=1.25, color=black] (16) to  (10);
\draw [line width=1.25, color=black] (15) to  (9);
\draw [line width=1.25, color=black] (21) to  (22);
\draw [line width=1.25, dotted, color=black] (21) to  (19);
\draw [line width=1.25, dotted, color=black] (22) to  (20);
\draw [line width=1.25, color=black] (20) to  (19);
\draw [line width=1.25, color=black] (10) to  (9);
\draw [line width=1.25, color=black] (10) to  (19);
\draw [line width=1.25, color=black] (20) to  (9);
\draw [line width=1.25, color=black] (4) to  (2);
\draw [line width=1.25, color=black] (5) to  (4);
\draw [line width=1.25, color=black] (2) to  (21);
\draw [line width=1.25, color=black] (4) to  (22);
\draw [line width=1.25, color=black] (14) to  [in=162, out=18] (5);
\draw [line width=1.25, color=black] (6) to  [in=342, out=198] (12);
\node at (189.125pt, -160.375pt) {\textcolor{black}{$\mathbf{L_n}$}};
\end{tikzpicture}}
\hspace{2.5cc}
\psscalebox{0.8 0.8} {\begin{tikzpicture}
\node (14) [circle, minimum size=3pt, fill=black, line width=0.625pt, draw=black] at (75.0pt, -62.5pt)  {};
\node (12) [circle, minimum size=3pt, fill=black, line width=0.625pt, draw=black] at (75.0pt, -100.0pt)  {};
\node (15) [circle, minimum size=3pt, fill=black, line width=0.625pt, draw=black] at (112.5pt, -62.5pt)  {};
\node (16) [circle, minimum size=3pt, fill=black, line width=0.625pt, draw=black] at (112.5pt, -100.0pt)  {};
\node (10) [circle, minimum size=3pt, fill=black, line width=0.625pt, draw=black] at (150.0pt, -100.0pt)  {};
\node (9) [circle, minimum size=3pt, fill=black, line width=0.625pt, draw=black] at (150.0pt, -62.5pt)  {};
\node (19) [circle, minimum size=3pt, fill=black, line width=0.625pt, draw=black] at (187.5pt, -100.0pt)  {};
\node (20) [circle, minimum size=3pt, fill=black, line width=0.625pt, draw=black] at (187.5pt, -62.5pt)  {};
\node (22) [circle, minimum size=3pt, fill=black, line width=0.625pt, draw=black] at (225.0pt, -62.5pt)  {};
\node (21) [circle, minimum size=3pt, fill=black, line width=0.625pt, draw=black] at (225.0pt, -100.0pt)  {};
\node (2) [circle, minimum size=3pt, fill=black, line width=0.625pt, draw=black] at (262.5pt, -100.0pt)  {};
\node (4) [circle, minimum size=3pt, fill=black, line width=0.625pt, draw=black] at (262.5pt, -62.5pt)  {};
\node (6) [circle, minimum size=3pt, fill=black, line width=0.625pt, draw=black] at (300.0pt, -100.0pt)  {};
\node (5) [circle, minimum size=3pt, fill=black, line width=0.625pt, draw=black] at (300.0pt, -62.5pt)  {};
\draw [line width=1.25, color=black] (14) to  (15);
\draw [line width=1.25, color=black] (16) to  (15);
\draw [line width=1.25, color=black] (12) to  (16);
\draw [line width=1.25, color=black] (12) to  (14);
\draw [line width=1.25, color=black] (2) to  (6);
\draw [line width=1.25, color=black] (5) to  (6);
\draw [line width=1.25, color=black] (16) to  (10);
\draw [line width=1.25, color=black] (15) to  (9);
\draw [line width=1.25, color=black] (21) to  (22);
\draw [line width=1.25, dotted, color=black] (21) to  (19);
\draw [line width=1.25, dotted, color=black] (22) to  (20);
\draw [line width=1.25, color=black] (20) to  (19);
\draw [line width=1.25, color=black] (10) to  (9);
\draw [line width=1.25, color=black] (10) to  (19);
\draw [line width=1.25, color=black] (20) to  (9);
\draw [line width=1.25, color=black] (4) to  (2);
\draw [line width=1.25, color=black] (5) to  (4);
\draw [line width=1.25, color=black] (2) to  (21);
\draw [line width=1.25, color=black] (4) to  (22);
\draw [line width=1.25, color=black] (5) to  [in=66, out=135] (12);
\draw [line width=1.25, color=black] (6) to  [in=294, out=225] (14);
\node at (189.125pt, -160.375pt) {\textcolor{black}{$\mathbf{M_n}$}};
\end{tikzpicture}}
\caption{Two graphs which are not isomorphic to $K(\alpha,\beta)$}\label{f-3}
\end{figure}
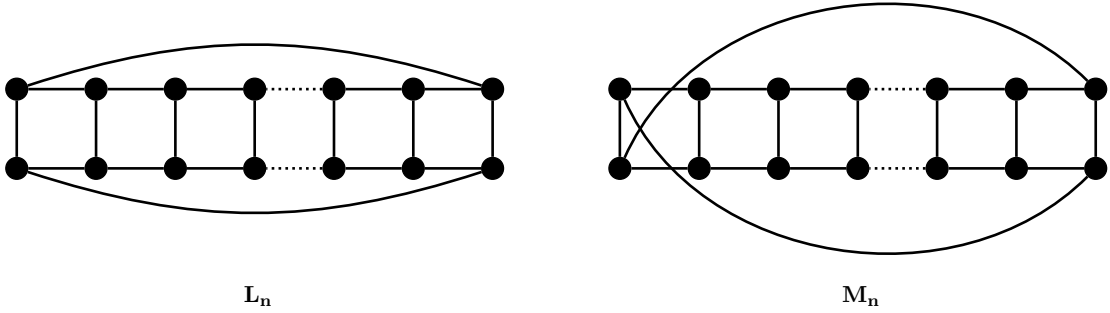

\begin{thm}\label{L-M}
The Kaplansky graph $K(\alpha,\beta)$ is isomorphic to none of the graphs $L_n$ and $M_n$ from Figure \ref{f-3}.
\end{thm}
\begin{proof}
If $n$ is the number of vertices of the graph $L_n$ or $M_n$, then it can be seen in Figure \ref{f-3} that the number of cycles of length $4$ in $L_n$ or $M_n$ is equal to $\frac{n}{2}$, and each two consecutive $C_4$ cycles have a common edge. So, if $K(\alpha,\beta)$ contains no two $C_4$ cycles with exactly one common edge, then it cannot be isomorphic to the graphs $L_n$ and $M_n$. Suppose that the graph $K(\alpha,\beta)$ contains two cycles of length $4$ with exactly one edge in common. 
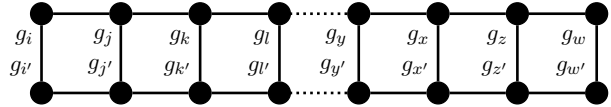
\begin{figure}[h]
\psscalebox{0.8 0.8} {\begin{tikzpicture}
\node (14) [circle, minimum size=3pt, fill=black, line width=0.625pt, draw=black] at (75.0pt, -62.5pt)  {};
\node (12) [circle, minimum size=3pt, fill=black, line width=0.625pt, draw=black] at (75.0pt, -100.0pt)  {};
\node (15) [circle, minimum size=3pt, fill=black, line width=0.625pt, draw=black] at (112.5pt, -62.5pt)  {};
\node (16) [circle, minimum size=3pt, fill=black, line width=0.625pt, draw=black] at (112.5pt, -100.0pt)  {};
\node (10) [circle, minimum size=3pt, fill=black, line width=0.625pt, draw=black] at (150.0pt, -100.0pt)  {};
\node (9) [circle, minimum size=3pt, fill=black, line width=0.625pt, draw=black] at (150.0pt, -62.5pt)  {};
\node (19) [circle, minimum size=3pt, fill=black, line width=0.625pt, draw=black] at (187.5pt, -100.0pt)  {};
\node (20) [circle, minimum size=3pt, fill=black, line width=0.625pt, draw=black] at (187.5pt, -62.5pt)  {};
\node (22) [circle, minimum size=3pt, fill=black, line width=0.625pt, draw=black] at (225.0pt, -62.5pt)  {};
\node (21) [circle, minimum size=3pt, fill=black, line width=0.625pt, draw=black] at (225.0pt, -100.0pt)  {};
\node (2) [circle, minimum size=3pt, fill=black, line width=0.625pt, draw=black] at (262.5pt, -100.0pt)  {};
\node (4) [circle, minimum size=3pt, fill=black, line width=0.625pt, draw=black] at (262.5pt, -62.5pt)  {};
\node (6) [circle, minimum size=3pt, fill=black, line width=0.625pt, draw=black] at (300.0pt, -100.0pt)  {};
\node (5) [circle, minimum size=3pt, fill=black, line width=0.625pt, draw=black] at (300.0pt, -62.5pt)  {};
\node (1) [circle, minimum size=3pt, fill=black, line width=0.625pt, draw=black] at (337.5pt, -62.5pt)  {};
\node (3) [circle, minimum size=3pt, fill=black, line width=0.625pt, draw=black] at (337.5pt, -100.0pt)  {};
\draw [line width=1.25, color=black] (14) to  (15);
\draw [line width=1.25, color=black] (16) to  (15);
\draw [line width=1.25, color=black] (12) to  (16);
\draw [line width=1.25, color=black] (12) to  (14);
\draw [line width=1.25, color=black] (2) to  (6);
\draw [line width=1.25, color=black] (5) to  (6);
\draw [line width=1.25, color=black] (16) to  (10);
\draw [line width=1.25, color=black] (15) to  (9);
\draw [line width=1.25, color=black] (21) to  (22);
\draw [line width=1.25, dotted, color=black] (21) to  (19);
\draw [line width=1.25, dotted, color=black] (22) to  (20);
\draw [line width=1.25, color=black] (20) to  (19);
\draw [line width=1.25, color=black] (10) to  (9);
\draw [line width=1.25, color=black] (10) to  (19);
\draw [line width=1.25, color=black] (20) to  (9);
\draw [line width=1.25, color=black] (4) to  (2);
\draw [line width=1.25, color=black] (5) to  (4);
\draw [line width=1.25, color=black] (2) to  (21);
\draw [line width=1.25, color=black] (4) to  (22);
\draw [line width=1.25, color=black] (5) to  (1);
\draw [line width=1.25, color=black] (1) to  (3);
\draw [line width=1.25, color=black] (6) to  (3);
\node at (66.875pt, -73.75pt) {\textcolor{black}{$g_i$}};
\node at (65.625pt, -88.125pt) {\textcolor{black}{$g_{i'}$}};
\node at (104.375pt, -73.75pt) {\textcolor{black}{$g_j$}};
\node at (103.125pt, -88.125pt) {\textcolor{black}{$g_{j'}$}};
\node at (139.375pt, -88.125pt) {\textcolor{black}{$g_{k'}$}};
\node at (140.625pt, -73.75pt) {\textcolor{black}{$g_k$}};
\node at (178.125pt, -88.125pt) {\textcolor{black}{$g_{l'}$}};
\node at (179.375pt, -73.75pt) {\textcolor{black}{$g_l$}};
\node at (214.375pt, -73.75pt) {\textcolor{black}{$g_y$}};
\node at (213.125pt, -88.125pt) {\textcolor{black}{$g_{y'}$}};
\node at (251.875pt, -88.125pt) {\textcolor{black}{$g_{x'}$}};
\node at (253.125pt, -73.75pt) {\textcolor{black}{$g_x$}};
\node at (289.375pt, -88.125pt) {\textcolor{black}{$g_{z'}$}};
\node at (290.625pt, -73.75pt) {\textcolor{black}{$g_z$}};
\node at (326.875pt, -73.75pt) {\textcolor{black}{$g_w$}};
\node at (325.625pt, -88.125pt) {\textcolor{black}{$g_{w'}$}};
\end{tikzpicture}}
\caption{Consecutive cycles of length $4$ in the graph $K(\alpha,\beta)$}\label{f-4}
\end{figure}

Suppose that $K(\alpha,\beta)$ contains a subgraph as like as Figure \ref{f-4}, whose number of consecutive $C_4$ cycles is denoted by $m$. We denote these cycles by $C_1,C_2,\ldots,C_{m-1}$ and $C_m$ from the left to the right. For all $i\in \{1,2,\ldots,m-1\}$, the $8$-tuples of $C_i$ and $C_{i+1}$ must be of exactly one of the $6$ cases in Remark \ref{r-1}, where the first two components are corresponded to their common edge. In the following, we prove the statement of the theorem for the case that $C_1=C$ and $C_2=C'$ in the first item of Remark \ref{r-1}. The other $5$ cases can be proven similarly.

Suppose that $C_1=C$, $C_2=C'$ are the same as the first item in Remark \ref{r-1}, i.e.  $T_{C_1}=[h_2,1,h_2,h_3,1,\\h_2,1,h_3]$ and $T_{C_2}=[h_2,1,h_3,h_2,1,h_2,h_3,1]$ such that 
\begin{equation}\label{e-11}
R(T_{C_1}):\left\{
\begin{array}{l}
h_2g_j=g_{j'}\\
h_2g_{j'}=h_3g_{i'}\\
g_{i'}=h_2g_{i}\\
g_{i}=h_3g_j\\
\end{array} \right.
\qquad \qquad
R(T_{C_2}):\left\{
\begin{array}{l}
h_2g_j=g_{j'}\\
h_3g_{j'}=h_2g_{k'}\\
g_{k'}=h_2g_{k}\\
h_3g_{k}=g_j\\
\end{array} \right.
\end{equation}
where $\{g_j,g_{j'},g_{i'},g_i\}$ and $\{g_j,g_{j'},g_{k'},g_k\}$ are the vertex sets of $C_1$ and $C_2$, respectively.\\
Using induction on $m$, we show that for all $i\in \{1,2,\ldots,m-1\}$, the corresponding relations of $C_i$ and $C_{i+1}$ are the same as the relations \ref{e-11}, where the first relation of each part is related to the common edge between $C_i$ and $C_{i+1}$ and the relations of these cycles are written clockwise and counter clockwise, respectively.

If $m=2$, then the above statement is obviously true. Suppose that this statement is true for $m-1$. So, $T_{C_{m-2}}=[h_2,1,h_2,h_3,1,h_2,1,h_3]$ and $T_{C_{m-1}}=[h_2,1,h_3,h_2,1,h_2,h_3,1]$, where the first two components are related to the common edge between these cycles. Suppose that $\{g_x,g_{x'},g_{y'},g_y\}$ and $\{g_x,g_{x'},g_{z'},g_z\}$ are the vertex sets of $C_{m-2}$ and $C_{m-1}$, respectively, and the corresponding relations of these cycles are as relations \ref{e-11}, by replacing $g_j,g_{j'},g_{i'},g_i,g_{k'}$ and $g_k$ with $g_x,g_{x'},g_{y'},g_y,g_{z'}$ and $g_z$, respectively. Obviously, we may rewrite the relations of $C_{m-1}$ such that the first relation being $g_{z'}=h_2g_{z}$. So, there is $T'_{C_{m-1}}\in \mathcal{T}(C_{m-1})$ equal to $T_{C_1}=T_C=[h_2,1,h_2,h_3,1,h_2,1,h_3]$ with the corresponding relations as follows:
\begin{equation*}
R(T'_{C_{m-1}}):\left\{
\begin{array}{l}
h_2g_z=g_{z'}\\
h_2g_{z'}=h_3g_{x'}\\
g_{x'}=h_2g_{x}\\
g_{x}=h_3g_z\\
\end{array} \right.
\end{equation*}
Therefore by Remark \ref{r-1}, there is $T_{C_{m}}\in \mathcal{T}(C_{m})$ equal to $[h_2,1,h_3,h_2,1,h_2,h_3,1]$ with the corresponding relations
\begin{equation*}
R(T_{C_m}):\left\{
\begin{array}{l}
h_2g_z=g_{z'}\\
h_3g_{z'}=h_2g_{w'}\\
g_{w'}=h_2g_{w}\\
h_3g_{w}=g_z\\
\end{array} \right.
\end{equation*}
where $\{g_z,g_{z'},g_{w'},g_w\}$ is the vertex set of $C_{m}$.\\
Hence, for all $i\in \{1,2,\ldots,m-1\}$, the corresponding relations of $C_i$ and $C_{i+1}$ are the same as the relations \ref{e-11}, where the first relation of each part is related to the common edge between $C_i$ and $C_{i+1}$ and the relations of these cycles are written clockwise and counter clockwise, respectively. 

Now suppose that $K(\alpha,\beta)$ is isomorphic to the graph $L_n$, where the vertex set of $K(\alpha,\beta)$ is equal to the set $B=\{g_1,g_2,\ldots,g_n\}$. Also suppose that the vertices of $L_n$ in the top row and the bottom row are $g_1',g_2',g_3',g_4',\ldots, g_{\frac{n}{2}-2}',g_{\frac{n}{2}-1}',g_{\frac{n}{2}}'$ and $g_n',g_{n-1}',g_{n-2}',g_{n-3}',\ldots, g_{\frac{n}{2}+3}',g_{\frac{n}{2}+2}',g_{\frac{n}{2}+1}'$, respectively from the left to the right, where $\{g_a' \vert a\in \{1,2,\ldots,n\}\}=\{g_1,g_2,\ldots,g_n\}$. Let $m=\frac{n}{2}$ be the number of cycles of length $4$ in  $L_n$ that each two of them which are consecutive have a common edge. We denote these cycles with $C_1,C_2,\ldots,C_{m-1}$ and $C_m$, where $C_1$ and $C_2$ are the cycles with vertex sets $\{g_1',g_2',g_{n-1}',g_n'\}$ and $\{g_2',g_3',g_{n-2}',g_{n-1}'\}$, respectively. In addition, Suppose that $C_1=C$, $C_2=C'$ are the same as the first item in Remark \ref{r-1}, i.e. the  $8$-tuples of $C_1$ and $C_2$ are  $[h_2,1,h_2,h_3,1,h_2,1,h_3]$ and $[h_2,1,h_3,h_2,1,h_2,h_3,1]$, respectively, and the corresponding relations of these cycles are as relations \ref{e-11}, by replacing $g_j,g_{j'},g_{i'},g_i,g_{k'}$ and $g_k$ with $g_2',g_{n-1}',g_{n}',g_1',g_{n-2}'$ and $g_3'$, respectively. \\
With the discussion above, for all $i\in \{1,2,\ldots,m-1\}$, the corresponding relations of $C_i$ and $C_{i+1}$ are the same as the relations of $C_1=C$, $C_2=C'$, when the first relation of each part is related to the common edge between $C_i$ and $C_{i+1}$ and the relations of these cycles are written clockwise and counter clockwise, respectively. Therefore we have $g_1'=h_3g_2', g_2'=h_3g_3', \ldots,g_{\frac{n}{2}-2}'=h_3g_{\frac{n}{2}-1}', g_{\frac{n}{2}-1}'=h_3g_{\frac{n}{2}}'$ and $g_{\frac{n}{2}}'=h_3g_1'$. So, $g_1'=h_3^{\frac{n}{2}}g_1'$ and hence $h_3^{\frac{n}{2}}=1$. Since $G$ is a torsion-free group, we have $h_3=1$ that is a contradiction with $\vert supp(\alpha)\vert=3$. Therefore, with above assumptions on $C_1$ and $C_2$, the graph $K(\alpha,\beta)$ cannot be isomorphic to the graph $L_n$. 

Now suppose that $K(\alpha,\beta)$ is isomorphic to the graph $M_n$, where the vertex set of $K(\alpha,\beta)$ is equal to the set $B=\{g_1,g_2,\ldots,g_n\}$. With the same assumptions such as above on the cycles of length $4$ in $M_n$ and by using the corresponding relations of these cycles, we have $g_1'=h_3g_2', g_2'=h_3g_3', \ldots,g_{\frac{n}{2}-2}'=h_3g_{\frac{n}{2}-1}', g_{\frac{n}{2}-1}'=h_3g_{\frac{n}{2}}',g_{\frac{n}{2}}'=h_3g_n'$ and $g_n'=h_2g_1'$. So, $g_1'=h_3^{\frac{n}{2}}g_n'$ and $g_n'=h_2g_1'$. Therefore, $g_1'=h_3^{\frac{n}{2}}h_2g_1'$ and hence $h_3^{\frac{n}{2}}h_2=1$. So, $h_2=h_3^{-\frac{n}{2}}$ and therefore $\langle h_2,h_3\rangle=\langle h_3\rangle$ is abelian, that is a contradiction since we know that abelian groups satisfy the zero divisor conjecture. Therefore, with above assumptions on $C_1$ and $C_2$, the graph $K(\alpha,\beta)$ cannot be isomorphic to the graph $M_n$. \\
Hence, the statement of this theorem is proven for the case that $T_{C_1}=[h_2,1,h_2,h_3,1,h_2,1,h_3]$ and $T_{C_2}=[h_2,1,h_3,h_2,1,h_2,h_3,1]$. Since all cycles of length $4$ in the graphs $L_n$ and $M_n$ are consecutive, it is easy to see that if we consider $T_{C_1}=[h_2,1,h_3,h_2,1,h_2,h_3,1]$ and $T_{C_2}=[h_2,1,h_2,h_3,1,h_2,1,h_3]$, then the above discussion is true for the latter case too. Therefore, the statement of this theorem is proven for the first item of Remark \ref{r-1}. 

With a similar discussion such as above, the statement of this theorem can be proven for the other $5$ cases in Remark \ref{r-1}. Hence, the graph $K(\alpha,\beta)$ is not isomorphic to the graphs $L_n$ and $M_n$ in Figure \ref{f-3}, where the number of vertices of $K(\alpha,\beta)$, $L_n$ and $M_n$ is equal to $n$.
\end{proof}


\section{\bf Forbidden subgraphs of Kaplansky graphs over $\mathbb{F}_2$}\label{S2}
In previous sections, we studied the existence of triangles and square and some subgraphs containing them in the graph $K_{\mathbb{F}}(\alpha,\beta)$. Also, we saw that $K_3$ and $K_{2,3}$ are two forbidden subgraphs of the graph $K(\alpha,\beta)$.  With the same discussion such as about $C_4$ cycles, we can study cycles of other lengths in the graph $K(\alpha,\beta)$ by using their relations. In this section, by using cycles up to lengths $7$ and their relations, we find another forbidden subgraphs of the graph $K(\alpha,\beta)$. The procedure of finding these forbidden subgraphs is as like as the procedure of finding previous examples. So, the frequent tedious details has been omitted. Such forbidden subgraphs of the Kaplansky graph are listed in Table \ref{tab-forbiddens}. In the following, we discuss about such subgraphs. Here, forbidden subgraphs are numbered from $1$ to $44$ such that the forbidden subgraph $K_{2,3}$ is numbered by $1$. 

\subsection{$\mathbf{K_{2,3}}$}

By Theorem \ref{K2,3}, $K(\alpha,\beta)$ contains no subgraph isomorphic to the complete bipartite graph $K_{2,3}$.

\subsection{$\mathbf{C_4--C_5}$}
With the same discussion such as about $K_{2,3}$, it can be seen that there are $121$ different cases for the relations of the cycles $C_4$ and $C_5$ in this structure. Using GAP \cite{gap}, we see that the groups with two generators $h_2$ and $h_3$ and two relations which are between $111$ cases of these $121$ cases are finite and solvable, that is a contradiction with the assumptions. So, there are just $10$ cases for the relations of the cycles $C_4$ and $C_5$ which may lead to the existence of a subgraph isomorphic to the graph $C_4--C_5$ in $K(\alpha,\beta)$. It can be seen that all of these $10$ cases lead to contradictions and so, the graph $K(\alpha,\beta)$ contains no subgraph isomorphic to the graph $C_4--C_5$.

Each group with two generators $h_2$ and $h_3$ and two relations which are between the latter $10$ cases is a quotient of $B(1,k)$, for some integer $k$, or has a torsion element.

\subsection{$\mathbf{C_4--C_6}$}
It can be seen that there are $658$ different cases for the relations of the cycles $C_4$ and $C_6$ in this structure. By considering all groups with two generators $h_2$ and $h_3$ and two relations which are between these cases and by using GAP \cite{gap}, we see that $632$ groups are finite and solvable, or just finite. So, there are just $20$ cases for the relations of the cycles $C_4$ and $C_6$ which may lead to the existence of a subgraph isomorphic to the graph $C_4--C_6$ in $K(\alpha,\beta)$. It can be seen that all of these $20$ cases lead to  contradictions and so, the graph $K(\alpha,\beta)$ contains no subgraph isomorphic to the graph $C_4--C_6$.

Each group with two generators $h_2$ and $h_3$ and two relations which are between the latter $20$ cases is a quotient of $B(1,k)$, for some integer $k$, has a torsion element or is a cyclic group.

\subsection{$\mathbf{C_4-C_5(-C_5-)}$}
It can be seen that there are $42$ cases for the relations of a cycle $C_4$ and two cycles $C_5$ in the graph $C_4-C_5(-C_5-)$. Using GAP \cite{gap}, we see that all groups with two generators $h_2$ and $h_3$ and three relations which are between $38$ cases of these $42$ cases are finite and solvable, that is a contradiction. So, there are just $4$ cases for the relations of these cycles which may lead to the existence of a subgraph isomorphic to the graph $C_4-C_5(-C_5-)$ in $K(\alpha,\beta)$. It can be seen that all of these $4$ cases lead to  contradictions and so, the graph $K(\alpha,\beta)$ contains no subgraph isomorphic to the graph $C_4-C_5(-C_5-)$.

Each group with two generators $h_2$ and $h_3$ and two relations which are between the latter $4$ cases is a quotient of $B(1,k)$, for some integer $k$, or has a torsion element.

\subsection{$\mathbf{C_4-C_5(-C_4-)}$}
It can be seen that there are $4$ cases for the relations of two cycles $C_4$ and a cycle $C_5$ in this structure. By considering all groups with two generators $h_2$ and $h_3$ and three relations which are between these cases and by using GAP \cite{gap}, we see that all of these groups are finite and solvable. So, the graph $K(\alpha,\beta)$ contains no subgraph isomorphic to the graph $C_4-C_5(-C_4-)$.

\subsection{$\mathbf{C_4-C_5(-C_6--)}$}
It can be seen that there are $126$ cases for the relations of a cycle $C_4$, a cycle $C_5$ and a cycle $C_6$ in the graph $C_4-C_5(-C_6--)$. Using GAP \cite{gap}, we see that all groups with two generators $h_2$ and $h_3$ and three relations which are between $122$ cases of these $126$ cases are finite and solvable, that is a contradiction with the assumptions. So, there are just $4$ cases for the relations of these cycles which may lead to the existence of a subgraph isomorphic to the graph $C_4-C_5(-C_6--)$ in $K(\alpha,\beta)$. It can be seen that all of these $4$ cases lead to contradictions and so, the graph $K(\alpha,\beta)$ contains no subgraph isomorphic to the graph $C_4-C_5(-C_6--)$.

Each group with two generators $h_2$ and $h_3$ and two relations which are between the latter $4$ cases is a quotient of $B(1,k)$, for some integer $k$.

\subsection{$\mathbf{C_4-C_5(-C_6-)}$}
It can be seen that there are $462$ cases for the relations of a cycle $C_4$, a cycle $C_5$ and a cycle $C_6$ in the graph $C_4-C_5(-C_6-)$. By considering all groups with two generators $h_2$ and $h_3$ and two relations which are between these cases and by using GAP \cite{gap}, we see that $436$ groups are finite and solvable, or just finite. So, there are just $22$ cases for the relations of these cycles which may lead to the existence of a subgraph isomorphic to the graph $C_4-C_5(-C_6-)$ in $K(\alpha,\beta)$.  It can be shown that all of these $22$ cases lead to contradictions and so, the graph $K(\alpha,\beta)$ contains no subgraph isomorphic to the graph $C_4-C_5(-C_6-)$.

Each group with two generators $h_2$ and $h_3$ and two relations which are between the latter $22$ cases is a quotient of $B(1,k)$, for some integer $k$, is a cyclic group or is a solvable group.

\subsection{$\mathbf{C_4-C_5(-C_7--)}$}
It can be seen that there are $648$ cases for the relations of a cycle $C_4$, a cycle $C_5$ and a cycle $C_7$ in the graph $C_4-C_5(-C_7--)$. Using GAP \cite{gap}, we see that all groups with two generators $h_2$ and $h_3$ and three relations which are between $608$ cases of these $648$ cases are finite and solvable, that is a contradiction with the assumptions. So, there are just $40$ cases for the relations of these cycles which may lead to the existence of a subgraph isomorphic to the graph $C_4-C_5(-C_7--)$ in $K(\alpha,\beta)$. It can be shown that all of these $40$ cases lead to contradictions and so, the graph $K(\alpha,\beta)$ contains no subgraph isomorphic to the graph $C_4-C_5(-C_7--)$.

Each group with two generators $h_2$ and $h_3$ and two relations which are between the latter $40$ cases is a quotient of $B(1,k)$, for some integer $k$, has a torsion element, is a cyclic group or is a solvable group.

\subsection{$\mathbf{C_5--C_5(--C_5)}$}
It can be seen that there are $192$ cases for the relations of three cycles $C_5$ in the graph $C_5--C_5(--C_5)$. Using GAP \cite{gap}, we see that all groups with two generators $h_2$ and $h_3$ and three relations which are between $188$ cases of these $192$ cases are finite and solvable, that is a contradiction with the assumptions. So, there are just $4$ cases for the relations of these cycles which may lead to the existence of a subgraph isomorphic to the graph $C_5--C_5(--C_5)$ in $K(\alpha,\beta)$. It can be shown that all of these $4$ cases lead to contradictions and so, the graph $K(\alpha,\beta)$ contains no subgraph isomorphic to the graph $C_5--C_5(--C_5)$.

Each group with two generators $h_2$ and $h_3$ and two relations which are between the latter $4$ cases is a quotient of $B(1,k)$, for some integer $k$.

\subsection{$\mathbf{C_5--C_5(--C_6)}$}
It can be seen that there are $1006$ cases for the relations of two cycles $C_5$ and a cycle $C_6$ in the graph $C_5--C_5(--C_6)$. Using GAP \cite{gap}, we see that all groups with two generators $h_2$ and $h_3$ and three relations which are between $986$ cases of these $1006$ cases are finite and solvable, that is a contradiction with the assumptions. So, there are just $20$ cases for the relations of these cycles which may lead to the existence of a subgraph isomorphic to the graph $C_5--C_5(--C_6)$ in $K(\alpha,\beta)$. It can be shown that all of these $20$ cases lead to contradictions and so, the graph $K(\alpha,\beta)$ contains no subgraph isomorphic to the graph $C_5--C_5(--C_6)$.

Each group with two generators $h_2$ and $h_3$ and two relations which are between the latter $20$ cases is a quotient of $B(1,k)$, for some integer $k$, is a cyclic group or is a solvable group.

\subsection{$\mathbf{C_4-C_6(--C_7--)(C_7-1)}$}
It can be seen that there are $176$ cases for the relations of a cycle $C_4$, two cycles $C_7$ and a cycle $C_6$ in this structure. By considering all groups with two generators $h_2$ and $h_3$ and four relations which are between these cases and by using GAP \cite{gap}, we see that all of these groups are finite and solvable. So, the graph $K(\alpha,\beta)$ contains no subgraph isomorphic to the graph $C_4-C_6(--C_7--)(C_7-1)$.

\subsection{$\mathbf{C_4-C_6(--C_7--)(--C_5-)}$}
It can be seen that there are $28$ cases for the relations of a cycle $C_4$, a cycle $C_6$, a cycle $C_7$ and a cycle $C_5$ in the graph $C_4-C_6(--C_7--)(--C_5-)$. Using GAP \cite{gap}, we see that all groups with two generators $h_2$ and $h_3$ and four relations which are between $24$ cases of these $28$ cases are finite and solvable, that is a contradiction with the assumptions. So, there are just $4$ cases for the relations of these cycles which may lead to the existence of a subgraph isomorphic to the graph $C_4-C_6(--C_7--)(--C_5-)$ in $K(\alpha,\beta)$. It can be shown that all of these $4$ cases lead to contradictions and so, the graph $K(\alpha,\beta)$ contains no subgraph isomorphic to the graph $C_4-C_6(--C_7--)(--C_5-)$.

Each group with two generators $h_2$ and $h_3$ and two relations which are between the latter $4$ cases is a cyclic group.

\subsection{$\mathbf{C_4-C_6(-C_6--)(-C_4-)}$}
It can be seen that there is no case for the relations of two cycles $C_4$ and two cycles $C_6$ in this structure. It means that the graph $K(\alpha,\beta)$ contains no subgraph isomorphic to the graph $C_4-C_6(-C_6--)(-C_4-)$.

\subsection{$\mathbf{C_4-C_6(-C_6--)(--C_5-)}$}
It can be seen that there are $22$ cases for the relations of a cycle $C_4$, two cycles $C_6$ and a cycle $C_5$ in this structure. By considering all groups with two generators $h_2$ and $h_3$ and four relations which are between these cases and by using GAP \cite{gap}, we see that all of these groups are finite and solvable. So, the graph $K(\alpha,\beta)$ contains no subgraph isomorphic to the graph $C_4-C_6(-C_6--)(--C_5-)$.

\subsection{$\mathbf{C_4-C_6(-C_6--)(C_6---)}$}
It can be seen that there are $66$ cases for the relations of a cycle $C_4$ and three cycles $C_6$ in the graph $C_4-C_6(-C_6--)(C_6---)$. Using GAP \cite{gap}, we see that all groups with two generators $h_2$ and $h_3$ and four relations which are between $62$ cases of these $66$ cases are finite and solvable, that is a contradiction with the assumptions. So, there are just $4$ cases for the relations of these cycles which may lead to the existence of a subgraph isomorphic to the graph $C_4-C_6(-C_6--)(C_6---)$ in $K(\alpha,\beta)$. It can be shown that all of these $4$ cases lead to contradictions and so, the graph $K(\alpha,\beta)$ contains no subgraph isomorphic to the graph $C_4-C_6(-C_6--)(C_6---)$.

Each group with two generators $h_2$ and $h_3$ and two relations which are between the latter $4$ cases is a quotient of $B(1,k)$, for some integer $k$.

\subsection{$\mathbf{C_4-C_6(-C_6--)(---C_4)}$}
It can be seen that there is no case for the relations of two cycles $C_4$ and two cycles $C_6$ in this structure. It means that the graph $K(\alpha,\beta)$ contains no subgraph isomorphic to the graph $C_4-C_6(-C_6--)(---C_4)$.

\subsection{$\mathbf{C_5-C_5(--C_6--)}$}
It can be seen that there are $440$ cases for the relations of two cycles $C_5$ and a cycle $C_6$ in the graph $C_5-C_5(--C_6--)$. Using GAP \cite{gap}, we see that all groups with two generators $h_2$ and $h_3$ and three relations which are between $404$ cases of these $440$ cases are finite and solvable, that is a contradiction with the assumptions. So, there are just $36$ cases for the relations of these cycles which may lead to the existence of a subgraph isomorphic to the graph $C_5-C_5(--C_6--)$ in $K(\alpha,\beta)$. It can be shown that all of these $36$ cases lead to contradictions and so, the graph $K(\alpha,\beta)$ contains no subgraph isomorphic to the graph $C_5-C_5(--C_6--)$.

Each group with two generators $h_2$ and $h_3$ and two relations which are between the latter $36$ cases is a quotient of $B(1,k)$, for some integer $k$, has a torsion element or is a cyclic group.

\subsection{$\mathbf{C_5-C_5(-C_6--)(C_6---)}$}
It can be seen that there are $56$ cases for the relations of two cycles $C_5$ and two cycles $C_6$ in the graph $C_5-C_5(-C_6--)(C_6---)$. By considering all groups with two generators $h_2$ and $h_3$ and four relations which are between $54$ cases of these $56$ cases are finite or solvable, that is a contradiction with the assumptions. So, there are just $2$ cases for the relations of these cycles which may lead to the existence of a subgraph isomorphic to the graph $C_5-C_5(-C_6--)(C_6---)$ in $K(\alpha,\beta)$. It can be shown that all of these $2$ cases lead to contradictions and so, the graph $K(\alpha,\beta)$ contains no subgraph isomorphic to the graph $C_5-C_5(-C_6--)(C_6---)$.

Each group with two generators $h_2$ and $h_3$ and two relations which are between the latter $2$ cases is a cyclic group.

\subsection{$\mathbf{C_5-C_5(-C_6--)(--C_6-1)}$}
It can be seen that there are $56$ cases for the relations of two cycles $C_5$ and two cycles $C_6$ in the graph $C_5-C_5(-C_6--)(--C_6-1)$. By considering all groups with two generators $h_2$ and $h_3$ and four relations which are between these cases and by using GAP \cite{gap}, we see that all of these groups are finite and solvable. So, the graph $K(\alpha,\beta)$ contains no subgraph isomorphic to the graph $C_5-C_5(-C_6--)(--C_6-1)$.

\subsection{$\mathbf{C_5-C_5(-C_6--)(-C_5--)}$}
It can be seen that there are $14$ cases for the relations of three cycles $C_5$ and one cycle $C_6$ in the graph $C_5-C_5(-C_6--)(-C_5--)$. By considering all groups with two generators $h_2$ and $h_3$ and four relations which are between these cases and by using GAP \cite{gap}, we see that all of these groups are finite and solvable. So, the graph $K(\alpha,\beta)$ contains no subgraph isomorphic to the graph $C_5-C_5(-C_6--)(-C_5--)$.

\subsection{$\mathbf{C_6---C_6(C_6---C_6)}$}
It can be seen that there are $46$ cases for the relations of four $C_6$ cycles in this structure. Using GAP \cite{gap}, we see that all groups with two generators $h_2$ and $h_3$ and four relations which are between $30$ cases of these $46$ cases are finite and solvable, or just finite, that is a contradiction with the assumptions. So, there are just $16$ cases for the relations of these cycles which may lead to the existence of a subgraph isomorphic to the graph $C_6---C_6(C_6---C_6)$ in $K(\alpha,\beta)$. It can be shown that these $16$ cases lead to contradictions and so, the graph $K(\alpha,\beta)$ contains no subgraph isomorphic to the graph $C_6---C_6(C_6---C_6)$.

Each group with two generators $h_2$ and $h_3$ and two relations which are between the latter $16$ cases is a quotient of $B(1,k)$, for some integer $k$, or has a torsion element.

\subsection{$\mathbf{C_6---C_6(C_6)(C_6)(C_6)}$}
It can be seen that there are $10$ cases for the relations of five $C_6$ cycles in this structure. Using GAP \cite{gap}, we see that all groups with two generators $h_2$ and $h_3$ and five relations which are between $6$ cases of these $10$ cases are finite and solvable, or just finite, that is a contradiction with the assumptions. So, there are just $4$ cases for the relations of these cycles which may lead to the existence of a subgraph isomorphic to the graph $C_6---C_6(C_6)(C_6)(C_6)$ in $K(\alpha,\beta)$. It can be shown that these $4$ cases lead to contradictions and so, the graph $K(\alpha,\beta)$ contains no subgraph isomorphic to the graph $C_6---C_6(C_6)(C_6)(C_6)$.

Each group with two generators $h_2$ and $h_3$ and two relations which are between the latter $4$ cases is a quotient of $B(1,k)$, for some integer $k$, or is a cyclic group.

\subsection{$\mathbf{C_5(--C_6--)C_5(---C_6)}$}
It can be seen that there are $134$ cases for the relations of two cycles $C_5$ and two cycles $C_6$ in this structure. Using GAP \cite{gap}, we see that all groups with two generators $h_2$ and $h_3$ and four relations which are between $130$ cases of these $134$ cases are finite and solvable, or just finite, that is a contradiction with the assumptions. So, there are just $4$ cases for the relations of these cycles which may lead to the existence of a subgraph isomorphic to the graph $C_5(--C_6--)C_5(---C_6)$ in $K(\alpha,\beta)$. It can be shown that these $4$ cases lead to contradictions and so, the graph $K(\alpha,\beta)$ contains no subgraph isomorphic to the graph $C_5(--C_6--)C_5(---C_6)$.

Each group with two generators $h_2$ and $h_3$ and two relations which are between the latter $4$ cases is a quotient of $B(1,k)$, for some integer $k$, or is a cyclic group.

\subsection{$\mathbf{C_6--C_6(C_6--C_6)}$}
It can be seen that there are $5119$ cases for the relations of four $C_6$ cycles in this structure. Using GAP \cite{gap}, we see that all groups with two generators $h_2$ and $h_3$ and four relations which are between $4983$ cases of these $5119$ cases are finite and solvable, or just finite, that is a contradiction with the assumptions. So, there are just $136$ cases for the relations of these cycles which may lead to the existence of a subgraph isomorphic to the graph $C_6--C_6(C_6--C_6)$ in $K(\alpha,\beta)$. It can be shown that these $136$ cases lead to contradictions and so, the graph $K(\alpha,\beta)$ contains no subgraph isomorphic to the graph $C_6--C_6(C_6--C_6)$.

Each group with two generators $h_2$ and $h_3$ and two relations which are between the latter $136$ cases is a quotient of $B(1,k)$, for some integer $k$, has a torsion element or is a cyclic group.

\subsection{$\mathbf{C_6---C_6(C_6--C_6)}$}
It can be seen that there are $1594$ cases for the relations of four $C_6$ cycles in this structure. Using GAP \cite{gap}, we see that all groups with two generators $h_2$ and $h_3$ and four relations which are between $1446$ cases of these $1594$ cases are finite and solvable, or just finite, that is a contradiction with the assumptions. So, there are just $148$ cases for the relations of these cycles which may lead to the existence of a subgraph isomorphic to the graph $C_6---C_6(C_6--C_6)$ in $K(\alpha,\beta)$. It can be shown that these $148$ cases lead to contradictions and so, the graph $K(\alpha,\beta)$ contains no subgraph isomorphic to the graph $C_6---C_6(C_6--C_6)$.

Each group with two generators $h_2$ and $h_3$ and two relations which are between the latter $148$ cases is a quotient of $B(1,k)$, for some integer $k$, has a torsion element or is a cyclic group.

\subsection{$\mathbf{C_6---C_6(-C_5-)}$}
It can be seen that there are $1482$ cases for the relations of two cycles $C_6$ and a cycle $C_5$ in the graph $C_6---C_6(-C_5-)$. Using GAP \cite{gap}, we see that all groups with two generators $h_2$ and $h_3$ and three relations which are between $1358$ cases of these $1482$ cases are finite and solvable, that is a contradiction with the assumptions. So, there are just $124$ cases for the relations of these cycles which may lead to the existence of a subgraph isomorphic to the graph $C_6---C_6(-C_5-)$ in $K(\alpha,\beta)$. It can be shown that all of these $124$ cases lead to contradictions and so, the graph $K(\alpha,\beta)$ contains no subgraph isomorphic to the graph $C_6---C_6(-C_5-)$.

Each group with two generators $h_2$ and $h_3$ and two relations which are between the latter $124$ cases is a quotient of $B(1,k)$, for some integer $k$, has a torsion element or is a cyclic group.

\subsection{$\mathbf{C_4-C_6(--C_7--)(---C_6)}$}
It can be seen that there are $124$ cases for the relations of a cycle $C_4$, two cycles $C_6$ and a cycle $C_7$ in the graph $C_4-C_6(--C_7--)(---C_6)$. Using GAP \cite{gap}, we see that all groups with two generators $h_2$ and $h_3$ and four relations which are between $112$ cases of these $124$ cases are finite and solvable, that is a contradiction with the assumptions. So, there are just $12$ cases for the relations of these cycles which may lead to the existence of a subgraph isomorphic to the graph $C_4-C_6(--C_7--)(---C_6)$ in $K(\alpha,\beta)$. It can be shown that all of these $12$ cases lead to contradictions and so, the graph $K(\alpha,\beta)$ contains no subgraph isomorphic to the graph $C_4-C_6(--C_7--)(---C_6)$.

Each group with two generators $h_2$ and $h_3$ and two relations which are between the latter $12$ cases is a quotient of $B(1,k)$, for some integer $k$.

\subsection{$\mathbf{C_4-C_6(--C_7--)(C_4)(C_4)}$}
It can be seen that there are $8$ cases for the relations of three cycles $C_4$, a cycle $C_7$ and a cycle $C_6$ in this structure. By considering all groups with two generators $h_2$ and $h_3$ and five relations which are between these cases and by using GAP \cite{gap}, we see that all of these groups are finite and solvable. So, the graph $K(\alpha,\beta)$ contains no subgraph isomorphic to the graph $C_4-C_6(--C_7--)(C_4)(C_4)$.

\subsection{$\mathbf{C_6---C_6(-C_5--)}$}
It can be seen that there are $418$ cases for the relations of two cycles $C_6$ and a cycle $C_5$ in the graph $C_6---C_6(-C_5--)$. Using GAP \cite{gap}, we see that all groups with two generators $h_2$ and $h_3$ and three relations which are between $358$ cases of these $418$ cases are finite and solvable, that is a contradiction with the assumptions. So, there are just $60$ cases for the relations of these cycles which may lead to the existence of a subgraph isomorphic to the graph $C_6---C_6(-C_5--)$ in $K(\alpha,\beta)$. It can be shown that all of these $60$ cases lead to contradictions and so, the graph $K(\alpha,\beta)$ contains no subgraph isomorphic to the graph $C_6---C_6(-C_5--)$.

Each group with two generators $h_2$ and $h_3$ and two relations which are between the latter $60$ cases is a quotient of $B(1,k)$, for some integer $k$, has a torsion element, is a cyclic group or is a solvable group.

\subsection{$\mathbf{C_6--C_6(--C_5-)(-C_5-)}$}
It can be seen that there are $62$ cases for the relations of two cycles $C_5$ and two cycles $C_6$ in this structure. Using GAP \cite{gap}, we see that all groups with two generators $h_2$ and $h_3$ and four relations which are between $56$ cases of these $62$ cases are finite and solvable, or just finite, that is a contradiction with the assumptions. So, there are just $6$ cases for the relations of these cycles which may lead to the existence of a subgraph isomorphic to the graph $C_6--C_6(--C_5-)(-C_5-)$ in $K(\alpha,\beta)$. It can be shown that these $6$ cases lead to contradictions and so, the graph $K(\alpha,\beta)$ contains no subgraph isomorphic to the graph $C_6--C_6(--C_5-)(-C_5-)$.

Each group with two generators $h_2$ and $h_3$ and two relations which are between the latter $6$ cases has a torsion element or is a cyclic group.

\subsection{$\mathbf{C_6--C_6(--C_5-)(C_6---)}$}
It can be seen that there are $76$ cases for the relations of a cycle $C_5$ and three cycles $C_6$ in the graph $C_6--C_6(--C_5-)(C_6---)$. Using GAP \cite{gap}, we see that all groups with two generators $h_2$ and $h_3$ and four relations which are between $64$ cases of these $76$ cases are finite or solvable, that is a contradiction with the assumptions. So, there are just $12$ cases for the relations of these cycles which may lead to the existence of a subgraph isomorphic to the graph $C_6--C_6(--C_5-)(C_6---)$ in $K(\alpha,\beta)$. It can be shown that all of these $12$ cases lead to contradictions and so, the graph $K(\alpha,\beta)$ contains no subgraph isomorphic to the graph $C_6--C_6(--C_5-)(C_6---)$.

Each group with two generators $h_2$ and $h_3$ and two relations which are between the latter $12$ cases is a quotient of $B(1,k)$, for some integer $k$, has a torsion element or is a cyclic group.

\subsection{$\mathbf{C_5(--C_6--)C_5(C_6)}$}
It can be seen that there are $120$ cases for the relations of two cycles $C_5$ and two cycles $C_6$ in this structure. Using GAP \cite{gap}, we see that all groups with two generators $h_2$ and $h_3$ and four relations which are between $104$ cases of these $120$ cases are finite and solvable, or just finite, that is a contradiction with the assumptions. So, there are just $16$ cases for the relations of these cycles which may lead to the existence of a subgraph isomorphic to the graph $C_5(--C_6--)C_5(C_6)$ in $K(\alpha,\beta)$. It can be shown that these $16$ cases lead to contradictions and so, the graph $K(\alpha,\beta)$ contains no subgraph isomorphic to the graph $C_5(--C_6--)C_5(C_6)$.

Each group with two generators $h_2$ and $h_3$ and two relations which are between the latter $16$ cases is a quotient of $B(1,k)$, for some integer $k$, has a torsion element or is a cyclic group.

\subsection{$\mathbf{C_5(--C_6--)C_5(C_7)}$}
It can be seen that there are $248$ cases for the relations of two cycles $C_5$, a cycle $C_6$ and a cycle $C_7$ in this structure. Using GAP \cite{gap}, we see that all groups with two generators $h_2$ and $h_3$ and four relations which are between $220$ cases of these $248$ cases are finite and solvable, or just finite, that is a contradiction with the assumptions. So, there are just $28$ cases for the relations of these cycles which may lead to the existence of a subgraph isomorphic to the graph $C_5(--C_6--)C_5(C_7)$ in $K(\alpha,\beta)$. It can be shown that these $28$ cases lead to contradictions and so, the graph $K(\alpha,\beta)$ contains no subgraph isomorphic to the graph $C_5(--C_6--)C_5(C_7)$.

Each group with two generators $h_2$ and $h_3$ and two relations which are between the latter $28$ cases is a quotient of $B(1,k)$, for some integer $k$, has a torsion element or is a cyclic group.

\subsection{$\mathbf{C_5-C_5(--C_7--)(--C_5)}$}
It can be seen that there are $394$ cases for the relations of a cycle $C_7$ and three cycles $C_5$ in the graph $C_5-C_5(--C_7--)(--C_5)$. Using GAP \cite{gap}, we see that all groups with two generators $h_2$ and $h_3$ and four relations which are between $352$ cases of these $394$ cases are finite or solvable, that is a contradiction with the assumptions. So, there are just $42$ cases for the relations of these cycles which may lead to the existence of a subgraph isomorphic to the graph $C_5-C_5(--C_7--)(--C_5)$ in $K(\alpha,\beta)$. It can be shown that all of these $42$ cases lead to contradictions and so, the graph $K(\alpha,\beta)$ contains no subgraph isomorphic to the graph $C_5-C_5(--C_7--)(--C_5)$.

Each group with two generators $h_2$ and $h_3$ and two relations which are between the latter $42$ cases is a quotient of $B(1,k)$, for some integer $k$, has a torsion element, is a cyclic group or is a solvable group.

\subsection{$\mathbf{C_5-C_5(--C_7--)(-C_5-)}$}
It can be seen that there are $138$ cases for the relations of a cycle $C_7$ and three cycles $C_5$ in the graph $C_5-C_5(--C_7--)(-C_5-)$. Using GAP \cite{gap}, we see that all groups with two generators $h_2$ and $h_3$ and four relations which are between $132$ cases of these $138$ cases are finite or solvable, that is a contradiction with the assumptions. So, there are just $6$ cases for the relations of these cycles which may lead to the existence of a subgraph isomorphic to the graph $C_5-C_5(--C_7--)(-C_5-)$ in $K(\alpha,\beta)$. It can be shown that all of these $6$ cases lead to contradictions and so, the graph $K(\alpha,\beta)$ contains no subgraph isomorphic to the graph $C_5-C_5(--C_7--)(-C_5-)$.

Each group with two generators $h_2$ and $h_3$ and two relations which are between the latter $6$ cases has a torsion element, is a cyclic group or is a solvable group.

\subsection{$\mathbf{C_5-C_5(-C_6--)(--C_6-2)}$}
It can be seen that there are $22$ cases for the relations of two cycles $C_5$ and two cycles $C_6$ in the graph $C_5-C_5(-C_6--)(--C_6-2)$. By considering all groups with two generators $h_2$ and $h_3$ and four relations which are between these cases and by using GAP \cite{gap}, we see that all of these groups are finite and solvable. So, the graph $K(\alpha,\beta)$ contains no subgraph isomorphic to the graph $C_5-C_5(-C_6--)(--C_6-2)$.

\subsection{$\mathbf{C_4-C_4(-C_7-)(C_4)}$}
It can be seen that there are $32$ cases for the relations of a cycle $C_7$ and three cycles $C_4$ in the graph $C_4-C_4(-C_7-)(C_4)$. By considering all groups with two generators $h_2$ and $h_3$ and four relations which are between these cases and by using GAP \cite{gap}, we see that all of these groups are finite and solvable. So, the graph $K(\alpha,\beta)$ contains no subgraph isomorphic to the graph $C_4-C_4(-C_7-)(C_4)$.

\subsection{$\mathbf{C_5--C_5(-C_5--)}$}
It can be seen that there are $64$ cases for the relations of three cycles $C_5$ in the graph $C_5--C_5(-C_5--)$. Using GAP \cite{gap}, we see that all groups with two generators $h_2$ and $h_3$ and three relations which are between $58$ cases of these $64$ cases are finite and solvable, that is a contradiction with the assumptions. So, there are just $6$ cases for the relations of these cycles which may lead to the existence of a subgraph isomorphic to the graph $C_5--C_5(-C_5--)$ in $K(\alpha,\beta)$. It can be shown that all of these $6$ cases lead to contradictions and so, the graph $K(\alpha,\beta)$ contains no subgraph isomorphic to the graph $C_5--C_5(-C_5--)$.

Each group with two generators $h_2$ and $h_3$ and two relations which are between the latter $6$ cases is a quotient of $B(1,k)$, for some integer $k$ or is a cyclic group.

\subsection{$\mathbf{C_6---C_6(-C_4)}$}
It can be seen that there are $420$ cases for the relations of two cycles $C_6$ and a cycle $C_4$ in the graph $C_6---C_6(-C_4)$. Using GAP \cite{gap}, we see that all groups with two generators $h_2$ and $h_3$ and three relations which are between $398$ cases of these $420$ cases are finite or solvable, that is a contradiction with the assumptions. So, there are just $22$ cases for the relations of these cycles which may lead to the existence of a subgraph isomorphic to the graph $C_6---C_6(-C_4)$ in $K(\alpha,\beta)$. It can be shown that all of these $22$ cases lead to contradictions and so, the graph $K(\alpha,\beta)$ contains no subgraph isomorphic to the graph $C_6---C_6(-C_4)$.

Each group with two generators $h_2$ and $h_3$ and two relations which are between the latter $22$ cases is a cyclic group.

\subsection{$\mathbf{C_6--C_6(C_4)}$}
It can be seen that there are $279$ cases for the relations of two cycles $C_6$ and a cycle $C_4$ in the graph $C_6--C_6(C_4)$. Using GAP \cite{gap}, we see that all groups with two generators $h_2$ and $h_3$ and three relations which are between $268$ cases of these $279$ cases are finite or solvable, that is a contradiction with the assumptions. So, there are just $11$ cases for the relations of these cycles which may lead to the existence of a subgraph isomorphic to the graph $C_6--C_6(C_4)$ in $K(\alpha,\beta)$. It can be shown that all of these $11$ cases lead to contradictions and so, the graph $K(\alpha,\beta)$ contains no subgraph isomorphic to the graph $C_6--C_6(C_4)$.

Each group with two generators $h_2$ and $h_3$ and two relations which are between the latter $11$ cases is a quotient of $B(1,k)$, for some integer $k$, or is a cyclic group.

\subsection{$\mathbf{C_4-C_6(-C_4)(-C_4)}$}
It can be seen that there are $36$ cases for the relations of a cycle $C_6$ and three cycles $C_4$ in the graph $C_4-C_6(-C_4)(-C_4)$. By considering all groups with two generators $h_2$ and $h_3$ and four relations which are between these cases and by using GAP \cite{gap}, we see that all of these groups are solvable. So, the graph $K(\alpha,\beta)$ contains no subgraph isomorphic to the graph $C_4-C_6(-C_4)(-C_4)$.

\subsection{$\mathbf{C_4-C_6(--C_7--)(-C_5-)}$}
It can be seen that there are $62$ cases for the relations of a cycle $C_4$, a cycle $C_6$, a cycle $C_7$ and a cycle $C_5$ in the graph $C_4-C_6(--C_7--)(-C_5-)$. Using GAP \cite{gap}, we see that all groups with two generators $h_2$ and $h_3$ and four relations which are between $58$ cases of these $62$ cases are solvable, that is a contradiction with the assumptions. So, there are just $4$ cases for the relations of these cycles which may lead to the existence of a subgraph isomorphic to the graph $C_4-C_6(--C_7--)(-C_5-)$ in $K(\alpha,\beta)$. It can be shown that all of these $4$ cases lead to contradictions and so, the graph $K(\alpha,\beta)$ contains no subgraph isomorphic to the graph $C_4-C_6(--C_7--)(-C_5-)$.

Each group with two generators $h_2$ and $h_3$ and two relations which are between the latter $4$ cases is a cyclic group.

\subsection{$\mathbf{C_5-C_5(-C_6--)(--C_5-)}$}
It can be seen that there are $14$ cases for the relations of three cycles $C_5$ and a cycle $C_6$ in the graph $C_5-C_5(-C_6--)(--C_5-)$. By considering all groups with two generators $h_2$ and $h_3$ and four relations which are between these cases and by using GAP \cite{gap}, we see that all of these groups are solvable. So, the graph $K(\alpha,\beta)$ contains no subgraph isomorphic to the graph $C_5-C_5(-C_6--)(--C_5-)$.

\subsection{$\mathbf{C_4-C_6(--C_7--)(C_7-2)}$}
It can be seen that there are $168$ cases for the relations of a cycle $C_4$, a cycle $C_6$ and two cycles $C_7$ in the graph $C_4-C_6(--C_7--)(C_7-2)$. Using GAP \cite{gap}, we see that all groups with two generators $h_2$ and $h_3$ and four relations which are between $152$ cases of these $168$ cases are solvable, that is a contradiction with the assumptions. So, there are just $16$ cases for the relations of these cycles which may lead to the existence of a subgraph isomorphic to the graph $C_4-C_6(--C_7--)(C_7-2)$ in $K(\alpha,\beta)$. It can be shown that all of these $16$ cases lead to contradictions and so, the graph $K(\alpha,\beta)$ contains no subgraph isomorphic to the graph $C_4-C_6(--C_7--)(C_7-2)$.

Each group with two generators $h_2$ and $h_3$ and two relations which are between the latter $16$ cases is a quotient of $B(1,k)$, for some integer $k$, or is a cyclic group.

In the following theorem, we summarize our results about forbidden subgraphs of the graph $K(\alpha,\beta)$. 
\begin{thm}\label{forbiddens}
The Kaplansky graph $K(\alpha,\beta)$ is a triangle-free graph which contains no subgraph isomorphic to one of the $46$ other graphs in Table \ref{tab-forbiddens}.
\end{thm}

\begin{rem}\label{r-label}\rm
In $K_{\mathbb{F}}(\alpha,\beta)$, for all $g_i,g_j\in supp(\beta)$, $g_i\sim g_j$ if and only if $ag_i=bg_j$ for some $a,b\in supp(\alpha)$. Let $e$ be the edge between such latter vertices. Suppose that we give an orientation and a label in $L=\{(1,h_2),(1,h_3),(h_2,h_3)\}$ to $e$ such as follows:
\begin{enumerate}
\item
If $(a,b)\in L$, then $e$ is labelled by $(a,b)$ and oriented from $g_i$ to $g_j$.
\item
If $(b,a)\in L$, then $e$ is labelled by $(b,a)$ and oriented from $g_j$ to $g_i$.
\end{enumerate}
Therefore, $K_{\mathbb{F}}(\alpha,\beta)$ can be considered as a directed graph in which all edges are labelled by $L$ with the above method. 

\begin{figure}[ht]
\psscalebox{0.9 0.9} 
{
\begin{tikzpicture}
[every node/.style={inner sep=0pt}]
\node (1) [circle, minimum size=12.5pt, fill=black, line width=1.25pt, draw=black] at (112.5pt, -75.0pt)  {};
\node (2) [circle, minimum size=12.5pt, fill=black, line width=1.25pt, draw=black] at (187.5pt, -75.0pt)  {};
\draw [line width=1.25, ->, color=black] (1) to  (2);
\node at (112.5pt, -60.0pt) {\textcolor{black}{$g_i$}};
\node at (187.5pt, -60.0pt) {\textcolor{black}{$g_j$}};
\node at (148.75pt, -65.625pt) {\textcolor{black}{$(a,b)$}};
\end{tikzpicture}
}
\caption{An edge in $K_{\mathbb{F}}(\alpha,\beta)$}\label{f-edge}
\end{figure}
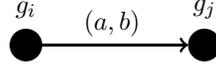

Let $e$ be an edge in $K_{\mathbb{F}}(\alpha,\beta)$ as Figure \ref{f-edge}, where $(a,b)\in L$. Suppose that when we traverse $e$ from $g_i$ to $g_j$ (or from $g_j$ to $g_i$), the word $a^{-1}b$ (or $b^{-1}a$,  respectively) is corresponded to this edge. Then $G$ is the group $\langle h_2,h_3 | \text{ words of closed paths in } K_{\mathbb{F}}(\alpha,\beta)\rangle$.

Let $\mathbb{F}=\mathbb{F}_2$. It is easy to see that for each vertex $g_i$ of $K_{\mathbb{F}_2}(\alpha,\beta)$, none of the cases in Figure \ref{f-label} can be happened.

\begin{figure}[h]
\psscalebox{0.7 0.7} 
{
\begin{tikzpicture}
[every node/.style={inner sep=0pt}]
\node (1) [circle, minimum size=12.5pt, fill=black, line width=1.25pt, draw=black] at (112.5pt, -50.0pt)  {};
\node (2) [circle, minimum size=12.5pt, fill=black, line width=1.25pt, draw=black] at (62.5pt, -100.0pt)  {};
\node (3) [circle, minimum size=12.5pt, fill=black, line width=1.25pt, draw=black] at (162.5pt, -100.0pt)  {};
\draw [line width=1.25, ->, color=black] (2) to  (1);
\draw [line width=1.25, ->, color=black] (3) to  (1);
\node at (120.375pt, -38.125pt) {\textcolor{black}{$g_i$}};
\node at (81.875pt, -67.5pt) [rotate=45] {\textcolor{black}{$(1,h_2)$}};
\node at (145.625pt, -70.0pt) [rotate=315] {\textcolor{black}{$(1,h_2)$}};
\end{tikzpicture}
\quad
\begin{tikzpicture}
[every node/.style={inner sep=0pt}]
\node (1) [circle, minimum size=12.5pt, fill=black, line width=1.25pt, draw=black] at (112.5pt, -50.0pt)  {};
\node (2) [circle, minimum size=12.5pt, fill=black, line width=1.25pt, draw=black] at (62.5pt, -100.0pt)  {};
\node (3) [circle, minimum size=12.5pt, fill=black, line width=1.25pt, draw=black] at (162.5pt, -100.0pt)  {};
\draw [line width=1.25, ->, color=black] (2) to  (1);
\draw [line width=1.25, ->, color=black] (3) to  (1);
\node at (120.375pt, -38.125pt) {\textcolor{black}{$g_i$}};
\node at (81.875pt, -67.5pt) [rotate=45] {\textcolor{black}{$(1,h_3)$}};
\node at (145.625pt, -70.0pt) [rotate=315] {\textcolor{black}{$(1,h_3)$}};
\end{tikzpicture}
\quad
\begin{tikzpicture}
[every node/.style={inner sep=0pt}]
\node (1) [circle, minimum size=12.5pt, fill=black, line width=1.25pt, draw=black] at (112.5pt, -50.0pt)  {};
\node (2) [circle, minimum size=12.5pt, fill=black, line width=1.25pt, draw=black] at (62.5pt, -100.0pt)  {};
\node (3) [circle, minimum size=12.5pt, fill=black, line width=1.25pt, draw=black] at (162.5pt, -100.0pt)  {};
\draw [line width=1.25, ->, color=black] (2) to  (1);
\draw [line width=1.25, ->, color=black] (3) to  (1);
\node at (120.375pt, -38.125pt) {\textcolor{black}{$g_i$}};
\node at (81.875pt, -67.5pt) [rotate=45] {\textcolor{black}{$(1,h_3)$}};
\node at (144.375pt, -68.75pt) [rotate=315] {\textcolor{black}{$(h_2,h_3)$}};
\end{tikzpicture}
\quad
\begin{tikzpicture}
[every node/.style={inner sep=0pt}]
\node (1) [circle, minimum size=12.5pt, fill=black, line width=1.25pt, draw=black] at (112.5pt, -50.0pt)  {};
\node (2) [circle, minimum size=12.5pt, fill=black, line width=1.25pt, draw=black] at (62.5pt, -100.0pt)  {};
\node (3) [circle, minimum size=12.5pt, fill=black, line width=1.25pt, draw=black] at (162.5pt, -100.0pt)  {};
\draw [line width=1.25, ->, color=black] (2) to  (1);
\draw [line width=1.25, ->, color=black] (3) to  (1);
\node at (120.375pt, -38.125pt) {\textcolor{black}{$g_i$}};
\node at (81.25pt, -67.5pt) [rotate=45] {\textcolor{black}{$(h_2,h_3)$}};
\node at (144.375pt, -68.75pt) [rotate=315] {\textcolor{black}{$(h_2,h_3)$}};
\end{tikzpicture}
}
\end{figure}
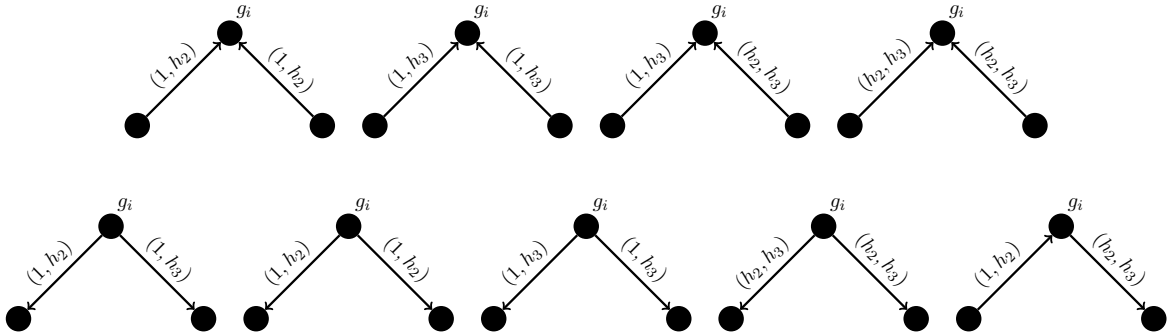
\begin{figure}[h]
\psscalebox{0.7 0.7} 
{
\begin{tikzpicture}
[every node/.style={inner sep=0pt}]
\node (1) [circle, minimum size=12.5pt, fill=black, line width=1.25pt, draw=black] at (112.5pt, -50.0pt)  {};
\node (2) [circle, minimum size=12.5pt, fill=black, line width=1.25pt, draw=black] at (62.5pt, -100.0pt)  {};
\node (3) [circle, minimum size=12.5pt, fill=black, line width=1.25pt, draw=black] at (162.5pt, -100.0pt)  {};
\draw [line width=1.25, ->, color=black] (1) to  (2);
\draw [line width=1.25, ->, color=black] (1) to  (3);
\node at (120.375pt, -38.125pt) {\textcolor{black}{$g_i$}};
\node at (79.375pt, -70.0pt) [rotate=45] {\textcolor{black}{$(1,h_2)$}};
\node at (143.75pt, -68.125pt) [rotate=315] {\textcolor{black}{$(1,h_3)$}};
\end{tikzpicture}
\quad
\begin{tikzpicture}
[every node/.style={inner sep=0pt}]
\node (1) [circle, minimum size=12.5pt, fill=black, line width=1.25pt, draw=black] at (112.5pt, -50.0pt)  {};
\node (2) [circle, minimum size=12.5pt, fill=black, line width=1.25pt, draw=black] at (62.5pt, -100.0pt)  {};
\node (3) [circle, minimum size=12.5pt, fill=black, line width=1.25pt, draw=black] at (162.5pt, -100.0pt)  {};
\draw [line width=1.25, ->, color=black] (1) to  (2);
\draw [line width=1.25, ->, color=black] (1) to  (3);
\node at (120.375pt, -38.125pt) {\textcolor{black}{$g_i$}};
\node at (79.375pt, -70.0pt) [rotate=45] {\textcolor{black}{$(1,h_2)$}};
\node at (143.75pt, -68.125pt) [rotate=315] {\textcolor{black}{$(1,h_2)$}};
\end{tikzpicture}
\quad
\begin{tikzpicture}
[every node/.style={inner sep=0pt}]
\node (1) [circle, minimum size=12.5pt, fill=black, line width=1.25pt, draw=black] at (112.5pt, -50.0pt)  {};
\node (2) [circle, minimum size=12.5pt, fill=black, line width=1.25pt, draw=black] at (62.5pt, -100.0pt)  {};
\node (3) [circle, minimum size=12.5pt, fill=black, line width=1.25pt, draw=black] at (162.5pt, -100.0pt)  {};
\draw [line width=1.25, ->, color=black] (1) to  (2);
\draw [line width=1.25, ->, color=black] (1) to  (3);
\node at (120.375pt, -38.125pt) {\textcolor{black}{$g_i$}};
\node at (79.375pt, -70.0pt) [rotate=45] {\textcolor{black}{$(1,h_3)$}};
\node at (143.75pt, -68.125pt) [rotate=315] {\textcolor{black}{$(1,h_3)$}};
\end{tikzpicture}
\quad
\begin{tikzpicture}
[every node/.style={inner sep=0pt}]
\node (1) [circle, minimum size=12.5pt, fill=black, line width=1.25pt, draw=black] at (112.5pt, -50.0pt)  {};
\node (2) [circle, minimum size=12.5pt, fill=black, line width=1.25pt, draw=black] at (62.5pt, -100.0pt)  {};
\node (3) [circle, minimum size=12.5pt, fill=black, line width=1.25pt, draw=black] at (162.5pt, -100.0pt)  {};
\draw [line width=1.25, ->, color=black] (1) to  (2);
\draw [line width=1.25, ->, color=black] (1) to  (3);
\node at (120.375pt, -38.125pt) {\textcolor{black}{$g_i$}};
\node at (79.375pt, -70.0pt) [rotate=45] {\textcolor{black}{$(h_2,h_3)$}};
\node at (143.75pt, -68.125pt) [rotate=315] {\textcolor{black}{$(h_2,h_3)$}};
\end{tikzpicture}
\quad
\begin{tikzpicture}
[every node/.style={inner sep=0pt}]
\node (1) [circle, minimum size=12.5pt, fill=black, line width=1.25pt, draw=black] at (112.5pt, -50.0pt)  {};
\node (2) [circle, minimum size=12.5pt, fill=black, line width=1.25pt, draw=black] at (62.5pt, -100.0pt)  {};
\node (3) [circle, minimum size=12.5pt, fill=black, line width=1.25pt, draw=black] at (162.5pt, -100.0pt)  {};
\draw [line width=1.25, ->, color=black] (2) to  (1);
\draw [line width=1.25, ->, color=black] (1) to  (3);
\node at (120.375pt, -38.125pt) {\textcolor{black}{$g_i$}};
\node at (78.125pt, -70.625pt) [rotate=45] {\textcolor{black}{$(1,h_2)$}};
\node at (143.75pt, -68.125pt) [rotate=315] {\textcolor{black}{$(h_2,h_3)$}};
\end{tikzpicture}
}
\caption{Impossible cases for labelling the edges of each vertex $g_i$ in $K_{\mathbb{F}_2}(\alpha,\beta)$}\label{f-label}
\end{figure}
\end{rem}

\begin{rem}{\rm
Suppose that $\Gamma$ is a directed graph in which all edges are labelled by $L$. Let $e$ be an edge in $\Gamma$ as Figure \ref{f-edge}, where $(a,b)\in L$. Suppose that when we traverse $e$ from $g_i$ to $g_j$ (or from $g_j$ to $g_i$), the word $a^{-1}b$ (or $b^{-1}a$,  respectively) is corresponded to this edge. Let $G(\Gamma):=\langle h_2,h_3 | \text{ words of closed paths in } \Gamma\rangle$. 
Suppose that the ground graph of $\Gamma$ contains a subgraph isomorphic to one of the graphs in Table \ref{tab-forbiddens} and for each vertex $g_i$ of such latter subgraph, none of the cases in Figure \ref{f-label} are happened. By Theorem \ref{forbiddens} and Remark \ref{r-label},  $G(\Gamma)$ has at least one of the following properties:
\begin{enumerate}
\item
It is a finite group,
\item 
It is an abelian group,
\item
It is a quotient of $BS(1,k)$ or $BS(k,1)$, where $k$ is an integer,
\item
It has a non-trivial torsion element,
\item 
It is a solvable group.
\end{enumerate}
}
\end{rem}

\section{\bf The possible number of vertices of $K(\alpha,\beta)$}\label{S3}
By Remark \ref{r-F2}, the number of vertices of $K(\alpha,\beta)$ must be an even positive integer $n\geq4$. Also by Theorem \ref{thm-graph}, the graph $K(\alpha,\beta)$ is a connected cubic triangle-free graph and we found $44$ other forbidden subgraphs of such graph in Section \ref{S2}. Furthermore in Section \ref{S-C4}, we found two graphs, namely $L_n$ and $M_n$, with $n$ vertices which are not isomorphic to the graph $K(\alpha,\beta)$. 

Using Sage Mathematics Software \cite{sage} and its package \textit{Nauty-geng}, all non-isomorphic connected cubic triangle-free graphs with the size of the vertex sets $n$ can be found. In this section by using Sage Mathematics Software \cite{sage}, we give some results about checking each of the mentioned forbidden subgraphs in all of the non-isomorphic connected cubic triangle-free graphs with the size of vertex sets $n \leq 20$. By using these results we show that $n$ must be greater than or equal to $20$. Also, some results about the case $n=20$ is given.

Table \ref{tab-graphs} lists all results about the number of non-isomorphic connected cubic triangle-free graphs with the size of vertex sets $n \leq 20$ which contain each of the forbidden subgraphs. The results in this table from top to bottom are presented in such a way that by checking each of the forbidden subgraphs in a row, the number of graphs containing these subgraph are omitted from the total number and the existence of the next forbidden subgraph is checked among the remaining ones.

\begin{scriptsize}
\begin{longtable}{|rl|l|l|l|l|l|l|l|l|l|}
\caption{Existence of the forbidden subgraphs in non-isomorphic connected cubic triangle-free graphs with the size of vertex sets $n \leq 20$}\label{tab-graphs}\\
\hline \multicolumn{1}{|r}{$ $} & \multicolumn{1}{l|}{$ $} & \multicolumn{1}{l|}{$n=4$} & \multicolumn{1}{l|}{$n=6$}
& \multicolumn{1}{l|}{$n=8$} & \multicolumn{1}{l|}{$n=10$} & \multicolumn{1}{l|}{$n=12$} & \multicolumn{1}{l|}{$n=14$}
& \multicolumn{1}{l|}{$n=16$} & \multicolumn{1}{l|}{$n=18$} & \multicolumn{1}{l|}{$n=20$}
\\\hline
\endfirsthead
\multicolumn{11}{c}%
{{\tablename\ \thetable{} -- continued from previous page}} \\
\hline \multicolumn{1}{|r}{$ $} & \multicolumn{1}{l|}{$ $} & \multicolumn{1}{l|}{$n=4$} & \multicolumn{1}{l|}{$n=6$}
& \multicolumn{1}{l|}{$n=8$} & \multicolumn{1}{l|}{$n=10$} & \multicolumn{1}{l|}{$n=12$} & \multicolumn{1}{l|}{$n=14$}
& \multicolumn{1}{l|}{$n=16$} & \multicolumn{1}{l|}{$n=18$} & \multicolumn{1}{l|}{$n=20$}
\\\hline
\endhead
\hline \multicolumn{11}{|r|}{{Continued on next page}} \\\hline
\endfoot
\hline 
\endlastfoot
$ $&$\text{Total}$&$0$&$1$&$2$&$6$&$22$&$110$&$792$&$7805$&$97546$\\\hline
$1)$&$K_{2,3}$&$0$&$1$&$0$&$1$&$4$&$22$&$144$&$1222$&$12991$\\
$2)$&$C_4--C_5$&$0$&$0$&$1$&$2$&$6$&$30$&$223$&$2161$&$25427$\\
$3)$&$C_4--C_6$&$0$&$0$&$1$&$1$&$6$&$31$&$223$&$2228$&$28080$\\
$4)$&$C_4-C_5(-C_5-)$&$0$&$0$&$0$&$0$&$2$&$6$&$40$&$319$&$3396$\\
$5)$&$C_4-C_5(-C_4-)$&$0$&$0$&$0$&$1$&$0$&$3$&$12$&$88$&$1123$\\
$6)$&$C_4-C_5(-C_6--)$&$0$&$0$&$0$&$0$&$0$&$4$&$42$&$389$&$4548$\\
$7)$&$C_4-C_5(-C_6-)$&$0$&$0$&$0$&$0$&$0$&$1$&$20$&$382$&$5661$\\
$8)$&$C_4-C_5(-C_7--)$&$0$&$0$&$0$&$0$&$0$&$1$&$10$&$176$&$3172$\\
$9)$&$C_5--C_5(--C_5)$&$0$&$0$&$0$&$1$&$2$&$3$&$18$&$157$&$1617$\\
$10)$&$C_5--C_5(--C_6)$&$0$&$0$&$0$&$0$&$0$&$4$&$32$&$291$&$4289$\\
$11)$&$C_4-C_6(--C_7--)(C_7-1)$&$0$&$0$&$0$&$0$&$1$&$0$&$9$&$64$&$446$\\
$12)$&$C_4-C_6(--C_7--)(--C_5-)$&$0$&$0$&$0$&$0$&$0$&$0$&$0$&$2$&$51$\\
$13)$&$C_4-C_6(-C_6--)(-C_4-)$&$0$&$0$&$0$&$0$&$1$&$0$&$1$&$5$&$35$\\
$14)$&$C_4-C_6(-C_6--)(--C_5-)$&$0$&$0$&$0$&$0$&$0$&$0$&$1$&$1$&$149$\\
$15)$&$C_4-C_6(-C_6--)(C_6---)$&$0$&$0$&$0$&$0$&$0$&$0$&$1$&$30$&$404$\\
$16)$&$C_4-C_6(-C_6--)(---C_4)$&$0$&$0$&$0$&$0$&$0$&$0$&$1$&$4$&$12$\\
$17)$&$C_5-C_5(--C_6--)$&$0$&$0$&$0$&$0$&$0$&$2$&$0$&$41$&$352$\\
$18)$&$C_5-C_5(-C_6--)(C_6---)$&$0$&$0$&$0$&$0$&$0$&$0$&$7$&$47$&$529$\\
$19)$&$C_5-C_5(-C_6--)(--C_6- 1)$&$0$&$0$&$0$&$0$&$0$&$0$&$1$&$31$&$249$\\
$20)$&$C_5-C_5(-C_6--)(-C_5--)$&$0$&$0$&$0$&$0$&$0$&$0$&$1$&$1$&$69$\\
$21)$&$C_6---C_6(C_6---C_6)$&$0$&$0$&$0$&$0$&$0$&$1$&$1$&$8$&$43$\\
$22)$&$C_6---C_6(C_6)(C_6)(C_6)$&$0$&$0$&$0$&$0$&$0$&$0$&$2$&$6$&$25$\\
$23)$&$C_5(--C_6--)C_5(---C_6)$&$0$&$0$&$0$&$0$&$0$&$0$&$1$&$16$&$374$\\
$24)$&$C_6--C_6(C_6--C_6)$&$0$&$0$&$0$&$0$&$0$&$0$&$0$&$29$&$438$\\
$25)$&$C_6---C_6(C_6--C_6)$&$0$&$0$&$0$&$0$&$0$&$0$&$0$&$20$&$505$\\
$26)$&$C_6---C_6(-C_5-)$&$0$&$0$&$0$&$0$&$0$&$0$&$0$&$27$&$721$\\
$27)$&$C_4-C_6(--C_7--)(---C_6)$&$0$&$0$&$0$&$0$&$0$&$0$&$0$&$8$&$257$\\
$28)$&$C_4-C_6(--C_7--)(C_4)(C_4)$&$0$&$0$&$0$&$0$&$0$&$0$&$0$&$2$&$2$\\
$29)$&$C_6---C_6(-C_5--)$&$0$&$0$&$0$&$0$&$0$&$0$&$0$&$5$&$66$\\
$30)$&$C_6--C_6(--C_5-)(-C_5-)$&$0$&$0$&$0$&$0$&$0$&$0$&$0$&$12$&$293$\\
$31)$&$C_6--C_6(--C_5-)(C_6---)$&$0$&$0$&$0$&$0$&$0$&$0$&$0$&$2$&$267$\\
$32)$&$C_5(--C_6--)C_5(C_6)$&$0$&$0$&$0$&$0$&$0$&$0$&$0$&$1$&$43$\\
$33)$&$C_5(--C_6--)C_5(C_7)$&$0$&$0$&$0$&$0$&$0$&$0$&$0$&$2$&$50$\\
$34)$&$C_5-C_5(--C_7--)(--C_5)$&$0$&$0$&$0$&$0$&$0$&$0$&$0$&$6$&$199$\\
$35)$&$C_5-C_5(--C_7--)(-C_5-)$&$0$&$0$&$0$&$0$&$0$&$0$&$0$&$2$&$69$\\
$36)$&$C_5-C_5(-C_6--)(--C_6-2)$&$0$&$0$&$0$&$0$&$0$&$0$&$0$&$2$&$114$\\
$37)$&$C_4-C_4(-C_7-)(C_4)$&$0$&$0$&$0$&$0$&$0$&$1$&$0$&$6$&$72$\\
$38)$&$C_5--C_5(-C_5--)$&$0$&$0$&$0$&$0$&$0$&$0$&$0$&$1$&$11$\\
$39)$&$C_6---C_6(-C_4)$&$0$&$0$&$0$&$0$&$0$&$0$&$0$&$3$&$94$\\
$40)$&$C_6--C_6(C_4)$&$0$&$0$&$0$&$0$&$0$&$0$&$0$&$2$&$67$\\
$41)$&$C_4-C_6(-C_4)(-C_4)$&$0$&$0$&$0$&$0$&$0$&$0$&$0$&$1$&$30$\\
$42)$&$C_4-C_6(--C_7--)(-C_5-)$&$0$&$0$&$0$&$0$&$0$&$0$&$0$&$1$&$26$\\
$43)$&$C_5-C_5(-C_6--)(--C_5-)$&$0$&$0$&$0$&$0$&$0$&$0$&$0$&$1$&$17$\\
$44)$&$C_4-C_6(--C_7--)(C_7-2)$&$0$&$0$&$0$&$0$&$0$&$0$&$0$&$1$&$41$\\
\hline
$ $&$\text{Isomorphic to }L_n$&$0$&$0$&$0$&$0$&$0$&$0$&$1$&$1$&$1$\\
$ $&$\text{Isomorphic to }M_n$&$0$&$0$&$0$&$0$&$0$&$1$&$1$&$1$&$1$\\
$ $&$\text{Remains}$&$0$&$0$&$0$&$0$&$0$&$0$&$0$&$0$&$1120$
\end{longtable}
\end{scriptsize}

The discussion above and the results of Table \ref{tab-graphs} are summarized in the following theorem. 
\begin{thm}\label{T2}
The vertex set size of $K(\alpha,\beta)$ must be greater than or equal to $20$. Furthermore, there are just $1120$ graphs with vertex set size equal to $20$ which may be isomorphic to $K(\alpha,\beta)$.
\end{thm}
\begin{cor}\label{maincor}
Let $\alpha$ and $\beta$ be non-zero elements of the group algebra of any torsion-free group over the field with two elements. If $|supp(\alpha)|=3$ and $\alpha \beta=0$ then $|supp(\beta)|\geq 20$. 
\end{cor}

\section{\bf Some results on the unit conjecture}\label{S-unit}
Throughout this section let $\mathbb{F}$ be an arbitrary field and $G$ be a torsion-free group and  $\gamma=\gamma_1 h_1+\gamma_2h_2+\gamma_3h_3 \in \mathbb{F}[G]$ such that $\vert supp(\gamma)\vert=3$. Suppose further that  $\gamma \delta =1$ for some  $\delta \in \mathbb{F}[G]$ and assume that  $n:=|supp(\delta)|$ is minimum with respect to the latter property. Let $\delta=\delta_1 g_1+\delta_2g_2+\cdots+\delta_ng_n$ and $A=\{1,2,3\}\times \{1,2,\ldots,n\}$. Since $\gamma \delta =1$, there must be at least one $(i,j)\in A$ such that $h_ig_j=1$. By renumbering, we may assume that $(i,j)=(1,1)$. Replacing $\gamma$ by $h_1^{-1}\gamma$ and $\delta$ by $\delta g_1^{-1}$ we may assume that $h_1=g_1=1$. So we may suppose that $1\in supp(\gamma)$ and $1\in supp(\delta)$. 

There is a partition $\pi$ of $A$ such that $(i, j)$ and $(i',j')$ belong to the same set of $\pi$ if and only if $h_ig_j=h_{i'}g_{j'}$ and because of the relation $\gamma \delta =1$, for all $E\in \pi$ we have 
\begin{equation}\label{e-u}
\displaystyle\sum_{(i,j)\in E}{\gamma_i\delta_j}=\left\{
\begin{array}{lr}
1 &(1,1)\in E\\
0 &(1,1)\notin E\\
\end{array} \right.
\end{equation}
Let $E_1$ be the set in $\pi$ which contains $(1,1)$. Obviously, $n\geq 2$. Abelian groups satisfy Conjecture \ref{conj-unit}. So, $G$ must be a nonabelian torsion-free group. Firstly in this section we show that $n\geq 4$. Then, we examine some small positive integers greater that $3$ as the possible values of $n$ and show that $n$ must be at least $8$. 

\subsection{The support of $\delta$ is of size at least $4$}
\begin{thm}[Corollary of \cite{kemp}]\label{kemperman}
Let $G$ be an arbitrary group and let $B$ and $C$ be finite non-empty subsets of $G$. Suppose that each non-identity element $g$ of $G$ has a finite or infinite order greater than or equal to $|B|+|C|-1$. Then $|BC|\geq |B|+|C|-1$.
\end{thm}
By Theorem \ref{kemperman}, $|supp(\gamma) supp(\delta)|\geq |supp(\gamma)|+|supp(\delta)|-1$ because $G$ is torsion-free. Hence, $3n\geq |supp(\gamma) supp(\delta)|\geq 2+n$.
\begin{enumerate}
\item
Let $n=2$. Then with the above discussion, $6\geq |supp(\gamma) supp(\delta)|\geq 4$ and so, there are at least $3$ sets different from $E_1$ in $\pi$, namely $E_2,E_3,E_4$. Since $\gamma \delta =1$, each of such sets must have at least two elements such that $\sum_{(i,j)\in E_k}{\gamma_i\delta_j}=0$ for all $k\in \{2,3,4\}$, but since $6-4=2$, $|E_k|\leq 1$ for some $k\in \{2,3,4\}$, a contradiction. Therefore, $n\not=2$.
\item
Let $n=3$. Then with the above discussion, $9\geq |supp(\gamma) supp(\delta)|\geq 5$ and so, there are at least $4$ sets different from $E_1$ in $\pi$, namely $E_2,E_3,E_4,E_5$. Since $\gamma \delta =1$, each of such sets must have at least two elements such that $\sum_{(i,j)\in E_k}{\gamma_i\delta_j}=0$ for all $k\in \{2,3,4,5\}$, but since $9-5=4$, $|E_1|=1$ and $|E_k|=2$ for all $k\in \{2,3,4,5\}$. Therefore $E_1=\{(1,1)\}$ and for all $k\in \{2,3,4,5\}$, $E_k=\{(i,j),(i',j')\}$ where $h_ig_j=h_{i'}g_{j'}$ for some $(i,j),(i',j')\in A$ such that $i\not=i'$ and $j\not=j'$. Let $\gamma'=\sum_{a\in supp(\gamma)}{a}$ and $\delta'=\sum_{b\in supp(\delta)}{b}$. So, $\gamma', \delta' \in \mathbb{F}_2[G]$, $|supp(\gamma')|=3$ and $|supp(\delta')|=3$ and with the above discussion we have $\gamma'\delta'=1$, that is a contradiction with Theorem \ref{dyk-rem}. Therefore, $n\not=3$.
\end{enumerate}

\subsection{The support of $\delta$ must be of size greater than or equal to $8$}\label{unitsub-2}
Without loss of generality we may assume that $G$ is generated by $supp(\gamma)\cup supp(\delta)$, since otherwise we replace $G$ by the subgroup generated by this set. Also, $1\in supp(\gamma)$ and $1\in supp(\delta)$ and with the same discussion such as in Lemma \ref{supp}, $\langle supp(\delta) \rangle=\langle  supp(\gamma) \rangle$. Therefore  by Theorem \ref{hamidoune}, if $n\geq 4$, then $|supp(\gamma) supp(\delta)|\geq |supp(\gamma)|+|supp(\delta)|+1$. Also, it is easy to see that $|supp(\gamma)| |supp(\delta)|\geq|supp(\gamma) supp(\delta)|$. Hence, $3n\geq |supp(\gamma) supp(\delta)|\geq 4+n$.
\begin{enumerate}
\item
Let $n=4$. Then with the above discussion, $12\geq |supp(\gamma) supp(\delta)|\geq 8$ and so, there are at least $7$ sets different from $E_1$ in $\pi$, namely $E_2,E_3,\ldots,E_8$. Since $\gamma \delta =1$, each of such sets must have at least two elements such that $\sum_{(i,j)\in E_k}{\gamma_i\delta_j}=0$ for all $k\in \{2,3,\ldots,8\}$, but since $12-8=4$, $|E_k|\leq 1$ for some $k\in \{2,3,\ldots,8\}$, a contradiction. Therefore, $n\not=4$.
\item
Let $n=5$. Then with the discussion above $15\geq |supp(\gamma) supp(\delta)|\geq 9$ and so, there are at least $8$ sets different from $E_1$ in $\pi$, namely $E_2,E_3,\ldots,E_9$. Since $\gamma \delta =1$, each of such sets must have at least two elements such that $\sum_{(i,j)\in E_k}{\gamma_i\delta_j}=0$ for all $k\in \{2,3,\ldots,9\}$, but since $15-9=6$, $|E_k|\leq 1$ for some $k\in \{2,3,\ldots,9\}$, a contradiction. Therefore, $n\not=5$.
\item
Let $n=6$. Then with the discussion above $18\geq |supp(\gamma) supp(\delta)|\geq 10$ and so, there are at least $9$ sets different from $E_1$ in $\pi$, namely $E_2,E_3,\ldots,E_{10}$. Since $\gamma \delta =1$, each of such sets must have at least two elements such that $\sum_{(i,j)\in E_k}{\gamma_i\delta_j}=0$ for all $k\in \{2,3,\ldots,10\}$, but since $18-10=8$, $|E_k|\leq 1$ for some $k\in \{2,3,\ldots,10\}$, a contradiction. Therefore, $n\not=6$.
\item
Let $n=7$. Then with the discussion above $21\geq |supp(\gamma) supp(\delta)|\geq 11$ and so, there are at least $10$ sets different from $E_1$ in $\pi$, namely $E_2,E_3,\ldots,E_{11}$. Since $\gamma \delta =1$, each of such sets must have at least two elements such that $\sum_{(i,j)\in E_k}{\gamma_i\delta_j}=0$ for all $k\in \{2,3,\ldots,11\}$, but since $21-11=10$, $|E_1|=1$ and $|E_k|=2$ for all $k\in \{2,3,\ldots,11\}$. Therefore $E_1=\{(1,1)\}$ and for all $k\in \{2,3,\ldots,11\}$, $E_k=\{(i,j),(i',j')\}$ where $h_ig_j=h_{i'}g_{j'}$ for some $(i,j),(i',j')\in A$ such that $i\not=i'$ and $j\not=j'$. Let $\gamma'=\sum_{a\in supp(\gamma)}{a}$ and $\delta'=\sum_{b\in supp(\delta)}{b}$. So, $\gamma', \delta' \in \mathbb{F}_2[G]$, $|supp(\gamma')|=3$ and $|supp(\delta')|=7$ and with the above discussion we have $\gamma'\delta'=1$, that is a contradiction with Theorem \ref{dyk-rem}. Therefore, $n\not=7$.
\item
Let $n=8$. Then with the discussion above $24\geq |supp(\gamma) supp(\delta)|\geq 12$. Let $|supp(\gamma) supp(\delta)|> 12$. Then $|supp(\gamma) supp(\delta)|\geq 13$ and so, there are at least $12$ sets different from $E_1$ in $\pi$, namely $E_2,E_3,\ldots,E_{13}$. Since $\gamma \delta =1$, each of such sets must have at least two elements such that $\sum_{(i,j)\in E_k}{\gamma_i\delta_j}=0$ for all $k\in \{2,3,\ldots,13\}$, but since $24-13=11$, $|E_k|\leq 1$ for some $k\in \{2,3,\ldots,13\}$, a contradiction. So, $|supp(\gamma) supp(\delta)|=12$. Therefore, there are $11$ sets different from $E_1$ in $\pi$, namely $E_2,E_3,\ldots,E_{12}$. Since $\gamma \delta =1$, there are two cases for the number of elements in such sets. 
\begin{enumerate}
\item
$|E_1|=2$ and for all $k\in \{2,3,\ldots,12\}$, $E_k=\{(i,j),(i',j')\}$ where $h_ig_j=h_{i'}g_{j'}$ for some $(i,j),(i',j')\in A$ such that $i\not=i'$ and $j\not=j'$. Let $\gamma'=\sum_{a\in supp(\gamma)}{a}$ and $\delta'=\sum_{b\in supp(\delta)}{b}$. So, $\gamma', \delta' \in \mathbb{F}_2[G]$, $|supp(\gamma')|=3$ and $|supp(\delta')|=8$ and with the above discussion we have $\gamma'\delta'=0$, that is a contradiction (see Corollary \ref{maincor}). 
\item
$E_1=\{(1,1)\}$ and $|E_l|=3$ for exactly one $l\in \{2,3,\ldots,12\}$ and for all $k\in \{2,3,\ldots,12\}\setminus \{l\}$, $E_k=\{(i,j),(i',j')\}$ where $h_ig_j=h_{i'}g_{j'}$ for some $(i,j),(i',j')\in A$ such that $i\not=i'$ and $j\not=j'$.
\end{enumerate}
\end{enumerate}

\begin{thm}\label{main-unit1}
Let $\gamma$ and $\delta$ be elements of the group algebra of any torsion-free group over an arbitrary field. If $|supp(\gamma)|=3$ and $\gamma \delta =1$ then $|supp(\delta)|\geq 8$. 
\end{thm}
\subsection{Kaplansky unit graphs over $\mathbb{F}$}
By the discussion from the beginning of the section, similar to the case of zero divisors, it can be associated a graph to $\gamma$ and $\delta$ with the vertex set $supp(\delta)$ such that two vertices $g_i$ and $g_j$ are adjacent whenever $h_{i'}g_i=h_{j'} g_j$ for some distinct $i',j' \in \{1,2,3\}$. We call the graph Kaplansky unit graph of $(\gamma,\delta)$ over $\mathbb{F}$ and it is denoted  by ${Ku}_{\mathbb{F}}(\gamma,\delta)$. The connectedness follows from the way we have chosen $\delta$ of minimum support size with respect to the property $\gamma\delta=1$.  Also by Lemma \ref{size of S}, $|S|=6$ where $S=\{h^{-1} h' \;|\; h,h'\in supp(\gamma), h\not=h' \}$. Therefore, ${Ku}_{\mathbb{F}}(\gamma,\delta)$ is a simple graph. Furthermore by Theorems \ref{K3-K3} and \ref{K2,3}, with a similar discussion as about  Kaplansky graphs, ${Ku}_{\mathbb{F}}(\gamma,\delta)$ contains no $K_3-K_3$ or $K_{2,3}$ as a subgraph, too.

By item $(5)$ of Subsection \ref{unitsub-2}, if $n=8$, $|supp(\gamma) supp(\delta)|=12$ and there are exactly $11$ sets different from $E_1$ in $\pi$, namely $E_2,E_3,\ldots,E_{12}$. Also, $E_1=\{(1,1)\}$ and $|E_l|=3$ for exactly one $l\in \{2,3,\ldots,12\}$ and for all $k\in \{2,3,\ldots,12\}\setminus \{l\}$, $E_k=\{(i,j),(i',j')\}$ such that $i\not=i'$, $j\not=j'$ and $h_ig_j=h_{i'}g_{j'}$.

\begin{figure}[ht]
\centering
\psscalebox{1.0 1.0} 
{
\begin{pspicture}(0,-2.305)(5.52,2.305)
\psdots[linecolor=black, dotsize=0.4](2.8,0.495)
\psdots[linecolor=black, dotsize=0.4](2.0,-0.705)
\psdots[linecolor=black, dotsize=0.4](3.6,-0.705)
\psdots[linecolor=black, dotsize=0.4](2.8,1.695)
\psdots[linecolor=black, dotsize=0.4](0.8,-0.705)
\psdots[linecolor=black, dotsize=0.4](1.6,-1.905)
\psdots[linecolor=black, dotsize=0.4](4.0,-1.905)
\psdots[linecolor=black, dotsize=0.4](4.8,-0.705)
\psline[linecolor=black, linewidth=0.04](2.8,0.495)(2.0,-0.705)(3.6,-0.705)(2.8,0.495)(2.8,1.695)
\psline[linecolor=black, linewidth=0.04](2.0,-0.705)(1.6,-1.905)
\psline[linecolor=black, linewidth=0.04](0.8,-0.705)(2.0,-0.705)
\psline[linecolor=black, linewidth=0.04](4.0,-1.905)(3.6,-0.705)(4.8,-0.705)
\rput[bl](3.04,0.495){$g_1$}
\rput[bl](2.8,1.92){$g_2$}
\rput[bl](3.6,-0.5){$g_3$}
\rput[bl](1.6,-0.5){$g_4$}
\rput[bl](0.2,-0.705){$g_5$}
\rput[bl](1.2,-2.305){$g_6$}
\rput[bl](4.1,-2.305){$g_7$}
\rput[bl](5.03,-0.705){$g_8$}
\end{pspicture}
}
\caption{The subgraph $H$ of ${Ku}_{\mathbb{F}}(\gamma,\delta)$ for the case that $n=8$}\label{f-H}
\end{figure}
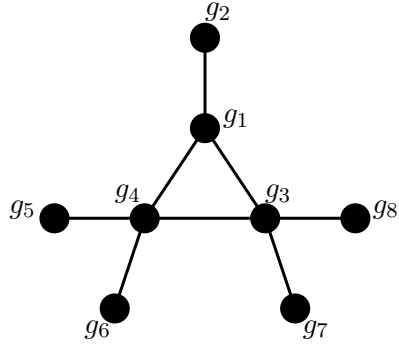

Let $E_l=\{(i_1,j_1),(i_2,j_2),(i_3,j_3)\}$. Therefore, $h_{i_1}g_{j_1}=h_{i_2}g_{j_2}=h_{i_3}g_{j_3}$ and so there is a triangle in ${Ku}_{\mathbb{F}}(\gamma,\delta)$ with the vertex set $\{g_{j_1},g_{j_2},g_{j_3}\}$ and there is no other triangle in the latter graph. Let $(2,1),(3,1)\notin E_l$. Then by the way we have chosen $E_k$ for $k\in \{1,2,3,\ldots,12\}$, the degree of $g_{j_1}$, $g_{j_2}$ and $g_{j_3}$ are equal to $4$, denoted by $deg(g_{j_1})=deg(g_{j_2})=deg(g_{j_3})=4$. So, there must be $6$ other vertices different from $g_{j_1}$, $g_{j_2}$ and $g_{j_3}$ in the vertex set of ${Ku}_{\mathbb{F}}(\gamma,\delta)$ because there is no other triangle in the latter graph. This leads us to a contradiction because the size of the vertex set of ${Ku}_{\mathbb{F}}(\gamma,\delta)$ is $n=8$. Hence, $E_l=\{(a,1),(i,j),(i',j')\}$ where $a\in\{2,3\}$ and $\{h_a,h_i,h_{i'}\}=supp(\gamma)$. Since $|E_1|=1$, $1=h_1g_1\not=h_mg_n$ for all $(m,n)\in A\setminus E_1$. So, $deg(g_1)=3$ and by renumbering, we may assume that ${Ku}_{\mathbb{F}}(\gamma,\delta)$ has the graph $H$ in Figure \ref{f-H} as a subgraph and there is no other vertices in ${Ku}_{\mathbb{F}}(\gamma,\delta)$. In $H$, $g_5\sim g_2$ or $g_5\not\sim g_2$.

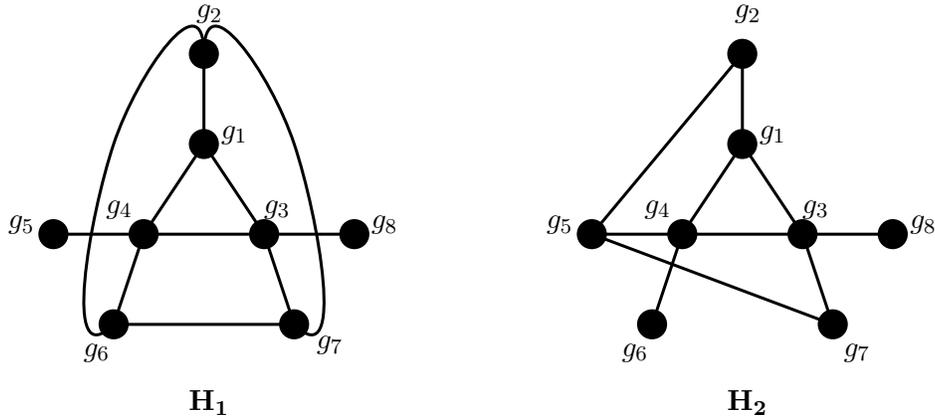
\begin{figure}[h]
\centering
\psscalebox{1.0 1.0} 
{
\begin{pspicture}(0,-2.705)(5.52,2.705)
\psdots[linecolor=black, dotsize=0.4](2.8,0.895)
\psdots[linecolor=black, dotsize=0.4](2.0,-0.305)
\psdots[linecolor=black, dotsize=0.4](3.6,-0.305)
\psdots[linecolor=black, dotsize=0.4](2.8,2.095)
\psdots[linecolor=black, dotsize=0.4](0.8,-0.305)
\psdots[linecolor=black, dotsize=0.4](1.6,-1.505)
\psdots[linecolor=black, dotsize=0.4](4.0,-1.505)
\psdots[linecolor=black, dotsize=0.4](4.8,-0.305)
\psline[linecolor=black, linewidth=0.04](2.8,0.895)(2.0,-0.305)(3.6,-0.305)(2.8,0.895)(2.8,2.095)
\psline[linecolor=black, linewidth=0.04](2.0,-0.305)(1.6,-1.505)
\psline[linecolor=black, linewidth=0.04](0.8,-0.305)(2.0,-0.305)
\psline[linecolor=black, linewidth=0.04](4.0,-1.505)(3.6,-0.305)(4.8,-0.305)
\rput[bl](3.03,0.895){$g_1$}
\rput[bl](2.7,2.5){$g_2$}
\rput[bl](3.6,-0.095){$g_3$}
\rput[bl](1.5,-0.095){$g_4$}
\rput[bl](0.2,-0.305){$g_5$}
\rput[bl](1.2,-2){$g_6$}
\rput[bl](4.3,-1.905){$g_7$}
\rput[bl](5.03,-0.305){$g_8$}
\psbezier[linecolor=black, linewidth=0.04](4.0,-1.505)(4.7071066,-2.2121067)(4.316228,-0.0536833)(4.0,0.895)(3.6837723,1.8436832)(2.8,3.095)(2.8,2.095)
\psbezier[linecolor=black, linewidth=0.04](1.6,-1.505)(0.8928932,-2.2121067)(1.2837722,-0.0536833)(1.6,0.895)(1.9162278,1.8436832)(2.8,3.095)(2.8,2.095)
\psline[linecolor=black, linewidth=0.04](1.6,-1.505)(4.0,-1.505)
\rput[bl](2.6,-2.705){$\mathbf{H_1}$}
\end{pspicture}
} 
\hspace{2.5cc}
\psscalebox{1.0 1.0} 
{
\begin{pspicture}(0,-2.705)(5.52,2.705)
\psdots[linecolor=black, dotsize=0.4](2.8,0.895)
\psdots[linecolor=black, dotsize=0.4](2.0,-0.305)
\psdots[linecolor=black, dotsize=0.4](3.6,-0.305)
\psdots[linecolor=black, dotsize=0.4](2.8,2.095)
\psdots[linecolor=black, dotsize=0.4](0.8,-0.305)
\psdots[linecolor=black, dotsize=0.4](1.6,-1.505)
\psdots[linecolor=black, dotsize=0.4](4.0,-1.505)
\psdots[linecolor=black, dotsize=0.4](4.8,-0.305)
\psline[linecolor=black, linewidth=0.04](2.8,0.895)(2.0,-0.305)(3.6,-0.305)(2.8,0.895)(2.8,2.095)
\psline[linecolor=black, linewidth=0.04](2.0,-0.305)(1.6,-1.505)
\psline[linecolor=black, linewidth=0.04](0.8,-0.305)(2.0,-0.305)
\psline[linecolor=black, linewidth=0.04](4.0,-1.505)(3.6,-0.305)(4.8,-0.305)
\rput[bl](3.03,0.895){$g_1$}
\rput[bl](2.7,2.5){$g_2$}
\rput[bl](3.6,-0.095){$g_3$}
\rput[bl](1.5,-0.095){$g_4$}
\rput[bl](0.2,-0.305){$g_5$}
\rput[bl](1.2,-2){$g_6$}
\rput[bl](4.15,-2){$g_7$}
\rput[bl](5.03,-0.305){$g_8$}
\psline[linecolor=black, linewidth=0.04](2.8,2.095)(0.8,-0.305)(4.0,-1.505)
\rput[bl](2.6,-2.705){$\mathbf{H_2}$}
\end{pspicture}
}
\caption{Two possible subgraphs of ${Ku}_{\mathbb{F}}(\gamma,\delta)$ for the case that $n=8$}\label{f-H1-H2}
\end{figure}

Let $g_5\not\sim g_2$. Since $(a,2),(a,5),(a,6),(a,7),(a,8) \notin E_l,E_1$ for all $a\in \{1,2,3\}$ and $|E_k|=2$ for all  $k\in \{2,3,\ldots,12\}\setminus \{l\}$, $deg(g_2)=deg(g_5)=deg(g_6)=deg(g_7)=deg(g_8)=3$. Since ${Ku}_{\mathbb{F}}(\gamma,\delta)$ is a simple graph which contains exactly one triangle, $g_6\not\sim g_i$ for all $i\in \{1,3,4,5,6\}$. So, without loss of generality we may assume that $g_6\sim g_7$ since $deg(g_6)=3$.  Suppose, for a contradiction, that $g_6\sim g_8$. Then there is a subgraph isomorphic to $K_{2,3}$ in ${Ku}_{\mathbb{F}}(\gamma,\delta)$ with the vertex set $\{g_3,g_4,g_6,g_7,g_8\}$, a contradiction. So, $g_6\sim g_2$. With a same discussion we have $g_7\sim g_2$ and ${Ku}_{\mathbb{F}}(\gamma,\delta)$ has the graph $H_1$ in Figure \ref{f-H1-H2} as a subgraph. Since $deg(g_2)=deg(g_5)=deg(g_6)=deg(g_7)=deg(g_8)=3$ and $deg(g_1)=deg(g_3)=deg(g_4)=4$, we must have $g_5\sim g_8$ with a double edge, a contradiction because ${Ku}_{\mathbb{F}}(\gamma,\delta)$ is simple. Therefore, $g_5\sim g_2$.

Since ${Ku}_{\mathbb{F}}(\gamma,\delta)$ is a simple graph which contains exactly one triangle, $g_5\not\sim g_i$ for all $i\in \{1,2,3,4,5,6\}$. So, $g_5\sim g_7$ or $g_5\sim g_8$ because $deg(g_5)=3$. As we can see in Figure \ref{f-H}, without loss of generality we may assume that $g_5\sim g_7$ and ${Ku}_{\mathbb{F}}(\gamma,\delta)$ has the graph $H_2$ in Figure \ref{f-H1-H2} as a subgraph. Since ${Ku}_{\mathbb{F}}(\gamma,\delta)$ is a simple graph which contains exactly one triangle, $g_6\not\sim g_i$ for all $i\in \{1,3,4,5,6\}$. So, $g_6\sim g_2$, $g_6\sim g_7$ or $g_6\sim g_8$ because $deg(g_6)=3$. If $g_6\sim g_2$ or $g_6\sim g_7$ then there is a subgraph isomorphic to $K_{2,3}$ in ${Ku}_{\mathbb{F}}(\gamma,\delta)$ as we can see in the graph $H_2$ of Figure \ref{f-H1-H2}, a contradiction. Therefore we must have $g_6\sim g_8$ with a double edge, a contradiction because ${Ku}_{\mathbb{F}}(\gamma,\delta)$ is simple. 

Hence with the above discussion, $n\not=8$ and by Theorem \ref{main-unit1}, we have the following result.
\begin{thm}\label{main-unit2}
Let $\gamma$ and $\delta$ be elements of the group algebra of any torsion-free group over an arbitrary field. If $|supp(\gamma)|=3$ and $\gamma \delta =1$ then $|supp(\delta)|\geq 9$. 
\end{thm}

\newpage

\section{\bf Appendix}\label{S2-app}
In Section \ref{S2}, we found $44$ forbidden subgraphs for the Kaplansky graphs over $\mathbb{F}_2$ without details. In this section, we give more details about finding such subgraph of $K(\alpha,\beta)$. Firstly with the same discussion such as about $C_3$ and $C_4$ cycles, we study cycles of lengths $5$ and $6$ with their relations. Then we use such relations and the relations of $C_4$ and $C_7$ cycles to show that Kaplansky graphs do not contain the latter $44$ graphs. Such forbidden subgraphs are numbered from $1$ to $44$ such that the forbidden subgraph $K_{2,3}$ is numbered by $1$. 

$\mathbf{C_5}$ \textbf{cycles:} With the same discussion such as about $C_4$ cycles, there are $105$ non-equivalent cases for the relations of a $C_5$ cycle in the graph $K(\alpha,\beta)$. Such relations are listed in table \ref{tab-C5}. In the following, we show that some of these relations lead to a contradiction. Each of such relations is marked by a $*$ in the Table \ref{tab-C5}.

\begin{enumerate}
\item[(1)]$h_2^5=1$:\\
$h_2^5=1 \text{ and } G \text{ is torsion-free}\Rightarrow h_2=1\text{, a contradiction}$.
\item[(2)]$h_2^4h_3=1$:\\
$h_3=h_2^{-4}\Rightarrow \langle h_2,h_3\rangle=\langle h_2\rangle \text{ is abelian, a contradiction}$.
\item[(3)]$h_2^4h_3^{-1}h_2=1$:\\
$h_3=h_2^5\Rightarrow \langle h_2,h_3\rangle=\langle h_2\rangle \text{ is abelian, a contradiction}$.
\item[(6)]$h_2^3h_3^{-1}h_2^{-1}h_3=1$:\\
$h_3^{-1}h_2h_3=h_2^3\Rightarrow \langle h_2,h_3\rangle\cong BS(1,3) \text{ is solvable, a contradiction}$.
\item[(8)]$h_2^3h_3^{-1}h_2h_3=1$:\\
$h_3^{-1}h_2h_3=h_2^{-3}\Rightarrow \langle h_2,h_3\rangle\cong BS(1,-3) \text{ is solvable, a contradiction}$.
\item[(10)]$h_2(h_2h_3)^2=1$:\\
$h_2=(h_2h_3)^{-2} \text{ and } \langle h_2,h_3\rangle=\langle h_2,h_2h_3\rangle\Rightarrow \langle h_2,h_3\rangle=\langle h_2h_3\rangle \text{ is abelian, a contradiction}$.
\item[(11)]$h_2^2h_3h_2h_3^{-1}h_2=1$:\\
$h_3^{-1}h_2^{-3}h_3=h_2\Rightarrow \langle h_2,h_3\rangle\cong BS(-3,1) \text{ is solvable, a contradiction}$.
\item[(15)]$h_2^2h_3h_2^{-1}h_3^{-1}h_2=1$:\\
$h_3^{-1}h_2^{3}h_3=h_2\Rightarrow \langle h_2,h_3\rangle\cong BS(3,1) \text{ is solvable, a contradiction}$.
\item[(27)]$(h_2^2h_3^{-1})^2h_2=1$:\\
$h_2=(h_3h_2^{-2})^{2} \text{ and } \langle h_2,h_3\rangle=\langle h_2,h_3h_2^{-2}\rangle\Rightarrow \langle h_2,h_3\rangle=\langle h_3h_2^{-2}\rangle \text{ is abelian, a contradiction}$.
\item[(34)]$(h_2h_3)^2h_3=1$:\\
By interchanging $h_2$ and $h_3$ in (10) and with the same discussion, there is a contradiction.
\item[(36)]$h_2h_3h_2h_3^{-1}h_2^{-1}h_3=1$:\\
Using Tietze transformation where $h_2\mapsto h_2h_3^{-1}$ and $h_3\mapsto h_3$, we have:\\ $h_2^2h_3^{-1}h_2^{-1}h_3=1\Rightarrow h_3^{-1}h_2h_3=h_2^2\Rightarrow \langle h_2,h_3\rangle\cong BS(1,2) \text{ is solvable, a contradiction}$.
\item[(41)]$h_2h_3^4=1$:\\
By interchanging $h_2$ and $h_3$ in (2) and with the same discussion, there is a contradiction.

\begin{center}
\begin{longtable}{|c|l||c|l||c|l|}
\caption{{\small All non-equivalent cases for the relations of a $C_5$ cycle}}\label{tab-C5}\\
\hline \multicolumn{1}{|c|}{$n$} & \multicolumn{1}{l||}{$R$} & \multicolumn{1}{c|}{$n$} & \multicolumn{1}{l||}{$R$} & \multicolumn{1}{c|}{$n$} & \multicolumn{1}{l|}{$R$}\\\hline
\endfirsthead
\multicolumn{6}{c}%
{{\tablename\ \thetable{} -- continued from previous page}} \\
\hline \multicolumn{1}{|c|}{$n$} & \multicolumn{1}{l||}{$R$} & \multicolumn{1}{c|}{$n$} & \multicolumn{1}{l||}{$R$} & \multicolumn{1}{c|}{$n$} & \multicolumn{1}{l|}{$R$}\\\hline
\endhead
\hline \multicolumn{6}{|r|}{{Continued on next page}} \\\hline
\endfoot
\hline
\endlastfoot
$1$&$h_2^5=1 \ *$&$36$&$h_2h_3h_2h_3^{-1}h_2^{-1}h_3=1 \ *$&$71$&$h_2h_3^{-3}h_2h_3=1$\\
$2$&$h_2^4h_3=1 \ *$&$37$&$h_2h_3h_2h_3^{-2}h_2=1$&$72$&$h_2h_3^{-2}(h_3^{-1}h_2)^2=1$\\
$3$&$h_2^4h_3^{-1}h_2=1 \ *$&$38$&$h_2h_3h_2h_3^{-1}h_2h_3=1$&$73$&$h_2h_3^{-2}h_2^2h_3^{-1}h_2=1 \ *$\\
$4$&$h_2^3h_3^2=1$&$39$&$h_2h_3(h_2h_3^{-1})^2h_2=1$&$74$&$h_2h_3^{-2}h_2h_3^2=1$\\
$5$&$h_2^3h_3h_2^{-1}h_3=1$&$40$&$h_2h_3^2h_2h_3^{-1}h_2=1$&$75$&$h_2h_3^{-2}h_2h_3h_2^{-1}h_3=1$\\
$6$&$h_2^3h_3^{-1}h_2^{-1}h_3=1 \ *$&$41$&$h_2h_3^4=1 \ *$&$76$&$h_2h_3^{-2}h_2h_3^{-1}h_2^{-1}h_3=1 \ *$\\
$7$&$h_2^3h_3^{-2}h_2=1$&$42$&$h_2h_3^3h_2^{-1}h_3=1 \ *$&$77$&$(h_2h_3^{-2})^2h_2=1 \ *$\\
$8$&$h_2^3h_3^{-1}h_2h_3=1 \ *$&$43$&$h_2h_3^2h_2^{-2}h_3=1$&$78$&$h_2h_3^{-1}(h_3^{-1}h_2)^2h_3=1$\\
$9$&$h_2^2(h_2h_3^{-1})^2h_2=1$&$44$&$h_2h_3^2h_2^{-1}h_3^{-1}h_2=1$&$79$&$h_2h_3^{-1}(h_3^{-1}h_2)^3=1 \ *$\\
$10$&$h_2(h_2h_3)^2=1 \ *$&$45$&$h_2h_3^2h_2^{-1}h_3^2=1$&$80$&$(h_2h_3^{-1}h_2)^2h_3=1$\\
$11$&$h_2^2h_3h_2h_3^{-1}h_2=1 \ *$&$46$&$h_2h_3(h_3h_2^{-1})^2h_3=1$&$81$&$(h_2h_3^{-1}h_2)^2h_3^{-1}h_2=1 \ *$\\
$12$&$h_2^2h_3^3=1$&$47$&$h_2h_3h_2^{-2}h_3^2=1$&$82$&$h_2h_3^{-1}h_2h_3^3=1$\\
$13$&$h_2^2h_3^2h_2^{-1}h_3=1$&$48$&$h_2h_3h_2^{-1}(h_2^{-1}h_3)^2=1$&$83$&$h_2h_3^{-1}h_2h_3^2h_2^{-1}h_3=1$\\
$14$&$h_2^2h_3h_2^{-2}h_3=1$&$49$&$h_2h_3h_2^{-1}h_3^{-1}h_2^{-1}h_3=1$&$84$&$h_2h_3^{-1}h_2h_3h_2^{-1}h_3^2=1$\\
$15$&$h_2^2h_3h_2^{-1}h_3^{-1}h_2=1 \ *$&$50$&$h_2h_3h_2^{-1}h_3^{-2}h_2=1$&$85$&$h_2h_3^{-1}h_2(h_3h_2^{-1})^2h_3=1$\\
$16$&$h_2^2h_3h_2^{-1}h_3^2=1$&$51$&$h_2h_3h_2^{-1}h_3^{-1}h_2h_3=1 \ *$&$86$&$(h_2h_3^{-1})^2h_2^{-1}h_3^2=1$\\
$17$&$h_2^2(h_3h_2^{-1})^2h_3=1$&$52$&$h_2h_3h_2^{-1}(h_3^{-1}h_2)^2=1$&$87$&$(h_2h_3^{-1})^2(h_2^{-1}h_3)^2=1$\\
$18$&$h_2^2h_3^{-1}h_2^{-2}h_3=1$&$53$&$h_2h_3h_2^{-1}h_3h_2h_3^{-1}h_2=1$&$88$&$(h_2h_3^{-1})^2h_3^{-1}h_2^{-1}h_3=1$\\
$19$&$h_2^2h_3^{-1}h_2^{-1}h_3^{-1}h_2=1$&$54$&$h_2h_3h_2^{-1}h_3^3=1 \ *$&$89$&$(h_2h_3^{-1})^2h_3^{-2}h_2=1$\\
$20$&$h_2^2h_3^{-1}h_2^{-1}h_3^2=1$&$55$&$h_2(h_3h_2^{-1}h_3)^2=1$&$90$&$(h_2h_3^{-1})^2h_3^{-1}h_2h_3=1$\\
$21$&$h_2^2h_3^{-1}(h_2^{-1}h_3)^2=1$&$56$&$h_2(h_3h_2^{-1})^2h_2^{-1}h_3=1$&$91$&$h_2(h_3^{-1}h_2h_3^{-1})^2h_2=1$\\
$22$&$h_2^2h_3^{-2}h_2^{-1}h_3=1$&$57$&$h_2(h_3h_2^{-1})^2h_3^{-1}h_2=1$&$92$&$(h_2h_3^{-1})^2h_2h_3^2=1$\\
$23$&$h_2^2h_3^{-3}h_2=1$&$58$&$h_2(h_3h_2^{-1})^2h_3^2=1$&$93$&$(h_2h_3^{-1})^2h_2h_3h_2^{-1}h_3=1$\\
$24$&$h_2^2h_3^{-2}h_2h_3=1$&$59$&$h_2(h_3h_2^{-1})^3h_3=1$&$94$&$(h_2h_3^{-1})^3h_2^{-1}h_3=1 \ *$\\
$25$&$h_2^2h_3^{-1}(h_3^{-1}h_2)^2=1$&$60$&$h_2h_3^{-1}h_2^{-1}h_3h_2h_3^{-1}h_2=1 \ *$&$95$&$(h_2h_3^{-1})^3h_3^{-1}h_2=1 \ *$\\
$26$&$h_2^2h_3^{-1}h_2^2h_3=1$&$61$&$h_2h_3^{-1}h_2^{-1}h_3^3=1 \ *$&$96$&$(h_2h_3^{-1})^3h_2h_3=1$\\
$27$&$(h_2^2h_3^{-1})^2h_2=1 \ *$&$62$&$h_2h_3^{-1}h_2^{-1}h_3^2h_2^{-1}h_3=1 \ *$&$97$&$(h_2h_3^{-1})^4h_2=1 \ *$\\
$28$&$h_2^2h_3^{-1}h_2h_3^2=1$&$63$&$h_2h_3^{-1}h_2^{-1}h_3h_2^{-2}h_3=1 \ *$&$98$&$h_3^5=1 \ *$\\
$29$&$h_2^2h_3^{-1}h_2h_3h_2^{-1}h_3=1$&$64$&$h_2h_3^{-1}h_2^{-1}h_3h_2^{-1}h_3^{-1}h_2=1$&$99$&$h_3^4h_2^{-1}h_3=1 \ *$\\
$30$&$h_2(h_2h_3^{-1})^2h_2^{-1}h_3=1$&$65$&$h_2h_3^{-1}(h_2^{-1}h_3)^2h_3=1$&$100$&$h_3^2(h_3h_2^{-1})^2h_3=1$\\
$31$&$h_2(h_2h_3^{-1})^2h_3^{-1}h_2=1$&$66$&$h_2h_3^{-1}(h_2^{-1}h_3)^3=1 \ *$&$101$&$(h_3^2h_2^{-1})^2h_3=1 \ *$\\
$32$&$h_2(h_2h_3^{-1})^2h_2h_3=1$&$67$&$h_2h_3^{-2}h_2^{-1}h_3^2=1$&$102$&$h_3(h_3h_2^{-1})^3h_3=1$\\
$33$&$h_2(h_2h_3^{-1})^3h_2=1$&$68$&$h_2h_3^{-2}(h_2^{-1}h_3)^2=1$&$103$&$(h_3h_2^{-1}h_3)^2h_2^{-1}h_3=1 \ *$\\
$34$&$(h_2h_3)^2h_3=1 \ *$&$69$&$h_2h_3^{-3}h_2^{-1}h_3=1 \ *$&$104$&$(h_3h_2^{-1})^4h_3=1 \ *$\\
$35$&$(h_2h_3)^2h_2^{-1}h_3=1$&$70$&$h_2h_3^{-4}h_2=1$&$105$&$(h_2^{-1}h_3)^5=1 \ *$\\
\end{longtable}
\end{center}
\item[(42)]$h_2h_3^3h_2^{-1}h_3=1$:\\
By interchanging $h_2$ and $h_3$ in (8) and with the same discussion, there is a contradiction.
\item[(51)]$h_2h_3h_2^{-1}h_3^{-1}h_2h_3=1$:\\
By interchanging $h_2$ and $h_3$ in (36) and with the same discussion, there is a contradiction.
\item[(54)]$h_2h_3h_2^{-1}h_3^3=1$:\\
By interchanging $h_2$ and $h_3$ in (11) and with the same discussion, there is a contradiction.
\item[(60)]$h_2h_3^{-1}h_2^{-1}h_3h_2h_3^{-1}h_2=1$:\\
Using Tietze transformation where $h_3\mapsto h_3h_2$ and $h_2\mapsto h_2$, we have $h_3^{-1}h_2^{-1}h_3h_2h_3^{-1}h_2=1$. Using Tietze transformation again where $h_3\mapsto h_3h_2$ and $h_2\mapsto h_2$, we have:\\ $h_3^{-2}h_2^{-1}h_3h_2=1\Rightarrow h_2^{-1}h_3h_2=h_3^2\Rightarrow \langle h_2,h_3\rangle\cong BS(1,2) \text{ is solvable, a contradiction}$.
\item[(61)]$h_2h_3^{-1}h_2^{-1}h_3^3=1$:\\
By interchanging $h_2$ and $h_3$ in (15) and with the same discussion, there is a contradiction
\item[(62)]$h_2h_3^{-1}h_2^{-1}h_3^2h_2^{-1}h_3=1$:\\
Using Tietze transformation where $h_2\mapsto h_3h_2$ and $h_3\mapsto h_3$, we have $h_3h_2h_3^{-1}h_2^{-1}h_3h_2^{-1}=1$. Using Tietze transformation again where $h_2\mapsto h_2h_3$ and $h_3\mapsto h_3$, we have:\\ $h_3h_2h_3^{-1}h_2^{-2}=1\Rightarrow h_3^{-1}h_2^2h_3=h_2\Rightarrow \langle h_2,h_3\rangle\cong BS(2,1) \text{ is solvable, a contradiction}$.
\item[(63)]$h_2h_3^{-1}h_2^{-1}h_3h_2^{-2}h_3=1$:\\
By interchanging $h_2$ and $h_3$ in (62) and with the same discussion, there is a contradiction
\item[(66)]$h_2h_3^{-1}(h_2^{-1}h_3)^3=1$:\\
Using Tietze transformation where $h_2\mapsto h_3h_2$ and $h_3\mapsto h_3$, we have:\\ $h_3h_2h_3^{-1}h_2^{-3}=1\Rightarrow h_3^{-1}h_2^{3}h_3=h_2\Rightarrow \langle h_2,h_3\rangle\cong BS(3,1) \text{ is solvable, a contradiction}$.
\item[(69)]$h_2h_3^{-3}h_2^{-1}h_3=1$:\\
By interchanging $h_2$ and $h_3$ in (6) and with the same discussion, there is a contradiction
\item[(73)]$h_2h_3^{-2}h_2^2h_3^{-1}h_2=1$:\\
Using Tietze transformation where $h_3\mapsto h_2h_3$ and $h_2\mapsto h_2$, we have:\\ $h_2h_3^{-1}h_2^{-1}h_3^{-2}=1\Rightarrow h_2^{-1}h_3^{-2}h_2=h_3\Rightarrow \langle h_2,h_3\rangle\cong BS(-2,1) \text{ is solvable, a contradiction}$.
\item[(76)]$h_2h_3^{-2}h_2h_3^{-1}h_2^{-1}h_3=1$:\\
By interchanging $h_2$ and $h_3$ in (60) and with the same discussion, there is a contradiction
\item[(77)]$(h_2h_3^{-2})^2h_2=1$:\\
By interchanging $h_2$ and $h_3$ in (73) and with the same discussion, there is a contradiction
\item[(79)]$h_2h_3^{-1}(h_3^{-1}h_2)^3=1$:\\
Using Tietze transformation where $h_2\mapsto h_3h_2$ and $h_3\mapsto h_3$, we have:\\ $h_3h_2h_3^{-1}h_2^{3}=1\Rightarrow h_3^{-1}h_2^{-3}h_3=h_2\Rightarrow \langle h_2,h_3\rangle\cong BS(-3,1) \text{ is solvable, a contradiction}$.
\item[(81)]$(h_2h_3^{-1}h_2)^2h_3^{-1}h_2=1$:\\
$(h_2h_3^{-1}h_2)^2=h_2^{-1}h_3 \text{ and } \langle h_2,h_3\rangle=\langle h_2h_3^{-1}h_2,h_2^{-1}h_3\rangle\Rightarrow \langle h_2,h_3\rangle=\langle h_2h_3^{-1}h_2\rangle$  is abelian, \\ a contradiction.
\item[(94)]$(h_2h_3^{-1})^3h_2^{-1}h_3=1$:\\
Using Tietze transformation where $h_3\mapsto h_3h_2$ and $h_2\mapsto h_2$, we have:\\ $h_3^{-3}h_2^{-1}h_3h_2=1\Rightarrow h_2^{-1}h_3h_2=h_3^3\Rightarrow \langle h_2,h_3\rangle\cong BS(1,3) \text{ is solvable, a contradiction}$.
\item[(95)]$(h_2h_3^{-1})^3h_3^{-1}h_2=1$:\\
Using Tietze transformation where $h_2\mapsto h_2h_3$ and $h_3\mapsto h_3$, we have:\\ $h_2^{3}h_3^{-1}h_2h_3=1\Rightarrow h_3^{-1}h_2h_3=h_2^{-3}\Rightarrow \langle h_2,h_3\rangle\cong BS(1,-3) \text{ is solvable, a contradiction}$.
\item[(97)]$(h_2h_3^{-1})^4h_2=1$:\\
$h_2=(h_3h_2^{-1})^4\text{ and } \langle h_2,h_3\rangle=\langle h_2,h_3h_2^{-1}\rangle\Rightarrow \langle h_2,h_3\rangle=\langle h_3h_2^{-1}\rangle\text{ is abelian, a contradiction}$.
\item[(98)]$h_3^5=1$:\\
By interchanging $h_2$ and $h_3$ in (1) and with the same discussion, there is a contradiction.
\item[(99)]$h_3^4h_2^{-1}h_3=1$:\\
By interchanging $h_2$ and $h_3$ in (3) and with the same discussion, there is a contradiction.
\item[(101)]$(h_3^2h_2^{-1})^2h_3=1$:\\
By interchanging $h_2$ and $h_3$ in (27) and with the same discussion, there is a contradiction.
\item[(103)]$(h_3h_2^{-1}h_3)^2h_2^{-1}h_3=1$:\\
By interchanging $h_2$ and $h_3$ in (81) and with the same discussion, there is a contradiction.
\item[(104)]$(h_3h_2^{-1})^4h_3=1$:\\
By interchanging $h_2$ and $h_3$ in (97) and with the same discussion, there is a contradiction.
\item[(105)]$(h_2^{-1}h_3)^5=1$:\\
$(h_2^{-1}h_3)^5=1 \text{ and } G \text{ is torsion-free}\Rightarrow h_2=h_3\Rightarrow \langle h_2,h_3\rangle=\langle h_2\rangle\text{ is abelian, a contradiction}$.
\end{enumerate}

$\mathbf{C_6}$ \textbf{cycles:} With the same discussion such as about $C_4$ cycles and by considering the relations corresponding to some $C_6$ cycles equivalent, there are $351$ non-equivalent cases for the relations of a  $C_6$ cycle in the graph $K(\alpha,\beta)$. These relations are listed in table \ref{tab-C6}. In the following, we show that some of these relations lead to a contradiction. Each of such relations is marked by a $*$ in the Table \ref{tab-C6}.

\begin{enumerate}
\item[(1)]$h_2^6=1$:\\
$h_2^6=1 \text{ and } G \text{ is torsion-free}\Rightarrow h_2=1\text{, a contradiction}$.
\item[(2)]$h_2^5h_3=1$:\\
$h_3=h_2^{-5}\Rightarrow \langle h_2,h_3\rangle=\langle h_2\rangle \text{ is abelian, a contradiction}$.
\item[(3)]$h_2^5h_3^{-1}h_2=1$:\\
$h_3=h_2^6\Rightarrow \langle h_2,h_3\rangle=\langle h_2\rangle \text{ is abelian, a contradiction}$.
\item[(6)]$h_2^4h_3^{-1}h_2^{-1}h_3=1$:\\
$h_3^{-1}h_2h_3=h_2^4\Rightarrow \langle h_2,h_3\rangle\cong BS(1,4) \text{ is solvable, a contradiction}$.
\item[(8)]$h_2^4h_3^{-1}h_2h_3=1$:\\
$h_3^{-1}h_2h_3=h_2^{-4}\Rightarrow \langle h_2,h_3\rangle\cong BS(1,-4) \text{ is solvable, a contradiction}$.
\item[(10)]$h_2^2(h_2h_3)^2=1$:\\
$h_2^2=(h_3^{-1}h_2^{-1})^{2} \Rightarrow \langle h_2,h_3\rangle=\langle h_2,h_3^{-1}h_2^{-1}\rangle\cong BS(1,-1) \text{ is solvable, a contradiction}$.
\item[(11)]$h_2^3h_3h_2h_3^{-1}h_2=1$:\\
$h_3^{-1}h_2^{-4}h_3=h_2\Rightarrow \langle h_2,h_3\rangle\cong BS(-4,1) \text{ is solvable, a contradiction}$.
\item[(15)]$h_2^3h_3h_2^{-1}h_3^{-1}h_2=1$:\\
$h_3^{-1}h_2^{4}h_3=h_2\Rightarrow \langle h_2,h_3\rangle\cong BS(4,1) \text{ is solvable, a contradiction}$.
\item[(27)]$h_2(h_2^2h_3^{-1})^2h_2=1$:\\
$h_2^2=(h_3h_2^{-2})^{2} \Rightarrow \langle h_2,h_3\rangle=\langle h_2,h_3h_2^{-2}\rangle\cong BS(1,-1) \text{ is solvable, a contradiction}$.

\begin{center}
\small
\begin{longtable}{|c|l||c|l||c|l|}
\caption{All non-equivalent cases for the relations of a $C_6$ cycle}\label{tab-C6}\\
\hline \multicolumn{1}{|c|}{$n$} & \multicolumn{1}{l||}{$R$} & \multicolumn{1}{c|}{$n$} & \multicolumn{1}{l||}{$R$} & \multicolumn{1}{c|}{$n$} & \multicolumn{1}{l|}{$R$}\\\hline
\endfirsthead
\multicolumn{6}{c}%
{{\tablename\ \thetable{} -- continued from previous page}} \\
\hline \multicolumn{1}{|c|}{$n$} & \multicolumn{1}{l||}{$R$} & \multicolumn{1}{c|}{$n$} & \multicolumn{1}{l||}{$R$} & \multicolumn{1}{c|}{$n$} & \multicolumn{1}{l|}{$R$}\\\hline
\endhead
\hline \multicolumn{6}{|r|}{{Continued on next page}} \\\hline
\endfoot
\hline
\endlastfoot
$1$&$h_2^6=1 \ *$&$40$&$h_2^2h_3h_2h_3^{-1}h_2h_3=1$&$79$&$h_2^2h_3^{-1}(h_2^{-1}h_3)^3=1$\\
$2$&$h_2^5h_3=1 \ *$&$41$&$h_2^2h_3(h_2h_3^{-1})^2h_2=1$&$80$&$h_2^2h_3^{-2}h_2^{-2}h_3=1$\\
$3$&$h_2^5h_3^{-1}h_2=1 \ *$&$42$&$h_2^2h_3^2h_2h_3=1$&$81$&$h_2^2h_3^{-2}h_2^{-1}h_3^{-1}h_2=1$\\
$4$&$h_2^4h_3^2=1$&$43$&$h_2^2h_3^2h_2h_3^{-1}h_2=1$&$82$&$h_2^2h_3^{-2}h_2^{-1}h_3^2=1$\\
$5$&$h_2^4h_3h_2^{-1}h_3=1$&$44$&$h_2^2h_3^4=1$&$83$&$h_2^2h_3^{-2}(h_2^{-1}h_3)^2=1$\\
$6$&$h_2^4h_3^{-1}h_2^{-1}h_3=1 \ *$&$45$&$h_2^2h_3^3h_2^{-1}h_3=1$&$84$&$h_2^2h_3^{-3}h_2^{-1}h_3=1$\\
$7$&$h_2^4h_3^{-2}h_2=1$&$46$&$h_2^2h_3^2h_2^{-2}h_3=1$&$85$&$h_2^2h_3^{-4}h_2=1$\\
$8$&$h_2^4h_3^{-1}h_2h_3=1 \ *$&$47$&$h_2^2h_3^2h_2^{-1}h_3^{-1}h_2=1$&$86$&$h_2^2h_3^{-3}h_2h_3=1$\\
$9$&$h_2^3(h_2h_3^{-1})^2h_2=1$&$48$&$h_2^2h_3^2h_2^{-1}h_3^2=1$&$87$&$h_2^2h_3^{-2}(h_3^{-1}h_2)^2=1$\\
$10$&$h_2^2(h_2h_3)^2=1 \ *$&$49$&$h_2^2h_3(h_3h_2^{-1})^2h_3=1$&$88$&$h_2^2h_3^{-2}h_2^2h_3=1$\\
$11$&$h_2^3h_3h_2h_3^{-1}h_2=1 \ *$&$50$&$h_2^2h_3h_2^{-3}h_3=1$&$89$&$h_2^2h_3^{-2}h_2^2h_3^{-1}h_2=1$\\
$12$&$h_2^3h_3^3=1$&$51$&$h_2^2h_3h_2^{-2}h_3^{-1}h_2=1$&$90$&$h_2^2h_3^{-2}h_2h_3^2=1$\\
$13$&$h_2^3h_3^2h_2^{-1}h_3=1$&$52$&$h_2^2h_3h_2^{-2}h_3^2=1$&$91$&$h_2^2h_3^{-2}h_2h_3h_2^{-1}h_3=1$\\
$14$&$h_2^3h_3h_2^{-2}h_3=1$&$53$&$h_2^2h_3h_2^{-1}(h_2^{-1}h_3)^2=1$&$92$&$h_2^2h_3^{-2}h_2h_3^{-1}h_2^{-1}h_3=1$\\
$15$&$h_2^3h_3h_2^{-1}h_3^{-1}h_2=1 \ *$&$54$&$h_2^2h_3h_2^{-1}h_3^{-1}h_2^{-1}h_3=1$&$93$&$h_2(h_2h_3^{-2})^2h_2=1$\\
$16$&$h_2^3h_3h_2^{-1}h_3^2=1$&$55$&$h_2^2h_3h_2^{-1}h_3^{-2}h_2=1$&$94$&$h_2^2h_3^{-1}(h_3^{-1}h_2)^2h_3=1$\\
$17$&$h_2^3(h_3h_2^{-1})^2h_3=1$&$56$&$h_2^2h_3h_2^{-1}h_3^{-1}h_2h_3=1 \ *$&$95$&$h_2^2h_3^{-1}(h_3^{-1}h_2)^3=1$\\
$18$&$h_2^3h_3^{-1}h_2^{-2}h_3=1$&$57$&$h_2^2h_3h_2^{-1}(h_3^{-1}h_2)^2=1$&$96$&$(h_2^2h_3^{-1}h_2)^2=1 \ *$\\
$19$&$h_2^3h_3^{-1}h_2^{-1}h_3^{-1}h_2=1$&$58$&$h_2^2h_3h_2^{-1}h_3h_2h_3=1$&$97$&$h_2^2h_3^{-1}h_2^2h_3^2=1$\\
$20$&$h_2^3h_3^{-1}h_2^{-1}h_3^2=1$&$59$&$h_2^2h_3h_2^{-1}h_3h_2h_3^{-1}h_2=1$&$98$&$h_2^2h_3^{-1}h_2^2h_3h_2^{-1}h_3=1$\\
$21$&$h_2^3h_3^{-1}(h_2^{-1}h_3)^2=1$&$60$&$h_2^2h_3h_2^{-1}h_3^3=1$&$99$&$(h_2^2h_3^{-1})^2h_2^{-1}h_3=1$\\
$22$&$h_2^3h_3^{-2}h_2^{-1}h_3=1$&$61$&$h_2^2(h_3h_2^{-1}h_3)^2=1$&$100$&$(h_2^2h_3^{-1})^2h_3^{-1}h_2=1$\\
$23$&$h_2^3h_3^{-3}h_2=1$&$62$&$h_2^2(h_3h_2^{-1})^2h_2^{-1}h_3=1$&$101$&$(h_2^2h_3^{-1})^2h_2h_3=1$\\
$24$&$h_2^3h_3^{-2}h_2h_3=1$&$63$&$h_2^2(h_3h_2^{-1})^2h_3^{-1}h_2=1$&$102$&$(h_2^2h_3^{-1})^2h_2h_3^{-1}h_2=1 \ *$\\
$25$&$h_2^3h_3^{-1}(h_3^{-1}h_2)^2=1$&$64$&$h_2^2(h_3h_2^{-1})^2h_3^2=1$&$103$&$h_2^2h_3^{-1}(h_2h_3)^2=1$\\
$26$&$h_2^3h_3^{-1}h_2^2h_3=1$&$65$&$h_2^2(h_3h_2^{-1})^3h_3=1$&$104$&$h_2^2h_3^{-1}h_2h_3h_2h_3^{-1}h_2=1$\\
$27$&$h_2(h_2^2h_3^{-1})^2h_2=1 \ *$&$66$&$h_2^2h_3^{-1}h_2^{-2}h_3^2=1$&$105$&$h_2^2h_3^{-1}h_2h_3^3=1$\\
$28$&$h_2^3h_3^{-1}h_2h_3^2=1$&$67$&$h_2^2h_3^{-1}h_2^{-1}(h_2^{-1}h_3)^2=1$&$106$&$h_2^2h_3^{-1}h_2h_3^2h_2^{-1}h_3=1$\\
$29$&$h_2^3h_3^{-1}h_2h_3h_2^{-1}h_3=1$&$68$&$h_2^2(h_3^{-1}h_2^{-1})^2h_3=1$&$107$&$h_2^2h_3^{-1}h_2h_3h_2^{-2}h_3=1$\\
$30$&$h_2^2(h_2h_3^{-1})^2h_2^{-1}h_3=1$&$69$&$h_2^2h_3^{-1}h_2^{-1}h_3^{-2}h_2=1$&$108$&$h_2^2h_3^{-1}h_2h_3h_2^{-1}h_3^{-1}h_2=1 \ *$\\
$31$&$h_2^2(h_2h_3^{-1})^2h_3^{-1}h_2=1$&$70$&$h_2^2h_3^{-1}h_2^{-1}h_3^{-1}h_2h_3=1$&$109$&$h_2^2h_3^{-1}h_2h_3h_2^{-1}h_3^2=1$\\
$32$&$h_2^2(h_2h_3^{-1})^2h_2h_3=1$&$71$&$h_2^2h_3^{-1}h_2^{-1}(h_3^{-1}h_2)^2=1$&$110$&$h_2^2h_3^{-1}h_2(h_3h_2^{-1})^2h_3=1$\\
$33$&$h_2^2(h_2h_3^{-1})^3h_2=1$&$72$&$h_2^2h_3^{-1}h_2^{-1}h_3h_2h_3=1$&$111$&$h_2(h_2h_3^{-1})^2h_2^{-2}h_3=1$\\
$34$&$(h_2^2h_3)^2=1 \ *$&$73$&$h_2^2h_3^{-1}h_2^{-1}h_3h_2h_3^{-1}h_2=1 \ *$&$112$&$h_2(h_2h_3^{-1})^2h_2^{-1}h_3^{-1}h_2=1$\\
$35$&$h_2^2h_3h_2^2h_3^{-1}h_2=1$&$74$&$h_2^2h_3^{-1}h_2^{-1}h_3^3=1$&$113$&$h_2(h_2h_3^{-1})^2h_2^{-1}h_3^2=1$\\
$36$&$h_2(h_2h_3)^2h_3=1$&$75$&$h_2^2h_3^{-1}h_2^{-1}h_3^2h_2^{-1}h_3=1$&$114$&$h_2(h_2h_3^{-1})^2(h_2^{-1}h_3)^2=1$\\
$37$&$h_2(h_2h_3)^2h_2^{-1}h_3=1$&$76$&$h_2^2h_3^{-1}h_2^{-1}h_3h_2^{-2}h_3=1$&$115$&$h_2(h_2h_3^{-1})^2h_3^{-1}h_2^{-1}h_3=1$\\
$38$&$h_2^2h_3h_2h_3^{-1}h_2^{-1}h_3=1 \ *$&$77$&$h_2^2h_3^{-1}h_2^{-1}h_3h_2^{-1}h_3^{-1}h_2=1$&$116$&$h_2(h_2h_3^{-1})^2h_3^{-2}h_2=1$\\
$39$&$h_2^2h_3h_2h_3^{-2}h_2=1$&$78$&$h_2^2h_3^{-1}(h_2^{-1}h_3)^2h_3=1$&$117$&$h_2(h_2h_3^{-1})^2h_3^{-1}h_2h_3=1$\\
$118$&$h_2^2(h_3^{-1}h_2h_3^{-1})^2h_2=1$&$157$&$h_2h_3^5=1 \ *$&$196$&$(h_2h_3h_2^{-1}h_3)^2=1 \ *$\\
$119$&$h_2(h_2h_3^{-1})^2h_2^2h_3=1$&$158$&$h_2h_3^4h_2^{-1}h_3=1 \ *$&$197$&$h_2h_3h_2^{-1}h_3h_2h_3^{-1}h_2^{-1}h_3=1$\\
$120$&$h_2^2h_3^{-1}(h_2h_3^{-1}h_2)^2=1$&$159$&$h_2h_3^3h_2^{-2}h_3=1$&$198$&$h_2h_3h_2^{-1}h_3h_2h_3^{-2}h_2=1$\\
$121$&$h_2(h_2h_3^{-1})^2h_2h_3^2=1$&$160$&$h_2h_3^3h_2^{-1}h_3^{-1}h_2=1$&$199$&$h_2h_3h_2^{-1}h_3h_2h_3^{-1}h_2h_3=1$\\
$122$&$h_2(h_2h_3^{-1})^2h_2h_3h_2^{-1}h_3=1$&$161$&$h_2h_3^3h_2^{-1}h_3^2=1$&$200$&$h_2h_3h_2^{-1}h_3(h_2h_3^{-1})^2h_2=1$\\
$123$&$h_2(h_2h_3^{-1})^3h_2^{-1}h_3=1$&$162$&$h_2h_3^2(h_3h_2^{-1})^2h_3=1$&$201$&$h_2h_3h_2^{-1}h_3^2h_2h_3^{-1}h_2=1$\\
$124$&$h_2(h_2h_3^{-1})^3h_3^{-1}h_2=1$&$163$&$h_2h_3^2h_2^{-2}h_3^2=1$&$202$&$h_2h_3h_2^{-1}h_3^4=1 \ *$\\
$125$&$h_2(h_2h_3^{-1})^3h_2h_3=1$&$164$&$h_2h_3^2h_2^{-1}(h_2^{-1}h_3)^2=1$&$203$&$h_2h_3h_2^{-1}h_3^3h_2^{-1}h_3=1$\\
$126$&$h_2(h_2h_3^{-1})^4h_2=1$&$165$&$h_2h_3^2h_2^{-1}h_3^{-1}h_2^{-1}h_3=1$&$204$&$h_2h_3h_2^{-1}h_3^2h_2^{-2}h_3=1$\\
$127$&$(h_2h_3)^3=1 \ *$&$166$&$h_2h_3^2h_2^{-1}h_3^{-2}h_2=1$&$205$&$h_2h_3h_2^{-1}h_3^2h_2^{-1}h_3^{-1}h_2=1$\\
$128$&$(h_2h_3)^2h_2h_3^{-1}h_2=1$&$167$&$h_2h_3^2h_2^{-1}h_3^{-1}h_2h_3=1$&$206$&$h_2(h_3h_2^{-1}h_3)^2h_3=1$\\
$129$&$(h_2h_3)^2h_3^2=1 \ *$&$168$&$h_2h_3^2h_2^{-1}(h_3^{-1}h_2)^2=1$&$207$&$h_2(h_3h_2^{-1}h_3)^2h_2^{-1}h_3=1$\\
$130$&$(h_2h_3)^2h_3h_2^{-1}h_3=1$&$169$&$h_2h_3^2h_2^{-1}h_3h_2h_3^{-1}h_2=1$&$208$&$h_2(h_3h_2^{-1})^2h_2^{-1}h_3^2=1$\\
$131$&$(h_2h_3)^2h_2^{-2}h_3=1$&$170$&$h_2h_3^2h_2^{-1}h_3^3=1$&$209$&$h_2h_3(h_2^{-1}h_3h_2^{-1})^2h_3=1$\\
$132$&$(h_2h_3)^2h_2^{-1}h_3^{-1}h_2=1$&$171$&$h_2(h_3^2h_2^{-1})^2h_3=1$&$210$&$h_2(h_3h_2^{-1})^2h_3^{-1}h_2^{-1}h_3=1$\\
$133$&$(h_2h_3)^2h_2^{-1}h_3^2=1$&$172$&$h_2h_3(h_3h_2^{-1})^2h_2^{-1}h_3=1$&$211$&$h_2(h_3h_2^{-1})^2h_3^{-2}h_2=1$\\
$134$&$(h_2h_3)^2(h_2^{-1}h_3)^2=1$&$173$&$h_2h_3(h_3h_2^{-1})^2h_3^{-1}h_2=1$&$212$&$h_2(h_3h_2^{-1})^2h_3^{-1}h_2h_3=1$\\
$135$&$h_2h_3h_2h_3^{-1}h_2^{-2}h_3=1$&$174$&$h_2h_3(h_3h_2^{-1})^2h_3^2=1$&$213$&$h_2(h_3h_2^{-1})^2(h_3^{-1}h_2)^2=1$\\
$136$&$h_2h_3h_2h_3^{-1}h_2^{-1}h_3^{-1}h_2=1$&$175$&$h_2h_3(h_3h_2^{-1})^3h_3=1$&$214$&$h_2(h_3h_2^{-1})^2h_3h_2h_3^{-1}h_2=1$\\
$137$&$h_2h_3h_2h_3^{-1}h_2^{-1}h_3^2=1$&$176$&$h_2h_3h_2^{-2}h_3h_2h_3^{-1}h_2=1$&$215$&$h_2(h_3h_2^{-1})^2h_3^3=1$\\
$138$&$h_2h_3h_2h_3^{-1}(h_2^{-1}h_3)^2=1$&$177$&$h_2h_3h_2^{-2}h_3^3=1$&$216$&$h_2h_3h_2^{-1}(h_3h_2^{-1}h_3)^2=1$\\
$139$&$h_2h_3h_2h_3^{-2}h_2^{-1}h_3=1$&$178$&$h_2h_3h_2^{-2}h_3^2h_2^{-1}h_3=1$&$217$&$h_2(h_3h_2^{-1})^3h_2^{-1}h_3=1$\\
$140$&$h_2h_3h_2h_3^{-3}h_2=1$&$179$&$h_2(h_3h_2^{-2})^2h_3=1$&$218$&$h_2(h_3h_2^{-1})^3h_3^{-1}h_2=1$\\
$141$&$h_2h_3h_2h_3^{-2}h_2h_3=1$&$180$&$h_2h_3h_2^{-2}h_3h_2^{-1}h_3^{-1}h_2=1$&$219$&$h_2(h_3h_2^{-1})^3h_3^2=1$\\
$142$&$h_2h_3h_2h_3^{-1}(h_3^{-1}h_2)^2=1$&$181$&$h_2h_3h_2^{-1}(h_2^{-1}h_3)^2h_3=1$&$220$&$h_2(h_3h_2^{-1})^4h_3=1$\\
$143$&$h_2h_3(h_2h_3^{-1}h_2)^2=1$&$182$&$h_2h_3h_2^{-1}(h_2^{-1}h_3)^3=1$&$221$&$(h_2h_3^{-1}h_2^{-1}h_3)^2=1 \ *$\\
$144$&$h_2h_3h_2h_3^{-1}h_2h_3^2=1$&$183$&$h_2h_3h_2^{-1}h_3^{-1}h_2^{-1}h_3^2=1$&$222$&$h_2h_3^{-1}h_2^{-1}h_3h_2h_3^{-2}h_2=1$\\
$145$&$h_2h_3h_2h_3^{-1}h_2h_3h_2^{-1}h_3=1$&$184$&$h_2h_3h_2^{-1}h_3^{-1}(h_2^{-1}h_3)^2=1$&$223$&$h_2h_3^{-1}h_2^{-1}h_3h_2h_3^{-1}h_2h_3=1$\\
$146$&$h_2h_3(h_2h_3^{-1})^2h_2^{-1}h_3=1$&$185$&$h_2h_3h_2^{-1}h_3^{-2}h_2^{-1}h_3=1$&$224$&$h_2h_3^{-1}h_2^{-1}h_3(h_2h_3^{-1})^2h_2=1 \ *$\\
$147$&$h_2h_3(h_2h_3^{-1})^2h_3^{-1}h_2=1$&$186$&$h_2h_3h_2^{-1}h_3^{-3}h_2=1$&$225$&$h_2h_3^{-1}h_2^{-1}h_3^2h_2h_3^{-1}h_2=1$\\
$148$&$h_2h_3(h_2h_3^{-1})^2h_2h_3=1$&$187$&$h_2h_3h_2^{-1}h_3^{-2}h_2h_3=1$&$226$&$h_2h_3^{-1}h_2^{-1}h_3^4=1 \ *$\\
$149$&$h_2h_3(h_2h_3^{-1})^3h_2=1$&$188$&$h_2h_3h_2^{-1}h_3^{-1}(h_3^{-1}h_2)^2=1$&$227$&$h_2h_3^{-1}h_2^{-1}h_3^3h_2^{-1}h_3=1 \ *$\\
$150$&$(h_2h_3^2)^2=1 \ *$&$189$&$h_2h_3h_2^{-1}h_3^{-1}h_2^2h_3^{-1}h_2=1$&$228$&$h_2h_3^{-1}h_2^{-1}h_3^2h_2^{-2}h_3=1$\\
$151$&$h_2h_3^2h_2h_3h_2^{-1}h_3=1$&$190$&$h_2h_3h_2^{-1}h_3^{-1}h_2h_3^2=1 \ *$&$229$&$h_2h_3^{-1}h_2^{-1}h_3^2h_2^{-1}h_3^{-1}h_2=1$\\
$152$&$h_2h_3^2h_2h_3^{-1}h_2^{-1}h_3=1 \ *$&$191$&$h_2h_3h_2^{-1}h_3^{-1}h_2h_3h_2^{-1}h_3=1$&$230$&$h_2h_3^{-1}(h_2^{-1}h_3^2)^2=1$\\
$153$&$h_2h_3^2h_2h_3^{-2}h_2=1$&$192$&$h_2h_3h_2^{-1}h_3^{-1}h_2h_3^{-1}h_2^{-1}h_3=1$&$231$&$h_2h_3^{-1}h_2^{-1}h_3(h_3h_2^{-1})^2h_3=1 \ *$\\
$154$&$h_2h_3^2h_2h_3^{-1}h_2h_3=1$&$193$&$h_2h_3h_2^{-1}h_3^{-1}h_2h_3^{-2}h_2=1$&$232$&$h_2h_3^{-1}h_2^{-1}h_3h_2^{-2}h_3^2=1$\\
$155$&$h_2h_3^2(h_2h_3^{-1})^2h_2=1$&$194$&$h_2h_3h_2^{-1}(h_3^{-1}h_2)^2h_3=1$&$233$&$h_2h_3^{-1}(h_2^{-1}h_3h_2^{-1})^2h_3=1$\\
$156$&$h_2h_3^3h_2h_3^{-1}h_2=1$&$195$&$h_2h_3h_2^{-1}(h_3^{-1}h_2)^3=1$&$234$&$h_2h_3^{-1}h_2^{-1}h_3h_2^{-1}h_3^{-1}h_2^{-1}h_3=1$\\
$$&$$&$ $&$$&$$&$ $\\
$235$&$h_2h_3^{-1}h_2^{-1}h_3h_2^{-1}h_3^{-2}h_2=1$&$274$&$h_2h_3^{-2}h_2h_3^{-1}(h_2^{-1}h_3)^2=1$&$313$&$(h_2h_3^{-1})^2h_3^{-2}h_2h_3=1$\\
$236$&$h_2h_3^{-1}h_2^{-1}h_3h_2^{-1}h_3^{-1}h_2h_3=1$&$275$&$(h_2h_3^{-2})^2h_2^{-1}h_3=1$&$314$&$(h_2h_3^{-1})^2h_3^{-1}(h_3^{-1}h_2)^2=1$\\
$237$&$h_2h_3^{-1}h_2^{-1}h_3h_2^{-1}(h_3^{-1}h_2)^2=1$&$276$&$(h_2h_3^{-2})^2h_3^{-1}h_2=1$&$315$&$(h_2h_3^{-1})^2h_3^{-1}h_2h_3^2=1$\\
$238$&$h_2h_3^{-1}(h_2^{-1}h_3)^2h_2h_3^{-1}h_2=1$&$277$&$(h_2h_3^{-2})^2h_2h_3=1$&$316$&$(h_2h_3^{-1})^2h_3^{-1}h_2h_3h_2^{-1}h_3=1$\\
$239$&$h_2h_3^{-1}(h_2^{-1}h_3)^2h_3^2=1$&$278$&$(h_2h_3^{-2})^2h_2h_3^{-1}h_2=1$&$317$&$h_2(h_3^{-1}h_2h_3^{-1})^2h_2^{-1}h_3=1$\\
$240$&$h_2h_3^{-1}h_2^{-1}(h_3h_2^{-1}h_3)^2=1$&$279$&$h_2h_3^{-2}(h_2h_3^{-1}h_2)^2=1$&$318$&$h_2(h_3^{-1}h_2h_3^{-1})^2h_3^{-1}h_2=1$\\
$241$&$h_2h_3^{-1}(h_2^{-1}h_3)^2h_2^{-2}h_3=1 \ *$&$280$&$h_2h_3^{-1}(h_3^{-1}h_2)^2h_3^2=1$&$319$&$h_2(h_3^{-1}h_2h_3^{-1})^2h_2h_3=1$\\
$242$&$h_2h_3^{-1}(h_2^{-1}h_3)^2h_2^{-1}h_3^{-1}h_2=1$&$281$&$h_2h_3^{-1}(h_3^{-1}h_2)^2h_3h_2^{-1}h_3=1$&$320$&$h_2(h_3^{-1}h_2h_3^{-1})^2h_2h_3^{-1}h_2=1$\\
$243$&$h_2h_3^{-1}(h_2^{-1}h_3)^3h_3=1$&$282$&$h_2h_3^{-1}(h_3^{-1}h_2)^2h_3^{-1}h_2^{-1}h_3=1 \ *$&$321$&$((h_2h_3^{-1})^2h_2)^2=1 \ *$\\
$244$&$h_2h_3^{-1}(h_2^{-1}h_3)^4=1 \ *$&$283$&$h_2h_3^{-1}(h_3^{-1}h_2)^2h_3^{-2}h_2=1$&$322$&$(h_2h_3^{-1})^2h_2h_3^3=1$\\
$245$&$h_2h_3^{-2}h_2^{-1}h_3h_2h_3^{-1}h_2=1$&$284$&$h_2h_3^{-1}(h_3^{-1}h_2)^3h_3=1$&$323$&$(h_2h_3^{-1})^2h_2h_3^2h_2^{-1}h_3=1$\\
$246$&$h_2h_3^{-2}h_2^{-1}h_3^3=1$&$285$&$h_2h_3^{-1}(h_3^{-1}h_2)^4=1 \ *$&$324$&$(h_2h_3^{-1})^2h_2h_3h_2^{-1}h_3^2=1$\\
$247$&$h_2h_3^{-2}h_2^{-1}h_3^2h_2^{-1}h_3=1$&$286$&$(h_2h_3^{-1}h_2)^3=1 \ *$&$325$&$(h_2h_3^{-1})^2h_2(h_3h_2^{-1})^2h_3=1$\\
$248$&$h_2h_3^{-2}h_2^{-1}h_3h_2^{-2}h_3=1$&$287$&$(h_2h_3^{-1}h_2)^2h_3^2=1$&$326$&$(h_2h_3^{-1})^3h_2^{-1}h_3^2=1$\\
$249$&$h_2h_3^{-2}h_2^{-1}h_3h_2^{-1}h_3^{-1}h_2=1$&$288$&$(h_2h_3^{-1}h_2)^2h_3h_2^{-1}h_3=1$&$327$&$(h_2h_3^{-1})^3(h_2^{-1}h_3)^2=1$\\
$250$&$h_2h_3^{-2}(h_2^{-1}h_3)^2h_3=1$&$289$&$(h_2h_3^{-1}h_2)^2h_3^{-1}h_2^{-1}h_3=1$&$328$&$(h_2h_3^{-1})^3h_3^{-1}h_2^{-1}h_3=1$\\
$251$&$h_2h_3^{-2}(h_2^{-1}h_3)^3=1$&$290$&$(h_2h_3^{-1}h_2)^2h_3^{-2}h_2=1$&$329$&$(h_2h_3^{-1})^3h_3^{-2}h_2=1$\\
$252$&$h_2h_3^{-3}h_2^{-1}h_3^2=1$&$291$&$(h_2h_3^{-1}h_2)^2h_3^{-1}h_2h_3=1$&$330$&$(h_2h_3^{-1})^3h_3^{-1}h_2h_3=1$\\
$253$&$h_2h_3^{-3}(h_2^{-1}h_3)^2=1$&$292$&$(h_2h_3^{-1}h_2)^2h_3^{-1}h_2h_3^{-1}h_2=1 \ *$&$331$&$h_2h_3^{-1}h_2(h_3^{-1}h_2h_3^{-1})^2h_2=1$\\
$254$&$h_2h_3^{-4}h_2^{-1}h_3=1 \ *$&$293$&$(h_2h_3^{-1}h_2h_3)^2=1 \ *$&$332$&$(h_2h_3^{-1})^3h_2h_3^2=1$\\
$255$&$h_2h_3^{-5}h_2=1$&$294$&$h_2h_3^{-1}h_2h_3(h_2h_3^{-1})^2h_2=1$&$333$&$(h_2h_3^{-1})^3h_2h_3h_2^{-1}h_3=1$\\
$256$&$h_2h_3^{-4}h_2h_3=1$&$295$&$h_2h_3^{-1}h_2h_3^4=1$&$334$&$(h_2h_3^{-1})^4h_2^{-1}h_3=1 \ *$\\
$257$&$h_2h_3^{-3}(h_3^{-1}h_2)^2=1$&$296$&$h_2h_3^{-1}h_2h_3^3h_2^{-1}h_3=1$&$335$&$(h_2h_3^{-1})^4h_3^{-1}h_2=1 \ *$\\
$258$&$h_2h_3^{-3}h_2^2h_3^{-1}h_2=1$&$297$&$h_2h_3^{-1}h_2h_3^2h_2^{-1}h_3^2=1$&$336$&$(h_2h_3^{-1})^4h_2h_3=1$\\
$259$&$h_2h_3^{-3}h_2h_3^2=1$&$298$&$h_2h_3^{-1}h_2h_3(h_3h_2^{-1})^2h_3=1$&$337$&$(h_2h_3^{-1})^5h_2=1 \ *$\\
$260$&$h_2h_3^{-3}h_2h_3h_2^{-1}h_3=1$&$299$&$h_2h_3^{-1}h_2h_3h_2^{-1}h_3^3=1$&$338$&$h_3^6=1 \ *$\\
$261$&$h_2h_3^{-3}h_2h_3^{-1}h_2^{-1}h_3=1 \ *$&$300$&$h_2h_3^{-1}h_2(h_3h_2^{-1}h_3)^2=1$&$339$&$h_3^5h_2^{-1}h_3=1 \ *$\\
$262$&$h_2h_3^{-1}(h_3^{-2}h_2)^2=1$&$301$&$h_2h_3^{-1}h_2(h_3h_2^{-1})^2h_2^{-1}h_3=1$&$340$&$h_3^3(h_3h_2^{-1})^2h_3=1$\\
$263$&$h_2h_3^{-2}(h_3^{-1}h_2)^2h_3=1$&$302$&$h_2h_3^{-1}h_2(h_3h_2^{-1})^2h_3^{-1}h_2=1$&$341$&$h_3(h_3^2h_2^{-1})^2h_3=1 \ *$\\
$264$&$h_2h_3^{-2}(h_3^{-1}h_2)^3=1$&$303$&$h_2h_3^{-1}h_2(h_3h_2^{-1})^2h_3^2=1$&$342$&$h_3^2(h_3h_2^{-1})^3h_3=1$\\
$265$&$(h_2h_3^{-2}h_2)^2=1 \ *$&$304$&$h_2h_3^{-1}h_2(h_3h_2^{-1})^3h_3=1$&$343$&$(h_3^2h_2^{-1}h_3)^2=1 \ *$\\
$266$&$h_2h_3^{-2}h_2^2h_3^{-1}h_2h_3=1$&$305$&$(h_2h_3^{-1})^2h_2^{-1}h_3^3=1$&$344$&$(h_3^2h_2^{-1})^2h_3h_2^{-1}h_3=1 \ *$\\
$267$&$h_2h_3^{-2}h_2(h_2h_3^{-1})^2h_2=1$&$306$&$(h_2h_3^{-1})^2h_2^{-1}h_3^2h_2^{-1}h_3=1$&$345$&$h_3^2h_2^{-1}(h_3h_2^{-1}h_3)^2=1$\\
$268$&$h_2h_3^{-2}h_2h_3h_2h_3^{-1}h_2=1$&$307$&$(h_2h_3^{-1})^2(h_2^{-1}h_3)^2h_3=1$&$346$&$h_3(h_3h_2^{-1})^4h_3=1$\\
$269$&$h_2h_3^{-2}h_2h_3^3=1$&$308$&$(h_2h_3^{-1})^2(h_2^{-1}h_3)^3=1$&$347$&$(h_3h_2^{-1}h_3)^3=1 \ *$\\
$270$&$h_2h_3^{-2}h_2h_3^2h_2^{-1}h_3=1$&$309$&$(h_2h_3^{-1})^2h_3^{-1}h_2^{-1}h_3^2=1$&$348$&$(h_3h_2^{-1}h_3)^2h_2^{-1}h_3h_2^{-1}h_3=1 \ *$\\
$271$&$h_2h_3^{-2}h_2h_3h_2^{-1}h_3^2=1$&$310$&$(h_2h_3^{-1})^2h_3^{-1}(h_2^{-1}h_3)^2=1$&$349$&$((h_3h_2^{-1})^2h_3)^2=1 \ *$\\
$272$&$h_2h_3^{-2}h_2(h_3h_2^{-1})^2h_3=1$&$311$&$(h_2h_3^{-1})^2h_3^{-2}h_2^{-1}h_3=1$&$350$&$(h_3h_2^{-1})^5h_3=1 \ *$\\
$273$&$h_2h_3^{-2}h_2h_3^{-1}h_2^{-1}h_3^2=1$&$312$&$(h_2h_3^{-1})^2h_3^{-3}h_2=1$&$351$&$(h_2^{-1}h_3)^6=1 \ *$\\
\end{longtable}
\end{center}

\item[(34)]$(h_2^2h_3)^2=1$:\\
$(h_2^2h_3)^2=1 \text{ and } G \text{ is torsion-free}\Rightarrow h_2^2h_3=1\Rightarrow h_3=h_2^{-2}\Rightarrow \langle h_2,h_3\rangle=\langle h_2\rangle \text{ is abelian,}$ a contradiction.
\item[(38)]$h_2^2h_3h_2h_3^{-1}h_2^{-1}h_3=1$:\\
Using Tietze transformation where $h_3\mapsto h_2^{-1}h_3$ and $h_2\mapsto h_2$, we have $h_2h_3h_2h_3^{-1}h_2^{-1}h_3=1$. Using Tietze transformation again where $h_3\mapsto h_2^{-1}h_3$ and $h_2\mapsto h_2$, we have:\\ $h_3h_2h_3^{-1}h_2^{-1}h_3=1\Rightarrow h_2^{-1}h_3^2h_2=h_3\Rightarrow \langle h_2,h_3\rangle\cong BS(2,1) \text{ is solvable, a contradiction}$.
\item[(56)]$h_2^2h_3h_2^{-1}h_3^{-1}h_2h_3=1$:\\
Using Tietze transformation where $h_3\mapsto h_2^{-1}h_3$ and $h_2\mapsto h_2$, we have $h_2h_3h_2^{-1}h_3^{-1}h_2h_3=1$. Using Tietze transformation again where $h_3\mapsto h_2^{-1}h_3$ and $h_2\mapsto h_2$, we have:\\ $h_3h_2^{-1}h_3^{-1}h_2h_3=1\Rightarrow h_2^{-1}h_3h_2=h_3^2\Rightarrow \langle h_2,h_3\rangle\cong BS(1,2) \text{ is solvable, a contradiction}$.
\item[(73)]$h_2^2h_3^{-1}h_2^{-1}h_3h_2h_3^{-1}h_2=1$:\\
Using Tietze transformation where $h_3\mapsto h_3h_2$ and $h_2\mapsto h_2$, we have $h_2^2h_3^{-1}h_2^{-1}h_3h_2h_3^{-1}=1$. Using Tietze transformation again where $h_3\mapsto h_3h_2$ and $h_2\mapsto h_2$, we have $h_2h_3^{-1}h_2^{-1}h_3h_2h_3^{-1}=1$.  Now by using Tietze transformation again where $h_3\mapsto h_3h_2$ and $h_2\mapsto h_2$, we have:\\ $h_3^{-1}h_2^{-1}h_3h_2h_3^{-1}=1\Rightarrow h_2^{-1}h_3h_2=h_3^2\Rightarrow \langle h_2,h_3\rangle\cong BS(1,2) \text{ is solvable, a contradiction}$.
\item[(96)]$(h_2^2h_3^{-1}h_2)^2=1$:\\
$(h_2^2h_3^{-1}h_2)^2=1 \text{ and } G \text{ is torsion-free}\Rightarrow h_2^2h_3^{-1}h_2=1\Rightarrow h_3=h_2^3\Rightarrow \langle h_2,h_3\rangle=\langle h_2\rangle \text{ is abelian,}$ a contradiction.
\item[(102)]$(h_2^2h_3^{-1})^2h_2h_3^{-1}h_2=1$:\\
Using Tietze transformation where $h_3\mapsto h_3h_2$ and $h_2\mapsto h_2$, we have $(h_2h_3^{-1})^2h_3^{-1}h_2=1$. Using Tietze transformation again where $h_3\mapsto h_3h_2$ and $h_2\mapsto h_2$, we have:\\ $(h_3^{-1})^2h_2^{-1}h_3^{-1}h_2=1\Rightarrow h_2^{-1}h_3h_2=h_3^{-2}\Rightarrow \langle h_2,h_3\rangle\cong BS(1,-2) \text{ is solvable, a contradiction}$.
\item[(108)]$h_2^2h_3^{-1}h_2h_3h_2^{-1}h_3^{-1}h_2=1$:\\
Using Tietze transformation where $h_3\mapsto h_3h_2$ and $h_2\mapsto h_2$, we have $h_2^2h_3^{-1}h_2h_3h_2^{-1}h_3^{-1}=1$. Using Tietze transformation again where $h_3\mapsto h_3h_2$ and $h_2\mapsto h_2$, we have $h_2h_3^{-1}h_2h_3h_2^{-1}h_3^{-1}=1$.  Now by using Tietze transformation again where $h_3\mapsto h_3h_2$ and $h_2\mapsto h_2$, we have:\\ $h_3^{-1}h_2h_3h_2^{-1}h_3^{-1}=1\Rightarrow h_2^{-1}h_3^2h_2=h_3\Rightarrow \langle h_2,h_3\rangle\cong BS(2,1) \text{ is solvable, a contradiction}$.
\item[(127)]$(h_2h_3)^3=1$:\\
$(h_2h_3)^3=1 \text{ and } G \text{ is torsion-free}\Rightarrow h_2h_3=1\Rightarrow h_3=h_2^{-1}\Rightarrow \langle h_2,h_3\rangle=\langle h_2\rangle \text{ is abelian,}$ a contradiction.
\item[(129)]$(h_2h_3)^2h_3^2=1$:\\
By interchanging $h_2$ and $h_3$ in (10) and with the same discussion, there is a contradiction.
\item[(150)]$(h_2h_3^2)^2=1$:\\
By interchanging $h_2$ and $h_3$ in (34) and with the same discussion, there is a contradiction.
\item[(152)]$h_2h_3^2h_2h_3^{-1}h_2^{-1}h_3=1$:\\
By interchanging $h_2$ and $h_3$ in (56) and with the same discussion, there is a contradiction.
\item[(157)]$h_2h_3^5=1$:\\
By interchanging $h_2$ and $h_3$ in (2) and with the same discussion, there is a contradiction.
\item[(158)]$h_2h_3^4h_2^{-1}h_3=1$:\\
By interchanging $h_2$ and $h_3$ in (8) and with the same discussion, there is a contradiction.
\item[(190)]$h_2h_3h_2^{-1}h_3^{-1}h_2h_3^2=1$:\\
By interchanging $h_2$ and $h_3$ in (38) and with the same discussion, there is a contradiction.
\item[(196)]$(h_2h_3h_2^{-1}h_3)^2=1$:\\
$(h_2h_3h_2^{-1}h_3)^2=1 \text{ and } G \text{ is torsion-free}\Rightarrow h_2h_3h_2^{-1}h_3=1\Rightarrow h_2^{-1}h_3h_2=h_3^{-1} \Rightarrow \langle h_2,h_3\rangle\cong BS(1,-1) \text{ is solvable, a contradiction}$.
\item[(202)]$h_2h_3h_2^{-1}h_3^4=1$:\\
By interchanging $h_2$ and $h_3$ in (11) and with the same discussion, there is a contradiction.
\item[(221)]$(h_2h_3^{-1}h_2^{-1}h_3)^2=1$:\\
$(h_2h_3^{-1}h_2^{-1}h_3)^2=1 \text{ and } G \text{ is torsion-free}\Rightarrow h_2h_3^{-1}h_2^{-1}h_3=1\Rightarrow h_3^{-1}h_2h_3=h_2 \Rightarrow \langle h_2,h_3\rangle\cong BS(1,1) \text{ is solvable, a contradiction}$.
\item[(224)]$h_2h_3^{-1}h_2^{-1}h_3(h_2h_3^{-1})^2h_2=1$:\\
Using Tietze transformation where $h_3\mapsto h_3h_2$ and $h_2\mapsto h_2$, we have $h_2h_3^{-1}h_2^{-1}h_3h_2(h_3^{-1})^2=1$. Using Tietze transformation again where $h_2\mapsto h_3h_2$ and $h_3\mapsto h_3$, we have $h_3^{-1}h_2h_3^{-1}h_2^{-1}h_3h_2=1$.  Now by using Tietze transformation again where $h_2\mapsto h_3h_2$ and $h_3\mapsto h_3$, we have:\\ $h_2h_3^{-1}h_2^{-1}h_3h_2=1\Rightarrow h_3^{-1}h_2h_3=h_2^2\Rightarrow \langle h_2,h_3\rangle\cong BS(1,2) \text{ is solvable, a contradiction}$.
\item[(226)]$h_2h_3^{-1}h_2^{-1}h_3^4=1$:\\
By interchanging $h_2$ and $h_3$ in (15) and with the same discussion, there is a contradiction.
\item[(227)]$h_2h_3^{-1}h_2^{-1}h_3^3h_2^{-1}h_3=1$:\\
By interchanging $h_2$ and $h_3$ in (108) and with the same discussion, there is a contradiction.
\item[(231)]$h_2h_3^{-1}h_2^{-1}h_3(h_3h_2^{-1})^2h_3=1$:\\
By using Tietze transformation where $h_3\mapsto h_3h_2$ and $h_2\mapsto h_2$, we have $h_3^{-1}h_2^{-1}h_3h_2h_3^3h_2=1$. Using Tietze transformation again where $h_2\mapsto h_3^{-1}h_2$ and $h_3\mapsto h_3$, we have $h_3^{-1}h_2^{-1}h_3h_2h_3^2h_2=1$.  Now by using Tietze transformation again where $h_2\mapsto h_3^{-1}h_2$ and $h_3\mapsto h_3$, we have $h_3^{-1}h_2^{-1}h_3h_2h_3h_2=1$. By using Tietze transformation again where $h_2\mapsto h_3^{-1}h_2$ and $h_3\mapsto h_3$, we have $h_3^{-1}h_2^{-1}h_3h_2^2=1\Rightarrow h_3^{-1}h_2h_3=h_2^2\Rightarrow \langle h_2,h_3\rangle\cong BS(1,2) \text{ is solvable, a contradiction}$.
\item[(241)]$h_2h_3^{-1}(h_2^{-1}h_3)^2h_2^{-2}h_3=1$:\\
By interchanging $h_2$ and $h_3$ in (231) and with the same discussion, there is a contradiction.
\item[(244)]$h_2h_3^{-1}(h_2^{-1}h_3)^4=1$:\\
Using Tietze transformation where $h_3\mapsto h_2h_3$ and $h_2\mapsto h_2$, we have:\\ $h_2h_3^{-1}h_2^{-1}h_3^4=1\Rightarrow h_2^{-1}h_3^4h_2=h_3\Rightarrow \langle h_2,h_3\rangle\cong BS(4,1) \text{ is solvable, a contradiction}$.
\item[(254)]$h_2h_3^{-4}h_2^{-1}h_3=1$:\\
By interchanging $h_2$ and $h_3$ in (6) and with the same discussion, there is a contradiction.
\item[(261)]$h_2h_3^{-3}h_2h_3^{-1}h_2^{-1}h_3=1$:\\
By interchanging $h_2$ and $h_3$ in (73) and with the same discussion, there is a contradiction.
\item[(265)]$(h_2h_3^{-2}h_2)^2=1$:\\
$(h_2h_3^{-2}h_2)^2=1 \text{ and } G \text{ is torsion-free}\Rightarrow h_2h_3^{-2}h_2=1\Rightarrow h_2^2=h_3^2 \Rightarrow \langle h_2,h_3\rangle\cong BS(1,-1)$ is solvable, a contradiction.
\item[(282)]$h_2h_3^{-1}(h_3^{-1}h_2)^2h_3^{-1}h_2^{-1}h_3=1$:\\
By interchanging $h_2$ and $h_3$ in (224) and with the same discussion, there is a contradiction.
\item[(285)]$h_2h_3^{-1}(h_3^{-1}h_2)^4=1$:\\
Using Tietze transformation where $h_2\mapsto h_3h_2$ and $h_3\mapsto h_3$, we have:\\ $h_3h_2h_3^{-1}h_2^4=1\Rightarrow h_3^{-1}h_2^{-4}h_3=h_2\Rightarrow \langle h_2,h_3\rangle\cong BS(-4,1) \text{ is solvable, a contradiction}$.
\item[(286)]$(h_2h_3^{-1}h_2)^3=1$:\\
$(h_2h_3^{-1}h_2)^3=1 \text{ and } G \text{ is torsion-free}\Rightarrow h_2h_3^{-1}h_2=1\Rightarrow h_3=h_2^2 \Rightarrow \langle h_2,h_3\rangle=\langle h_2\rangle$ is abelian, a contradiction.
\item[(292)]$(h_2h_3^{-1}h_2)^2h_3^{-1}h_2h_3^{-1}h_2=1$:\\
Using Tietze transformation where $h_3\mapsto h_2h_3$ and $h_2\mapsto h_2$, we have $(h_2h_3^{-1})^2h_3^{-2}=1$. Using Tietze transformation again where $h_2\mapsto h_2h_3$ and $h_3\mapsto h_3$, we have  $h_2^2=h_3^2\Rightarrow \langle h_2,h_3\rangle\cong BS(1,-1) \text{ is solvable, a contradiction}$.
\item[(293)]$(h_2h_3^{-1}h_2h_3)^2=1$:\\
By interchanging $h_2$ and $h_3$ in (196) and with the same discussion, there is a contradiction.
\item[(321)]$((h_2h_3^{-1})^2h_2)^2=1$:\\
$((h_2h_3^{-1})^2h_2)^2=1 \text{ and } G \text{ is torsion-free}\Rightarrow (h_2h_3^{-1})^2h_2=1\Rightarrow h_2=(h_3h_2^{-1})^2 \Rightarrow \langle h_2,h_3\rangle=\langle h_2,h_3h_2^{-1}\rangle=\langle h_3h_2^{-1}\rangle$ is abelian, a contradiction.
\item[(334)]$(h_2h_3^{-1})^4h_2^{-1}h_3=1$:\\
Using Tietze transformation where $h_2\mapsto h_2h_3$ and $h_3\mapsto h_3$, we have:\\ $h_2^4h_3^{-1}h_2^{-1}h_3=1\Rightarrow h_3^{-1}h_2h_3=h_2^4\Rightarrow \langle h_2,h_3\rangle\cong BS(1,4) \text{ is solvable, a contradiction}$.
\item[(335)]$(h_2h_3^{-1})^4h_3^{-1}h_2=1$:\\
Using Tietze transformation where $h_2\mapsto h_2h_3$ and $h_3\mapsto h_3$, we have:\\ $h_2^4h_3^{-1}h_2h_3=1\Rightarrow h_3^{-1}h_2h_3=h_2^{-4}\Rightarrow \langle h_2,h_3\rangle\cong BS(1,-4) \text{ is solvable, a contradiction}$.
\item[(337)]$(h_2h_3^{-1})^5h_2=1$:\\
$h_2=(h_3h_2^{-1})^5\Rightarrow \langle h_2,h_3\rangle=\langle h_2,h_3h_2^{-1}\rangle=\langle h_3h_2^{-1}\rangle \text{ is abelian, a contradiction}$.
\item[(338)]$h_3^6=1$:\\
By interchanging $h_2$ and $h_3$ in (1) and with the same discussion, there is a contradiction.
\item[(339)]$h_3^5h_2^{-1}h_3=1$:\\
By interchanging $h_2$ and $h_3$ in (3) and with the same discussion, there is a contradiction.
\item[(341)]$h_3(h_3^2h_2^{-1})^2h_3=1$:\\
By interchanging $h_2$ and $h_3$ in (27) and with the same discussion, there is a contradiction.
\item[(343)]$(h_3^2h_2^{-1}h_3)^2=1$:\\
By interchanging $h_2$ and $h_3$ in (96) and with the same discussion, there is a contradiction.
\item[(344)]$(h_3^2h_2^{-1})^2h_3h_2^{-1}h_3=1$:\\
By interchanging $h_2$ and $h_3$ in (102) and with the same discussion, there is a contradiction.
\item[(347)]$(h_3h_2^{-1}h_3)^3=1$:\\
By interchanging $h_2$ and $h_3$ in (286) and with the same discussion, there is a contradiction.
\item[(348)]$(h_3h_2^{-1}h_3)^2h_2^{-1}h_3h_2^{-1}h_3=1$:\\
By interchanging $h_2$ and $h_3$ in (292) and with the same discussion, there is a contradiction.
\item[(349)]$((h_3h_2^{-1})^2h_3)^2=1$:\\
By interchanging $h_2$ and $h_3$ in (321) and with the same discussion, there is a contradiction.
\item[(350)]$(h_3h_2^{-1})^5h_3=1$:\\
By interchanging $h_2$ and $h_3$ in (337) and with the same discussion, there is a contradiction.
\item[(351)]$(h_2^{-1}h_3)^6=1$:\\
$(h_2^{-1}h_3)^6=1 \text{ and } G \text{ is torsion-free}\Rightarrow h_2^{-1}h_3=1\Rightarrow h_2=h_3 \Rightarrow \langle h_2,h_3\rangle=\langle h_2\rangle$ is abelian, a contradiction.
\end{enumerate}


\subsection{$\mathbf{K_{2,3}}$}
By Theorem \ref{K2,3}, the Kaplansky graph $K(\alpha,\beta)$ contains no subgraph isomorphic to the complete bipartite graph $K_{2,3}$.
\begin{center}
\begin{figure}[t]
\psscalebox{0.85 0.85} 
{
\begin{pspicture}(0,-2.2985578)(9.194231,2.2985578)
\psdots[linecolor=black, dotsize=0.4](1.7971153,-1.0985577)
\psdots[linecolor=black, dotsize=0.4](2.5971153,-0.2985577)
\psdots[linecolor=black, dotsize=0.4](3.3971152,0.5014423)
\psdots[linecolor=black, dotsize=0.4](1.7971153,2.1014423)
\psdots[linecolor=black, dotsize=0.4](0.19711533,0.5014423)
\psdots[linecolor=black, dotsize=0.4](1.7971153,0.5014423)
\psdots[linecolor=black, dotsize=0.4](7.397115,2.1014423)
\psdots[linecolor=black, dotsize=0.4](7.397115,0.5014423)
\psdots[linecolor=black, dotsize=0.4](7.397115,-1.0985577)
\psdots[linecolor=black, dotsize=0.4](8.197115,-0.2985577)
\psdots[linecolor=black, dotsize=0.4](8.997115,0.5014423)
\psdots[linecolor=black, dotsize=0.4](8.197115,1.3014423)
\psdots[linecolor=black, dotsize=0.4](5.7971153,0.5014423)
\psline[linecolor=black, linewidth=0.04](1.7971153,2.1014423)(3.3971152,0.5014423)(2.5971153,-0.2985577)(1.7971153,-1.0985577)(0.19711533,0.5014423)(1.7971153,2.1014423)
\psline[linecolor=black, linewidth=0.04](1.7971153,2.1014423)(1.7971153,0.5014423)(1.7971153,-1.0985577)(1.7971153,-1.0985577)
\psline[linecolor=black, linewidth=0.04](5.7971153,0.5014423)(7.397115,2.1014423)(8.197115,1.3014423)(8.997115,0.5014423)(8.197115,-0.2985577)(7.397115,-1.0985577)(5.7971153,0.5014423)
\psline[linecolor=black, linewidth=0.04](7.397115,2.1014423)(7.397115,0.5014423)(7.397115,-1.0985577)(7.397115,-1.0985577)
\rput[bl](0.9971153,-2.2985578){$\mathbf{C_4--C_5}$}
\rput[bl](6.5971155,-2.2985578){$\mathbf{C_4--C_6}$}
\end{pspicture}
}
\caption{($C_4--C_5$) A $C_4$ and a $C_5$ cycle with two common edges, and ($C_4--C_6$) A $C_4$ and a $C_6$ cycle with two common edges in $K(\alpha,\beta)$}\label{f-5}
\end{figure}
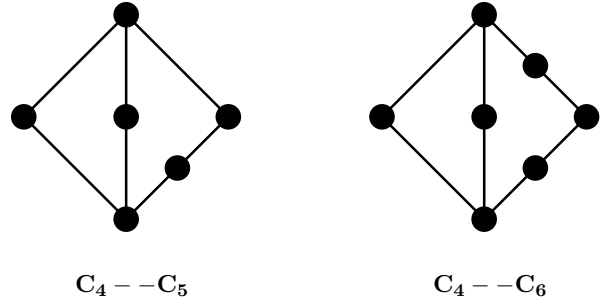
\end{center}
\subsection{$\mathbf{C_4--C_5}$}
Suppose that $[h_1',h_1'',h_2',h_2'',h_3',h_3'',h_4',h_4'']$ is the $8-$tuple related to the cycle $C_4$ and $[h_1',h_1'',h_2',h_2'',h_5',h_5'',h_6',h_6'',h_7',h_7'']$ is the $10-$tuple related to the cycle $C_5$ in the graph $C_4--C_5$, where the first four  components of these tuples are related to the common edges of $C_4$ and $C_5$. Without loss of generality we may assume that $h_1'=1$ and  $\alpha=1+h_2+h_3$. It can be seen that there are $121$ different cases for the relations of the cycles $C_4$ and $C_5$ in this structure. Using GAP \cite{gap}, we see that the groups with two generators $h_2$ and $h_3$ and two relations which are between $111$ cases of these $121$ cases are finite and solvable, that is a contradiction with the assumptions. So, there are just $10$ cases for the relations of the cycles $C_4$ and $C_5$ which may lead to the existence of a subgraph isomorphic to the graph $C_4--C_5$ in $K(\alpha,\beta)$. Such cases are listed in table \ref{tab-C4--C5}.  In the following, we show that all of these $10$ cases lead to contradictions and so, the graph $K(\alpha,\beta)$ contains no subgraph isomorphic to the graph $C_4--C_5$.

\begin{enumerate}
\item[(1)]$R_1:h_2^2h_3h_2^{-1}h_3=1$, $R_2: h_2^2(h_2h_3^{-1})^2h_2=1$:\\
Using Tietze transformation where $h_3\mapsto h_3h_2$ and $h_2\mapsto h_2$, we have:\\ $R_1:h_2^3=h_3^{-2} \text{ and } R_2: h_2^3=h_3^2 \Rightarrow h_3^4=1 \text{ and } G \text{ is torsion-free}\Rightarrow h_3=1$, a contradiction.
\item[(2)]$ R_1: h_2^2h_3h_2^{-1}h_3=1$, $R_2: h_2(h_3^{-1}h_2h_3^{-1})^2h_2=1$:\\
Using Tietze transformation where $h_3\mapsto h_3h_2$ and $h_2\mapsto h_2$, we have:\\ $R_1:h_2^3=h_3^{-2} \text{ and } R_2: h_2h_3^{-2}h_2^{-1}h_3^{-2}=1 \Rightarrow h_2h_2^3h_2^{-1}h_2^3=1 \Rightarrow h_2^6=1\text{ and } G \text{ is torsion-free}\Rightarrow h_2=1$, a contradiction.
\item[(3)]$ R_1: h_2^2h_3^{-2}h_2=1$, $R_2: h_2^3h_3^2=1$:\\
$R_1:h_2^3=h_3^{2} \text{ and } R_2: h_2^3=h_3^{-2} \Rightarrow h_3^4=1 \text{ and } G \text{ is torsion-free}\Rightarrow h_3=1$, a contradiction.
\item[(4)]$ R_1: h_2h_3h_2^{-2}h_3=1$, $R_2: h_2h_3^{-2}h_2^{-1}h_3^2=1$:\\
Using Tietze transformation where $h_2\mapsto h_3h_2$ and $h_3\mapsto h_3$, we have $R_1: h_3h_2h_3h_2^{-1}h_3^{-1}h_2^{-1}=1 \text{ and } R_2: h_2h_3^{-2}h_2^{-1}h_3^2=1$. Using Tietze transformation again where $h_3\mapsto h_3h_2^{-1}$ and $h_2\mapsto h_2$, we have $R_1: h_3h_2h_3^{-1}=h_2^{-1}h_3 \text{ and } R_2: (h_2^{-1}h_3)^2(h_2h_3^{-1})^2=1$. Using $R_1 \text{ and } R_2$,  we have $(h_3h_2h_3^{-1})^2=(h_3h_2^{-1})^2\Rightarrow \langle h_2,h_3\rangle=\langle h_3h_2h_3^{-1},h_3h_2^{-1}\rangle\cong BS(1,-1)$ is solvable, a contradiction.
\item[(5)]$ R_1: h_2(h_3h_2^{-1})^2h_3=1$, $R_2: h_2^2h_3^{-1}h_2^2h_3=1$:\\
Using Tietze transformation where $h_3\mapsto h_3h_2$ and $h_2\mapsto h_2$, we have:\\ $R_1:h_2^2=h_3^{-3} \text{ and } R_2: h_2^2h_3^{-1}h_2^2h_3=1 \Rightarrow h_3^6=1 \text{ and } G \text{ is torsion-free}\Rightarrow h_3=1$, a contradiction.
\item[(6)]$ R_1: h_2h_3^{-3}h_2=1$, $R_2: h_2^2h_3^3=1$:\\
By interchanging $h_2$ and $h_3$ in (3) and with the same discussion, there is a contradiction.
\item[(7)]$ R_1: h_2h_3^{-2}h_2h_3=1$, $R_2: h_2^2h_3^{-1}h_2^{-2}h_3=1$:\\
By interchanging $h_2$ and $h_3$ in (4) and with the same discussion, there is a contradiction.
\item[(8)]$ R_1: h_2h_3^{-1}h_2h_3^2=1$, $R_2: h_2(h_3^{-1}h_2h_3^{-1})^2h_2=1$:\\
By interchanging $h_2$ and $h_3$ in (2) and with the same discussion, there is a contradiction.
\item[(9)]$ R_1: h_2h_3^{-1}h_2h_3^2=1$, $R_2: h_3^2(h_3h_2^{-1})^2h_3=1$:\\
By interchanging $h_2$ and $h_3$ in (1) and with the same discussion, there is a contradiction.
\item[(10)]$ R_1: (h_2h_3^{-1})^2h_2h_3=1$, $R_2: h_2h_3^2h_2^{-1}h_3^2=1$:\\
By interchanging $h_2$ and $h_3$ in (5) and with the same discussion, there is a contradiction.
\end{enumerate}
\begin{table}[t]
\centering
\caption{The relations of a $C_4--C_5$ in the Kaplansky graph}\label{tab-C4--C5}
\begin{tabular}{|c|l|l|}\hline
$n$&$R_1$&$R_2$\\\hline
$1$&$h_2^2h_3h_2^{-1}h_3=1$&$h_2^2(h_2h_3^{-1})^2h_2=1$\\
$2$&$h_2^2h_3h_2^{-1}h_3=1$&$h_2(h_3^{-1}h_2h_3^{-1})^2h_2=1$\\
$3$&$h_2^2h_3^{-2}h_2=1$&$h_2^3h_3^2=1$\\
$4$&$h_2h_3h_2^{-2}h_3=1$&$h_2h_3^{-2}h_2^{-1}h_3^2=1$\\
$5$&$h_2(h_3h_2^{-1})^2h_3=1$&$h_2^2h_3^{-1}h_2^2h_3=1$\\
$6$&$h_2h_3^{-3}h_2=1$&$h_2^2h_3^3=1$\\
$7$&$h_2h_3^{-2}h_2h_3=1$&$h_2^2h_3^{-1}h_2^{-2}h_3=1$\\
$8$&$h_2h_3^{-1}h_2h_3^2=1$&$h_2(h_3^{-1}h_2h_3^{-1})^2h_2=1$\\
$9$&$h_2h_3^{-1}h_2h_3^2=1$&$h_3^2(h_3h_2^{-1})^2h_3=1$\\
$10$&$(h_2h_3^{-1})^2h_2h_3=1$&$h_2h_3^2h_2^{-1}h_3^2=1$\\\hline
\end{tabular}
\end{table}
\subsection{$\mathbf{C_4--C_6}$}
Suppose that $[h_1',h_1'',h_2',h_2'',h_3',h_3'',h_4',h_4'']$ is the $8-$tuple related to the cycle $C_4$ and $[h_1',h_1'',h_2',h_2'',h_5',h_5'',h_6',h_6'',h_7',h_7'',h_8',h_8'']$ is the $12-$tuple related to the cycle $C_6$ in the graph $C_4--C_6$, where the first four components of these tuples are related to the common edges of $C_4$ and $C_6$. Without loss of generality we may assume that $h_1'=1$ and $\alpha=1+h_2+h_3$. Also it is easy to see that $h_3' \neq h_5'$ and $h_4'' \neq h_8''$. With these assumptions and by considering the relations from Tables \ref{tab-C4} and \ref{tab-C6} which are not disproved, it can be seen that there are $658$ different cases for the relations of the cycles $C_4$ and $C_6$ in this structure. By considering all groups with two generators $h_2$ and $h_3$ and two relations which are between these cases and by using GAP \cite{gap}, we see that $632$ groups are finite and solvable, or just finite. So, there are just $20$ cases for the relations of the cycles $C_4$ and $C_6$ which may lead to the existence of a subgraph isomorphic to the graph $C_4--C_6$ in $K(\alpha,\beta)$. These cases are listed in table \ref{tab-C4--C6}.  In the following, we show that all of such $20$ cases lead to  contradictions and so, the graph $K(\alpha,\beta)$ contains no subgraph isomorphic to the graph $C_4--C_6$.

\begin{table}[ht]
\centering
\caption{The relations of a $C_4--C_6$ in the Kaplansky graph}\label{tab-C4--C6}
\begin{tabular}{|c|l|l|}\hline
$n$&$R_1$&$R_2$\\\hline
$1$&$h_2^2h_3h_2^{-1}h_3=1$&$h_2h_3h_2h_3^{-1}h_2^{-1}h_3^2=1$\\
$2$&$h_2^2h_3h_2^{-1}h_3=1$&$h_2h_3^2h_2^{-1}h_3^{-1}h_2h_3=1$\\
$3$&$h_2^2h_3^{-2}h_2=1$&$h_2h_3h_2^{-1}h_3^{-1}h_2h_3^{-1}h_2^{-1}h_3=1$\\
$4$&$h_2^2h_3^{-2}h_2=1$&$h_2h_3h_2^{-1}(h_3^{-1}h_2)^3=1$\\
$5$&$h_2^2h_3^{-2}h_2=1$&$h_2h_3^{-1}(h_2^{-1}h_3h_2^{-1})^2h_3=1$\\
$6$&$h_2^2h_3^{-2}h_2=1$&$h_2h_3^{-1}h_2^{-1}h_3h_2^{-1}h_3^{-1}h_2h_3=1$\\
$7$&$h_2h_3h_2^{-2}h_3=1$&$h_2h_3^2h_2^{-1}h_3^{-1}h_2^{-1}h_3=1$\\
$8$&$h_2h_3h_2^{-2}h_3=1$&$h_2h_3^{-1}h_2^{-1}h_3^2h_2^{-1}h_3^{-1}h_2=1$\\
$9$&$h_2(h_3h_2^{-1})^2h_3=1$&$h_2h_3^{-1}(h_2^{-1}h_3^2)^2=1$\\
$10$&$h_2(h_3h_2^{-1})^2h_3=1$&$(h_2h_3^{-2})^2h_2^{-1}h_3=1$\\
$11$&$h_2h_3^{-3}h_2=1$&$h_2h_3h_2^{-1}h_3^{-1}h_2h_3^{-1}h_2^{-1}h_3=1$\\
$12$&$h_2h_3^{-3}h_2=1$&$h_2h_3^{-1}h_2^{-1}h_3h_2^{-1}h_3^{-1}h_2h_3=1$\\
$13$&$h_2h_3^{-3}h_2=1$&$h_2h_3^{-1}h_2^{-1}(h_3h_2^{-1}h_3)^2=1$\\
$14$&$h_2h_3^{-3}h_2=1$&$h_2h_3^{-1}(h_2^{-1}h_3)^3h_3=1$\\
$15$&$h_2h_3^{-2}h_2h_3=1$&$h_2^2h_3^{-1}h_2^{-1}h_3^{-1}h_2h_3=1$\\
$16$&$h_2h_3^{-2}h_2h_3=1$&$h_2h_3^{-1}h_2^{-1}h_3^2h_2^{-1}h_3^{-1}h_2=1$\\
$17$&$h_2h_3^{-1}h_2h_3^2=1$&$h_2^2h_3^{-1}h_2^{-1}h_3h_2h_3=1$\\
$18$&$h_2h_3^{-1}h_2h_3^2=1$&$(h_2h_3)^2h_2^{-1}h_3^{-1}h_2=1$\\
$19$&$(h_2h_3^{-1})^2h_2h_3=1$&$(h_2^2h_3^{-1})^2h_2^{-1}h_3=1$\\
$20$&$(h_2h_3^{-1})^2h_2h_3=1$&$h_2h_3h_2^{-1}h_3^{-1}h_2^2h_3^{-1}h_2=1$\\\hline
\end{tabular}
\end{table}
\begin{enumerate}
\item[(1)]$R_1: h_2^2h_3h_2^{-1}h_3=1$, $R_2:h_2h_3h_2h_3^{-1}h_2^{-1}h_3^2=1$:\\
$R_2:h_2h_3h_2h_3^{-1}h_2^{-1}h_3^2=1\Rightarrow h_2h_3h_2^{-1}=h_3^2h_2h_3 (*)$. Using $R_1$ and $(*)$ we have $(h_3^2h_2)^2=1 \text{ and } G \text{ is torsion-free}\Rightarrow h_3^2h_2=1\Rightarrow h_2=h_3^{-2}\Rightarrow \langle h_2,h_3\rangle=\langle h_3\rangle \text{ is abelian, a contradiction}$.
\item[(2)]$R_1: h_2^2h_3h_2^{-1}h_3=1$, $R_2:h_2h_3^2h_2^{-1}h_3^{-1}h_2h_3=1$:\\
$R_2:h_2h_3^2h_2^{-1}h_3^{-1}h_2h_3=1\Rightarrow h_2^{-1}h_3h_2=h_3h_2h_3^2 (*)$. Using $R_1$ and $(*)$ we have $(h_2h_3^2)^2=1 \text{ and } G \text{ is torsion-free}\Rightarrow h_2h_3^2=1\Rightarrow h_2=h_3^{-2}\Rightarrow \langle h_2,h_3\rangle=\langle h_3\rangle \text{ is abelian, a contradiction}$.
\item[(3)]$R_1: h_2^2h_3^{-2}h_2=1$, $R_2:h_2h_3h_2^{-1}h_3^{-1}h_2h_3^{-1}h_2^{-1}h_3=1$:\\
$R_1: h_2^2h_3^{-2}h_2=1\Rightarrow h_2^{2}h_3^{-1}=h_2^{-1}h_3 (*) \text{ and } h_3^{2}=h_2^{3} (**)$. $R_2:h_2h_3h_2^{-1}h_3^{-1}h_2h_3^{-1}h_2^{-1}h_3=1\Rightarrow h_3h_2^{-1}h_3^{-1}h_2h_3^{-1}h_2^{-1}h_3h_2=1\Rightarrow h_3^2h_2^{-1}h_3^{-1}h_2h_3^{-1}h_2^{-1}h_3h_2h_3^{-1}=1 (***)$. Using $(**)$ and $(***)$ we have $h_2^2h_3^{-1}h_2h_3^{-1}h_2^{-1}h_3h_2h_3^{-1}=1\Rightarrow \text{By } (*), h_2^{-1}h_3h_2h_3^{-1}h_2^{-1}h_3h_2h_3^{-1}=1\Rightarrow (h_2^{-1}h_3h_2h_3^{-1})^2=1 \text{ and } G \text{ is torsion-free}\Rightarrow h_2^{-1}h_3h_2=h_3\Rightarrow \langle h_2,h_3\rangle\cong BS(1,1)$ is solvable, a contradiction.
\item[(4)]$R_1: h_2^2h_3^{-2}h_2=1$, $R_2:h_2h_3h_2^{-1}(h_3^{-1}h_2)^3=1$:\\
$R_1: h_2^2h_3^{-2}h_2=1\Rightarrow h_2^{-2}=h_2h_3^{-2} (*) \text{ and } h_3^{-2}=h_2^{-3} (**)$. Using $R_2$ and $(*)$ we have $(h_2h_3^{-1})^2(h_2^{-1}h_3)^2=1\Rightarrow (h_2h_3^{-1})^2=(h_3^{-1}h_2)^2$. By $(**)$, $(h_2h_3^{-1})(h_3^{-1}h_2)=h_2h_3^{-2}h_2=h_2^{-1}\Rightarrow \langle h_2,h_3\rangle=\langle h_2h_3^{-1},h_3^{-1}h_2\rangle\cong BS(1,-1)$ is solvable, a contradiction.
\item[(5)]$R_1: h_2^2h_3^{-2}h_2=1$, $R_2:h_2h_3^{-1}(h_2^{-1}h_3h_2^{-1})^2h_3=1$:\\
$R_1: h_2^2h_3^{-2}h_2=1\Rightarrow h_2^{-2}=h_2h_3^{-2} (*)$. Using $R_2$ and $(*)$ we have $h_2h_3^{-1}h_2^{-1}h_3h_2h_3^{-2}h_3h_2^{-1}h_3=1\Rightarrow (h_2h_3^{-1}h_2^{-1}h_3)^2=1\text{ and } G \text{ is torsion-free}\Rightarrow h_2h_3^{-1}h_2^{-1}h_3=1\Rightarrow h_2^{-1}h_3h_2=h_3\Rightarrow \langle h_2,h_3\rangle\cong BS(1,1)$ is solvable, a contradiction.
\item[(6)]$R_1: h_2^2h_3^{-2}h_2=1$, $R_2:h_2h_3^{-1}h_2^{-1}h_3h_2^{-1}h_3^{-1}h_2h_3=1$:\\
$R_1: h_2^2h_3^{-2}h_2=1\Rightarrow h_2^{2}=h_3^{2}h_2^{-1} (*) \text{ and } h_3^{2}=h_2^{3} (**)$. Using $R_2$, $(*)$ and $(**)$ we have $(h_2h_3^{-1})^2(h_2^{-1}h_3)^2=1\Rightarrow (h_2h_3^{-1})^2=(h_3^{-1}h_2)^2$. By $(**)$, $(h_2h_3^{-1})(h_3^{-1}h_2)=h_2h_3^{-2}h_2=h_2^{-1}\Rightarrow \langle h_2,h_3\rangle=\langle h_2h_3^{-1},h_3^{-1}h_2\rangle\cong BS(1,-1)$ is solvable, a contradiction.
\item[(7)]$R_1: h_2h_3h_2^{-2}h_3=1$, $R_2:h_2h_3^2h_2^{-1}h_3^{-1}h_2^{-1}h_3=1$:\\
$R_1: h_2h_3h_2^{-2}h_3=1\Rightarrow h_3h_2h_3=h_2^2 \text{ and } R_2:h_2h_3^2h_2^{-1}h_3^{-1}h_2^{-1}h_3=1\Rightarrow h_3h_2^{-1}h_3^{-1}h_2=1\Rightarrow h_3^{-1}h_2h_3=h_2\Rightarrow \langle h_2,h_3\rangle\cong BS(1,1)$ is solvable, a contradiction.
\item[(8)]$R_1: h_2h_3h_2^{-2}h_3=1$, $R_2:h_2h_3^{-1}h_2^{-1}h_3^2h_2^{-1}h_3^{-1}h_2=1$:\\
$R_1: h_2h_3h_2^{-2}h_3=1\Rightarrow h_2^{-1}h_3^{-1}h_2^{2}h_3^{-1}=1 (*)$. Using $R_2$ and $(*)$ we have $h_2^{-1}h_3^2=1\Rightarrow h_2=h_3^{2}\Rightarrow \langle h_2,h_3\rangle=\langle h_3\rangle \text{ is abelian, a contradiction}$.
\item[(9)]$R_1: h_2(h_3h_2^{-1})^2h_3=1$, $R_2:h_2h_3^{-1}(h_2^{-1}h_3^2)^2=1$:\\
$R_2:h_2h_3^{-1}(h_2^{-1}h_3^2)^2=1\Rightarrow h_2h_3^{-1}h_2^{-1}h_3^2h_2^{-1}h_3^2=1\Rightarrow h_2h_3h_2^{-1}h_3^{-2}h_2h_3^{-2}=1\Rightarrow h_2h_3h_2^{-1}=h_3^{2}h_2^{-1}h_3^{2} (*)$. By using $R_1$ and $(*)$ we have $(h_3^3h_2^{-1})^2=1 \text{ and } G \text{ is torsion-free}\Rightarrow h_3^3h_2^{-1}=1\Rightarrow h_2=h_3^{3}\Rightarrow \langle h_2,h_3\rangle=\langle h_3\rangle \text{ is abelian, a contradiction}$.
\item[(10)]$R_1: h_2(h_3h_2^{-1})^2h_3=1$, $R_2:(h_2h_3^{-2})^2h_2^{-1}h_3=1$:\\
$R_2:(h_2h_3^{-2})^2h_2^{-1}h_3=1\Rightarrow h_2^{-1}h_3h_2=h_3^{2}h_2^{-1}h_3^{2} (*)$. By using $R_1$ and $(*)$ we have $(h_2^{-1}h_3^3)^2=1 \text{ and } G \text{ is torsion-free}\Rightarrow h_2^{-1}h_3^3=1\Rightarrow h_2=h_3^{3}\Rightarrow \langle h_2,h_3\rangle=\langle h_3\rangle \text{ is abelian, a contradiction}$.
\item[(11)]$R_1: h_2h_3^{-3}h_2=1$, $R_2:h_2h_3h_2^{-1}h_3^{-1}h_2h_3^{-1}h_2^{-1}h_3=1$:\\
By interchanging $h_2$ and $h_3$ in (6) and with the same discussion, there is a contradiction.
\item[(12)]$R_1: h_2h_3^{-3}h_2=1$, $R_2:h_2h_3^{-1}h_2^{-1}h_3h_2^{-1}h_3^{-1}h_2h_3=1$:\\
By interchanging $h_2$ and $h_3$ in (3) and with the same discussion, there is a contradiction.
\item[(13)]$R_1: h_2h_3^{-3}h_2=1$, $R_2:h_2h_3^{-1}h_2^{-1}(h_3h_2^{-1}h_3)^2=1$:\\
By interchanging $h_2$ and $h_3$ in (5) and with the same discussion, there is a contradiction.
\item[(14)]$R_1: h_2h_3^{-3}h_2=1$, $R_2:h_2h_3^{-1}(h_2^{-1}h_3)^3h_3=1$:\\
By interchanging $h_2$ and $h_3$ in (4) and with the same discussion, there is a contradiction.
\item[(15)]$R_1: h_2h_3^{-2}h_2h_3=1$, $R_2:h_2^2h_3^{-1}h_2^{-1}h_3^{-1}h_2h_3=1$:\\
By interchanging $h_2$ and $h_3$ in (7) and with the same discussion, there is a contradiction.
\item[(16)]$R_1: h_2h_3^{-2}h_2h_3=1$, $R_2:h_2h_3^{-1}h_2^{-1}h_3^2h_2^{-1}h_3^{-1}h_2=1$:\\
By interchanging $h_2$ and $h_3$ in (8) and with the same discussion, there is a contradiction.
\item[(17)]$R_1: h_2h_3^{-1}h_2h_3^2=1$, $R_2:h_2^2h_3^{-1}h_2^{-1}h_3h_2h_3=1$:\\
By interchanging $h_2$ and $h_3$ in (2) and with the same discussion, there is a contradiction.
\item[(18)]$R_1: h_2h_3^{-1}h_2h_3^2=1$, $R_2:(h_2h_3)^2h_2^{-1}h_3^{-1}h_2=1$:\\
By interchanging $h_2$ and $h_3$ in (1) and with the same discussion, there is a contradiction.
\item[(19)]$R_1: (h_2h_3^{-1})^2h_2h_3=1$, $R_2:(h_2^2h_3^{-1})^2h_2^{-1}h_3=1$:\\
By interchanging $h_2$ and $h_3$ in (10) and with the same discussion, there is a contradiction.
\item[(20)]$R_1: (h_2h_3^{-1})^2h_2h_3=1$, $R_2:h_2h_3h_2^{-1}h_3^{-1}h_2^2h_3^{-1}h_2=1$:\\
By interchanging $h_2$ and $h_3$ in (9) and with the same discussion, there is a contradiction.
\end{enumerate}

$\mathbf{C_4-C_5}$ \textbf{subgraph:} Suppose that $[h_1',h_1'',h_2',h_2'',h_3',h_3'',h_4',h_4'']$ is the $8-$tuple related to the cycle $C_4$ and $[h_1',h_1'',h_5',h_5'',h_6',h_6'',h_7',h_7'',h_8',h_8'']$ is the $10-$tuple related to the cycle $C_5$ in the graph $C_4-C_5$, where the first two components of these tuples are related to the common edge of $C_4$ and $C_5$. With the same argument such as about $C_4--C_5$, without loss of generality we may assume that $h_1'=1$, where $h_1'',h_2',h_2'',h_3',h_3'',h_4',h_4'',h_5',h_5'',h_6',h_6'',h_7',h_7'',h_8',h_8'' \in supp(\alpha)$ and $\alpha=1+h_2+h_3$. Also it is easy to see that $h_2' \neq h_5'$ and $h_4'' \neq h_8''$.  With these assumptions and by considering the relations from Tables \ref{tab-C4} and \ref{tab-C5} which are not disproved, it can be seen that there are $482$ cases for the relations of the cycles $C_4$ and $C_5$ in this structure. Using Gap \cite{gap}, we see that all groups with two generators $h_2$ and $h_3$ and two relations which are between $426$ cases of these $482$ cases are finite and solvable, that is a contradiction with the assumptions. So there are just $56$ cases for the relations of the cycles $C_4$ and $C_5$ which may lead to the existence of a subgraph isomorphic to the graph $C_4-C_5$ in the graph $K(\alpha,\beta)$. These cases are listed in table \ref{tab-C4-C5}.  In the following, we show that $50$ cases of these relations lead to a contradiction and just $6$ cases of them may lead to the  existence of a subgraph isomorphic to the graph $C_4-C_5$ in the graph $K(\alpha,\beta)$. Cases which are not disproved are marked by $*$s in the Table \ref{tab-C4-C5}.
\begin{figure}[h]
\psscalebox{1.0 1.0} 
{
\begin{pspicture}(0,-4.0985575)(11.594231,4.0985575)
\psdots[linecolor=black, dotsize=0.3](0.9971153,3.1014423)
\psdots[linecolor=black, dotsize=0.3](1.7971153,3.5014422)
\psdots[linecolor=black, dotsize=0.3](2.5971153,3.1014423)
\psdots[linecolor=black, dotsize=0.3](2.1971154,2.3014424)
\psdots[linecolor=black, dotsize=0.3](1.3971153,2.3014424)
\psdots[linecolor=black, dotsize=0.3](0.59711534,1.9014423)
\psdots[linecolor=black, dotsize=0.3](0.19711533,2.7014422)
\psline[linecolor=black, linewidth=0.04](0.19711533,2.7014422)(0.9971153,3.1014423)(1.7971153,3.5014422)(2.5971153,3.1014423)(2.1971154,2.3014424)(1.3971153,2.3014424)(0.9971153,3.1014423)(0.9971153,3.1014423)
\psline[linecolor=black, linewidth=0.04](1.3971153,2.3014424)(0.59711534,1.9014423)(0.19711533,2.7014422)(0.19711533,2.7014422)
\rput[bl](0.59711534,0.7014423){$\mathbf{C_4-C_5}$}
\psdots[linecolor=black, dotsize=0.3](4.5971155,2.7014422)
\psdots[linecolor=black, dotsize=0.3](4.5971155,1.9014423)
\psdots[linecolor=black, dotsize=0.3](5.397115,1.9014423)
\psdots[linecolor=black, dotsize=0.3](5.397115,2.7014422)
\psdots[linecolor=black, dotsize=0.3](6.1971154,2.7014422)
\psdots[linecolor=black, dotsize=0.3](6.1971154,1.9014423)
\psdots[linecolor=black, dotsize=0.3](4.997115,3.5014422)
\psdots[linecolor=black, dotsize=0.3](5.7971153,3.5014422)
\psline[linecolor=black, linewidth=0.04](6.1971154,2.7014422)(6.1971154,1.9014423)(5.397115,1.9014423)(4.5971155,1.9014423)(4.5971155,2.7014422)(5.397115,2.7014422)(6.1971154,2.7014422)(6.1971154,2.7014422)
\psline[linecolor=black, linewidth=0.04](5.397115,1.9014423)(5.397115,2.7014422)(5.397115,2.7014422)
\psline[linecolor=black, linewidth=0.04](6.1971154,2.7014422)(5.7971153,3.5014422)(4.997115,3.5014422)(4.5971155,2.7014422)
\rput[bl](4.1971154,0.7014423){$\mathbf{C_4-C_5(-C_4-)}$}
\psdots[linecolor=black, dotsize=0.3](9.797115,3.5014422)
\psdots[linecolor=black, dotsize=0.3](9.797115,2.7014422)
\psdots[linecolor=black, dotsize=0.3](8.997115,3.9014423)
\psdots[linecolor=black, dotsize=0.3](8.197115,3.1014423)
\psdots[linecolor=black, dotsize=0.3](8.5971155,2.3014424)
\psdots[linecolor=black, dotsize=0.3](10.997115,2.3014424)
\psdots[linecolor=black, dotsize=0.3](11.397116,3.1014423)
\psdots[linecolor=black, dotsize=0.3](10.5971155,3.9014423)
\psdots[linecolor=black, dotsize=0.3](9.797115,1.5014423)
\psline[linecolor=black, linewidth=0.04](9.797115,1.5014423)(10.997115,2.3014424)(9.797115,2.7014422)(8.5971155,2.3014424)(9.797115,1.5014423)
\psline[linecolor=black, linewidth=0.04](8.5971155,2.3014424)(8.197115,3.1014423)(8.997115,3.9014423)(9.797115,3.5014422)(9.797115,2.7014422)(9.797115,2.7014422)
\psline[linecolor=black, linewidth=0.04](10.997115,2.3014424)(11.397116,3.1014423)(10.5971155,3.9014423)(9.797115,3.5014422)(9.797115,3.5014422)
\rput[bl](8.5971155,0.7014423){$\mathbf{C_4-C_5(-C_5-)}$}
\rput[bl](0.19711533,-4.0985575){$\mathbf{C_4-C_5(-C_6--)}$}
\rput[bl](4.5971155,-4.0985575){$\mathbf{C_4-C_5(-C_6-)}$}
\rput[bl](8.197115,-4.0985575){$\mathbf{C_4-C_5(-C_7--)}$}
\psdots[linecolor=black, dotsize=0.3](1.7971153,-0.4985577)
\psdots[linecolor=black, dotsize=0.3](2.5971153,-0.89855766)
\psdots[linecolor=black, dotsize=0.3](0.9971153,-0.89855766)
\psdots[linecolor=black, dotsize=0.3](1.3971153,-1.6985577)
\psdots[linecolor=black, dotsize=0.3](2.1971154,-1.6985577)
\psdots[linecolor=black, dotsize=0.3](0.59711534,-2.0985577)
\psdots[linecolor=black, dotsize=0.3](0.19711533,-1.2985576)
\psdots[linecolor=black, dotsize=0.3](0.9971153,-2.8985577)
\psdots[linecolor=black, dotsize=0.3](2.9971154,-2.0985577)
\psline[linecolor=black, linewidth=0.04](0.19711533,-1.2985576)(0.9971153,-0.89855766)(1.7971153,-0.4985577)(2.5971153,-0.89855766)(2.1971154,-1.6985577)(1.3971153,-1.6985577)(0.9971153,-0.89855766)(0.9971153,-0.89855766)
\psline[linecolor=black, linewidth=0.04](1.3971153,-1.6985577)(0.59711534,-2.0985577)(0.19711533,-1.2985576)(0.19711533,-1.2985576)
\psline[linecolor=black, linewidth=0.04](0.59711534,-2.0985577)(0.9971153,-2.8985577)(2.9971154,-2.0985577)(2.5971153,-0.89855766)
\psdots[linecolor=black, dotsize=0.3](6.1971154,-0.4985577)
\psdots[linecolor=black, dotsize=0.3](6.997115,-0.89855766)
\psdots[linecolor=black, dotsize=0.3](5.397115,-0.89855766)
\psdots[linecolor=black, dotsize=0.3](5.7971153,-1.6985577)
\psdots[linecolor=black, dotsize=0.3](6.5971155,-1.6985577)
\psdots[linecolor=black, dotsize=0.3](4.997115,-2.0985577)
\psdots[linecolor=black, dotsize=0.3](4.5971155,-1.2985576)
\psdots[linecolor=black, dotsize=0.3](4.997115,-2.8985577)
\psdots[linecolor=black, dotsize=0.3](5.7971153,-2.8985577)
\psdots[linecolor=black, dotsize=0.3](6.5971155,-2.4985578)
\psline[linecolor=black, linewidth=0.04](4.5971155,-1.2985576)(5.397115,-0.89855766)(6.1971154,-0.4985577)(6.997115,-0.89855766)(6.5971155,-1.6985577)(5.7971153,-1.6985577)(5.397115,-0.89855766)(5.397115,-0.89855766)
\psline[linecolor=black, linewidth=0.04](5.7971153,-1.6985577)(4.997115,-2.0985577)(4.5971155,-1.2985576)(4.5971155,-1.2985576)
\psline[linecolor=black, linewidth=0.04](4.997115,-2.0985577)(4.997115,-2.8985577)(5.7971153,-2.8985577)(6.5971155,-2.4985578)(6.5971155,-1.6985577)
\psdots[linecolor=black, dotsize=0.3](9.797115,-0.4985577)
\psdots[linecolor=black, dotsize=0.3](10.5971155,-0.89855766)
\psdots[linecolor=black, dotsize=0.3](8.997115,-0.89855766)
\psdots[linecolor=black, dotsize=0.3](9.397116,-1.6985577)
\psdots[linecolor=black, dotsize=0.3](10.197115,-1.6985577)
\psdots[linecolor=black, dotsize=0.3](8.5971155,-2.0985577)
\psdots[linecolor=black, dotsize=0.3](8.197115,-1.2985576)
\psdots[linecolor=black, dotsize=0.3](8.997115,-2.8985577)
\psdots[linecolor=black, dotsize=0.3](10.197115,-2.4985578)
\psdots[linecolor=black, dotsize=0.3](11.397116,-1.6985577)
\psline[linecolor=black, linewidth=0.04](8.5971155,-2.0985577)(8.197115,-1.2985576)(8.997115,-0.89855766)(9.797115,-0.4985577)(10.5971155,-0.89855766)(11.397116,-1.6985577)(10.197115,-2.4985578)(8.997115,-2.8985577)(8.5971155,-2.0985577)(9.397116,-1.6985577)(10.197115,-1.6985577)(10.5971155,-0.89855766)(10.5971155,-0.89855766)
\psline[linecolor=black, linewidth=0.04](9.397116,-1.6985577)(8.997115,-0.89855766)
\end{pspicture}
}
\caption{The graph $C_4-C_5$ and some forbidden subgraphs which contain such graph}\label{f-6}
\end{figure}
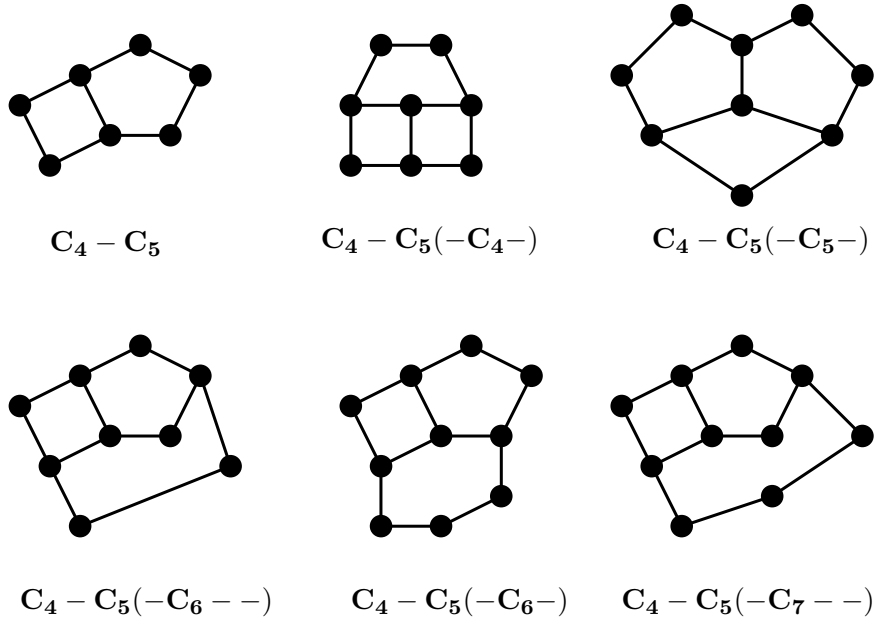
\begin{enumerate}
\item[(1)]$ R_1: h_2^2h_3h_2^{-1}h_3=1$, $R_2: h_2^2(h_2h_3^{-1})^2h_2=1$:\\
$ R_1: h_2^2h_3h_2^{-1}h_3=1\Rightarrow h_2^2=h_3^{-1}h_2h_3^{-1} \ (*)\Rightarrow R_2: h_2h_3^{-1}h_2h_3^{-1}h_2h_3^{-1}h_2h_3^{-1}=1\Rightarrow (h_2h_3^{-1})^4=1 \text{ and } G \text{ is torsion-free}\Rightarrow h_2h_3^{-1}=1\Rightarrow h_2=h_3\Rightarrow \langle h_2,h_3\rangle=\langle h_2\rangle \text{ is abelian, a contradiction}$.

\item[(2)]$ R_1: h_2^2h_3h_2^{-1}h_3=1$, $R_2: h_2^2h_3^{-1}h_2^{-2}h_3=1$:\\
$R_2: h_2^2h_3^{-1}h_2^{-2}h_3=1\Rightarrow h_2^2h_3=h_3h_2^2 \ (*)\Rightarrow R_1: h_3h_2^2h_2^{-1}h_3=1\Rightarrow h_2=h_3^{-2}\Rightarrow \langle h_2,h_3\rangle=\langle h_3\rangle \text{ is abelian, a contradiction}$.

\item[(3)]$ R_1: h_2^2h_3h_2^{-1}h_3=1$, $R_2: h_2(h_3^{-1}h_2h_3^{-1})^2h_2=1$:\\
$ R_1: h_2^2h_3h_2^{-1}h_3=1\Rightarrow h_2^2=h_3^{-1}h_2h_3^{-1} \ (*)\Rightarrow R_2: h_2^6=1\Rightarrow h_2=1$, a contradiction.

\item[(4)]$ R_1: h_2^2h_3h_2^{-1}h_3=1$, $R_2: h_2h_3^{-2}h_2^{-1}h_3^2=1$:\\
$R_2: h_2h_3^{-2}h_2^{-1}h_3^2=1\Rightarrow  h_3^2h_2^2=h_2^2h_3^2  \ (*) \text{ and } R_1: h_2^2h_3h_2^{-1}h_3=1\Rightarrow (h_3h_2^2h_3h_2^{-1})(h_2^2h_3h_2^{-1}h_3)=1  \ (**).$ By $(*)$ and $(**)$ we have $h_2h_3^3h_2h_3=1\Rightarrow h_3^2(h_2h_3)^2=1\Rightarrow h_3^2=(h_3^{-1}h_2^{-1})^2 \Rightarrow \langle h_2,h_3\rangle=\langle h_3,h_3^{-1}h_2^{-1}\rangle\cong BS(1,-1)$ is solvable, a contradiction.

\item[(5)]$ R_1: h_2^2h_3h_2^{-1}h_3=1$, $R_2: h_2h_3^{-1}h_2h_3h_2^{-1}h_3^2=1$:\\
$ R_1: h_2^2h_3h_2^{-1}h_3=1\Rightarrow h_2h_3h_2^{-1}h_3=h_2^{-1} \ (*) \Rightarrow R_2: h_2h_3^{-1}h_2^{-1}h_3=1\Rightarrow h_2^{-1}h_3h_2=h_3\Rightarrow \langle h_2,h_3\rangle\cong BS(1,1)$ is solvable, a contradiction.

\item[(6)]$ R_1: h_2^2h_3h_2^{-1}h_3=1$, $R_2: h_2h_3^{-1}h_2h_3^2h_2^{-1}h_3=1$:\\
$ R_1: h_2^2h_3h_2^{-1}h_3=1\Rightarrow h_3h_2^{-1}h_3h_2=h_2^{-1} \ (*) \Rightarrow R_2:h_3^{-1}h_2h_3h_2^{-1}=1\Rightarrow h_3^{-1}h_2h_3=h_2\Rightarrow \langle h_2,h_3\rangle\cong BS(1,1)$ is solvable, a contradiction.

%

\item[(9)]$ R_1: h_2^2h_3^{-2}h_2=1$, $R_2: h_2h_3^{-1}(h_3^{-1}h_2)^2h_3=1$:\\
$ R_1: h_2^2h_3^{-2}h_2=1\Rightarrow h_3^{-2}=h_2^{-3} \ (*)  \text{ and } R_2: h_2h_3^{-2}h_2h_3^{-1}h_2h_3=1 \ (**)$. By $(*)$ and $(**)$ we have $h_2^{-1}h_3^{-1}h_2h_3=1\Rightarrow h_2^{-1}h_3h_2=h_3\Rightarrow \langle h_2,h_3\rangle\cong BS(1,1)$ is solvable, a contradiction.

\item[(10)]$ R_1: h_2^2h_3^{-2}h_2=1$, $R_2: (h_2h_3^{-1})^2(h_2^{-1}h_3)^2=1$:\\
Using Tietze transformation where $h_3\mapsto h_3h_2$ and $h_2\mapsto h_2$, we have $ R_1: h_2^2h_3^{-1}h_2^{-1}h_3^{-1}=1$ and $R_2:h_2h_3^{-2}h_2^{-1}h_3^2=1$. $R_2\Rightarrow h_2^2h_3^{-2}=h_3^{-2}h_2^2 \ (*)   \text{ and } R_1\Rightarrow (h_3^{-1}h_2^2h_3^{-1}h_2^{-1})(h_2^2h_3^{-1}h_2^{-1}h_3^{-1})=1 \ (**)$. By $(*)$ and $(**)$ we have $h_2h_3^{-3}h_2h_3^{-1}=1\Rightarrow h_3^{-2}(h_2h_3^{-1})^2=1\Rightarrow h_3^2=(h_2h_3^{-1})^2 \Rightarrow \langle h_2,h_3\rangle=\langle h_3,h_2h_3^{-1}\rangle\cong BS(1,-1)$ is solvable, a contradiction.

\item[(11)]$ R_1: h_2^2h_3^{-2}h_2=1$, $R_2: (h_2h_3^{-1})^2h_3^{-1}h_2h_3=1$:\\
$ R_1\Rightarrow h_3^{-2}=h_2^{-3} \ (*)   \text{ and } R_2: h_2h_3^{-1}h_2h_3^{-2}h_2h_3=1 \ (**)$. By $(*)$ and $(**)$ we have $h_2h_3^{-1}h_2h_2^{-3}h_2h_3=1\Rightarrow h_2h_3^{-1}h_2^{-1}h_3=1 \Rightarrow h_2^{-1}h_3h_2=h_3\Rightarrow \langle h_2,h_3\rangle\cong BS(1,1)$ is solvable, a contradiction.

\item[(12)]$ R_1: h_2^2h_3^{-2}h_2=1$, $R_2: h_2^3h_3^2=1$:\\
$ R_1\Rightarrow h_2^3=h_3^2 \text{ and } R_2\Rightarrow h_2^3=h_3^{-2}. \ \Rightarrow h_3^4=1\Rightarrow h_3=1$, a contradiction.
\begin{table}[h]
\centering
\caption{The relations of a $C_4-C_5$ in the Kaplansky graph}\label{tab-C4-C5}
\begin{tabular}{|c|l|l||c|l|l|}\hline
$n$&$R_1$&$R_2$&$n$&$R_1$&$R_2$\\\hline
$1$&$h_2^2h_3h_2^{-1}h_3=1$&$h_2^2(h_2h_3^{-1})^2h_2=1$&$29$&$h_2h_3^{-3}h_2=1$&$h_2h_3h_2^{-1}(h_2^{-1}h_3)^2=1$\\
$2$&$h_2^2h_3h_2^{-1}h_3=1$&$h_2^2h_3^{-1}h_2^{-2}h_3=1$&$30$&$h_2h_3^{-3}h_2=1$&$h_2(h_3h_2^{-1})^2h_2^{-1}h_3=1$\\
$3$&$h_2^2h_3h_2^{-1}h_3=1$&$h_2(h_3^{-1}h_2h_3^{-1})^2h_2=1$&$31$&$h_2h_3^{-3}h_2=1$&$(h_2h_3^{-1})^2(h_2^{-1}h_3)^2=1$\\
$4$&$h_2^2h_3h_2^{-1}h_3=1$&$h_2h_3^{-2}h_2^{-1}h_3^2=1$&$32$&$h_2h_3^{-3}h_2=1$&$h_2h_3^{-2}h_2^{-1}h_3^2=1$\\
$5$&$h_2^2h_3h_2^{-1}h_3=1$&$h_2h_3^{-1}h_2h_3h_2^{-1}h_3^2=1$&$33$&$h_2h_3^{-2}h_2h_3=1$&$h_2^3h_3^{-2}h_2=1$\\
$6$&$h_2^2h_3h_2^{-1}h_3=1$&$h_2h_3^{-1}h_2h_3^2h_2^{-1}h_3=1$&$34$&$h_2h_3^{-2}h_2h_3=1$&$h_2^2h_3^{-1}h_2^{-2}h_3=1$\\
$7*$&$h_2^2h_3h_2^{-1}h_3=1$&$(h_2h_3^{-1})^2(h_2^{-1}h_3)^2=1$&$35$&$h_2h_3^{-2}h_2h_3=1$&$h_2h_3h_2^{-1}(h_3^{-1}h_2)^2=1$\\
$8*$&$h_2^2h_3^{-2}h_2=1$&$h_2h_3^{-2}h_2^{-1}h_3^2=1$&$36$&$h_2h_3^{-2}h_2h_3=1$&$h_2h_3^{-2}h_2^{-1}h_3^2=1$\\
$9$&$h_2^2h_3^{-2}h_2=1$&$h_2h_3^{-1}(h_3^{-1}h_2)^2h_3=1$&$37$&$h_2h_3^{-2}h_2h_3=1$&$(h_2h_3^{-1})^3h_2h_3=1$\\
$10$&$h_2^2h_3^{-2}h_2=1$&$(h_2h_3^{-1})^2(h_2^{-1}h_3)^2=1$&$38$&$h_2h_3^{-1}h_2h_3^2=1$&$h_2^2h_3^{-1}h_2^{-2}h_3=1$\\
$11$&$h_2^2h_3^{-2}h_2=1$&$(h_2h_3^{-1})^2h_3^{-1}h_2h_3=1$&$39$&$h_2h_3^{-1}h_2h_3^2=1$&$h_2h_3h_2^{-1}h_3h_2h_3^{-1}h_2=1$\\
$12$&$h_2^2h_3^{-2}h_2=1$&$h_2^3h_3^2=1$&$40$&$h_2h_3^{-1}h_2h_3^2=1$&$h_2^2h_3^{-1}h_2h_3h_2^{-1}h_3=1$\\
$13$&$h_2^2h_3^{-2}h_2=1$&$h_2^2h_3^{-1}h_2^{-2}h_3=1$&$41*$&$h_2h_3^{-1}h_2h_3^2=1$&$(h_2h_3^{-1})^2(h_2^{-1}h_3)^2=1$\\
$14$&$h_2^2h_3^{-2}h_2=1$&$h_2h_3^2h_2^{-1}h_3^2=1$&$42$&$h_2h_3^{-1}h_2h_3^2=1$&$h_2h_3^{-2}h_2^{-1}h_3^2=1$\\
$15$&$h_2h_3h_2^{-2}h_3=1$&$h_2^2h_3^{-1}h_2^{-2}h_3=1$&$43$&$h_2h_3^{-1}h_2h_3^2=1$&$h_2(h_3^{-1}h_2h_3^{-1})^2h_2=1$\\
$16$&$h_2h_3h_2^{-2}h_3=1$&$h_2h_3^{-4}h_2=1$&$44$&$h_2h_3^{-1}h_2h_3^2=1$&$h_3^2(h_3h_2^{-1})^2h_3=1$\\
$17$&$h_2h_3h_2^{-2}h_3=1$&$h_2h_3^{-1}(h_2^{-1}h_3)^2h_3=1$&$45$&$h_2h_3^{-1}h_2h_3h_2^{-1}h_3=1$&$h_2^3h_3h_2^{-1}h_3=1$\\
$18$&$h_2h_3h_2^{-2}h_3=1$&$h_2h_3^{-2}h_2^{-1}h_3^2=1$&$46$&$h_2h_3^{-1}h_2h_3h_2^{-1}h_3=1$&$h_2^2h_3^{-1}h_2^{-2}h_3=1$\\
$19$&$h_2h_3h_2^{-2}h_3=1$&$h_2(h_3h_2^{-1})^3h_3=1$&$47$&$h_2h_3^{-1}h_2h_3h_2^{-1}h_3=1$&$h_2h_3^2h_2^{-1}h_3^{-1}h_2=1$\\
$20*$&$h_2(h_3h_2^{-1})^2h_3=1$&$h_2^2h_3^{-1}h_2^{-2}h_3=1$&$48$&$h_2h_3^{-1}h_2h_3h_2^{-1}h_3=1$&$h_2^2h_3^{-1}h_2^{-1}h_3^2=1$\\
$21$&$h_2(h_3h_2^{-1})^2h_3=1$&$h_2^2h_3^{-1}h_2^2h_3=1$&$49$&$h_2h_3^{-1}h_2h_3h_2^{-1}h_3=1$&$h_2h_3^{-2}h_2^{-1}h_3^2=1$\\
$22$&$h_2(h_3h_2^{-1})^2h_3=1$&$h_2(h_2h_3^{-1})^3h_2=1$&$50$&$h_2h_3^{-1}h_2h_3h_2^{-1}h_3=1$&$h_2h_3^{-1}h_2h_3^3=1$\\
$23$&$h_2(h_3h_2^{-1})^2h_3=1$&$h_2h_3^2h_2^{-2}h_3=1$&$51$&$(h_2h_3^{-1})^2h_2h_3=1$&$h_2^2h_3^{-1}h_2^{-2}h_3=1$\\
$24$&$h_2(h_3h_2^{-1})^2h_3=1$&$h_2h_3h_2^{-2}h_3^2=1$&$52$&$(h_2h_3^{-1})^2h_2h_3=1$&$h_2h_3h_2h_3^{-2}h_2=1$\\
$25$&$h_2(h_3h_2^{-1})^2h_3=1$&$h_2h_3^{-2}h_2^{-1}h_3^2=1$&$53$&$(h_2h_3^{-1})^2h_2h_3=1$&$h_2^2h_3^{-2}h_2h_3=1$\\
$26$&$h_2h_3^{-3}h_2=1$&$h_2^2h_3^3=1$&$54$&$(h_2h_3^{-1})^2h_2h_3=1$&$h_2h_3^2h_2^{-1}h_3^2=1$\\
$27*$&$h_2h_3^{-3}h_2=1$&$h_2^2h_3^{-1}h_2^{-2}h_3=1$&$55*$&$(h_2h_3^{-1})^2h_2h_3=1$&$h_2h_3^{-2}h_2^{-1}h_3^2=1$\\
$28$&$h_2h_3^{-3}h_2=1$&$h_2^2h_3^{-1}h_2^2h_3=1$&$56$&$(h_2h_3^{-1})^2h_2h_3=1$&$h_3(h_3h_2^{-1})^3h_3=1$\\
\hline
\end{tabular}
\end{table}

\item[(13)]$ R_1: h_2^2h_3^{-2}h_2=1$, $R_2: h_2^2h_3^{-1}h_2^{-2}h_3=1$:\\
$ R_1\Rightarrow h_2^2h_3^{-1}=h_2^{-1}h_3 \ (*)  \Rightarrow R_2:h_2^{-1}h_3h_2^{-2}h_3=1\Rightarrow h_2=(h_2^{-1}h_3)^2\Rightarrow \langle h_2,h_3\rangle=\langle h_2,h_2^{-1}h_3\rangle=\langle h_2^{-1}h_3\rangle  \text{ is abelian}$, a contradiction.

\item[(14)]$ R_1: h_2^2h_3^{-2}h_2=1$, $R_2: h_2h_3^2h_2^{-1}h_3^2=1$:\\
$ R_1\Rightarrow h_2^3=h_3^2 \ (*)  \Rightarrow R_2:h_2^6=1\Rightarrow h_2=1$, a contradiction.

\item[(15)]$ R_1: h_2h_3h_2^{-2}h_3=1$, $R_2: h_2^2h_3^{-1}h_2^{-2}h_3=1$:\\
$ R_2\Rightarrow h_2^{-2}h_3=h_3h_2^{-2} \ (*)  \Rightarrow R_1:h_2=h_3^2\Rightarrow \langle h_2,h_3\rangle=\langle h_3\rangle \text{ is abelian, a contradiction}$.

\item[(16)]$ R_1: h_2h_3h_2^{-2}h_3=1$, $R_2: h_2h_3^{-4}h_2=1$:\\
$ R_2\Rightarrow h_2^{-2}=h_3^{-4} \ (*)  \Rightarrow R_1:h_2h_3h_3^{-4}h_3=1\Rightarrow h_2=h_3^2\Rightarrow \langle h_2,h_3\rangle=\langle h_3\rangle  \text{ is abelian}$, a contradiction.

\item[(17)]$ R_1: h_2h_3h_2^{-2}h_3=1$, $R_2: h_2h_3^{-1}(h_2^{-1}h_3)^2h_3=1$:\\
$ R_1\Rightarrow h_2h_3^{-1}h_2^{-1}=h_2^{-1}h_3 \ (*)   \text{ and } R_2: h_2h_3^{-1}h_2^{-1}h_3h_2^{-1}h_3^2=1 \ (**)$. By $(*)$ and $(**)$ we have $(h_2^{-1}h_3^2)^2=1\text{ and } G \text{ is torsion-free}\Rightarrow h_2^{-1}h_3^2=1\Rightarrow h_2=h_3^2\Rightarrow \langle h_2,h_3\rangle=\langle h_3\rangle \text{ is abelian}$, a contradiction.

\item[(18)]$ R_1: h_2h_3h_2^{-2}h_3=1$, $R_2: h_2h_3^{-2}h_2^{-1}h_3^2=1$:\\
$R_2\Rightarrow h_2h_3^2=h_3^2h_2 \ (*)   \text{ and } R_1\Rightarrow (h_3h_2h_3h_2^{-2})(h_2h_3h_2^{-2}h_3)=1 \ (**)$. By $(*)$ and $(**)$ we have $h_2^{-1}h_3^3h_2^{-1}h_3=1\Rightarrow h_3^2=(h_3^{-1}h_2)^2 \Rightarrow \langle h_2,h_3\rangle=\langle h_3,h_3^{-1}h_2\rangle\cong BS(1,-1)$ is solvable, a contradiction.

\item[(19)]$ R_1: h_2h_3h_2^{-2}h_3=1$, $R_2: h_2(h_3h_2^{-1})^3h_3=1$:\\
$ R_1\Rightarrow h_3h_2h_3h_2^{-1}=h_2 \ (*)   \text{ and } R_2: h_2h_3h_2^{-1}h_3h_2^{-1}h_3h_2^{-1}h_3=1 \ (**)$. By $(*)$ and $(**)$ we have $h_2=(h_2h_3^{-1})^2\Rightarrow \langle h_2,h_3\rangle=\langle h_2,h_2h_3^{-1}\rangle=\langle h_2h_3^{-1}\rangle  \text{ is abelian}$, a contradiction.


\item[(21)]$ R_1: h_2(h_3h_2^{-1})^2h_3=1$, $R_2: h_2^2h_3^{-1}h_2^2h_3=1$:\\
Using Tietze transformation where $h_3\mapsto h_3h_2$ and $h_2\mapsto h_2$, we have $ R_1: h_2^2h_3^3=1  \text{ and } R_2: h_2^2h_3^{-1}h_2^2h_3=1$. $ R_1\Rightarrow h_2^2=h_3^{-3} \ (*)  \Rightarrow R_2:h_3^6=1\Rightarrow h_3=1$, a contradiction.

\item[(22)]$ R_1: h_2(h_3h_2^{-1})^2h_3=1$, $R_2: h_2(h_2h_3^{-1})^3h_2=1$:\\
$ R_1\Rightarrow (h_3^{-1}h_2)^2h_3^{-1}=h_2 \ (*)   \text{ and } R_2:h_2^3(h_3^{-1}h_2)^2h_3^{-1}=1\ (**)$. By $(*)$ and $(**)$ we have $h_2^4=1\Rightarrow h_2=1$, a contradiction.

\item[(23)]$ R_1: h_2(h_3h_2^{-1})^2h_3=1$, $R_2: h_2h_3^2h_2^{-2}h_3=1$:\\
$ R_1\Rightarrow h_2^{-1}h_3h_2h_3=h_3^{-1}h_2 \ (*)  \Rightarrow R_2:h_2^{-1}h_3^{-1}h_2h_3=1 \Rightarrow h_2^{-1}h_3h_2=h_3\Rightarrow \langle h_2,h_3\rangle\cong BS(1,1)$ is solvable, a contradiction.

\item[(24)]$ R_1: h_2(h_3h_2^{-1})^2h_3=1$, $R_2: h_2h_3h_2^{-2}h_3^2=1$:\\
$ R_1\Rightarrow h_3h_2h_3h_2^{-1}=h_2h_3^{-1} \ (*)  \Rightarrow R_2:h_2^{-1}h_3h_2h_3^{-1}=1 \Rightarrow h_2^{-1}h_3h_2=h_3\Rightarrow \langle h_2,h_3\rangle\cong BS(1,1)$ is solvable, a contradiction.

\item[(25)]$ R_1: h_2(h_3h_2^{-1})^2h_3=1$, $R_2: h_2h_3^{-2}h_2^{-1}h_3^2=1$:\\
Using Tietze transformation where $h_2\mapsto h_2h_3$ and $h_2\mapsto h_2$, we have $ R_1: h_2h_3h_2^{-2}h_3=1$ and $R_2: h_2h_3^{-2}h_2^{-1}h_3^2=1$. With the same discussion such as (18), there is a contradiction.

\item[(26)]$ R_1: h_2h_3^{-3}h_2=1$, $R_2: h_2^2h_3^3=1$:\\
$ R_1\Rightarrow h_2^2=h_3^3 \text{ and } R_2\Rightarrow h_2^2=h_3^{-3}. \ \Rightarrow h_3^6=1\Rightarrow h_3=1$, a contradiction.


\item[(28)]$ R_1: h_2h_3^{-3}h_2=1$, $R_2: h_2^2h_3^{-1}h_2^2h_3=1$:\\
$ R_1\Rightarrow h_2^2=h_3^3  \ (*)  \Rightarrow R_2:h_3^6=1\Rightarrow h_3=1$, a contradiction.

\item[(29)]$ R_1: h_2h_3^{-3}h_2=1$, $R_2: h_2h_3h_2^{-1}(h_2^{-1}h_3)^2=1$:\\
By interchanging $h_2$ and $h_3$ in (9) and with the same discussion, there is a contradiction.

\item[(30)]$ R_1: h_2h_3^{-3}h_2=1$, $R_2: h_2(h_3h_2^{-1})^2h_2^{-1}h_3=1$:\\
By interchanging $h_2$ and $h_3$ in (11) and with the same discussion, there is a contradiction.

\item[(31)]$ R_1: h_2h_3^{-3}h_2=1$, $R_2: (h_2h_3^{-1})^2(h_2^{-1}h_3)^2=1$:\\
By interchanging $h_2$ and $h_3$ in (10) and with the same discussion, there is a contradiction.

\item[(32)]$ R_1: h_2h_3^{-3}h_2=1$, $R_2: h_2h_3^{-2}h_2^{-1}h_3^2=1$:\\
By interchanging $h_2$ and $h_3$ in (13) and with the same discussion, there is a contradiction.

\item[(33)]$ R_1: h_2h_3^{-2}h_2h_3=1$, $R_2: h_2^3h_3^{-2}h_2=1$:\\
By interchanging $h_2$ and $h_3$ in (16) and with the same discussion, there is a contradiction.

\item[(34)]$ R_1: h_2h_3^{-2}h_2h_3=1$, $R_2: h_2^2h_3^{-1}h_2^{-2}h_3=1$:\\
By interchanging $h_2$ and $h_3$ in (18) and with the same discussion, there is a contradiction.

\item[(35)]$ R_1: h_2h_3^{-2}h_2h_3=1$, $R_2: h_2h_3h_2^{-1}(h_3^{-1}h_2)^2=1$:\\
By interchanging $h_2$ and $h_3$ in (17) and with the same discussion, there is a contradiction.

\item[(36)]$ R_1: h_2h_3^{-2}h_2h_3=1$, $R_2: h_2h_3^{-2}h_2^{-1}h_3^2=1$:\\
By interchanging $h_2$ and $h_3$ in (15) and with the same discussion, there is a contradiction.

\item[(37)]$ R_1: h_2h_3^{-2}h_2h_3=1$, $R_2: (h_2h_3^{-1})^3h_2h_3=1$:\\
By interchanging $h_2$ and $h_3$ in (19) and with the same discussion, there is a contradiction.

\item[(38)]$ R_1: h_2h_3^{-1}h_2h_3^2=1$, $R_2: h_2^2h_3^{-1}h_2^{-2}h_3=1$:\\
By interchanging $h_2$ and $h_3$ in (4) and with the same discussion, there is a contradiction.

\item[(39)]$ R_1: h_2h_3^{-1}h_2h_3^2=1$, $R_2: h_2h_3h_2^{-1}h_3h_2h_3^{-1}h_2=1$:\\
By interchanging $h_2$ and $h_3$ in (5) and with the same discussion, there is a contradiction.

\item[(40)]$ R_1: h_2h_3^{-1}h_2h_3^2=1$, $R_2: h_2^2h_3^{-1}h_2h_3h_2^{-1}h_3=1$:\\
By interchanging $h_2$ and $h_3$ in (6) and with the same discussion, there is a contradiction.


\item[(42)]$ R_1: h_2h_3^{-1}h_2h_3^2=1$, $R_2: h_2h_3^{-2}h_2^{-1}h_3^2=1$:\\
By interchanging $h_2$ and $h_3$ in (2) and with the same discussion, there is a contradiction.

\item[(43)]$ R_1: h_2h_3^{-1}h_2h_3^2=1$, $R_2: h_2(h_3^{-1}h_2h_3^{-1})^2h_2=1$:\\
By interchanging $h_2$ and $h_3$ in (3) and with the same discussion, there is a contradiction.

\item[(44)]$ R_1: h_2h_3^{-1}h_2h_3^2=1$, $R_2: h_3^2(h_3h_2^{-1})^2h_3=1$:\\
By interchanging $h_2$ and $h_3$ in (1) and with the same discussion, there is a contradiction.

\item[(45)]$ R_1: h_2h_3^{-1}h_2h_3h_2^{-1}h_3=1$, $R_2: h_2^3h_3h_2^{-1}h_3=1$:\\
$ R_2\Rightarrow h_2h_3h_2^{-1}h_3=h_2^{-2} \ (*)  \Rightarrow R_1: h_3^{-1}h_2^{-1}=1\Rightarrow h_3=h_2^{-1}\Rightarrow \langle h_2,h_3\rangle=\langle h_2\rangle \text{ is abelian}$, a contradiction.

\item[(46)]$ R_1: h_2h_3^{-1}h_2h_3h_2^{-1}h_3=1$, $R_2: h_2^2h_3^{-1}h_2^{-2}h_3=1$:\\
$ R_2\Rightarrow h_2^{-1}h_3h_2=h_2h_3h_2^{-1} \ (*) \text{ and } h_3^{-1}h_2^{-2}=h_2^{-2}h_3^{-1} \ (**)$. $(*) \Rightarrow R_1:h_3^{-1}h_2h_3h_2h_3h_2^{-1}=1\Rightarrow (h_2^{-1}h_3^{-1}h_2h_3h_2h_3)(h_3^{-1}h_2h_3h_2h_3h_2^{-1})=1\ (***)$.  By $(**)$ and $(***)$ we have $h_2^2=h_3^{-2}\Rightarrow \langle h_2,h_3\rangle=\langle h_2,h_3^{-1}\rangle\cong BS(1,-1)$ is solvable, a contradiction.

\item[(47)]$ R_1: h_2h_3^{-1}h_2h_3h_2^{-1}h_3=1$, $R_2: h_2h_3^2h_2^{-1}h_3^{-1}h_2=1$:\\
$ R_1\Rightarrow h_3h_2^{-1}h_3^{-1}h_2=h_2h_3 \ (*)  \Rightarrow R_2:(h_2h_3)^2=1\text{ and } G \text{ is torsion-free}\Rightarrow h_3=h_2^{-1}\Rightarrow \langle h_2,h_3\rangle=\langle h_2\rangle \text{ is abelian}$, a contradiction.

\item[(48)]$ R_1: h_2h_3^{-1}h_2h_3h_2^{-1}h_3=1$, $R_2: h_2^2h_3^{-1}h_2^{-1}h_3^2=1$:\\
By interchanging $h_2$ and $h_3$ in (47) and with the same discussion, there is a contradiction.

\item[(49)]$ R_1: h_2h_3^{-1}h_2h_3h_2^{-1}h_3=1$, $R_2: h_2h_3^{-2}h_2^{-1}h_3^2=1$:\\
By interchanging $h_2$ and $h_3$ in (46) and with the same discussion, there is a contradiction.

\item[(50)]$ R_1: h_2h_3^{-1}h_2h_3h_2^{-1}h_3=1$, $R_2: h_2h_3^{-1}h_2h_3^3=1$:\\
By interchanging $h_2$ and $h_3$ in (45) and with the same discussion, there is a contradiction.

\item[(51)]$ R_1: (h_2h_3^{-1})^2h_2h_3=1$, $R_2: h_2^2h_3^{-1}h_2^{-2}h_3=1$:\\
By interchanging $h_2$ and $h_3$ in (25) and with the same discussion, there is a contradiction.

\item[(52)]$ R_1: (h_2h_3^{-1})^2h_2h_3=1$, $R_2: h_2h_3h_2h_3^{-2}h_2=1$:\\
By interchanging $h_2$ and $h_3$ in (24) and with the same discussion, there is a contradiction.

\item[(53)]$ R_1: (h_2h_3^{-1})^2h_2h_3=1$, $R_2: h_2^2h_3^{-2}h_2h_3=1$:\\
By interchanging $h_2$ and $h_3$ in (23) and with the same discussion, there is a contradiction.

\item[(54)]$ R_1: (h_2h_3^{-1})^2h_2h_3=1$, $R_2: h_2h_3^2h_2^{-1}h_3^2=1$:\\
By interchanging $h_2$ and $h_3$ in (21) and with the same discussion, there is a contradiction.


\item[(56)]$ R_1: (h_2h_3^{-1})^2h_2h_3=1$, $R_2: h_3(h_3h_2^{-1})^3h_3=1$:\\
By interchanging $h_2$ and $h_3$ in (22) and with the same discussion, there is a contradiction.
\end{enumerate} 

\subsection{$\mathbf{C_4-C_5(-C_5-)}$}
It can be seen that there are $42$ cases for the relations of a cycle $C_4$ and two cycles $C_5$ in the graph $C_4-C_5(-C_5-)$. Using GAP \cite{gap}, we see that all groups with two generators $h_2$ and $h_3$ and three relations which are between $38$ cases of these $42$ cases are finite and solvable, that is a contradiction. So, there are just $4$ cases for the relations of these cycles which may lead to the existence of a subgraph isomorphic to the graph $C_4-C_5(-C_5-)$ in $K(\alpha,\beta)$. Such cases are listed in table \ref{tab-C4-C5(-C5-)}. In the following it can be seen that all of these $4$ cases lead to  contradictions and so, the graph $K(\alpha,\beta)$ contains no subgraph isomorphic to the graph $C_4-C_5(-C_5-)$.
\begin{enumerate}
\item[(1)]$ R_1:h_2^2h_3h_2^{-1}h_3=1$, $R_2:(h_2h_3^{-1})^2(h_2^{-1}h_3)^2=1$, $R_3:h_2h_3^{-2}h_2^{-1}h_3^2=1$:\\
$R_3: h_2h_3^{-2}h_2^{-1}h_3^2=1\Rightarrow  h_3^2h_2^2=h_2^2h_3^2  \ (*) \text{ and } R_1: h_2^2h_3h_2^{-1}h_3=1\Rightarrow (h_3h_2^2h_3h_2^{-1})(h_2^2h_3h_2^{-1}h_3)=1  \ (**).$ By $(*)$ and $(**)$ we have $h_2h_3^3h_2h_3=1\Rightarrow h_3^2(h_2h_3)^2=1\Rightarrow h_3^2=(h_3^{-1}h_2^{-1})^2 \Rightarrow \langle h_2,h_3\rangle=\langle h_3,h_3^{-1}h_2^{-1}\rangle\cong BS(1,-1)$ is solvable, a contradiction.

\item[(2)]$ R_1:h_2^2h_3^{-2}h_2=1$, $R_2:h_2h_3^{-2}h_2^{-1}h_3^2=1$, $R_3:h_2h_3^2h_2^{-1}h_3^2=1$:\\
$ R_1\Rightarrow h_2^3=h_3^2 \ (*)  \Rightarrow R_3:h_2^6=1\Rightarrow h_2=1$, a contradiction.

\item[(3)]$ R_1:h_2h_3^{-3}h_2=1$, $R_2:h_2^2h_3^{-1}h_2^{-2}h_3=1$, $R_3:h_2^2h_3^{-1}h_2^2h_3=1$:\\
By interchanging $h_2$ and $h_3$ in (2) and with the same discussion, there is a contradiction.

\item[(4)]$ R_1:h_2h_3^{-1}h_2h_3^2=1$, $R_2:(h_2h_3^{-1})^2(h_2^{-1}h_3)^2=1$, $R_3:h_2^2h_3^{-1}h_2^{-2}h_3=1$:\\
By interchanging $h_2$ and $h_3$ in (1) and with the same discussion, there is a contradiction.
\end{enumerate}

\begin{table}[h]
\centering
\caption{All the relations related to the existence of a $C_4-C_5(-C_5-)$ in $K(\alpha,\beta)$}\label{tab-C4-C5(-C5-)}
\begin{tabular}{|c|l|l|l|}\hline
$n$&$R_1$&$R_2$&$R_3$\\\hline
$1$&$h_2^2h_3h_2^{-1}h_3=1$&$(h_2h_3^{-1})^2(h_2^{-1}h_3)^2=1$&$h_2h_3^{-2}h_2^{-1}h_3^2=1$\\
$2$&$h_2^2h_3^{-2}h_2=1$&$h_2h_3^{-2}h_2^{-1}h_3^2=1$&$h_2h_3^2h_2^{-1}h_3^2=1$\\
$3$&$h_2h_3^{-3}h_2=1$&$h_2^2h_3^{-1}h_2^{-2}h_3=1$&$h_2^2h_3^{-1}h_2^2h_3=1$\\
$4$&$h_2h_3^{-1}h_2h_3^2=1$&$(h_2h_3^{-1})^2(h_2^{-1}h_3)^2=1$&$h_2^2h_3^{-1}h_2^{-2}h_3=1$\\
\hline
\end{tabular}
\end{table}

\subsection{$\mathbf{C_4-C_5(-C_4-)}$}
It can be seen that there are $4$ cases for the relations of two cycles $C_4$ and a cycle $C_5$ in this structure. By considering all groups with two generators $h_2$ and $h_3$ and three relations which are between these cases and by using GAP \cite{gap}, we see that all of these groups are finite and solvable. So, the graph $K(\alpha,\beta)$ contains no subgraph isomorphic to the graph $C_4-C_5(-C_4-)$.

\subsection{$\mathbf{C_4-C_5(-C_6--)}$}
It can be seen that there are $126$ cases for the relations of a cycle $C_4$, a cycle $C_5$ and a cycle $C_6$ in the graph $C_4-C_5(-C_6--)$. Using GAP \cite{gap}, we see that all groups with two generators $h_2$ and $h_3$ and three relations which are between $122$ cases of these $126$ cases are finite and solvable, that is a contradiction with the assumptions. So, there are just $4$ cases for the relations of these cycles which may lead to the existence of a subgraph isomorphic to the graph $C_4-C_5(-C_6--)$ in $K(\alpha,\beta)$. These cases are listed in table \ref{tab-C4-C5(-C_6--)}.  In the following, it can be seen that all of these $4$ cases lead to contradictions and so, the graph $K(\alpha,\beta)$ contains no subgraph isomorphic to the graph $C_4-C_5(-C_6--)$.

\begin{enumerate}
\item[(1)]$ R_1:h_2^2h_3^{-2}h_2=1$, $R_2:h_2h_3^{-2}h_2^{-1}h_3^2=1$, $R_3:h_2^3h_3^{-1}h_2^{-1}h_3^{-1}h_2=1$:\\
$R_1\Rightarrow h_2^3h_3^{-1}=h_3\Rightarrow R_3:h_3h_2^{-1}h_3^{-1}h_2=1\Rightarrow h_2^{-1}h_3h_2=h_3\Rightarrow \langle h_2,h_3\rangle\cong BS(1,1)$ is solvable, a contradiction.

\item[(2)]$ R_1:h_2(h_3h_2^{-1})^2h_3=1$, $R_2:h_2^2h_3^{-1}h_2^{-2}h_3=1$, $R_3:h_2h_3^{-1}h_2(h_3h_2^{-1})^3h_3=1$:\\
$ R_1\Rightarrow (h_3h_2^{-1})^2h_3h_2=1\Rightarrow R_3:h_3^{-1}h_2h_3h_2^{-1}=1\Rightarrow h_2^{-1}h_3h_2=h_3\Rightarrow \langle h_2,h_3\rangle\cong BS(1,1)$ is solvable, a contradiction.

\item[(3)]$ R_1:h_2h_3^{-3}h_2=1$, $R_2:h_2^2h_3^{-1}h_2^{-2}h_3=1$, $R_3:h_2h_3^{-4}h_2h_3=1$:\\
By interchanging $h_2$ and $h_3$ in (1) and with the same discussion, there is a contradiction.

\item[(4)]$ R_1:(h_2h_3^{-1})^2h_2h_3=1$, $R_2:h_2h_3^{-2}h_2^{-1}h_3^2=1$, $R_3:(h_2h_3^{-1})^3h_2h_3h_2^{-1}h_3=1$:\\
By interchanging $h_2$ and $h_3$ in (2) and with the same discussion, there is a contradiction.
\end{enumerate}

\begin{table}[h]
\centering
\caption{All the relations related to the existence of a $C_4-C_5(-C_6--)$ in $K(\alpha,\beta)$}\label{tab-C4-C5(-C_6--)}
\begin{tabular}{|c|l|l|l|}\hline
$n$&$R_1$&$R_2$&$R_3$\\\hline
$1$&$h_2^2h_3^{-2}h_2=1$&$h_2h_3^{-2}h_2^{-1}h_3^2=1$&$h_2^3h_3^{-1}h_2^{-1}h_3^{-1}h_2=1$\\
$2$&$h_2(h_3h_2^{-1})^2h_3=1$&$h_2^2h_3^{-1}h_2^{-2}h_3=1$&$h_2h_3^{-1}h_2(h_3h_2^{-1})^3h_3=1$\\
$3$&$h_2h_3^{-3}h_2=1$&$h_2^2h_3^{-1}h_2^{-2}h_3=1$&$h_2h_3^{-4}h_2h_3=1$\\
$4$&$(h_2h_3^{-1})^2h_2h_3=1$&$h_2h_3^{-2}h_2^{-1}h_3^2=1$&$(h_2h_3^{-1})^3h_2h_3h_2^{-1}h_3=1$\\
\hline
\end{tabular}
\end{table}
Before finding other forbidden subgraphs, we give a useful lemma.
\begin{lem}\label{n-sub1}
Suppose that $G$ is a group with the generating set $S$ and $H$ is a group with the generating set $T$. Also suppose that $H\le G$. If for every $g\in S\cup S^{-1}$ and $h\in T$, $h^g=g^{-1} h g \in T$, then $H \trianglelefteq G$.
\end{lem}
\begin{proof}
The proof is straightforward.
\end{proof}
\begin{cor}\label{n-sub2}
Suppose that $G=\langle x,y\rangle$ and $y^2\in Z(G)$ and $H=\langle x,x^y\rangle$. Then $H \trianglelefteq G$.
\end{cor}

\subsection{$\mathbf{C_4-C_5(-C_6-)}$}
It can be seen that there are $462$ cases for the relations of a cycle $C_4$, a cycle $C_5$ and a cycle $C_6$ in the graph $C_4-C_5(-C_6-)$. By considering all groups with two generators $h_2$ and $h_3$ and two relations which are between these cases and by using GAP \cite{gap}, we see that $436$ groups are finite and solvable, or just finite. So, there are just $22$ cases for the relations of these cycles which may lead to the existence of a subgraph isomorphic to the graph $C_4-C_5(-C_6-)$ in $K(\alpha,\beta)$.  These cases are listed in table \ref{tab-C4-C5(-C6-)}.  In the following, it can be seen that all of these $22$ cases lead to contradictions and so, the graph $K(\alpha,\beta)$ contains no subgraph isomorphic to the graph $C_4-C_5(-C_6-)$.

\begin{enumerate}
\item[(1)]$ R_1:h_2^2h_3h_2^{-1}h_3=1$, $R_2:(h_2h_3^{-1})^2(h_2^{-1}h_3)^2=1$, $R_3:h_2^2h_3^{-1}h_2^{-1}h_3^3=1$:\\
$R_1\Rightarrow h_3h_2^2=h_2h_3^{-1} \ (*)  \Rightarrow R_3:h_2h_3^{-2}h_2^{-1}h_3^2=1\Rightarrow h_3^2\in Z(G)$ where $G=\langle h_2,h_3\rangle$. Let $x=h_2h_3^{-1}$. So $h_3^{-1}h_2=x^{h_3}$. \\
$R_2:(h_2h_3^{-1})^2(h_2^{-1}h_3)^2=1\Rightarrow (h_2h_3^{-1})^2=(h_3^{-1}h_2)^2 \text{ so if } H=\langle x,x^{h_3}\rangle\Rightarrow H\cong BS(1,-1)$ is solvable. By Corollary \ref{n-sub2} $H \trianglelefteq G$ since $G=\langle x,h_3\rangle$. Since $\frac{G}{H}=\frac{\langle h_3\rangle H}{H}$ is a cyclic group, it is solvable. So $G$ is solvable, a contradiction.

\item[(2)]$ R_1:h_2^2h_3h_2^{-1}h_3=1$, $R_2:(h_2h_3^{-1})^2(h_2^{-1}h_3)^2=1$, $R_3:h_2^2h_3^{-1}h_2(h_3h_2^{-1})^2h_3=1$:\\
$R_1\Rightarrow h_3h_2^2=h_2h_3^{-1} \ (*)  \Rightarrow R_3:h_3^{-2}h_2h_3h_2^{-1}h_3=1\Rightarrow h_2^{-1}h_3h_2=h_3\Rightarrow \langle h_2,h_3\rangle\cong BS(1,1)$ is solvable, a contradiction.

\item[(3)]$ R_1:h_2^2h_3h_2^{-1}h_3=1$, $R_2:(h_2h_3^{-1})^2(h_2^{-1}h_3)^2=1$, $R_3:h_2h_3h_2h_3^{-1}h_2^{-1}h_3^2=1$:\\
Using Tietze transformation where $h_3\mapsto h_3h_2$ and $h_2\mapsto h_2$, we have $ R_1:h_2^3h_3^2=1\Rightarrow h_3h_2^2=h_3^{-1}h_2^{-1} \ (*)$, and $R_3:h_3h_2^2h_3h_2h_3^{-1}h_2^{-1}h_3h_2=1$. By $(*)$ and $R_3$ we have $(h_3^{-1}h_2^{-1}h_3h_2)^2=1 \text{ and } G \text{ is torsion-free}\Rightarrow h_3^{-1}h_2^{-1}h_3h_2=1\Rightarrow h_2^{-1}h_3h_2=h_3\Rightarrow \langle h_2,h_3\rangle\cong BS(1,1)$ is solvable, a contradiction.

\item[(4)]$ R_1:h_2^2h_3h_2^{-1}h_3=1$, $R_2:(h_2h_3^{-1})^2(h_2^{-1}h_3)^2=1$, $R_3:h_2(h_3h_2^{-1})^2h_3^{-2}h_2=1$:\\
$ R_1\Rightarrow h_2^2h_3h_2^{-1}h_3=1  \ (*) \text{ and } R_3:h_2^2h_3h_2^{-1}h_3h_2^{-1}h_3^{-2}=1\Rightarrow h_2=h_3^{-2}\Rightarrow \langle h_2,h_3\rangle=\langle h_3\rangle \text{ is abelian}$, a contradiction.

\item[(5)]$ R_1:h_2^2h_3^{-2}h_2=1$, $R_2:h_2h_3^{-2}h_2^{-1}h_3^2=1$, $R_3:h_2^2h_3^{-3}h_2h_3=1$:\\
$ R_1\Rightarrow h_2^2h_3^{-2}=h_2^{-1} \ (*)\Rightarrow R_3:h_2^{-1}h_3^{-1}h_2h_3=1\Rightarrow h_2^{-1}h_3h_2=h_3\Rightarrow \langle h_2,h_3\rangle\cong BS(1,1)$ is solvable, a contradiction.

\item[(6)]$ R_1:h_2^2h_3^{-2}h_2=1$, $R_2:h_2h_3^{-2}h_2^{-1}h_3^2=1$, $R_3:h_2(h_2h_3^{-1})^3h_2^{-1}h_3=1$:\\
$ R_1\Rightarrow h_3^2\in Z(G) \text{, where } G=\langle h_2,h_3\rangle, \text{ and } h_2^2h_3^{-1}=h_2^{-1}h_3\ (*), \text{ and } R_3:h_2^2h_3^{-1}h_2h_3^{-1}h_2h_3^{-1}h_2^{-1}h_3=1$. By $(*)$ and $R_3$ we have $(h_2^{-1}h_3)^2(h_2h_3^{-1})^2=1\Rightarrow (h_2h_3^{-1})^2=(h_3^{-1}h_2)^2\ (**)$. Let $x=h_2h_3^{-1}$. So $h_3^{-1}h_2=x^{h_3}$. If $H=\langle x,x^{h_3}\rangle$, by $(**)$ we have $H\cong BS(1,-1)$ is solvable. By Corollary \ref{n-sub2} $H \trianglelefteq G$ since $G=\langle x,h_3\rangle$. Since $\frac{G}{H}=\frac{\langle h_3\rangle H}{H}$ is a cyclic group, it is solvable. So $G$ is solvable, a contradiction.

\item[(7)]$ R_1:h_2^2h_3^{-2}h_2=1$, $R_2:h_2h_3^{-2}h_2^{-1}h_3^2=1$, $R_3:h_2(h_3h_2^{-1})^2h_3^{-2}h_2=1$:\\
$ R_1\Rightarrow h_3^{-2}h_2^2=h_2^{-1} \ (*) \text{ and } R_3:h_3h_2^{-1}h_3h_2^{-1}h_3^{-2}h_2^2=1$. By $(*)$ and $R_3$ we have $h_2=(h_3h_2^{-1})^2\Rightarrow \langle h_2,h_3\rangle=\langle h_2,h_3h_2^{-1}\rangle=\langle h_3h_2^{-1}\rangle \text{ is abelian}$, a contradiction.

\item[(8)]$ R_1:h_2^2h_3^{-2}h_2=1$, $R_2:h_2h_3^{-2}h_2^{-1}h_3^2=1$, $R_3:(h_2h_3^{-1}h_2)^2h_3^{-1}h_2^{-1}h_3=1$:\\
$ R_1\Rightarrow h_2^2h_3^{-1}=h_2^{-1}h_3 \ (*)$, and $R_3:h_2h_3^{-1}h_2^2h_3^{-1}h_2h_3^{-1}h_2^{-1}h_3=1$. By $(*)$ and $R_3$ we have $(h_2h_3^{-1}h_2^{-1}h_3)^2=1 \text{ and } G \text{ is torsion-free}\Rightarrow h_2h_3^{-1}h_2^{-1}h_3=1\Rightarrow h_2^{-1}h_3h_2=h_3\Rightarrow \langle h_2,h_3\rangle\cong BS(1,1)$ is solvable, a contradiction.

\begin{table}[ht]
\centering
\caption{The relations of a $C_4-C_5(-C_6-)$ in $K(\alpha,\beta)$}\label{tab-C4-C5(-C6-)}
\begin{tabular}{|c|l|l|l|}\hline
$n$&$R_1$&$R_2$&$R_3$\\\hline
$1$&$h_2^2h_3h_2^{-1}h_3=1$&$(h_2h_3^{-1})^2(h_2^{-1}h_3)^2=1$&$h_2^2h_3^{-1}h_2^{-1}h_3^3=1$\\
$2$&$h_2^2h_3h_2^{-1}h_3=1$&$(h_2h_3^{-1})^2(h_2^{-1}h_3)^2=1$&$h_2^2h_3^{-1}h_2(h_3h_2^{-1})^2h_3=1$\\
$3$&$h_2^2h_3h_2^{-1}h_3=1$&$(h_2h_3^{-1})^2(h_2^{-1}h_3)^2=1$&$h_2h_3h_2h_3^{-1}h_2^{-1}h_3^2=1$\\
$4$&$h_2^2h_3h_2^{-1}h_3=1$&$(h_2h_3^{-1})^2(h_2^{-1}h_3)^2=1$&$h_2(h_3h_2^{-1})^2h_3^{-2}h_2=1$\\
$5$&$h_2^2h_3^{-2}h_2=1$&$h_2h_3^{-2}h_2^{-1}h_3^2=1$&$h_2^2h_3^{-3}h_2h_3=1$\\
$6$&$h_2^2h_3^{-2}h_2=1$&$h_2h_3^{-2}h_2^{-1}h_3^2=1$&$h_2(h_2h_3^{-1})^3h_2^{-1}h_3=1$\\
$7$&$h_2^2h_3^{-2}h_2=1$&$h_2h_3^{-2}h_2^{-1}h_3^2=1$&$h_2(h_3h_2^{-1})^2h_3^{-2}h_2=1$\\
$8$&$h_2^2h_3^{-2}h_2=1$&$h_2h_3^{-2}h_2^{-1}h_3^2=1$&$(h_2h_3^{-1}h_2)^2h_3^{-1}h_2^{-1}h_3=1$\\
$9$&$h_2^2h_3^{-2}h_2=1$&$h_2h_3^{-2}h_2^{-1}h_3^2=1$&$h_2(h_2h_3^{-1})^4h_2=1$\\
$10$&$h_2(h_3h_2^{-1})^2h_3=1$&$h_2^2h_3^{-1}h_2^{-2}h_3=1$&$h_2h_3^{-1}(h_2^{-1}h_3)^2h_3^2=1$\\
$11$&$h_2(h_3h_2^{-1})^2h_3=1$&$h_2^2h_3^{-1}h_2^{-2}h_3=1$&$h_3^3(h_3h_2^{-1})^2h_3=1$\\
$12$&$h_2h_3^{-3}h_2=1$&$h_2^2h_3^{-1}h_2^{-2}h_3=1$&$(h_2h_3^{-1})^3h_3^{-1}h_2^{-1}h_3=1$\\
$13$&$h_2h_3^{-3}h_2=1$&$h_2^2h_3^{-1}h_2^{-2}h_3=1$&$h_3(h_3h_2^{-1})^4h_3=1$\\
$14$&$h_2h_3^{-3}h_2=1$&$h_2^2h_3^{-1}h_2^{-2}h_3=1$&$h_2^2h_3^{-2}h_2^{-1}h_3^{-1}h_2=1$\\
$15$&$h_2h_3^{-3}h_2=1$&$h_2^2h_3^{-1}h_2^{-2}h_3=1$&$h_2(h_3h_2^{-1})^2h_3^{-2}h_2=1$\\
$16$&$h_2h_3^{-3}h_2=1$&$h_2^2h_3^{-1}h_2^{-2}h_3=1$&$h_2(h_3^{-1}h_2h_3^{-1})^2h_2^{-1}h_3=1$\\
$17$&$h_2h_3^{-1}h_2h_3^2=1$&$(h_2h_3^{-1})^2(h_2^{-1}h_3)^2=1$&$h_2^2h_3^2h_2^{-1}h_3^{-1}h_2=1$\\
$18$&$h_2h_3^{-1}h_2h_3^2=1$&$(h_2h_3^{-1})^2(h_2^{-1}h_3)^2=1$&$(h_2h_3)^2h_2^{-1}h_3^{-1}h_2=1$\\
$19$&$h_2h_3^{-1}h_2h_3^2=1$&$(h_2h_3^{-1})^2(h_2^{-1}h_3)^2=1$&$h_2(h_3h_2^{-1})^2h_3^{-2}h_2=1$\\
$20$&$h_2h_3^{-1}h_2h_3^2=1$&$(h_2h_3^{-1})^2(h_2^{-1}h_3)^2=1$&$(h_2h_3^{-1})^2h_2h_3^2h_2^{-1}h_3=1$\\
$21$&$(h_2h_3^{-1})^2h_2h_3=1$&$h_2h_3^{-2}h_2^{-1}h_3^2=1$&$h_2^3(h_2h_3^{-1})^2h_2=1$\\
$22$&$(h_2h_3^{-1})^2h_2h_3=1$&$h_2h_3^{-2}h_2^{-1}h_3^2=1$&$h_2^2h_3h_2^{-1}(h_3^{-1}h_2)^2=1$\\
\hline
\end{tabular}
\end{table}
\item[(9)]$ R_1:h_2^2h_3^{-2}h_2=1$, $R_2:h_2h_3^{-2}h_2^{-1}h_3^2=1$, $R_3:h_2(h_2h_3^{-1})^4h_2=1$:\\
$ R_1\Rightarrow h_2^3h_3^{-1}=h_3 \ (*)\Rightarrow R_3:h_2h_3^{-1}h_2h_3^{-1}h_2=1\Rightarrow h_2=(h_3h_2^{-1})^2\Rightarrow \langle h_2,h_3\rangle=\langle h_2,h_3h_2^{-1}\rangle=\langle h_3h_2^{-1}\rangle \text{ is abelian}$, a contradiction.

\item[(10)]$ R_1:h_2(h_3h_2^{-1})^2h_3=1$, $R_2:h_2^2h_3^{-1}h_2^{-2}h_3=1$, $R_3:h_2h_3^{-1}(h_2^{-1}h_3)^2h_3^2=1$:\\
$ R_1\Rightarrow (h_2^{-1}h_3)^2=h_3^{-1}h_2^{-1}\ (*)$. By $(*)$ and $R_3$ we have  $h_2h_3^{-2}h_2^{-1}h_3^2=1\ (**)$. Using Tietze transformation where $h_2\mapsto h_2h_3$ and $h_3\mapsto h_3$, we have $R_2:(h_3^{-1}h_2^{-1})^2=(h_2^{-1}h_3^{-1})^2\ (***)$ and $(**): h_2h_3^{-2}h_2^{-1}h_3^2=1$. $(**)\Rightarrow h_3^2\in Z(G) \text{, where } G=\langle h_2,h_3\rangle$. Let $x=h_2^{-1}h_3^{-1}$. So $h_3^{-1}h_2^{-1}=x^{h_3}$. If $H=\langle x,x^{h_3}\rangle$, by $(***)$ we have $H\cong BS(1,-1)$ is solvable. By Corollary \ref{n-sub2} $H \trianglelefteq G$ since $G=\langle x,h_3\rangle$. Since $\frac{G}{H}=\frac{\langle h_3\rangle H}{H}$ is a cyclic group, it is solvable. So $G$ is solvable, a contradiction.

\item[(11)]$ R_1:h_2(h_3h_2^{-1})^2h_3=1$, $R_2:h_2^2h_3^{-1}h_2^{-2}h_3=1$, $R_3:h_3^3(h_3h_2^{-1})^2h_3=1$:\\
$ R_1\Rightarrow (h_3h_2^{-1})^2h_3=h_2^{-1} \ (*)\Rightarrow R_3:h_3^3h_2^{-1}=1\Rightarrow h_2=h_3^3\Rightarrow \langle h_2,h_3\rangle=\langle h_3\rangle \text{ is abelian}$, a contradiction.

\item[(12)]$ R_1:h_2h_3^{-3}h_2=1$, $R_2:h_2^2h_3^{-1}h_2^{-2}h_3=1$, $R_3:(h_2h_3^{-1})^3h_3^{-1}h_2^{-1}h_3=1$:\\
By interchanging $h_2$ and $h_3$ in (6) and with the same discussion, there is a contradiction.

\item[(13)]$ R_1:h_2h_3^{-3}h_2=1$, $R_2:h_2^2h_3^{-1}h_2^{-2}h_3=1$, $R_3:h_3(h_3h_2^{-1})^4h_3=1$:\\
By interchanging $h_2$ and $h_3$ in (9) and with the same discussion, there is a contradiction.

\item[(14)]$ R_1:h_2h_3^{-3}h_2=1$, $R_2:h_2^2h_3^{-1}h_2^{-2}h_3=1$, $R_3:h_2^2h_3^{-2}h_2^{-1}h_3^{-1}h_2=1$:\\
By interchanging $h_2$ and $h_3$ in (5) and with the same discussion, there is a contradiction.

\item[(15)]$ R_1:h_2h_3^{-3}h_2=1$, $R_2:h_2^2h_3^{-1}h_2^{-2}h_3=1$, $R_3:h_2(h_3h_2^{-1})^2h_3^{-2}h_2=1$:\\
By interchanging $h_2$ and $h_3$ in (7) and with the same discussion, there is a contradiction.

\item[(16)]$ R_1:h_2h_3^{-3}h_2=1$, $R_2:h_2^2h_3^{-1}h_2^{-2}h_3=1$, $R_3:h_2(h_3^{-1}h_2h_3^{-1})^2h_2^{-1}h_3=1$:\\
By interchanging $h_2$ and $h_3$ in (8) and with the same discussion, there is a contradiction.

\item[(17)]$ R_1:h_2h_3^{-1}h_2h_3^2=1$, $R_2:(h_2h_3^{-1})^2(h_2^{-1}h_3)^2=1$, $R_3:h_2^2h_3^2h_2^{-1}h_3^{-1}h_2=1$:\\
By interchanging $h_2$ and $h_3$ in (1) and with the same discussion, there is a contradiction.

\item[(18)]$ R_1:h_2h_3^{-1}h_2h_3^2=1$, $R_2:(h_2h_3^{-1})^2(h_2^{-1}h_3)^2=1$, $R_3:(h_2h_3)^2h_2^{-1}h_3^{-1}h_2=1$:\\
By interchanging $h_2$ and $h_3$ in (3) and with the same discussion, there is a contradiction.

\item[(19)]$ R_1:h_2h_3^{-1}h_2h_3^2=1$, $R_2:(h_2h_3^{-1})^2(h_2^{-1}h_3)^2=1$, $R_3:h_2(h_3h_2^{-1})^2h_3^{-2}h_2=1$:\\
By interchanging $h_2$ and $h_3$ in (4) and with the same discussion, there is a contradiction.

\item[(20)]$ R_1:h_2h_3^{-1}h_2h_3^2=1$, $R_2:(h_2h_3^{-1})^2(h_2^{-1}h_3)^2=1$, $R_3:(h_2h_3^{-1})^2h_2h_3^2h_2^{-1}h_3=1$:\\
By interchanging $h_2$ and $h_3$ in (2) and with the same discussion, there is a contradiction.

\item[(21)]$ R_1:(h_2h_3^{-1})^2h_2h_3=1$, $R_2:h_2h_3^{-2}h_2^{-1}h_3^2=1$, $R_3:h_2^3(h_2h_3^{-1})^2h_2=1$:\\
By interchanging $h_2$ and $h_3$ in (11) and with the same discussion, there is a contradiction.

\item[(22)]$ R_1:(h_2h_3^{-1})^2h_2h_3=1$, $R_2:h_2h_3^{-2}h_2^{-1}h_3^2=1$, $R_3:h_2^2h_3h_2^{-1}(h_3^{-1}h_2)^2=1$:\\
By interchanging $h_2$ and $h_3$ in (10) and with the same discussion, there is a contradiction.
\end{enumerate}

\subsection{$\mathbf{C_4-C_5(-C_7--)}$}
With the same discussion such as about $C_4$ cycles and considering the relations of $C_7$ cycles equivalent, there are $1173$ non-equivalent cases for the relations of a $C_7$ cycle in the graph $K(\alpha,\beta)$. By considering these relations, it can be seen that there are $648$ cases for the relations of a cycle $C_4$, a cycle $C_5$ and a cycle $C_7$ in the graph $C_4-C_5(-C_7--)$. Using GAP \cite{gap}, we see that all groups with two generators $h_2$ and $h_3$ and three relations which are between $608$ cases of these $648$ cases are finite and solvable, that is a contradiction with the assumptions. So, there are just $40$ cases for the relations of these cycles which may lead to the existence of a subgraph isomorphic to the graph $C_4-C_5(-C_7--)$ in $K(\alpha,\beta)$. These cases are listed in table \ref{tab-C4-C5(-C7--)}.  In the following, we show that all of these $40$ cases lead to contradictions and so, the graph $K(\alpha,\beta)$ contains no subgraph isomorphic to the graph $C_4-C_5(-C_7--)$.

\begin{enumerate}
\item[(1)]$ R_1:h_2^2h_3h_2^{-1}h_3=1$, $R_2:(h_2h_3^{-1})^2(h_2^{-1}h_3)^2=1$, $R_3:h_2^3h_3^2h_2^{-1}h_3^2=1$:\\
$\Rightarrow\langle h_2,h_3\rangle\cong BS(1,1)$ is solvable, a contradiction.

\item[(2)]$ R_1:h_2^2h_3h_2^{-1}h_3=1$, $R_2:(h_2h_3^{-1})^2(h_2^{-1}h_3)^2=1$, $R_3:h_2^2(h_2h_3^{-1})^4h_2=1$:\\
$\Rightarrow h_3=1$ that is a contradiction.

\item[(3)]$ R_1:h_2^2h_3h_2^{-1}h_3=1$, $R_2:(h_2h_3^{-1})^2(h_2^{-1}h_3)^2=1$, $R_3:h_2h_3^{-2}h_2h_3^2h_2^{-1}h_3^2=1$:\\
$\Rightarrow\langle h_2,h_3\rangle\cong BS(1,1)$ is solvable, a contradiction.

\item[(4)]$ R_1:h_2^2h_3h_2^{-1}h_3=1$, $R_2:(h_2h_3^{-1})^2(h_2^{-1}h_3)^2=1$, $R_3:(h_2h_3^{-1})^4h_3^{-1}h_2h_3^{-1}h_2=1$:\\
$\Rightarrow h_3=1$ that is a contradiction.

\item[(5)]$ R_1:h_2^2h_3^{-2}h_2=1$, $R_2:h_2h_3^{-2}h_2^{-1}h_3^2=1$, $R_3:h_2^3h_3^4=1$:\\
$\Rightarrow h_3=1$ that is a contradiction.

\item[(6)]$ R_1:h_2^2h_3^{-2}h_2=1$, $R_2:h_2h_3^{-2}h_2^{-1}h_3^2=1$, $R_3:h_2^3(h_3^{-1}h_2h_3^{-1})^2h_2=1$:\\
$\Rightarrow\langle h_2,h_3\rangle\cong BS(1,1)$ is solvable, a contradiction.

\item[(7)]$ R_1:h_2^2h_3^{-2}h_2=1$, $R_2:h_2h_3^{-2}h_2^{-1}h_3^2=1$, $R_3:h_2h_3^2h_2^{-1}h_3^4=1$:\\
$\Rightarrow h_3=1$ that is a contradiction.

\item[(8)]$ R_1:h_2^2h_3^{-2}h_2=1$, $R_2:h_2h_3^{-2}h_2^{-1}h_3^2=1$, $R_3:h_2(h_3^{-1}h_2h_3^{-1})^2h_2h_3h_2^{-1}h_3=1$:\\
$\Rightarrow\langle h_2,h_3\rangle\cong BS(1,1)$ is solvable, a contradiction.

\item[(9)]$ R_1:h_2^2h_3^{-2}h_2=1$, $R_2:h_2h_3^{-2}h_2^{-1}h_3^2=1$, $R_3:h_2^3h_3^{-1}(h_3^{-1}h_2)^3=1$:\\
$\Rightarrow \langle h_2,h_3\rangle=\langle h_2,h_3h_2^{-1}\rangle=\langle h_3h_2^{-1}\rangle$ is abelian, a contradiction.

\item[(10)]$ R_1:h_2^2h_3^{-2}h_2=1$, $R_2:h_2h_3^{-2}h_2^{-1}h_3^2=1$, $R_3:h_2^2(h_3^{-1}h_2^{-1})^2h_3^2=1$:\\
$\Rightarrow\langle h_2,h_3\rangle\cong BS(1,1)$ is solvable, a contradiction.

\begin{center}
\footnotesize
\begin{longtable}{|c|l|l|l|}
\caption{{\tiny The relations of a $C_4-C_5(-C_7--)$ in $K(\alpha,\beta)$}}\label{tab-C4-C5(-C7--)}\\
\hline \multicolumn{1}{|c|}{$n$} & \multicolumn{1}{l|}{$R_1$} & \multicolumn{1}{l|}{$R_2$} & \multicolumn{1}{l|}{$R_3$} \\\hline
\endfirsthead
\multicolumn{4}{c}%
{{\tablename\ \thetable{} -- continued from previous page}} \\
\hline  \multicolumn{1}{|c|}{$n$} & \multicolumn{1}{l|}{$R_1$} & \multicolumn{1}{l|}{$R_2$} & \multicolumn{1}{l|}{$R_3$}\\\hline
\endhead
\hline \multicolumn{4}{|r|}{{Continued on next page}} \\\hline
\endfoot
\hline
\endlastfoot
$1$&$h_2^2h_3h_2^{-1}h_3=1$&$(h_2h_3^{-1})^2(h_2^{-1}h_3)^2=1$&$h_2^3h_3^2h_2^{-1}h_3^2=1$\\
$2$&$h_2^2h_3h_2^{-1}h_3=1$&$(h_2h_3^{-1})^2(h_2^{-1}h_3)^2=1$&$h_2^2(h_2h_3^{-1})^4h_2=1$\\
$3$&$h_2^2h_3h_2^{-1}h_3=1$&$(h_2h_3^{-1})^2(h_2^{-1}h_3)^2=1$&$h_2h_3^{-2}h_2h_3^2h_2^{-1}h_3^2=1$\\
$4$&$h_2^2h_3h_2^{-1}h_3=1$&$(h_2h_3^{-1})^2(h_2^{-1}h_3)^2=1$&$(h_2h_3^{-1})^4h_3^{-1}h_2h_3^{-1}h_2=1$\\
$5$&$h_2^2h_3^{-2}h_2=1$&$h_2h_3^{-2}h_2^{-1}h_3^2=1$&$h_2^3h_3^4=1$\\
$6$&$h_2^2h_3^{-2}h_2=1$&$h_2h_3^{-2}h_2^{-1}h_3^2=1$&$h_2^3(h_3^{-1}h_2h_3^{-1})^2h_2=1$\\
$7$&$h_2^2h_3^{-2}h_2=1$&$h_2h_3^{-2}h_2^{-1}h_3^2=1$&$h_2h_3^2h_2^{-1}h_3^4=1$\\
$8$&$h_2^2h_3^{-2}h_2=1$&$h_2h_3^{-2}h_2^{-1}h_3^2=1$&$h_2(h_3^{-1}h_2h_3^{-1})^2h_2h_3h_2^{-1}h_3=1$\\
$9$&$h_2^2h_3^{-2}h_2=1$&$h_2h_3^{-2}h_2^{-1}h_3^2=1$&$h_2^3h_3^{-1}(h_3^{-1}h_2)^3=1$\\
$10$&$h_2^2h_3^{-2}h_2=1$&$h_2h_3^{-2}h_2^{-1}h_3^2=1$&$h_2^2(h_3^{-1}h_2^{-1})^2h_3^2=1$\\
$11$&$h_2^2h_3^{-2}h_2=1$&$h_2h_3^{-2}h_2^{-1}h_3^2=1$&$h_2^2h_3^{-1}h_2^{-1}h_3^{-1}(h_3^{-1}h_2)^2=1$\\
$12$&$h_2^2h_3^{-2}h_2=1$&$h_2h_3^{-2}h_2^{-1}h_3^2=1$&$h_2^2h_3^{-1}h_2^{-1}(h_3^{-1}h_2)^2h_3=1$\\
$13$&$h_2^2h_3^{-2}h_2=1$&$h_2h_3^{-2}h_2^{-1}h_3^2=1$&$h_2h_3h_2^{-2}h_3h_2^{-1}h_3^{-1}h_2^{-1}h_3=1$\\
$14$&$h_2^2h_3^{-2}h_2=1$&$h_2h_3^{-2}h_2^{-1}h_3^2=1$&$h_2h_3^{-1}(h_3^{-1}h_2)^3h_3h_2^{-1}h_3=1$\\
$15$&$h_2(h_3h_2^{-1})^2h_3=1$&$h_2^2h_3^{-1}h_2^{-2}h_3=1$&$h_2h_3^3h_2^{-1}h_3^{-1}h_2^{-1}h_3=1$\\
$16$&$h_2(h_3h_2^{-1})^2h_3=1$&$h_2^2h_3^{-1}h_2^{-2}h_3=1$&$h_2h_3^2(h_3h_2^{-1})^3h_3=1$\\
$17$&$h_2(h_3h_2^{-1})^2h_3=1$&$h_2^2h_3^{-1}h_2^{-2}h_3=1$&$h_2h_3(h_3h_2^{-1})^2h_3h_2h_3^{-1}h_2=1$\\
$18$&$h_2(h_3h_2^{-1})^2h_3=1$&$h_2^2h_3^{-1}h_2^{-2}h_3=1$&$h_2h_3h_2^{-2}h_3h_2^{-1}(h_3^{-1}h_2)^2=1$\\
$19$&$h_2(h_3h_2^{-1})^2h_3=1$&$h_2^2h_3^{-1}h_2^{-2}h_3=1$&$h_2h_3^{-1}h_2^{-1}h_3h_2^{-1}h_3^{-1}h_2h_3^{-2}h_2=1$\\
$20$&$h_2(h_3h_2^{-1})^2h_3=1$&$h_2^2h_3^{-1}h_2^{-2}h_3=1$&$h_2h_3^{-2}h_2^{-1}h_3h_2^{-1}(h_3^{-1}h_2)^2=1$\\
$21$&$h_2h_3^{-3}h_2=1$&$h_2^2h_3^{-1}h_2^{-2}h_3=1$&$h_2h_3h_2h_3^{-2}h_2^{-2}h_3=1$\\
$22$&$h_2h_3^{-3}h_2=1$&$h_2^2h_3^{-1}h_2^{-2}h_3=1$&$h_2h_3h_2h_3^{-2}(h_3^{-1}h_2)^2=1$\\
$23$&$h_2h_3^{-3}h_2=1$&$h_2^2h_3^{-1}h_2^{-2}h_3=1$&$h_2h_3^2h_2^{-1}h_3^{-1}(h_2^{-1}h_3)^2=1$\\
$24$&$h_2h_3^{-3}h_2=1$&$h_2^2h_3^{-1}h_2^{-2}h_3=1$&$h_2h_3h_2^{-1}h_3^2h_2^{-1}h_3^{-1}h_2^{-1}h_3=1$\\
$25$&$h_2h_3^{-3}h_2=1$&$h_2^2h_3^{-1}h_2^{-2}h_3=1$&$h_2h_3^{-1}h_2^{-1}h_3h_2^{-1}(h_3^{-1}h_2)^3=1$\\
$26$&$h_2h_3^{-3}h_2=1$&$h_2^2h_3^{-1}h_2^{-2}h_3=1$&$h_2h_3^{-3}(h_3^{-1}h_2)^3=1$\\
$27$&$h_2h_3^{-3}h_2=1$&$h_2^2h_3^{-1}h_2^{-2}h_3=1$&$h_2^4h_3^3=1$\\
$28$&$h_2h_3^{-3}h_2=1$&$h_2^2h_3^{-1}h_2^{-2}h_3=1$&$h_2^3h_3h_2^2h_3^{-1}h_2=1$\\
$29$&$h_2h_3^{-3}h_2=1$&$h_2^2h_3^{-1}h_2^{-2}h_3=1$&$h_2h_3^{-1}h_2h_3(h_2^{-1}h_3h_2^{-1})^2h_3=1$\\
$30$&$h_2h_3^{-3}h_2=1$&$h_2^2h_3^{-1}h_2^{-2}h_3=1$&$(h_2h_3^{-1})^2h_3^{-2}(h_3^{-1}h_2)^2=1$\\
$31$&$h_2h_3^{-1}h_2h_3^2=1$&$(h_2h_3^{-1})^2(h_2^{-1}h_3)^2=1$&$h_2^2h_3^{-1}h_2^2h_3^3=1$\\
$32$&$h_2h_3^{-1}h_2h_3^2=1$&$(h_2h_3^{-1})^2(h_2^{-1}h_3)^2=1$&$h_2^2h_3^{-1}h_2^2h_3h_2^{-2}h_3=1$\\
$33$&$h_2h_3^{-1}h_2h_3^2=1$&$(h_2h_3^{-1})^2(h_2^{-1}h_3)^2=1$&$(h_2h_3^{-1})^4h_3^{-1}h_2h_3^{-1}h_2=1$\\
$34$&$h_2h_3^{-1}h_2h_3^2=1$&$(h_2h_3^{-1})^2(h_2^{-1}h_3)^2=1$&$h_3^2(h_3h_2^{-1})^4h_3=1$\\
$35$&$(h_2h_3^{-1})^2h_2h_3=1$&$h_2h_3^{-2}h_2^{-1}h_3^2=1$&$h_2^3h_3^{-1}h_2^{-1}h_3^{-1}h_2h_3=1$\\
$36$&$(h_2h_3^{-1})^2h_2h_3=1$&$h_2h_3^{-2}h_2^{-1}h_3^2=1$&$h_2^2(h_2h_3^{-1})^3h_2h_3=1$\\
$37$&$(h_2h_3^{-1})^2h_2h_3=1$&$h_2h_3^{-2}h_2^{-1}h_3^2=1$&$h_2^2h_3^{-1}(h_3^{-1}h_2)^2h_3h_2^{-1}h_3=1$\\
$38$&$(h_2h_3^{-1})^2h_2h_3=1$&$h_2h_3^{-2}h_2^{-1}h_3^2=1$&$h_2(h_2h_3^{-1})^2h_2h_3h_2^{-1}h_3^2=1$\\
$39$&$(h_2h_3^{-1})^2h_2h_3=1$&$h_2h_3^{-2}h_2^{-1}h_3^2=1$&$h_2h_3^{-2}h_2^2h_3^{-1}h_2h_3h_2^{-1}h_3=1$\\
$40$&$(h_2h_3^{-1})^2h_2h_3=1$&$h_2h_3^{-2}h_2^{-1}h_3^2=1$&$h_2h_3^{-2}h_2h_3^{-1}(h_2^{-1}h_3)^2h_3=1$\\
\end{longtable}
\end{center}

\item[(11)]$ R_1:h_2^2h_3^{-2}h_2=1$, $R_2:h_2h_3^{-2}h_2^{-1}h_3^2=1$, $R_3:h_2^2h_3^{-1}h_2^{-1}h_3^{-1}(h_3^{-1}h_2)^2=1$:\\
$\Rightarrow h_2=1$ that is a contradiction.

\item[(12)]$ R_1:h_2^2h_3^{-2}h_2=1$, $R_2:h_2h_3^{-2}h_2^{-1}h_3^2=1$, $R_3:h_2^2h_3^{-1}h_2^{-1}(h_3^{-1}h_2)^2h_3=1$:\\
$R_2:h_2h_3^{-2}h_2^{-1}h_3^2=1\Rightarrow h_3^2\in Z(G)$ where $G=\langle h_2,h_3\rangle$. Let $x=h_2h_3^{-1}$. So $h_3^{-1}h_2=x^{h_3}$. \\
$R_1 \text{ and } R_3 \Rightarrow (h_2h_3^{-1})^2=(h_3^{-1}h_2)^2 \text{ so if } H=\langle x,x^{h_3}\rangle\Rightarrow H\cong BS(1,-1)$ is solvable. By Corollary \ref{n-sub2} $H \trianglelefteq G$ since $G=\langle x,h_3\rangle$. Since $\frac{G}{H}=\frac{\langle h_3\rangle H}{H}$ is a cyclic group, it is solvable. So $G$ is solvable, a contradiction.

\item[(13)]$ R_1:h_2^2h_3^{-2}h_2=1$, $R_2:h_2h_3^{-2}h_2^{-1}h_3^2=1$, $R_3:h_2h_3h_2^{-2}h_3h_2^{-1}h_3^{-1}h_2^{-1}h_3=1$:\\
$R_2:h_2h_3^{-2}h_2^{-1}h_3^2=1\Rightarrow h_3^2\in Z(G)$ where $G=\langle h_2,h_3\rangle$. Let $x=h_2^{-1}h_3^{-1}$. So $h_3^{-1}h_2^{-1}=x^{h_3}$. \\
$R_1 \text{ and } R_3 \Rightarrow (h_3^{-1}h_2^{-1})^2=(h_2^{-1}h_3^{-1})^2 \text{ so if } H=\langle x,x^{h_3}\rangle\Rightarrow H\cong BS(1,-1)$ is solvable. By Corollary \ref{n-sub2} $H \trianglelefteq G$ since $G=\langle x,h_3\rangle$. Since $\frac{G}{H}=\frac{\langle h_3\rangle H}{H}$ is a cyclic group, it is solvable. So $G$ is solvable, a contradiction.

\item[(14)]$ R_1:h_2^2h_3^{-2}h_2=1$, $R_2:h_2h_3^{-2}h_2^{-1}h_3^2=1$, $R_3:h_2h_3^{-1}(h_3^{-1}h_2)^3h_3h_2^{-1}h_3=1$:\\
$R_2:h_2h_3^{-2}h_2^{-1}h_3^2=1\Rightarrow h_3^2\in Z(G)$ where $G=\langle h_2,h_3\rangle$. Let $x=h_2h_3^{-1}$. So $h_3^{-1}h_2=x^{h_3}$. \\
$R_1 \text{ and } R_3 \Rightarrow (h_2h_3^{-1})^2=(h_3^{-1}h_2)^2 \text{ so if } H=\langle x,x^{h_3}\rangle\Rightarrow H\cong BS(1,-1)$ is solvable. By Corollary \ref{n-sub2} $H \trianglelefteq G$ since $G=\langle x,h_3\rangle$. Since $\frac{G}{H}=\frac{\langle h_3\rangle H}{H}$ is a cyclic group, it is solvable. So $G$ is solvable, a contradiction.

\item[(15)]$ R_1:h_2(h_3h_2^{-1})^2h_3=1$, $R_2:h_2^2h_3^{-1}h_2^{-2}h_3=1$, $R_3:h_2h_3^3h_2^{-1}h_3^{-1}h_2^{-1}h_3=1$:\\
$R_2\Rightarrow h_2^2\in Z(G)$ where $G=\langle h_2,h_3\rangle$. Let $x=h_3^{-1}h_2^{-1}$. So $h_2^{-1}h_3^{-1}=x^{h_2}$. \\
$R_1 \text{ and } R_3 \Rightarrow (h_3^{-1}h_2^{-1})^2=(h_2^{-1}h_3^{-1})^2 \text{ so if } H=\langle x,x^{h_2}\rangle\Rightarrow H\cong BS(1,-1)$ is solvable. By Corollary \ref{n-sub2} $H \trianglelefteq G$ since $G=\langle x,h_2\rangle$. Since $\frac{G}{H}=\frac{\langle h_2\rangle H}{H}$ is a cyclic group, it is solvable. So $G$ is solvable, a contradiction.

\item[(16)]$ R_1:h_2(h_3h_2^{-1})^2h_3=1$, $R_2:h_2^2h_3^{-1}h_2^{-2}h_3=1$, $R_3:h_2h_3^2(h_3h_2^{-1})^3h_3=1$:\\
$\Rightarrow \langle h_2,h_3\rangle=\langle h_3\rangle$ is abelian, a contradiction.

\item[(17)]$ R_1:h_2(h_3h_2^{-1})^2h_3=1$, $R_2:h_2^2h_3^{-1}h_2^{-2}h_3=1$, $R_3:h_2h_3(h_3h_2^{-1})^2h_3h_2h_3^{-1}h_2=1$:\\
$\Rightarrow h_2=1$ that is a contradiction.

\item[(18)]$ R_1:h_2(h_3h_2^{-1})^2h_3=1$, $R_2:h_2^2h_3^{-1}h_2^{-2}h_3=1$, $R_3:h_2h_3h_2^{-2}h_3h_2^{-1}(h_3^{-1}h_2)^2=1$:\\
$\Rightarrow\langle h_2,h_3\rangle\cong BS(1,1)$ is solvable, a contradiction.

\item[(19)]$ R_1:h_2(h_3h_2^{-1})^2h_3=1$, $R_2:h_2^2h_3^{-1}h_2^{-2}h_3=1$, $R_3:h_2h_3^{-1}h_2^{-1}h_3h_2^{-1}h_3^{-1}h_2h_3^{-2}h_2=1$:\\
$R_2\Rightarrow h_2^2\in Z(G)$ where $G=\langle h_2,h_3\rangle$. Let $x=h_3h_2^{-1}$. So $h_2^{-1}h_3=x^{h_2}$. \\
$R_1 \text{ and } R_3 \Rightarrow (h_3h_2^{-1})^2=(h_2^{-1}h_3)^2 \text{ so if } H=\langle x,x^{h_2}\rangle\Rightarrow H\cong BS(1,-1)$ is solvable. By Corollary \ref{n-sub2} $H \trianglelefteq G$ since $G=\langle x,h_2\rangle$. Since $\frac{G}{H}=\frac{\langle h_2\rangle H}{H}$ is a cyclic group, it is solvable. So $G$ is solvable, a contradiction.

\item[(20)]$ R_1:h_2(h_3h_2^{-1})^2h_3=1$, $R_2:h_2^2h_3^{-1}h_2^{-2}h_3=1$, $R_3:h_2h_3^{-2}h_2^{-1}h_3h_2^{-1}(h_3^{-1}h_2)^2=1$:\\
$R_2\Rightarrow h_2^2\in Z(G)$ where $G=\langle h_2,h_3\rangle$. Let $x=h_3h_2^{-1}$. So $h_2^{-1}h_3=x^{h_2}$. \\
$R_1 \text{ and } R_3 \Rightarrow (h_3h_2^{-1})^2=(h_2^{-1}h_3)^2 \text{ so if } H=\langle x,x^{h_2}\rangle\Rightarrow H\cong BS(1,-1)$ is solvable. By Corollary \ref{n-sub2} $H \trianglelefteq G$ since $G=\langle x,h_2\rangle$. Since $\frac{G}{H}=\frac{\langle h_2\rangle H}{H}$ is a cyclic group, it is solvable. So $G$ is solvable, a contradiction.

\item[(21)]$ R_1:h_2h_3^{-3}h_2=1$, $R_2:h_2^2h_3^{-1}h_2^{-2}h_3=1$, $R_3:h_2h_3h_2h_3^{-2}h_2^{-2}h_3=1$:\\
By interchanging $h_2$ and $h_3$ in (10) and with the same discussion, there is a contradiction.

\item[(22)]$ R_1:h_2h_3^{-3}h_2=1$, $R_2:h_2^2h_3^{-1}h_2^{-2}h_3=1$, $R_3:h_2h_3h_2h_3^{-2}(h_3^{-1}h_2)^2=1$:\\
By interchanging $h_2$ and $h_3$ in (11) and with the same discussion, there is a contradiction.

\item[(23)]$ R_1:h_2h_3^{-3}h_2=1$, $R_2:h_2^2h_3^{-1}h_2^{-2}h_3=1$, $R_3:h_2h_3^2h_2^{-1}h_3^{-1}(h_2^{-1}h_3)^2=1$:\\
By interchanging $h_2$ and $h_3$ in (12) and with the same discussion, there is a contradiction.

\item[(24)]$ R_1:h_2h_3^{-3}h_2=1$, $R_2:h_2^2h_3^{-1}h_2^{-2}h_3=1$, $R_3:h_2h_3h_2^{-1}h_3^2h_2^{-1}h_3^{-1}h_2^{-1}h_3=1$:\\
By interchanging $h_2$ and $h_3$ in (13) and with the same discussion, there is a contradiction.

\item[(25)]$ R_1:h_2h_3^{-3}h_2=1$, $R_2:h_2^2h_3^{-1}h_2^{-2}h_3=1$, $R_3:h_2h_3^{-1}h_2^{-1}h_3h_2^{-1}(h_3^{-1}h_2)^3=1$:\\
By interchanging $h_2$ and $h_3$ in (14) and with the same discussion, there is a contradiction.

\item[(26)]$ R_1:h_2h_3^{-3}h_2=1$, $R_2:h_2^2h_3^{-1}h_2^{-2}h_3=1$, $R_3:h_2h_3^{-3}(h_3^{-1}h_2)^3=1$:\\
By interchanging $h_2$ and $h_3$ in (9) and with the same discussion, there is a contradiction.

\item[(27)]$ R_1:h_2h_3^{-3}h_2=1$, $R_2:h_2^2h_3^{-1}h_2^{-2}h_3=1$, $R_3:h_2^4h_3^3=1$:\\
By interchanging $h_2$ and $h_3$ in (5) and with the same discussion, there is a contradiction.

\item[(28)]$ R_1:h_2h_3^{-3}h_2=1$, $R_2:h_2^2h_3^{-1}h_2^{-2}h_3=1$, $R_3:h_2^3h_3h_2^2h_3^{-1}h_2=1$:\\
By interchanging $h_2$ and $h_3$ in (7) and with the same discussion, there is a contradiction.

\item[(29)]$ R_1:h_2h_3^{-3}h_2=1$, $R_2:h_2^2h_3^{-1}h_2^{-2}h_3=1$, $R_3:h_2h_3^{-1}h_2h_3(h_2^{-1}h_3h_2^{-1})^2h_3=1$:\\
By interchanging $h_2$ and $h_3$ in (8) and with the same discussion, there is a contradiction.

\item[(30)]$ R_1:h_2h_3^{-3}h_2=1$, $R_2:h_2^2h_3^{-1}h_2^{-2}h_3=1$, $R_3:(h_2h_3^{-1})^2h_3^{-2}(h_3^{-1}h_2)^2=1$:\\
By interchanging $h_2$ and $h_3$ in (6) and with the same discussion, there is a contradiction.

\item[(31)]$ R_1:h_2h_3^{-1}h_2h_3^2=1$, $R_2:(h_2h_3^{-1})^2(h_2^{-1}h_3)^2=1$, $R_3:h_2^2h_3^{-1}h_2^2h_3^3=1$:\\
By interchanging $h_2$ and $h_3$ in (1) and with the same discussion, there is a contradiction.

\item[(32)]$ R_1:h_2h_3^{-1}h_2h_3^2=1$, $R_2:(h_2h_3^{-1})^2(h_2^{-1}h_3)^2=1$, $R_3:h_2^2h_3^{-1}h_2^2h_3h_2^{-2}h_3=1$:\\
By interchanging $h_2$ and $h_3$ in (3) and with the same discussion, there is a contradiction.

\item[(33)]$ R_1:h_2h_3^{-1}h_2h_3^2=1$, $R_2:(h_2h_3^{-1})^2(h_2^{-1}h_3)^2=1$, $R_3:(h_2h_3^{-1})^4h_3^{-1}h_2h_3^{-1}h_2=1$:\\
By interchanging $h_2$ and $h_3$ in (4) and with the same discussion, there is a contradiction.

\item[(34)]$ R_1:h_2h_3^{-1}h_2h_3^2=1$, $R_2:(h_2h_3^{-1})^2(h_2^{-1}h_3)^2=1$, $R_3:h_3^2(h_3h_2^{-1})^4h_3=1$:\\
By interchanging $h_2$ and $h_3$ in (2) and with the same discussion, there is a contradiction.

\item[(35)]$ R_1:(h_2h_3^{-1})^2h_2h_3=1$, $R_2:h_2h_3^{-2}h_2^{-1}h_3^2=1$, $R_3:h_2^3h_3^{-1}h_2^{-1}h_3^{-1}h_2h_3=1$:\\
By interchanging $h_2$ and $h_3$ in (15) and with the same discussion, there is a contradiction.

\item[(36)]$ R_1:(h_2h_3^{-1})^2h_2h_3=1$, $R_2:h_2h_3^{-2}h_2^{-1}h_3^2=1$, $R_3:h_2^2(h_2h_3^{-1})^3h_2h_3=1$:\\
By interchanging $h_2$ and $h_3$ in (16) and with the same discussion, there is a contradiction.

\item[(37)]$ R_1:(h_2h_3^{-1})^2h_2h_3=1$, $R_2:h_2h_3^{-2}h_2^{-1}h_3^2=1$, $R_3:h_2^2h_3^{-1}(h_3^{-1}h_2)^2h_3h_2^{-1}h_3=1$:\\
By interchanging $h_2$ and $h_3$ in (20) and with the same discussion, there is a contradiction.

\item[(38)]$ R_1:(h_2h_3^{-1})^2h_2h_3=1$, $R_2:h_2h_3^{-2}h_2^{-1}h_3^2=1$, $R_3:h_2(h_2h_3^{-1})^2h_2h_3h_2^{-1}h_3^2=1$:\\
By interchanging $h_2$ and $h_3$ in (17) and with the same discussion, there is a contradiction.

\item[(39)]$ R_1:(h_2h_3^{-1})^2h_2h_3=1$, $R_2:h_2h_3^{-2}h_2^{-1}h_3^2=1$, $R_3:h_2h_3^{-2}h_2^2h_3^{-1}h_2h_3h_2^{-1}h_3=1$:\\
By interchanging $h_2$ and $h_3$ in (19) and with the same discussion, there is a contradiction.

\item[(40)]$ R_1:(h_2h_3^{-1})^2h_2h_3=1$, $R_2:h_2h_3^{-2}h_2^{-1}h_3^2=1$, $R_3:h_2h_3^{-2}h_2h_3^{-1}(h_2^{-1}h_3)^2h_3=1$:\\
By interchanging $h_2$ and $h_3$ in (18) and with the same discussion, there is a contradiction.
\end{enumerate}

$\mathbf{C_5--C_5}$ \textbf{subgraph:} Suppose that there are two cycles of length $5$ in the graph $K(\alpha,\beta)$. By considering the relations from Table \ref{tab-C5} which are not disproved, it can be seen that there are $2485$ cases for existing two cycles $C_5$ in the graph $K(\alpha,\beta)$. Using Gap \cite{gap}, we see that all groups with two generators $h_2$ and $h_3$ and two relations which are between $2038$ cases of these $2485$ cases are solvable or have the same ``structure description'' $SL(2,5)$ according to the function StructureDescription of GAP, that is finite. So there are $447$ cases for the relations of two cycles of length $5$ in the graph $K(\alpha,\beta)$. Now suppose that the graph $K(\alpha,\beta)$ has a subgraph isomorphic to the graph $C_5--C_5$. Since this structure has two cycles of length $5$, the relations of these $C_5$ cycles must be between $447$ cases that have mentioned above.

Suppose that $[h_1',h_1'',h_2',h_2'',h_3',h_3'',h_4',h_4'',h_5',h_5'']$ and $[h_1',h_1'',h_2',h_2'',h_6',h_6'',h_7',h_7'',h_8',h_8'']$ are $10-$tuples related to the cycles $C_5$ in the graph $C_5--C_5$, where the first four components of these tuples are related to the common edges of $C_5$ and $C_5$. Without loss of generality we may assume that $h_1'=1$, where $h_1'',h_2',h_2'',h_3',h_3'',h_4',h_4'',h_5',h_5'',h_6',h_6'',h_7',h_7'',h_8',h_8'' \in supp(\alpha)$ and $\alpha=1+h_2+h_3$. With the same discussion such as about $K_{2,3}$, it is easy to see that $h_3' \neq h_6'$ and $h_5'' \neq h_8''$. 

With such assumptions and by the discussion above, it can be seen that there are $99$ cases which may lead to the existence of a subgraph isomorphic to the graph $C_5--C_5$ in the graph $K(\alpha,\beta)$. These cases are listed in table \ref{tab-C5--C5}. In the following, we show that $87$ cases of these relations lead to a contradiction and just $12$ cases of them may lead to the  existence of a subgraph isomorphic to the graph $C_5--C_5$ in the graph $K(\alpha,\beta)$. Cases which are not disproved are marked by $*$s in the Table \ref{tab-C5--C5}.
\begin{center}
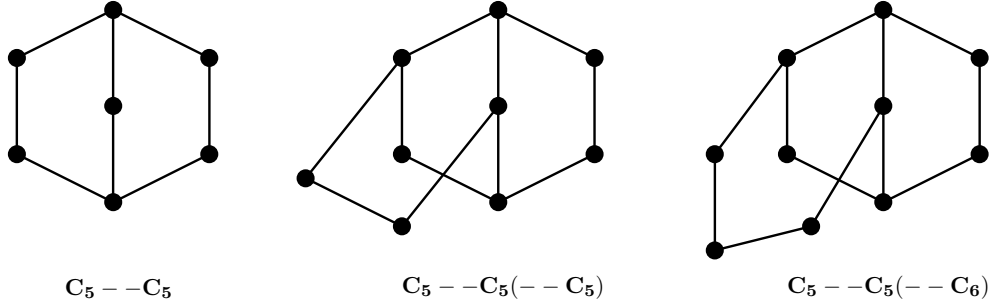
\begin{figure}[t]
\psscalebox{0.8 0.8} 
{
\begin{pspicture}(0,-2.4985578)(16.394232,2.4985578)
\psdots[linecolor=black, dotsize=0.3](1.7971154,2.3014424)
\psdots[linecolor=black, dotsize=0.3](1.7971154,0.7014423)
\psdots[linecolor=black, dotsize=0.3](1.7971154,-0.8985577)
\psdots[linecolor=black, dotsize=0.3](8.197115,2.3014424)
\psdots[linecolor=black, dotsize=0.3](8.197115,0.7014423)
\psdots[linecolor=black, dotsize=0.3](8.197115,-0.8985577)
\psdots[linecolor=black, dotsize=0.3](9.797115,1.5014423)
\psdots[linecolor=black, dotsize=0.3](9.797115,-0.0985577)
\psdots[linecolor=black, dotsize=0.3](6.5971155,1.5014423)
\psdots[linecolor=black, dotsize=0.3](6.5971155,-0.0985577)
\psdots[linecolor=black, dotsize=0.3](3.3971155,1.5014423)
\psdots[linecolor=black, dotsize=0.3](3.3971155,-0.0985577)
\psdots[linecolor=black, dotsize=0.3](0.1971154,-0.0985577)
\psdots[linecolor=black, dotsize=0.3](0.1971154,1.5014423)
\psline[linecolor=black, linewidth=0.04](1.7971154,2.3014424)(3.3971155,1.5014423)(3.3971155,-0.0985577)(1.7971154,-0.8985577)(0.1971154,-0.0985577)(0.1971154,1.5014423)(1.7971154,2.3014424)
\psline[linecolor=black, linewidth=0.04](1.7971154,2.3014424)(1.7971154,0.7014423)(1.7971154,-0.8985577)
\psline[linecolor=black, linewidth=0.04](8.197115,2.3014424)(9.797115,1.5014423)(9.797115,-0.0985577)(8.197115,-0.8985577)(6.5971155,-0.0985577)(6.5971155,1.5014423)(8.197115,2.3014424)(8.197115,0.7014423)(8.197115,-0.8985577)(8.197115,-0.8985577)
\psdots[linecolor=black, dotsize=0.3](6.5971155,-1.2985578)
\psline[linecolor=black, linewidth=0.04](8.197115,0.7014423)(6.5971155,-1.2985578)(4.9971156,-0.4985577)(6.5971155,1.5014423)(6.5971155,1.5014423)
\psdots[linecolor=black, dotsize=0.3](4.9971156,-0.4985577)
\psdots[linecolor=black, dotsize=0.3](14.5971155,2.3014424)
\psdots[linecolor=black, dotsize=0.3](14.5971155,0.7014423)
\psdots[linecolor=black, dotsize=0.3](14.5971155,-0.8985577)
\psdots[linecolor=black, dotsize=0.3](16.197115,-0.0985577)
\psdots[linecolor=black, dotsize=0.3](16.197115,1.5014423)
\psdots[linecolor=black, dotsize=0.3](12.997115,1.5014423)
\psdots[linecolor=black, dotsize=0.3](12.997115,-0.0985577)
\psdots[linecolor=black, dotsize=0.3](13.397116,-1.2985578)
\psdots[linecolor=black, dotsize=0.3](11.797115,-1.6985577)
\psdots[linecolor=black, dotsize=0.3](11.797115,-0.0985577)
\psline[linecolor=black, linewidth=0.04](12.997115,1.5014423)(14.5971155,2.3014424)(16.197115,1.5014423)(16.197115,-0.0985577)(14.5971155,-0.8985577)(12.997115,-0.0985577)(12.997115,1.5014423)(12.997115,1.5014423)
\psline[linecolor=black, linewidth=0.04](14.5971155,2.3014424)(14.5971155,0.7014423)(14.5971155,-0.8985577)(14.5971155,-0.8985577)
\psline[linecolor=black, linewidth=0.04](14.5971155,0.7014423)(13.397116,-1.2985578)(11.797115,-1.6985577)(11.797115,-0.0985577)(12.997115,1.5014423)
\rput[bl](0.9971154,-2.4985578){$\mathbf{C_5--C_5}$}
\rput[bl](6.5971155,-2.4985578){$\mathbf{C_5--C_5(--C_5)}$}
\rput[bl](12.997115,-2.4985578){$\mathbf{C_5--C_5(--C_6)}$}
\end{pspicture}
}
\vspace{1cc}
\caption{The graph $C_5--C_5$  and some forbidden subgraphs which contain it}\label{f-8}
\end{figure}
\end{center}
\begin{enumerate}
\item[(1)]$ R_1: h_2^3h_3^2=1$, $R_2: h_2^2h_3^{-1}h_2^{-2}h_3=1$:\\
$\Rightarrow\langle h_2,h_3\rangle\cong BS(1,1)$ is solvable, a contradiction.

\item[(2)]$ R_1: h_2^3h_3h_2^{-1}h_3=1$, $R_2: h_2^2h_3^{-1}h_2^{-2}h_3=1$:\\
$\Rightarrow \langle h_2,h_3\rangle=\langle h_2\rangle$ is abelian, a contradiction.

\item[(3)]$ R_1: h_2^3h_3h_2^{-1}h_3=1$, $R_2: h_2^2h_3^{-1}h_2^{-1}h_3^{-1}h_2=1$:\\
$\Rightarrow h_3=1$ that is a contradiction.

\item[(4)]$ R_1: h_2^3h_3h_2^{-1}h_3=1$, $R_2: h_2^2h_3^{-1}h_2^{-1}h_3^2=1$:\\
$\Rightarrow\langle h_2,h_3\rangle\cong BS(1,2)$ is solvable, a contradiction.

\item[(5)]$ R_1: h_2^3h_3h_2^{-1}h_3=1$, $R_2: h_2^2h_3^{-1}h_2^2h_3=1$:\\
$\Rightarrow h_2=1$ that is a contradiction.

\begin{center}
\tiny
\begin{longtable}{|c|l|l||c|l|l|}
\caption{{\footnotesize The relations of a $C_5--C_5$}}\label{tab-C5--C5}\\
\hline \multicolumn{1}{|c|}{$n$} & \multicolumn{1}{l|}{$R_1$} & \multicolumn{1}{l||}{$R_2$} & \multicolumn{1}{c|}{$n$} & \multicolumn{1}{l|}{$R_1$} & \multicolumn{1}{l|}{$R_2$}\\\hline
\endfirsthead
\multicolumn{6}{c}%
{{\tablename\ \thetable{} -- continued from previous page}} \\
\hline \multicolumn{1}{|c|}{$n$} & \multicolumn{1}{l|}{$R_1$} & \multicolumn{1}{l||}{$R_2$} & \multicolumn{1}{c|}{$n$} & \multicolumn{1}{l|}{$R_1$} & \multicolumn{1}{l|}{$R_2$}\\\hline
\endhead
\hline \multicolumn{6}{|r|}{{Continued on next page}} \\\hline
\endfoot
\hline
\endlastfoot
$1$&$h_2^3h_3^2=1$&$h_2^2h_3^{-1}h_2^{-2}h_3=1$&$51$&$h_2^2h_3^{-1}h_2^2h_3=1$&$h_2h_3^{-2}h_2h_3^2=1$\\
$2$&$h_2^3h_3h_2^{-1}h_3=1$&$h_2^2h_3^{-1}h_2^{-2}h_3=1$&$52$&$h_2^2h_3^{-1}h_2^2h_3=1$&$h_2h_3^{-1}h_2h_3^3=1$\\
$3$&$h_2^3h_3h_2^{-1}h_3=1$&$h_2^2h_3^{-1}h_2^{-1}h_3^{-1}h_2=1$&$53$&$h_2^2h_3^{-1}h_2^2h_3=1$&$(h_2h_3^{-1})^2h_2^{-1}h_3^2=1$\\
$4$&$h_2^3h_3h_2^{-1}h_3=1$&$h_2^2h_3^{-1}h_2^{-1}h_3^2=1$&$54$&$h_2^2h_3^{-1}h_2^2h_3=1$&$(h_2h_3^{-1})^2(h_2^{-1}h_3)^2=1$\\
$5$&$h_2^3h_3h_2^{-1}h_3=1$&$h_2^2h_3^{-1}h_2^2h_3=1$&$55$&$h_2^2h_3^{-1}h_2^2h_3=1$&$(h_2h_3^{-1})^2h_2h_3^2=1$\\
$6$&$h_2^3h_3h_2^{-1}h_3=1$&$h_2h_3^2h_2^{-1}h_3^2=1$&$56$&$h_2^2h_3^{-1}h_2^2h_3=1$&$(h_2h_3^{-1})^2h_2h_3h_2^{-1}h_3=1$\\
$7$&$h_2^3h_3h_2^{-1}h_3=1$&$h_2(h_3^{-1}h_2h_3^{-1})^2h_2=1$&$57$&$h_2^2h_3^{-1}h_2^2h_3=1$&$(h_2h_3^{-1})^3h_2h_3=1$\\
$8$&$h_2^3h_3^{-2}h_2=1$&$h_2^2h_3h_2^{-2}h_3=1$&$58$&$h_2^2h_3^{-1}h_2h_3^2=1$&$h_2h_3h_2h_3^{-1}h_2h_3=1$\\
$9$&$h_2^3h_3^{-2}h_2=1$&$h_2^2h_3^{-1}h_2^{-2}h_3=1$&$59$&$h_2^2h_3^{-1}h_2h_3^2=1$&$(h_2h_3^{-1})^2(h_2^{-1}h_3)^2=1$\\
$10$&$h_2^3h_3^{-2}h_2=1$&$h_2^2h_3^{-1}h_2^2h_3=1$&$60$&$h_2^2h_3^{-1}h_2h_3h_2^{-1}h_3=1$&$(h_2h_3^{-1})^2(h_2^{-1}h_3)^2=1$\\
$11$&$h_2^3h_3^{-2}h_2=1$&$h_2(h_2h_3^{-1})^2h_2^{-1}h_3=1$&$61$&$h_2(h_2h_3^{-1})^2h_2^{-1}h_3=1$&$h_2h_3^{-1}h_2^{-1}h_3h_2h_3^{-1}h_2=1$\\
$12$&$h_2^3h_3^{-2}h_2=1$&$h_2h_3h_2^{-1}(h_3^{-1}h_2)^2=1$&$62$&$h_2(h_2h_3^{-1})^2h_2h_3=1$&$h_2h_3^{-2}h_2^{-1}h_3^2=1$\\
$13$&$h_2^2(h_2h_3^{-1})^2h_2=1$&$h_2^2h_3^{-1}h_2^{-2}h_3=1$&$63$&$h_2(h_2h_3^{-1})^2h_2h_3=1$&$(h_2h_3^{-1}h_2)^2h_3=1$\\
$14$&$h_2^2h_3^3=1$&$h_2h_3^{-2}h_2^{-1}h_3^2=1$&$64$&$(h_2h_3)^2h_2^{-1}h_3=1$&$h_2h_3^{-2}h_2^{-1}h_3^2=1$\\
$15$&$h_2^2h_3^2h_2^{-1}h_3=1$&$(h_2h_3)^2h_2^{-1}h_3=1$&$65$&$h_2h_3h_2h_3^{-1}h_2h_3=1$&$h_2h_3^2h_2h_3^{-1}h_2=1$\\
$16$&$h_2^2h_3^2h_2^{-1}h_3=1$&$(h_2h_3^{-1})^2(h_2^{-1}h_3)^2=1$&$66$&$h_2h_3h_2h_3^{-1}h_2h_3=1$&$h_2h_3^2h_2^{-1}h_3^2=1$\\
$17$&$h_2^2h_3h_2^{-2}h_3=1$&$h_2h_3^2h_2^{-1}h_3^2=1$&$67 *$&$h_2h_3h_2h_3^{-1}h_2h_3=1$&$h_2h_3^{-2}h_2^{-1}h_3^2=1$\\
$18$&$h_2^2h_3h_2^{-2}h_3=1$&$h_2h_3^{-2}h_2^{-1}h_3^2=1$&$68$&$h_2h_3(h_2h_3^{-1})^2h_2=1$&$(h_2h_3^{-1}h_2)^2h_3=1$\\
$19$&$h_2^2h_3h_2^{-2}h_3=1$&$h_2h_3^{-2}h_2h_3^2=1$&$69$&$h_2h_3^2h_2^{-1}h_3^{-1}h_2=1$&$h_2h_3^{-1}h_2h_3^3=1$\\
$20$&$h_2^2h_3h_2^{-1}h_3^2=1$&$(h_2h_3)^2h_2^{-1}h_3=1$&$70$&$h_2h_3^2h_2^{-1}h_3^2=1$&$h_2(h_3h_2^{-1})^2h_3^{-1}h_2=1$\\
$21 *$&$h_2^2(h_3h_2^{-1})^2h_3=1$&$h_2^2h_3^{-3}h_2=1$&$71$&$h_2h_3^2h_2^{-1}h_3^2=1$&$h_2(h_3h_2^{-1})^3h_3=1$\\
$22$&$h_2^2(h_3h_2^{-1})^2h_3=1$&$h_2h_3^2h_2^{-1}h_3^2=1$&$72$&$h_2h_3^2h_2^{-1}h_3^2=1$&$h_2h_3^{-4}h_2=1$\\
$23 *$&$ h_2^2h_3^{-1}h_2^{-2}h_3=1$&$(h_2h_3)^2h_2^{-1}h_3=1$&$73$&$h_2h_3^2h_2^{-1}h_3^2=1$&$h_2h_3^{-3}h_2h_3=1$\\
$24$&$h_2^2h_3^{-1}h_2^{-2}h_3=1$&$h_2h_3h_2h_3^{-2}h_2=1$&$74$&$h_2h_3^2h_2^{-1}h_3^2=1$&$h_2h_3^{-1}h_2h_3^3=1$\\
$25$&$h_2^2h_3^{-1}h_2^{-2}h_3=1$&$h_2h_3h_2h_3^{-1}h_2h_3=1$&$75$&$h_2h_3^2h_2^{-1}h_3^2=1$&$h_2h_3^{-1}h_2(h_3h_2^{-1})^2h_3=1$\\
$26$&$h_2^2h_3^{-1}h_2^{-2}h_3=1$&$h_2h_3^2h_2^{-2}h_3=1$&$76$&$h_2h_3^2h_2^{-1}h_3^2=1$&$(h_2h_3^{-1})^2(h_2^{-1}h_3)^2=1$\\
$27$&$h_2^2h_3^{-1}h_2^{-2}h_3=1$&$h_2h_3^2h_2^{-1}h_3^{-1}h_2=1$&$77$&$h_2h_3^2h_2^{-1}h_3^2=1$&$(h_2h_3^{-1})^3h_2h_3=1$\\
$28$&$h_2^2h_3^{-1}h_2^{-2}h_3=1$&$h_2h_3^2h_2^{-1}h_3^2=1$&$78$&$h_2h_3(h_3h_2^{-1})^2h_3=1$&$h_2(h_3h_2^{-1}h_3)^2=1$\\
$29$&$h_2^2h_3^{-1}h_2^{-2}h_3=1$&$h_2h_3(h_3h_2^{-1})^2h_3=1$&$79$&$h_2h_3h_2^{-2}h_3^2=1$&$h_2h_3^{-2}h_2^{-1}h_3^2=1$\\
$30$&$h_2^2h_3^{-1}h_2^{-2}h_3=1$&$h_2h_3^{-2}h_2^{-1}h_3^2=1$&$80 *$&$h_2h_3h_2^{-1}(h_2^{-1}h_3)^2=1$&$h_2h_3^{-2}h_2^{-1}h_3^2=1$\\
$31$&$h_2^2h_3^{-1}h_2^{-2}h_3=1$&$h_2h_3^{-2}(h_2^{-1}h_3)^2=1$&$81$&$h_2h_3h_2^{-1}h_3^{-2}h_2=1$&$h_2h_3^{-1}(h_3^{-1}h_2)^3=1$\\
$32 *$&$h_2^2h_3^{-1}h_2^{-2}h_3=1$&$h_2h_3^{-3}h_2h_3=1$&$82$&$h_2(h_3h_2^{-1}h_3)^2=1$&$h_2(h_3h_2^{-1})^2h_3^2=1$\\
$33$&$h_2^2h_3^{-1}h_2^{-2}h_3=1$&$h_2h_3^{-2}h_2h_3^2=1$&$83 *$&$h_2(h_3h_2^{-1}h_3)^2=1$&$(h_2h_3^{-1})^2(h_2^{-1}h_3)^2=1$\\
$34$&$h_2^2h_3^{-1}h_2^{-2}h_3=1$&$h_2h_3^{-2}h_2h_3h_2^{-1}h_3=1$&$84$&$h_2(h_3h_2^{-1}h_3)^2=1$&$h_2(h_3^{-1}h_2h_3^{-1})^2h_2=1$\\
$35 *$&$h_2^2h_3^{-1}h_2^{-2}h_3=1$&$h_2h_3^{-1}(h_3^{-1}h_2)^2h_3=1$&$85$&$h_2h_3^{-1}h_2^{-1}h_3h_2^{-1}h_3^{-1}h_2=1$&$h_2h_3^{-2}h_2^{-1}h_3^2=1$\\
$36$&$h_2^2h_3^{-1}h_2^{-1}h_3^{-1}h_2=1$&$h_2^2h_3^{-1}h_2^2h_3=1$&$86$&$h_2h_3^{-1}(h_2^{-1}h_3)^2h_3=1$&$h_2h_3^{-4}h_2=1$\\
$37 *$&$h_2^2h_3^{-1}h_2^{-1}h_3^{-1}h_2=1$&$h_2h_3^{-2}h_2^{-1}h_3^2=1$&$87$&$h_2h_3^{-2}h_2^{-1}h_3^2=1$&$h_2h_3^{-4}h_2=1$\\
$38 *$&$h_2^2h_3^{-1}h_2^{-1}h_3^{-1}h_2=1$&$h_2h_3^{-3}h_2h_3=1$&$88$&$h_2h_3^{-2}h_2^{-1}h_3^2=1$&$h_2h_3^{-1}h_2h_3^3=1$\\
$39$&$h_2^2h_3^{-1}h_2^{-1}h_3^2=1$&$h_2h_3^2h_2^{-1}h_3^{-1}h_2=1$&$89$&$h_2h_3^{-2}h_2^{-1}h_3^2=1$&$h_3^2(h_3h_2^{-1})^2h_3=1$\\
$40$&$h_2^2h_3^{-1}h_2^{-1}h_3^2=1$&$h_2h_3^{-2}h_2^{-1}h_3^2=1$&$90$&$h_2h_3^{-4}h_2=1$&$h_2h_3^{-2}h_2h_3^2=1$\\
$41$&$h_2^2h_3^{-1}(h_2^{-1}h_3)^2=1$&$h_2h_3^{-2}h_2^{-1}h_3^2=1$&$91$&$h_2h_3^{-4}h_2=1$&$(h_2h_3^{-1})^2h_3^{-1}h_2^{-1}h_3=1$\\
$42$&$h_2^2h_3^{-2}h_2^{-1}h_3=1$&$h_2h_3h_2^{-1}h_3^{-1}h_2^{-1}h_3=1$&$92$&$h_2h_3^{-3}h_2h_3=1$&$h_2h_3^{-1}h_2h_3^3=1$\\
$43 *$&$h_2^2h_3^{-2}h_2^{-1}h_3=1$&$h_2h_3h_2^{-1}h_3^{-2}h_2=1$&$93$&$h_2h_3^{-2}h_2h_3^{-1}h_2^{-1}h_3=1$&$(h_2h_3^{-1})^2h_3^{-1}h_2^{-1}h_3=1$\\
$44$&$h_2^2h_3^{-3}h_2=1$&$h_2h_3h_2^{-1}h_3^{-1}h_2^{-1}h_3=1$&$94 *$&$(h_2h_3^{-1}h_2)^2h_3=1$&$(h_2h_3^{-1})^2(h_2^{-1}h_3)^2=1$\\
$45$&$h_2^2h_3^{-3}h_2=1$&$h_2h_3h_2^{-1}h_3^{-2}h_2=1$&$95$&$(h_2h_3^{-1}h_2)^2h_3=1$&$h_2(h_3^{-1}h_2h_3^{-1})^2h_2=1$\\
$46 *$&$h_2^2h_3^{-3}h_2=1$&$(h_2h_3^{-1})^2h_2h_3^2=1$&$96$&$h_2h_3^{-1}h_2h_3^3=1$&$h_2(h_3^{-1}h_2h_3^{-1})^2h_2=1$\\
$47$&$h_2^2h_3^{-2}h_2h_3=1$&$h_2h_3^{-2}h_2^{-1}h_3^2=1$&$97$&$h_2h_3^{-1}h_2h_3^2h_2^{-1}h_3=1$&$(h_2h_3^{-1})^2(h_2^{-1}h_3)^2=1$\\
$48$&$h_2^2h_3^{-1}h_2^2h_3=1$&$(h_2h_3)^2h_2^{-1}h_3=1$&$98$&$h_2h_3^{-1}h_2(h_3h_2^{-1})^2h_3=1$&$h_2(h_3^{-1}h_2h_3^{-1})^2h_2=1$\\
$49$&$h_2^2h_3^{-1}h_2^2h_3=1$&$h_2(h_3h_2^{-1})^3h_3=1$&$99$&$h_2(h_3^{-1}h_2h_3^{-1})^2h_2=1$&$(h_2h_3^{-1})^2h_2h_3h_2^{-1}h_3=1$\\
$50$&$h_2^2h_3^{-1}h_2^2h_3=1$&$h_2h_3^{-2}h_2^{-1}h_3^2=1$&$$&$$&$ $\\
\end{longtable}
\end{center}

\item[(6)]$ R_1: h_2^3h_3h_2^{-1}h_3=1$, $R_2: h_2h_3^2h_2^{-1}h_3^2=1$:\\
$\Rightarrow h_3=1$ that is a contradiction.

\item[(7)]$ R_1: h_2^3h_3h_2^{-1}h_3=1$, $R_2: h_2(h_3^{-1}h_2h_3^{-1})^2h_2=1$:\\
$\Rightarrow h_2=1$ that is a contradiction.

\item[(8)]$ R_1: h_2^3h_3^{-2}h_2=1$, $R_2: h_2^2h_3h_2^{-2}h_3=1$:\\
$\Rightarrow h_2=1$ that is a contradiction.

\item[(9)]$ R_1: h_2^3h_3^{-2}h_2=1$, $R_2: h_2^2h_3^{-1}h_2^{-2}h_3=1$:\\
$\Rightarrow \langle h_2,h_3\rangle=\langle h_2\rangle$ is abelian, a contradiction.

\item[(10)]$ R_1: h_2^3h_3^{-2}h_2=1$, $R_2: h_2^2h_3^{-1}h_2^2h_3=1$:\\
$\Rightarrow h_2=1$ that is a contradiction.

\item[(11)]$ R_1: h_2^3h_3^{-2}h_2=1$, $R_2: h_2(h_2h_3^{-1})^2h_2^{-1}h_3=1$:\\
$\Rightarrow\langle h_2,h_3\rangle\cong BS(2,1)$ is solvable, a contradiction.

\item[(12)]$ R_1: h_2^3h_3^{-2}h_2=1$, $R_2: h_2h_3h_2^{-1}(h_3^{-1}h_2)^2=1$:\\
$\Rightarrow\langle h_2,h_3\rangle\cong BS(1,2)$ is solvable, a contradiction.

\item[(13)]$ R_1: h_2^2(h_2h_3^{-1})^2h_2=1$, $R_2: h_2^2h_3^{-1}h_2^{-2}h_3=1$:\\
$\Rightarrow \langle h_2,h_3\rangle=\langle h_2,h_3h_2^{-2}\rangle=\langle h_3h_2^{-2}\rangle$ is abelian, a contradiction.

\item[(14)]$ R_1: h_2^2h_3^3=1$, $R_2: h_2h_3^{-2}h_2^{-1}h_3^2=1$:\\
By interchanging $h_2$ and $h_3$ in (1) and with the same discussion, there is a contradiction.

\item[(15)]$ R_1: h_2^2h_3^2h_2^{-1}h_3=1$, $R_2: (h_2h_3)^2h_2^{-1}h_3=1$:\\
$\Rightarrow\langle h_2,h_3\rangle\cong BS(1,1)$ is solvable, a contradiction.

\item[(16)]$ R_1: h_2^2h_3^2h_2^{-1}h_3=1$, $R_2: (h_2h_3^{-1})^2(h_2^{-1}h_3)^2=1$:\\
$\Rightarrow\langle h_2,h_3\rangle\cong BS(1,1)$ is solvable, a contradiction.

\item[(17)]$ R_1: h_2^2h_3h_2^{-2}h_3=1$, $R_2: h_2h_3^2h_2^{-1}h_3^2=1$:\\
$\Rightarrow h_3=1$ that is a contradiction.

\item[(18)]$ R_1: h_2^2h_3h_2^{-2}h_3=1$, $R_2: h_2h_3^{-2}h_2^{-1}h_3^2=1$:\\
$\Rightarrow h_3=1$ that is a contradiction.

\item[(19)]$ R_1: h_2^2h_3h_2^{-2}h_3=1$, $R_2: h_2h_3^{-2}h_2h_3^2=1$:\\
$\Rightarrow h_3=1$ that is a contradiction.

\item[(20)]$ R_1: h_2^2h_3h_2^{-1}h_3^2=1$, $R_2: (h_2h_3)^2h_2^{-1}h_3=1$:\\
$\Rightarrow\langle h_2,h_3\rangle\cong BS(1,1)$ is solvable, a contradiction.


\item[(22)]$ R_1: h_2^2(h_3h_2^{-1})^2h_3=1$, $R_2: h_2h_3^2h_2^{-1}h_3^2=1$:\\
$\Rightarrow h_3=1$ that is a contradiction.


\item[(24)]$ R_1: h_2^2h_3^{-1}h_2^{-2}h_3=1$, $R_2: h_2h_3h_2h_3^{-2}h_2=1$:\\
$\Rightarrow \langle h_2,h_3\rangle=\langle h_2\rangle$ is abelian, a contradiction.

\item[(25)]$ R_1: h_2^2h_3^{-1}h_2^{-2}h_3=1$, $R_2: h_2h_3h_2h_3^{-1}h_2h_3=1$:\\
$\Rightarrow\langle h_2,h_3\rangle\cong BS(1,-1)$ is solvable, a contradiction.

\item[(26)]$ R_1: h_2^2h_3^{-1}h_2^{-2}h_3=1$, $R_2: h_2h_3^2h_2^{-2}h_3=1$:\\
$\Rightarrow \langle h_2,h_3\rangle=\langle h_3\rangle$ is abelian, a contradiction.

\item[(27)]$ R_1: h_2^2h_3^{-1}h_2^{-2}h_3=1$, $R_2: h_2h_3^2h_2^{-1}h_3^{-1}h_2=1$:\\
$\Rightarrow \langle h_2,h_3\rangle=\langle h_2\rangle$ is abelian, a contradiction.

\item[(28)]$ R_1: h_2^2h_3^{-1}h_2^{-2}h_3=1$, $R_2: h_2h_3^2h_2^{-1}h_3^2=1$:\\
Using Tietze transformation where $h_3\mapsto h_3h_2$ and $h_2\mapsto h_2$, we have: $R_1\Rightarrow h_2^2\in Z(G)$ where $G=\langle h_2,h_3\rangle$. Let $x=h_2h_3$. So $h_2^{-1}h_3^{-1}=(x^{h_2})^{-1}$. \\
$R_2 \Rightarrow (h_2h_3)^2=(h_2^{-1}h_3^{-1})^2 \text{ so } H=\langle x,(x^{h_2})^{-1}\rangle=\langle x,x^{h_2}\rangle\cong BS(1,-1)$ is solvable. By Corollary \ref{n-sub2} $H \trianglelefteq G$ since $G=\langle x,h_2\rangle$. Since $\frac{G}{H}=\frac{\langle h_2\rangle H}{H}$ is a cyclic group, it is solvable. So $G$ is solvable, a contradiction.

\item[(29)]$ R_1: h_2^2h_3^{-1}h_2^{-2}h_3=1$, $R_2: h_2h_3(h_3h_2^{-1})^2h_3=1$:\\
$\Rightarrow\langle h_2,h_3\rangle\cong BS(1,1)$ is solvable, a contradiction.

\item[(30)]$ R_1: h_2^2h_3^{-1}h_2^{-2}h_3=1$, $R_2: h_2h_3^{-2}h_2^{-1}h_3^2=1$:\\
Using Tietze transformation where $h_3\mapsto h_3h_2^{-1}$ and $h_2\mapsto h_2$, we have: $R_1\Rightarrow h_2^2\in Z(G)$ where $G=\langle h_2,h_3\rangle$. Let $x=h_3h_2^{-1}$. So $h_2^{-1}h_3=x^{h_2}$. \\
$R_2 \Rightarrow (h_3h_2^{-1})^2=(h_2^{-1}h_3)^2 \text{ so if } H=\langle x,x^{h_2}\rangle\Rightarrow H\cong BS(1,-1)$ is solvable. By Corollary \ref{n-sub2} $H \trianglelefteq G$ since $G=\langle x,h_2\rangle$. Since $\frac{G}{H}=\frac{\langle h_2\rangle H}{H}$ is a cyclic group, it is solvable. So $G$ is solvable, a contradiction.

\item[(31)]$ R_1: h_2^2h_3^{-1}h_2^{-2}h_3=1$, $R_2: h_2h_3^{-2}(h_2^{-1}h_3)^2=1$:\\
$\Rightarrow h_3=1$ that is a contradiction.


\item[(33)]$ R_1: h_2^2h_3^{-1}h_2^{-2}h_3=1$, $R_2: h_2h_3^{-2}h_2h_3^2=1$:\\
By interchanging $h_2$ and $h_3$ in (18) and with the same discussion, there is a contradiction.

\item[(34)]$ R_1: h_2^2h_3^{-1}h_2^{-2}h_3=1$, $R_2: h_2h_3^{-2}h_2h_3h_2^{-1}h_3=1$:\\
$\Rightarrow h_2=1$ that is a contradiction.



\item[(36)]$ R_1: h_2^2h_3^{-1}h_2^{-1}h_3^{-1}h_2=1$, $R_2: h_2^2h_3^{-1}h_2^2h_3=1$:\\
$\Rightarrow h_3=1$ that is a contradiction.



\item[(39)]$ R_1: h_2^2h_3^{-1}h_2^{-1}h_3^2=1$, $R_2: h_2h_3^2h_2^{-1}h_3^{-1}h_2=1$:\\
$\Rightarrow h_3=1$ that is a contradiction.

\item[(40)]$ R_1: h_2^2h_3^{-1}h_2^{-1}h_3^2=1$, $R_2: h_2h_3^{-2}h_2^{-1}h_3^2=1$:\\
By interchanging $h_2$ and $h_3$ in (27) and with the same discussion, there is a contradiction.

\item[(41)]$ R_1: h_2^2h_3^{-1}(h_2^{-1}h_3)^2=1$, $R_2: h_2h_3^{-2}h_2^{-1}h_3^2=1$:\\
By interchanging $h_2$ and $h_3$ in (31) and with the same discussion, there is a contradiction.

\item[(42)]$ R_1: h_2^2h_3^{-2}h_2^{-1}h_3=1$, $R_2: h_2h_3h_2^{-1}h_3^{-1}h_2^{-1}h_3=1$:\\
$\Rightarrow \langle h_2,h_3\rangle=\langle h_3\rangle$ is abelian, a contradiction.


\item[(44)]$ R_1: h_2^2h_3^{-3}h_2=1$, $R_2: h_2h_3h_2^{-1}h_3^{-1}h_2^{-1}h_3=1$:\\
$\Rightarrow \langle h_2,h_3\rangle=\langle h_2,h_2^{-1}h_3\rangle=\langle h_2^{-1}h_3\rangle$ is abelian, a contradiction.

\item[(45)]$ R_1: h_2^2h_3^{-3}h_2=1$, $R_2: h_2h_3h_2^{-1}h_3^{-2}h_2=1$:\\
$\Rightarrow \langle h_2,h_3\rangle=\langle h_2\rangle$ is abelian, a contradiction.


\item[(47)]$ R_1: h_2^2h_3^{-2}h_2h_3=1$, $R_2: h_2h_3^{-2}h_2^{-1}h_3^2=1$:\\
By interchanging $h_2$ and $h_3$ in (26) and with the same discussion, there is a contradiction.

\item[(48)]$ R_1: h_2^2h_3^{-1}h_2^2h_3=1$, $R_2: (h_2h_3)^2h_2^{-1}h_3=1$:\\
$\Rightarrow h_3=1$ that is a contradiction.

\item[(49)]$ R_1: h_2^2h_3^{-1}h_2^2h_3=1$, $R_2: h_2(h_3h_2^{-1})^3h_3=1$:\\
$\Rightarrow h_2=1$ that is a contradiction.

\item[(50)]$ R_1: h_2^2h_3^{-1}h_2^2h_3=1$, $R_2: h_2h_3^{-2}h_2^{-1}h_3^2=1$:\\
By interchanging $h_2$ and $h_3$ in (28) and with the same discussion, there is a contradiction.

\item[(51)]$ R_1: h_2^2h_3^{-1}h_2^2h_3=1$, $R_2: h_2h_3^{-2}h_2h_3^2=1$:\\
By interchanging $h_2$ and $h_3$ in (17) and with the same discussion, there is a contradiction.

\item[(52)]$ R_1: h_2^2h_3^{-1}h_2^2h_3=1$, $R_2: h_2h_3^{-1}h_2h_3^3=1$:\\
By interchanging $h_2$ and $h_3$ in (6) and with the same discussion, there is a contradiction.

\item[(53)]$ R_1: h_2^2h_3^{-1}h_2^2h_3=1$, $R_2: (h_2h_3^{-1})^2h_2^{-1}h_3^2=1$:\\
$\Rightarrow h_3=1$ that is a contradiction.

\item[(54)]$ R_1: h_2^2h_3^{-1}h_2^2h_3=1$, $R_2: (h_2h_3^{-1})^2(h_2^{-1}h_3)^2=1$:\\
Using Tietze transformation where $h_3\mapsto h_2h_3$ and $h_2\mapsto h_2$, we have $ R_1: h_2^2h_3^{-1}h_2^2h_3=1$, $R_2: h_2h_3^{-2}h_2^{-1}h_3^2=1$. So by the same discussion such as item 50, there is a contradiction.

\item[(55)]$ R_1: h_2^2h_3^{-1}h_2^2h_3=1$, $R_2: (h_2h_3^{-1})^2h_2h_3^2=1$:\\
By interchanging $h_2$ and $h_3$ in (22) and with the same discussion, there is a contradiction.

\item[(56)]$ R_1: h_2^2h_3^{-1}h_2^2h_3=1$, $R_2: (h_2h_3^{-1})^2h_2h_3h_2^{-1}h_3=1$:\\
$\Rightarrow h_3=1$ that is a contradiction.

\item[(57)]$ R_1: h_2^2h_3^{-1}h_2^2h_3=1$, $R_2: (h_2h_3^{-1})^3h_2h_3=1$:\\
$\Rightarrow h_3=1$ that is a contradiction.

\item[(58)]$ R_1: h_2^2h_3^{-1}h_2h_3^2=1$, $R_2: h_2h_3h_2h_3^{-1}h_2h_3=1$:\\
By interchanging $h_2$ and $h_3$ in (15) and with the same discussion, there is a contradiction.

\item[(59)]$ R_1: h_2^2h_3^{-1}h_2h_3^2=1$, $R_2: (h_2h_3^{-1})^2(h_2^{-1}h_3)^2=1$:\\
By interchanging $h_2$ and $h_3$ in (16) and with the same discussion, there is a contradiction.

\item[(60)]$ R_1: h_2^2h_3^{-1}h_2h_3h_2^{-1}h_3=1$, $R_2: (h_2h_3^{-1})^2(h_2^{-1}h_3)^2=1$:\\
$\Rightarrow \langle h_2,h_3\rangle=\langle h_2\rangle$ is abelian, a contradiction.

\item[(61)]$ R_1: h_2(h_2h_3^{-1})^2h_2^{-1}h_3=1$, $R_2: h_2h_3^{-1}h_2^{-1}h_3h_2h_3^{-1}h_2=1$:\\
$\Rightarrow\langle h_2,h_3\rangle\cong BS(1,1)$ is solvable, a contradiction.

\item[(62)]$ R_1: h_2(h_2h_3^{-1})^2h_2h_3=1$, $R_2: h_2h_3^{-2}h_2^{-1}h_3^2=1$:\\
By interchanging $h_2$ and $h_3$ in (29) and with the same discussion, there is a contradiction.

\item[(63)]$ R_1: h_2(h_2h_3^{-1})^2h_2h_3=1$, $R_2: (h_2h_3^{-1}h_2)^2h_3=1$:\\
$\Rightarrow\langle h_2,h_3\rangle\cong BS(1,1)$ is solvable, a contradiction.

\item[(64)]$ R_1: (h_2h_3)^2h_2^{-1}h_3=1$, $R_2: h_2h_3^{-2}h_2^{-1}h_3^2=1$:\\
By interchanging $h_2$ and $h_3$ in (25) and with the same discussion, there is a contradiction.

\item[(65)]$ R_1: h_2h_3h_2h_3^{-1}h_2h_3=1$, $R_2: h_2h_3^2h_2h_3^{-1}h_2=1$:\\
By interchanging $h_2$ and $h_3$ in (20) and with the same discussion, there is a contradiction.

\item[(66)]$ R_1: h_2h_3h_2h_3^{-1}h_2h_3=1$, $R_2: h_2h_3^2h_2^{-1}h_3^2=1$:\\
By interchanging $h_2$ and $h_3$ in (48) and with the same discussion, there is a contradiction.


\item[(68)]$ R_1: h_2h_3(h_2h_3^{-1})^2h_2=1$, $R_2: (h_2h_3^{-1}h_2)^2h_3=1$:\\
$\Rightarrow\langle h_2,h_3\rangle\cong BS(1,1)$ is solvable, a contradiction.

\item[(69)]$ R_1: h_2h_3^2h_2^{-1}h_3^{-1}h_2=1$, $R_2: h_2h_3^{-1}h_2h_3^3=1$:\\
By interchanging $h_2$ and $h_3$ in (4) and with the same discussion, there is a contradiction.

\item[(70)]$ R_1: h_2h_3^2h_2^{-1}h_3^2=1$, $R_2: h_2(h_3h_2^{-1})^2h_3^{-1}h_2=1$:\\
By interchanging $h_2$ and $h_3$ in (53) and with the same discussion, there is a contradiction.

\item[(71)]$ R_1: h_2h_3^2h_2^{-1}h_3^2=1$, $R_2: h_2(h_3h_2^{-1})^3h_3=1$:\\
By interchanging $h_2$ and $h_3$ in (57) and with the same discussion, there is a contradiction.

\item[(72)]$ R_1: h_2h_3^2h_2^{-1}h_3^2=1$, $R_2: h_2h_3^{-4}h_2=1$:\\
By interchanging $h_2$ and $h_3$ in (10) and with the same discussion, there is a contradiction.

\item[(73)]$ R_1: h_2h_3^2h_2^{-1}h_3^2=1$, $R_2: h_2h_3^{-3}h_2h_3=1$:\\
By interchanging $h_2$ and $h_3$ in (36) and with the same discussion, there is a contradiction.

\item[(74)]$ R_1: h_2h_3^2h_2^{-1}h_3^2=1$, $R_2: h_2h_3^{-1}h_2h_3^3=1$:\\
By interchanging $h_2$ and $h_3$ in (5) and with the same discussion, there is a contradiction.

\item[(75)]$ R_1: h_2h_3^2h_2^{-1}h_3^2=1$, $R_2: h_2h_3^{-1}h_2(h_3h_2^{-1})^2h_3=1$:\\
By interchanging $h_2$ and $h_3$ in (56) and with the same discussion, there is a contradiction.

\item[(76)]$ R_1: h_2h_3^2h_2^{-1}h_3^2=1$, $R_2: (h_2h_3^{-1})^2(h_2^{-1}h_3)^2=1$:\\
By interchanging $h_2$ and $h_3$ in (54) and with the same discussion, there is a contradiction.

\item[(77)]$ R_1: h_2h_3^2h_2^{-1}h_3^2=1$, $R_2: (h_2h_3^{-1})^3h_2h_3=1$:\\
By interchanging $h_2$ and $h_3$ in (49) and with the same discussion, there is a contradiction.

\item[(78)]$ R_1: h_2h_3(h_3h_2^{-1})^2h_3=1$, $R_2: h_2(h_3h_2^{-1}h_3)^2=1$:\\
By interchanging $h_2$ and $h_3$ in (63) and with the same discussion, there is a contradiction.

\item[(79)]$ R_1: h_2h_3h_2^{-2}h_3^2=1$, $R_2: h_2h_3^{-2}h_2^{-1}h_3^2=1$:\\
By interchanging $h_2$ and $h_3$ in (24) and with the same discussion, there is a contradiction.



\item[(81)]$ R_1: h_2h_3h_2^{-1}h_3^{-2}h_2=1$, $R_2: h_2h_3^{-1}(h_3^{-1}h_2)^3=1$:\\
$\Rightarrow\langle h_2,h_3\rangle\cong BS(-3,1)$ is solvable, a contradiction.

\item[(82)]$ R_1: h_2(h_3h_2^{-1}h_3)^2=1$, $R_2: h_2(h_3h_2^{-1})^2h_3^2=1$:\\
By interchanging $h_2$ and $h_3$ in (68) and with the same discussion, there is a contradiction.


\item[(84)]$ R_1: h_2(h_3h_2^{-1}h_3)^2=1$, $R_2: h_2(h_3^{-1}h_2h_3^{-1})^2h_2=1$:\\
$\Rightarrow h_2=1$ that is a contradiction.


\item[(85)]$ R_1: h_2h_3^{-1}h_2^{-1}h_3h_2^{-1}h_3^{-1}h_2=1$, $R_2: h_2h_3^{-2}h_2^{-1}h_3^2=1$:\\
By interchanging $h_2$ and $h_3$ in (34) and with the same discussion, there is a contradiction.

\item[(86)]$ R_1: h_2h_3^{-1}(h_2^{-1}h_3)^2h_3=1$, $R_2: h_2h_3^{-4}h_2=1$:\\
By interchanging $h_2$ and $h_3$ in (12) and with the same discussion, there is a contradiction.

\item[(87)]$ R_1: h_2h_3^{-2}h_2^{-1}h_3^2=1$, $R_2: h_2h_3^{-4}h_2=1$:\\
By interchanging $h_2$ and $h_3$ in (9) and with the same discussion, there is a contradiction.

\item[(88)]$ R_1: h_2h_3^{-2}h_2^{-1}h_3^2=1$, $R_2: h_2h_3^{-1}h_2h_3^3=1$:\\
By interchanging $h_2$ and $h_3$ in (2) and with the same discussion, there is a contradiction.

\item[(89)]$ R_1: h_2h_3^{-2}h_2^{-1}h_3^2=1$, $R_2: h_3^2(h_3h_2^{-1})^2h_3=1$:\\
By interchanging $h_2$ and $h_3$ in (13) and with the same discussion, there is a contradiction.

\item[(90)]$ R_1: h_2h_3^{-4}h_2=1$, $R_2: h_2h_3^{-2}h_2h_3^2=1$:\\
By interchanging $h_2$ and $h_3$ in (8) and with the same discussion, there is a contradiction.

\item[(91)]$ R_1: h_2h_3^{-4}h_2=1$, $R_2: (h_2h_3^{-1})^2h_3^{-1}h_2^{-1}h_3=1$:\\
By interchanging $h_2$ and $h_3$ in (11) and with the same discussion, there is a contradiction.

\item[(92)]$ R_1: h_2h_3^{-3}h_2h_3=1$, $R_2: h_2h_3^{-1}h_2h_3^3=1$:\\
By interchanging $h_2$ and $h_3$ in (3) and with the same discussion, there is a contradiction.

\item[(93)]$ R_1: h_2h_3^{-2}h_2h_3^{-1}h_2^{-1}h_3=1$, $R_2: (h_2h_3^{-1})^2h_3^{-1}h_2^{-1}h_3=1$:\\
By interchanging $h_2$ and $h_3$ in (61) and with the same discussion, there is a contradiction.


\item[(95)]$ R_1: (h_2h_3^{-1}h_2)^2h_3=1$, $R_2: h_2(h_3^{-1}h_2h_3^{-1})^2h_2=1$:\\
By interchanging $h_2$ and $h_3$ in (84) and with the same discussion, there is a contradiction.

\item[(96)]$ R_1: h_2h_3^{-1}h_2h_3^3=1$, $R_2: h_2(h_3^{-1}h_2h_3^{-1})^2h_2=1$:\\
By interchanging $h_2$ and $h_3$ in (7) and with the same discussion, there is a contradiction.

\item[(97)]$ R_1: h_2h_3^{-1}h_2h_3^2h_2^{-1}h_3=1$, $R_2: (h_2h_3^{-1})^2(h_2^{-1}h_3)^2=1$:\\
By interchanging $h_2$ and $h_3$ in (60) and with the same discussion, there is a contradiction.

\item[(98)]$ R_1: h_2h_3^{-1}h_2(h_3h_2^{-1})^2h_3=1$, $R_2: h_2(h_3^{-1}h_2h_3^{-1})^2h_2=1$:\\
$\Rightarrow h_3=1$ that is a contradiction.

\item[(99)]$ R_1: h_2(h_3^{-1}h_2h_3^{-1})^2h_2=1$, $R_2: (h_2h_3^{-1})^2h_2h_3h_2^{-1}h_3=1$:\\
By interchanging $h_2$ and $h_3$ in (98) and with the same discussion, there is a contradiction.
\end{enumerate}

\subsection{$\mathbf{C_5--C_5(--C_5)}$}
It can be seen that there are $192$ cases for the relations of three cycles $C_5$ in the graph $C_5--C_5(--C_5)$. Using GAP \cite{gap}, we see that all groups with two generators $h_2$ and $h_3$ and three relations which are between $188$ cases of these $192$ cases are finite and solvable, that is a contradiction with the assumptions. So, there are just $4$ cases for the relations of these cycles which may lead to the existence of a subgraph isomorphic to the graph $C_5--C_5(--C_5)$ in $K(\alpha,\beta)$. These cases are listed in table \ref{tab-C5--C5(--C5)}.  In the following, we how that all of these $4$ cases lead to contradictions and so, the graph $K(\alpha,\beta)$ contains no subgraph isomorphic to the graph $C_5--C_5(--C_5)$.

\begin{enumerate}
\item[(1)]$ R_1:h_2^2h_3^{-1}h_2^{-2}h_3=1$, $R_2:h_2h_3^{-1}(h_3^{-1}h_2)^2h_3=1$, $R_3:h_2h_3^{-2}h_2^{-1}h_3^2=1$:\\
$\Rightarrow\langle h_2,h_3\rangle\cong BS(1,1)$ is solvable, a contradiction.

\item[(2)]$ R_1:h_2h_3h_2^{-1}(h_2^{-1}h_3)^2=1$, $R_2:h_2h_3^{-2}h_2^{-1}h_3^2=1$, $R_3:h_2^2h_3^{-1}h_2^{-2}h_3=1$:\\
By interchanging $h_2$ and $h_3$ in (1) and with the same discussion, there is a contradiction.

\item[(3)]$ R_1:h_2(h_3h_2^{-1}h_3)^2=1$, $R_2:(h_2h_3^{-1})^2(h_2^{-1}h_3)^2=1$, $R_3:h_2h_3^{-2}h_2^{-1}h_3^2=1$:\\
$\Rightarrow\langle h_2,h_3\rangle\cong BS(1,-1)$ is solvable, a contradiction.

\item[(4)]$ R_1:(h_2h_3^{-1}h_2)^2h_3=1$, $R_2:(h_2h_3^{-1})^2(h_2^{-1}h_3)^2=1$, $R_3:h_2^2h_3^{-1}h_2^{-2}h_3=1$:\\
By interchanging $h_2$ and $h_3$ in (3) and with the same discussion, there is a contradiction.

\end{enumerate}
\begin{table}[h]
\centering
\caption{The relations of a $C_5--C_5(--C_5)$ in $K(\alpha,\beta)$}\label{tab-C5--C5(--C5)}
\begin{tabular}{|c|l|l|l|}\hline
$n$&$R_1$&$R_2$&$R_3$\\\hline
$1$&$h_2^2h_3^{-1}h_2^{-2}h_3=1$&$h_2h_3^{-1}(h_3^{-1}h_2)^2h_3=1$&$h_2h_3^{-2}h_2^{-1}h_3^2=1$\\
$2$&$h_2h_3h_2^{-1}(h_2^{-1}h_3)^2=1$&$h_2h_3^{-2}h_2^{-1}h_3^2=1$&$h_2^2h_3^{-1}h_2^{-2}h_3=1$\\
$3$&$h_2(h_3h_2^{-1}h_3)^2=1$&$(h_2h_3^{-1})^2(h_2^{-1}h_3)^2=1$&$h_2h_3^{-2}h_2^{-1}h_3^2=1$\\
$4$&$(h_2h_3^{-1}h_2)^2h_3=1$&$(h_2h_3^{-1})^2(h_2^{-1}h_3)^2=1$&$h_2^2h_3^{-1}h_2^{-2}h_3=1$\\
\hline
\end{tabular}
\end{table}

\subsection{$\mathbf{C_5--C_5(--C_6)}$}
It can be seen that there are $1006$ cases for the relations of two cycles $C_5$ and a cycle $C_6$ in the graph $C_5--C_5(--C_6)$. Using GAP \cite{gap}, we see that all groups with two generators $h_2$ and $h_3$ and three relations which are between $986$ cases of these $1006$ cases are finite and solvable, that is a contradiction with the assumptions. So, there are just $20$ cases for the relations of these cycles which may lead to the existence of a subgraph isomorphic to the graph $C_5--C_5(--C_6)$ in $K(\alpha,\beta)$. These cases are listed in table \ref{tab-C5--C5(--C6)}.  In the following, we show that all of these $20$ cases lead to contradictions and so, the graph $K(\alpha,\beta)$ contains no subgraph isomorphic to the graph $C_5--C_5(--C_6)$.

\begin{table}[h]
\centering
\caption{The relations of a $C5--C5(--C6)$ in $K(\alpha,\beta)$}\label{tab-C5--C5(--C6)}
\begin{tabular}{|c|l|l|l|}\hline
$n$&$R_1$&$R_2$&$R_3$\\\hline
$1$&$h_2^2h_3^{-1}h_2^{-2}h_3=1$&$(h_2h_3)^2h_2^{-1}h_3=1$&$h_2h_3^{-1}h_2^{-1}h_3^2h_2^{-1}h_3^{-1}h_2=1$\\
$2$&$h_2^2h_3^{-1}h_2^{-2}h_3=1$&$(h_2h_3)^2h_2^{-1}h_3=1$&$h_2h_3^{-1}h_2^{-1}h_3h_2^{-1}h_3^{-1}h_2h_3=1$\\
$3$&$h_2^2h_3^{-1}h_2^{-2}h_3=1$&$h_2h_3^{-1}(h_3^{-1}h_2)^2h_3=1$&$h_2^3h_3^{-1}h_2^{-1}h_3^{-1}h_2=1$\\
$4$&$h_2^2h_3^{-1}h_2^{-2}h_3=1$&$h_2h_3^{-1}(h_3^{-1}h_2)^2h_3=1$&$h_2^2h_3^{-3}h_2h_3=1$\\
$5$&$h_2^2h_3^{-1}h_2^{-2}h_3=1$&$h_2h_3^{-1}(h_3^{-1}h_2)^2h_3=1$&$h_2h_3h_2h_3^{-3}h_2=1$\\
$6$&$h_2^2h_3^{-1}h_2^{-2}h_3=1$&$h_2h_3^{-1}(h_3^{-1}h_2)^2h_3=1$&$h_2(h_3h_2^{-1})^2h_3^{-2}h_2=1$\\
$7$&$h_2h_3h_2h_3^{-1}h_2h_3=1$&$h_2h_3^{-2}h_2^{-1}h_3^2=1$&$h_2h_3h_2^{-1}h_3^{-1}h_2h_3^{-1}h_2^{-1}h_3=1$\\
$8$&$h_2h_3h_2h_3^{-1}h_2h_3=1$&$h_2h_3^{-2}h_2^{-1}h_3^2=1$&$h_2h_3^{-1}h_2^{-1}h_3^2h_2^{-1}h_3^{-1}h_2=1$\\
$9$&$h_2h_3h_2^{-1}(h_2^{-1}h_3)^2=1$&$h_2h_3^{-2}h_2^{-1}h_3^2=1$&$h_2h_3^{-4}h_2h_3=1$\\
$10$&$h_2h_3h_2^{-1}(h_2^{-1}h_3)^2=1$&$h_2h_3^{-2}h_2^{-1}h_3^2=1$&$h_2^2h_3^{-1}h_2^{-1}h_3^{-2}h_2=1$\\
$11$&$h_2h_3h_2^{-1}(h_2^{-1}h_3)^2=1$&$h_2h_3^{-2}h_2^{-1}h_3^2=1$&$h_2^2h_3^{-2}h_2^{-1}h_3^{-1}h_2=1$\\
$12$&$h_2h_3h_2^{-1}(h_2^{-1}h_3)^2=1$&$h_2h_3^{-2}h_2^{-1}h_3^2=1$&$h_2(h_3h_2^{-1})^2h_3^{-2}h_2=1$\\
$13$&$h_2(h_3h_2^{-1}h_3)^2=1$&$(h_2h_3^{-1})^2(h_2^{-1}h_3)^2=1$&$h_2h_3^2h_2^{-2}h_3^2=1$\\
$14$&$h_2(h_3h_2^{-1}h_3)^2=1$&$(h_2h_3^{-1})^2(h_2^{-1}h_3)^2=1$&$h_2h_3h_2^{-2}h_3^3=1$\\
$15$&$h_2(h_3h_2^{-1}h_3)^2=1$&$(h_2h_3^{-1})^2(h_2^{-1}h_3)^2=1$&$h_2(h_3h_2^{-1})^2h_3^{-2}h_2=1$\\
$16$&$h_2(h_3h_2^{-1}h_3)^2=1$&$(h_2h_3^{-1})^2(h_2^{-1}h_3)^2=1$&$h_2h_3^{-1}h_2^{-1}h_3h_2^{-1}h_3^{-1}h_2h_3=1$\\
$17$&$(h_2h_3^{-1}h_2)^2h_3=1$&$(h_2h_3^{-1})^2(h_2^{-1}h_3)^2=1$&$h_2^2h_3h_2h_3^{-2}h_2=1$\\
$18$&$(h_2h_3^{-1}h_2)^2h_3=1$&$(h_2h_3^{-1})^2(h_2^{-1}h_3)^2=1$&$h_2^2h_3^{-2}h_2^2h_3=1$\\
$19$&$(h_2h_3^{-1}h_2)^2h_3=1$&$(h_2h_3^{-1})^2(h_2^{-1}h_3)^2=1$&$h_2h_3h_2^{-1}h_3^{-1}h_2h_3^{-1}h_2^{-1}h_3=1$\\
$20$&$(h_2h_3^{-1}h_2)^2h_3=1$&$(h_2h_3^{-1})^2(h_2^{-1}h_3)^2=1$&$h_2(h_3h_2^{-1})^2h_3^{-2}h_2=1$\\
\hline
\end{tabular}
\end{table}

\begin{enumerate}
\item[(1)]$ R_1:h_2^2h_3^{-1}h_2^{-2}h_3=1$, $R_2:(h_2h_3)^2h_2^{-1}h_3=1$, $R_3:h_2h_3^{-1}h_2^{-1}h_3^2h_2^{-1}h_3^{-1}h_2=1$:\\
$\Rightarrow\langle h_2,h_3\rangle\cong BS(1,-1)$ is solvable, a contradiction.

\item[(2)]$ R_1:h_2^2h_3^{-1}h_2^{-2}h_3=1$, $R_2:(h_2h_3)^2h_2^{-1}h_3=1$, $R_3:h_2h_3^{-1}h_2^{-1}h_3h_2^{-1}h_3^{-1}h_2h_3=1$:\\
$\Rightarrow \langle h_2,h_3\rangle=\langle h_3\rangle$ is abelian, a contradiction.

\item[(3)]$ R_1:h_2^2h_3^{-1}h_2^{-2}h_3=1$, $R_2:h_2h_3^{-1}(h_3^{-1}h_2)^2h_3=1$, $R_3:h_2^3h_3^{-1}h_2^{-1}h_3^{-1}h_2=1$:\\
$\Rightarrow\langle h_2,h_3\rangle\cong BS(1,1)$ is solvable, a contradiction.

\item[(4)]$ R_1:h_2^2h_3^{-1}h_2^{-2}h_3=1$, $R_2:h_2h_3^{-1}(h_3^{-1}h_2)^2h_3=1$, $R_3:h_2^2h_3^{-3}h_2h_3=1$:\\
$\Rightarrow\langle h_2,h_3\rangle\cong BS(1,1)$ is solvable, a contradiction.

\item[(5)]$ R_1:h_2^2h_3^{-1}h_2^{-2}h_3=1$, $R_2:h_2h_3^{-1}(h_3^{-1}h_2)^2h_3=1$, $R_3:h_2h_3h_2h_3^{-3}h_2=1$:\\
$\Rightarrow\langle h_2,h_3\rangle\cong BS(1,1)$ is solvable, a contradiction.

\item[(6)]$ R_1:h_2^2h_3^{-1}h_2^{-2}h_3=1$, $R_2:h_2h_3^{-1}(h_3^{-1}h_2)^2h_3=1$, $R_3:h_2(h_3h_2^{-1})^2h_3^{-2}h_2=1$:\\
$\Rightarrow\langle h_2,h_3\rangle\cong BS(1,1)$ is solvable, a contradiction.

\item[(7)]$ R_1:h_2h_3h_2h_3^{-1}h_2h_3=1$, $R_2:h_2h_3^{-2}h_2^{-1}h_3^2=1$, $R_3:h_2h_3h_2^{-1}h_3^{-1}h_2h_3^{-1}h_2^{-1}h_3=1$:\\
By interchanging $h_2$ and $h_3$ in (2) and with the same discussion, there is a contradiction.

\item[(8)]$ R_1:h_2h_3h_2h_3^{-1}h_2h_3=1$, $R_2:h_2h_3^{-2}h_2^{-1}h_3^2=1$, $R_3:h_2h_3^{-1}h_2^{-1}h_3^2h_2^{-1}h_3^{-1}h_2=1$:\\
By interchanging $h_2$ and $h_3$ in (1) and with the same discussion, there is a contradiction.

\item[(9)]$ R_1:h_2h_3h_2^{-1}(h_2^{-1}h_3)^2=1$, $R_2:h_2h_3^{-2}h_2^{-1}h_3^2=1$, $R_3:h_2h_3^{-4}h_2h_3=1$:\\
By interchanging $h_2$ and $h_3$ in (3) and with the same discussion, there is a contradiction.

\item[(10)]$ R_1:h_2h_3h_2^{-1}(h_2^{-1}h_3)^2=1$, $R_2:h_2h_3^{-2}h_2^{-1}h_3^2=1$, $R_3:h_2^2h_3^{-1}h_2^{-1}h_3^{-2}h_2=1$:\\
By interchanging $h_2$ and $h_3$ in (5) and with the same discussion, there is a contradiction.

\item[(11)]$ R_1:h_2h_3h_2^{-1}(h_2^{-1}h_3)^2=1$, $R_2:h_2h_3^{-2}h_2^{-1}h_3^2=1$, $R_3:h_2^2h_3^{-2}h_2^{-1}h_3^{-1}h_2=1$:\\
By interchanging $h_2$ and $h_3$ in (4) and with the same discussion, there is a contradiction.

\item[(12)]$ R_1:h_2h_3h_2^{-1}(h_2^{-1}h_3)^2=1$, $R_2:h_2h_3^{-2}h_2^{-1}h_3^2=1$, $R_3:h_2(h_3h_2^{-1})^2h_3^{-2}h_2=1$:\\
By interchanging $h_2$ and $h_3$ in (6) and with the same discussion, there is a contradiction.

\item[(13)]$ R_1:h_2(h_3h_2^{-1}h_3)^2=1$, $R_2:(h_2h_3^{-1})^2(h_2^{-1}h_3)^2=1$, $R_3:h_2h_3^2h_2^{-2}h_3^2=1$:\\
$\Rightarrow\langle h_2,h_3\rangle\cong BS(1,1)$ is solvable, a contradiction.

\item[(14)]$ R_1:h_2(h_3h_2^{-1}h_3)^2=1$, $R_2:(h_2h_3^{-1})^2(h_2^{-1}h_3)^2=1$, $R_3:h_2h_3h_2^{-2}h_3^3=1$:\\
Using Tietze transformation where $h_3\mapsto h_2h_3$ and $h_2\mapsto h_2$, we have:\\
$R_2 \Rightarrow h_3^2\in Z(G)$ where $G=\langle h_2,h_3\rangle$. Let $x=h_2^{-1}h_3^{-1}$. So $h_3^{-1}h_2^{-1}=x^{h_3}$. \\
$R_1 \text{ and } R_3 \Rightarrow (h_3^{-1}h_2^{-1})^2=(h_2^{-1}h_3^{-1})^2 \text{ so if } H=\langle x,x^{h_3}\rangle\Rightarrow H\cong BS(1,-1)$ is solvable. By Corollary \ref{n-sub2} $H \trianglelefteq G$ since $G=\langle x,h_3\rangle$. Since $\frac{G}{H}=\frac{\langle h_3\rangle H}{H}$ is a cyclic group, it is solvable. So $G$ is solvable, a contradiction.

\item[(15)]$ R_1:h_2(h_3h_2^{-1}h_3)^2=1$, $R_2:(h_2h_3^{-1})^2(h_2^{-1}h_3)^2=1$, $R_3:h_2(h_3h_2^{-1})^2h_3^{-2}h_2=1$:\\
$\Rightarrow\langle h_2,h_3\rangle\cong BS(1,1)$ is solvable, a contradiction.

\item[(16)]$ R_1:h_2(h_3h_2^{-1}h_3)^2=1$, $R_2:(h_2h_3^{-1})^2(h_2^{-1}h_3)^2=1$, $R_3:h_2h_3^{-1}h_2^{-1}h_3h_2^{-1}h_3^{-1}h_2h_3=1$:\\
$\Rightarrow h_3^2\in Z(G)$ where $G=\langle h_2,h_3\rangle$. Let $x=h_2h_3^{-1}$. So $h_3^{-1}h_2=x^{h_3}$. Furthermore we have $(h_2h_3^{-1})^2=(h_3^{-1}h_2)^2 \text{ so if } H=\langle x,x^{h_3}\rangle\Rightarrow H\cong BS(1,-1)$ is solvable. By Corollary \ref{n-sub2} $H \trianglelefteq G$ since $G=\langle x,h_3\rangle$. Since $\frac{G}{H}=\frac{\langle h_3\rangle H}{H}$ is a cyclic group, it is solvable. So $G$ is solvable, a contradiction.

\item[(17)]$ R_1:(h_2h_3^{-1}h_2)^2h_3=1$, $R_2:(h_2h_3^{-1})^2(h_2^{-1}h_3)^2=1$, $R_3:h_2^2h_3h_2h_3^{-2}h_2=1$:\\
By interchanging $h_2$ and $h_3$ in (14) and with the same discussion, there is a contradiction.

\item[(18)]$ R_1:(h_2h_3^{-1}h_2)^2h_3=1$, $R_2:(h_2h_3^{-1})^2(h_2^{-1}h_3)^2=1$, $R_3:h_2^2h_3^{-2}h_2^2h_3=1$:\\
By interchanging $h_2$ and $h_3$ in (13) and with the same discussion, there is a contradiction.

\item[(19)]$ R_1:(h_2h_3^{-1}h_2)^2h_3=1$, $R_2:(h_2h_3^{-1})^2(h_2^{-1}h_3)^2=1$, $R_3:h_2h_3h_2^{-1}h_3^{-1}h_2h_3^{-1}h_2^{-1}h_3=1$:\\
By interchanging $h_2$ and $h_3$ in (16) and with the same discussion, there is a contradiction.

\item[(20)]$ R_1:(h_2h_3^{-1}h_2)^2h_3=1$, $R_2:(h_2h_3^{-1})^2(h_2^{-1}h_3)^2=1$, $R_3:h_2(h_3h_2^{-1})^2h_3^{-2}h_2=1$:\\
By interchanging $h_2$ and $h_3$ in (15) and with the same discussion, there is a contradiction.
\end{enumerate}

$\mathbf{C_4-C_6}$ \textbf{subgraph:} Suppose that $[h_1',h_1'',h_2',h_2'',h_3',h_3'',h_4',h_4'']$ is the $8-$tuple related to the cycle $C_4$ and $[h_1',h_1'',h_5',h_5'',h_6',h_6'',h_7',h_7'',h_8',h_8'',h_9',h_9'']$ is the $12-$tuple related to the cycle $C_6$ in the graph $C_4-C_6$, where the first two components of these tuples are related to the common edge of $C_4$ and $C_6$. With the same argument such as about $C_4-C_5$, without loss of generality we may assume that $h_1'=1$, where $h_1'',h_2',h_2'',h_3',h_3'',h_4',h_4'',h_5',h_5'',h_6',h_6'',h_7',h_7'',h_8',h_8'',h_9',h_9'' \in supp(\alpha)$ and $\alpha=1+h_2+h_3$. Also it is easy to see that $h_2' \neq h_5'$ and $h_4'' \neq h_9''$.  With these assumptions and by considering the relations from Tables \ref{tab-C4} and \ref{tab-C6} which are not disproved, it can be seen that there are $2188$ cases for the relations of the cycles $C_4$ and $C_6$ in this structure. Using Gap \cite{gap}, we see that all groups with two generators $h_2$ and $h_3$ and two relations which are between $2035$ cases of these $2188$ cases are finite and solvable and  $14$ groups have the same ``structure description'' $SL(2,5)$ according to the function StructureDescription of GAP, that is finite. So there are just $139$ cases for the relations of the cycles $C_4$ and $C_6$ which may lead to the existence of a subgraph isomorphic to the graph $C_4-C_6$ in the graph $K(\alpha,\beta)$. These cases are listed in table \ref{tab-C4-C6}.  In the following, we show that $127$ cases of these relations lead to a contradiction and just $12$ cases of them may lead to the  existence of a subgraph isomorphic to the graph $C_4-C_6$ in the graph $K(\alpha,\beta)$. Cases which are not disproved are marked by $*$s in the Table \ref{tab-C4-C6}.

\begin{center}
\begin{figure}[t]
\psscalebox{0.9 0.9} 
{
\begin{pspicture}(0,-6.2985578)(16.86,6.2985578)
\psdots[linecolor=black, dotsize=0.3](3.2,6.1014423)
\psdots[linecolor=black, dotsize=0.3](2.4,5.7014422)
\psdots[linecolor=black, dotsize=0.3](2.4,4.9014425)
\psdots[linecolor=black, dotsize=0.3](3.2,4.5014424)
\psdots[linecolor=black, dotsize=0.3](4.0,4.9014425)
\psdots[linecolor=black, dotsize=0.3](4.0,5.7014422)
\psdots[linecolor=black, dotsize=0.3](1.6,4.9014425)
\psdots[linecolor=black, dotsize=0.3](1.6,5.7014422)
\psline[linecolor=black, linewidth=0.04](2.4,5.7014422)(3.2,6.1014423)(4.0,5.7014422)(4.0,4.9014425)(3.2,4.5014424)(2.4,4.9014425)(2.4,5.7014422)(1.6,5.7014422)(1.6,4.9014425)(2.4,4.9014425)
\psdots[linecolor=black, dotsize=0.3](8.4,6.1014423)
\psdots[linecolor=black, dotsize=0.3](9.2,5.7014422)
\psdots[linecolor=black, dotsize=0.3](9.2,4.9014425)
\psdots[linecolor=black, dotsize=0.3](8.4,4.5014424)
\psdots[linecolor=black, dotsize=0.3](7.6,4.9014425)
\psdots[linecolor=black, dotsize=0.3](7.6,5.7014422)
\psdots[linecolor=black, dotsize=0.3](6.8,5.7014422)
\psdots[linecolor=black, dotsize=0.3](6.8,4.9014425)
\psdots[linecolor=black, dotsize=0.3](7.2,3.7014422)
\psdots[linecolor=black, dotsize=0.3](8.8,3.7014422)
\psdots[linecolor=black, dotsize=0.3](14.4,6.1014423)
\psdots[linecolor=black, dotsize=0.3](15.2,5.7014422)
\psdots[linecolor=black, dotsize=0.3](15.2,4.9014425)
\psdots[linecolor=black, dotsize=0.3](14.4,4.5014424)
\psdots[linecolor=black, dotsize=0.3](13.6,4.9014425)
\psdots[linecolor=black, dotsize=0.3](13.6,5.7014422)
\psdots[linecolor=black, dotsize=0.3](12.8,5.7014422)
\psdots[linecolor=black, dotsize=0.3](12.8,4.9014425)
\psdots[linecolor=black, dotsize=0.3](14.8,3.7014422)
\psline[linecolor=black, linewidth=0.04](6.8,5.7014422)(7.6,5.7014422)(8.4,6.1014423)(9.2,5.7014422)(9.2,4.9014425)(8.4,4.5014424)(7.6,4.9014425)(6.8,4.9014425)(6.8,5.7014422)
\psline[linecolor=black, linewidth=0.04](7.6,5.7014422)(7.6,4.9014425)
\psline[linecolor=black, linewidth=0.04](9.2,4.9014425)(8.8,3.7014422)(7.2,3.7014422)(6.8,4.9014425)
\psline[linecolor=black, linewidth=0.04](12.8,5.7014422)(13.6,5.7014422)(14.4,6.1014423)(15.2,5.7014422)(15.2,4.9014425)(14.4,4.5014424)(13.6,4.9014425)(12.8,4.9014425)(12.8,5.7014422)
\psline[linecolor=black, linewidth=0.04](13.6,5.7014422)(13.6,4.9014425)
\psdots[linecolor=black, dotsize=0.3](12.0,3.7014422)
\psline[linecolor=black, linewidth=0.04](15.2,4.9014425)(14.8,3.7014422)(12.0,3.7014422)(12.8,5.7014422)
\psdots[linecolor=black, dotsize=0.3](3.2,1.7014422)
\psdots[linecolor=black, dotsize=0.3](4.0,1.3014423)
\psdots[linecolor=black, dotsize=0.3](4.0,0.5014423)
\psdots[linecolor=black, dotsize=0.3](3.2,0.1014423)
\psdots[linecolor=black, dotsize=0.3](2.4,0.5014423)
\psdots[linecolor=black, dotsize=0.3](2.4,1.3014423)
\psdots[linecolor=black, dotsize=0.3](1.6,1.3014423)
\psdots[linecolor=black, dotsize=0.3](1.6,0.5014423)
\psdots[linecolor=black, dotsize=0.3](4.8,1.3014423)
\psdots[linecolor=black, dotsize=0.3](4.8,0.5014423)
\psdots[linecolor=black, dotsize=0.3](3.2,-1.0985577)
\psdots[linecolor=black, dotsize=0.3](8.4,1.7014422)
\psdots[linecolor=black, dotsize=0.3](9.2,1.3014423)
\psdots[linecolor=black, dotsize=0.3](9.2,0.5014423)
\psdots[linecolor=black, dotsize=0.3](8.4,0.1014423)
\psdots[linecolor=black, dotsize=0.3](7.6,0.5014423)
\psdots[linecolor=black, dotsize=0.3](7.6,1.3014423)
\psdots[linecolor=black, dotsize=0.3](6.8,1.3014423)
\psdots[linecolor=black, dotsize=0.3](6.8,0.5014423)
\psdots[linecolor=black, dotsize=0.3](7.2,-0.6985577)
\psdots[linecolor=black, dotsize=0.3](8.8,-0.6985577)
\psline[linecolor=black, linewidth=0.04](1.6,1.3014423)(2.4,1.3014423)(3.2,1.7014422)(4.0,1.3014423)(4.8,1.3014423)(4.8,0.5014423)(4.0,0.5014423)(3.2,0.1014423)(2.4,0.5014423)(1.6,0.5014423)(1.6,1.3014423)
\psline[linecolor=black, linewidth=0.04](2.4,1.3014423)(2.4,0.5014423)
\psline[linecolor=black, linewidth=0.04](4.0,1.3014423)(4.0,0.5014423)
\psline[linecolor=black, linewidth=0.04](4.8,0.5014423)(3.2,-1.0985577)(1.6,0.5014423)
\psline[linecolor=black, linewidth=0.04](6.8,1.3014423)(7.6,1.3014423)(8.4,1.7014422)(9.2,1.3014423)(9.2,0.5014423)(8.4,0.1014423)(7.6,0.5014423)(6.8,0.5014423)(6.8,1.3014423)
\psline[linecolor=black, linewidth=0.04](7.6,1.3014423)(7.6,0.5014423)
\psline[linecolor=black, linewidth=0.04](9.2,0.5014423)(8.8,-0.6985577)(7.2,-0.6985577)(6.8,0.5014423)(6.8,0.5014423)
\psline[linecolor=black, linewidth=0.04](8.4,1.7014422)(10.4,1.7014422)(8.8,-0.6985577)(8.8,-0.6985577)
\psdots[linecolor=black, dotsize=0.3](10.4,1.7014422)
\psdots[linecolor=black, dotsize=0.3](12.8,1.3014423)
\psdots[linecolor=black, dotsize=0.3](13.6,1.3014423)
\psdots[linecolor=black, dotsize=0.3](14.4,1.3014423)
\psdots[linecolor=black, dotsize=0.3](15.2,1.3014423)
\psdots[linecolor=black, dotsize=0.3](15.2,0.5014423)
\psdots[linecolor=black, dotsize=0.3](14.4,0.5014423)
\psdots[linecolor=black, dotsize=0.3](13.6,0.5014423)
\psdots[linecolor=black, dotsize=0.3](12.8,0.5014423)
\psdots[linecolor=black, dotsize=0.3](13.6,-0.6985577)
\psdots[linecolor=black, dotsize=0.3](14.4,-0.6985577)
\psline[linecolor=black, linewidth=0.04](12.8,1.3014423)(13.6,1.3014423)(14.4,1.3014423)(15.2,1.3014423)(15.2,0.5014423)(14.4,0.5014423)(13.6,0.5014423)(12.8,0.5014423)(12.8,1.3014423)
\psline[linecolor=black, linewidth=0.04](13.6,1.3014423)(13.6,0.5014423)
\psline[linecolor=black, linewidth=0.04](14.4,1.3014423)(14.4,0.5014423)
\psline[linecolor=black, linewidth=0.04](15.2,0.5014423)(14.4,-0.6985577)(13.6,-0.6985577)(12.8,0.5014423)
\psdots[linecolor=black, dotsize=0.3](3.2,-2.6985576)
\psdots[linecolor=black, dotsize=0.3](4.0,-3.0985577)
\psdots[linecolor=black, dotsize=0.3](4.0,-3.8985577)
\psdots[linecolor=black, dotsize=0.3](3.2,-4.2985578)
\psdots[linecolor=black, dotsize=0.3](2.4,-3.8985577)
\psdots[linecolor=black, dotsize=0.3](2.4,-3.0985577)
\psdots[linecolor=black, dotsize=0.3](2.8,-4.698558)
\psdots[linecolor=black, dotsize=0.3](2.4,-5.0985575)
\psdots[linecolor=black, dotsize=0.3](2.0,-5.4985576)
\psdots[linecolor=black, dotsize=0.3](1.6,-3.8985577)
\psdots[linecolor=black, dotsize=0.3](1.6,-3.0985577)
\psdots[linecolor=black, dotsize=0.3](3.6,-5.4985576)
\psline[linecolor=black, linewidth=0.04](1.6,-3.0985577)(2.4,-3.0985577)(3.2,-2.6985576)(4.0,-3.0985577)(4.0,-3.8985577)(3.2,-4.2985578)(2.4,-3.8985577)(2.4,-3.0985577)
\psline[linecolor=black, linewidth=0.04](2.4,-3.8985577)(1.6,-3.8985577)(1.6,-3.0985577)
\psline[linecolor=black, linewidth=0.04](1.6,-3.8985577)(2.0,-5.4985576)(2.4,-5.0985575)(2.8,-4.698558)(3.2,-4.2985578)
\psline[linecolor=black, linewidth=0.04](2.0,-5.4985576)(3.6,-5.4985576)(4.0,-3.8985577)
\psdots[linecolor=black, dotsize=0.3](8.4,-2.6985576)
\psdots[linecolor=black, dotsize=0.3](9.2,-3.0985577)
\psdots[linecolor=black, dotsize=0.3](9.2,-3.8985577)
\psdots[linecolor=black, dotsize=0.3](8.4,-4.2985578)
\psdots[linecolor=black, dotsize=0.3](7.6,-3.8985577)
\psdots[linecolor=black, dotsize=0.3](7.6,-3.0985577)
\psdots[linecolor=black, dotsize=0.3](6.8,-3.0985577)
\psdots[linecolor=black, dotsize=0.3](6.8,-3.8985577)
\psdots[linecolor=black, dotsize=0.3](6.0,-5.4985576)
\psdots[linecolor=black, dotsize=0.3](8.8,-5.4985576)
\psdots[linecolor=black, dotsize=0.3](10.4,-2.6985576)
\psline[linecolor=black, linewidth=0.04](6.8,-3.0985577)(7.6,-3.0985577)(8.4,-2.6985576)(9.2,-3.0985577)(9.2,-3.8985577)(8.4,-4.2985578)(7.6,-3.8985577)(6.8,-3.8985577)(6.8,-3.0985577)
\psline[linecolor=black, linewidth=0.04](7.6,-3.0985577)(7.6,-3.8985577)
\psline[linecolor=black, linewidth=0.04](6.8,-3.0985577)(6.0,-5.4985576)
\psline[linecolor=black, linewidth=0.04](6.0,-5.4985576)(8.8,-5.4985576)(9.2,-3.8985577)
\psline[linecolor=black, linewidth=0.04](8.4,-2.6985576)(10.4,-2.6985576)(8.8,-5.4985576)
\psdots[linecolor=black, dotsize=0.3](14.4,-3.0985577)
\psdots[linecolor=black, dotsize=0.3](15.2,-3.4985578)
\psdots[linecolor=black, dotsize=0.3](15.2,-4.2985578)
\psdots[linecolor=black, dotsize=0.3](14.4,-4.698558)
\psdots[linecolor=black, dotsize=0.3](13.6,-4.2985578)
\psdots[linecolor=black, dotsize=0.3](13.6,-3.4985578)
\psdots[linecolor=black, dotsize=0.3](12.8,-3.4985578)
\psdots[linecolor=black, dotsize=0.3](12.8,-4.2985578)
\psdots[linecolor=black, dotsize=0.3](12.0,-2.6985576)
\psdots[linecolor=black, dotsize=0.3](14.8,-2.6985576)
\psdots[linecolor=black, dotsize=0.3](12.0,-5.4985576)
\psdots[linecolor=black, dotsize=0.3](14.8,-5.4985576)
\psline[linecolor=black, linewidth=0.04](14.4,-3.0985577)(15.2,-3.4985578)(15.2,-4.2985578)(14.4,-4.698558)(13.6,-4.2985578)(12.8,-4.2985578)(12.8,-3.4985578)(13.6,-3.4985578)(14.4,-3.0985577)(14.4,-3.0985577)
\psline[linecolor=black, linewidth=0.04](13.6,-3.4985578)(13.6,-4.2985578)
\psline[linecolor=black, linewidth=0.04](15.2,-3.4985578)(14.8,-2.6985576)(12.0,-2.6985576)(12.8,-4.2985578)
\psline[linecolor=black, linewidth=0.04](12.8,-3.4985578)(12.0,-5.4985576)
\psline[linecolor=black, linewidth=0.04](12.0,-5.4985576)(14.8,-5.4985576)(15.2,-4.2985578)
\rput[bl](2.4,2.9014423){$\mathbf{C_4-C_6}$}
\rput[bl](6.4,2.9014423){$\mathbf{C_4-C_6(-C_6--)}$}
\rput[bl](12.0,2.9014423){$\mathbf{C_4-C_6(--C_7--)}$}
\rput[bl](0.1,-1.8985577){$\mathbf{C_4-C_6(-C_6--)(-C_4-)}$}
\rput[bl](5.6,-1.8985577){$\mathbf{C_4-C_6(-C_6--)(--C_5-)}$}
\rput[bl](11.9,-1.8985577){$\mathbf{C_4-C_6(-C_6--)(---C_4)}$}
\rput[bl](-0.6,-6.2985578){$\mathbf{C_4-C_6(-C_6--)(C_6---)}$}
\rput[bl](5.2,-6.2985578){$\mathbf{C_4-C_6(--C_7--)(--C_5-)}$}
\rput[bl](11.9,-6.2985578){$\mathbf{C_4-C_6(--C_7--)(C_7-1)}$}
\end{pspicture}
}
\vspace{1cc}
\caption{The graphs $C_4-C_6$, $C_4-C_6(-C_6--)$ and $C_4-C_6(--C_7--)$, and some forbidden subgraphs which contain them}\label{f-7}
\end{figure}
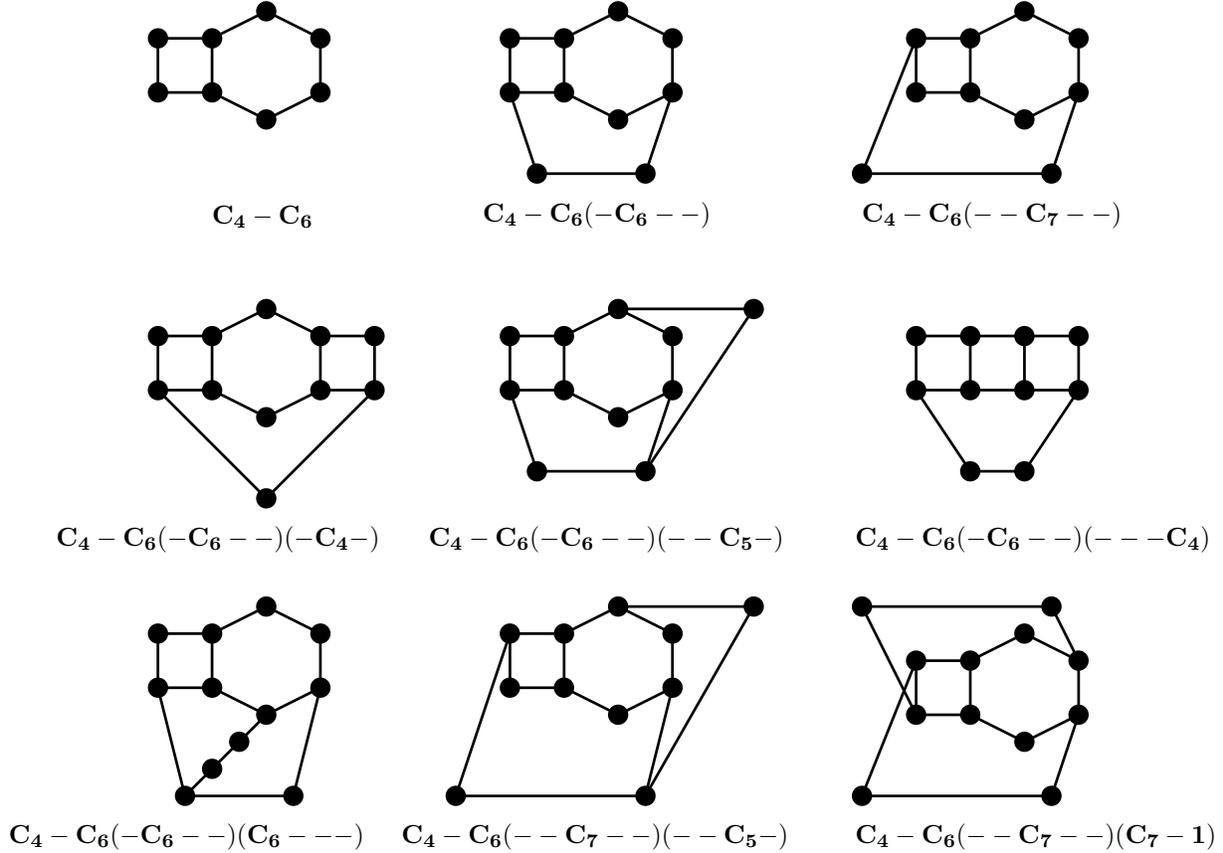
\end{center}
\begin{enumerate}
\item[(1)]$ R_1: h_2^2h_3h_2^{-1}h_3=1$, $R_2: h_2^3h_3h_2^{-2}h_3=1$:\\
$\Rightarrow\langle h_2,h_3\rangle\cong BS(1,1)$ is solvable, a contradiction.

\begin{center}
\small
\begin{longtable}{|c|l|l||c|l|l|}
\caption{{\footnotesize The relations of a $C_4-C_6$ in $K(\alpha,\beta)$}}\label{tab-C4-C6}\\
\hline \multicolumn{1}{|c|}{$n$} & \multicolumn{1}{l|}{$R_1$} & \multicolumn{1}{l||}{$R_2$} & \multicolumn{1}{c|}{$n$} & \multicolumn{1}{l|}{$R_1$} & \multicolumn{1}{l|}{$R_2$}\\\hline
\endfirsthead
\multicolumn{6}{c}%
{{\tablename\ \thetable{} -- continued from previous page}} \\
\hline \multicolumn{1}{|c|}{$n$} & \multicolumn{1}{l|}{$R_1$} & \multicolumn{1}{l||}{$R_2$} & \multicolumn{1}{c|}{$n$} & \multicolumn{1}{l|}{$R_1$} & \multicolumn{1}{l|}{$R_2$}\\\hline
\endhead
\hline \multicolumn{6}{|r|}{{Continued on next page}} \\\hline
\endfoot
\hline
\endlastfoot
$1$&$h_2^2h_3h_2^{-1}h_3=1$&$h_2^3h_3h_2^{-2}h_3=1$&$
40$&$h_2h_3h_2^{-2}h_3=1$&$h_2h_3^{-1}h_2^{-1}h_3^2h_2^{-1}h_3^{-1}h_2=1$\\
$2$&$h_2^2h_3h_2^{-1}h_3=1$&$h_2^2h_3^{-1}h_2^{-1}h_3^3=1$&$
41$&$h_2h_3h_2^{-2}h_3=1$&$h_2h_3^{-1}h_2^{-1}h_3h_2^{-1}(h_3^{-1}h_2)^2=1$\\
$3$&$h_2^2h_3h_2^{-1}h_3=1$&$h_2^2h_3^{-2}(h_2^{-1}h_3)^2=1$&$
42$&$h_2h_3h_2^{-2}h_3=1$&$h_2h_3^{-2}h_2(h_2h_3^{-1})^2h_2=1$\\
$4$&$h_2^2h_3h_2^{-1}h_3=1$&$h_2^2h_3^{-1}h_2(h_3h_2^{-1})^2h_3=1$&$
43$&$h_2h_3h_2^{-2}h_3=1$&$h_2h_3^{-1}h_2(h_3h_2^{-1})^2h_2^{-1}h_3=1$\\
$5$&$h_2^2h_3h_2^{-1}h_3=1$&$h_2(h_2h_3^{-1})^2h_2h_3h_2^{-1}h_3=1$&$
44$&$h_2h_3h_2^{-2}h_3=1$&$h_2h_3^2h_2h_3h_2^{-1}h_3=1$\\
$6$&$h_2^2h_3h_2^{-1}h_3=1$&$h_2h_3^{-1}h_2^{-1}h_3^2h_2^{-1}h_3^{-1}h_2=1$&$
45$&$h_2h_3h_2^{-2}h_3=1$&$h_2h_3h_2^{-1}h_3^{-1}h_2^{-1}h_3^2=1$\\
$7$&$h_2^2h_3h_2^{-1}h_3=1$&$h_2^2h_3^4=1$&$46$&$h_2h_3h_2^{-2}h_3=1$&$h_2h_3h_2^{-1}h_3^{-3}h_2=1$\\
$8$&$h_2^2h_3h_2^{-1}h_3=1$&$h_2h_3h_2h_3^{-1}h_2^{-1}h_3^2=1$&$
47$&$h_2h_3h_2^{-2}h_3=1$&$h_2(h_3h_2^{-1})^2h_3^{-2}h_2=1$\\
$9$&$h_2^2h_3h_2^{-1}h_3=1$&$h_2h_3^3h_2^{-1}h_3^{-1}h_2=1$&$
48 *$&$h_2h_3h_2^{-2}h_3=1$&$h_3^2(h_3h_2^{-1})^3h_3=1$\\
$10$&$h_2^2h_3h_2^{-1}h_3=1$&$h_2(h_3h_2^{-1})^2h_3^{-2}h_2=1$&$
49 *$&$h_2h_3h_2^{-2}h_3=1$&$h_3^2h_2^{-1}(h_3h_2^{-1}h_3)^2=1$\\
$11$&$h_2^2h_3h_2^{-1}h_3=1$&$h_2h_3^2h_2^{-1}h_3^{-1}h_2h_3=1$&$
50$&$h_2h_3h_2^{-2}h_3=1$&$h_2h_3^2h_2^{-1}h_3^{-1}h_2^{-1}h_3=1$\\
$12$&$h_2^2h_3h_2^{-1}h_3=1$&$h_2h_3h_2^{-1}h_3^{-1}h_2h_3^{-1}h_2^{-1}h_3=1$&$
51$&$h_2h_3h_2^{-2}h_3=1$&$h_2h_3h_2^{-1}h_3^{-1}h_2h_3^{-1}h_2^{-1}h_3=1$\\
$13$&$h_2^2h_3h_2^{-1}h_3=1$&$h_2h_3h_2^{-1}h_3(h_2h_3^{-1})^2h_2=1$&$
52$&$h_2(h_3h_2^{-1})^2h_3=1$&$h_2^2h_3h_2^{-1}(h_2^{-1}h_3)^2=1$\\
$14$&$h_2^2h_3h_2^{-1}h_3=1$&$h_2(h_3h_2^{-1})^2h_3h_2h_3^{-1}h_2=1$&$
53$&$h_2(h_3h_2^{-1})^2h_3=1$&$h_2^2(h_3h_2^{-1})^2h_2^{-1}h_3=1$\\
$15$&$h_2^2h_3h_2^{-1}h_3=1$&$h_2h_3^{-1}h_2^{-1}h_3h_2^{-1}h_3^{-1}h_2h_3=1$&$
54$&$h_2(h_3h_2^{-1})^2h_3=1$&$h_2^2h_3^{-1}h_2^{-1}(h_3^{-1}h_2)^2=1$\\
$16$&$h_2^2h_3^{-2}h_2=1$&$h_2h_3h_2h_3^{-3}h_2=1$&$55$&$h_2(h_3h_2^{-1})^2h_3=1$&$h_2^2h_3^{-2}(h_2^{-1}h_3)^2=1$\\
$17$&$h_2^2h_3^{-2}h_2=1$&$h_2h_3^2h_2h_3^{-2}h_2=1$&$
56$&$h_2(h_3h_2^{-1})^2h_3=1$&$h_2(h_2h_3^{-1})^2h_2^{-1}h_3^{-1}h_2=1$\\
$18$&$h_2^2h_3^{-2}h_2=1$&$h_2h_3h_2^{-1}h_3^{-1}h_2h_3^{-1}h_2^{-1}h_3=1$&$
57$&$h_2(h_3h_2^{-1})^2h_3=1$&$h_2h_3^{-1}h_2^{-1}h_3^2h_2^{-1}h_3^{-1}h_2=1$\\
$19$&$h_2^2h_3^{-2}h_2=1$&$h_2h_3h_2^{-1}(h_3^{-1}h_2)^3=1$&$
58$&$h_2(h_3h_2^{-1})^2h_3=1$&$h_2(h_3h_2^{-1})^2h_3^{-2}h_2=1$\\
$20$&$h_2^2h_3^{-2}h_2=1$&$h_2(h_3h_2^{-1})^2h_3^{-2}h_2=1$&$
59$&$h_2(h_3h_2^{-1})^2h_3=1$&$h_2h_3^{-1}(h_2^{-1}h_3^2)^2=1$\\
$21$&$h_2^2h_3^{-2}h_2=1$&$h_2h_3^{-1}h_2^{-1}h_3^2h_2^{-1}h_3^{-1}h_2=1$&$
60$&$h_2(h_3h_2^{-1})^2h_3=1$&$h_2h_3^{-1}(h_2^{-1}h_3)^2h_3^2=1$\\
$22$&$h_2^2h_3^{-2}h_2=1$&$h_2h_3^{-1}(h_2^{-1}h_3h_2^{-1})^2h_3=1$&$
61$&$h_2(h_3h_2^{-1})^2h_3=1$&$(h_2h_3^{-1})^2h_3^{-2}h_2^{-1}h_3=1$\\
$23$&$h_2^2h_3^{-2}h_2=1$&$h_2h_3^{-1}h_2^{-1}h_3h_2^{-1}h_3^{-1}h_2h_3=1$&$
62$&$h_2(h_3h_2^{-1})^2h_3=1$&$h_3^3(h_3h_2^{-1})^2h_3=1$\\
$24$&$h_2^2h_3^{-2}h_2=1$&$(h_2h_3^{-1}h_2)^2h_3^{-1}h_2^{-1}h_3=1$&$
63$&$h_2h_3^{-3}h_2=1$&$h_2^2h_3^2h_2^{-2}h_3=1$\\
$25$&$h_2^2h_3^{-2}h_2=1$&$h_2^3h_3^{-1}h_2^{-1}h_3^{-1}h_2=1$&$64$&$h_2h_3^{-3}h_2=1$&$h_2^2h_3h_2^{-2}h_3^2=1$\\
$26$&$h_2^2h_3^{-2}h_2=1$&$h_2^2h_3^{-2}(h_2^{-1}h_3)^2=1$&$65$&$h_2h_3^{-3}h_2=1$&$h_2^2h_3^{-1}h_2^{-1}h_3^{-2}h_2=1$\\
$27$&$h_2^2h_3^{-2}h_2=1$&$h_2^2h_3^{-3}h_2h_3=1$&$66$&$h_2h_3^{-3}h_2=1$&$h_2^2h_3^{-2}h_2^{-1}h_3^{-1}h_2=1$\\
$28$&$h_2^2h_3^{-2}h_2=1$&$h_2^2h_3^{-2}h_2h_3^2=1$&$67$&$h_2h_3^{-3}h_2=1$&$h_2^2h_3^{-2}(h_2^{-1}h_3)^2=1$\\
$29$&$h_2^2h_3^{-2}h_2=1$&$h_2(h_2h_3^{-1})^3h_2^{-1}h_3=1$&$
68$&$h_2h_3^{-3}h_2=1$&$h_2h_3h_2^{-1}h_3^{-1}h_2h_3^{-1}h_2^{-1}h_3=1$\\
$30$&$h_2h_3h_2^{-2}h_3=1$&$h_2^2h_3h_2^{-3}h_3=1$&$69$&$h_2h_3^{-3}h_2=1$&$h_2(h_3h_2^{-1})^2h_3^{-2}h_2=1$\\
$31$&$h_2h_3h_2^{-2}h_3=1$&$h_2^2h_3^{-1}(h_2^{-1}h_3)^3=1$&$
70$&$h_2h_3^{-3}h_2=1$&$h_2h_3^{-1}h_2^{-1}h_3^2h_2^{-1}h_3^{-1}h_2=1$\\
$32$&$h_2h_3h_2^{-2}h_3=1$&$h_2^2h_3^{-2}(h_2^{-1}h_3)^2=1$&$
71$&$h_2h_3^{-3}h_2=1$&$h_2h_3^{-1}h_2^{-1}h_3h_2^{-1}h_3^{-1}h_2h_3=1$\\
$33$&$h_2h_3h_2^{-2}h_3=1$&$h_2^2h_3^{-3}h_2^{-1}h_3=1$&$
72$&$h_2h_3^{-3}h_2=1$&$h_2h_3^{-1}h_2^{-1}(h_3h_2^{-1}h_3)^2=1$\\
$34 *$&$h_2h_3h_2^{-2}h_3=1$&$h_2h_3^2h_2^{-1}(h_2^{-1}h_3)^2=1$&$
73$&$h_2h_3^{-3}h_2=1$&$h_2h_3^{-1}(h_2^{-1}h_3)^3h_3=1$\\
$35$&$h_2h_3h_2^{-2}h_3=1$&$h_2h_3(h_3h_2^{-1})^2h_2^{-1}h_3=1$&$
74$&$h_2h_3^{-3}h_2=1$&$h_2(h_3^{-1}h_2h_3^{-1})^2h_2^{-1}h_3=1$\\
$36$&$h_2h_3h_2^{-2}h_3=1$&$h_2h_3h_2^{-2}h_3^2h_2^{-1}h_3=1$&$75$&$h_2h_3^{-3}h_2=1$&$h_2h_3^{-4}h_2h_3=1$\\
$37$&$h_2h_3h_2^{-2}h_3=1$&$h_2h_3h_2^{-1}(h_2^{-1}h_3)^2h_3=1$&$
76$&$h_2h_3^{-3}h_2=1$&$(h_2h_3^{-1})^3h_3^{-1}h_2^{-1}h_3=1$\\
$38$&$h_2h_3h_2^{-2}h_3=1$&$h_2h_3h_2^{-1}h_3^2h_2^{-2}h_3=1$&$
77$&$h_2h_3^{-2}h_2h_3=1$&$h_2^3h_3^{-2}h_2^{-1}h_3=1$\\
$39 *$&$h_2h_3h_2^{-2}h_3=1$&$h_2(h_3h_2^{-1})^2h_2^{-1}h_3^2=1$&$
78 *$&$h_2h_3^{-2}h_2h_3=1$&$h_2^2(h_2h_3^{-1})^3h_2=1$\\
$79$&$h_2h_3^{-2}h_2h_3=1$&$h_2^2h_3h_2h_3^{-1}h_2h_3=1$&$
110$&$h_2h_3^{-1}h_2h_3^2=1$&$h_2^2h_3^{-2}(h_2^{-1}h_3)^2=1$\\
$80$&$h_2h_3^{-2}h_2h_3=1$&$h_2^2h_3h_2^{-1}h_3^{-2}h_2=1$&$111$&$h_2h_3^{-1}h_2h_3^2=1$&$h_2h_3^{-2}h_2h_3^3=1$\\
$81 *$&$h_2h_3^{-2}h_2h_3=1$&$h_2^2h_3^{-1}(h_2h_3^{-1}h_2)^2=1$&$
112$&$h_2h_3^{-1}h_2h_3^2=1$&$h_2h_3^{-1}h_2h_3(h_3h_2^{-1})^2h_3=1$\\
$82$&$h_2h_3^{-2}h_2h_3=1$&$h_2h_3h_2h_3^{-1}h_2^{-1}h_3^{-1}h_2=1$&$
113$&$h_2h_3^{-1}h_2h_3^2=1$&$(h_2h_3^{-1})^2h_2h_3^2h_2^{-1}h_3=1$\\
$83$&$h_2h_3^{-2}h_2h_3=1$&$h_2h_3h_2h_3^{-1}(h_3^{-1}h_2)^2=1$&$
114 *$&$h_2h_3^{-1}h_2h_3h_2^{-1}h_3=1$&$h_2^3h_3^3=1$\\
$84 *$&$h_2h_3^{-2}h_2h_3=1$&$h_2h_3(h_2h_3^{-1})^2h_3^{-1}h_2=1$&$
115$&$h_2h_3^{-1}h_2h_3h_2^{-1}h_3=1$&$h_2^3h_3^{-1}(h_2^{-1}h_3)^2=1$\\
$85$&$h_2h_3^{-2}h_2h_3=1$&$h_2(h_3h_2^{-1})^2h_3^{-2}h_2=1$&$
116 *$&$h_2h_3^{-1}h_2h_3h_2^{-1}h_3=1$&$h_2(h_2h_3)^2h_3=1$\\
$86$&$h_2h_3^{-2}h_2h_3=1$&$h_2^2h_3^{-2}(h_2^{-1}h_3)^2=1$&$
117 *$&$h_2h_3^{-1}h_2h_3h_2^{-1}h_3=1$&$h_2^2h_3^2h_2h_3=1$\\
$87 *$&$h_2h_3^{-2}h_2h_3=1$&$h_2^2h_3^{-1}(h_3^{-1}h_2)^2h_3=1$&$
118$&$h_2h_3^{-1}h_2h_3h_2^{-1}h_3=1$&$h_2^2(h_3h_2^{-1})^2h_3^{-1}h_2=1$\\
$88$&$h_2h_3^{-2}h_2h_3=1$&$h_2(h_2h_3^{-1})^2h_3^{-1}h_2h_3=1$&$
119$&$h_2h_3^{-1}h_2h_3h_2^{-1}h_3=1$&$h_2^2h_3^{-1}h_2h_3h_2h_3^{-1}h_2=1$\\
$89$&$h_2h_3^{-2}h_2h_3=1$&$h_2h_3^{-1}h_2^{-1}h_3^2h_2^{-1}h_3^{-1}h_2=1$&$
120$&$h_2h_3^{-1}h_2h_3h_2^{-1}h_3=1$&$h_2h_3^2h_2^{-1}h_3h_2h_3^{-1}h_2=1$\\
$90$&$h_2h_3^{-2}h_2h_3=1$&$h_2h_3^{-2}(h_2^{-1}h_3)^3=1$&$
121$&$h_2h_3^{-1}h_2h_3h_2^{-1}h_3=1$&$h_2h_3h_2^{-2}h_3h_2h_3^{-1}h_2=1$\\
$91$&$h_2h_3^{-2}h_2h_3=1$&$h_2h_3^{-3}h_2h_3^2=1$&$
122 *$&$h_2h_3^{-1}h_2h_3h_2^{-1}h_3=1$&$h_2h_3h_2^{-1}h_3^2h_2h_3^{-1}h_2=1$\\
$92$&$h_2h_3^{-2}h_2h_3=1$&$h_2h_3^{-2}h_2^2h_3^{-1}h_2h_3=1$&$
123$&$h_2h_3^{-1}h_2h_3h_2^{-1}h_3=1$&$h_2(h_3h_2^{-1})^2h_3^{-2}h_2=1$\\
$93$&$h_2h_3^{-2}h_2h_3=1$&$h_2h_3^{-2}h_2h_3h_2h_3^{-1}h_2=1$&$
124$&$h_2h_3^{-1}h_2h_3h_2^{-1}h_3=1$&$h_2^2h_3^{-1}h_2h_3h_2^{-1}h_3^2=1$\\
$94$&$h_2h_3^{-2}h_2h_3=1$&$h_2h_3^{-1}(h_3^{-1}h_2)^2h_3h_2^{-1}h_3=1$&$
125$&$h_2h_3^{-1}h_2h_3h_2^{-1}h_3=1$&$h_2h_3h_2^{-1}h_3^3h_2^{-1}h_3=1$\\
$95$&$h_2h_3^{-2}h_2h_3=1$&$h_2h_3^{-1}(h_3^{-1}h_2)^2h_3^{-2}h_2=1$&$
126$&$h_2h_3^{-1}h_2h_3h_2^{-1}h_3=1$&$h_2h_3^{-3}(h_2^{-1}h_3)^2=1$\\
$96$&$h_2h_3^{-2}h_2h_3=1$&$(h_2h_3^{-1})^2h_3^{-1}h_2h_3h_2^{-1}h_3=1$&$
127$&$h_2h_3^{-1}h_2h_3h_2^{-1}h_3=1$&$h_2h_3^{-2}h_2h_3h_2^{-1}h_3^2=1$\\
$97$&$h_2h_3^{-2}h_2h_3=1$&$h_2^2h_3^{-1}h_2^{-1}h_3^{-1}h_2h_3=1$&$
128$&$h_2h_3^{-1}h_2h_3h_2^{-1}h_3=1$&$(h_2h_3^{-1})^2h_2^{-1}h_3^3=1$\\
$98$&$h_2h_3^{-2}h_2h_3=1$&$h_2h_3^{-1}h_2^{-1}h_3h_2^{-1}h_3^{-1}h_2h_3=1$&$
129$&$(h_2h_3^{-1})^2h_2h_3=1$&$h_2^3(h_2h_3^{-1})^2h_2=1$\\
$99$&$h_2h_3^{-1}h_2h_3^2=1$&$h_2^4h_3^2=1$&$130$&$(h_2h_3^{-1})^2h_2h_3=1$&$h_2^2(h_2h_3^{-1})^2h_2^{-1}h_3=1$\\
$100$&$h_2h_3^{-1}h_2h_3^2=1$&$h_2^3h_3^{-1}h_2^{-1}h_3^2=1$&$
131$&$(h_2h_3^{-1})^2h_2h_3=1$&$h_2^2h_3h_2^{-1}(h_3^{-1}h_2)^2=1$\\
$101$&$h_2h_3^{-1}h_2h_3^2=1$&$h_2^2h_3^2h_2^{-1}h_3^{-1}h_2=1$&$
132$&$(h_2h_3^{-1})^2h_2h_3=1$&$h_2h_3h_2^{-1}h_3^{-1}h_2^2h_3^{-1}h_2=1$\\
$102$&$h_2h_3^{-1}h_2h_3^2=1$&$(h_2h_3)^2h_2^{-1}h_3^{-1}h_2=1$&$
133$&$(h_2h_3^{-1})^2h_2h_3=1$&$h_2(h_3h_2^{-1})^2h_3^{-2}h_2=1$\\
$103$&$h_2h_3^{-1}h_2h_3^2=1$&$h_2(h_3h_2^{-1})^2h_3^{-2}h_2=1$&$
134$&$(h_2h_3^{-1})^2h_2h_3=1$&$h_2^2h_3^{-2}(h_2^{-1}h_3)^2=1$\\
$104$&$h_2h_3^{-1}h_2h_3^2=1$&$h_2^2h_3^{-1}h_2^{-1}h_3h_2h_3=1$&$
135$&$(h_2h_3^{-1})^2h_2h_3=1$&$h_2h_3^{-1}h_2^{-1}h_3^2h_2^{-1}h_3^{-1}h_2=1$\\
$105$&$h_2h_3^{-1}h_2h_3^2=1$&$h_2h_3h_2^{-1}h_3^{-1}h_2h_3^{-1}h_2^{-1}h_3=1$&$
136$&$(h_2h_3^{-1})^2h_2h_3=1$&$h_2h_3^{-2}(h_3^{-1}h_2)^2h_3=1$\\
$106$&$h_2h_3^{-1}h_2h_3^2=1$&$h_2h_3^{-1}h_2^{-1}h_3^2h_2^{-1}h_3^{-1}h_2=1$&$
137$&$(h_2h_3^{-1})^2h_2h_3=1$&$h_2h_3^{-1}(h_3^{-1}h_2)^2h_3^2=1$\\
$107$&$h_2h_3^{-1}h_2h_3^2=1$&$h_2h_3^{-1}h_2^{-1}h_3h_2^{-1}h_3^{-1}h_2h_3=1$&$
138$&$(h_2h_3^{-1})^2h_2h_3=1$&$(h_2h_3^{-1})^2h_3^{-2}h_2h_3=1$\\
$108$&$h_2h_3^{-1}h_2h_3^2=1$&$h_2h_3^{-1}h_2(h_3h_2^{-1})^2h_3^2=1$&$
139$&$(h_2h_3^{-1})^2h_2h_3=1$&$(h_2h_3^{-1})^2h_3^{-1}h_2h_3^2=1$\\
$109$&$h_2h_3^{-1}h_2h_3^2=1$&$(h_2h_3^{-1})^2h_2h_3h_2^{-1}h_3^2=1$&$$&$$&$ $\\
\hline
\end{longtable}
\end{center}

\item[(2)]$ R_1: h_2^2h_3h_2^{-1}h_3=1$, $R_2: h_2^2h_3^{-1}h_2^{-1}h_3^3=1$:\\
$\Rightarrow\langle h_2,h_3\rangle\cong BS(1,-1)$ is solvable, a contradiction.

\item[(3)]$ R_1: h_2^2h_3h_2^{-1}h_3=1$, $R_2: h_2^2h_3^{-2}(h_2^{-1}h_3)^2=1$:\\
$\Rightarrow \langle h_2,h_3\rangle=\langle h_3\rangle$ is abelian, a contradiction.

\item[(4)]$ R_1: h_2^2h_3h_2^{-1}h_3=1$, $R_2: h_2^2h_3^{-1}h_2(h_3h_2^{-1})^2h_3=1$:\\
$\Rightarrow\langle h_2,h_3\rangle\cong BS(1,1)$ is solvable, a contradiction.

\item[(5)]$ R_1: h_2^2h_3h_2^{-1}h_3=1$, $R_2: h_2(h_2h_3^{-1})^2h_2h_3h_2^{-1}h_3=1$:\\
$\Rightarrow \langle h_2,h_3\rangle=\langle h_2\rangle$ is abelian, a contradiction.

\item[(6)]$ R_1: h_2^2h_3h_2^{-1}h_3=1$, $R_2: h_2h_3^{-1}h_2^{-1}h_3^2h_2^{-1}h_3^{-1}h_2=1$:\\
$\Rightarrow \langle h_2,h_3\rangle=\langle h_3\rangle$ is abelian, a contradiction.

\item[(7)]$ R_1: h_2^2h_3h_2^{-1}h_3=1$, $R_2: h_2^2h_3^4=1$:\\
$\Rightarrow \langle h_2,h_3\rangle=\langle h_3\rangle$ is abelian, a contradiction.

\item[(8)]$ R_1: h_2^2h_3h_2^{-1}h_3=1$, $R_2: h_2h_3h_2h_3^{-1}h_2^{-1}h_3^2=1$:\\
$\Rightarrow \langle h_2,h_3\rangle=\langle h_3\rangle$ is abelian, a contradiction.

\item[(9)]$ R_1: h_2^2h_3h_2^{-1}h_3=1$, $R_2: h_2h_3^3h_2^{-1}h_3^{-1}h_2=1$:\\
$\Rightarrow\langle h_2,h_3\rangle\cong BS(1,-1)$ is solvable, a contradiction.

\item[(10)]$ R_1: h_2^2h_3h_2^{-1}h_3=1$, $R_2: h_2(h_3h_2^{-1})^2h_3^{-2}h_2=1$:\\
$\Rightarrow \langle h_2,h_3\rangle=\langle h_3\rangle$ is abelian, a contradiction.

\item[(11)]$ R_1: h_2^2h_3h_2^{-1}h_3=1$, $R_2: h_2h_3^2h_2^{-1}h_3^{-1}h_2h_3=1$:\\
$\Rightarrow \langle h_2,h_3\rangle=\langle h_3\rangle$ is abelian, a contradiction.

\item[(12)]$ R_1: h_2^2h_3h_2^{-1}h_3=1$, $R_2: h_2h_3h_2^{-1}h_3^{-1}h_2h_3^{-1}h_2^{-1}h_3=1$:\\
$\Rightarrow \langle h_2,h_3\rangle=\langle h_3\rangle$ is abelian, a contradiction.

\item[(13)]$ R_1: h_2^2h_3h_2^{-1}h_3=1$, $R_2: h_2h_3h_2^{-1}h_3(h_2h_3^{-1})^2h_2=1$:\\
$\Rightarrow \langle h_2,h_3\rangle=\langle h_2\rangle$ is abelian, a contradiction.

\item[(14)]$ R_1: h_2^2h_3h_2^{-1}h_3=1$, $R_2: h_2(h_3h_2^{-1})^2h_3h_2h_3^{-1}h_2=1$:\\
$\Rightarrow\langle h_2,h_3\rangle\cong BS(1,1)$ is solvable, a contradiction.

\item[(15)]$ R_1: h_2^2h_3h_2^{-1}h_3=1$, $R_2: h_2h_3^{-1}h_2^{-1}h_3h_2^{-1}h_3^{-1}h_2h_3=1$:\\
$\Rightarrow\langle h_2,h_3\rangle\cong BS(1,-1)$ is solvable, a contradiction.

\item[(16)]$ R_1: h_2^2h_3^{-2}h_2=1$, $R_2: h_2h_3h_2h_3^{-3}h_2=1$:\\
$\Rightarrow\langle h_2,h_3\rangle\cong BS(1,1)$ is solvable, a contradiction.

\item[(17)]$ R_1: h_2^2h_3^{-2}h_2=1$, $R_2: h_2h_3^2h_2h_3^{-2}h_2=1$:\\
$\Rightarrow h_2=1$ that is a contradiction.

\item[(18)]$ R_1: h_2^2h_3^{-2}h_2=1$, $R_2: h_2h_3h_2^{-1}h_3^{-1}h_2h_3^{-1}h_2^{-1}h_3=1$:\\
$\Rightarrow \langle h_2,h_3\rangle=\langle h_3\rangle$ is abelian, a contradiction.

\item[(19)]$ R_1: h_2^2h_3^{-2}h_2=1$, $R_2: h_2h_3h_2^{-1}(h_3^{-1}h_2)^3=1$:\\
$\Rightarrow h_3^2\in Z(G)$ where $G=\langle h_2,h_3\rangle$. Let $x=h_2h_3^{-1}$. So $h_3^{-1}h_2=x^{h_3}$. Also $(h_2h_3^{-1})^2=(h_3^{-1}h_2)^2 \text{ so if } H=\langle x,x^{h_3}\rangle\Rightarrow H\cong BS(1,-1)$ is solvable. By Corollary \ref{n-sub2} $H \trianglelefteq G$ since $G=\langle x,h_3\rangle$. Since $\frac{G}{H}=\frac{\langle h_3\rangle H}{H}$ is a cyclic group, it is solvable. So $G$ is solvable, a contradiction.

\item[(20)]$ R_1: h_2^2h_3^{-2}h_2=1$, $R_2: h_2(h_3h_2^{-1})^2h_3^{-2}h_2=1$:\\
$\Rightarrow \langle h_2,h_3\rangle=\langle h_3h_2^{-1}\rangle$ is abelian, a contradiction.

\item[(21)]$ R_1: h_2^2h_3^{-2}h_2=1$, $R_2: h_2h_3^{-1}h_2^{-1}h_3^2h_2^{-1}h_3^{-1}h_2=1$:\\
$\Rightarrow \langle h_2,h_3\rangle=\langle h_3h_2^{-1}\rangle$ is abelian, a contradiction.

\item[(22)]$ R_1: h_2^2h_3^{-2}h_2=1$, $R_2: h_2h_3^{-1}(h_2^{-1}h_3h_2^{-1})^2h_3=1$:\\
$\Rightarrow\langle h_2,h_3\rangle\cong BS(1,1)$ is solvable, a contradiction.

\item[(23)]$ R_1: h_2^2h_3^{-2}h_2=1$, $R_2: h_2h_3^{-1}h_2^{-1}h_3h_2^{-1}h_3^{-1}h_2h_3=1$:\\
With the same discussion such as the item 19, there is a contradiction.

\item[(24)]$ R_1: h_2^2h_3^{-2}h_2=1$, $R_2: (h_2h_3^{-1}h_2)^2h_3^{-1}h_2^{-1}h_3=1$:\\
$\Rightarrow\langle h_2,h_3\rangle\cong BS(1,1)$ is solvable, a contradiction.

\item[(25)]$ R_1: h_2^2h_3^{-2}h_2=1$, $R_2: h_2^3h_3^{-1}h_2^{-1}h_3^{-1}h_2=1$:\\
$\Rightarrow\langle h_2,h_3\rangle\cong BS(1,1)$ is solvable, a contradiction.

\item[(26)]$ R_1: h_2^2h_3^{-2}h_2=1$, $R_2: h_2^2h_3^{-2}(h_2^{-1}h_3)^2=1$:\\
$\Rightarrow \langle h_2,h_3\rangle=\langle h_2^{-1}h_3\rangle$ is abelian, a contradiction.

\item[(27)]$ R_1: h_2^2h_3^{-2}h_2=1$, $R_2: h_2^2h_3^{-3}h_2h_3=1$:\\
$\Rightarrow\langle h_2,h_3\rangle\cong BS(1,1)$ is solvable, a contradiction.

\item[(28)]$ R_1: h_2^2h_3^{-2}h_2=1$, $R_2: h_2^2h_3^{-2}h_2h_3^2=1$:\\
$\Rightarrow h_3=1$ that is a contradiction.

\item[(29)]$ R_1: h_2^2h_3^{-2}h_2=1$, $R_2: h_2(h_2h_3^{-1})^3h_2^{-1}h_3=1$:\\
With the same discussion such as the item 19, there is a contradiction.

\item[(30)]$ R_1: h_2h_3h_2^{-2}h_3=1$, $R_2: h_2^2h_3h_2^{-3}h_3=1$:\\
$\Rightarrow\langle h_2,h_3\rangle\cong BS(1,1)$ is solvable, a contradiction.

\item[(31)]$ R_1: h_2h_3h_2^{-2}h_3=1$, $R_2: h_2^2h_3^{-1}(h_2^{-1}h_3)^3=1$:\\
$\Rightarrow\langle h_2,h_3\rangle\cong BS(1,-1)$ is solvable, a contradiction.

\item[(32)]$ R_1: h_2h_3h_2^{-2}h_3=1$, $R_2: h_2^2h_3^{-2}(h_2^{-1}h_3)^2=1$:\\
$\Rightarrow \langle h_2,h_3\rangle=\langle h_3\rangle$ is abelian, a contradiction.

\item[(33)]$ R_1: h_2h_3h_2^{-2}h_3=1$, $R_2: h_2^2h_3^{-3}h_2^{-1}h_3=1$:\\
$\Rightarrow\langle h_2,h_3\rangle\cong BS(1,-1)$ is solvable, a contradiction.


\item[(35)]$ R_1: h_2h_3h_2^{-2}h_3=1$, $R_2: h_2h_3(h_3h_2^{-1})^2h_2^{-1}h_3=1$:\\
$\Rightarrow \langle h_2,h_3\rangle=\langle h_3\rangle$ is abelian, a contradiction.

\item[(36)]$ R_1: h_2h_3h_2^{-2}h_3=1$, $R_2: h_2h_3h_2^{-2}h_3^2h_2^{-1}h_3=1$:\\
$\Rightarrow \langle h_2,h_3\rangle=\langle h_3\rangle$ is abelian, a contradiction.

\item[(37)]$ R_1: h_2h_3h_2^{-2}h_3=1$, $R_2: h_2h_3h_2^{-1}(h_2^{-1}h_3)^2h_3=1$:\\
$\Rightarrow \langle h_2,h_3\rangle=\langle h_3\rangle$ is abelian, a contradiction.

\item[(38)]$ R_1: h_2h_3h_2^{-2}h_3=1$, $R_2: h_2h_3h_2^{-1}h_3^2h_2^{-2}h_3=1$:\\
$\Rightarrow \langle h_2,h_3\rangle=\langle h_3\rangle$ is abelian, a contradiction.


\item[(40)]$ R_1: h_2h_3h_2^{-2}h_3=1$, $R_2: h_2h_3^{-1}h_2^{-1}h_3^2h_2^{-1}h_3^{-1}h_2=1$:\\
$\Rightarrow \langle h_2,h_3\rangle=\langle h_3\rangle$ is abelian, a contradiction.

\item[(41)]$ R_1: h_2h_3h_2^{-2}h_3=1$, $R_2: h_2h_3^{-1}h_2^{-1}h_3h_2^{-1}(h_3^{-1}h_2)^2=1$:\\
$\Rightarrow\langle h_2,h_3\rangle\cong BS(1,1)$ is solvable, a contradiction.

\item[(42)]$ R_1: h_2h_3h_2^{-2}h_3=1$, $R_2: h_2h_3^{-2}h_2(h_2h_3^{-1})^2h_2=1$:\\
$\Rightarrow h_2=1$ that is a contradiction.

\item[(43)]$ R_1: h_2h_3h_2^{-2}h_3=1$, $R_2: h_2h_3^{-1}h_2(h_3h_2^{-1})^2h_2^{-1}h_3=1$:\\
$\Rightarrow\langle h_2,h_3\rangle\cong BS(1,1)$ is solvable, a contradiction.

\item[(44)]$ R_1: h_2h_3h_2^{-2}h_3=1$, $R_2: h_2h_3^2h_2h_3h_2^{-1}h_3=1$:\\
$\Rightarrow \langle h_2,h_3\rangle=\langle h_2\rangle$ is abelian, a contradiction.

\item[(45)]$ R_1: h_2h_3h_2^{-2}h_3=1$, $R_2: h_2h_3h_2^{-1}h_3^{-1}h_2^{-1}h_3^2=1$:\\
$\Rightarrow\langle h_2,h_3\rangle\cong BS(1,1)$ is solvable, a contradiction.

\item[(46)]$ R_1: h_2h_3h_2^{-2}h_3=1$, $R_2: h_2h_3h_2^{-1}h_3^{-3}h_2=1$:\\
$\Rightarrow\langle h_2,h_3\rangle\cong BS(1,-1)$ is solvable, a contradiction.

\item[(47)]$ R_1: h_2h_3h_2^{-2}h_3=1$, $R_2: h_2(h_3h_2^{-1})^2h_3^{-2}h_2=1$:\\
$\Rightarrow\langle h_2,h_3\rangle\cong BS(1,1)$ is solvable, a contradiction.



\item[(50)]$ R_1: h_2h_3h_2^{-2}h_3=1$, $R_2: h_2h_3^2h_2^{-1}h_3^{-1}h_2^{-1}h_3=1$:\\
$\Rightarrow\langle h_2,h_3\rangle\cong BS(1,1)$ is solvable, a contradiction.

\item[(51)]$ R_1: h_2h_3h_2^{-2}h_3=1$, $R_2: h_2h_3h_2^{-1}h_3^{-1}h_2h_3^{-1}h_2^{-1}h_3=1$:\\
$\Rightarrow \langle h_2,h_3\rangle=\langle h_3\rangle$ is abelian, a contradiction.

\item[(52)]$ R_1: h_2(h_3h_2^{-1})^2h_3=1$, $R_2: h_2^2h_3h_2^{-1}(h_2^{-1}h_3)^2=1$:\\
$\Rightarrow\langle h_2,h_3\rangle\cong BS(1,1)$ is solvable, a contradiction.

\item[(53)]$ R_1: h_2(h_3h_2^{-1})^2h_3=1$, $R_2: h_2^2(h_3h_2^{-1})^2h_2^{-1}h_3=1$:\\
$\Rightarrow\langle h_2,h_3\rangle\cong BS(1,1)$ is solvable, a contradiction.

\item[(54)]$ R_1: h_2(h_3h_2^{-1})^2h_3=1$, $R_2: h_2^2h_3^{-1}h_2^{-1}(h_3^{-1}h_2)^2=1$:\\
$\Rightarrow h_2=1$ that is a contradiction.

\item[(55)]$ R_1: h_2(h_3h_2^{-1})^2h_3=1$, $R_2: h_2^2h_3^{-2}(h_2^{-1}h_3)^2=1$:\\
$\Rightarrow \langle h_2,h_3\rangle=\langle h_3\rangle$ is abelian, a contradiction.

\item[(56)]$ R_1: h_2(h_3h_2^{-1})^2h_3=1$, $R_2: h_2(h_2h_3^{-1})^2h_2^{-1}h_3^{-1}h_2=1$:\\
$\Rightarrow h_2=1$ that is a contradiction.

\item[(57)]$ R_1: h_2(h_3h_2^{-1})^2h_3=1$, $R_2: h_2h_3^{-1}h_2^{-1}h_3^2h_2^{-1}h_3^{-1}h_2=1$:\\
$\Rightarrow\langle h_2,h_3\rangle\cong BS(1,-1)$ is solvable, a contradiction.

\item[(58)]$ R_1: h_2(h_3h_2^{-1})^2h_3=1$, $R_2: h_2(h_3h_2^{-1})^2h_3^{-2}h_2=1$:\\
$\Rightarrow \langle h_2,h_3\rangle=\langle h_3\rangle$ is abelian, a contradiction.

\item[(59)]$ R_1: h_2(h_3h_2^{-1})^2h_3=1$, $R_2: h_2h_3^{-1}(h_2^{-1}h_3^2)^2=1$:\\
$\Rightarrow \langle h_2,h_3\rangle=\langle h_3\rangle$ is abelian, a contradiction.

\item[(60)]$ R_1: h_2(h_3h_2^{-1})^2h_3=1$, $R_2: h_2h_3^{-1}(h_2^{-1}h_3)^2h_3^2=1$:\\
$\Rightarrow\langle h_2,h_3\rangle\cong BS(1,-1)$ is solvable, a contradiction.

\item[(61)]$ R_1: h_2(h_3h_2^{-1})^2h_3=1$, $R_2: (h_2h_3^{-1})^2h_3^{-2}h_2^{-1}h_3=1$:\\
$\Rightarrow\langle h_2,h_3\rangle\cong BS(1,-1)$ is solvable, a contradiction.

\item[(62)]$ R_1: h_2(h_3h_2^{-1})^2h_3=1$, $R_2: h_3^3(h_3h_2^{-1})^2h_3=1$:\\
$\Rightarrow \langle h_2,h_3\rangle=\langle h_3\rangle$ is abelian, a contradiction.

\item[(63)]$ R_1: h_2h_3^{-3}h_2=1$, $R_2: h_2^2h_3^2h_2^{-2}h_3=1$:\\
By interchanging $h_2$ and $h_3$ in (28) and with the same discussion, there is a contradiction.

\item[(64)]$ R_1: h_2h_3^{-3}h_2=1$, $R_2: h_2^2h_3h_2^{-2}h_3^2=1$:\\
By interchanging $h_2$ and $h_3$ in (17) and with the same discussion, there is a contradiction.

\item[(65)]$ R_1: h_2h_3^{-3}h_2=1$, $R_2: h_2^2h_3^{-1}h_2^{-1}h_3^{-2}h_2=1$:\\
By interchanging $h_2$ and $h_3$ in (16) and with the same discussion, there is a contradiction.

\item[(66)]$ R_1: h_2h_3^{-3}h_2=1$, $R_2: h_2^2h_3^{-2}h_2^{-1}h_3^{-1}h_2=1$:\\
By interchanging $h_2$ and $h_3$ in (27) and with the same discussion, there is a contradiction.

\item[(67)]$ R_1: h_2h_3^{-3}h_2=1$, $R_2: h_2^2h_3^{-2}(h_2^{-1}h_3)^2=1$:\\
By interchanging $h_2$ and $h_3$ in (26) and with the same discussion, there is a contradiction.

\item[(68)]$ R_1: h_2h_3^{-3}h_2=1$, $R_2: h_2h_3h_2^{-1}h_3^{-1}h_2h_3^{-1}h_2^{-1}h_3=1$:\\
By interchanging $h_2$ and $h_3$ in (23) and with the same discussion, there is a contradiction.

\item[(69)]$ R_1: h_2h_3^{-3}h_2=1$, $R_2: h_2(h_3h_2^{-1})^2h_3^{-2}h_2=1$:\\
By interchanging $h_2$ and $h_3$ in (20) and with the same discussion, there is a contradiction.

\item[(70)]$ R_1: h_2h_3^{-3}h_2=1$, $R_2: h_2h_3^{-1}h_2^{-1}h_3^2h_2^{-1}h_3^{-1}h_2=1$:\\
By interchanging $h_2$ and $h_3$ in (21) and with the same discussion, there is a contradiction.

\item[(71)]$ R_1: h_2h_3^{-3}h_2=1$, $R_2: h_2h_3^{-1}h_2^{-1}h_3h_2^{-1}h_3^{-1}h_2h_3=1$:\\
By interchanging $h_2$ and $h_3$ in (18) and with the same discussion, there is a contradiction.

\item[(72)]$ R_1: h_2h_3^{-3}h_2=1$, $R_2: h_2h_3^{-1}h_2^{-1}(h_3h_2^{-1}h_3)^2=1$:\\
By interchanging $h_2$ and $h_3$ in (22) and with the same discussion, there is a contradiction.

\item[(73)]$ R_1: h_2h_3^{-3}h_2=1$, $R_2: h_2h_3^{-1}(h_2^{-1}h_3)^3h_3=1$:\\
By interchanging $h_2$ and $h_3$ in (19) and with the same discussion, there is a contradiction.

\item[(74)]$ R_1: h_2h_3^{-3}h_2=1$, $R_2: h_2(h_3^{-1}h_2h_3^{-1})^2h_2^{-1}h_3=1$:\\
By interchanging $h_2$ and $h_3$ in (24) and with the same discussion, there is a contradiction.

\item[(75)]$ R_1: h_2h_3^{-3}h_2=1$, $R_2: h_2h_3^{-4}h_2h_3=1$:\\
By interchanging $h_2$ and $h_3$ in (25) and with the same discussion, there is a contradiction.

\item[(76)]$ R_1: h_2h_3^{-3}h_2=1$, $R_2: (h_2h_3^{-1})^3h_3^{-1}h_2^{-1}h_3=1$:\\
By interchanging $h_2$ and $h_3$ in (29) and with the same discussion, there is a contradiction.

\item[(77)]$ R_1: h_2h_3^{-2}h_2h_3=1$, $R_2: h_2^3h_3^{-2}h_2^{-1}h_3=1$:\\
By interchanging $h_2$ and $h_3$ in (33) and with the same discussion, there is a contradiction.


\item[(79)]$ R_1: h_2h_3^{-2}h_2h_3=1$, $R_2: h_2^2h_3h_2h_3^{-1}h_2h_3=1$:\\
By interchanging $h_2$ and $h_3$ in (44) and with the same discussion, there is a contradiction.

\item[(80)]$ R_1: h_2h_3^{-2}h_2h_3=1$, $R_2: h_2^2h_3h_2^{-1}h_3^{-2}h_2=1$:\\
By interchanging $h_2$ and $h_3$ in (46) and with the same discussion, there is a contradiction.


\item[(82)]$ R_1: h_2h_3^{-2}h_2h_3=1$, $R_2: h_2h_3h_2h_3^{-1}h_2^{-1}h_3^{-1}h_2=1$:\\
By interchanging $h_2$ and $h_3$ in (45) and with the same discussion, there is a contradiction.

\item[(83)]$ R_1: h_2h_3^{-2}h_2h_3=1$, $R_2: h_2h_3h_2h_3^{-1}(h_3^{-1}h_2)^2=1$:\\
By interchanging $h_2$ and $h_3$ in (37) and with the same discussion, there is a contradiction.


\item[(85)]$ R_1: h_2h_3^{-2}h_2h_3=1$, $R_2: h_2(h_3h_2^{-1})^2h_3^{-2}h_2=1$:\\
By interchanging $h_2$ and $h_3$ in (47) and with the same discussion, there is a contradiction.

\item[(86)]$ R_1: h_2h_3^{-2}h_2h_3=1$, $R_2: h_2^2h_3^{-2}(h_2^{-1}h_3)^2=1$:\\
By interchanging $h_2$ and $h_3$ in (32) and with the same discussion, there is a contradiction.


\item[(88)]$ R_1: h_2h_3^{-2}h_2h_3=1$, $R_2: h_2(h_2h_3^{-1})^2h_3^{-1}h_2h_3=1$:\\
By interchanging $h_2$ and $h_3$ in (35) and with the same discussion, there is a contradiction.

\item[(89)]$ R_1: h_2h_3^{-2}h_2h_3=1$, $R_2: h_2h_3^{-1}h_2^{-1}h_3^2h_2^{-1}h_3^{-1}h_2=1$:\\
By interchanging $h_2$ and $h_3$ in (40) and with the same discussion, there is a contradiction.

\item[(90)]$ R_1: h_2h_3^{-2}h_2h_3=1$, $R_2: h_2h_3^{-2}(h_2^{-1}h_3)^3=1$:\\
By interchanging $h_2$ and $h_3$ in (31) and with the same discussion, there is a contradiction.

\item[(91)]$ R_1: h_2h_3^{-2}h_2h_3=1$, $R_2: h_2h_3^{-3}h_2h_3^2=1$:\\
By interchanging $h_2$ and $h_3$ in (30) and with the same discussion, there is a contradiction.

\item[(92)]$ R_1: h_2h_3^{-2}h_2h_3=1$, $R_2: h_2h_3^{-2}h_2^2h_3^{-1}h_2h_3=1$:\\
By interchanging $h_2$ and $h_3$ in (36) and with the same discussion, there is a contradiction.

\item[(93)]$ R_1: h_2h_3^{-2}h_2h_3=1$, $R_2: h_2h_3^{-2}h_2h_3h_2h_3^{-1}h_2=1$:\\
By interchanging $h_2$ and $h_3$ in (38) and with the same discussion, there is a contradiction.

\item[(94)]$ R_1: h_2h_3^{-2}h_2h_3=1$, $R_2: h_2h_3^{-1}(h_3^{-1}h_2)^2h_3h_2^{-1}h_3=1$:\\
By interchanging $h_2$ and $h_3$ in (41) and with the same discussion, there is a contradiction.

\item[(95)]$ R_1: h_2h_3^{-2}h_2h_3=1$, $R_2: h_2h_3^{-1}(h_3^{-1}h_2)^2h_3^{-2}h_2=1$:\\
By interchanging $h_2$ and $h_3$ in (42) and with the same discussion, there is a contradiction.

\item[(96)]$ R_1: h_2h_3^{-2}h_2h_3=1$, $R_2: (h_2h_3^{-1})^2h_3^{-1}h_2h_3h_2^{-1}h_3=1$:\\
By interchanging $h_2$ and $h_3$ in (43) and with the same discussion, there is a contradiction.

\item[(97)]$ R_1: h_2h_3^{-2}h_2h_3=1$, $R_2: h_2^2h_3^{-1}h_2^{-1}h_3^{-1}h_2h_3=1$:\\
By interchanging $h_2$ and $h_3$ in (50) and with the same discussion, there is a contradiction.

\item[(98)]$ R_1: h_2h_3^{-2}h_2h_3=1$, $R_2: h_2h_3^{-1}h_2^{-1}h_3h_2^{-1}h_3^{-1}h_2h_3=1$:\\
By interchanging $h_2$ and $h_3$ in (51) and with the same discussion, there is a contradiction.

\item[(99)]$ R_1: h_2h_3^{-1}h_2h_3^2=1$, $R_2: h_2^4h_3^2=1$:\\
By interchanging $h_2$ and $h_3$ in (7) and with the same discussion, there is a contradiction.

\item[(100)]$ R_1: h_2h_3^{-1}h_2h_3^2=1$, $R_2: h_2^3h_3^{-1}h_2^{-1}h_3^2=1$:\\
By interchanging $h_2$ and $h_3$ in (9) and with the same discussion, there is a contradiction.

\item[(101)]$ R_1: h_2h_3^{-1}h_2h_3^2=1$, $R_2: h_2^2h_3^2h_2^{-1}h_3^{-1}h_2=1$:\\
By interchanging $h_2$ and $h_3$ in (2) and with the same discussion, there is a contradiction.

\item[(102)]$ R_1: h_2h_3^{-1}h_2h_3^2=1$, $R_2: (h_2h_3)^2h_2^{-1}h_3^{-1}h_2=1$:\\
By interchanging $h_2$ and $h_3$ in (8) and with the same discussion, there is a contradiction.

\item[(103)]$ R_1: h_2h_3^{-1}h_2h_3^2=1$, $R_2: h_2(h_3h_2^{-1})^2h_3^{-2}h_2=1$:\\
By interchanging $h_2$ and $h_3$ in (10) and with the same discussion, there is a contradiction.

\item[(104)]$ R_1: h_2h_3^{-1}h_2h_3^2=1$, $R_2: h_2^2h_3^{-1}h_2^{-1}h_3h_2h_3=1$:\\
By interchanging $h_2$ and $h_3$ in (11) and with the same discussion, there is a contradiction.

\item[(105)]$ R_1: h_2h_3^{-1}h_2h_3^2=1$, $R_2: h_2h_3h_2^{-1}h_3^{-1}h_2h_3^{-1}h_2^{-1}h_3=1$:\\
By interchanging $h_2$ and $h_3$ in (15) and with the same discussion, there is a contradiction.

\item[(106)]$ R_1: h_2h_3^{-1}h_2h_3^2=1$, $R_2: h_2h_3^{-1}h_2^{-1}h_3^2h_2^{-1}h_3^{-1}h_2=1$:\\
By interchanging $h_2$ and $h_3$ in (6) and with the same discussion, there is a contradiction.

\item[(107)]$ R_1: h_2h_3^{-1}h_2h_3^2=1$, $R_2: h_2h_3^{-1}h_2^{-1}h_3h_2^{-1}h_3^{-1}h_2h_3=1$:\\
By interchanging $h_2$ and $h_3$ in (12) and with the same discussion, there is a contradiction.

\item[(108)]$ R_1: h_2h_3^{-1}h_2h_3^2=1$, $R_2: h_2h_3^{-1}h_2(h_3h_2^{-1})^2h_3^2=1$:\\
By interchanging $h_2$ and $h_3$ in (13) and with the same discussion, there is a contradiction.

\item[(109)]$ R_1: h_2h_3^{-1}h_2h_3^2=1$, $R_2: (h_2h_3^{-1})^2h_2h_3h_2^{-1}h_3^2=1$:\\
By interchanging $h_2$ and $h_3$ in (14) and with the same discussion, there is a contradiction.

\item[(110)]$ R_1: h_2h_3^{-1}h_2h_3^2=1$, $R_2: h_2^2h_3^{-2}(h_2^{-1}h_3)^2=1$:\\
By interchanging $h_2$ and $h_3$ in (3) and with the same discussion, there is a contradiction.

\item[(111)]$ R_1: h_2h_3^{-1}h_2h_3^2=1$, $R_2: h_2h_3^{-2}h_2h_3^3=1$:\\
By interchanging $h_2$ and $h_3$ in (1) and with the same discussion, there is a contradiction.

\item[(112)]$ R_1: h_2h_3^{-1}h_2h_3^2=1$, $R_2: h_2h_3^{-1}h_2h_3(h_3h_2^{-1})^2h_3=1$:\\
By interchanging $h_2$ and $h_3$ in (5) and with the same discussion, there is a contradiction.

\item[(113)]$ R_1: h_2h_3^{-1}h_2h_3^2=1$, $R_2: (h_2h_3^{-1})^2h_2h_3^2h_2^{-1}h_3=1$:\\
By interchanging $h_2$ and $h_3$ in (4) and with the same discussion, there is a contradiction.


\item[(115)]$ R_1: h_2h_3^{-1}h_2h_3h_2^{-1}h_3=1$, $R_2: h_2^3h_3^{-1}(h_2^{-1}h_3)^2=1$:\\
$\Rightarrow\langle h_2,h_3\rangle\cong BS(1,-1)$ is solvable, a contradiction.



\item[(118)]$ R_1: h_2h_3^{-1}h_2h_3h_2^{-1}h_3=1$, $R_2: h_2^2(h_3h_2^{-1})^2h_3^{-1}h_2=1$:\\
$\Rightarrow\langle h_2,h_3\rangle\cong BS(1,-1)$ is solvable, a contradiction.

\item[(119)]$ R_1: h_2h_3^{-1}h_2h_3h_2^{-1}h_3=1$, $R_2: h_2^2h_3^{-1}h_2h_3h_2h_3^{-1}h_2=1$:\\
$\Rightarrow \langle h_2,h_3\rangle=\langle h_2\rangle$ is abelian, a contradiction.

\item[(120)]$ R_1: h_2h_3^{-1}h_2h_3h_2^{-1}h_3=1$, $R_2: h_2h_3^2h_2^{-1}h_3h_2h_3^{-1}h_2=1$:\\
$\Rightarrow \langle h_2,h_3\rangle=\langle h_2\rangle$ is abelian, a contradiction.

\item[(121)]$ R_1: h_2h_3^{-1}h_2h_3h_2^{-1}h_3=1$, $R_2: h_2h_3h_2^{-2}h_3h_2h_3^{-1}h_2=1$:\\
$\Rightarrow\langle h_2,h_3\rangle\cong BS(1,1)$ is solvable, a contradiction.


\item[(123)]$ R_1: h_2h_3^{-1}h_2h_3h_2^{-1}h_3=1$, $R_2: h_2(h_3h_2^{-1})^2h_3^{-2}h_2=1$:\\
$\Rightarrow \langle h_2,h_3\rangle=\langle h_2\rangle$ is abelian, a contradiction.

\item[(124)]$ R_1: h_2h_3^{-1}h_2h_3h_2^{-1}h_3=1$, $R_2: h_2^2h_3^{-1}h_2h_3h_2^{-1}h_3^2=1$:\\
By interchanging $h_2$ and $h_3$ in (120) and with the same discussion, there is a contradiction.

\item[(125)]$ R_1: h_2h_3^{-1}h_2h_3h_2^{-1}h_3=1$, $R_2: h_2h_3h_2^{-1}h_3^3h_2^{-1}h_3=1$:\\
By interchanging $h_2$ and $h_3$ in (119) and with the same discussion, there is a contradiction.

\item[(126)]$ R_1: h_2h_3^{-1}h_2h_3h_2^{-1}h_3=1$, $R_2: h_2h_3^{-3}(h_2^{-1}h_3)^2=1$:\\
By interchanging $h_2$ and $h_3$ in (115) and with the same discussion, there is a contradiction.

\item[(127)]$ R_1: h_2h_3^{-1}h_2h_3h_2^{-1}h_3=1$, $R_2: h_2h_3^{-2}h_2h_3h_2^{-1}h_3^2=1$:\\
By interchanging $h_2$ and $h_3$ in (121) and with the same discussion, there is a contradiction.

\item[(128)]$ R_1: h_2h_3^{-1}h_2h_3h_2^{-1}h_3=1$, $R_2: (h_2h_3^{-1})^2h_2^{-1}h_3^3=1$:\\
By interchanging $h_2$ and $h_3$ in (118) and with the same discussion, there is a contradiction.

\item[(129)]$ R_1: (h_2h_3^{-1})^2h_2h_3=1$, $R_2: h_2^3(h_2h_3^{-1})^2h_2=1$:\\
By interchanging $h_2$ and $h_3$ in (62) and with the same discussion, there is a contradiction.

\item[(130)]$ R_1: (h_2h_3^{-1})^2h_2h_3=1$, $R_2: h_2^2(h_2h_3^{-1})^2h_2^{-1}h_3=1$:\\
By interchanging $h_2$ and $h_3$ in (61) and with the same discussion, there is a contradiction.

\item[(131)]$ R_1: (h_2h_3^{-1})^2h_2h_3=1$, $R_2: h_2^2h_3h_2^{-1}(h_3^{-1}h_2)^2=1$:\\
By interchanging $h_2$ and $h_3$ in (60) and with the same discussion, there is a contradiction.

\item[(132)]$ R_1: (h_2h_3^{-1})^2h_2h_3=1$, $R_2: h_2h_3h_2^{-1}h_3^{-1}h_2^2h_3^{-1}h_2=1$:\\
By interchanging $h_2$ and $h_3$ in (59) and with the same discussion, there is a contradiction.

\item[(133)]$ R_1: (h_2h_3^{-1})^2h_2h_3=1$, $R_2: h_2(h_3h_2^{-1})^2h_3^{-2}h_2=1$:\\
By interchanging $h_2$ and $h_3$ in (58) and with the same discussion, there is a contradiction.

\item[(134)]$ R_1: (h_2h_3^{-1})^2h_2h_3=1$, $R_2: h_2^2h_3^{-2}(h_2^{-1}h_3)^2=1$:\\
By interchanging $h_2$ and $h_3$ in (55) and with the same discussion, there is a contradiction.

\item[(135)]$ R_1: (h_2h_3^{-1})^2h_2h_3=1$, $R_2: h_2h_3^{-1}h_2^{-1}h_3^2h_2^{-1}h_3^{-1}h_2=1$:\\
By interchanging $h_2$ and $h_3$ in (57) and with the same discussion, there is a contradiction.

\item[(136)]$ R_1: (h_2h_3^{-1})^2h_2h_3=1$, $R_2: h_2h_3^{-2}(h_3^{-1}h_2)^2h_3=1$:\\
By interchanging $h_2$ and $h_3$ in (54) and with the same discussion, there is a contradiction.

\item[(137)]$ R_1: (h_2h_3^{-1})^2h_2h_3=1$, $R_2: h_2h_3^{-1}(h_3^{-1}h_2)^2h_3^2=1$:\\
By interchanging $h_2$ and $h_3$ in (52) and with the same discussion, there is a contradiction.

\item[(138)]$ R_1: (h_2h_3^{-1})^2h_2h_3=1$, $R_2: (h_2h_3^{-1})^2h_3^{-2}h_2h_3=1$:\\
By interchanging $h_2$ and $h_3$ in (56) and with the same discussion, there is a contradiction.

\item[(139)]$ R_1: (h_2h_3^{-1})^2h_2h_3=1$, $R_2: (h_2h_3^{-1})^2h_3^{-1}h_2h_3^2=1$:\\
By interchanging $h_2$ and $h_3$ in (53) and with the same discussion, there is a contradiction.

\end{enumerate}

$\mathbf{C_4-C_6(-C_6--)}$ \textbf{subgraph:} By considering the relations from Tables \ref{tab-C6} and \ref{tab-C4-C6} which are not disproved, it can be seen that there are $342$ cases for the relations of a cycle $C_4$ and two cycles $C_6$ in the graph $C_4-C_6(-C_6--)$. Using Gap \cite{gap}, we see that all groups with two generators $h_2$ and $h_3$ and three relations which are between $324$ cases of these $342$ cases are finite and solvable, that is a contradiction with the assumptions. So there are just $18$ cases for the relations of these cycles which may lead to the existence of a subgraph isomorphic to the graph $C_4-C_6(-C_6--)$ in the graph $K(\alpha,\beta)$. These cases are listed in table \ref{tab-C4-C6(-C6--)}.  In the following, we show that $12$ cases of these relations lead to a contradiction and just $6$ cases of them may lead to the  existence of a subgraph isomorphic to the graph $C_4-C_6(-C_6--)$ in the graph $K(\alpha,\beta)$. Cases which are not disproved are marked by $*$s in the Table \ref{tab-C4-C6(-C6--)}.

\begin{table}[h]
\centering
\caption{The relations of a $C_4-C_6(-C_6--)$ in $K(\alpha,\beta)$}\label{tab-C4-C6(-C6--)}
\begin{tabular}{|c|l|l|l|}\hline
$n$&$R_1$&$R_2$&$R_3$\\\hline
$1 *$&$h_2h_3h_2^{-2}h_3=1$&$h_2(h_3h_2^{-1})^2h_2^{-1}h_3^2=1$&$h_3^2h_2^{-1}(h_3h_2^{-1}h_3)^2=1$\\
$2$&$h_2h_3h_2^{-2}h_3=1$&$h_3^2(h_3h_2^{-1})^3h_3=1$&$h_2h_3h_2^{-1}h_3^2h_2^{-2}h_3=1$\\
$3$&$h_2h_3h_2^{-2}h_3=1$&$h_3^2h_2^{-1}(h_3h_2^{-1}h_3)^2=1$&$h_2h_3h_2^{-1}h_3^2h_2^{-2}h_3=1$\\
$4 *$&$h_2h_3h_2^{-2}h_3=1$&$h_3^2h_2^{-1}(h_3h_2^{-1}h_3)^2=1$&$h_2(h_3h_2^{-1})^2h_2^{-1}h_3^2=1$\\
$5$&$h_2h_3h_2^{-2}h_3=1$&$h_3^2h_2^{-1}(h_3h_2^{-1}h_3)^2=1$&$h_2h_3^{-1}h_2(h_3h_2^{-1})^2h_2^{-1}h_3=1$\\
$6$&$h_2h_3^{-2}h_2h_3=1$&$h_2^2(h_2h_3^{-1})^3h_2=1$&$h_2h_3^{-2}h_2h_3h_2h_3^{-1}h_2=1$\\
$7 *$&$h_2h_3^{-2}h_2h_3=1$&$h_2^2h_3^{-1}(h_2h_3^{-1}h_2)^2=1$&$h_2h_3(h_2h_3^{-1})^2h_3^{-1}h_2=1$\\
$8$&$h_2h_3^{-2}h_2h_3=1$&$h_2^2h_3^{-1}(h_2h_3^{-1}h_2)^2=1$&$h_2h_3^{-2}h_2h_3h_2h_3^{-1}h_2=1$\\
$9$&$h_2h_3^{-2}h_2h_3=1$&$h_2^2h_3^{-1}(h_2h_3^{-1}h_2)^2=1$&$(h_2h_3^{-1})^2h_3^{-1}h_2h_3h_2^{-1}h_3=1$\\
$10 *$&$h_2h_3^{-2}h_2h_3=1$&$h_2h_3(h_2h_3^{-1})^2h_3^{-1}h_2=1$&$h_2^2h_3^{-1}(h_2h_3^{-1}h_2)^2=1$\\
$11$&$h_2h_3^{-1}h_2h_3h_2^{-1}h_3=1$&$h_2^3h_3^3=1$&$h_2h_3h_2h_3^{-1}h_2h_3h_2^{-1}h_3=1$\\
$12$&$h_2h_3^{-1}h_2h_3h_2^{-1}h_3=1$&$h_2^3h_3^3=1$&$h_2h_3h_2^{-1}h_3h_2h_3^{-1}h_2h_3=1$\\
$13$&$h_2h_3^{-1}h_2h_3h_2^{-1}h_3=1$&$h_2(h_2h_3)^2h_3=1$&$h_2h_3h_2h_3^{-1}h_2h_3h_2^{-1}h_3=1$\\
$14 *$&$h_2h_3^{-1}h_2h_3h_2^{-1}h_3=1$&$h_2(h_2h_3)^2h_3=1$&$h_2h_3h_2^{-1}h_3^2h_2h_3^{-1}h_2=1$\\
$15$&$h_2h_3^{-1}h_2h_3h_2^{-1}h_3=1$&$h_2(h_2h_3)^2h_3=1$&$h_2h_3^{-2}h_2h_3h_2^{-1}h_3^2=1$\\
$16 *$&$h_2h_3^{-1}h_2h_3h_2^{-1}h_3=1$&$h_2^2h_3^2h_2h_3=1$&$h_2h_3h_2^{-1}h_3^2h_2h_3^{-1}h_2=1$\\
$17$&$h_2h_3^{-1}h_2h_3h_2^{-1}h_3=1$&$h_2^2h_3^2h_2h_3=1$&$h_2h_3h_2^{-2}h_3h_2h_3^{-1}h_2=1$\\
$18$&$h_2h_3^{-1}h_2h_3h_2^{-1}h_3=1$&$h_2^2h_3^2h_2h_3=1$&$h_2h_3h_2^{-1}h_3h_2h_3^{-1}h_2h_3=1$\\
\hline
\end{tabular}
\end{table}

\begin{enumerate}

\item[(2)]$ R_1:h_2h_3h_2^{-2}h_3=1$, $R_2:h_3^2(h_3h_2^{-1})^3h_3=1$, $R_3:h_2h_3h_2^{-1}h_3^2h_2^{-2}h_3=1$:\\
$\Rightarrow \langle h_2,h_3\rangle=\langle h_3\rangle$ is abelian, a contradiction.

\item[(3)]$ R_1:h_2h_3h_2^{-2}h_3=1$, $R_2:h_3^2h_2^{-1}(h_3h_2^{-1}h_3)^2=1$, $R_3:h_2h_3h_2^{-1}h_3^2h_2^{-2}h_3=1$:\\
$\Rightarrow \langle h_2,h_3\rangle=\langle h_3\rangle$ is abelian, a contradiction.


\item[(5)]$ R_1:h_2h_3h_2^{-2}h_3=1$, $R_2:h_3^2h_2^{-1}(h_3h_2^{-1}h_3)^2=1$, $R_3:h_2h_3^{-1}h_2(h_3h_2^{-1})^2h_2^{-1}h_3=1$:\\
$\Rightarrow\langle h_2,h_3\rangle\cong BS(1,1)$ is solvable, a contradiction.

\item[(6)]$ R_1:h_2h_3^{-2}h_2h_3=1$, $R_2:h_2^2(h_2h_3^{-1})^3h_2=1$, $R_3:h_2h_3^{-2}h_2h_3h_2h_3^{-1}h_2=1$:\\
By interchanging $h_2$ and $h_3$ in (2) and with the same discussion, there is a contradiction.


\item[(8)]$ R_1:h_2h_3^{-2}h_2h_3=1$, $R_2:h_2^2h_3^{-1}(h_2h_3^{-1}h_2)^2=1$, $R_3:h_2h_3^{-2}h_2h_3h_2h_3^{-1}h_2=1$:\\
By interchanging $h_2$ and $h_3$ in (3) and with the same discussion, there is a contradiction.

\item[(9)]$ R_1:h_2h_3^{-2}h_2h_3=1$, $R_2:h_2^2h_3^{-1}(h_2h_3^{-1}h_2)^2=1$, $R_3:(h_2h_3^{-1})^2h_3^{-1}h_2h_3h_2^{-1}h_3=1$:\\
By interchanging $h_2$ and $h_3$ in (5) and with the same discussion, there is a contradiction.


\item[(11)]$ R_1:h_2h_3^{-1}h_2h_3h_2^{-1}h_3=1$, $R_2:h_2^3h_3^3=1$, $R_3:h_2h_3h_2h_3^{-1}h_2h_3h_2^{-1}h_3=1$:\\
$\Rightarrow \langle h_2,h_3\rangle=\langle h_2\rangle$ is abelian, a contradiction.

\item[(12)]$ R_1:h_2h_3^{-1}h_2h_3h_2^{-1}h_3=1$, $R_2:h_2^3h_3^3=1$, $R_3:h_2h_3h_2^{-1}h_3h_2h_3^{-1}h_2h_3=1$:\\
By interchanging $h_2$ and $h_3$ in (11) and with the same discussion, there is a contradiction.

\item[(13)]$ R_1:h_2h_3^{-1}h_2h_3h_2^{-1}h_3=1$, $R_2:h_2(h_2h_3)^2h_3=1$, $R_3:h_2h_3h_2h_3^{-1}h_2h_3h_2^{-1}h_3=1$:\\
$\Rightarrow \langle h_2,h_3\rangle=\langle h_2\rangle$ is abelian, a contradiction.


\item[(15)]$ R_1:h_2h_3^{-1}h_2h_3h_2^{-1}h_3=1$, $R_2:h_2(h_2h_3)^2h_3=1$, $R_3:h_2h_3^{-2}h_2h_3h_2^{-1}h_3^2=1$:\\
$\Rightarrow\langle h_2,h_3\rangle\cong BS(1,1)$ is solvable, a contradiction.


\item[(17)]$ R_1:h_2h_3^{-1}h_2h_3h_2^{-1}h_3=1$, $R_2:h_2^2h_3^2h_2h_3=1$, $R_3:h_2h_3h_2^{-2}h_3h_2h_3^{-1}h_2=1$:\\
By interchanging $h_2$ and $h_3$ in (15) and with the same discussion, there is a contradiction.

\item[(18)]$ R_1:h_2h_3^{-1}h_2h_3h_2^{-1}h_3=1$, $R_2:h_2^2h_3^2h_2h_3=1$, $R_3:h_2h_3h_2^{-1}h_3h_2h_3^{-1}h_2h_3=1$:\\
By interchanging $h_2$ and $h_3$ in (13) and with the same discussion, there is a contradiction.

\end{enumerate}

$\mathbf{C_4-C_6(--C_7--)}$ \textbf{subgraph:} With the same discussion such as about $C_4$ cycles, by considering the relations of $C_7$ cycles and the relations from Table \ref{tab-C4-C6} which are not disproved, it can be seen that there are $521$ cases for the relations of a cycle $C_4$, a cycle $C_6$ and a cycle $C_7$ in the graph $C_4-C_6(--C_7--)$. By considering all groups with two generators $h_2$ and $h_3$ and three relations which are between these cases and by using Gap \cite{gap}, we see that $484$ groups are finite and solvable and $4$ groups have the same ``structure description'' $C_5\times SL(2,5)$ according to the function StructureDescription of GAP, that is finite. So there are just $33$ cases for the relations of these cycles which may lead to the existence of a subgraph isomorphic to the graph $C_4-C_6(--C_7--)$ in the graph $K(\alpha,\beta)$. These cases are listed in table \ref{tab-C4-C6(--C7--)}.  In the following, we show that $17$ cases of these relations lead to a contradiction and $16$ cases of them may lead to the existence of a subgraph isomorphic to the graph $C_4-C_6(--C_7--)$ in the graph $K(\alpha,\beta)$. Cases which are not disproved are marked by $*$s in the Table \ref{tab-C4-C6(--C7--)}.
\begin{table}[h]
\centering
\caption{The relations of a $C_4-C_6(--C_7--)$ in $K(\alpha,\beta)$}\label{tab-C4-C6(--C7--)}
\begin{tabular}{|c|l|l|l|}\hline
$n$&$R_1$&$R_2$&$R_3$\\\hline
$1*$&$h_2h_3h_2^{-2}h_3=1$&$h_2h_3^2h_2^{-1}(h_2^{-1}h_3)^2=1$&$h_2h_3^{-5}h_2h_3=1$\\
$2*$&$h_2h_3h_2^{-2}h_3=1$&$h_2h_3^2h_2^{-1}(h_2^{-1}h_3)^2=1$&$h_2h_3^{-4}h_2h_3h_2^{-1}h_3=1$\\
$3*$&$h_2h_3h_2^{-2}h_3=1$&$h_2h_3^2h_2^{-1}(h_2^{-1}h_3)^2=1$&$h_2h_3^{-2}h_2(h_3h_2^{-1})^3h_3=1$\\
$4*$&$h_2h_3h_2^{-2}h_3=1$&$h_2h_3^2h_2^{-1}(h_2^{-1}h_3)^2=1$&$h_2h_3^{-1}h_2(h_3h_2^{-1})^4h_3=1$\\
$5*$&$h_2h_3h_2^{-2}h_3=1$&$h_2(h_3h_2^{-1})^2h_2^{-1}h_3^2=1$&$(h_2h_3^{-2})^2h_2^{-1}h_3^2=1$\\
$6*$&$h_2h_3h_2^{-2}h_3=1$&$h_2(h_3h_2^{-1})^2h_2^{-1}h_3^2=1$&$(h_2h_3^{-2})^2h_2h_3^{-1}h_2^{-1}h_3=1$\\
$7$&$h_2h_3h_2^{-2}h_3=1$&$h_2(h_3h_2^{-1})^2h_2^{-1}h_3^2=1$&$h_2^3h_3^{-1}h_2^3h_3=1$\\
$8*$&$h_2h_3h_2^{-2}h_3=1$&$h_2(h_3h_2^{-1})^2h_2^{-1}h_3^2=1$&$h_2h_3^{-1}(h_2^{-1}h_3^2)^2h_2^{-1}h_3=1$\\
$9*$&$h_2h_3h_2^{-2}h_3=1$&$h_2(h_3h_2^{-1})^2h_2^{-1}h_3^2=1$&$h_2h_3^{-1}h_2h_3^{-1}h_2^{-1}(h_3h_2^{-1}h_3)^2=1$\\
$10$&$h_2h_3h_2^{-2}h_3=1$&$h_3^2(h_3h_2^{-1})^3h_3=1$&$h_2h_3^{-2}h_2(h_3h_2^{-1})^3h_3=1$\\
$11$&$h_2h_3h_2^{-2}h_3=1$&$h_3^2h_2^{-1}(h_3h_2^{-1}h_3)^2=1$&$(h_2h_3^{-1})^2h_2^{-1}h_3(h_3h_2^{-1})^2h_3=1$\\
$12$&$h_2h_3h_2^{-2}h_3=1$&$h_3^2h_2^{-1}(h_3h_2^{-1}h_3)^2=1$&$h_2h_3^{-1}h_2h_3^{-1}h_2^{-1}(h_3h_2^{-1}h_3)^2=1$\\
$13$&$h_2h_3^{-2}h_2h_3=1$&$h_2^2(h_2h_3^{-1})^3h_2=1$&$h_2h_3^{-1}(h_2^{-1}h_3)^3h_2^{-1}h_3^{-1}h_2=1$\\
$14$&$h_2h_3^{-2}h_2h_3=1$&$h_2^2h_3^{-1}(h_2h_3^{-1}h_2)^2=1$&$h_2h_3^{-1}h_2(h_3h_2^{-1})^2(h_3^{-1}h_2)^2=1$\\
$15$&$h_2h_3^{-2}h_2h_3=1$&$h_2^2h_3^{-1}(h_2h_3^{-1}h_2)^2=1$&$(h_2h_3^{-1})^2(h_2^{-1}h_3)^2h_2^{-2}h_3=1$\\
$16$&$h_2h_3^{-2}h_2h_3=1$&$h_2h_3(h_2h_3^{-1})^2h_3^{-1}h_2=1$&$h_2h_3^3h_2^{-1}h_3^3=1$\\
$17*$&$h_2h_3^{-2}h_2h_3=1$&$h_2h_3(h_2h_3^{-1})^2h_3^{-1}h_2=1$&$(h_2^2h_3^{-1})^2h_2^{-2}h_3=1$\\
$18*$&$h_2h_3^{-2}h_2h_3=1$&$h_2h_3(h_2h_3^{-1})^2h_3^{-1}h_2=1$&$h_2h_3^{-1}h_2^{-1}h_3(h_2h_3^{-1}h_2)^2=1$\\
$19*$&$h_2h_3^{-2}h_2h_3=1$&$h_2h_3(h_2h_3^{-1})^2h_3^{-1}h_2=1$&$h_2h_3^{-1}(h_2^{-1}h_3h_2^{-1})^2h_2^{-1}h_3=1$\\
$20*$&$h_2h_3^{-2}h_2h_3=1$&$h_2h_3(h_2h_3^{-1})^2h_3^{-1}h_2=1$&$h_2h_3^{-1}h_2(h_3h_2^{-1})^2(h_3^{-1}h_2)^2=1$\\
$21*$&$h_2h_3^{-2}h_2h_3=1$&$h_2^2h_3^{-1}(h_3^{-1}h_2)^2h_3=1$&$h_2^4h_3^{-1}h_2^{-1}h_3^{-1}h_2=1$\\
$22*$&$h_2h_3^{-2}h_2h_3=1$&$h_2^2h_3^{-1}(h_3^{-1}h_2)^2h_3=1$&$h_2^3h_3^{-1}h_2^{-1}h_3h_2^{-1}h_3^{-1}h_2=1$\\
$23*$&$h_2h_3^{-2}h_2h_3=1$&$h_2^2h_3^{-1}(h_3^{-1}h_2)^2h_3=1$&$h_2h_3^{-1}(h_2^{-1}h_3)^3h_2^{-1}h_3^{-1}h_2=1$\\
$24*$&$h_2h_3^{-2}h_2h_3=1$&$h_2^2h_3^{-1}(h_3^{-1}h_2)^2h_3=1$&$(h_2h_3^{-1})^4h_2h_3h_2^{-1}h_3=1$\\
$25$&$h_2h_3^{-1}h_2h_3h_2^{-1}h_3=1$&$h_2^3h_3^3=1$&$h_2h_3h_2^{-1}h_3^{-3}h_2^{-1}h_3=1$\\
$26$&$h_2h_3^{-1}h_2h_3h_2^{-1}h_3=1$&$h_2^3h_3^3=1$&$h_2^3h_3h_2^{-1}h_3^{-1}h_2^{-1}h_3=1$\\
$27$&$h_2h_3^{-1}h_2h_3h_2^{-1}h_3=1$&$h_2(h_2h_3)^2h_3=1$&$h_2h_3^2h_2^{-1}h_3^{-2}h_2h_3=1$\\
$28$&$h_2h_3^{-1}h_2h_3h_2^{-1}h_3=1$&$h_2(h_2h_3)^2h_3=1$&$h_2h_3h_2^{-1}h_3^{-2}h_2h_3^2=1$\\
$29$&$h_2h_3^{-1}h_2h_3h_2^{-1}h_3=1$&$h_2(h_2h_3)^2h_3=1$&$h_2(h_2h_3^{-1})^2h_2^{-1}(h_2^{-1}h_3)^2=1$\\
$30$&$h_2h_3^{-1}h_2h_3h_2^{-1}h_3=1$&$h_2^2h_3^2h_2h_3=1$&$(h_2h_3^{-1})^2h_3^{-1}(h_2^{-1}h_3)^2h_3=1$\\
$31$&$h_2h_3^{-1}h_2h_3h_2^{-1}h_3=1$&$h_2^2h_3^2h_2h_3=1$&$h_2^2h_3h_2h_3^{-1}h_2^{-2}h_3=1$\\
$32$&$h_2h_3^{-1}h_2h_3h_2^{-1}h_3=1$&$h_2^2h_3^2h_2h_3=1$&$h_2^2h_3^{-1}h_2^{-2}h_3h_2h_3=1$\\
$33$&$h_2h_3^{-1}h_2h_3h_2^{-1}h_3=1$&$h_2h_3h_2^{-1}h_3^2h_2h_3^{-1}h_2=1$&$h_2((h_3^{-1}h_2)^2h_3^{-1})^2h_2=1$\\
\hline
\end{tabular}
\end{table}

\begin{enumerate}






\item[(7)]$ R_1:h_2h_3h_2^{-2}h_3=1$, $R_2:h_2(h_3h_2^{-1})^2h_2^{-1}h_3^2=1$, $R_3:h_2^3h_3^{-1}h_2^3h_3=1$:\\
$\Rightarrow h_2=1$ that is a contradiction.



\item[(10)]$ R_1:h_2h_3h_2^{-2}h_3=1$, $R_2:h_3^2(h_3h_2^{-1})^3h_3=1$, $R_3:h_2h_3^{-2}h_2(h_3h_2^{-1})^3h_3=1$:\\
$\Rightarrow \langle h_2,h_3\rangle=\langle h_3\rangle$ is abelian, a contradiction.

\item[(11)]$ R_1:h_2h_3h_2^{-2}h_3=1$, $R_2:h_3^2h_2^{-1}(h_3h_2^{-1}h_3)^2=1$, $R_3:(h_2h_3^{-1})^2h_2^{-1}h_3(h_3h_2^{-1})^2h_3=1$:\\
$\Rightarrow \langle h_2,h_3\rangle=\langle h_3\rangle$ is abelian, a contradiction.

\item[(12)]$ R_1:h_2h_3h_2^{-2}h_3=1$, $R_2:h_3^2h_2^{-1}(h_3h_2^{-1}h_3)^2=1$, $R_3:h_2h_3^{-1}h_2h_3^{-1}h_2^{-1}(h_3h_2^{-1}h_3)^2=1$:\\
$\Rightarrow\langle h_2,h_3\rangle\cong BS(1,-1)$ is solvable, a contradiction.

\item[(13)]$ R_1:h_2h_3^{-2}h_2h_3=1$, $R_2:h_2^2(h_2h_3^{-1})^3h_2=1$, $R_3:h_2h_3^{-1}(h_2^{-1}h_3)^3h_2^{-1}h_3^{-1}h_2=1$:\\
By interchanging $h_2$ and $h_3$ in (10) and with the same discussion, there is a contradiction.

\item[(14)]$ R_1:h_2h_3^{-2}h_2h_3=1$, $R_2:h_2^2h_3^{-1}(h_2h_3^{-1}h_2)^2=1$, $R_3:h_2h_3^{-1}h_2(h_3h_2^{-1})^2(h_3^{-1}h_2)^2=1$:\\
By interchanging $h_2$ and $h_3$ in (12) and with the same discussion, there is a contradiction.

\item[(15)]$ R_1:h_2h_3^{-2}h_2h_3=1$, $R_2:h_2^2h_3^{-1}(h_2h_3^{-1}h_2)^2=1$, $R_3:(h_2h_3^{-1})^2(h_2^{-1}h_3)^2h_2^{-2}h_3=1$:\\
By interchanging $h_2$ and $h_3$ in (11) and with the same discussion, there is a contradiction.

\item[(16)]$ R_1:h_2h_3^{-2}h_2h_3=1$, $R_2:h_2h_3(h_2h_3^{-1})^2h_3^{-1}h_2=1$, $R_3:h_2h_3^3h_2^{-1}h_3^3=1$:\\
By interchanging $h_2$ and $h_3$ in (7) and with the same discussion, there is a contradiction.









\item[(25)]$ R_1:h_2h_3^{-1}h_2h_3h_2^{-1}h_3=1$, $R_2:h_2^3h_3^3=1$, $R_3:h_2h_3h_2^{-1}h_3^{-3}h_2^{-1}h_3=1$:\\
$\Rightarrow \langle h_2,h_3\rangle=\langle h_2\rangle$ is abelian, a contradiction.

\item[(26)]$ R_1:h_2h_3^{-1}h_2h_3h_2^{-1}h_3=1$, $R_2:h_2^3h_3^3=1$, $R_3:h_2^3h_3h_2^{-1}h_3^{-1}h_2^{-1}h_3=1$:\\
By interchanging $h_2$ and $h_3$ in (25) and with the same discussion, there is a contradiction.

\item[(27)]$ R_1:h_2h_3^{-1}h_2h_3h_2^{-1}h_3=1$, $R_2:h_2(h_2h_3)^2h_3=1$, $R_3:h_2h_3^2h_2^{-1}h_3^{-2}h_2h_3=1$:\\
$\Rightarrow\langle h_2,h_3\rangle\cong BS(1,-1)$ is solvable, a contradiction.

\item[(28)]$ R_1:h_2h_3^{-1}h_2h_3h_2^{-1}h_3=1$, $R_2:h_2(h_2h_3)^2h_3=1$, $R_3:h_2h_3h_2^{-1}h_3^{-2}h_2h_3^2=1$:\\
$\Rightarrow \langle h_2,h_3\rangle=\langle h_2\rangle$ is abelian, a contradiction.

\item[(29)]$ R_1:h_2h_3^{-1}h_2h_3h_2^{-1}h_3=1$, $R_2:h_2(h_2h_3)^2h_3=1$, $R_3:h_2(h_2h_3^{-1})^2h_2^{-1}(h_2^{-1}h_3)^2=1$:\\
$\Rightarrow\langle h_2,h_3\rangle\cong BS(1,-1)$ is solvable, a contradiction.

\item[(30)]$ R_1:h_2h_3^{-1}h_2h_3h_2^{-1}h_3=1$, $R_2:h_2^2h_3^2h_2h_3=1$, $R_3:(h_2h_3^{-1})^2h_3^{-1}(h_2^{-1}h_3)^2h_3=1$:\\
By interchanging $h_2$ and $h_3$ in (29) and with the same discussion, there is a contradiction.

\item[(31)]$ R_1:h_2h_3^{-1}h_2h_3h_2^{-1}h_3=1$, $R_2:h_2^2h_3^2h_2h_3=1$, $R_3:h_2^2h_3h_2h_3^{-1}h_2^{-2}h_3=1$:\\
By interchanging $h_2$ and $h_3$ in (28) and with the same discussion, there is a contradiction.

\item[(32)]$ R_1:h_2h_3^{-1}h_2h_3h_2^{-1}h_3=1$, $R_2:h_2^2h_3^2h_2h_3=1$, $R_3:h_2^2h_3^{-1}h_2^{-2}h_3h_2h_3=1$:\\
By interchanging $h_2$ and $h_3$ in (27) and with the same discussion, there is a contradiction.

\item[(33)]$ R_1:h_2h_3^{-1}h_2h_3h_2^{-1}h_3=1$, $R_2:h_2h_3h_2^{-1}h_3^2h_2h_3^{-1}h_2=1$, $R_3:h_2((h_3^{-1}h_2)^2h_3^{-1})^2h_2=1$:\\
$\Rightarrow \langle h_2,h_3\rangle=\langle h_2\rangle$ is abelian, a contradiction.

\end{enumerate}

\subsection{$\mathbf{C_4-C_6(--C_7--)(C_7-1)}$}
By considering the relations of a $C_7$ cycle in the graph $K(\alpha,\beta)$ and relations from Table \ref{tab-C4-C6(--C7--)} which are not disproved, it can be seen that there are $176$ cases for the relations of a cycle $C_4$, two cycles $C_7$ and a cycle $C_6$ in this structure. By considering all groups with two generators $h_2$ and $h_3$ and four relations which are between these cases and by using GAP \cite{gap}, we see that all of these groups are finite and solvable. So, the graph $K(\alpha,\beta)$ contains no subgraph isomorphic to the graph $C_4-C_6(--C_7--)(C_7-1)$.

\subsection{$\mathbf{C_4-C_6(--C_7--)(--C_5-)}$}
By considering the relations from Tables \ref{tab-C5} and \ref{tab-C4-C6(--C7--)} which are not disproved, it can be seen that there are $28$ cases for the relations of a cycle $C_4$, a cycle $C_6$, a cycle $C_7$ and a cycle $C_5$ in the graph $C_4-C_6(--C_7--)(--C_5-)$. Using GAP \cite{gap}, we see that all groups with two generators $h_2$ and $h_3$ and four relations which are between $24$ cases of these $28$ cases are finite and solvable, that is a contradiction with the assumptions. So, there are just $4$ cases for the relations of these cycles which may lead to the existence of a subgraph isomorphic to the graph $C_4-C_6(--C_7--)(--C_5-)$ in $K(\alpha,\beta)$. These cases are listed in table \ref{tab-C4-C6(--C7--)(--C5-)}.  In the following, we show that all of these $4$ cases lead to contradictions and so, the graph $K(\alpha,\beta)$ contains no subgraph isomorphic to the graph $C_4-C_6(--C_7--)(--C_5-)$.

\begin{table}[h]
\centering
\caption{The relations of a $C_4-C_6(--C_7--)(--C_5-)$ in $K(\alpha,\beta)$}\label{tab-C4-C6(--C7--)(--C5-)}
\begin{tabular}{|c|l|l|l|l|}\hline
$n$&$R_1$&$R_2$&$R_3$&$R_4$\\\hline
$1$&$h_2h_3h_2^{{-2}}h_3=1$&$h_2h_3^2h_2^{-1}(h_2^{-1}h_3)^2=1$&$h_2h_3^{-4}h_2h_3h_2^{-1}h_3=1$&$(h_2h_3^{-1})^2h_3^{-1}h_2^{-1}h_3=1$\\
$2$&$h_2h_3h_2^{{-2}}h_3=1$&$h_2h_3^2h_2^{-1}(h_2^{-1}h_3)^2=1$&$h_2h_3^{{-2}}h_2(h_3h_2^{-1})^3h_3=1$&$h_2h_3^{-1}(h_2^{-1}h_3)^2h_3=1$\\
$3$&$h_2h_3^{{-2}}h_2h_3=1$&$h_2^2h_3^{-1}(h_3^{-1}h_2)^2h_3=1$&$h_2^3h_3^{-1}h_2^{-1}h_3h_2^{-1}h_3^{-1}h_2=1$&$h_2(h_2h_3^{-1})^2h_2^{-1}h_3=1$\\
$4$&$h_2h_3^{{-2}}h_2h_3=1$&$h_2^2h_3^{-1}(h_3^{-1}h_2)^2h_3=1$&$h_2h_3^{-1}(h_2^{-1}h_3)^3h_2^{-1}h_3^{-1}h_2=1$&$h_2h_3h_2^{-1}(h_3^{-1}h_2)^2=1$\\
\hline
\end{tabular}
\end{table}

\begin{enumerate}
\item[(1)]$ R_1:h_2h_3h_2^{-2}h_3=1$, $R_2:h_2h_3^2h_2^{-1}(h_2^{-1}h_3)^2=1$, $R_3:h_2h_3^{-4}h_2h_3h_2^{-1}h_3=1$,\\ $R_4:(h_2h_3^{-1})^2h_3^{-1}h_2^{-1}h_3=1$:\\
$\Rightarrow \langle h_2,h_3\rangle=\langle h_3\rangle$ is abelian, a contradiction.

\item[(2)]$ R_1:h_2h_3h_2^{-2}h_3=1$, $R_2:h_2h_3^2h_2^{-1}(h_2^{-1}h_3)^2=1$, $R_3:h_2h_3^{-2}h_2(h_3h_2^{-1})^3h_3=1$,\\ $R_4:h_2h_3^{-1}(h_2^{-1}h_3)^2h_3=1$:\\
$\Rightarrow \langle h_2,h_3\rangle=\langle h_3\rangle$ is abelian, a contradiction.

\item[(3)]$ R_1:h_2h_3^{-2}h_2h_3=1$, $R_2:h_2^2h_3^{-1}(h_3^{-1}h_2)^2h_3=1$, $R_3:h_2^3h_3^{-1}h_2^{-1}h_3h_2^{-1}h_3^{-1}h_2=1$, \\$R_4:h_2(h_2h_3^{-1})^2h_2^{-1}h_3=1$:\\
By interchanging $h_2$ and $h_3$ in (1) and with the same discussion, there is a contradiction.

\item[(4)]$ R_1:h_2h_3^{-2}h_2h_3=1$, $R_2:h_2^2h_3^{-1}(h_3^{-1}h_2)^2h_3=1$, $R_3:h_2h_3^{-1}(h_2^{-1}h_3)^3h_2^{-1}h_3^{-1}h_2=1$, \\$R_4:h_2h_3h_2^{-1}(h_3^{-1}h_2)^2=1$:\\
By interchanging $h_2$ and $h_3$ in (2) and with the same discussion, there is a contradiction.

\end{enumerate}

\subsection{$\mathbf{C_4-C_6(-C_6--)(-C_4-)}$}
By considering the relations from Tables \ref{tab-C4} and \ref{tab-C4-C6(-C6--)} which are not disproved, it can be seen that there is no case for the relations of two cycles $C_4$ and two cycles $C_6$ in this structure. It means that the graph $K(\alpha,\beta)$ contains no subgraph isomorphic to the graph $C_4-C_6(-C_6--)(-C_4-)$.

\subsection{$\mathbf{C_4-C_6(-C_6--)(--C_5-)}$}
By considering the relations from Tables \ref{tab-C5} and \ref{tab-C4-C6(-C6--)} which are not disproved, it can be seen that there are $22$ cases for the relations of a cycle $C_4$, two cycles $C_6$ and a cycle $C_5$ in this structure. By considering all groups with two generators $h_2$ and $h_3$ and four relations which are between these cases and by using GAP \cite{gap}, we see that all of these groups are finite and solvable. So, the graph $K(\alpha,\beta)$ contains no subgraph isomorphic to the graph $C_4-C_6(-C_6--)(--C_5-)$.

\subsection{$\mathbf{C_4-C_6(-C_6--)(C_6---)}$}
By considering the relations from Tables \ref{tab-C6} and \ref{tab-C4-C6(-C6--)} which are not disproved, it can be seen that there are $66$ cases for the relations of a cycle $C_4$ and three cycles $C_6$ in the graph $C_4-C_6(-C_6--)(C_6---)$. Using GAP \cite{gap}, we see that all groups with two generators $h_2$ and $h_3$ and four relations which are between $62$ cases of these $66$ cases are finite and solvable, that is a contradiction with the assumptions. So, there are just $4$ cases for the relations of these cycles which may lead to the existence of a subgraph isomorphic to the graph $C_4-C_6(-C_6--)(C_6---)$ in $K(\alpha,\beta)$. These cases are listed in table \ref{tab-C4-C6(-C6--)(C6---)}.  In the following, we show that all of these $4$ cases lead to contradictions and so, the graph $K(\alpha,\beta)$ contains no subgraph isomorphic to the graph $C_4-C_6(-C_6--)(C_6---)$.

\begin{enumerate}
\item[(1)]$ R_1:h_2h_3h_2^{-2}h_3=1$, $R_2:h_2(h_3h_2^{-1})^2h_2^{-1}h_3^2=1$, $R_3:h_3^2h_2^{-1}(h_3h_2^{-1}h_3)^2=1$,\\ $R_4:h_2(h_3h_2^{-1})^3h_3^{-1}h_2=1$:\\
$\Rightarrow\langle h_2,h_3\rangle\cong BS(1,-1)$ is solvable, a contradiction.

\item[(2)]$ R_1:h_2h_3h_2^{-2}h_3=1$, $R_2:h_2(h_3h_2^{-1})^2h_2^{-1}h_3^2=1$, $R_3:h_3^2h_2^{-1}(h_3h_2^{-1}h_3)^2=1$,\\ $R_4:h_2h_3h_2^{-1}h_3^{-3}h_2=1$:\\
$\Rightarrow\langle h_2,h_3\rangle\cong BS(1,-1)$ is solvable, a contradiction.

\item[(3)]$ R_1:h_2h_3^{-2}h_2h_3=1$, $R_2:h_2h_3(h_2h_3^{-1})^2h_3^{-1}h_2=1$, $R_3:h_2^2h_3^{-1}(h_2h_3^{-1}h_2)^2=1$,\\ $R_4:(h_2h_3^{-1})^3h_2^{-1}h_3^2=1$:\\
By interchanging $h_2$ and $h_3$ in (1) and with the same discussion, there is a contradiction.

\item[(4)]$ R_1:h_2h_3^{-2}h_2h_3=1$, $R_2:h_2h_3(h_2h_3^{-1})^2h_3^{-1}h_2=1$, $R_3:h_2^2h_3^{-1}(h_2h_3^{-1}h_2)^2=1$,\\ $R_4:h_2^2h_3h_2^{-1}h_3^{-2}h_2=1$:\\
By interchanging $h_2$ and $h_3$ in (2) and with the same discussion, there is a contradiction.

\end{enumerate}
\begin{table}[h]
\centering
\caption{The relations of a $C_4-C_6(-C_6--)(C_6---)$ in $K(\alpha,\beta)$}\label{tab-C4-C6(-C6--)(C6---)}
\begin{tabular}{|c|l|l|l|l|}\hline
$n$&$R_1$&$R_2$&$R_3$&$R_4$\\\hline
$1 $&$h_2h_3h_2^{-2}h_3=1$&$h_2(h_3h_2^{-1})^2h_2^{-1}h_3^2=1$&$h_3^2h_2^{-1}(h_3h_2^{-1}h_3)^2=1$&$h_2(h_3h_2^{-1})^3h_3^{-1}h_2=1$\\
$2 $&$h_2h_3h_2^{-2}h_3=1$&$h_2(h_3h_2^{-1})^2h_2^{-1}h_3^2=1$&$h_3^2h_2^{-1}(h_3h_2^{-1}h_3)^2=1$&$h_2h_3h_2^{-1}h_3^{-3}h_2=1$\\
$3$&$h_2h_3^{-2}h_2h_3=1$&$h_2h_3(h_2h_3^{-1})^2h_3^{-1}h_2=1$&$h_2^2h_3^{-1}(h_2h_3^{-1}h_2)^2=1$&$(h_2h_3^{-1})^3h_2^{-1}h_3^2=1$\\
$4$&$h_2h_3^{-2}h_2h_3=1$&$h_2h_3(h_2h_3^{-1})^2h_3^{-1}h_2=1$&$h_2^2h_3^{-1}(h_2h_3^{-1}h_2)^2=1$&$h_2^2h_3h_2^{-1}h_3^{-2}h_2=1$\\
\hline
\end{tabular}
\end{table}

\subsection{$\mathbf{C_4-C_6(-C_6--)(---C_4)}$}
By considering the relations from Tables \ref{tab-C4} and \ref{tab-C4-C6(-C6--)} which are not disproved, it can be seen that there is no case for the relations of two cycles $C_4$ and two cycles $C_6$ in this structure. It means that the graph $K(\alpha,\beta)$ contains no subgraph isomorphic to the graph $C_4-C_6(-C_6--)(---C_4)$.

$\mathbf{C_5-C_5}$ \textbf{subgraph:} With the same discussion such as about $C_5--C_5$, there are $447$ cases for the relations of the existence of two cycles of length $5$ in the graph $K(\alpha,\beta)$. Now suppose that there are two $C_5$ with one common edge in the graph $K(\alpha,\beta)$. Since this structure has two cycles of length $5$, the relations of these $C_5$ cycles must be between $447$ cases that have mentioned above.

Suppose that $[h_1',h_1'',h_2',h_2'',h_3',h_3'',h_4',h_4'',h_5',h_5'']$ and $[h_1',h_1'',h_6',h_6'',h_7',h_7'',h_8',h_8'',h_9',h_9'']$ are $10-$tuples related to the cycles $C_5$ in the graph $C_5-C_5$, where the first two components of these tuples are related to the common edge of $C_5$ and $C_5$. Without loss of generality we may assume that $h_1'=1$, where $h_1'',h_2',h_2'',h_3',h_3'',h_4',h_4'',h_5',h_5'',h_6',h_6'',h_7',h_7'',h_8',h_8'',h_9',h_9''\in supp(\alpha)$ and $\alpha=1+h_2+h_3$. With the same discussion such as about $K_{2,3}$, it is easy to see that $h_2' \neq h_6'$ and $h_5'' \neq h_9''$. 

With these assumptions and by the discussion above, it can be seen that there are $355$ cases which may lead to the existence of a subgraph isomorphic to the graph $C_5-C_5$ in the graph $K(\alpha,\beta)$. These cases are listed in table \ref{tab-C5-C5}. In the following, we show that $225$ cases of these relations lead to a contradiction and $130$ cases of them may lead to the  existence of a subgraph isomorphic to the graph $C_5-C_5$ in the graph $K(\alpha,\beta)$. Cases which are not disproved are marked by $*$s in the Table \ref{tab-C5-C5}.

\begin{center}
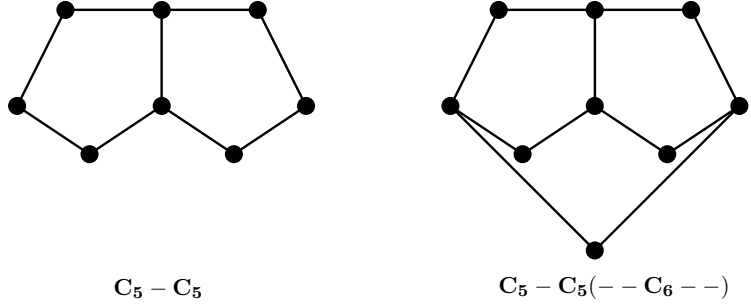
\begin{figure}
\psscalebox{0.8 0.8} 
{
\begin{pspicture}(0,-2.4985578)(12.394231,2.4985578)
\psdots[linecolor=black, dotsize=0.3](2.5971153,2.3014424)
\psdots[linecolor=black, dotsize=0.3](2.5971153,0.7014423)
\psdots[linecolor=black, dotsize=0.3](3.7971153,-0.0985577)
\psdots[linecolor=black, dotsize=0.3](4.997115,0.7014423)
\psdots[linecolor=black, dotsize=0.3](4.1971154,2.3014424)
\psdots[linecolor=black, dotsize=0.3](0.9971153,2.3014424)
\psdots[linecolor=black, dotsize=0.3](0.19711533,0.7014423)
\psdots[linecolor=black, dotsize=0.3](1.3971153,-0.0985577)
\psdots[linecolor=black, dotsize=0.3](9.797115,2.3014424)
\psdots[linecolor=black, dotsize=0.3](9.797115,0.7014423)
\psdots[linecolor=black, dotsize=0.3](10.997115,-0.0985577)
\psdots[linecolor=black, dotsize=0.3](12.197115,0.7014423)
\psdots[linecolor=black, dotsize=0.3](11.397116,2.3014424)
\psdots[linecolor=black, dotsize=0.3](8.197115,2.3014424)
\psdots[linecolor=black, dotsize=0.3](7.397115,0.7014423)
\psdots[linecolor=black, dotsize=0.3](8.5971155,-0.0985577)
\psdots[linecolor=black, dotsize=0.3](9.797115,-1.6985577)
\psline[linecolor=black, linewidth=0.04](2.5971153,0.7014423)(2.5971153,2.3014424)(4.1971154,2.3014424)(4.997115,0.7014423)(3.7971153,-0.0985577)(2.5971153,0.7014423)(1.3971153,-0.0985577)(0.19711533,0.7014423)(0.9971153,2.3014424)(2.5971153,2.3014424)
\psline[linecolor=black, linewidth=0.04](9.797115,2.3014424)(11.397116,2.3014424)(12.197115,0.7014423)(10.997115,-0.0985577)(9.797115,0.7014423)(8.5971155,-0.0985577)(7.397115,0.7014423)(8.197115,2.3014424)(9.797115,2.3014424)(9.797115,0.7014423)
\psline[linecolor=black, linewidth=0.04](12.197115,0.7014423)(9.797115,-1.6985577)(7.397115,0.7014423)
\rput[bl](1.7971153,-2.4985578){$\mathbf{C_5-C_5}$}
\rput[bl](8.197115,-2.4985578){$\mathbf{C_5-C_5(--C_6--)}$}
\end{pspicture}
}
\vspace{1cc}
\caption{The graph $C_5-C_5$  and a forbidden subgraph which contains it}\label{f-9}
\end{figure}
\end{center}

\begin{enumerate}

\item[(1)]$ R_1: h_2^3h_3^2=1$, $R_2: h_2^2h_3^{-1}h_2^{-2}h_3=1$:\\
$\Rightarrow\langle h_2,h_3\rangle\cong BS(1,1)$ is solvable, a contradiction.

\item[(2)]$ R_1: h_2^3h_3^2=1$, $R_2: h_2h_3^2h_2^{-1}h_3^2=1$:\\
$\Rightarrow h_2=1$ that is a contradiction.


\item[(4)]$ R_1: h_2^3h_3^2=1$, $R_2: (h_2h_3^{-1})^2(h_2^{-1}h_3)^2=1$:\\
$R_1 \Rightarrow h_3^2\in Z(G)$ where $G=\langle h_2,h_3\rangle$. Let $x=h_2h_3^{-1}$. So $h_3^{-1}h_2=x^{h_3}$. Also 
$ (h_2h_3^{-1})^2=(h_3^{-1}h_2)^2 \text{ so if } H=\langle x,x^{h_3}\rangle\Rightarrow H\cong BS(1,-1)$ is solvable. By Corollary \ref{n-sub2} $H \trianglelefteq G$ since $G=\langle x,h_3\rangle$. Since $\frac{G}{H}=\frac{\langle h_3\rangle H}{H}$ is a cyclic group, it is solvable. So $G$ is solvable, a contradiction.

\item[(5)]$ R_1: h_2^3h_3h_2^{-1}h_3=1$, $R_2: h_2^2h_3^{-1}h_2^{-2}h_3=1$:\\
$\Rightarrow \langle h_2,h_3\rangle=\langle h_2\rangle$ is abelian, a contradiction.

\item[(6)]$ R_1: h_2^3h_3h_2^{-1}h_3=1$, $R_2: h_2^2h_3^{-1}h_2^{-1}h_3^{-1}h_2=1$:\\
$\Rightarrow h_3=1$ that is a contradiction.

\item[(7)]$ R_1: h_2^3h_3h_2^{-1}h_3=1$, $R_2: h_2^2h_3^{-1}h_2^{-1}h_3^2=1$:\\
$\Rightarrow\langle h_2,h_3\rangle\cong BS(1,2)$ is solvable, a contradiction.

\item[(8)]$ R_1: h_2^3h_3h_2^{-1}h_3=1$, $R_2: h_2^2h_3^{-1}h_2^2h_3=1$:\\
$\Rightarrow h_2=1$ that is a contradiction.

\item[(9)]$ R_1: h_2^3h_3h_2^{-1}h_3=1$, $R_2: h_2h_3^2h_2^{-1}h_3^{-1}h_2=1$:\\
$\Rightarrow\langle h_2,h_3\rangle\cong BS(1,-2)$ is solvable, a contradiction.

\item[(10)]$ R_1: h_2^3h_3h_2^{-1}h_3=1$, $R_2: h_2h_3^2h_2^{-1}h_3^2=1$:\\
$\Rightarrow h_3=1$ that is a contradiction.

\item[(11)]$ R_1: h_2^3h_3h_2^{-1}h_3=1$, $R_2: h_2h_3^{-2}h_2^{-1}h_3^2=1$:\\
$\Rightarrow h_3^2\in Z(G)$ where $G=\langle h_2,h_3\rangle$. Let $x=h_2^{-1}h_3^{-1}$. So $h_3^{-1}h_2^{-1}=x^{h_3}$. Also 
$(h_3^{-1}h_2^{-1})^2=(h_2^{-1}h_3^{-1})^2 \text{ so if } H=\langle x,x^{h_3}\rangle\Rightarrow H\cong BS(1,-1)$ is solvable. By Corollary \ref{n-sub2} $H \trianglelefteq G$ since $G=\langle x,h_3\rangle$. Since $\frac{G}{H}=\frac{\langle h_3\rangle H}{H}$ is a cyclic group, it is solvable. So $G$ is solvable, a contradiction.



\item[(14)]$ R_1: h_2^3h_3h_2^{-1}h_3=1$, $R_2: h_2(h_3^{-1}h_2h_3^{-1})^2h_2=1$:\\
$\Rightarrow h_2=1$ that is a contradiction.

\item[(15)]$ R_1: h_2^3h_3h_2^{-1}h_3=1$, $R_2: (h_2h_3^{-1})^2h_2h_3h_2^{-1}h_3=1$:\\
$\Rightarrow\langle h_2,h_3\rangle\cong BS(1,-1)$ is solvable, a contradiction.

\item[(16)]$ R_1: h_2^3h_3^{-2}h_2=1$, $R_2: h_2^2h_3h_2^{-2}h_3=1$:\\
$\Rightarrow h_2=1$ that is a contradiction.

\item[(17)]$ R_1: h_2^3h_3^{-2}h_2=1$, $R_2: h_2^2h_3^{-1}h_2^{-2}h_3=1$:\\
$\Rightarrow \langle h_2,h_3\rangle=\langle h_2\rangle$ is abelian, a contradiction.

\item[(18)]$ R_1: h_2^3h_3^{-2}h_2=1$, $R_2: h_2^2h_3^{-1}h_2^2h_3=1$:\\
$\Rightarrow h_2=1$ that is a contradiction. 

\item[(19)]$ R_1: h_2^3h_3^{-2}h_2=1$, $R_2: h_2(h_2h_3^{-1})^2h_2^{-1}h_3=1$:\\
$\Rightarrow\langle h_2,h_3\rangle\cong BS(2,1)$ is solvable, a contradiction.

\item[(20)]$ R_1: h_2^3h_3^{-2}h_2=1$, $R_2: h_2h_3^2h_2^{-1}h_3^2=1$:\\
$\Rightarrow h_2=1$ that is a contradiction.

\item[(21)]$ R_1: h_2^3h_3^{-2}h_2=1$, $R_2: h_2h_3h_2^{-1}(h_3^{-1}h_2)^2=1$:\\
$\Rightarrow\langle h_2,h_3\rangle\cong BS(1,2)$ is solvable, a contradiction.


\item[(23)]$ R_1: h_2^3h_3^{-2}h_2=1$, $R_2: h_2h_3^{-2}h_2h_3^2=1$:\\
$\Rightarrow\langle h_2,h_3\rangle\cong BS(1,-1)$ is solvable, a contradiction.

\item[(24)]$ R_1: h_2^3h_3^{-2}h_2=1$, $R_2: (h_2h_3^{-1})^2(h_2^{-1}h_3)^2=1$:\\
$R_1 \Rightarrow h_3^2\in Z(G)$ where $G=\langle h_2,h_3\rangle$. Let $x=h_2h_3^{-1}$. So $h_3^{-1}h_2=x^{h_3}$. Also 
$ (h_2h_3^{-1})^2=(h_3^{-1}h_2)^2 \text{ so if } H=\langle x,x^{h_3}\rangle\Rightarrow H\cong BS(1,-1)$ is solvable. By Corollary \ref{n-sub2} $H \trianglelefteq G$ since $G=\langle x,h_3\rangle$. Since $\frac{G}{H}=\frac{\langle h_3\rangle H}{H}$ is a cyclic group, it is solvable. So $G$ is solvable, a contradiction.


\item[(26)]$ R_1: h_2^2(h_2h_3^{-1})^2h_2=1$, $R_2: h_2^2h_3^{-1}h_2^{-2}h_3=1$:\\
$\Rightarrow \langle h_2,h_3\rangle=\langle h_3h_2^{-2}\rangle$ is abelian, a contradiction.

\item[(27)]$ R_1: h_2^2(h_2h_3^{-1})^2h_2=1$, $R_2: h_2(h_3h_2^{-1}h_3)^2=1$:\\
$\Rightarrow h_2=1$ that is a contradiction.

\item[(28)]$ R_1: h_2^2(h_2h_3^{-1})^2h_2=1$, $R_2: h_2h_3^{-2}h_2^{-1}h_3^2=1$:\\
$\Rightarrow h_3^2\in Z(G)$ where $G=\langle h_2,h_3\rangle$. Let $x=h_2^{-1}h_3^{-1}$. So $h_3^{-1}h_2^{-1}=x^{h_3}$. Also 
$(h_3^{-1}h_2^{-1})^2=(h_2^{-1}h_3^{-1})^2 \text{ so if } H=\langle x,x^{h_3}\rangle\Rightarrow H\cong BS(1,-1)$ is solvable. By Corollary \ref{n-sub2} $H \trianglelefteq G$ since $G=\langle x,h_3\rangle$. Since $\frac{G}{H}=\frac{\langle h_3\rangle H}{H}$ is a cyclic group, it is solvable. So $G$ is solvable, a contradiction.



\item[(31)]$ R_1: h_2^2h_3^3=1$, $R_2: h_2^2h_3^{-1}h_2^2h_3=1$:\\
By interchanging $h_2$ and $h_3$ in (2) and with the same discussion, there is a contradiction.

\item[(32)]$ R_1: h_2^2h_3^3=1$, $R_2: h_2h_3^{-2}h_2^{-1}h_3^2=1$:\\
By interchanging $h_2$ and $h_3$ in (1) and with the same discussion, there is a contradiction.

\item[(33)]$ R_1: h_2^2h_3^3=1$, $R_2: (h_2h_3^{-1})^2(h_2^{-1}h_3)^2=1$:\\
By interchanging $h_2$ and $h_3$ in (4) and with the same discussion, there is a contradiction.


\item[(35)]$ R_1: h_2^2h_3^2h_2^{-1}h_3=1$, $R_2: h_2^2h_3^{-1}h_2^{-2}h_3=1$:\\
$\Rightarrow \langle h_2,h_3\rangle=\langle h_3\rangle$ is abelian, a contradiction.

\item[(36)]$ R_1: h_2^2h_3^2h_2^{-1}h_3=1$, $R_2: (h_2h_3)^2h_2^{-1}h_3=1$:\\
$\Rightarrow\langle h_2,h_3\rangle\cong BS(1,1)$ is solvable, a contradiction.

\item[(37)]$ R_1: h_2^2h_3^2h_2^{-1}h_3=1$, $R_2: h_2h_3^{-2}h_2^{-1}h_3^2=1$:\\
$\Rightarrow \langle h_2,h_3\rangle=\langle h_3\rangle$ is abelian, a contradiction.

\begin{center}
\tiny
\begin{longtable}{|c|l|l||c|l|l|}
\caption{{\footnotesize The relations of a $C_5-C_5$}}\label{tab-C5-C5}\\
\hline \multicolumn{1}{|c|}{$n$} & \multicolumn{1}{l|}{$R_1$} & \multicolumn{1}{l||}{$R_2$} & \multicolumn{1}{c|}{$n$} & \multicolumn{1}{l|}{$R_1$} & \multicolumn{1}{l|}{$R_2$}\\\hline
\endfirsthead
\multicolumn{6}{c}%
{{\tablename\ \thetable{} -- continued from previous page}} \\
\hline \multicolumn{1}{|c|}{$n$} & \multicolumn{1}{l|}{$R_1$} & \multicolumn{1}{l||}{$R_2$} & \multicolumn{1}{c|}{$n$} & \multicolumn{1}{l|}{$R_1$} & \multicolumn{1}{l|}{$R_2$}\\\hline
\endhead
\hline \multicolumn{6}{|r|}{{Continued on next page}} \\\hline
\endfoot
\hline
\endlastfoot
$1$&$h_2^3h_3^2=1$&$h_2^2h_3^{-1}h_2^{-2}h_3=1$&$61$&$h_2^2(h_3h_2^{-1})^2h_3=1$&$h_2h_3^2h_2^{-1}h_3^2=1$\\
$2$&$h_2^3h_3^2=1$&$h_2h_3^2h_2^{-1}h_3^2=1$&$62$&$h_2^2(h_3h_2^{-1})^2h_3=1$&$h_2(h_3h_2^{-1})^2h_3^{-1}h_2=1$\\
$3*$&$h_2^3h_3^2=1$&$h_2h_3^{-2}h_2^{-1}h_3^2=1$&$63*$&$h_2^2(h_3h_2^{-1})^2h_3=1$&$h_2h_3^{-1}h_2^{-1}h_3h_2^{-1}h_3^{-1}h_2=1$\\
$4$&$h_2^3h_3^2=1$&$(h_2h_3^{-1})^2(h_2^{-1}h_3)^2=1$&$64*$&$h_2^2(h_3h_2^{-1})^2h_3=1$&$h_2h_3^{-2}h_2^{-1}h_3^2=1$\\
$5$&$h_2^3h_3h_2^{-1}h_3=1$&$h_2^2h_3^{-1}h_2^{-2}h_3=1$&$65$&$h_2^2(h_3h_2^{-1})^2h_3=1$&$h_2h_3^{-1}h_2(h_3h_2^{-1})^2h_3=1$\\
$6$&$h_2^3h_3h_2^{-1}h_3=1$&$h_2^2h_3^{-1}h_2^{-1}h_3^{-1}h_2=1$&$66$&$h_2^2(h_3h_2^{-1})^2h_3=1$&$(h_2h_3^{-1})^2(h_2^{-1}h_3)^2=1$\\
$7$&$h_2^3h_3h_2^{-1}h_3=1$&$h_2^2h_3^{-1}h_2^{-1}h_3^2=1$&$67*$&$h_2^2(h_3h_2^{-1})^2h_3=1$&$(h_2h_3^{-1})^2h_2h_3^2=1$\\
$8$&$h_2^3h_3h_2^{-1}h_3=1$&$h_2^2h_3^{-1}h_2^2h_3=1$&$68*$&$h_2^2h_3^{-1}h_2^{-2}h_3=1$&$h_2^2h_3^{-1}h_2^{-2}h_3=1$\\
\end{longtable}
\end{center}

\item[(38)]$ R_1: h_2^2h_3^2h_2^{-1}h_3=1$, $R_2: (h_2h_3^{-1})^2(h_2^{-1}h_3)^2=1$:\\
$\Rightarrow\langle h_2,h_3\rangle\cong BS(1,1)$ is solvable, a contradiction.


\item[(40)]$ R_1: h_2^2h_3h_2^{-2}h_3=1$, $R_2: h_2^2(h_3h_2^{-1})^2h_3=1$:\\
$\Rightarrow h_3=1$ that is a contradiction.

\item[(41)]$ R_1: h_2^2h_3h_2^{-2}h_3=1$, $R_2: h_2^2h_3^{-1}h_2^{-2}h_3=1$:\\
$\Rightarrow h_3=1$ that is a contradiction.

\item[(42)]$ R_1: h_2^2h_3h_2^{-2}h_3=1$, $R_2: h_2^2h_3^{-1}(h_2^{-1}h_3)^2=1$:\\
$\Rightarrow\langle h_2,h_3\rangle\cong BS(1,2)$ is solvable, a contradiction.

\item[(43)]$ R_1: h_2^2h_3h_2^{-2}h_3=1$, $R_2: h_2^2h_3^{-1}h_2^2h_3=1$:\\
$\Rightarrow h_2=1$ that is a contradiction.

\item[(44)]$ R_1: h_2^2h_3h_2^{-2}h_3=1$, $R_2: h_2h_3^2h_2^{-1}h_3^2=1$:\\
$\Rightarrow h_3=1$ that is a contradiction.

\item[(45)]$ R_1: h_2^2h_3h_2^{-2}h_3=1$, $R_2: h_2(h_3h_2^{-1})^2h_3^{-1}h_2=1$:\\
$\Rightarrow\langle h_2,h_3\rangle\cong BS(2,1)$ is solvable, a contradiction.

\item[(46)]$ R_1: h_2^2h_3h_2^{-2}h_3=1$, $R_2: h_2h_3^{-1}h_2^{-1}h_3h_2^{-1}h_3^{-1}h_2=1$:\\
$\Rightarrow h_3=1$ that is a contradiction.

\item[(47)]$ R_1: h_2^2h_3h_2^{-2}h_3=1$, $R_2: h_2h_3^{-2}h_2^{-1}h_3^2=1$:\\
$\Rightarrow h_3=1$ that is a contradiction.

\item[(48)]$ R_1: h_2^2h_3h_2^{-2}h_3=1$, $R_2: h_2h_3^{-4}h_2=1$:\\
By interchanging $h_2$ and $h_3$ in (23) and with the same discussion, there is a contradiction.

\item[(49)]$ R_1: h_2^2h_3h_2^{-2}h_3=1$, $R_2: h_2h_3^{-2}h_2h_3^2=1$:\\
$\Rightarrow h_3=1$ that is a contradiction.



\item[(52)]$ R_1: h_2^2h_3h_2^{-2}h_3=1$, $R_2: h_2(h_3^{-1}h_2h_3^{-1})^2h_2=1$:\\
$\Rightarrow h_3=1$ that is a contradiction.


\item[(54)]$ R_1: h_2^2h_3h_2^{-1}h_3^2=1$, $R_2: h_2^2h_3^{-1}h_2^{-2}h_3=1$:\\
$\Rightarrow \langle h_2,h_3\rangle=\langle h_3\rangle$ is abelian, a contradiction.

\item[(55)]$ R_1: h_2^2h_3h_2^{-1}h_3^2=1$, $R_2: (h_2h_3)^2h_2^{-1}h_3=1$:\\
$\Rightarrow\langle h_2,h_3\rangle\cong BS(1,1)$ is solvable, a contradiction.

\item[(56)]$ R_1: h_2^2h_3h_2^{-1}h_3^2=1$, $R_2: h_2h_3^{-2}h_2^{-1}h_3^2=1$:\\
$\Rightarrow \langle h_2,h_3\rangle=\langle h_3\rangle$ is abelian, a contradiction.

\item[(57)]$ R_1: h_2^2h_3h_2^{-1}h_3^2=1$, $R_2: (h_2h_3^{-1})^2(h_2^{-1}h_3)^2=1$:\\
$\Rightarrow\langle h_2,h_3\rangle\cong BS(1,1)$ is solvable, a contradiction.

\item[(58)]$ R_1: h_2^2(h_3h_2^{-1})^2h_3=1$, $R_2: h_2^2h_3^{-1}h_2^{-2}h_3=1$:\\
$\Rightarrow\langle h_2,h_3\rangle\cong BS(-2,1)$ is solvable, a contradiction.

\item[(59)]$ R_1: h_2^2(h_3h_2^{-1})^2h_3=1$, $R_2: h_2^2h_3^{-1}(h_2^{-1}h_3)^2=1$:\\
$\Rightarrow h_3=1$ that is a contradiction.


\item[(61)]$ R_1: h_2^2(h_3h_2^{-1})^2h_3=1$, $R_2: h_2h_3^2h_2^{-1}h_3^2=1$:\\
$\Rightarrow h_3=1$ that is a contradiction.

\item[(62)]$ R_1: h_2^2(h_3h_2^{-1})^2h_3=1$, $R_2: h_2(h_3h_2^{-1})^2h_3^{-1}h_2=1$:\\
$\Rightarrow h_3=1$ that is a contradiction.



\item[(65)]$ R_1: h_2^2(h_3h_2^{-1})^2h_3=1$, $R_2: h_2h_3^{-1}h_2(h_3h_2^{-1})^2h_3=1$:\\
$\Rightarrow h_3=1$ that is a contradiction.

\item[(66)]$ R_1: h_2^2(h_3h_2^{-1})^2h_3=1$, $R_2: (h_2h_3^{-1})^2(h_2^{-1}h_3)^2=1$:\\
$\Rightarrow\langle h_2,h_3\rangle\cong BS(1,-2)$ is solvable, a contradiction.



\item[(69)]$ R_1: h_2^2h_3^{-1}h_2^{-2}h_3=1$, $R_2: h_2^2h_3^{-1}h_2^{-1}h_3^{-1}h_2=1$:\\
$\Rightarrow \langle h_2,h_3\rangle=\langle h_2\rangle$ is abelian, a contradiction.

\item[(70)]$ R_1: h_2^2h_3^{-1}h_2^{-2}h_3=1$, $R_2: h_2^2h_3^{-1}h_2^{-1}h_3^2=1$:\\
$\Rightarrow \langle h_2,h_3\rangle=\langle h_2\rangle$ is abelian, a contradiction.

\item[(71)]$ R_1: h_2^2h_3^{-1}h_2^{-2}h_3=1$, $R_2: h_2^2h_3^{-1}(h_2^{-1}h_3)^2=1$:\\
$\Rightarrow h_3=1$ that is a contradiction.

\item[(72)]$ R_1: h_2^2h_3^{-1}h_2^{-2}h_3=1$, $R_2: h_2^2h_3^{-2}h_2^{-1}h_3=1$:\\
$\Rightarrow \langle h_2,h_3\rangle=\langle h_2\rangle$ is abelian, a contradiction.

\item[(73)]$ R_1: h_2^2h_3^{-1}h_2^{-2}h_3=1$, $R_2: h_2^2h_3^{-3}h_2=1$:\\
$\Rightarrow\langle h_2,h_3\rangle\cong BS(1,1)$ is solvable, a contradiction.

\item[(74)]$ R_1: h_2^2h_3^{-1}h_2^{-2}h_3=1$, $R_2: h_2^2h_3^{-2}h_2h_3=1$:\\
$\Rightarrow \langle h_2,h_3\rangle=\langle h_2\rangle$ is abelian, a contradiction.

\item[(75)]$ R_1: h_2^2h_3^{-1}h_2^{-2}h_3=1$, $R_2: h_2^2h_3^{-1}(h_3^{-1}h_2)^2=1$:\\
$\Rightarrow \langle h_2,h_3\rangle=\langle h_2^{-1}h_3\rangle$ is abelian, a contradiction.

\item[(76)]$ R_1: h_2^2h_3^{-1}h_2^{-2}h_3=1$, $R_2: h_2^2h_3^{-1}h_2^2h_3=1$:\\
$\Rightarrow h_2=1$ that is a contradiction.

\item[(77)]$ R_1: h_2^2h_3^{-1}h_2^{-2}h_3=1$, $R_2: h_2^2h_3^{-1}h_2h_3^2=1$:\\
By interchanging $h_2$ and $h_3$ in (37) and with the same discussion, there is a contradiction.

\item[(78)]$ R_1: h_2^2h_3^{-1}h_2^{-2}h_3=1$, $R_2: h_2^2h_3^{-1}h_2h_3h_2^{-1}h_3=1$:\\
$\Rightarrow \langle h_2,h_3\rangle=\langle h_2\rangle$ is abelian, a contradiction.

\item[(79)]$ R_1: h_2^2h_3^{-1}h_2^{-2}h_3=1$, $R_2: h_2(h_2h_3^{-1})^2h_2^{-1}h_3=1$:\\
$\Rightarrow \langle h_2,h_3\rangle=\langle h_2\rangle$ is abelian, a contradiction.

\item[(80)]$ R_1: h_2^2h_3^{-1}h_2^{-2}h_3=1$, $R_2: h_2(h_2h_3^{-1})^2h_3^{-1}h_2=1$:\\
$\Rightarrow \langle h_2,h_3\rangle=\langle h_3\rangle$ is abelian, a contradiction.

\item[(81)]$ R_1: h_2^2h_3^{-1}h_2^{-2}h_3=1$, $R_2: h_2(h_2h_3^{-1})^2h_2h_3=1$:\\
$\Rightarrow \langle h_2,h_3\rangle=\langle h_2\rangle$ is abelian, a contradiction.



\item[(84)]$ R_1: h_2^2h_3^{-1}h_2^{-2}h_3=1$, $R_2: h_2h_3h_2h_3^{-2}h_2=1$:\\
$\Rightarrow \langle h_2,h_3\rangle=\langle h_2\rangle$ is abelian, a contradiction.

\item[(85)]$ R_1: h_2^2h_3^{-1}h_2^{-2}h_3=1$, $R_2: h_2h_3h_2h_3^{-1}h_2h_3=1$:\\
$\Rightarrow\langle h_2,h_3\rangle\cong BS(1,-1)$ is solvable, a contradiction.

\item[(86)]$ R_1: h_2^2h_3^{-1}h_2^{-2}h_3=1$, $R_2: h_2h_3(h_2h_3^{-1})^2h_2=1$:\\
$\Rightarrow \langle h_2,h_3\rangle=\langle h_2\rangle$ is abelian, a contradiction.

\item[(87)]$ R_1: h_2^2h_3^{-1}h_2^{-2}h_3=1$, $R_2: h_2h_3^2h_2h_3^{-1}h_2=1$:\\
By interchanging $h_2$ and $h_3$ in (56) and with the same discussion, there is a contradiction.

\item[(88)]$ R_1: h_2^2h_3^{-1}h_2^{-2}h_3=1$, $R_2: h_2h_3^2h_2^{-2}h_3=1$:\\
$\Rightarrow \langle h_2,h_3\rangle=\langle h_3\rangle$ is abelian, a contradiction.

\item[(89)]$ R_1: h_2^2h_3^{-1}h_2^{-2}h_3=1$, $R_2: h_2h_3^2h_2^{-1}h_3^{-1}h_2=1$:\\
$\Rightarrow \langle h_2,h_3\rangle=\langle h_2\rangle$ is abelian, a contradiction.

\item[(90)]$ R_1: h_2^2h_3^{-1}h_2^{-2}h_3=1$, $R_2: h_2h_3^2h_2^{-1}h_3^2=1$:\\
Using Tietze transformation where $h_3\mapsto h_3h_2$ and $h_2\mapsto h_2$, we have: $R_1\Rightarrow h_2^2\in Z(G)$ where $G=\langle h_2,h_3\rangle$. Let $x=h_2h_3$. So $h_2^{-1}h_3^{-1}=(x^{h_2})^{-1}$. Also 
$(h_2h_3)^2=(h_2^{-1}h_3^{-1})^2 \text{ so } H=\langle x,(x^{h_2})^{-1}\rangle=\langle x,x^{h_2}\rangle\cong BS(1,-1)$ is solvable. By Corollary \ref{n-sub2} $H \trianglelefteq G$ since $G=\langle x,h_2\rangle$. Since $\frac{G}{H}=\frac{\langle h_2\rangle H}{H}$ is a cyclic group, it is solvable. So $G$ is solvable, a contradiction.

\item[(91)]$ R_1: h_2^2h_3^{-1}h_2^{-2}h_3=1$, $R_2: h_2h_3(h_3h_2^{-1})^2h_3=1$:\\
$\Rightarrow\langle h_2,h_3\rangle\cong BS(1,1)$ is solvable, a contradiction.

\item[(92)]$ R_1: h_2^2h_3^{-1}h_2^{-2}h_3=1$, $R_2: h_2h_3h_2^{-2}h_3^2=1$:\\
$\Rightarrow \langle h_2,h_3\rangle=\langle h_3\rangle$ is abelian, a contradiction.

\item[(93)]$ R_1: h_2^2h_3^{-1}h_2^{-2}h_3=1$, $R_2: h_2h_3h_2^{-1}(h_2^{-1}h_3)^2=1$:\\
$\Rightarrow \langle h_2,h_3\rangle=\langle h_3^{-1}h_2\rangle$ is abelian, a contradiction.

\item[(94)]$ R_1: h_2^2h_3^{-1}h_2^{-2}h_3=1$, $R_2: h_2h_3h_2^{-1}h_3^{-1}h_2^{-1}h_3=1$:\\
$\Rightarrow\langle h_2,h_3\rangle\cong BS(1,-1)$ is solvable, a contradiction.

\item[(95)]$ R_1: h_2^2h_3^{-1}h_2^{-2}h_3=1$, $R_2: h_2h_3h_2^{-1}h_3^{-2}h_2=1$:\\
$\Rightarrow \langle h_2,h_3\rangle=\langle h_2\rangle$ is abelian, a contradiction.

\item[(96)]$ R_1: h_2^2h_3^{-1}h_2^{-2}h_3=1$, $R_2: h_2h_3h_2^{-1}(h_3^{-1}h_2)^2=1$:\\
$\Rightarrow \langle h_2,h_3\rangle=\langle h_2\rangle$ is abelian, a contradiction.

\item[(97)]$ R_1: h_2^2h_3^{-1}h_2^{-2}h_3=1$, $R_2: h_2h_3h_2^{-1}h_3h_2h_3^{-1}h_2=1$:\\
$\Rightarrow \langle h_2,h_3\rangle=\langle h_2\rangle$ is abelian, a contradiction.


\item[(99)]$ R_1: h_2^2h_3^{-1}h_2^{-2}h_3=1$, $R_2: h_2(h_3h_2^{-1})^2h_2^{-1}h_3=1$:\\
$\Rightarrow \langle h_2,h_3\rangle=\langle h_2h_3^{-1}\rangle$ is abelian, a contradiction.

\item[(100)]$ R_1: h_2^2h_3^{-1}h_2^{-2}h_3=1$, $R_2: h_2(h_3h_2^{-1})^2h_3^{-1}h_2=1$:\\
$\Rightarrow h_3=1$ that is a contradiction.



\item[(103)]$ R_1: h_2^2h_3^{-1}h_2^{-2}h_3=1$, $R_2: h_2h_3^{-1}h_2^{-1}h_3h_2^{-1}h_3^{-1}h_2=1$:\\
$\Rightarrow\langle h_2,h_3\rangle\cong BS(1,-1)$ is solvable, a contradiction.


\item[(105)]$ R_1: h_2^2h_3^{-1}h_2^{-2}h_3=1$, $R_2: h_2h_3^{-2}h_2^{-1}h_3^2=1$:\\
Using Tietze transformation where $h_3\mapsto h_3h_2^{-1}$ and $h_2\mapsto h_2$, we have: $R_1\Rightarrow h_2^2\in Z(G)$ where $G=\langle h_2,h_3\rangle$. Let $x=h_3h_2^{-1}$. So $h_2^{-1}h_3=x^{h_2}$. Also 
$(h_3h_2^{-1})^2=(h_2^{-1}h_3)^2 \text{ so if } H=\langle x,x^{h_2}\rangle\Rightarrow H\cong BS(1,-1)$ is solvable. By Corollary \ref{n-sub2} $H \trianglelefteq G$ since $G=\langle x,h_2\rangle$. Since $\frac{G}{H}=\frac{\langle h_2\rangle H}{H}$ is a cyclic group, it is solvable. So $G$ is solvable, a contradiction.

\item[(106)]$ R_1: h_2^2h_3^{-1}h_2^{-2}h_3=1$, $R_2: h_2h_3^{-2}(h_2^{-1}h_3)^2=1$:\\
$\Rightarrow h_3=1$ that is a contradiction.




\item[(110)]$ R_1: h_2^2h_3^{-1}h_2^{-2}h_3=1$, $R_2: h_2h_3^{-2}h_2h_3^2=1$:\\
By interchanging $h_2$ and $h_3$ in (47) and with the same discussion, there is a contradiction.

\item[(111)]$ R_1: h_2^2h_3^{-1}h_2^{-2}h_3=1$, $R_2: h_2h_3^{-2}h_2h_3h_2^{-1}h_3=1$:\\
$\Rightarrow h_2=1$ that is a contradiction.


\item[(113)]$ R_1: h_2^2h_3^{-1}h_2^{-2}h_3=1$, $R_2: (h_2h_3^{-1}h_2)^2h_3=1$:\\
$\Rightarrow\langle h_2,h_3\rangle\cong BS(1,-1)$ is solvable, a contradiction.

\item[(114)]$ R_1: h_2^2h_3^{-1}h_2^{-2}h_3=1$, $R_2: h_2h_3^{-1}h_2h_3^3=1$:\\
By interchanging $h_2$ and $h_3$ in (11) and with the same discussion, there is a contradiction.





\item[(119)]$ R_1: h_2^2h_3^{-1}h_2^{-2}h_3=1$, $R_2: (h_2h_3^{-1})^2(h_2^{-1}h_3)^2=1$:\\
$\Rightarrow h_2^2\in Z(G)$ where $G=\langle h_2,h_3\rangle$. Let $x=h_3h_2^{-1}$. So $h_2^{-1}h_3=x^{h_2}$. Also 
$(h_3h_2^{-1})^2=(h_2^{-1}h_3)^2 \text{ so if } H=\langle x,x^{h_2}\rangle\Rightarrow H\cong BS(1,-1)$ is solvable. By Corollary \ref{n-sub2} $H \trianglelefteq G$ since $G=\langle x,h_2\rangle$. Since $\frac{G}{H}=\frac{\langle h_2\rangle H}{H}$ is a cyclic group, it is solvable. So $G$ is solvable, a contradiction.







\item[(126)]$ R_1: h_2^2h_3^{-1}h_2^{-2}h_3=1$, $R_2: (h_2h_3^{-1})^3h_2h_3=1$:\\

\item[(127)]$ R_1: h_2^2h_3^{-1}h_2^{-2}h_3=1$, $R_2: h_3^2(h_3h_2^{-1})^2h_3=1$:\\
By interchanging $h_2$ and $h_3$ in (28) and with the same discussion, there is a contradiction.

\item[(128)]$ R_1: h_2^2h_3^{-1}h_2^{-2}h_3=1$, $R_2: h_3(h_3h_2^{-1})^3h_3=1$:\\
$\Rightarrow h_3^2\in Z(G)$ where $G=\langle h_2,h_3\rangle$. Let $x=h_2^{-1}h_3^{-1}$. So $h_3^{-1}h_2^{-1}=x^{h_3}$. Also 
$(h_3^{-1}h_2^{-1})^2=(h_2^{-1}h_3^{-1})^2 \text{ so if } H=\langle x,x^{h_3}\rangle\Rightarrow H\cong BS(1,-1)$ is solvable. By Corollary \ref{n-sub2} $H \trianglelefteq G$ since $G=\langle x,h_3\rangle$. Since $\frac{G}{H}=\frac{\langle h_3\rangle H}{H}$ is a cyclic group, it is solvable. So $G$ is solvable, a contradiction.



\item[(131)]$ R_1: h_2^2h_3^{-1}h_2^{-1}h_3^{-1}h_2=1$, $R_2: h_2^2h_3^{-3}h_2=1$:\\
$\Rightarrow \langle h_2,h_3\rangle=\langle h_2\rangle$ is abelian, a contradiction.

\item[(132)]$ R_1: h_2^2h_3^{-1}h_2^{-1}h_3^{-1}h_2=1$, $R_2: h_2^2h_3^{-1}h_2^2h_3=1$:\\
$\Rightarrow h_3=1$ that is a contradiction.

\item[(133)]$ R_1: h_2^2h_3^{-1}h_2^{-1}h_3^{-1}h_2=1$, $R_2: h_2h_3^2h_2^{-1}h_3^2=1$:\\
$\Rightarrow h_3=1$ that is a contradiction.

\item[(134)]$ R_1: h_2^2h_3^{-1}h_2^{-1}h_3^{-1}h_2=1$, $R_2: h_2h_3h_2^{-1}h_3^{-1}h_2^{-1}h_3=1$:\\
$\Rightarrow \langle h_2,h_3\rangle=\langle h_2\rangle$ is abelian, a contradiction.


\item[(136)]$ R_1: h_2^2h_3^{-1}h_2^{-1}h_3^{-1}h_2=1$, $R_2: h_2(h_3h_2^{-1})^3h_3=1$:\\
$\Rightarrow h_3=1$ that is a contradiction.



\item[(139)]$ R_1: h_2^2h_3^{-1}h_2^{-1}h_3^{-1}h_2=1$, $R_2: (h_2h_3^{-1})^2(h_2^{-1}h_3)^2=1$:\\
$\Rightarrow h_2=1$ that is a contradiction.




\item[(143)]$ R_1: h_2^2h_3^{-1}h_2^{-1}h_3^2=1$, $R_2: h_2h_3^{-2}h_2^{-1}h_3^2=1$:\\
By interchanging $h_2$ and $h_3$ in (89) and with the same discussion, there is a contradiction.

\item[(144)]$ R_1: h_2^2h_3^{-1}h_2^{-1}h_3^2=1$, $R_2: h_2h_3^{-1}h_2h_3^3=1$:\\
By interchanging $h_2$ and $h_3$ in (9) and with the same discussion, there is a contradiction.




\item[(148)]$ R_1: h_2^2h_3^{-1}(h_2^{-1}h_3)^2=1$, $R_2: h_2h_3^{-1}h_2^{-1}h_3h_2^{-1}h_3^{-1}h_2=1$:\\
$\Rightarrow h_3=1$ that is a contradiction.

\item[(149)]$ R_1: h_2^2h_3^{-1}(h_2^{-1}h_3)^2=1$, $R_2: h_2h_3^{-2}h_2^{-1}h_3^2=1$:\\
By interchanging $h_2$ and $h_3$ in (106) and with the same discussion, there is a contradiction.


\item[(151)]$ R_1: h_2^2h_3^{-2}h_2^{-1}h_3=1$, $R_2: h_2^2h_3^{-3}h_2=1$:\\
$\Rightarrow \langle h_2,h_3\rangle=\langle h_2\rangle$ is abelian, a contradiction.

\item[(152)]$ R_1: h_2^2h_3^{-2}h_2^{-1}h_3=1$, $R_2: h_2h_3h_2^{-1}h_3^{-1}h_2^{-1}h_3=1$:\\
$\Rightarrow \langle h_2,h_3\rangle=\langle h_3\rangle$ is abelian, a contradiction.

\item[(153)]$ R_1: h_2^2h_3^{-2}h_2^{-1}h_3=1$, $R_2: h_2h_3^{-2}h_2^{-1}h_3^2=1$:\\
By interchanging $h_2$ and $h_3$ in (72) and with the same discussion, there is a contradiction.



\item[(156)]$ R_1: h_2^2h_3^{-3}h_2=1$, $R_2: h_2h_3h_2^{-1}h_3^{-1}h_2^{-1}h_3=1$:\\
$\Rightarrow \langle h_2,h_3\rangle=\langle h_2^{-1}h_3\rangle$ is abelian, a contradiction.

\item[(157)]$ R_1: h_2^2h_3^{-3}h_2=1$, $R_2: h_2h_3h_2^{-1}h_3^{-2}h_2=1$:\\
$\Rightarrow \langle h_2,h_3\rangle=\langle h_2\rangle$ is abelian, a contradiction.

\item[(158)]$ R_1: h_2^2h_3^{-3}h_2=1$, $R_2: h_2h_3^{-2}h_2^{-1}h_3^2=1$:\\
By interchanging $h_2$ and $h_3$ in (73) and with the same discussion, there is a contradiction.

\item[(159)]$ R_1: h_2^2h_3^{-3}h_2=1$, $R_2: h_2h_3^{-3}h_2h_3=1$:\\
By interchanging $h_2$ and $h_3$ in (131) and with the same discussion, there is a contradiction.




\item[(163)]$ R_1: h_2^2h_3^{-2}h_2h_3=1$, $R_2: h_2h_3h_2h_3^{-2}h_2=1$:\\
$\Rightarrow\langle h_2,h_3\rangle\cong BS(1,1)$ is solvable, a contradiction.

\item[(164)]$ R_1: h_2^2h_3^{-2}h_2h_3=1$, $R_2: h_2h_3^{-2}h_2^{-1}h_3^2=1$:\\
By interchanging $h_2$ and $h_3$ in (88) and with the same discussion, there is a contradiction.


\item[(166)]$ R_1: h_2^2h_3^{-1}(h_3^{-1}h_2)^2=1$, $R_2: (h_2h_3^{-1})^2(h_2^{-1}h_3)^2=1$:\\
$\Rightarrow \langle h_2,h_3\rangle=\langle h_3h_2^{-1}\rangle$ is abelian, a contradiction.


\item[(168)]$ R_1: h_2^2h_3^{-1}h_2^2h_3=1$, $R_2: h_2(h_2h_3^{-1})^3h_2=1$:\\
$\Rightarrow h_3=1$ that is a contradiction.

\item[(169)]$ R_1: h_2^2h_3^{-1}h_2^2h_3=1$, $R_2: (h_2h_3)^2h_2^{-1}h_3=1$:\\
$\Rightarrow h_3=1$ that is a contradiction.


\item[(171)]$ R_1: h_2^2h_3^{-1}h_2^2h_3=1$, $R_2: h_2(h_3h_2^{-1})^3h_3=1$:\\
$\Rightarrow h_2=1$ that is a contradiction.

\item[(172)]$ R_1: h_2^2h_3^{-1}h_2^2h_3=1$, $R_2: h_2h_3^{-2}h_2^{-1}h_3^2=1$:\\
By interchanging $h_2$ and $h_3$ in (90) and with the same discussion, there is a contradiction.


\item[(174)]$ R_1: h_2^2h_3^{-1}h_2^2h_3=1$, $R_2: h_2h_3^{-4}h_2=1$:\\
By interchanging $h_2$ and $h_3$ in (20) and with the same discussion, there is a contradiction.

\item[(175)]$ R_1: h_2^2h_3^{-1}h_2^2h_3=1$, $R_2: h_2h_3^{-3}h_2h_3=1$:\\
By interchanging $h_2$ and $h_3$ in (133) and with the same discussion, there is a contradiction.

\item[(176)]$ R_1: h_2^2h_3^{-1}h_2^2h_3=1$, $R_2: h_2h_3^{-2}h_2h_3^2=1$:\\
By interchanging $h_2$ and $h_3$ in (44) and with the same discussion, there is a contradiction.


\item[(178)]$ R_1: h_2^2h_3^{-1}h_2^2h_3=1$, $R_2: h_2h_3^{-1}h_2h_3^3=1$:\\
By interchanging $h_2$ and $h_3$ in (10) and with the same discussion, there is a contradiction.


\item[(180)]$ R_1: h_2^2h_3^{-1}h_2^2h_3=1$, $R_2: (h_2h_3^{-1})^2h_2^{-1}h_3^2=1$:\\
$\Rightarrow h_3=1$ that is a contradiction.

\item[(181)]$ R_1: h_2^2h_3^{-1}h_2^2h_3=1$, $R_2: (h_2h_3^{-1})^2(h_2^{-1}h_3)^2=1$:\\
Using Tietze transformation where $h_3\mapsto h_2h_3$ and $h_2\mapsto h_2$, we have $ R_1: h_2^2h_3^{-1}h_2^2h_3=1$ and $R_2: h_2h_3^{-2}h_2^{-1}h_3^2=1$. So with the same discussion such as item (172), there is a contradiction.


\item[(183)]$ R_1: h_2^2h_3^{-1}h_2^2h_3=1$, $R_2: (h_2h_3^{-1})^2h_2h_3^2=1$:\\
By interchanging $h_2$ and $h_3$ in (61) and with the same discussion, there is a contradiction.

\item[(184)]$ R_1: h_2^2h_3^{-1}h_2^2h_3=1$, $R_2: (h_2h_3^{-1})^2h_2h_3h_2^{-1}h_3=1$:\\
$\Rightarrow h_3=1$ that is a contradiction.

\item[(185)]$ R_1: h_2^2h_3^{-1}h_2^2h_3=1$, $R_2: (h_2h_3^{-1})^3h_2h_3=1$:\\
$\Rightarrow h_3=1$ that is a contradiction.


\item[(187)]$ R_1: h_2^2h_3^{-1}h_2h_3^2=1$, $R_2: h_2h_3h_2h_3^{-1}h_2h_3=1$:\\
By interchanging $h_2$ and $h_3$ in (36) and with the same discussion, there is a contradiction.

\item[(188)]$ R_1: h_2^2h_3^{-1}h_2h_3^2=1$, $R_2: h_2h_3^{-2}h_2^{-1}h_3^2=1$:\\
By interchanging $h_2$ and $h_3$ in (35) and with the same discussion, there is a contradiction.

\item[(189)]$ R_1: h_2^2h_3^{-1}h_2h_3^2=1$, $R_2: (h_2h_3^{-1})^2(h_2^{-1}h_3)^2=1$:\\
By interchanging $h_2$ and $h_3$ in (38) and with the same discussion, there is a contradiction.




\item[(193)]$ R_1: h_2(h_2h_3^{-1})^2h_2^{-1}h_3=1$, $R_2: h_2h_3^{-1}h_2^{-1}h_3h_2h_3^{-1}h_2=1$:\\
$\Rightarrow\langle h_2,h_3\rangle\cong BS(1,1)$ is solvable, a contradiction.



\item[(196)]$ R_1: h_2(h_2h_3^{-1})^2h_3^{-1}h_2=1$, $R_2: (h_2h_3^{-1})^2(h_2^{-1}h_3)^2=1$:\\
$\Rightarrow \langle h_2,h_3\rangle=\langle h_2^{-1} h_3\rangle$ is abelian, a contradiction.


\item[(198)]$ R_1: h_2(h_2h_3^{-1})^2h_2h_3=1$, $R_2: h_2h_3^{-2}h_2^{-1}h_3^2=1$:\\
By interchanging $h_2$ and $h_3$ in (91) and with the same discussion, there is a contradiction.

\item[(199)]$ R_1: h_2(h_2h_3^{-1})^2h_2h_3=1$, $R_2: (h_2h_3^{-1}h_2)^2h_3=1$:\\
$\Rightarrow\langle h_2,h_3\rangle\cong BS(1,1)$ is solvable, a contradiction.

\item[(200)]$ R_1: h_2(h_2h_3^{-1})^3h_2=1$, $R_2: (h_2h_3)^2h_2^{-1}h_3=1$:\\
$\Rightarrow h_3=1$ that is a contradiction.

\item[(201)]$ R_1: h_2(h_2h_3^{-1})^3h_2=1$, $R_2: h_2h_3^{-2}h_2^{-1}h_3^2=1$:\\
By interchanging $h_2$ and $h_3$ in (128) and with the same discussion, there is a contradiction.

\item[(202)]$ R_1: h_2(h_2h_3^{-1})^3h_2=1$, $R_2: (h_2h_3^{-1})^2(h_2^{-1}h_3)^2=1$:\\
$\Rightarrow \langle h_2,h_3\rangle=\langle h_3^{-1} h_2\rangle$ is abelian, a contradiction.

\item[(203)]$ R_1: (h_2h_3)^2h_2^{-1}h_3=1$, $R_2: h_2h_3^{-2}h_2^{-1}h_3^2=1$:\\
By interchanging $h_2$ and $h_3$ in (85) and with the same discussion, there is a contradiction.


\item[(205)]$ R_1: h_2h_3h_2h_3^{-2}h_2=1$, $R_2: h_2h_3^{-2}h_2^{-1}h_3^2=1$:\\
By interchanging $h_2$ and $h_3$ in (92) and with the same discussion, there is a contradiction.


\item[(207)]$ R_1: h_2h_3h_2h_3^{-1}h_2h_3=1$, $R_2: h_2h_3^2h_2h_3^{-1}h_2=1$:\\
By interchanging $h_2$ and $h_3$ in (55) and with the same discussion, there is a contradiction.

\item[(208)]$ R_1: h_2h_3h_2h_3^{-1}h_2h_3=1$, $R_2: h_2h_3^2h_2^{-1}h_3^2=1$:\\
By interchanging $h_2$ and $h_3$ in (169) and with the same discussion, there is a contradiction.


\item[(210)]$ R_1: h_2h_3h_2h_3^{-1}h_2h_3=1$, $R_2: h_3(h_3h_2^{-1})^3h_3=1$:\\
By interchanging $h_2$ and $h_3$ in (200) and with the same discussion, there is a contradiction.



\item[(213)]$ R_1: h_2h_3(h_2h_3^{-1})^2h_2=1$, $R_2: (h_2h_3^{-1}h_2)^2h_3=1$:\\
$\Rightarrow\langle h_2,h_3\rangle\cong BS(1,1)$ is solvable, a contradiction.

\item[(214)]$ R_1: h_2h_3(h_2h_3^{-1})^2h_2=1$, $R_2: (h_2h_3^{-1})^2(h_2^{-1}h_3)^2=1$:\\
$\Rightarrow \langle h_2,h_3\rangle=\langle h_2\rangle$ is abelian, a contradiction.


\item[(216)]$ R_1: h_2h_3^2h_2h_3^{-1}h_2=1$, $R_2: h_2h_3^{-2}h_2^{-1}h_3^2=1$:\\
By interchanging $h_2$ and $h_3$ in (54) and with the same discussion, there is a contradiction.

\item[(217)]$ R_1: h_2h_3^2h_2h_3^{-1}h_2=1$, $R_2: (h_2h_3^{-1})^2(h_2^{-1}h_3)^2=1$:\\
By interchanging $h_2$ and $h_3$ in (57) and with the same discussion, there is a contradiction.


\item[(219)]$ R_1: h_2h_3^2h_2^{-2}h_3=1$, $R_2: h_2h_3h_2^{-2}h_3^2=1$:\\
By interchanging $h_2$ and $h_3$ in (163) and with the same discussion, there is a contradiction.

\item[(220)]$ R_1: h_2h_3^2h_2^{-2}h_3=1$, $R_2: h_2h_3^{-2}h_2^{-1}h_3^2=1$:\\
By interchanging $h_2$ and $h_3$ in (74) and with the same discussion, there is a contradiction.


\item[(222)]$ R_1: h_2h_3^2h_2^{-1}h_3^{-1}h_2=1$, $R_2: h_2h_3^{-2}h_2^{-1}h_3^2=1$:\\
By interchanging $h_2$ and $h_3$ in (70) and with the same discussion, there is a contradiction.

\item[(223)]$ R_1: h_2h_3^2h_2^{-1}h_3^{-1}h_2=1$, $R_2: h_2h_3^{-1}h_2h_3^3=1$:\\
By interchanging $h_2$ and $h_3$ in (7) and with the same discussion, there is a contradiction.



\item[(226)]$ R_1: h_2h_3^2h_2^{-1}h_3^2=1$, $R_2: h_2(h_3h_2^{-1})^2h_3^{-1}h_2=1$:\\
By interchanging $h_2$ and $h_3$ in (180) and with the same discussion, there is a contradiction.

\item[(227)]$ R_1: h_2h_3^2h_2^{-1}h_3^2=1$, $R_2: h_2(h_3h_2^{-1})^3h_3=1$:\\
By interchanging $h_2$ and $h_3$ in (185) and with the same discussion, there is a contradiction.


\item[(229)]$ R_1: h_2h_3^2h_2^{-1}h_3^2=1$, $R_2: h_2h_3^{-2}h_2^{-1}h_3^2=1$:\\
By interchanging $h_2$ and $h_3$ in (76) and with the same discussion, there is a contradiction.

\item[(230)]$ R_1: h_2h_3^2h_2^{-1}h_3^2=1$, $R_2: h_2h_3^{-4}h_2=1$:\\
By interchanging $h_2$ and $h_3$ in (18) and with the same discussion, there is a contradiction.

\item[(231)]$ R_1: h_2h_3^2h_2^{-1}h_3^2=1$, $R_2: h_2h_3^{-3}h_2h_3=1$:\\
By interchanging $h_2$ and $h_3$ in (132) and with the same discussion, there is a contradiction.

\item[(232)]$ R_1: h_2h_3^2h_2^{-1}h_3^2=1$, $R_2: h_2h_3^{-2}h_2h_3^2=1$:\\
By interchanging $h_2$ and $h_3$ in (43) and with the same discussion, there is a contradiction.

\item[(233)]$ R_1: h_2h_3^2h_2^{-1}h_3^2=1$, $R_2: h_2h_3^{-1}h_2h_3^3=1$:\\
By interchanging $h_2$ and $h_3$ in (8) and with the same discussion, there is a contradiction.

\item[(234)]$ R_1: h_2h_3^2h_2^{-1}h_3^2=1$, $R_2: h_2h_3^{-1}h_2(h_3h_2^{-1})^2h_3=1$:\\
By interchanging $h_2$ and $h_3$ in (184) and with the same discussion, there is a contradiction.

\item[(235)]$ R_1: h_2h_3^2h_2^{-1}h_3^2=1$, $R_2: (h_2h_3^{-1})^2(h_2^{-1}h_3)^2=1$:\\
By interchanging $h_2$ and $h_3$ in (181) and with the same discussion, there is a contradiction.



\item[(238)]$ R_1: h_2h_3^2h_2^{-1}h_3^2=1$, $R_2: (h_2h_3^{-1})^3h_2h_3=1$:\\
By interchanging $h_2$ and $h_3$ in (171) and with the same discussion, there is a contradiction.

\item[(239)]$ R_1: h_2h_3^2h_2^{-1}h_3^2=1$, $R_2: h_3(h_3h_2^{-1})^3h_3=1$:\\
By interchanging $h_2$ and $h_3$ in (168) and with the same discussion, there is a contradiction.


\item[(241)]$ R_1: h_2h_3(h_3h_2^{-1})^2h_3=1$, $R_2: h_2(h_3h_2^{-1}h_3)^2=1$:\\
By interchanging $h_2$ and $h_3$ in (199) and with the same discussion, there is a contradiction.

\item[(242)]$ R_1: h_2h_3(h_3h_2^{-1})^2h_3=1$, $R_2: h_2h_3^{-2}h_2^{-1}h_3^2=1$:\\
By interchanging $h_2$ and $h_3$ in (81) and with the same discussion, there is a contradiction.


\item[(244)]$ R_1: h_2h_3h_2^{-2}h_3^2=1$, $R_2: h_2h_3^{-2}h_2^{-1}h_3^2=1$:\\
By interchanging $h_2$ and $h_3$ in (84) and with the same discussion, there is a contradiction.



\item[(247)]$ R_1: h_2h_3h_2^{-1}(h_2^{-1}h_3)^2=1$, $R_2: h_2(h_3h_2^{-1})^2h_2^{-1}h_3=1$:\\
$\Rightarrow \langle h_2,h_3\rangle=\langle h_3^{-1} h_2\rangle$ is abelian, a contradiction.


\item[(249)]$ R_1: h_2h_3h_2^{-1}h_3^{-1}h_2^{-1}h_3=1$, $R_2: h_2h_3h_2^{-1}h_3^{-2}h_2=1$:\\
$\Rightarrow \langle h_2,h_3\rangle=\langle h_2\rangle$ is abelian, a contradiction.

\item[(250)]$ R_1: h_2h_3h_2^{-1}h_3^{-1}h_2^{-1}h_3=1$, $R_2: h_2h_3^{-2}h_2^{-1}h_3^2=1$:\\
By interchanging $h_2$ and $h_3$ in (94) and with the same discussion, there is a contradiction.

\item[(251)]$ R_1: h_2h_3h_2^{-1}h_3^{-1}h_2^{-1}h_3=1$, $R_2: h_2h_3^{-3}h_2h_3=1$:\\
By interchanging $h_2$ and $h_3$ in (134) and with the same discussion, there is a contradiction.



\item[(254)]$ R_1: h_2h_3h_2^{-1}h_3^{-2}h_2=1$, $R_2: h_2h_3^{-2}h_2^{-1}h_3^2=1$:\\
By interchanging $h_2$ and $h_3$ in (95) and with the same discussion, there is a contradiction.





\item[(259)]$ R_1: h_2h_3h_2^{-1}(h_3^{-1}h_2)^2=1$, $R_2: (h_2h_3^{-1})^2(h_2^{-1}h_3)^2=1$:\\
$\Rightarrow \langle h_2,h_3\rangle=\langle h_2\rangle$ is abelian, a contradiction.





\item[(264)]$ R_1: h_2(h_3h_2^{-1}h_3)^2=1$, $R_2: h_2(h_3h_2^{-1})^2h_3^2=1$:\\
By interchanging $h_2$ and $h_3$ in (213) and with the same discussion, there is a contradiction.

\item[(265)]$ R_1: h_2(h_3h_2^{-1}h_3)^2=1$, $R_2: h_2h_3^{-2}h_2^{-1}h_3^2=1$:\\
By interchanging $h_2$ and $h_3$ in (113) and with the same discussion, there is a contradiction.

\item[(266)]$ R_1: h_2(h_3h_2^{-1}h_3)^2=1$, $R_2: h_2(h_3^{-1}h_2h_3^{-1})^2h_2=1$:\\
$\Rightarrow h_2=1$ that is a contradiction.




\item[(270)]$ R_1: h_2(h_3h_2^{-1})^2h_3^{-1}h_2=1$, $R_2: h_2h_3^{-1}h_2^{-1}h_3h_2^{-1}h_3^{-1}h_2=1$:\\
$\Rightarrow h_3=1$ that is a contradiction.


\item[(272)]$ R_1: h_2(h_3h_2^{-1})^2h_3^{-1}h_2=1$, $R_2: h_2h_3^{-1}h_2(h_3h_2^{-1})^2h_3=1$:\\
$\Rightarrow\langle h_2,h_3\rangle\cong BS(2,1)$ is solvable, a contradiction.

\item[(273)]$ R_1: h_2(h_3h_2^{-1})^2h_3^{-1}h_2=1$, $R_2: (h_2h_3^{-1})^2(h_2^{-1}h_3)^2=1$:\\
$\Rightarrow h_3=1$ that is a contradiction.


\item[(275)]$ R_1: h_2(h_3h_2^{-1})^2h_3^2=1$, $R_2: h_2h_3^{-2}h_2^{-1}h_3^2=1$:\\
By interchanging $h_2$ and $h_3$ in (86) and with the same discussion, there is a contradiction.

\item[(276)]$ R_1: h_2(h_3h_2^{-1})^2h_3^2=1$, $R_2: (h_2h_3^{-1})^2(h_2^{-1}h_3)^2=1$:\\
By interchanging $h_2$ and $h_3$ in (214) and with the same discussion, there is a contradiction.


\item[(278)]$ R_1: h_2(h_3h_2^{-1})^3h_3=1$, $R_2: h_2h_3^{-2}h_2^{-1}h_3^2=1$:\\
By interchanging $h_2$ and $h_3$ in (126) and with the same discussion, there is a contradiction.



\item[(281)]$ R_1: h_2h_3^{-1}h_2^{-1}h_3h_2^{-1}h_3^{-1}h_2=1$, $R_2: h_2h_3^{-2}h_2^{-1}h_3^2=1$:\\
By interchanging $h_2$ and $h_3$ in (111) and with the same discussion, there is a contradiction.


\item[(283)]$ R_1: h_2h_3^{-1}(h_2^{-1}h_3)^2h_3=1$, $R_2: h_2h_3^{-2}h_2^{-1}h_3^2=1$:\\
By interchanging $h_2$ and $h_3$ in (96) and with the same discussion, there is a contradiction.

\item[(284)]$ R_1: h_2h_3^{-1}(h_2^{-1}h_3)^2h_3=1$, $R_2: h_2h_3^{-4}h_2=1$:\\
By interchanging $h_2$ and $h_3$ in (21) and with the same discussion, there is a contradiction.

\item[(285)]$ R_1: h_2h_3^{-1}(h_2^{-1}h_3)^2h_3=1$, $R_2: (h_2h_3^{-1})^2(h_2^{-1}h_3)^2=1$:\\
By interchanging $h_2$ and $h_3$ in (259) and with the same discussion, there is a contradiction.


\item[(287)]$ R_1: h_2h_3^{-2}h_2^{-1}h_3^2=1$, $R_2: h_2h_3^{-2}(h_2^{-1}h_3)^2=1$:\\
By interchanging $h_2$ and $h_3$ in (71) and with the same discussion, there is a contradiction.

\item[(288)]$ R_1: h_2h_3^{-2}h_2^{-1}h_3^2=1$, $R_2: h_2h_3^{-4}h_2=1$:\\
By interchanging $h_2$ and $h_3$ in (17) and with the same discussion, there is a contradiction.

\item[(289)]$ R_1: h_2h_3^{-2}h_2^{-1}h_3^2=1$, $R_2: h_2h_3^{-3}h_2h_3=1$:\\
By interchanging $h_2$ and $h_3$ in (69) and with the same discussion, there is a contradiction.

\item[(290)]$ R_1: h_2h_3^{-2}h_2^{-1}h_3^2=1$, $R_2: h_2h_3^{-2}(h_3^{-1}h_2)^2=1$:\\
By interchanging $h_2$ and $h_3$ in (75) and with the same discussion, there is a contradiction.

\item[(291)]$ R_1: h_2h_3^{-2}h_2^{-1}h_3^2=1$, $R_2: h_2h_3^{-2}h_2h_3^2=1$:\\
By interchanging $h_2$ and $h_3$ in (41) and with the same discussion, there is a contradiction.

\item[(292)]$ R_1: h_2h_3^{-2}h_2^{-1}h_3^2=1$, $R_2: h_2h_3^{-2}h_2h_3h_2^{-1}h_3=1$:\\
By interchanging $h_2$ and $h_3$ in (103) and with the same discussion, there is a contradiction.

\item[(293)]$ R_1: h_2h_3^{-2}h_2^{-1}h_3^2=1$, $R_2: h_2h_3^{-1}(h_3^{-1}h_2)^2h_3=1$:\\
By interchanging $h_2$ and $h_3$ in (93) and with the same discussion, there is a contradiction.


\item[(295)]$ R_1: h_2h_3^{-2}h_2^{-1}h_3^2=1$, $R_2: h_2h_3^{-1}h_2h_3^3=1$:\\
By interchanging $h_2$ and $h_3$ in (5) and with the same discussion, there is a contradiction.

\item[(296)]$ R_1: h_2h_3^{-2}h_2^{-1}h_3^2=1$, $R_2: h_2h_3^{-1}h_2h_3^2h_2^{-1}h_3=1$:\\
By interchanging $h_2$ and $h_3$ in (78) and with the same discussion, there is a contradiction.

\item[(297)]$ R_1: h_2h_3^{-2}h_2^{-1}h_3^2=1$, $R_2: h_2h_3^{-1}h_2h_3h_2^{-1}h_3^2=1$:\\
By interchanging $h_2$ and $h_3$ in (97) and with the same discussion, there is a contradiction.


\item[(299)]$ R_1: h_2h_3^{-2}h_2^{-1}h_3^2=1$, $R_2: (h_2h_3^{-1})^2h_2^{-1}h_3^2=1$:\\
By interchanging $h_2$ and $h_3$ in (100) and with the same discussion, there is a contradiction.

\item[(300)]$ R_1: h_2h_3^{-2}h_2^{-1}h_3^2=1$, $R_2: (h_2h_3^{-1})^2(h_2^{-1}h_3)^2=1$:\\
By interchanging $h_2$ and $h_3$ in (119) and with the same discussion, there is a contradiction.

\item[(301)]$ R_1: h_2h_3^{-2}h_2^{-1}h_3^2=1$, $R_2: (h_2h_3^{-1})^2h_3^{-1}h_2^{-1}h_3=1$:\\
By interchanging $h_2$ and $h_3$ in (79) and with the same discussion, there is a contradiction.

\item[(302)]$ R_1: h_2h_3^{-2}h_2^{-1}h_3^2=1$, $R_2: (h_2h_3^{-1})^2h_3^{-2}h_2=1$:\\
By interchanging $h_2$ and $h_3$ in (80) and with the same discussion, there is a contradiction.

\item[(303)]$ R_1: h_2h_3^{-2}h_2^{-1}h_3^2=1$, $R_2: (h_2h_3^{-1})^2h_3^{-1}h_2h_3=1$:\\
By interchanging $h_2$ and $h_3$ in (99) and with the same discussion, there is a contradiction.


\item[(305)]$ R_1: h_2h_3^{-2}h_2^{-1}h_3^2=1$, $R_2: (h_2h_3^{-1})^2h_2h_3^2=1$:\\
By interchanging $h_2$ and $h_3$ in (58) and with the same discussion, there is a contradiction.



\item[(308)]$ R_1: h_2h_3^{-2}h_2^{-1}h_3^2=1$, $R_2: h_3^2(h_3h_2^{-1})^2h_3=1$:\\
By interchanging $h_2$ and $h_3$ in (26) and with the same discussion, there is a contradiction.



\item[(311)]$ R_1: h_2h_3^{-2}(h_2^{-1}h_3)^2=1$, $R_2: h_2h_3^{-2}h_2h_3^2=1$:\\
By interchanging $h_2$ and $h_3$ in (42) and with the same discussion, there is a contradiction.

\item[(312)]$ R_1: h_2h_3^{-2}(h_2^{-1}h_3)^2=1$, $R_2: h_2h_3^{-2}h_2h_3h_2^{-1}h_3=1$:\\
By interchanging $h_2$ and $h_3$ in (148) and with the same discussion, there is a contradiction.

\item[(313)]$ R_1: h_2h_3^{-2}(h_2^{-1}h_3)^2=1$, $R_2: (h_2h_3^{-1})^2h_2h_3^2=1$:\\
By interchanging $h_2$ and $h_3$ in (59) and with the same discussion, there is a contradiction.

\item[(314)]$ R_1: h_2h_3^{-4}h_2=1$, $R_2: h_2h_3^{-2}h_2h_3^2=1$:\\
By interchanging $h_2$ and $h_3$ in (16) and with the same discussion, there is a contradiction.

\item[(315)]$ R_1: h_2h_3^{-4}h_2=1$, $R_2: (h_2h_3^{-1})^2(h_2^{-1}h_3)^2=1$:\\
By interchanging $h_2$ and $h_3$ in (24) and with the same discussion, there is a contradiction.

\item[(316)]$ R_1: h_2h_3^{-4}h_2=1$, $R_2: (h_2h_3^{-1})^2h_3^{-1}h_2^{-1}h_3=1$:\\
By interchanging $h_2$ and $h_3$ in (19) and with the same discussion, there is a contradiction.


\item[(318)]$ R_1: h_2h_3^{-3}h_2h_3=1$, $R_2: h_2h_3^{-1}h_2h_3^3=1$:\\
By interchanging $h_2$ and $h_3$ in (6) and with the same discussion, there is a contradiction.


\item[(320)]$ R_1: h_2h_3^{-3}h_2h_3=1$, $R_2: (h_2h_3^{-1})^2(h_2^{-1}h_3)^2=1$:\\
By interchanging $h_2$ and $h_3$ in (139) and with the same discussion, there is a contradiction.


\item[(322)]$ R_1: h_2h_3^{-3}h_2h_3=1$, $R_2: (h_2h_3^{-1})^3h_2h_3=1$:\\
By interchanging $h_2$ and $h_3$ in (136) and with the same discussion, there is a contradiction.

\item[(323)]$ R_1: h_2h_3^{-2}(h_3^{-1}h_2)^2=1$, $R_2: (h_2h_3^{-1})^2(h_2^{-1}h_3)^2=1$:\\
By interchanging $h_2$ and $h_3$ in (166) and with the same discussion, there is a contradiction.


\item[(325)]$ R_1: h_2h_3^{-2}h_2h_3^2=1$, $R_2: h_2h_3^{-2}h_2h_3h_2^{-1}h_3=1$:\\
By interchanging $h_2$ and $h_3$ in (46) and with the same discussion, there is a contradiction.

\item[(326)]$ R_1: h_2h_3^{-2}h_2h_3^2=1$, $R_2: (h_2h_3^{-1})^2h_2^{-1}h_3^2=1$:\\
By interchanging $h_2$ and $h_3$ in (45) and with the same discussion, there is a contradiction.


\item[(328)]$ R_1: h_2h_3^{-2}h_2h_3^2=1$, $R_2: h_2(h_3^{-1}h_2h_3^{-1})^2h_2=1$:\\
By interchanging $h_2$ and $h_3$ in (52) and with the same discussion, there is a contradiction.

\item[(329)]$ R_1: h_2h_3^{-2}h_2h_3^2=1$, $R_2: (h_2h_3^{-1})^2h_2h_3^2=1$:\\
By interchanging $h_2$ and $h_3$ in (40) and with the same discussion, there is a contradiction.


\item[(331)]$ R_1: h_2h_3^{-2}h_2h_3h_2^{-1}h_3=1$, $R_2: (h_2h_3^{-1})^2h_2^{-1}h_3^2=1$:\\
By interchanging $h_2$ and $h_3$ in (270) and with the same discussion, there is a contradiction.


\item[(333)]$ R_1: h_2h_3^{-2}h_2h_3^{-1}h_2^{-1}h_3=1$, $R_2: (h_2h_3^{-1})^2h_3^{-1}h_2^{-1}h_3=1$:\\
By interchanging $h_2$ and $h_3$ in (193) and with the same discussion, there is a contradiction.


\item[(335)]$ R_1: h_2h_3^{-1}(h_3^{-1}h_2)^2h_3=1$, $R_2: (h_2h_3^{-1})^2h_3^{-1}h_2h_3=1$:\\
By interchanging $h_2$ and $h_3$ in (247) and with the same discussion, there is a contradiction.

\item[(336)]$ R_1: (h_2h_3^{-1}h_2)^2h_3=1$, $R_2: h_2(h_3^{-1}h_2h_3^{-1})^2h_2=1$:\\
By interchanging $h_2$ and $h_3$ in (266) and with the same discussion, there is a contradiction.

\item[(337)]$ R_1: (h_2h_3^{-1}h_2)^2h_3=1$, $R_2: h_3^2(h_3h_2^{-1})^2h_3=1$:\\
By interchanging $h_2$ and $h_3$ in (27) and with the same discussion, there is a contradiction.

\item[(338)]$ R_1: h_2h_3^{-1}h_2h_3^3=1$, $R_2: h_2h_3^{-1}h_2(h_3h_2^{-1})^2h_3=1$:\\
By interchanging $h_2$ and $h_3$ in (15) and with the same discussion, there is a contradiction.


\item[(340)]$ R_1: h_2h_3^{-1}h_2h_3^3=1$, $R_2: h_2(h_3^{-1}h_2h_3^{-1})^2h_2=1$:\\
By interchanging $h_2$ and $h_3$ in (14) and with the same discussion, there is a contradiction.





\item[(345)]$ R_1: (h_2h_3^{-1})^2h_2^{-1}h_3^2=1$, $R_2: (h_2h_3^{-1})^2(h_2^{-1}h_3)^2=1$:\\
By interchanging $h_2$ and $h_3$ in (273) and with the same discussion, there is a contradiction.

\item[(346)]$ R_1: (h_2h_3^{-1})^2h_2^{-1}h_3^2=1$, $R_2: (h_2h_3^{-1})^2h_2h_3^2=1$:\\
By interchanging $h_2$ and $h_3$ in (62) and with the same discussion, there is a contradiction.

\item[(347)]$ R_1: (h_2h_3^{-1})^2h_2^{-1}h_3^2=1$, $R_2: (h_2h_3^{-1})^2h_2h_3h_2^{-1}h_3=1$:\\
By interchanging $h_2$ and $h_3$ in (272) and with the same discussion, there is a contradiction.

\item[(348)]$ R_1: (h_2h_3^{-1})^2(h_2^{-1}h_3)^2=1$, $R_2: (h_2h_3^{-1})^2h_3^{-2}h_2=1$:\\
By interchanging $h_2$ and $h_3$ in (196) and with the same discussion, there is a contradiction.

\item[(349)]$ R_1: (h_2h_3^{-1})^2(h_2^{-1}h_3)^2=1$, $R_2: (h_2h_3^{-1})^2h_2h_3^2=1$:\\
By interchanging $h_2$ and $h_3$ in (66) and with the same discussion, there is a contradiction.


\item[(351)]$ R_1: (h_2h_3^{-1})^2(h_2^{-1}h_3)^2=1$, $R_2: h_3(h_3h_2^{-1})^3h_3=1$:\\
By interchanging $h_2$ and $h_3$ in (202) and with the same discussion, there is a contradiction.




\item[(355)]$ R_1: (h_2h_3^{-1})^2h_2h_3^2=1$, $R_2: (h_2h_3^{-1})^2h_2h_3h_2^{-1}h_3=1$:\\
By interchanging $h_2$ and $h_3$ in (65) and with the same discussion, there is a contradiction.

\end{enumerate}

\subsection{$\mathbf{C_5-C_5(--C_6--)}$}
By considering the relations from Tables \ref{tab-C6} and \ref{tab-C5-C5} which are not disproved, it can be seen that there are $440$ cases for the relations of two cycles $C_5$ and a cycle $C_6$ in the graph $C_5-C_5(--C_6--)$. Using GAP \cite{gap}, we see that all groups with two generators $h_2$ and $h_3$ and three relations which are between $404$ cases of these $440$ cases are finite and solvable, that is a contradiction with the assumptions. So, there are just $36$ cases for the relations of these cycles which may lead to the existence of a subgraph isomorphic to the graph $C_5-C_5(--C_6--)$ in $K(\alpha,\beta)$. These cases are listed in table \ref{tab-C5-C5(--C6--)}.  In the following, we show that all of these $36$ cases lead to contradictions and so, the graph $K(\alpha,\beta)$ contains no subgraph isomorphic to the graph $C_5-C_5(--C_6--)$.

\begin{table}[h]
\small
\centering
\caption{The relations of a $C_5-C_5(--C_6--)$ in $K(\alpha,\beta)$}\label{tab-C5-C5(--C6--)}
\begin{tabular}{|c|l|l|l|}\hline
$n$&$R_1$&$R_2$&$R_3$\\\hline
$1$&$h_2^3h_3^{-2}h_2=1$&$h_2h_3^{-2}h_2^{-1}h_3^2=1$&$h_2^2h_3^4=1$\\
$2$&$h_2^2h_3^{-1}h_2^{-2}h_3=1$&$h_2^2h_3^{-1}h_2^{-2}h_3=1$&$h_2h_3h_2h_3^{-2}h_2^{-1}h_3=1$\\
$3$&$h_2^2h_3^{-1}h_2^{-2}h_3=1$&$h_2^2h_3^{-1}h_2^{-2}h_3=1$&$h_2h_3h_2h_3^{-2}h_2h_3=1$\\
$4$&$h_2^2h_3^{-1}h_2^{-2}h_3=1$&$h_2h_3^{-4}h_2=1$&$h_2^4h_3^2=1$\\
$5$&$h_2^2h_3^{-1}h_2^{-2}h_3=1$&$h_2h_3^{-1}(h_3^{-1}h_2)^2h_3=1$&$h_2h_3h_2h_3^{-3}h_2=1$\\
$6$&$h_2^2h_3^{-1}h_2^{-2}h_3=1$&$h_2h_3^{-1}h_2h_3^2h_2^{-1}h_3=1$&$h_2h_3^{-1}h_2^{-1}h_3^2h_2^{-1}h_3^{-1}h_2=1$\\
$7$&$h_2^2h_3^{-1}h_2^{-2}h_3=1$&$h_2(h_3^{-1}h_2h_3^{-1})^2h_2=1$&$h_2h_3^{-3}h_2h_3h_2^{-1}h_3=1$\\
$8$&$h_2^2h_3^{-1}h_2^{-1}h_3^{-1}h_2=1$&$h_2(h_3^{-1}h_2h_3^{-1})^2h_2=1$&$h_2h_3^{-2}h_2h_3^{-1}h_2^{-1}h_3^2=1$\\
$9$&$h_2^2h_3^{-2}h_2^{-1}h_3=1$&$h_2(h_3^{-1}h_2h_3^{-1})^2h_2=1$&$h_2h_3h_2^{-1}h_3^{-1}h_2h_3^{-1}h_2^{-1}h_3=1$\\
$10$&$h_2^2h_3^{-2}h_2^{-1}h_3=1$&$h_2(h_3^{-1}h_2h_3^{-1})^2h_2=1$&$h_2h_3^{-1}h_2^{-1}h_3h_2^{-1}h_3^{-1}h_2h_3=1$\\
$11$&$h_2^2h_3^{-3}h_2=1$&$(h_2h_3^{-1})^2(h_2^{-1}h_3)^2=1$&$h_2h_3h_2^{-1}h_3^{-1}h_2h_3^{-2}h_2=1$\\
$12$&$h_2^2h_3^{-3}h_2=1$&$(h_2h_3^{-1})^2(h_2^{-1}h_3)^2=1$&$h_2h_3^{-1}h_2^{-1}h_3h_2^{-2}h_3^2=1$\\
$13$&$h_2^2h_3^{-2}h_2h_3=1$&$h_2^2h_3^{-2}h_2h_3=1$&$h_2h_3^{-1}(h_2^{-1}h_3)^3h_3=1$\\
$14$&$h_2^2h_3^{-1}h_2^2h_3=1$&$h_2^2h_3^{-1}h_2^2h_3=1$&$h_2h_3^{-1}(h_2^{-1}h_3)^2h_2^{-1}h_3^{-1}h_2=1$\\
$15$&$h_2^2h_3^{-1}h_2^2h_3=1$&$h_2^2h_3^{-1}h_2^2h_3=1$&$h_2h_3^{-1}h_2(h_3h_2^{-1})^2h_3^{-1}h_2=1$\\
$16$&$h_2^2h_3^{-1}h_2h_3h_2^{-1}h_3=1$&$h_2^2h_3^{-1}h_2h_3h_2^{-1}h_3=1$&$(h_2h_3^{-1})^2h_3^{-2}h_2^{-1}h_3=1$\\
$17$&$h_2^2h_3^{-1}h_2h_3h_2^{-1}h_3=1$&$h_2h_3^{-2}h_2^{-1}h_3^2=1$&$h_2h_3^{-1}h_2^{-1}h_3^2h_2^{-1}h_3^{-1}h_2=1$\\
$18$&$h_2h_3^2h_2^{-2}h_3=1$&$h_2h_3^2h_2^{-2}h_3=1$&$h_2h_3h_2^{-1}(h_3^{-1}h_2)^3=1$\\
$19$&$h_2h_3^2h_2^{-1}h_3^2=1$&$h_2h_3^2h_2^{-1}h_3^2=1$&$h_2h_3^{-2}h_2(h_3h_2^{-1})^2h_3=1$\\
$20$&$h_2h_3^2h_2^{-1}h_3^2=1$&$h_2h_3^2h_2^{-1}h_3^2=1$&$(h_2h_3^{-1})^2h_2^{-1}h_3^2h_2^{-1}h_3=1$\\
$21$&$h_2h_3h_2^{-1}(h_2^{-1}h_3)^2=1$&$h_2h_3^{-2}h_2^{-1}h_3^2=1$&$h_2^2h_3^{-1}h_2^{-1}h_3^{-2}h_2=1$\\
$22$&$h_2h_3h_2^{-1}h_3^{-1}h_2^{-1}h_3=1$&$h_2(h_3^{-1}h_2h_3^{-1})^2h_2=1$&$h_2^2h_3^{-2}h_2h_3^{-1}h_2^{-1}h_3=1$\\
$23$&$h_2h_3h_2^{-1}h_3^{-1}h_2^{-1}h_3=1$&$h_2(h_3^{-1}h_2h_3^{-1})^2h_2=1$&$h_2h_3^{-2}h_2^{-1}h_3h_2h_3^{-1}h_2=1$\\
$24$&$h_2h_3h_2^{-1}h_3^{-2}h_2=1$&$(h_2h_3^{-1})^2(h_2^{-1}h_3)^2=1$&$h_2(h_2h_3^{-2})^2h_2=1$\\
$25$&$h_2h_3h_2^{-1}h_3^{-2}h_2=1$&$(h_2h_3^{-1})^2(h_2^{-1}h_3)^2=1$&$h_2h_3^{-3}h_2^2h_3^{-1}h_2=1$\\
$26$&$h_2h_3h_2^{-1}(h_3^{-1}h_2)^2=1$&$h_2h_3h_2^{-1}(h_3^{-1}h_2)^2=1$&$h_2h_3^2h_2^{-1}(h_2^{-1}h_3)^2=1$\\
$27$&$h_2h_3h_2^{-1}(h_3^{-1}h_2)^2=1$&$(h_2h_3^{-1})^3h_2h_3=1$&$h_2^2h_3^{-1}(h_3^{-1}h_2)^2h_3=1$\\
$28$&$h_2(h_3h_2^{-1})^3h_3=1$&$h_2h_3^{-1}(h_2^{-1}h_3)^2h_3=1$&$h_2h_3^2h_2^{-1}(h_2^{-1}h_3)^2=1$\\
$29$&$h_2(h_3h_2^{-1})^3h_3=1$&$h_2(h_3^{-1}h_2h_3^{-1})^2h_2=1$&$h_2^2(h_3h_2^{-1}h_3)^2=1$\\
$30$&$h_2h_3^{-1}(h_2^{-1}h_3)^2h_3=1$&$h_2h_3^{-1}(h_2^{-1}h_3)^2h_3=1$&$h_2^2h_3^{-1}(h_3^{-1}h_2)^2h_3=1$\\
$31$&$h_2h_3^{-2}h_2^{-1}h_3^2=1$&$h_2h_3^{-2}h_2^{-1}h_3^2=1$&$h_2^2(h_3^{-1}h_2^{-1})^2h_3=1$\\
$32$&$h_2h_3^{-2}h_2^{-1}h_3^2=1$&$h_2h_3^{-2}h_2^{-1}h_3^2=1$&$(h_2h_3)^2h_2^{-2}h_3=1$\\
$33$&$h_2h_3^{-2}h_2^{-1}h_3^2=1$&$h_2(h_3^{-1}h_2h_3^{-1})^2h_2=1$&$h_2^2h_3^{-1}h_2^{-1}h_3h_2^{-1}h_3^{-1}h_2=1$\\
$34$&$h_2h_3^{-3}h_2h_3=1$&$h_2(h_3^{-1}h_2h_3^{-1})^2h_2=1$&$h_2h_3h_2^{-2}h_3h_2^{-1}h_3^{-1}h_2=1$\\
$35$&$h_2h_3^{-1}h_2h_3^2h_2^{-1}h_3=1$&$h_2h_3^{-1}h_2h_3^2h_2^{-1}h_3=1$&$h_2^2(h_2h_3^{-1})^2h_2^{-1}h_3=1$\\
$36$&$h_2(h_3^{-1}h_2h_3^{-1})^2h_2=1$&$(h_2h_3^{-1})^3h_2h_3=1$&$(h_2h_3^{-1}h_2)^2h_3^2=1$\\
\hline
\end{tabular}
\end{table}

\begin{enumerate}
\item[(1)]$ R_1:h_2^3h_3^{-2}h_2=1$, $R_2:h_2h_3^{-2}h_2^{-1}h_3^2=1$, $R_3:h_2^2h_3^4=1$:\\
$\Rightarrow h_2=1$ that is a contradiction.

\item[(2)]$ R_1:h_2^2h_3^{-1}h_2^{-2}h_3=1$, $R_2:h_2^2h_3^{-1}h_2^{-2}h_3=1$, $R_3:h_2h_3h_2h_3^{-2}h_2^{-1}h_3=1$:\\
$\Rightarrow h_2=1$ that is a contradiction.

\item[(3)]$ R_1:h_2^2h_3^{-1}h_2^{-2}h_3=1$, $R_2:h_2^2h_3^{-1}h_2^{-2}h_3=1$, $R_3:h_2h_3h_2h_3^{-2}h_2h_3=1$:\\
$\Rightarrow h_2=1$ that is a contradiction.

\item[(4)]$ R_1:h_2^2h_3^{-1}h_2^{-2}h_3=1$, $R_2:h_2h_3^{-4}h_2=1$, $R_3:h_2^4h_3^2=1$:\\
By interchanging $h_2$ and $h_3$ in (1) and with the same discussion, there is a contradiction.

\item[(5)]$ R_1:h_2^2h_3^{-1}h_2^{-2}h_3=1$, $R_2:h_2h_3^{-1}(h_3^{-1}h_2)^2h_3=1$, $R_3:h_2h_3h_2h_3^{-3}h_2=1$:\\
$\Rightarrow\langle h_2,h_3\rangle\cong BS(1,1)$ is solvable, a contradiction.

\item[(6)]$ R_1:h_2^2h_3^{-1}h_2^{-2}h_3=1$, $R_2:h_2h_3^{-1}h_2h_3^2h_2^{-1}h_3=1$, $R_3:h_2h_3^{-1}h_2^{-1}h_3^2h_2^{-1}h_3^{-1}h_2=1$:\\
$\Rightarrow \langle h_2,h_3\rangle=\langle h_3\rangle$ is abelian, a contradiction.

\item[(7)]$ R_1:h_2^2h_3^{-1}h_2^{-2}h_3=1$, $R_2:h_2(h_3^{-1}h_2h_3^{-1})^2h_2=1$, $R_3:h_2h_3^{-3}h_2h_3h_2^{-1}h_3=1$:\\
$\Rightarrow\langle h_2,h_3\rangle\cong BS(1,-1)$ is solvable, a contradiction.

\item[(8)]$ R_1:h_2^2h_3^{-1}h_2^{-1}h_3^{-1}h_2=1$, $R_2:h_2(h_3^{-1}h_2h_3^{-1})^2h_2=1$, $R_3:h_2h_3^{-2}h_2h_3^{-1}h_2^{-1}h_3^2=1$:\\
$\Rightarrow \langle h_2,h_3\rangle=\langle h_2\rangle$ is abelian, a contradiction.

\item[(9)]$ R_1:h_2^2h_3^{-2}h_2^{-1}h_3=1$, $R_2:h_2(h_3^{-1}h_2h_3^{-1})^2h_2=1$, $R_3:h_2h_3h_2^{-1}h_3^{-1}h_2h_3^{-1}h_2^{-1}h_3=1$:\\
$\Rightarrow\langle h_2,h_3\rangle\cong BS(1,1)$ is solvable, a contradiction.

\item[(10)]$ R_1:h_2^2h_3^{-2}h_2^{-1}h_3=1$, $R_2:h_2(h_3^{-1}h_2h_3^{-1})^2h_2=1$, $R_3:h_2h_3^{-1}h_2^{-1}h_3h_2^{-1}h_3^{-1}h_2h_3=1$:\\
By interchanging $h_2$ and $h_3$ in (9) and with the same discussion, there is a contradiction.

\item[(11)]$ R_1:h_2^2h_3^{-3}h_2=1$, $R_2:(h_2h_3^{-1})^2(h_2^{-1}h_3)^2=1$, $R_3:h_2h_3h_2^{-1}h_3^{-1}h_2h_3^{-2}h_2=1$:\\
$\Rightarrow \langle h_2,h_3\rangle=\langle h_2\rangle$ is abelian, a contradiction.

\item[(12)]$ R_1:h_2^2h_3^{-3}h_2=1$, $R_2:(h_2h_3^{-1})^2(h_2^{-1}h_3)^2=1$, $R_3:h_2h_3^{-1}h_2^{-1}h_3h_2^{-2}h_3^2=1$:\\
By interchanging $h_2$ and $h_3$ in (11) and with the same discussion, there is a contradiction.

\item[(13)]$ R_1:h_2^2h_3^{-2}h_2h_3=1$, $R_2:h_2^2h_3^{-2}h_2h_3=1$, $R_3:h_2h_3^{-1}(h_2^{-1}h_3)^3h_3=1$:\\
$\Rightarrow h_2=1$ that is a contradiction.

\item[(14)]$ R_1:h_2^2h_3^{-1}h_2^2h_3=1$, $R_2:h_2^2h_3^{-1}h_2^2h_3=1$, $R_3:h_2h_3^{-1}(h_2^{-1}h_3)^2h_2^{-1}h_3^{-1}h_2=1$:\\
$\Rightarrow h_2=1$ that is a contradiction.

\item[(15)]$ R_1:h_2^2h_3^{-1}h_2^2h_3=1$, $R_2:h_2^2h_3^{-1}h_2^2h_3=1$, $R_3:h_2h_3^{-1}h_2(h_3h_2^{-1})^2h_3^{-1}h_2=1$:\\
$\Rightarrow h_2=1$ that is a contradiction.

\item[(16)]$ R_1:h_2^2h_3^{-1}h_2h_3h_2^{-1}h_3=1$, $R_2:h_2^2h_3^{-1}h_2h_3h_2^{-1}h_3=1$, $R_3:(h_2h_3^{-1})^2h_3^{-2}h_2^{-1}h_3=1$:\\
$\Rightarrow h_2=1$ that is a contradiction.

\item[(17)]$ R_1:h_2^2h_3^{-1}h_2h_3h_2^{-1}h_3=1$, $R_2:h_2h_3^{-2}h_2^{-1}h_3^2=1$, $R_3:h_2h_3^{-1}h_2^{-1}h_3^2h_2^{-1}h_3^{-1}h_2=1$:\\
By interchanging $h_2$ and $h_3$ in (6) and with the same discussion, there is a contradiction.

\item[(18)]$ R_1:h_2h_3^2h_2^{-2}h_3=1$, $R_2:h_2h_3^2h_2^{-2}h_3=1$, $R_3:h_2h_3h_2^{-1}(h_3^{-1}h_2)^3=1$:\\
By interchanging $h_2$ and $h_3$ in (13) and with the same discussion, there is a contradiction.

\item[(19)]$ R_1:h_2h_3^2h_2^{-1}h_3^2=1$, $R_2:h_2h_3^2h_2^{-1}h_3^2=1$, $R_3:h_2h_3^{-2}h_2(h_3h_2^{-1})^2h_3=1$:\\
By interchanging $h_2$ and $h_3$ in (14) and with the same discussion, there is a contradiction.

\item[(20)]$ R_1:h_2h_3^2h_2^{-1}h_3^2=1$, $R_2:h_2h_3^2h_2^{-1}h_3^2=1$, $R_3:(h_2h_3^{-1})^2h_2^{-1}h_3^2h_2^{-1}h_3=1$:\\
By interchanging $h_2$ and $h_3$ in (15) and with the same discussion, there is a contradiction.

\item[(21)]$ R_1:h_2h_3h_2^{-1}(h_2^{-1}h_3)^2=1$, $R_2:h_2h_3^{-2}h_2^{-1}h_3^2=1$, $R_3:h_2^2h_3^{-1}h_2^{-1}h_3^{-2}h_2=1$:\\
By interchanging $h_2$ and $h_3$ in (5) and with the same discussion, there is a contradiction.

\item[(22)]$ R_1:h_2h_3h_2^{-1}h_3^{-1}h_2^{-1}h_3=1$, $R_2:h_2(h_3^{-1}h_2h_3^{-1})^2h_2=1$, $R_3:h_2^2h_3^{-2}h_2h_3^{-1}h_2^{-1}h_3=1$:\\
$\Rightarrow \langle h_2,h_3\rangle=\langle h_2\rangle$ is abelian, a contradiction.

\item[(23)]$ R_1:h_2h_3h_2^{-1}h_3^{-1}h_2^{-1}h_3=1$, $R_2:h_2(h_3^{-1}h_2h_3^{-1})^2h_2=1$, $R_3:h_2h_3^{-2}h_2^{-1}h_3h_2h_3^{-1}h_2=1$:\\
By interchanging $h_2$ and $h_3$ in (22) and with the same discussion, there is a contradiction.

\item[(24)]$ R_1:h_2h_3h_2^{-1}h_3^{-2}h_2=1$, $R_2:(h_2h_3^{-1})^2(h_2^{-1}h_3)^2=1$, $R_3:h_2(h_2h_3^{-2})^2h_2=1$:\\
$\Rightarrow\langle h_2,h_3\rangle\cong BS(1,-1)$ is solvable, a contradiction.

\item[(25)]$ R_1:h_2h_3h_2^{-1}h_3^{-2}h_2=1$, $R_2:(h_2h_3^{-1})^2(h_2^{-1}h_3)^2=1$, $R_3:h_2h_3^{-3}h_2^2h_3^{-1}h_2=1$:\\
By interchanging $h_2$ and $h_3$ in (24) and with the same discussion, there is a contradiction.

\item[(26)]$ R_1:h_2h_3h_2^{-1}(h_3^{-1}h_2)^2=1$, $R_2:h_2h_3h_2^{-1}(h_3^{-1}h_2)^2=1$, $R_3:h_2h_3^2h_2^{-1}(h_2^{-1}h_3)^2=1$:\\
$\Rightarrow h_3=1$ that is a contradiction.

\item[(27)]$ R_1:h_2h_3h_2^{-1}(h_3^{-1}h_2)^2=1$, $R_2:(h_2h_3^{-1})^3h_2h_3=1$, $R_3:h_2^2h_3^{-1}(h_3^{-1}h_2)^2h_3=1$:\\
$\Rightarrow\langle h_2,h_3\rangle\cong BS(1,-1)$ is solvable, a contradiction.

\item[(28)]$ R_1:h_2(h_3h_2^{-1})^3h_3=1$, $R_2:h_2h_3^{-1}(h_2^{-1}h_3)^2h_3=1$, $R_3:h_2h_3^2h_2^{-1}(h_2^{-1}h_3)^2=1$:\\
By interchanging $h_2$ and $h_3$ in (27) and with the same discussion, there is a contradiction.

\item[(29)]$ R_1:h_2(h_3h_2^{-1})^3h_3=1$, $R_2:h_2(h_3^{-1}h_2h_3^{-1})^2h_2=1$, $R_3:h_2^2(h_3h_2^{-1}h_3)^2=1$:\\
$\Rightarrow h_2=1$ that is a contradiction.

\item[(30)]$ R_1:h_2h_3^{-1}(h_2^{-1}h_3)^2h_3=1$, $R_2:h_2h_3^{-1}(h_2^{-1}h_3)^2h_3=1$, $R_3:h_2^2h_3^{-1}(h_3^{-1}h_2)^2h_3=1$:\\
By interchanging $h_2$ and $h_3$ in (26) and with the same discussion, there is a contradiction.

\item[(31)]$ R_1:h_2h_3^{-2}h_2^{-1}h_3^2=1$, $R_2:h_2h_3^{-2}h_2^{-1}h_3^2=1$, $R_3:h_2^2(h_3^{-1}h_2^{-1})^2h_3=1$:\\
By interchanging $h_2$ and $h_3$ in (2) and with the same discussion, there is a contradiction.

\item[(32)]$ R_1:h_2h_3^{-2}h_2^{-1}h_3^2=1$, $R_2:h_2h_3^{-2}h_2^{-1}h_3^2=1$, $R_3:(h_2h_3)^2h_2^{-2}h_3=1$:\\
By interchanging $h_2$ and $h_3$ in (3) and with the same discussion, there is a contradiction.

\item[(33)]$ R_1:h_2h_3^{-2}h_2^{-1}h_3^2=1$, $R_2:h_2(h_3^{-1}h_2h_3^{-1})^2h_2=1$, $R_3:h_2^2h_3^{-1}h_2^{-1}h_3h_2^{-1}h_3^{-1}h_2=1$:\\
By interchanging $h_2$ and $h_3$ in (7) and with the same discussion, there is a contradiction.

\item[(34)]$ R_1:h_2h_3^{-3}h_2h_3=1$, $R_2:h_2(h_3^{-1}h_2h_3^{-1})^2h_2=1$, $R_3:h_2h_3h_2^{-2}h_3h_2^{-1}h_3^{-1}h_2=1$:\\
By interchanging $h_2$ and $h_3$ in (8) and with the same discussion, there is a contradiction.

\item[(35)]$ R_1:h_2h_3^{-1}h_2h_3^2h_2^{-1}h_3=1$, $R_2:h_2h_3^{-1}h_2h_3^2h_2^{-1}h_3=1$, $R_3:h_2^2(h_2h_3^{-1})^2h_2^{-1}h_3=1$:\\
By interchanging $h_2$ and $h_3$ in (16) and with the same discussion, there is a contradiction.

\item[(36)]$ R_1:h_2(h_3^{-1}h_2h_3^{-1})^2h_2=1$, $R_2:(h_2h_3^{-1})^3h_2h_3=1$, $R_3:(h_2h_3^{-1}h_2)^2h_3^2=1$:\\
By interchanging $h_2$ and $h_3$ in (29) and with the same discussion, there is a contradiction.

\end{enumerate}

$\mathbf{C_5-C_5(-C_6--)}$ \textbf{subgraph:} By considering the relations from Tables \ref{tab-C6} and \ref{tab-C5-C5} which are not disproved, it can be seen that there are $2121$ cases for the relations of two cycles $C_5$ and a cycle $C_6$ in the graph $C_5-C_5(-C_6--)$. Using Gap \cite{gap}, we see that all groups with two generators $h_2$ and $h_3$ and three relations which are between $1907$ cases of these $2121$ cases are finite and solvable, that is a contradiction with the assumptions. So there are $214$ cases for the relations of these cycles which may lead to the existence of a subgraph isomorphic to the graph $C_5-C_5(-C_6--)$ in the graph $K(\alpha,\beta)$. These cases are listed in table \ref{tab-C5-C5(-C6--)}. In the following, we show that $208$ cases of these relations lead to a contradiction and just $6$ cases of them may lead to the existence of a subgraph isomorphic to the graph $C_5-C_5(-C_6--)$ in the graph $K(\alpha,\beta)$. Cases which are not disproved are marked by $*$s in the Table \ref{tab-C5-C5(-C6--)}.

\begin{center}
\begin{figure}[h]
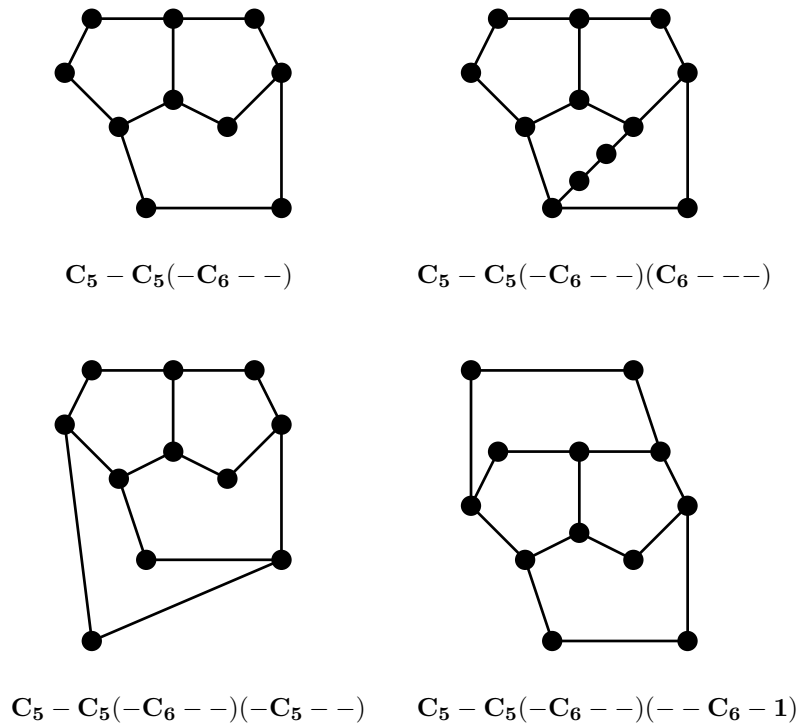

\psscalebox{0.9 0.9} 
{

}
\vspace{1cc}
\caption{The graph $C_5-C_5(-C_6--)$  and some forbidden subgraphs which contain it}\label{f-10}
\end{figure}
\end{center}

\begin{enumerate}
\item[(1)]$ R_1:h_2^3h_3^2=1$, $R_2:h_2h_3^{-2}h_2^{-1}h_3^2=1$, $R_3:h_2^2h_3^{-4}h_2=1$:\\
$\Rightarrow h_3=1$ that is a contradiction.

\begin{center}
\small
\begin{longtable}{|c|l|l|l|}
\caption{{\scriptsize The relations of a $C_5-C_5(-C_6--)$}}\label{tab-C5-C5(-C6--)}\\
\hline \multicolumn{1}{|c|}{$n$} & \multicolumn{1}{l|}{$R_1$} & \multicolumn{1}{l|}{$R_2$} & \multicolumn{1}{l|}{$R_3$} \\\hline
\endfirsthead
\multicolumn{4}{c}%
{{\tablename\ \thetable{} -- continued from previous page}} \\
\hline  \multicolumn{1}{|c|}{$n$} & \multicolumn{1}{l|}{$R_1$} & \multicolumn{1}{l|}{$R_2$} & \multicolumn{1}{l|}{$R_3$}\\\hline
\endhead
\hline \multicolumn{4}{|r|}{{Continued on next page}} \\\hline
\endfoot
\hline
\endlastfoot
$1$&$h_2^3h_3^2=1$&$h_2h_3^{-2}h_2^{-1}h_3^2=1$&$h_2^2h_3^{-4}h_2=1$\\
$2$&$h_2^3h_3^2=1$&$h_2h_3^{-2}h_2^{-1}h_3^2=1$&$h_2h_3^3h_2h_3^{-1}h_2=1$\\
$3$&$h_2^3h_3h_2^{-1}h_3=1$&$h_2h_3^{-1}h_2h_3^3=1$&$h_2h_3h_2^{-1}h_3h_2h_3^{-1}h_2h_3=1$\\
$4$&$h_2^3h_3h_2^{-1}h_3=1$&$h_2h_3^{-1}h_2h_3^3=1$&$(h_2h_3^{-1})^2h_2^{-1}h_3^3=1$\\
$5$&$h_2^3h_3h_2^{-1}h_3=1$&$(h_2h_3^{-1})^2(h_2^{-1}h_3)^2=1$&$h_2h_3^2h_2^{-1}h_3h_2h_3^{-1}h_2=1$\\
$6$&$h_2^3h_3h_2^{-1}h_3=1$&$(h_2h_3^{-1})^2(h_2^{-1}h_3)^2=1$&$h_2h_3h_2^{-1}h_3^{-1}h_2h_3^{-2}h_2=1$\\
$7$&$h_2^3h_3h_2^{-1}h_3=1$&$(h_2h_3^{-1})^2(h_2^{-1}h_3)^2=1$&$h_2h_3^{-2}h_2h_3^2h_2^{-1}h_3=1$\\
$8$&$h_2^3h_3^{-2}h_2=1$&$h_2h_3^{-2}h_2^{-1}h_3^2=1$&$h_2^2h_3^4=1$\\
\end{longtable}
\end{center}

\item[(2)]$ R_1:h_2^3h_3^2=1$, $R_2:h_2h_3^{-2}h_2^{-1}h_3^2=1$, $R_3:h_2h_3^3h_2h_3^{-1}h_2=1$:\\
$\Rightarrow\langle h_2,h_3\rangle\cong BS(1,1)$ is solvable, a contradiction.

\item[(3)]$ R_1:h_2^3h_3h_2^{-1}h_3=1$, $R_2:h_2h_3^{-1}h_2h_3^3=1$, $R_3:h_2h_3h_2^{-1}h_3h_2h_3^{-1}h_2h_3=1$:\\
$\Rightarrow\langle h_2,h_3\rangle\cong BS(1,1)$ is solvable, a contradiction.

\item[(4)]$ R_1:h_2^3h_3h_2^{-1}h_3=1$, $R_2:h_2h_3^{-1}h_2h_3^3=1$, $R_3:(h_2h_3^{-1})^2h_2^{-1}h_3^3=1$:\\
$\Rightarrow \langle h_2,h_3\rangle=\langle h_2\rangle$ is abelian, a contradiction.

\item[(5)]$ R_1:h_2^3h_3h_2^{-1}h_3=1$, $R_2:(h_2h_3^{-1})^2(h_2^{-1}h_3)^2=1$, $R_3:h_2h_3^2h_2^{-1}h_3h_2h_3^{-1}h_2=1$:\\
$\Rightarrow \langle h_2,h_3\rangle=\langle h_2\rangle$ is abelian, a contradiction.

\item[(6)]$ R_1:h_2^3h_3h_2^{-1}h_3=1$, $R_2:(h_2h_3^{-1})^2(h_2^{-1}h_3)^2=1$, $R_3:h_2h_3h_2^{-1}h_3^{-1}h_2h_3^{-2}h_2=1$:\\
$\Rightarrow h_3=1$ that is a contradiction.

\item[(7)]$ R_1:h_2^3h_3h_2^{-1}h_3=1$, $R_2:(h_2h_3^{-1})^2(h_2^{-1}h_3)^2=1$, $R_3:h_2h_3^{-2}h_2h_3^2h_2^{-1}h_3=1$:\\
$\Rightarrow\langle h_2,h_3\rangle\cong BS(-2,1)$ is solvable, a contradiction.

\item[(8)]$ R_1:h_2^3h_3^{-2}h_2=1$, $R_2:h_2h_3^{-2}h_2^{-1}h_3^2=1$, $R_3:h_2^2h_3^4=1$:\\
$\Rightarrow h_3=1$ that is a contradiction.

\item[(9)]$ R_1:h_2^3h_3^{-2}h_2=1$, $R_2:h_2h_3^{-2}h_2^{-1}h_3^2=1$, $R_3:h_2^2h_3^{-1}h_2^{-1}h_3^{-1}h_2h_3=1$:\\
$\Rightarrow\langle h_2,h_3\rangle\cong BS(1,2)$ is solvable, a contradiction.

\item[(10)]$ R_1:h_2^3h_3^{-2}h_2=1$, $R_2:h_2h_3^{-2}h_2^{-1}h_3^2=1$, $R_3:h_2h_3h_2^{-1}(h_2^{-1}h_3)^2h_3=1$:\\
$\Rightarrow h_3=1$ that is a contradiction.

\item[(11)]$ R_1:h_2^3h_3^{-2}h_2=1$, $R_2:(h_2h_3^{-1})^3h_2h_3=1$, $R_3:h_2^2h_3^{-2}(h_2^{-1}h_3)^2=1$:\\
$\Rightarrow\langle h_2,h_3\rangle\cong BS(1,1)$ is solvable, a contradiction.

\item[(12)]$ R_1:h_2^2h_3^3=1$, $R_2:h_2^2h_3^{-1}h_2^{-2}h_3=1$, $R_3:h_2^3h_3^{-3}h_2=1$:\\
By interchanging $h_2$ and $h_3$ in (1) and with the same discussion, there is a contradiction.

\item[(13)]$ R_1:h_2^2h_3^3=1$, $R_2:h_2^2h_3^{-1}h_2^{-2}h_3=1$, $R_3:h_2^3h_3h_2^{-1}h_3^2=1$:\\
By interchanging $h_2$ and $h_3$ in (2) and with the same discussion, there is a contradiction.


\item[(15)]$ R_1:h_2^2h_3h_2^{-2}h_3=1$, $R_2:h_2^2h_3h_2^{-2}h_3=1$, $R_3:h_2h_3h_2h_3^{-1}h_2^{-2}h_3=1$:\\
$\Rightarrow\langle h_2,h_3\rangle\cong BS(1,2)$ is solvable, a contradiction.

\item[(16)]$ R_1:h_2^2h_3h_2^{-2}h_3=1$, $R_2:h_2^2h_3h_2^{-2}h_3=1$, $R_3:h_2h_3h_2h_3^{-1}(h_2^{-1}h_3)^2=1$:\\
$\Rightarrow\langle h_2,h_3\rangle\cong BS(1,1)$ is solvable, a contradiction.

\item[(17)]$ R_1:h_2^2h_3h_2^{-2}h_3=1$, $R_2:h_2^2h_3h_2^{-2}h_3=1$, $R_3:h_2h_3(h_2h_3^{-1})^2h_2h_3=1$:\\
$\Rightarrow h_2=1$ that is a contradiction.

\item[(18)]$ R_1:h_2^2h_3h_2^{-2}h_3=1$, $R_2:h_2^2h_3h_2^{-2}h_3=1$, $R_3:h_2h_3^3h_2^{-1}h_3^2=1$:\\
$\Rightarrow\langle h_2,h_3\rangle\cong BS(-5,1)$ is solvable, a contradiction.

\item[(19)]$ R_1:h_2^2h_3h_2^{-2}h_3=1$, $R_2:h_2^2h_3h_2^{-2}h_3=1$, $R_3:h_2h_3^2h_2^{-1}h_3^3=1$:\\
$\Rightarrow\langle h_2,h_3\rangle\cong BS(1,-5)$ is solvable, a contradiction.

\item[(20)]$ R_1:h_2^2h_3h_2^{-2}h_3=1$, $R_2:h_2h_3^{-1}h_2(h_3h_2^{-1})^2h_3=1$, $R_3:h_2h_3^{-2}h_2^{-1}h_3^3=1$:\\
$\Rightarrow\langle h_2,h_3\rangle\cong BS(1,5)$ is solvable, a contradiction.

\item[(21)]$ R_1:h_2^2h_3h_2^{-2}h_3=1$, $R_2:h_2h_3^{-1}h_2(h_3h_2^{-1})^2h_3=1$, $R_3:h_2h_3^{-3}h_2^{-1}h_3^2=1$:\\
$\Rightarrow\langle h_2,h_3\rangle\cong BS(1,-5)$ is solvable, a contradiction.

\item[(22)]$ R_1:h_2^2h_3h_2^{-2}h_3=1$, $R_2:(h_2h_3^{-1})^2(h_2^{-1}h_3)^2=1$, $R_3:h_2h_3^3h_2^{-1}h_3^2=1$:\\
$\Rightarrow\langle h_2,h_3\rangle\cong BS(-5,1)$ is solvable, a contradiction.

\item[(23)]$ R_1:h_2^2h_3h_2^{-2}h_3=1$, $R_2:(h_2h_3^{-1})^2(h_2^{-1}h_3)^2=1$, $R_3:h_2h_3^2h_2^{-1}h_3^3=1$:\\
$\Rightarrow\langle h_2,h_3\rangle\cong BS(1,-5)$ is solvable, a contradiction.

\item[(24)]$ R_1:h_2^2(h_3h_2^{-1})^2h_3=1$, $R_2:h_2h_3^{-1}h_2^{-1}h_3h_2^{-1}h_3^{-1}h_2=1$, $R_3:h_2^2h_3h_2^{-1}h_3^{-1}h_2^{-1}h_3=1$:\\
$\Rightarrow h_3=1$ that is a contradiction.

\item[(25)]$ R_1:h_2^2(h_3h_2^{-1})^2h_3=1$, $R_2:h_2h_3^{-1}h_2^{-1}h_3h_2^{-1}h_3^{-1}h_2=1$, $R_3:h_2h_3^{-1}h_2^{-1}h_3h_2^{-1}h_3^{-1}h_2h_3=1$:\\
$\Rightarrow h_3=1$ that is a contradiction.

\item[(26)]$ R_1:h_2^2(h_3h_2^{-1})^2h_3=1$, $R_2:h_2h_3^{-1}h_2^{-1}h_3h_2^{-1}h_3^{-1}h_2=1$, $R_3:h_2h_3^{-2}h_2^{-1}h_3^3=1$:\\
$\Rightarrow h_3=1$ that is a contradiction.

\item[(27)]$ R_1:h_2^2(h_3h_2^{-1})^2h_3=1$, $R_2:h_2h_3^{-1}h_2^{-1}h_3h_2^{-1}h_3^{-1}h_2=1$, $R_3:h_2h_3^{-3}h_2^{-1}h_3^2=1$:\\
$\Rightarrow\langle h_2,h_3\rangle\cong BS(1,-5)$ is solvable, a contradiction.

\item[(28)]$ R_1:h_2^2(h_3h_2^{-1})^2h_3=1$, $R_2:h_2h_3^{-2}h_2^{-1}h_3^2=1$, $R_3:h_2^2h_3h_2^{-2}h_3^2=1$:\\
$\Rightarrow h_3=1$ that is a contradiction.

\item[(29)]$ R_1:h_2^2(h_3h_2^{-1})^2h_3=1$, $R_2:h_2h_3^{-2}h_2^{-1}h_3^2=1$, $R_3:(h_2h_3)^2h_2^{-2}h_3=1$:\\
$\Rightarrow h_3=1$ that is a contradiction.

\item[(30)]$ R_1:h_2^2(h_3h_2^{-1})^2h_3=1$, $R_2:h_2h_3^{-2}h_2^{-1}h_3^2=1$, $R_3:h_2h_3^{-1}h_2^{-1}h_3^2h_2^{-1}h_3^{-1}h_2=1$:\\
$R_2\Rightarrow R_3: h_2^2h_3^{-1}h_2^{-2}h_3=1$ and $R_2:h_2h_3^{-2}h_2^{-1}h_3^2=1$. Using Tietze transformation where $h_3\mapsto h_3h_2^{-1}$ and $h_2\mapsto h_2$, we have $R_3: h_2^2h_3^{-1}h_2^{-2}h_3=1$ and $R_2:(h_3h_2^{-1})^2=(h_2^{-1}h_3)^2$. So, $R_3\Rightarrow h_2^2\in Z(G)$ where $G=\langle h_2,h_3\rangle$. Let $x=h_3h_2^{-1}$. So $h_2^{-1}h_3=x^{h_2}$. Also 
$(h_3h_2^{-1})^2=(h_2^{-1}h_3)^2 \text{ so if } H=\langle x,x^{h_2}\rangle\Rightarrow H\cong BS(1,-1)$ is solvable. By Corollary \ref{n-sub2} $H \trianglelefteq G$ since $G=\langle x,h_2\rangle$. Since $\frac{G}{H}=\frac{\langle h_2\rangle H}{H}$ is a cyclic group, it is solvable. So $G$ is solvable, a contradiction.

\item[(31)]$ R_1:h_2^2(h_3h_2^{-1})^2h_3=1$, $R_2:h_2h_3^{-2}h_2^{-1}h_3^2=1$, $R_3:h_2h_3^{-1}h_2^{-1}h_3h_2^{-1}h_3^{-2}h_2=1$:\\
$\Rightarrow\langle h_2,h_3\rangle\cong BS(1,-1)$ is solvable, a contradiction.

\item[(32)]$ R_1:h_2^2(h_3h_2^{-1})^2h_3=1$, $R_2:(h_2h_3^{-1})^2h_2h_3^2=1$, $R_3:h_3^2(h_3h_2^{-1})^3h_3=1$:\\
$\Rightarrow h_3=1$ that is a contradiction.

\item[(33)]$ R_1:h_2^2h_3^{-1}h_2^{-2}h_3=1$, $R_2:h_2^2h_3^{-1}h_2^{-2}h_3=1$, $R_3:h_2^2h_3^2h_2h_3=1$:\\
$\Rightarrow\langle h_2,h_3\rangle\cong BS(1,-2)$ is solvable, a contradiction.

\item[(34)]$ R_1:h_2^2h_3^{-1}h_2^{-2}h_3=1$, $R_2:h_2^2h_3^{-1}h_2^{-2}h_3=1$, $R_3:h_2^2h_3h_2^{-1}h_3h_2h_3=1$:\\
$\Rightarrow \langle h_2,h_3\rangle=\langle h_2^{-1}h_3^{-1}\rangle$ is abelian, a contradiction.

\item[(35)]$ R_1:h_2^2h_3^{-1}h_2^{-2}h_3=1$, $R_2:h_2^2h_3^{-1}h_2^{-2}h_3=1$, $R_3:h_2^2h_3^{-1}h_2^{-1}h_3h_2h_3=1$:\\
$\Rightarrow \langle h_2,h_3\rangle=\langle h_2 \rangle$ is abelian, a contradiction.

\item[(36)]$ R_1:h_2^2h_3^{-1}h_2^{-2}h_3=1$, $R_2:h_2^2h_3^{-1}h_2^{-2}h_3=1$, $R_3:h_2^2h_3^{-1}(h_2h_3)^2=1$:\\
$\Rightarrow \langle h_2,h_3\rangle=\langle h_2 \rangle$ is abelian, a contradiction.

\item[(37)]$ R_1:h_2^2h_3^{-1}h_2^{-2}h_3=1$, $R_2:h_2^2h_3^{-1}h_2^{-2}h_3=1$, $R_3:(h_2h_3)^2h_3h_2^{-1}h_3=1$:\\
$\Rightarrow\langle h_2,h_3\rangle\cong BS(1,1)$ is solvable, a contradiction.

\item[(38)]$ R_1:h_2^2h_3^{-1}h_2^{-2}h_3=1$, $R_2:h_2^2h_3^{-1}h_2^{-2}h_3=1$, $R_3:(h_2h_3)^2h_2^{-2}h_3=1$:\\
$\Rightarrow\langle h_2,h_3\rangle\cong BS(-2,1)$ is solvable, a contradiction.

\item[(39)]$ R_1:h_2^2h_3^{-1}h_2^{-2}h_3=1$, $R_2:h_2^2h_3^{-1}h_2^{-2}h_3=1$, $R_3:(h_2h_3)^2h_2^{-1}h_3^2=1$:\\
$\Rightarrow\langle h_2,h_3\rangle\cong BS(1,1)$ is solvable, a contradiction.

\item[(40)]$ R_1:h_2^2h_3^{-1}h_2^{-2}h_3=1$, $R_2:h_2^2h_3^{-1}h_2^{-2}h_3=1$, $R_3:(h_2h_3)^2(h_2^{-1}h_3)^2=1$:\\
$\Rightarrow\langle h_2,h_3\rangle\cong BS(1,-1)$ is solvable, a contradiction.

\item[(41)]$ R_1:h_2^2h_3^{-1}h_2^{-2}h_3=1$, $R_2:h_2^2h_3^{-1}h_2^{-2}h_3=1$, $R_3:h_2h_3^2h_2h_3^{-2}h_2=1$:\\
$\Rightarrow h_2=1$ that is a contradiction.

\item[(42)]$ R_1:h_2^2h_3^{-1}h_2^{-2}h_3=1$, $R_2:h_2^2h_3^{-1}h_2^{-2}h_3=1$, $R_3:h_2h_3^2h_2^{-1}h_3^{-2}h_2=1$:\\
$\Rightarrow h_2=1$ that is a contradiction.

\item[(43)]$ R_1:h_2^2h_3^{-1}h_2^{-2}h_3=1$, $R_2:h_2^2h_3^{-1}h_2^{-2}h_3=1$, $R_3:h_2h_3^{-2}h_2^{-1}h_3^3=1$:\\
$\Rightarrow\langle h_2,h_3\rangle\cong BS(1,-1)$ is solvable, a contradiction.

\item[(44)]$ R_1:h_2^2h_3^{-1}h_2^{-2}h_3=1$, $R_2:h_2^2h_3^{-1}h_2^{-2}h_3=1$, $R_3:h_2h_3^{-2}(h_2^{-1}h_3)^2h_3=1$:\\
$\Rightarrow\langle h_2,h_3\rangle\cong BS(1,-1)$ is solvable, a contradiction.

\item[(45)]$ R_1:h_2^2h_3^{-1}h_2^{-2}h_3=1$, $R_2:h_2^2h_3^{-1}h_2^{-2}h_3=1$, $R_3:h_2h_3^{-3}h_2^{-1}h_3^2=1$:\\
$\Rightarrow\langle h_2,h_3\rangle\cong BS(1,-1)$ is solvable, a contradiction.

\item[(46)]$ R_1:h_2^2h_3^{-1}h_2^{-2}h_3=1$, $R_2:h_2^2h_3^{-1}h_2^{-2}h_3=1$, $R_3:h_2h_3^{-3}h_2h_3^2=1$:\\
$\Rightarrow\langle h_2,h_3\rangle\cong BS(1,-1)$ is solvable, a contradiction.

\item[(47)]$ R_1:h_2^2h_3^{-1}h_2^{-2}h_3=1$, $R_2:h_2^2h_3^{-1}h_2^{-2}h_3=1$, $R_3:h_2h_3^{-2}h_2h_3^3=1$:\\
$\Rightarrow\langle h_2,h_3\rangle\cong BS(1,-1)$ is solvable, a contradiction.

\item[(48)]$ R_1:h_2^2h_3^{-1}h_2^{-2}h_3=1$, $R_2:h_2^2h_3^{-1}h_2^{-2}h_3=1$, $R_3:h_2h_3^{-2}h_2h_3h_2^{-1}h_3^2=1$:\\
$\Rightarrow\langle h_2,h_3\rangle\cong BS(1,-1)$ is solvable, a contradiction.

\item[(49)]$ R_1:h_2^2h_3^{-1}h_2^{-2}h_3=1$, $R_2:h_2^2h_3^{-1}h_2^{-2}h_3=1$, $R_3:h_2h_3^{-2}h_2h_3^{-1}h_2^{-1}h_3^2=1$:\\
$\Rightarrow\langle h_2,h_3\rangle\cong BS(1,-1)$ is solvable, a contradiction.

\item[(50)]$ R_1:h_2^2h_3^{-1}h_2^{-2}h_3=1$, $R_2:h_2^2h_3^{-1}h_2^{-2}h_3=1$, $R_3:h_2h_3^{-1}(h_3^{-1}h_2)^2h_3^2=1$:\\
$\Rightarrow\langle h_2,h_3\rangle\cong BS(1,-1)$ is solvable, a contradiction.

\item[(51)]$ R_1:h_2^2h_3^{-1}h_2^{-2}h_3=1$, $R_2:h_2(h_2h_3^{-1})^3h_2=1$, $R_3:h_2^2h_3^{-1}h_2^{-1}(h_3^{-1}h_2)^2=1$:\\
$\Rightarrow \langle h_2,h_3\rangle=\langle h_2 \rangle$ is abelian, a contradiction.

\item[(52)]$ R_1:h_2^2h_3^{-1}h_2^{-2}h_3=1$, $R_2:h_2(h_2h_3^{-1})^3h_2=1$, $R_3:h_2^2h_3^{-2}(h_2^{-1}h_3)^2=1$:\\
$\Rightarrow\langle h_2,h_3\rangle\cong BS(1,1)$ is solvable, a contradiction.

\item[(53)]$ R_1:h_2^2h_3^{-1}h_2^{-2}h_3=1$, $R_2:(h_2h_3)^2h_2^{-1}h_3=1$, $R_3:h_2h_3^{-1}h_2h_3h_2^{-1}h_3^3=1$:\\
$R_2\Rightarrow R_3: h_3^2h_2^{-1}h_3^{-2}h_2=1$ and $R_1:h_2^2h_3^{-1}h_2^{-2}h_3=1$. Using Tietze transformation where $h_3\mapsto h_3h_2^{-1}$ and $h_2\mapsto h_2$, we have $R_1: h_2^2h_3^{-1}h_2^{-2}h_3=1$ and $R_3:(h_3h_2^{-1})^2=(h_2^{-1}h_3)^2$. So, $R_1\Rightarrow h_2^2\in Z(G)$ where $G=\langle h_2,h_3\rangle$. Let $x=h_3h_2^{-1}$. So $h_2^{-1}h_3=x^{h_2}$. Also 
$(h_3h_2^{-1})^2=(h_2^{-1}h_3)^2 \text{ so if } H=\langle x,x^{h_2}\rangle\Rightarrow H\cong BS(1,-1)$ is solvable. By Corollary \ref{n-sub2} $H \trianglelefteq G$ since $G=\langle x,h_2\rangle$. Since $\frac{G}{H}=\frac{\langle h_2\rangle H}{H}$ is a cyclic group, it is solvable. So $G$ is solvable, a contradiction.

\item[(54)]$ R_1:h_2^2h_3^{-1}h_2^{-2}h_3=1$, $R_2:h_2h_3^{-1}(h_2^{-1}h_3)^2h_3=1$, $R_3:h_2^2h_3^{-3}h_2^{-1}h_3=1$:\\
$\Rightarrow \langle h_2,h_3\rangle=\langle h_3 \rangle$ is abelian, a contradiction.

\item[(55)]$ R_1:h_2^2h_3^{-1}h_2^{-2}h_3=1$, $R_2:h_2h_3^{-1}(h_2^{-1}h_3)^2h_3=1$, $R_3:h_2^2h_3^{-2}(h_2^{-1}h_3)^2=1$:\\
$\Rightarrow\langle h_2,h_3\rangle\cong BS(1,1)$ is solvable, a contradiction.

\item[(56)]$ R_1:h_2^2h_3^{-1}h_2^{-2}h_3=1$, $R_2:h_2h_3^{-4}h_2=1$, $R_3:h_2^4h_3^2=1$:\\
By interchanging $h_2$ and $h_3$ in (8) and with the same discussion, there is a contradiction.

\item[(57)]$ R_1:h_2^2h_3^{-1}h_2^{-2}h_3=1$, $R_2:h_2h_3^{-4}h_2=1$, $R_3:h_2h_3h_2^{-1}h_3^{-1}h_2^{-1}h_3^2=1$:\\
$\Rightarrow\langle h_2,h_3\rangle\cong BS(2,1)$ is solvable, a contradiction.

\item[(58)]$ R_1:h_2^2h_3^{-1}h_2^{-2}h_3=1$, $R_2:h_2h_3^{-4}h_2=1$, $R_3:h_2h_3h_2^{-1}h_3^{-1}h_2h_3^{-1}h_2^{-1}h_3=1$:\\
$R_1\Rightarrow h_2^2\in Z(G)$ where $G=\langle h_2,h_3\rangle$. Let $x=h_3^{-1}h_2^{-1}$. So $h_2^{-1}h_3^{-1}=x^{h_2}$. \\
$R_1 \text{ and } R_3 \Rightarrow (h_3^{-1}h_2^{-1})^2=(h_2^{-1}h_3^{-1})^2 \text{ so if } H=\langle x,x^{h_2}\rangle\Rightarrow H\cong BS(1,-1)$ is solvable. By Corollary \ref{n-sub2} $H \trianglelefteq G$ since $G=\langle x,h_2\rangle$. Since $\frac{G}{H}=\frac{\langle h_2\rangle H}{H}$ is a cyclic group, it is solvable. So $G$ is solvable, a contradiction.

\item[(59)]$ R_1:h_2^2h_3^{-1}h_2^{-2}h_3=1$, $R_2:h_2h_3^{-3}h_2h_3=1$, $R_3:h_2^2h_3^{-1}h_2^{-1}(h_3^{-1}h_2)^2=1$:\\
$\Rightarrow \langle h_2,h_3\rangle=\langle h_2 \rangle$ is abelian, a contradiction.

\item[(60)]$ R_1:h_2^2h_3^{-1}h_2^{-2}h_3=1$, $R_2:h_2h_3^{-3}h_2h_3=1$, $R_3:h_2(h_3h_2^{-2})^2h_3=1$:\\
$\Rightarrow \langle h_2,h_3\rangle=\langle h_2 \rangle$ is abelian, a contradiction.

\item[(61)]$ R_1:h_2^2h_3^{-1}h_2^{-2}h_3=1$, $R_2:h_2h_3^{-2}(h_3^{-1}h_2)^2=1$, $R_3:h_2h_3h_2^{-1}h_3^{-1}h_2h_3^{-1}h_2^{-1}h_3=1$:\\
With the same discussion such as item (58), there is a contradiction.

\item[(62)]$ R_1:h_2^2h_3^{-1}h_2^{-2}h_3=1$, $R_2:h_2h_3^{-1}h_2h_3^2h_2^{-1}h_3=1$, $R_3:h_2^2h_3^{-1}h_2^{-1}h_3^3=1$:\\
$\Rightarrow \langle h_2,h_3\rangle=\langle h_3 \rangle$ is abelian, a contradiction.

\item[(63)]$ R_1:h_2^2h_3^{-1}h_2^{-2}h_3=1$, $R_2:h_2h_3^{-1}h_2h_3^2h_2^{-1}h_3=1$, $R_3:h_2h_3h_2h_3^{-1}h_2^{-1}h_3^2=1$:\\
$\Rightarrow\langle h_2,h_3\rangle\cong BS(-2,1)$ is solvable, a contradiction.

\item[(64)]$ R_1:h_2^2h_3^{-1}h_2^{-2}h_3=1$, $R_2:h_2h_3^{-1}h_2(h_3h_2^{-1})^2h_3=1$, $R_3:h_2h_3^{-2}h_2^{-1}h_3^3=1$:\\
With the same discussion such as item (43), there is a contradiction.

\item[(65)]$ R_1:h_2^2h_3^{-1}h_2^{-2}h_3=1$, $R_2:h_2h_3^{-1}h_2(h_3h_2^{-1})^2h_3=1$, $R_3:h_2h_3^{-3}h_2^{-1}h_3^2=1$:\\
With the same discussion such as item (45), there is a contradiction.

\item[(66)]$ R_1:h_2^2h_3^{-1}h_2^{-2}h_3=1$, $R_2:(h_2h_3^{-1})^2h_2^{-1}h_3^2=1$, $R_3:h_2^2h_3h_2^2h_3^{-1}h_2=1$:\\
$\Rightarrow h_2=1$ that is a contradiction.

\item[(67)]$ R_1:h_2^2h_3^{-1}h_2^{-2}h_3=1$, $R_2:(h_2h_3^{-1})^2h_2h_3^2=1$, $R_3:h_2^2h_3h_2^2h_3^{-1}h_2=1$:\\
With the same discussion such as item (66), there is a contradiction.

\item[(68)]$ R_1:h_2^2h_3^{-1}h_2^{-2}h_3=1$, $R_2:(h_2h_3^{-1})^2h_2h_3^2=1$, $R_3:h_2^2h_3^{-1}h_2^2h_3^2=1$:\\
$\Rightarrow \langle h_2,h_3\rangle=\langle h_2 \rangle$ is abelian, a contradiction.


\item[(70)]$ R_1:h_2^2h_3^{-1}h_2^{-1}h_3^{-1}h_2=1$, $R_2:h_2^2h_3^{-1}h_2^{-1}h_3^{-1}h_2=1$, $R_3:h_2h_3^{-3}h_2h_3h_2^{-1}h_3=1$:\\
$\Rightarrow\langle h_2,h_3\rangle\cong BS(1,-1)$ is solvable, a contradiction.

\item[(71)]$ R_1:h_2^2h_3^{-1}h_2^{-1}h_3^{-1}h_2=1$, $R_2:h_2^2h_3^{-1}h_2^{-1}h_3^{-1}h_2=1$, $R_3:h_2h_3^{-2}(h_3^{-1}h_2)^2h_3=1$:\\
$\Rightarrow h_3=1$ that is a contradiction.

\item[(72)]$ R_1:h_2^2h_3^{-1}h_2^{-1}h_3^{-1}h_2=1$, $R_2:h_2^2h_3^{-2}h_2^{-1}h_3=1$, $R_3:h_2h_3h_2^{-1}h_3^{-1}(h_3^{-1}h_2)^2=1$:\\
$\Rightarrow\langle h_2,h_3\rangle\cong BS(2,1)$ is solvable, a contradiction.

\item[(73)]$ R_1:h_2^2h_3^{-1}h_2^{-1}h_3^{-1}h_2=1$, $R_2:h_2^2h_3^{-2}h_2^{-1}h_3=1$, $R_3:h_2h_3^{-1}h_2^{-1}h_3h_2^{-2}h_3^2=1$:\\
$\Rightarrow \langle h_2,h_3\rangle=\langle h_2 \rangle$ is abelian, a contradiction.

\item[(74)]$ R_1:h_2^2h_3^{-1}h_2^{-1}h_3^{-1}h_2=1$, $R_2:h_2h_3h_2^{-1}h_3^{-2}h_2=1$, $R_3:h_2^2h_3^{-2}(h_2^{-1}h_3)^2=1$:\\
$\Rightarrow\langle h_2,h_3\rangle\cong BS(-2,1)$ is solvable, a contradiction.

\item[(75)]$ R_1:h_2^2h_3^{-1}h_2^{-1}h_3^{-1}h_2=1$, $R_2:h_2h_3h_2^{-1}h_3^{-2}h_2=1$, $R_3:h_2^2h_3^{-2}h_2h_3^{-1}h_2^{-1}h_3=1$:\\
$\Rightarrow \langle h_2,h_3\rangle=\langle h_2 \rangle$ is abelian, a contradiction.

\item[(76)]$ R_1:h_2^2h_3^{-1}h_2^{-1}h_3^{-1}h_2=1$, $R_2:h_2h_3^{-2}h_2^{-1}h_3^2=1$, $R_3:h_2^2h_3^{-1}h_2^{-1}(h_2^{-1}h_3)^2=1$:\\
$\Rightarrow \langle h_2,h_3\rangle=\langle h_2 \rangle$ is abelian, a contradiction.

\item[(77)]$ R_1:h_2^2h_3^{-1}h_2^{-1}h_3^{-1}h_2=1$, $R_2:h_2h_3^{-2}h_2^{-1}h_3^2=1$, $R_3:h_2^2h_3^{-1}h_2^{-1}h_3h_2^{-2}h_3=1$:\\
$\Rightarrow \langle h_2,h_3\rangle=\langle h_2 \rangle$ is abelian, a contradiction.

\item[(78)]$ R_1:h_2^2h_3^{-1}h_2^{-1}h_3^{-1}h_2=1$, $R_2:h_2h_3^{-3}h_2h_3=1$, $R_3:h_2^2h_3^{-1}h_2^{-1}(h_3^{-1}h_2)^2=1$:\\
$\Rightarrow \langle h_2,h_3\rangle=\langle h_2 \rangle$ is abelian, a contradiction.

\item[(79)]$ R_1:h_2^2h_3^{-1}h_2^{-1}h_3^{-1}h_2=1$, $R_2:h_2h_3^{-3}h_2h_3=1$, $R_3:h_2(h_3h_2^{-2})^2h_3=1$:\\
$\Rightarrow \langle h_2,h_3\rangle=\langle h_2 \rangle$ is abelian, a contradiction.

\item[(80)]$ R_1:h_2^2h_3^{-1}h_2^{-1}h_3^{-1}h_2=1$, $R_2:h_2(h_3^{-1}h_2h_3^{-1})^2h_2=1$, $R_3:h_2^2h_3^{-1}h_2h_3h_2^{-1}h_3^2=1$:\\
$\Rightarrow \langle h_2,h_3\rangle=\langle h_3 \rangle$ is abelian, a contradiction.

\item[(81)]$ R_1:h_2^2h_3^{-1}h_2^{-1}h_3^{-1}h_2=1$, $R_2:(h_2h_3^{-1})^2h_2h_3h_2^{-1}h_3=1$, $R_3:(h_2h_3)^2(h_2^{-1}h_3)^2=1$:\\
$\Rightarrow \langle h_2,h_3\rangle=\langle h_3 \rangle$ is abelian, a contradiction.

\item[(82)]$ R_1:h_2^2h_3^{-1}h_2^{-1}h_3^2=1$, $R_2:h_2^2h_3^{-1}h_2^{-1}h_3^2=1$, $R_3:h_2^2h_3^{-1}(h_3^{-1}h_2)^2h_3=1$:\\
$\Rightarrow \langle h_2,h_3\rangle=\langle h_3 \rangle$ is abelian, a contradiction.

\item[(83)]$ R_1:h_2^2h_3^{-1}h_2^{-1}h_3^2=1$, $R_2:h_2^2h_3^{-1}h_2^{-1}h_3^2=1$, $R_3:h_2h_3h_2^{-1}h_3h_2h_3^{-1}h_2h_3=1$:\\
$\Rightarrow\langle h_2,h_3\rangle\cong BS(-3,1)$ is solvable, a contradiction.

\item[(84)]$ R_1:h_2^2h_3^{-1}h_2^{-1}h_3^2=1$, $R_2:h_2^2h_3^{-1}h_2^{-1}h_3^2=1$, $R_3:h_2h_3^{-2}h_2h_3h_2^{-1}h_3^2=1$:\\
$\Rightarrow \langle h_2,h_3\rangle=\langle h_2 \rangle$ is abelian, a contradiction.

\item[(85)]$ R_1:h_2^2h_3^{-1}h_2^{-1}h_3^2=1$, $R_2:h_2^2h_3^{-1}h_2^{-1}h_3^2=1$, $R_3:h_2^3h_3^3=1$:\\
$\Rightarrow \langle h_2,h_3\rangle=\langle h_2 \rangle$ is abelian, a contradiction.

\item[(86)]$ R_1:h_2^2h_3^{-1}(h_2^{-1}h_3)^2=1$, $R_2:h_2^2h_3^{-1}(h_2^{-1}h_3)^2=1$, $R_3:h_2^2(h_3^{-1}h_2^{-1})^2h_3=1$:\\
$\Rightarrow h_3=1$ that is a contradiction.

\item[(87)]$ R_1:h_2^2h_3^{-1}(h_2^{-1}h_3)^2=1$, $R_2:h_2^2h_3^{-1}(h_2^{-1}h_3)^2=1$, $R_3:h_2(h_3h_2^{-1})^2h_3^{-1}h_2h_3=1$:\\
$\Rightarrow \langle h_2,h_3\rangle=\langle h_2 \rangle$ is abelian, a contradiction.

\item[(88)]$ R_1:h_2^2h_3^{-1}(h_2^{-1}h_3)^2=1$, $R_2:h_2^2h_3^{-1}(h_2^{-1}h_3)^2=1$, $R_3:h_2h_3^{-2}h_2^{-1}h_3^3=1$:\\
$\Rightarrow\langle h_2,h_3\rangle\cong BS(1,1)$ is solvable, a contradiction.

\item[(89)]$ R_1:h_2^2h_3^{-1}(h_2^{-1}h_3)^2=1$, $R_2:h_2^2h_3^{-1}(h_2^{-1}h_3)^2=1$, $R_3:h_2h_3^{-3}h_2^{-1}h_3^2=1$:\\
$\Rightarrow\langle h_2,h_3\rangle\cong BS(1,-4)$ is solvable, a contradiction.

\item[(90)]$ R_1:h_2^2h_3^{-1}(h_2^{-1}h_3)^2=1$, $R_2:h_2h_3^2h_2^{-1}h_3^2=1$, $R_3:h_2^2h_3h_2^{-2}h_3^2=1$:\\
$\Rightarrow\langle h_2,h_3\rangle\cong BS(2,1)$ is solvable, a contradiction.

\item[(91)]$ R_1:h_2^2h_3^{-1}(h_2^{-1}h_3)^2=1$, $R_2:h_2h_3^2h_2^{-1}h_3^2=1$, $R_3:h_2^2h_3h_2^{-1}h_3^{-1}h_2^{-1}h_3=1$:\\
$\Rightarrow h_3=1$ that is a contradiction.

\item[(92)]$ R_1:h_2^2h_3^{-1}(h_2^{-1}h_3)^2=1$, $R_2:h_2h_3^2h_2^{-1}h_3^2=1$, $R_3:h_2^2(h_3h_2^{-1}h_3)^2=1$:\\
$\Rightarrow\langle h_2,h_3\rangle\cong BS(1,-2)$ is solvable, a contradiction.

\item[(93)]$ R_1:h_2^2h_3^{-1}(h_2^{-1}h_3)^2=1$, $R_2:h_2h_3^2h_2^{-1}h_3^2=1$, $R_3:h_2^2(h_3h_2^{-1})^2h_3^2=1$:\\
$\Rightarrow\langle h_2,h_3\rangle\cong BS(1,1)$ is solvable, a contradiction.

\item[(94)]$ R_1:h_2^2h_3^{-2}h_2^{-1}h_3=1$, $R_2:h_2^2h_3^{-2}h_2^{-1}h_3=1$, $R_3:h_2h_3h_2^{-1}h_3^{-1}(h_2^{-1}h_3)^2=1$:\\
$\Rightarrow\langle h_2,h_3\rangle\cong BS(1,-3)$ is solvable, a contradiction.

\item[(95)]$ R_1:h_2^2h_3^{-2}h_2^{-1}h_3=1$, $R_2:h_2^2h_3^{-2}h_2^{-1}h_3=1$, $R_3:h_2h_3^{-1}h_2^{-1}h_3h_2^{-1}h_3^{-1}h_2h_3=1$:\\
$\Rightarrow\langle h_2,h_3\rangle\cong BS(1,1)$ is solvable, a contradiction.

\item[(96)]$ R_1:h_2^2h_3^{-2}h_2^{-1}h_3=1$, $R_2:h_2^2h_3^{-2}h_2^{-1}h_3=1$, $R_3:h_2h_3h_2^{-1}h_3^{-1}h_2h_3^{-1}h_2^{-1}h_3=1$:\\
By interchanging $h_2$ and $h_3$ in (95) and with the same discussion, there is a contradiction.

\item[(97)]$ R_1:h_2^2h_3^{-2}h_2^{-1}h_3=1$, $R_2:h_2h_3^{-3}h_2h_3=1$, $R_3:h_2^2h_3^{-1}h_2^{-1}(h_2^{-1}h_3)^2=1$:\\
$\Rightarrow\langle h_2,h_3\rangle\cong BS(2,1)$ is solvable, a contradiction.

\item[(98)]$ R_1:h_2^2h_3^{-2}h_2^{-1}h_3=1$, $R_2:h_2h_3^{-3}h_2h_3=1$, $R_3:h_2^2h_3^{-1}h_2^{-1}h_3h_2^{-2}h_3=1$:\\
$\Rightarrow\langle h_2,h_3\rangle\cong BS(2,1)$ is solvable, a contradiction.


\item[(100)]$ R_1:h_2^2h_3^{-2}h_2h_3=1$, $R_2:h_2^2h_3^{-2}h_2h_3=1$, $R_3:h_2h_3^{-1}h_2^{-1}h_3h_2^{-1}h_3^{-1}h_2h_3=1$:\\
$\Rightarrow\langle h_2,h_3\rangle\cong BS(1,1)$ is solvable, a contradiction.

\item[(101)]$ R_1:h_2^2h_3^{-1}h_2^2h_3=1$, $R_2:h_2^2h_3^{-1}h_2^2h_3=1$, $R_3:h_2^2h_3^2h_2^{-2}h_3=1$:\\
$\Rightarrow h_3=1$ that is a contradiction.

\item[(102)]$ R_1:h_2^2h_3^{-1}h_2^2h_3=1$, $R_2:h_2^2h_3^{-1}h_2^2h_3=1$, $R_3:h_2h_3(h_3h_2^{-1})^2h_2^{-1}h_3=1$:\\
$\Rightarrow \langle h_2,h_3\rangle=\langle h_3 \rangle$ is abelian, a contradiction.

\item[(103)]$ R_1:h_2^2h_3^{-1}h_2^2h_3=1$, $R_2:h_2^2h_3^{-1}h_2^2h_3=1$, $R_3:h_2h_3h_2^{-1}h_3(h_2h_3^{-1})^2h_2=1$:\\
$\Rightarrow\langle h_2,h_3\rangle\cong BS(2,1)$ is solvable, a contradiction.

\item[(104)]$ R_1:h_2^2h_3^{-1}h_2^2h_3=1$, $R_2:h_2^2h_3^{-1}h_2^2h_3=1$, $R_3:h_2(h_3h_2^{-1})^2(h_3^{-1}h_2)^2=1$:\\
$\Rightarrow\langle h_2,h_3\rangle\cong BS(-2,1)$ is solvable, a contradiction.

\item[(105)]$ R_1:h_2^2h_3^{-1}h_2^2h_3=1$, $R_2:h_2h_3^2h_2^{-1}h_3^2=1$, $R_3:h_3^3(h_3h_2^{-1})^2h_3=1$:\\
$\Rightarrow h_2=1$ that is a contradiction.

\item[(106)]$ R_1:h_2^2h_3^{-1}h_2^2h_3=1$, $R_2:h_2h_3^{-2}(h_2^{-1}h_3)^2=1$, $R_3:h_2^2h_3h_2^{-2}h_3^{-1}h_2=1$:\\
$\Rightarrow h_2=1$ that is a contradiction.

\item[(107)]$ R_1:h_2^2h_3^{-1}h_2^2h_3=1$, $R_2:h_2h_3^{-2}h_2h_3h_2^{-1}h_3=1$, $R_3:h_2^2h_3h_2^{-2}h_3^{-1}h_2=1$:\\
$\Rightarrow h_2=1$ that is a contradiction.

\item[(108)]$ R_1:h_2^2h_3^{-1}h_2^2h_3=1$, $R_2:h_2h_3^{-1}h_2(h_3h_2^{-1})^2h_3=1$, $R_3:h_2(h_2h_3^{-1})^2h_3^{-1}h_2h_3=1$:\\
$\Rightarrow h_3=1$ that is a contradiction.

\item[(109)]$ R_1:h_2^2h_3^{-1}h_2^2h_3=1$, $R_2:h_2(h_3^{-1}h_2h_3^{-1})^2h_2=1$, $R_3:h_2h_3^{-1}h_2^{-1}h_3^2h_2h_3^{-1}h_2=1$:\\
$\Rightarrow h_3=1$ that is a contradiction.

\item[(110)]$ R_1:h_2^2h_3^{-1}h_2^2h_3=1$, $R_2:h_2(h_3^{-1}h_2h_3^{-1})^2h_2=1$, $R_3:h_2h_3^{-2}h_2^{-1}h_3h_2h_3^{-1}h_2=1$:\\
$\Rightarrow h_3=1$ that is a contradiction.

\item[(111)]$ R_1:h_2^2h_3^{-1}h_2^2h_3=1$, $R_2:h_2(h_3^{-1}h_2h_3^{-1})^2h_2=1$, $R_3:(h_2h_3^{-1}h_2)^2h_3^2=1$:\\
$\Rightarrow h_3=1$ that is a contradiction.


\item[(113)]$ R_1:h_2^2h_3^{-1}h_2h_3h_2^{-1}h_3=1$, $R_2:h_2^2h_3^{-1}h_2h_3h_2^{-1}h_3=1$, $R_3:(h_2h_3^{-1})^2h_2h_3^2h_2^{-1}h_3=1$:\\
$\Rightarrow \langle h_2,h_3\rangle=\langle h_2 \rangle$ is abelian, a contradiction.

\item[(114)]$ R_1:h_2(h_2h_3^{-1})^2h_2^{-1}h_3=1$, $R_2:h_2(h_2h_3^{-1})^2h_2^{-1}h_3=1$, $R_3:h_2^2h_3^{-1}h_2h_3^2h_2^{-1}h_3=1$:\\
$\Rightarrow h_3=1$ that is a contradiction.

\item[(115)]$ R_1:h_2(h_2h_3^{-1})^2h_2^{-1}h_3=1$, $R_2:h_2(h_2h_3^{-1})^2h_2^{-1}h_3=1$, $R_3:h_2h_3^{-1}h_2^{-1}h_3h_2^{-1}h_3^{-1}h_2h_3=1$:\\
$\Rightarrow\langle h_2,h_3\rangle\cong BS(1,1)$ is solvable, a contradiction.

\item[(116)]$ R_1:h_2(h_2h_3^{-1})^2h_2^{-1}h_3=1$, $R_2:h_2h_3^{-2}h_2^{-1}h_3^2=1$, $R_3:h_2h_3^{-1}h_2^{-1}h_3^2h_2^{-1}h_3^{-1}h_2=1$:\\
$R_2\Rightarrow R_3: h_2^2h_3^{-1}h_2^{-2}h_3=1$ and $R_2:h_2h_3^{-2}h_2^{-1}h_3^2=1$. Using Tietze transformation where $h_3\mapsto h_3h_2^{-1}$ and $h_2\mapsto h_2$, we have $R_3: h_2^2h_3^{-1}h_2^{-2}h_3=1$ and $R_2:(h_3h_2^{-1})^2=(h_2^{-1}h_3)^2$. So, $R_3\Rightarrow h_2^2\in Z(G)$ where $G=\langle h_2,h_3\rangle$. Let $x=h_3h_2^{-1}$. So $h_2^{-1}h_3=x^{h_2}$. Also 
$(h_3h_2^{-1})^2=(h_2^{-1}h_3)^2 \text{ so if } H=\langle x,x^{h_2}\rangle\Rightarrow H\cong BS(1,-1)$ is solvable. By Corollary \ref{n-sub2} $H \trianglelefteq G$ since $G=\langle x,h_2\rangle$. Since $\frac{G}{H}=\frac{\langle h_2\rangle H}{H}$ is a cyclic group, it is solvable. So $G$ is solvable, a contradiction.

\item[(117)]$ R_1:h_2(h_2h_3^{-1})^2h_3^{-1}h_2=1$, $R_2:h_2h_3^{-2}h_2^{-1}h_3^2=1$, $R_3:h_2h_3^{-1}h_2^{-1}h_3^2h_2^{-1}h_3^{-1}h_2=1$:\\
With the same discussion such as item (116), there is a contradiction.

\item[(118)]$ R_1:h_2(h_2h_3^{-1})^2h_2h_3=1$, $R_2:h_2(h_2h_3^{-1})^2h_2h_3=1$, $R_3:h_2h_3^{-1}h_2^{-1}h_3h_2^{-1}h_3^{-1}h_2h_3=1$:\\
$\Rightarrow\langle h_2,h_3\rangle\cong BS(1,1)$ is solvable, a contradiction.

\item[(119)]$ R_1:h_2h_3(h_2h_3^{-1})^2h_2=1$, $R_2:h_2h_3^{-2}h_2^{-1}h_3^2=1$, $R_3:h_2^3h_3^{-2}h_2h_3=1$:\\
$\Rightarrow \langle h_2,h_3\rangle=\langle h_2 \rangle$ is abelian, a contradiction.


\item[(121)]$ R_1:h_2h_3^2h_2^{-2}h_3=1$, $R_2:h_2h_3^2h_2^{-2}h_3=1$, $R_3:h_2h_3h_2^{-1}h_3^{-1}h_2h_3^{-1}h_2^{-1}h_3=1$:\\
By interchanging $h_2$ and $h_3$ in (100) and with the same discussion, there is a contradiction.

\item[(122)]$ R_1:h_2h_3^2h_2^{-1}h_3^{-1}h_2=1$, $R_2:h_2h_3^2h_2^{-1}h_3^{-1}h_2=1$, $R_3:h_2h_3h_2^{-2}h_3h_2h_3^{-1}h_2=1$:\\
By interchanging $h_2$ and $h_3$ in (84) and with the same discussion, there is a contradiction.

\item[(123)]$ R_1:h_2h_3^2h_2^{-1}h_3^{-1}h_2=1$, $R_2:h_2h_3^2h_2^{-1}h_3^{-1}h_2=1$, $R_3:h_2h_3h_2h_3^{-1}h_2h_3h_2^{-1}h_3=1$:\\
By interchanging $h_2$ and $h_3$ in (83) and with the same discussion, there is a contradiction.

\item[(124)]$ R_1:h_2h_3^2h_2^{-1}h_3^{-1}h_2=1$, $R_2:h_2h_3^2h_2^{-1}h_3^{-1}h_2=1$, $R_3:h_2h_3^2h_2^{-1}(h_2^{-1}h_3)^2=1$:\\
By interchanging $h_2$ and $h_3$ in (82) and with the same discussion, there is a contradiction.

\item[(125)]$ R_1:h_2h_3^2h_2^{-1}h_3^{-1}h_2=1$, $R_2:h_2h_3^2h_2^{-1}h_3^{-1}h_2=1$, $R_3:h_2^3h_3^3=1$:\\
By interchanging $h_2$ and $h_3$ in (85) and with the same discussion, there is a contradiction.

\item[(126)]$ R_1:h_2h_3^2h_2^{-1}h_3^2=1$, $R_2:h_2h_3^2h_2^{-1}h_3^2=1$, $R_3:h_2h_3^{-1}h_2(h_3h_2^{-1})^2h_3^2=1$:\\
By interchanging $h_2$ and $h_3$ in (103) and with the same discussion, there is a contradiction.

\item[(127)]$ R_1:h_2h_3^2h_2^{-1}h_3^2=1$, $R_2:h_2h_3^2h_2^{-1}h_3^2=1$, $R_3:(h_2h_3^{-1})^2(h_2^{-1}h_3)^2h_3=1$:\\
By interchanging $h_2$ and $h_3$ in (104) and with the same discussion, there is a contradiction.

\item[(128)]$ R_1:h_2h_3^2h_2^{-1}h_3^2=1$, $R_2:h_2h_3^2h_2^{-1}h_3^2=1$, $R_3:h_2^2h_3^{-2}h_2h_3^2=1$:\\
By interchanging $h_2$ and $h_3$ in (101) and with the same discussion, there is a contradiction.

\item[(129)]$ R_1:h_2h_3^2h_2^{-1}h_3^2=1$, $R_2:h_2h_3^2h_2^{-1}h_3^2=1$, $R_3:h_2(h_2h_3^{-1})^2h_3^{-1}h_2h_3=1$:\\
By interchanging $h_2$ and $h_3$ in (102) and with the same discussion, there is a contradiction.

\item[(130)]$ R_1:h_2h_3^2h_2^{-1}h_3^2=1$, $R_2:h_2h_3^{-1}h_2^{-1}h_3h_2^{-1}h_3^{-1}h_2=1$, $R_3:h_2h_3^{-2}h_2^{-1}h_3^3=1$:\\
By interchanging $h_2$ and $h_3$ in (107) and with the same discussion, there is a contradiction.

\item[(131)]$ R_1:h_2h_3^2h_2^{-1}h_3^2=1$, $R_2:h_2(h_3^{-1}h_2h_3^{-1})^2h_2=1$, $R_3:h_2^2(h_3h_2^{-1}h_3)^2=1$:\\
By interchanging $h_2$ and $h_3$ in (111) and with the same discussion, there is a contradiction.

\item[(132)]$ R_1:h_2h_3^2h_2^{-1}h_3^2=1$, $R_2:h_2(h_3^{-1}h_2h_3^{-1})^2h_2=1$, $R_3:h_2^2h_3^{-2}h_2h_3^{-1}h_2^{-1}h_3=1$:\\
By interchanging $h_2$ and $h_3$ in (110) and with the same discussion, there is a contradiction.

\item[(133)]$ R_1:h_2h_3^2h_2^{-1}h_3^2=1$, $R_2:h_2(h_3^{-1}h_2h_3^{-1})^2h_2=1$, $R_3:h_2h_3h_2^{-1}h_3^2h_2^{-1}h_3^{-1}h_2=1$:\\
By interchanging $h_2$ and $h_3$ in (109) and with the same discussion, there is a contradiction.

\item[(134)]$ R_1:h_2h_3^2h_2^{-1}h_3^2=1$, $R_2:(h_2h_3^{-1})^2h_2h_3h_2^{-1}h_3=1$, $R_3:h_2h_3(h_3h_2^{-1})^2h_2^{-1}h_3=1$:\\
By interchanging $h_2$ and $h_3$ in (108) and with the same discussion, there is a contradiction.

\item[(135)]$ R_1:h_2h_3(h_3h_2^{-1})^2h_3=1$, $R_2:h_2h_3(h_3h_2^{-1})^2h_3=1$, $R_3:h_2h_3h_2^{-1}h_3^{-1}h_2h_3^{-1}h_2^{-1}h_3=1$:\\
By interchanging $h_2$ and $h_3$ in (118) and with the same discussion, there is a contradiction.

\item[(136)]$ R_1:h_2h_3h_2^{-1}(h_2^{-1}h_3)^2=1$, $R_2:h_2h_3h_2^{-1}(h_2^{-1}h_3)^2=1$, $R_3:h_2h_3^3h_2^{-1}h_3^{-1}h_2=1$:\\
$\Rightarrow h_3=1$ that is a contradiction.

\item[(137)]$ R_1:h_2h_3h_2^{-1}(h_2^{-1}h_3)^2=1$, $R_2:h_2h_3h_2^{-1}(h_2^{-1}h_3)^2=1$, $R_3:h_2h_3^{-4}h_2h_3=1$:\\
$\Rightarrow\langle h_2,h_3\rangle\cong BS(1,1)$ is solvable, a contradiction.

\item[(138)]$ R_1:h_2h_3h_2^{-1}h_3^{-2}h_2=1$, $R_2:h_2h_3h_2^{-1}h_3^{-2}h_2=1$, $R_3:h_2^2h_3^{-2}(h_3^{-1}h_2)^2=1$:\\
$\Rightarrow\langle h_2,h_3\rangle\cong BS(1,3)$ is solvable, a contradiction.

\item[(139)]$ R_1:h_2h_3h_2^{-1}h_3^{-2}h_2=1$, $R_2:h_2h_3h_2^{-1}h_3^{-2}h_2=1$, $R_3:h_2(h_2h_3^{-2})^2h_2=1$:\\
$\Rightarrow\langle h_2,h_3\rangle\cong BS(1,-1)$ is solvable, a contradiction.

\item[(140)]$ R_1:h_2h_3h_2^{-1}h_3^{-2}h_2=1$, $R_2:h_2h_3h_2^{-1}h_3^{-2}h_2=1$, $R_3:h_2h_3^{-3}h_2^2h_3^{-1}h_2=1$:\\
By interchanging $h_2$ and $h_3$ in (139) and with the same discussion, there is a contradiction.

\item[(141)]$ R_1:h_2h_3h_2^{-1}h_3^{-2}h_2=1$, $R_2:h_2h_3^{-3}h_2h_3=1$, $R_3:h_2^2h_3^{-1}h_2^{-1}(h_2^{-1}h_3)^2=1$:\\
$\Rightarrow\langle h_2,h_3\rangle\cong BS(1,1)$ is solvable, a contradiction.

\item[(142)]$ R_1:h_2h_3h_2^{-1}h_3^{-2}h_2=1$, $R_2:h_2h_3^{-3}h_2h_3=1$, $R_3:h_2h_3h_2^{-2}h_3h_2^{-1}h_3^{-1}h_2=1$:\\
$\Rightarrow \langle h_2,h_3\rangle=\langle h_2 \rangle$ is abelian, a contradiction.

\item[(143)]$ R_1:h_2(h_3h_2^{-1})^2h_2^{-1}h_3=1$, $R_2:h_2h_3^{-2}h_2^{-1}h_3^2=1$, $R_3:h_2h_3^{-1}h_2^{-1}h_3^2h_2^{-1}h_3^{-1}h_2=1$:\\
$R_2\Rightarrow R_3: h_2^2h_3^{-1}h_2^{-2}h_3=1$ and $R_2:h_2h_3^{-2}h_2^{-1}h_3^2=1$. Using Tietze transformation where $h_3\mapsto h_3h_2^{-1}$ and $h_2\mapsto h_2$, we have $R_3: h_2^2h_3^{-1}h_2^{-2}h_3=1$ and $R_2:(h_3h_2^{-1})^2=(h_2^{-1}h_3)^2$. So, $R_3\Rightarrow h_2^2\in Z(G)$ where $G=\langle h_2,h_3\rangle$. Let $x=h_3h_2^{-1}$. So $h_2^{-1}h_3=x^{h_2}$. Also 
$(h_3h_2^{-1})^2=(h_2^{-1}h_3)^2 \text{ so if } H=\langle x,x^{h_2}\rangle\Rightarrow H\cong BS(1,-1)$ is solvable. By Corollary \ref{n-sub2} $H \trianglelefteq G$ since $G=\langle x,h_2\rangle$. Since $\frac{G}{H}=\frac{\langle h_2\rangle H}{H}$ is a cyclic group, it is solvable. So $G$ is solvable, a contradiction.

\item[(144)]$ R_1:h_2(h_3h_2^{-1})^2h_3^{-1}h_2=1$, $R_2:h_2(h_3h_2^{-1})^2h_3^{-1}h_2=1$, $R_3:h_2h_3^3h_2^{-1}h_3^2=1$:\\
$\Rightarrow\langle h_2,h_3\rangle\cong BS(1,4)$ is solvable, a contradiction.

\item[(145)]$ R_1:h_2(h_3h_2^{-1})^2h_3^{-1}h_2=1$, $R_2:h_2(h_3h_2^{-1})^2h_3^{-1}h_2=1$, $R_3:h_2h_3^2h_2^{-1}h_3^3=1$:\\
$\Rightarrow\langle h_2,h_3\rangle\cong BS(-2,1)$ is solvable, a contradiction.

\item[(146)]$ R_1:h_2(h_3h_2^{-1})^2h_3^{-1}h_2=1$, $R_2:h_2(h_3h_2^{-1})^2h_3^{-1}h_2=1$, $R_3:h_2^2h_3^{-1}h_2^{-1}h_3^2h_2^{-1}h_3=1$:\\
$\Rightarrow h_3=1$ that is a contradiction.

\item[(147)]$ R_1:h_2(h_3h_2^{-1})^2h_3^{-1}h_2=1$, $R_2:h_2(h_3h_2^{-1})^2h_3^{-1}h_2=1$, $R_3:(h_2h_3^{-1})^2h_2^{-1}h_3^2h_2^{-1}h_3=1$:\\
$\Rightarrow h_3=1$ that is a contradiction.

\item[(148)]$ R_1:h_2(h_3h_2^{-1})^2h_3^{-1}h_2=1$, $R_2:h_2h_3^{-2}h_2^{-1}h_3^2=1$, $R_3:h_2^2h_3h_2^{-2}h_3^2=1$:\\
$\Rightarrow h_3=1$ that is a contradiction.

\item[(149)]$ R_1:h_2(h_3h_2^{-1})^2h_3^{-1}h_2=1$, $R_2:h_2h_3^{-2}h_2^{-1}h_3^2=1$, $R_3:(h_2h_3)^2h_2^{-2}h_3=1$:\\
$\Rightarrow h_3=1$ that is a contradiction.

\item[(150)]$ R_1:h_2(h_3h_2^{-1})^2h_3^{-1}h_2=1$, $R_2:h_2h_3^{-2}h_2^{-1}h_3^2=1$, $R_3:h_2h_3^{-1}h_2^{-1}h_3^2h_2^{-1}h_3^{-1}h_2=1$:\\
$\Rightarrow\langle h_2,h_3\rangle\cong BS(1,2)$ is solvable, a contradiction.

\item[(151)]$ R_1:h_2(h_3h_2^{-1})^2h_3^{-1}h_2=1$, $R_2:h_2h_3^{-2}h_2^{-1}h_3^2=1$, $R_3:h_2h_3^{-1}h_2^{-1}h_3h_2^{-1}h_3^{-2}h_2=1$:\\
$\Rightarrow\langle h_2,h_3\rangle\cong BS(1,-1)$ is solvable, a contradiction.

\item[(152)]$ R_1:h_2(h_3h_2^{-1})^3h_3=1$, $R_2:h_2h_3^{-1}(h_2^{-1}h_3)^2h_3=1$, $R_3:h_2h_3h_2^{-1}h_3^2h_2^{-2}h_3=1$:\\
$\Rightarrow \langle h_2,h_3\rangle=\langle h_3 \rangle$ is abelian, a contradiction.

\item[(153)]$ R_1:h_2(h_3h_2^{-1})^3h_3=1$, $R_2:h_2h_3^{-4}h_2=1$, $R_3:h_2^2h_3^4=1$:\\
$\Rightarrow h_3=1$ that is a contradiction.

\item[(154)]$ R_1:h_2(h_3h_2^{-1})^3h_3=1$, $R_2:h_2h_3^{-4}h_2=1$, $R_3:h_2^2h_3^{-3}h_2^{-1}h_3=1$:\\
$\Rightarrow \langle h_2,h_3\rangle=\langle h_3 \rangle$ is abelian, a contradiction.

\item[(155)]$ R_1:h_2(h_3h_2^{-1})^3h_3=1$, $R_2:h_2(h_3^{-1}h_2h_3^{-1})^2h_2=1$, $R_3:h_2^2h_3^{-1}h_2h_3h_2^{-1}h_3^2=1$:\\
$\Rightarrow h_3=1$ that is a contradiction.

\item[(156)]$ R_1:h_2h_3^{-2}h_2^{-1}h_3^2=1$, $R_2:h_2h_3^{-2}h_2^{-1}h_3^2=1$, $R_3:h_2(h_2h_3)^2h_3=1$:\\
By interchanging $h_2$ and $h_3$ in (33) and with the same discussion, there is a contradiction.

\item[(157)]$ R_1:h_2h_3^{-2}h_2^{-1}h_3^2=1$, $R_2:h_2h_3^{-2}h_2^{-1}h_3^2=1$, $R_3:h_2^2h_3^{-1}(h_2h_3)^2=1$:\\
By interchanging $h_2$ and $h_3$ in (37) and with the same discussion, there is a contradiction.

\item[(158)]$ R_1:h_2h_3^{-2}h_2^{-1}h_3^2=1$, $R_2:h_2h_3^{-2}h_2^{-1}h_3^2=1$, $R_3:(h_2h_3)^2h_2h_3^{-1}h_2=1$:\\
By interchanging $h_2$ and $h_3$ in (39) and with the same discussion, there is a contradiction.

\item[(159)]$ R_1:h_2h_3^{-2}h_2^{-1}h_3^2=1$, $R_2:h_2h_3^{-2}h_2^{-1}h_3^2=1$, $R_3:(h_2h_3)^2h_3h_2^{-1}h_3=1$:\\
By interchanging $h_2$ and $h_3$ in (36) and with the same discussion, there is a contradiction.

\item[(160)]$ R_1:h_2h_3^{-2}h_2^{-1}h_3^2=1$, $R_2:h_2h_3^{-2}h_2^{-1}h_3^2=1$, $R_3:h_2h_3h_2h_3^{-2}h_2h_3=1$:\\
By interchanging $h_2$ and $h_3$ in (38) and with the same discussion, there is a contradiction.

\item[(161)]$ R_1:h_2h_3^{-2}h_2^{-1}h_3^2=1$, $R_2:h_2h_3^{-2}h_2^{-1}h_3^2=1$, $R_3:h_2h_3(h_2h_3^{-1})^2h_2h_3=1$:\\
By interchanging $h_2$ and $h_3$ in (40) and with the same discussion, there is a contradiction.

\item[(162)]$ R_1:h_2h_3^{-2}h_2^{-1}h_3^2=1$, $R_2:h_2h_3^{-2}h_2^{-1}h_3^2=1$, $R_3:h_2h_3^2h_2h_3^{-1}h_2h_3=1$:\\
By interchanging $h_2$ and $h_3$ in (34) and with the same discussion, there is a contradiction.

\item[(163)]$ R_1:h_2h_3^{-2}h_2^{-1}h_3^2=1$, $R_2:h_2h_3^{-2}h_2^{-1}h_3^2=1$, $R_3:h_2h_3^2h_2^{-1}h_3^{-1}h_2h_3=1$:\\
By interchanging $h_2$ and $h_3$ in (35) and with the same discussion, there is a contradiction.

\item[(164)]$ R_1:h_2h_3^{-2}h_2^{-1}h_3^2=1$, $R_2:h_2h_3^{-2}h_2^{-1}h_3^2=1$, $R_3:h_2^3h_3h_2^{-2}h_3=1$:\\
By interchanging $h_2$ and $h_3$ in (47) and with the same discussion, there is a contradiction.

\item[(165)]$ R_1:h_2h_3^{-2}h_2^{-1}h_3^2=1$, $R_2:h_2h_3^{-2}h_2^{-1}h_3^2=1$, $R_3:h_2^3h_3^{-1}h_2^{-2}h_3=1$:\\
By interchanging $h_2$ and $h_3$ in (45) and with the same discussion, there is a contradiction.

\item[(166)]$ R_1:h_2h_3^{-2}h_2^{-1}h_3^2=1$, $R_2:h_2h_3^{-2}h_2^{-1}h_3^2=1$, $R_3:h_2^2h_3h_2^{-3}h_3=1$:\\
By interchanging $h_2$ and $h_3$ in (46) and with the same discussion, there is a contradiction.

\item[(167)]$ R_1:h_2h_3^{-2}h_2^{-1}h_3^2=1$, $R_2:h_2h_3^{-2}h_2^{-1}h_3^2=1$, $R_3:h_2^2h_3h_2^{-2}h_3^{-1}h_2=1$:\\
By interchanging $h_2$ and $h_3$ in (43) and with the same discussion, there is a contradiction.

\item[(168)]$ R_1:h_2h_3^{-2}h_2^{-1}h_3^2=1$, $R_2:h_2h_3^{-2}h_2^{-1}h_3^2=1$, $R_3:h_2^2h_3h_2^{-2}h_3^2=1$:\\
By interchanging $h_2$ and $h_3$ in (41) and with the same discussion, there is a contradiction.

\item[(169)]$ R_1:h_2h_3^{-2}h_2^{-1}h_3^2=1$, $R_2:h_2h_3^{-2}h_2^{-1}h_3^2=1$, $R_3:h_2^2h_3h_2^{-1}(h_2^{-1}h_3)^2=1$:\\
By interchanging $h_2$ and $h_3$ in (50) and with the same discussion, there is a contradiction.

\item[(170)]$ R_1:h_2h_3^{-2}h_2^{-1}h_3^2=1$, $R_2:h_2h_3^{-2}h_2^{-1}h_3^2=1$, $R_3:h_2^2h_3^{-1}h_2^{-2}h_3^2=1$:\\
By interchanging $h_2$ and $h_3$ in (42) and with the same discussion, there is a contradiction.

\item[(171)]$ R_1:h_2h_3^{-2}h_2^{-1}h_3^2=1$, $R_2:h_2h_3^{-2}h_2^{-1}h_3^2=1$, $R_3:h_2^2h_3^{-1}h_2^{-1}(h_2^{-1}h_3)^2=1$:\\
By interchanging $h_2$ and $h_3$ in (44) and with the same discussion, there is a contradiction.

\item[(172)]$ R_1:h_2h_3^{-2}h_2^{-1}h_3^2=1$, $R_2:h_2h_3^{-2}h_2^{-1}h_3^2=1$, $R_3:h_2h_3h_2^{-2}h_3h_2h_3^{-1}h_2=1$:\\
By interchanging $h_2$ and $h_3$ in (48) and with the same discussion, there is a contradiction.

\item[(173)]$ R_1:h_2h_3^{-2}h_2^{-1}h_3^2=1$, $R_2:h_2h_3^{-2}h_2^{-1}h_3^2=1$, $R_3:h_2h_3h_2^{-2}h_3h_2^{-1}h_3^{-1}h_2=1$:\\
By interchanging $h_2$ and $h_3$ in (49) and with the same discussion, there is a contradiction.

\item[(174)]$ R_1:h_2h_3^{-2}h_2^{-1}h_3^2=1$, $R_2:(h_2h_3^{-1})^2h_2h_3h_2^{-1}h_3=1$, $R_3:h_2^3h_3^{-1}h_2^{-2}h_3=1$:\\
By interchanging $h_2$ and $h_3$ in (65) and with the same discussion, there is a contradiction.

\item[(175)]$ R_1:h_2h_3^{-2}h_2^{-1}h_3^2=1$, $R_2:(h_2h_3^{-1})^2h_2h_3h_2^{-1}h_3=1$, $R_3:h_2^2h_3h_2^{-2}h_3^{-1}h_2=1$:\\
By interchanging $h_2$ and $h_3$ in (64) and with the same discussion, there is a contradiction.

\item[(176)]$ R_1:h_2h_3^{-2}h_2^{-1}h_3^2=1$, $R_2:h_3(h_3h_2^{-1})^3h_3=1$, $R_3:h_2^2h_3^{-2}(h_2^{-1}h_3)^2=1$:\\
By interchanging $h_2$ and $h_3$ in (52) and with the same discussion, there is a contradiction.

\item[(177)]$ R_1:h_2h_3^{-2}h_2^{-1}h_3^2=1$, $R_2:h_3(h_3h_2^{-1})^3h_3=1$, $R_3:h_2h_3^{-2}(h_3^{-1}h_2)^2h_3=1$:\\
By interchanging $h_2$ and $h_3$ in (51) and with the same discussion, there is a contradiction.

\item[(178)]$ R_1:h_2h_3^{-2}(h_2^{-1}h_3)^2=1$, $R_2:h_2h_3^{-2}(h_2^{-1}h_3)^2=1$, $R_3:h_2h_3h_2h_3^{-2}h_2^{-1}h_3=1$:\\
By interchanging $h_2$ and $h_3$ in (86) and with the same discussion, there is a contradiction.

\item[(179)]$ R_1:h_2h_3^{-2}(h_2^{-1}h_3)^2=1$, $R_2:h_2h_3^{-2}(h_2^{-1}h_3)^2=1$, $R_3:h_2h_3(h_2h_3^{-1})^2h_2^{-1}h_3=1$:\\
By interchanging $h_2$ and $h_3$ in (87) and with the same discussion, there is a contradiction.

\item[(180)]$ R_1:h_2h_3^{-2}(h_2^{-1}h_3)^2=1$, $R_2:h_2h_3^{-2}(h_2^{-1}h_3)^2=1$, $R_3:h_2^3h_3^{-1}h_2^{-2}h_3=1$:\\
By interchanging $h_2$ and $h_3$ in (89) and with the same discussion, there is a contradiction.

\item[(181)]$ R_1:h_2h_3^{-2}(h_2^{-1}h_3)^2=1$, $R_2:h_2h_3^{-2}(h_2^{-1}h_3)^2=1$, $R_3:h_2^2h_3h_2^{-2}h_3^{-1}h_2=1$:\\
By interchanging $h_2$ and $h_3$ in (88) and with the same discussion, there is a contradiction.


\item[(183)]$ R_1:h_2h_3^{-3}h_2h_3=1$, $R_2:h_2h_3^{-3}h_2h_3=1$, $R_3:h_2^2h_3^{-1}h_2^{-1}(h_3^{-1}h_2)^2=1$:\\
By interchanging $h_2$ and $h_3$ in (71) and with the same discussion, there is a contradiction.

\item[(184)]$ R_1:h_2h_3^{-3}h_2h_3=1$, $R_2:h_2h_3^{-3}h_2h_3=1$, $R_3:h_2^2h_3^{-1}h_2^{-1}h_3h_2^{-1}h_3^{-1}h_2=1$:\\
By interchanging $h_2$ and $h_3$ in (70) and with the same discussion, there is a contradiction.

\item[(185)]$ R_1:h_2h_3^{-3}h_2h_3=1$, $R_2:h_2h_3^{-1}h_2(h_3h_2^{-1})^2h_3=1$, $R_3:h_2h_3(h_2h_3^{-1})^2h_2h_3=1$:\\
By interchanging $h_2$ and $h_3$ in (81) and with the same discussion, there is a contradiction.

\item[(186)]$ R_1:h_2h_3^{-3}h_2h_3=1$, $R_2:h_2(h_3^{-1}h_2h_3^{-1})^2h_2=1$, $R_3:h_2h_3^2h_2^{-1}h_3h_2h_3^{-1}h_2=1$:\\
By interchanging $h_2$ and $h_3$ in (80) and with the same discussion, there is a contradiction.

\item[(187)]$ R_1:h_2h_3^{-2}h_2h_3^2=1$, $R_2:h_2h_3^{-2}h_2h_3^2=1$, $R_3:(h_2h_3)^2(h_2^{-1}h_3)^2=1$:\\
By interchanging $h_2$ and $h_3$ in (17) and with the same discussion, there is a contradiction.

\item[(188)]$ R_1:h_2h_3^{-2}h_2h_3^2=1$, $R_2:h_2h_3^{-2}h_2h_3^2=1$, $R_3:h_2h_3h_2^{-1}h_3^{-2}h_2h_3=1$:\\
By interchanging $h_2$ and $h_3$ in (15) and with the same discussion, there is a contradiction.

\item[(189)]$ R_1:h_2h_3^{-2}h_2h_3^2=1$, $R_2:h_2h_3^{-2}h_2h_3^2=1$, $R_3:h_2h_3h_2^{-1}(h_3^{-1}h_2)^2h_3=1$:\\
By interchanging $h_2$ and $h_3$ in (16) and with the same discussion, there is a contradiction.

\item[(190)]$ R_1:h_2h_3^{-2}h_2h_3^2=1$, $R_2:h_2h_3^{-2}h_2h_3^2=1$, $R_3:h_2^3h_3^{-1}h_2^2h_3=1$:\\
By interchanging $h_2$ and $h_3$ in (18) and with the same discussion, there is a contradiction.

\item[(191)]$ R_1:h_2h_3^{-2}h_2h_3^2=1$, $R_2:h_2h_3^{-2}h_2h_3^2=1$, $R_3:h_2^2h_3h_2^2h_3^{-1}h_2=1$:\\
By interchanging $h_2$ and $h_3$ in (19) and with the same discussion, there is a contradiction.

\item[(192)]$ R_1:h_2h_3^{-2}h_2h_3^2=1$, $R_2:(h_2h_3^{-1})^2(h_2^{-1}h_3)^2=1$, $R_3:h_2^3h_3^{-1}h_2^2h_3=1$:\\
By interchanging $h_2$ and $h_3$ in (22) and with the same discussion, there is a contradiction.

\item[(193)]$ R_1:h_2h_3^{-2}h_2h_3^2=1$, $R_2:(h_2h_3^{-1})^2(h_2^{-1}h_3)^2=1$, $R_3:h_2^2h_3h_2^2h_3^{-1}h_2=1$:\\
By interchanging $h_2$ and $h_3$ in (23) and with the same discussion, there is a contradiction.

\item[(194)]$ R_1:h_2h_3^{-2}h_2h_3^2=1$, $R_2:(h_2h_3^{-1})^2h_2h_3h_2^{-1}h_3=1$, $R_3:h_2^3h_3^{-1}h_2^{-2}h_3=1$:\\
By interchanging $h_2$ and $h_3$ in (21) and with the same discussion, there is a contradiction.

\item[(195)]$ R_1:h_2h_3^{-2}h_2h_3^2=1$, $R_2:(h_2h_3^{-1})^2h_2h_3h_2^{-1}h_3=1$, $R_3:h_2^2h_3h_2^{-2}h_3^{-1}h_2=1$:\\
By interchanging $h_2$ and $h_3$ in (20) and with the same discussion, there is a contradiction.

\item[(196)]$ R_1:h_2h_3^{-2}h_2h_3h_2^{-1}h_3=1$, $R_2:(h_2h_3^{-1})^2h_2h_3^2=1$, $R_3:(h_2h_3^{-1}h_2)^2h_3^2=1$:\\
$\Rightarrow h_2=1$ that is a contradiction.

\item[(197)]$ R_1:h_2h_3^{-2}h_2h_3h_2^{-1}h_3=1$, $R_2:(h_2h_3^{-1})^2h_2h_3^2=1$, $R_3:(h_2h_3^{-1}h_2)^2h_3h_2^{-1}h_3=1$:\\
$\Rightarrow h_2=1$ that is a contradiction.

\item[(198)]$ R_1:h_2h_3^{-2}h_2h_3h_2^{-1}h_3=1$, $R_2:(h_2h_3^{-1})^2h_2h_3^2=1$, $R_3:h_2^3h_3^{-1}h_2^2h_3=1$:\\
$\Rightarrow \langle h_2,h_3\rangle=\langle h_2 \rangle$ is abelian, a contradiction.

\item[(199)]$ R_1:h_2h_3^{-2}h_2h_3h_2^{-1}h_3=1$, $R_2:(h_2h_3^{-1})^2h_2h_3^2=1$, $R_3:h_2^2h_3h_2^2h_3^{-1}h_2=1$:\\
$\Rightarrow\langle h_2,h_3\rangle\cong BS(1,1)$ is solvable, a contradiction.

\item[(200)]$ R_1:h_2h_3^{-1}(h_3^{-1}h_2)^2h_3=1$, $R_2:h_2h_3^{-1}(h_3^{-1}h_2)^2h_3=1$, $R_3:h_2^3h_3^{-1}h_2^{-1}h_3^2=1$:\\
By interchanging $h_2$ and $h_3$ in (136) and with the same discussion, there is a contradiction.

\item[(201)]$ R_1:h_2h_3^{-1}(h_3^{-1}h_2)^2h_3=1$, $R_2:h_2h_3^{-1}(h_3^{-1}h_2)^2h_3=1$, $R_3:h_2^3h_3^{-1}h_2^{-1}h_3^{-1}h_2=1$:\\
By interchanging $h_2$ and $h_3$ in (137) and with the same discussion, there is a contradiction.

\item[(202)]$ R_1:h_2h_3^{-1}h_2h_3^3=1$, $R_2:(h_2h_3^{-1})^2(h_2^{-1}h_3)^2=1$, $R_3:h_2^2h_3^{-1}h_2h_3h_2^{-2}h_3=1$:\\
By interchanging $h_2$ and $h_3$ in (7) and with the same discussion, there is a contradiction.

\item[(203)]$ R_1:h_2h_3^{-1}h_2h_3^3=1$, $R_2:(h_2h_3^{-1})^2(h_2^{-1}h_3)^2=1$, $R_3:h_2^2h_3^{-1}h_2h_3h_2^{-1}h_3^2=1$:\\
By interchanging $h_2$ and $h_3$ in (5) and with the same discussion, there is a contradiction.

\item[(204)]$ R_1:h_2h_3^{-1}h_2h_3^3=1$, $R_2:(h_2h_3^{-1})^2(h_2^{-1}h_3)^2=1$, $R_3:h_2h_3^{-1}h_2^{-1}h_3h_2^{-2}h_3^2=1$:\\
By interchanging $h_2$ and $h_3$ in (6) and with the same discussion, there is a contradiction.

\item[(205)]$ R_1:h_2h_3^{-1}h_2h_3^2h_2^{-1}h_3=1$, $R_2:h_2h_3^{-1}h_2h_3^2h_2^{-1}h_3=1$, $R_3:h_2^2h_3^{-1}h_2(h_3h_2^{-1})^2h_3=1$:\\
By interchanging $h_2$ and $h_3$ in (113) and with the same discussion, there is a contradiction.

\item[(206)]$ R_1:(h_2h_3^{-1})^2h_2^{-1}h_3^2=1$, $R_2:(h_2h_3^{-1})^2h_2^{-1}h_3^2=1$, $R_3:h_2^3h_3^{-1}h_2^2h_3=1$:\\
By interchanging $h_2$ and $h_3$ in (144) and with the same discussion, there is a contradiction.

\item[(207)]$ R_1:(h_2h_3^{-1})^2h_2^{-1}h_3^2=1$, $R_2:(h_2h_3^{-1})^2h_2^{-1}h_3^2=1$, $R_3:h_2^2h_3h_2^2h_3^{-1}h_2=1$:\\
By interchanging $h_2$ and $h_3$ in (145) and with the same discussion, there is a contradiction.

\item[(208)]$ R_1:(h_2h_3^{-1})^2h_2^{-1}h_3^2=1$, $R_2:(h_2h_3^{-1})^2h_2^{-1}h_3^2=1$, $R_3:h_2h_3^{-2}h_2^{-1}h_3h_2^{-2}h_3=1$:\\
By interchanging $h_2$ and $h_3$ in (146) and with the same discussion, there is a contradiction.

\item[(209)]$ R_1:(h_2h_3^{-1})^2h_2^{-1}h_3^2=1$, $R_2:(h_2h_3^{-1})^2h_2^{-1}h_3^2=1$, $R_3:h_2h_3^{-1}h_2(h_3h_2^{-1})^2h_3^{-1}h_2=1$:\\
By interchanging $h_2$ and $h_3$ in (147) and with the same discussion, there is a contradiction.

\item[(210)]$ R_1:(h_2h_3^{-1})^2(h_2^{-1}h_3)^2=1$, $R_2:h_3^2(h_3h_2^{-1})^2h_3=1$, $R_3:h_2h_3h_2^{-1}h_3^{-1}h_2h_3^{-1}h_2^{-1}h_3=1$:\\
$R_1\Rightarrow R_3: h_2^2h_3^{-1}h_2^{-2}h_3=1$ and $R_2:h_2h_3^{-2}h_2^{-1}h_3^2=1$. So, $R_3\Rightarrow h_2^2\in Z(G)$ where $G=\langle h_2,h_3\rangle$. Let $x=h_3h_2^{-1}$. So $h_2^{-1}h_3=x^{h_2}$. Also 
$(h_3h_2^{-1})^2=(h_2^{-1}h_3)^2 \text{ so if } H=\langle x,x^{h_2}\rangle\Rightarrow H\cong BS(1,-1)$ is solvable. By Corollary \ref{n-sub2} $H \trianglelefteq G$ since $G=\langle x,h_2\rangle$. Since $\frac{G}{H}=\frac{\langle h_2\rangle H}{H}$ is a cyclic group, it is solvable. So $G$ is solvable, a contradiction.

\item[(211)]$ R_1:(h_2h_3^{-1})^2h_3^{-1}h_2^{-1}h_3=1$, $R_2:(h_2h_3^{-1})^2h_3^{-1}h_2^{-1}h_3=1$, $R_3:h_2^2h_3^{-1}h_2h_3^2h_2^{-1}h_3=1$:\\
By interchanging $h_2$ and $h_3$ in (114) and with the same discussion, there is a contradiction.

\item[(212)]$ R_1:(h_2h_3^{-1})^2h_3^{-1}h_2^{-1}h_3=1$, $R_2:(h_2h_3^{-1})^2h_3^{-1}h_2^{-1}h_3=1$, $R_3:h_2h_3h_2^{-1}h_3^{-1}h_2h_3^{-1}h_2^{-1}h_3=1$:\\
By interchanging $h_2$ and $h_3$ in (115) and with the same discussion, there is a contradiction.

\item[(213)]$ R_1:h_2(h_3^{-1}h_2h_3^{-1})^2h_2=1$, $R_2:(h_2h_3^{-1})^3h_2h_3=1$, $R_3:h_2h_3^{-1}h_2^{-1}h_3h_2^{-2}h_3^2=1$:\\
$\Rightarrow h_2=1$ that is a contradiction.

\item[(214)]$ R_1:h_2(h_3^{-1}h_2h_3^{-1})^2h_2=1$, $R_2:(h_2h_3^{-1})^3h_2h_3=1$, $R_3:h_2h_3^{-3}h_2^2h_3^{-1}h_2=1$:\\
$\Rightarrow h_2=1$ that is a contradiction.

\end{enumerate}

\subsection{$\mathbf{C_5-C_5(-C_6--)(C_6---)}$}
By considering the relations from Tables \ref{tab-C6} and \ref{tab-C5-C5(-C6--)} which are not disproved, it can be seen that there are $56$ cases for the relations of two cycles $C_5$ and two cycles $C_6$ in the graph $C_5-C_5(-C_6--)(C_6---)$. Using Gap \cite{gap}, we see that all groups with two generators $h_2$ and $h_3$ and four relations which are between $54$ cases of these $56$ cases are finite and solvable, that is a contradiction with the assumptions. So there are just $2$ cases for the relations of these cycles which may lead to the existence of a subgraph isomorphic to the graph $C_5-C_5(-C_6--)(C_6---)$ in the graph $K(\alpha,\beta)$. These cases are listed in table \ref{tab-C5-C5(-C6--)(C6---)}.  In the following, we show that these $2$ cases lead to contradictions and so, the graph $K(\alpha,\beta)$ contains no subgraph isomorphic to the graph $C_5-C_5(-C_6--)(C_6---)$.

\begin{enumerate}
\item[(1)]$ R_1:h_2^2h_3^{-2}h_2h_3=1$, $R_2:h_2^2h_3^{-2}h_2h_3=1$, $R_3:h_2^3(h_2h_3^{-1})^2h_2=1$, $R_4:(h_2^2h_3^{-1})^2h_2^{-1}h_3=1$:\\
$\Rightarrow \langle h_2,h_3\rangle=\langle h_2\rangle$ is abelian, a contradiction.

\item[(2)]$ R_1:h_2h_3^2h_2^{-2}h_3=1$, $R_2:h_2h_3^2h_2^{-2}h_3=1$, $R_3:h_3^3(h_3h_2^{-1})^2h_3=1$, $R_4:(h_2h_3^{-2})^2h_2^{-1}h_3=1$:\\
By interchanging $h_2$ and $h_3$ in (1) and with the same discussion, there is a contradiction.
\end{enumerate}

\begin{table}[h]
\centering
\caption{The relations of a $C_4-C_6(-C_6--)(C_6---)$ in $K(\alpha,\beta)$}\label{tab-C5-C5(-C6--)(C6---)}
\begin{tabular}{|c|l|l|l|l|}\hline
$n$&$R_1$&$R_2$&$R_3$&$R_4$\\\hline
$1$&$h_2^2h_3^{-2}h_2h_3=1$&$h_2^2h_3^{-2}h_2h_3=1$&$h_2^3(h_2h_3^{-1})^2h_2=1$&$(h_2^2h_3^{-1})^2h_2^{-1}h_3=1$\\
$2$&$h_2h_3^2h_2^{-2}h_3=1$&$h_2h_3^2h_2^{-2}h_3=1$&$h_3^3(h_3h_2^{-1})^2h_3=1$&$(h_2h_3^{-2})^2h_2^{-1}h_3=1$\\
\hline
\end{tabular}
\end{table}

\subsection{$\mathbf{C_5-C_5(-C_6--)(--C_6-1)}$}
By considering the relations from Tables \ref{tab-C6} and \ref{tab-C5-C5(-C6--)} which are not disproved, it can be seen that there are $56$ cases for the relations of two cycles $C_5$ and two cycles $C_6$ in the graph $C_5-C_5(-C_6--)(--C_6-1)$. By considering all groups with two generators $h_2$ and $h_3$ and four relations which are between these cases and by using GAP \cite{gap}, we see that all of these groups are finite and solvable. So, the graph $K(\alpha,\beta)$ contains no subgraph isomorphic to the graph $C_5-C_5(-C_6--)(--C_6-1)$.

\subsection{$\mathbf{C_5-C_5(-C_6--)(-C_5--)}$}
By considering the relations from Tables \ref{tab-C5} and \ref{tab-C5-C5(-C6--)} which are not disproved, it can be seen that there are $14$ cases for the relations of three cycles $C_5$ and one cycle $C_6$ in the graph $C_5-C_5(-C_6--)(-C_5--)$. By considering all groups with two generators $h_2$ and $h_3$ and four relations which are between these cases and by using GAP \cite{gap}, we see that all of these groups are finite and solvable. So, the graph $K(\alpha,\beta)$ contains no subgraph isomorphic to the graph $C_5-C_5(-C_6--)(-C_5--)$.

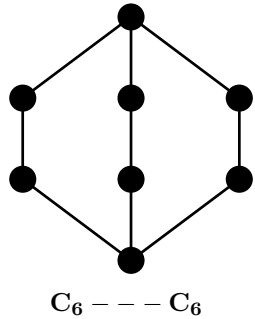
\begin{figure}[ht]
\psscalebox{0.9 0.9} 
{
\begin{pspicture}(0,-2.2985578)(3.594231,2.2985578)
\psdots[linecolor=black, dotsize=0.4](1.7971154,2.1014423)
\psdots[linecolor=black, dotsize=0.4](1.7971154,0.90144235)
\psdots[linecolor=black, dotsize=0.4](1.7971154,-0.29855767)
\psdots[linecolor=black, dotsize=0.4](1.7971154,-1.4985577)
\psdots[linecolor=black, dotsize=0.4](3.3971155,-0.29855767)
\psdots[linecolor=black, dotsize=0.4](3.3971155,0.90144235)
\psdots[linecolor=black, dotsize=0.4](0.19711548,0.90144235)
\psdots[linecolor=black, dotsize=0.4](0.19711548,-0.29855767)
\psline[linecolor=black, linewidth=0.04](1.7971154,2.1014423)(0.19711548,0.90144235)(0.19711548,-0.29855767)(1.7971154,-1.4985577)(1.7971154,-0.29855767)(1.7971154,0.90144235)(1.7971154,2.1014423)(3.3971155,0.90144235)(3.3971155,-0.29855767)(1.7971154,-1.4985577)
\rput[bl](0.59711546,-2.2985578){$\mathbf{C_6---C_6}$}
\end{pspicture}
}
\end{figure}
$\mathbf{C_6---C_6}$ \textbf{subgraph:} Suppose that there are two cycles of length $6$ with three successive common edges in the graph $K(\alpha,\beta)$ and denote such subgraph by $C_6---C_6$. \\
Let $[a_1,b_1,a_2,b_2,a_3,b_3,a_4,b_4,a_5,b_5,a_6,b_6]$ and $[a_1,b_1,a_2,b_2,a_3,b_3,a_4',b_4',a_5',b_5',a_6',b_6']$ be $12-$tuples related to the cycles $C_6$ in the graph $C_6---C_6$, where the first six components of these tuples are related to the three successive common edges of $C_6$ and $C_6$. Without loss of generality we may assume that $a_1=1$, where $a_1,b_1,a_2,b_2,a_3,b_3,a_4,b_4,a_5,b_5,a_6,b_6,a_4',b_4',a_5',b_5',a_6',b_6' \in supp(\alpha)$ and $\alpha=1+h_2+h_3$. With the same discussion such as about $K_{2,3}$, it is easy to see that $a_4 \neq a_4'$ and $b_6 \neq b_6'$. By considering the relations from Table \ref{tab-C6} which are not disproved and above assumptions, it can be seen that there are $4454$ cases for existing  two cycles of length $6$ with three successive common edges in the graph $K(\alpha,\beta)$. Using Gap \cite{gap}, we see that all groups with two generators $h_2$ and $h_3$ and two relations which are between $4022$ cases of these $4454$ cases are solvable or finite. So there are $432$ cases for the relations of the existence of $C_6---C_6$ in the graph $K(\alpha,\beta)$. Similar to the previous mentioned subgraphs it can be seen that $333$ cases of these relations lead to contradictions and  $99$ cases of them may lead to the  existence of a subgraph isomorphic to the graph $C_6---C_6$ in the graph $K(\alpha,\beta)$.

\subsection{$\mathbf{C_6---C_6(C_6---C_6)}$}
\begin{figure}[ht]
\psscalebox{0.9 0.9} 
{
\begin{pspicture}(0,-2.0985577)(4.394231,2.0985577)
\psdots[linecolor=black, dotsize=0.4](2.1971154,1.9014423)
\psdots[linecolor=black, dotsize=0.4](2.1971154,1.1014423)
\psdots[linecolor=black, dotsize=0.4](2.1971154,0.30144233)
\psdots[linecolor=black, dotsize=0.4](2.1971154,-0.49855766)
\psdots[linecolor=black, dotsize=0.4](3.3971155,1.1014423)
\psdots[linecolor=black, dotsize=0.4](3.3971155,0.30144233)
\psdots[linecolor=black, dotsize=0.4](0.9971155,1.1014423)
\psdots[linecolor=black, dotsize=0.4](0.9971155,0.30144233)
\psline[linecolor=black, linewidth=0.04](2.1971154,1.9014423)(0.9971155,1.1014423)(0.9971155,0.30144233)(2.1971154,-0.49855766)(3.3971155,0.30144233)(3.3971155,1.1014423)(2.1971154,1.9014423)(2.1971154,1.1014423)(2.1971154,0.30144233)(2.1971154,-0.49855766)
\psdots[linecolor=black, dotsize=0.4](1.3971155,-1.2985576)
\psdots[linecolor=black, dotsize=0.4](0.19711548,-0.49855766)
\psdots[linecolor=black, dotsize=0.4](4.1971154,-0.49855766)
\psline[linecolor=black, linewidth=0.04](0.9971155,1.1014423)(0.19711548,-0.49855766)(1.3971155,-1.2985576)(2.1971154,0.30144233)(2.1971154,0.30144233)
\psline[linecolor=black, linewidth=0.04](1.3971155,-1.2985576)(4.1971154,-0.49855766)(3.3971155,1.1014423)
\rput[bl](-0.6,-2.0985577){$\mathbf{21) \ C_6---C_6(C_6---C_6)}$}
\end{pspicture}
}
\end{figure}
By considering the $99$ cases related to the existence of $C_6---C_6$ in the graph $K(\alpha,\beta)$, it can be seen that there are $42$ cases for the relations of four $C_6$ cycles in this structure. Using GAP \cite{gap}, we see that all groups with two generators $h_2$ and $h_3$ and four relations which are between $26$ cases of these $42$ cases are finite and solvable, or just finite, that is a contradiction with the assumptions. So, there are just $16$ cases for the relations of these cycles which may lead to the existence of a subgraph isomorphic to the graph $C_6---C_6(C_6---C_6)$ in $K(\alpha,\beta)$. In the following, we show that these $16$ cases lead to contradictions and so, the graph $K(\alpha,\beta)$ contains no subgraph isomorphic to the graph $C_6---C_6(C_6---C_6)$.
\begin{enumerate}
\item[(1)]$ R_1:h_2^2h_3^{-1}h_2^{-1}(h_3^{-1}h_2)^2=1$, $R_2:h_2^2h_3^{-2}h_2h_3^{-1}h_2^{-1}h_3=1$, $R_3:h_2h_3h_2^{-1}h_3^{-1}h_2h_3^{-2}h_2=1$, $R_4:h_2h_3^{-1}h_2^{-1}h_3^2h_2^{-1}h_3^{-1}h_2=1$:\\
$\Rightarrow\langle h_2,h_3\rangle\cong BS(1,1)$ is solvable, a contradiction.

\item[(2)]$ R_1:h_2^2h_3^{-1}h_2^{-1}(h_3^{-1}h_2)^2=1$, $R_2:h_2h_3h_2^{-1}h_3^{-1}h_2h_3^{-2}h_2=1$, $R_3:h_2^2h_3^{-2}h_2h_3^{-1}h_2^{-1}h_3=1$, $R_4:h_2h_3^{-1}h_2^{-1}h_3h_2h_3^{-2}h_2=1$:\\
$\Rightarrow\langle h_2,h_3\rangle\cong BS(1,1)$ is solvable, a contradiction.

\item[(3)]$ R_1:h_2^2h_3^{-1}h_2^{-1}(h_3^{-1}h_2)^2=1$, $R_2:h_2h_3h_2^{-1}h_3^{-1}h_2h_3^{-2}h_2=1$, $R_3:h_2^2h_3^{-2}h_2h_3^{-1}h_2^{-1}h_3=1$, $R_4:h_2h_3^{-1}h_2^{-1}h_3^2h_2^{-1}h_3^{-1}h_2=1$:\\
$\Rightarrow\langle h_2,h_3\rangle\cong BS(1,1)$ is solvable, a contradiction.

\item[(4)]$ R_1:(h_2h_3)^2h_2h_3^{-1}h_2=1$, $R_2:(h_2h_3)^2h_2h_3^{-1}h_2=1$, $R_3:h_2(h_3h_2^{-1}h_3)^2h_3=1$,\\ $R_4:h_2(h_3h_2^{-1})^2h_3^{-2}h_2=1$:\\
$\Rightarrow$  $G$ has a torsion element, a contradiction.

\item[(5)]$ R_1:(h_2h_3)^2h_2h_3^{-1}h_2=1$, $R_2:(h_2h_3^{-1}h_2)^2h_3^{-2}h_2=1$, $R_3:h_2(h_3h_2^{-1})^2h_3^{-2}h_2=1$,\\ $R_4:(h_2h_3^{-1}h_2)^2h_3^{-2}h_2=1$:\\
$\Rightarrow$  $G$ has a torsion element, a contradiction.

\item[(6)]$ R_1:h_2h_3(h_2h_3^{-1}h_2)^2=1$, $R_2:h_2h_3(h_2h_3^{-1}h_2)^2=1$, $R_3:h_2(h_3h_2^{-1})^2h_3^{-2}h_2=1$,\\ $R_4:h_2(h_3^{-1}h_2h_3^{-1})^2h_3^{-1}h_2=1$:\\
$\Rightarrow$  $G$ has a torsion element, a contradiction.

\item[(7)]$ R_1:h_2h_3(h_2h_3^{-1}h_2)^2=1$, $R_2:h_2(h_3h_2^{-1})^2h_3^{-2}h_2=1$, $R_3:(h_2h_3)^2h_2^{-1}h_3^2=1$,\\ $R_4:(h_2h_3)^2h_2^{-1}h_3^2=1$:\\
By interchanging $h_2$ and $h_3$ in (4) and with the same discussion, there is a contradiction.

\item[(8)]$ R_1:h_2(h_3h_2^{-2})^2h_3=1$, $R_2:h_2h_3h_2^{-1}h_3^{-1}(h_3^{-1}h_2)^2=1$, $R_3:h_2h_3^{-1}h_2^{-1}h_3^2h_2^{-2}h_3=1$,\\ $R_4:(h_2h_3^{-2})^2h_2h_3=1$:\\
$\Rightarrow$  $G$ has a torsion element, a contradiction.

\item[(9)]$ R_1:h_2(h_3h_2^{-2})^2h_3=1$, $R_2:h_2h_3^{-1}h_2^{-1}h_3^2h_2^{-2}h_3=1$, $R_3:h_2h_3h_2^{-1}h_3^{-1}(h_3^{-1}h_2)^2=1$,\\\ $R_4:(h_2h_3^{-2})^2h_2h_3=1$:\\
By interchanging $h_2$ and $h_3$ in (8) and with the same discussion, there is a contradiction.

\item[(10)]$ R_1:h_2h_3h_2^{-2}h_3h_2^{-1}h_3^{-1}h_2=1$, $R_2:h_2h_3h_2^{-1}h_3^{-1}(h_2^{-1}h_3)^2=1$, $R_3:h_2h_3^{-1}h_2^{-1}h_3^2h_2^{-2}h_3=1$, $R_4:h_2h_3^{-2}h_2h_3^{-1}h_2^{-1}h_3^2=1$:\\
$\Rightarrow$  $G$ has a torsion element, a contradiction.

\item[(11)]$ R_1:h_2h_3h_2^{-2}h_3h_2^{-1}h_3^{-1}h_2=1$, $R_2:h_2h_3^{-1}h_2^{-1}h_3^2h_2^{-2}h_3=1$, $R_3:h_2h_3h_2^{-1}h_3^{-1}(h_2^{-1}h_3)^2=1$, $R_4:h_2h_3^{-2}h_2h_3^{-1}h_2^{-1}h_3^2=1$:\\
By interchanging $h_2$ and $h_3$ in (10) and with the same discussion, there is a contradiction.

\item[(12)]$ R_1:h_2(h_3h_2^{-1})^2h_3^{-2}h_2=1$, $R_2:(h_2h_3^{-1}h_2)^2h_3^{-2}h_2=1$, $R_3:h_2(h_3h_2^{-1}h_3)^2h_3=1$,\\ $R_4:h_2(h_3h_2^{-1}h_3)^2h_3=1$:\\
By interchanging $h_2$ and $h_3$ in (6) and with the same discussion, there is a contradiction.

\item[(13)]$ R_1:h_2h_3^{-2}h_2^{-1}h_3h_2h_3^{-1}h_2=1$, $R_2:h_2h_3^{-2}(h_3^{-1}h_2)^2h_3=1$, $R_3:h_2h_3^{-1}h_2^{-1}h_3^2h_2^{-1}h_3^{-1}h_2=1$, $R_4:h_2h_3^{-1}h_2^{-1}h_3h_2^{-2}h_3^2=1$:\\
By interchanging $h_2$ and $h_3$ in (1) and with the same discussion, there is a contradiction.

\item[(14)]$ R_1:h_2h_3^{-1}h_2^{-1}h_3h_2^{-2}h_3^2=1$, $R_2:h_2h_3^{-2}(h_3^{-1}h_2)^2h_3=1$, $R_3:h_2h_3^{-1}h_2^{-1}h_3h_2h_3^{-2}h_2=1$, \\$R_4:h_2h_3^{-2}h_2^{-1}h_3h_2h_3^{-1}h_2=1$:\\
By interchanging $h_2$ and $h_3$ in (2) and with the same discussion, there is a contradiction.

\item[(15)]$ R_1:h_2h_3^{-1}h_2^{-1}h_3h_2^{-2}h_3^2=1$, $R_2:h_2h_3^{-2}(h_3^{-1}h_2)^2h_3=1$, $R_3:h_2h_3^{-1}h_2^{-1}h_3^2h_2^{-1}h_3^{-1}h_2=1$, $R_4:h_2h_3^{-2}h_2^{-1}h_3h_2h_3^{-1}h_2=1$:\\
By interchanging $h_2$ and $h_3$ in (1) and with the same discussion, there is a contradiction.

\item[(16)]$ R_1:(h_2h_3)^2h_2^{-1}h_3^2=1$, $R_2:h_2(h_3^{-1}h_2h_3^{-1})^2h_3^{-1}h_2=1$, $R_3:h_2(h_3h_2^{-1})^2h_3^{-2}h_2=1$, \\$R_4:h_2(h_3^{-1}h_2h_3^{-1})^2h_3^{-1}h_2=1$:\\
By interchanging $h_2$ and $h_3$ in (5) and with the same discussion, there is a contradiction.

\end{enumerate}

$\mathbf{C_6---C_6(C_6)}$ \textbf{subgraph:} 
By considering the $99$ cases related to the existence of $C_6---C_6$ in the graph $K(\alpha,\beta)$ and the relations from Table \ref{tab-C6} which are not disproved, it can be seen that there are $2618$ cases for the relations of three $C_6$ cycles in this structure. Using GAP \cite{gap}, we see that all groups with two generators $h_2$ and $h_3$ and three relations which are between $2352$ cases of these $2618$ cases are finite and solvable, or just finite, that is a contradiction with the assumptions. So, there are just $266$ cases for the relations of these cycles which may lead to the existence of a subgraph isomorphic to the graph $C_6---C_6(C_6)$ in $K(\alpha,\beta)$. Similar to the previous mentioned subgraphs it can be seen that $228$ cases of these relations lead to contradictions and  $38$ cases of them may lead to the  existence of a subgraph isomorphic to the graph $C_6---C_6(C_6)$ in the graph $K(\alpha,\beta)$.

$\mathbf{C_6---C_6(C_6)(C_6)}$ \textbf{subgraph:} 
By considering the $38$ cases related to the existence of $C_6---C_6(C_6)$ in the graph $K(\alpha,\beta)$ and the relations from Table \ref{tab-C6} which are not disproved, it can be seen that there are $374$ cases for the relations of four $C_6$ cycles in this structure. Using GAP \cite{gap}, we see that all groups with two generators $h_2$ and $h_3$ and four relations which are between $340$ cases of these $374$ cases are finite and solvable, or just finite, that is a contradiction with the assumptions. So, there are just $34$ cases for the relations of these cycles which may lead to the existence of a subgraph isomorphic to the graph $C_6---C_6(C_6)(C_6)$ in $K(\alpha,\beta)$. Similar to the previous mentioned subgraphs it can be seen that $30$ cases of these relations lead to contradictions and  $4$ cases of them may lead to the  existence of a subgraph isomorphic to the graph $C_6---C_6(C_6)(C_6)$ in the graph $K(\alpha,\beta)$.

\subsection{$\mathbf{C_6---C_6(C_6)(C_6)(C_6)}$}
\begin{figure}[ht]
\psscalebox{0.9 0.9} 
{
\begin{pspicture}(0,-2.0)(5.81,2.2985578)
\psdots[linecolor=black, dotsize=0.4](2.8,2.1014423)
\psdots[linecolor=black, dotsize=0.4](2.8,1.3014424)
\psdots[linecolor=black, dotsize=0.4](2.8,0.5014423)
\psdots[linecolor=black, dotsize=0.4](2.8,-0.29855767)
\psdots[linecolor=black, dotsize=0.4](4.0,1.3014424)
\psdots[linecolor=black, dotsize=0.4](4.0,0.5014423)
\psdots[linecolor=black, dotsize=0.4](1.6,1.3014424)
\psdots[linecolor=black, dotsize=0.4](1.6,0.5014423)
\psline[linecolor=black, linewidth=0.04](2.8,2.1014423)(1.6,1.3014424)(1.6,0.5014423)(2.8,-0.29855767)(4.0,0.5014423)(4.0,1.3014424)(2.8,2.1014423)(2.8,1.3014424)(2.8,0.5014423)(2.8,-0.29855767)
\psdots[linecolor=black, dotsize=0.4](3.6,-0.6985577)
\psdots[linecolor=black, dotsize=0.4](4.4,-0.29855767)
\psline[linecolor=black, linewidth=0.04](4.0,0.5014423)(4.4,-0.29855767)(3.6,-0.6985577)(2.8,1.3014424)
\psdots[linecolor=black, dotsize=0.4](0.8,0.10144234)
\psdots[linecolor=black, dotsize=0.4](1.6,-0.6985577)
\psline[linecolor=black, linewidth=0.04](1.6,1.3014424)(0.8,0.10144234)(1.6,-0.6985577)(2.8,0.5014423)
\psdots[linecolor=black, dotsize=0.4](1.6,-1.4985577)
\psline[linecolor=black, linewidth=0.04](0.8,0.10144234)(1.6,-1.4985577)(3.6,-0.6985577)
\rput[bl](0.0,-2.){$\mathbf{22) \ C_6---C_6(C_6)(C_6)(C_6)}$}
\end{pspicture}
}
\end{figure}
By considering the $4$ cases related to the existence of $C_6---C_6(C_6)(C_6)$ in the graph $K(\alpha,\beta)$ and the relations from Table \ref{tab-C6} which are not disproved, it can be seen that there are $10$ cases for the relations of five $C_6$ cycles in this structure. Using GAP \cite{gap}, we see that all groups with two generators $h_2$ and $h_3$ and five relations which are between $6$ cases of these $10$ cases are finite and solvable, or just finite, that is a contradiction with the assumptions. So, there are just $4$ cases for the relations of these cycles which may lead to the existence of a subgraph isomorphic to the graph $C_6---C_6(C_6)(C_6)(C_6)$ in $K(\alpha,\beta)$. In the following, we show that these $4$ cases lead to contradictions and so, the graph $K(\alpha,\beta)$ contains no subgraph isomorphic to the graph $C_6---C_6(C_6)(C_6)(C_6)$.
\begin{enumerate}
\item[(1)]$ R_1:(h_2h_3)^2h_2h_3^{-1}h_2=1$, $R_2:(h_2h_3)^2h_2h_3^{-1}h_2=1$, $R_3:h_2(h_3h_2^{-1})^2h_3^{-2}h_2=1$,\\ $R_4:(h_2h_3)^2h_2h_3^{-1}h_2=1$, $R_5:h_2(h_3^{-1}h_2h_3^{-1})^2h_3^{-1}h_2=1$:\\
$\Rightarrow \langle h_2,h_3\rangle=\langle h_3\rangle$ is abelian, a contradiction.

\item[(2)]$ R_1:(h_2h_3)^2h_2^{-1}h_3^2=1$, $R_2:(h_2h_3)^2h_2^{-1}h_3^2=1$, $R_3:h_2(h_3h_2^{-1})^2h_3^{-2}h_2=1$,\\ $R_4:(h_2h_3)^2h_2^{-1}h_3^2=1$, $R_5:(h_2h_3^{-1}h_2)^2h_3^{-2}h_2=1$:\\
By interchanging $h_2$ and $h_3$ in (1) and with the same discussion, there is a contradiction.

\item[(3)]$ R_1:h_2h_3(h_2h_3^{-1}h_2)^2=1$, $R_2:h_2h_3(h_2h_3^{-1}h_2)^2=1$, $R_3:h_2(h_3h_2^{-1})^2h_3^{-2}h_2=1$,\\ $R_4:h_2h_3(h_2h_3^{-1}h_2)^2=1$, $R_5:h_2(h_3h_2^{-1}h_3)^2h_3=1$:\\
$\Rightarrow\langle h_2,h_3\rangle\cong BS(1,-2)$ is solvable, a contradiction.

\item[(4)]$ R_1:h_2(h_3h_2^{-1}h_3)^2h_3=1$, $R_2:h_2(h_3h_2^{-1}h_3)^2h_3=1$, $R_3:h_2(h_3h_2^{-1})^2h_3^{-2}h_2=1$,\\ $R_4:h_2(h_3h_2^{-1}h_3)^2h_3=1$, $R_5:h_2h_3(h_2h_3^{-1}h_2)^2=1$:\\
By interchanging $h_2$ and $h_3$ in (3) and with the same discussion, there is a contradiction.

\end{enumerate}

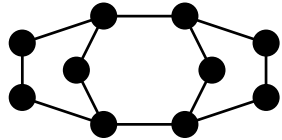
\begin{figure}[ht]
\psscalebox{0.9 0.9} 
{
\begin{pspicture}(0,-1.4985577)(3.994231,1.4985577)
\psdots[linecolor=black, dotsize=0.4](1.3971155,1.3014424)
\psdots[linecolor=black, dotsize=0.4](2.5971155,1.3014424)
\psdots[linecolor=black, dotsize=0.4](2.9971154,0.5014423)
\psdots[linecolor=black, dotsize=0.4](2.5971155,-0.29855767)
\psdots[linecolor=black, dotsize=0.4](1.3971155,-0.29855767)
\psdots[linecolor=black, dotsize=0.4](0.9971155,0.5014423)
\psdots[linecolor=black, dotsize=0.4](0.19711548,0.90144235)
\psdots[linecolor=black, dotsize=0.4](0.19711548,0.10144234)
\psdots[linecolor=black, dotsize=0.4](3.7971156,0.10144234)
\psdots[linecolor=black, dotsize=0.4](3.7971156,0.90144235)
\psline[linecolor=black, linewidth=0.04](1.3971155,1.3014424)(0.9971155,0.5014423)(1.3971155,-0.29855767)(2.5971155,-0.29855767)(2.9971154,0.5014423)(2.5971155,1.3014424)(1.3971155,1.3014424)(0.19711548,0.90144235)(0.19711548,0.10144234)(1.3971155,-0.29855767)
\psline[linecolor=black, linewidth=0.04](2.5971155,-0.29855767)(3.7971156,0.10144234)(3.7971156,0.90144235)(2.5971155,1.3014424)
\rput[bl](0.59711546,-1.4985577){$\mathbf{C_5(--C_6--)C_5}$}
\end{pspicture}
}
\end{figure}
$\mathbf{C_5(--C_6--)C_5}$ \textbf{subgraph:} 
It can be seen that there are $2215$ cases for the relations of two $C_5$ and one $C_6$ cycles in this structure. Using GAP \cite{gap}, we see that all groups with two generators $h_2$ and $h_3$ and three relations which are between $2046$ cases of these $2215$ cases are finite and solvable, or just finite, that is a contradiction with the assumptions. So, there are just $169$ cases for the relations of these cycles which may lead to the existence of a subgraph isomorphic to the graph $C_5(--C_6--)C_5$ in $K(\alpha,\beta)$. Similar to the previous mentioned subgraphs it can be seen that $163$ cases of these relations lead to contradictions and  $6$ cases of them may lead to the  existence of a subgraph isomorphic to the graph $C_5(--C_6--)C_5$ in the graph $K(\alpha,\beta)$.

\subsection{$\mathbf{C_5(--C_6--)C_5(---C_6)}$}
\begin{figure}[ht]
\psscalebox{0.9 0.9}  
{
\begin{pspicture}(0,-2.0985577)(4.56,2.0985577)
\psdots[linecolor=black, dotsize=0.4](1.6,1.9014423)
\psdots[linecolor=black, dotsize=0.4](2.8,1.9014423)
\psdots[linecolor=black, dotsize=0.4](3.2,1.1014423)
\psdots[linecolor=black, dotsize=0.4](2.8,0.3014423)
\psdots[linecolor=black, dotsize=0.4](1.6,0.3014423)
\psdots[linecolor=black, dotsize=0.4](1.2,1.1014423)
\psdots[linecolor=black, dotsize=0.4](0.4,1.5014423)
\psdots[linecolor=black, dotsize=0.4](0.4,0.7014423)
\psdots[linecolor=black, dotsize=0.4](4.0,0.7014423)
\psdots[linecolor=black, dotsize=0.4](4.0,1.5014423)
\psdots[linecolor=black, dotsize=0.4](0.8,-0.4985577)
\psdots[linecolor=black, dotsize=0.4](3.6,-0.4985577)
\psline[linecolor=black, linewidth=0.04](1.6,1.9014423)(0.4,1.5014423)(0.4,0.7014423)(1.6,0.3014423)(1.2,1.1014423)(1.6,1.9014423)(2.8,1.9014423)(3.2,1.1014423)(2.8,0.3014423)(1.6,0.3014423)
\psline[linecolor=black, linewidth=0.04](2.8,1.9014423)(4.0,1.5014423)(4.0,0.7014423)(2.8,0.3014423)
\psline[linecolor=black, linewidth=0.04](0.4,0.7014423)(0.8,-0.4985577)(3.6,-0.4985577)(3.2,1.1014423)
\rput[bl](-0.7,-2.0985577){$\mathbf{23) \ C_5(--C_6--)C_5(---C_6)}$}
\end{pspicture}
}
\end{figure}
By considering the $6$ cases related to the existence of $C_5(--C_6--)C_5$ in the graph $K(\alpha,\beta)$ and the relations from Table \ref{tab-C6} which are not disproved, it can be seen that there are $134$ cases for the relations of two cycles $C_5$ and two cycles $C_6$ in this structure. Using GAP \cite{gap}, we see that all groups with two generators $h_2$ and $h_3$ and four relations which are between $130$ cases of these $134$ cases are finite and solvable, or just finite, that is a contradiction with the assumptions. So, there are just $4$ cases for the relations of these cycles which may lead to the existence of a subgraph isomorphic to the graph $C_5(--C_6--)C_5(---C_6)$ in $K(\alpha,\beta)$. In the following, we show that these $4$ cases lead to contradictions and so, the graph $K(\alpha,\beta)$ contains no subgraph isomorphic to the graph $C_5(--C_6--)C_5(---C_6)$.
\begin{enumerate}
\item[(1)]$ R_1:h_2h_3h_2^{-1}(h_3^{-1}h_2)^2=1$, $R_2:h_2(h_2h_3^{-1})^2h_2^{-1}h_3=1$, $R_3:h_2^2h_3^{-1}(h_3^{-1}h_2)^2h_3=1$,\\ $R_4:(h_2h_3^{-1})^3h_2^{-1}h_3^2=1$:\\
$\Rightarrow\langle h_2,h_3\rangle\cong BS(1,1)$ is solvable, a contradiction.
\item[(2)]$ R_1:h_2h_3h_2^{-1}(h_3^{-1}h_2)^2=1$, $R_2:h_2h_3h_2^{-1}(h_3^{-1}h_2)^2=1$, $R_3:h_2^3h_3^{-2}h_2^{-1}h_3=1$,\\ $R_4:h_2h_3^{-1}h_2^{-1}h_3h_2^{-1}h_3^{-1}h_2h_3=1$:\\
$\Rightarrow \langle h_2,h_3\rangle=\langle h_2\rangle$ is abelian, a contradiction.
\item[(3)]$ R_1:h_2h_3^{-1}(h_2^{-1}h_3)^2h_3=1$, $R_2:h_2h_3^{-1}(h_2^{-1}h_3)^2h_3=1$, $R_3:h_2^2h_3^{-3}h_2^{-1}h_3=1$,\\ $R_4:h_2h_3h_2^{-1}h_3^{-1}h_2h_3^{-1}h_2^{-1}h_3=1$:\\
By interchanging $h_2$ and $h_3$ in (2) and with the same discussion, there is a contradiction.
\item[(4)]$ R_1:h_2h_3^{-1}(h_2^{-1}h_3)^2h_3=1$, $R_2:(h_2h_3^{-1})^2h_3^{-1}h_2^{-1}h_3=1$, $R_3:h_2h_3^2h_2^{-1}(h_2^{-1}h_3)^2=1$,\\ $R_4:h_2(h_3h_2^{-1})^3h_3^{-1}h_2=1$:\\
By interchanging $h_2$ and $h_3$ in (1) and with the same discussion, there is a contradiction.
\end{enumerate}

\begin{figure}[ht]
\psscalebox{0.9 0.9} 
{
\begin{pspicture}(0,-2.2985578)(3.594231,2.2985578)
\psdots[linecolor=black, dotsize=0.4](1.7971154,2.1014423)
\psdots[linecolor=black, dotsize=0.4](1.7971154,0.5014423)
\psdots[linecolor=black, dotsize=0.4](1.7971154,-1.0985577)
\psdots[linecolor=black, dotsize=0.4](3.3971155,-0.6985577)
\psdots[linecolor=black, dotsize=0.4](3.3971155,0.5014423)
\psdots[linecolor=black, dotsize=0.4](3.3971155,1.7014424)
\psdots[linecolor=black, dotsize=0.4](0.19711548,1.7014424)
\psdots[linecolor=black, dotsize=0.4](0.19711548,0.5014423)
\psdots[linecolor=black, dotsize=0.4](0.19711548,-0.6985577)
\psline[linecolor=black, linewidth=0.04](1.7971154,2.1014423)(0.19711548,1.7014424)(0.19711548,0.5014423)(0.19711548,-0.6985577)(1.7971154,-1.0985577)(3.3971155,-0.6985577)(3.3971155,0.5014423)(3.3971155,1.7014424)(1.7971154,2.1014423)(1.7971154,0.5014423)(1.7971154,-1.0985577)
\rput[bl](0.9971155,-2.2985578){$\mathbf{C_6--C_6}$}
\end{pspicture}
}
\end{figure}
$\mathbf{C_6--C_6}$ \textbf{subgraph:} Suppose that there are two cycles of length $6$ with two successive common edges in the graph $K(\alpha,\beta)$ and denote such subgraph by $C_6--C_6$. \\
Let $[a_1,b_1,a_2,b_2,a_3,b_3,a_4,b_4,a_5,b_5,a_6,b_6]$ and $[a_1,b_1,a_2,b_2,a_3',b_3',a_4',b_4',a_5',b_5',a_6',b_6']$ be $12-$tuples related to the cycles $C_6$ in the graph $C_6--C_6$, where the first four components of these tuples are related to the two successive common edges of $C_6$ and $C_6$. Without loss of generality we may assume that $a_1=1$, where $a_1,b_1,a_2,b_2,a_3,b_3,a_4,b_4,a_5,b_5,a_6,b_6,a_3',b_3',a_4',b_4',a_5',b_5',a_6',b_6' \in supp(\alpha)$ and $\alpha=1+h_2+h_3$. With the same discussion such as about $K_{2,3}$, it is easy to see that $a_3 \neq a_3'$ and $b_6 \neq b_6'$. By considering the relations from Table \ref{tab-C6} which are not disproved and above assumptions, it can be seen that there are $16462$ cases for existing  two cycles of length $6$ with two successive common edges in the graph $K(\alpha,\beta)$. Using Gap \cite{gap}, we see that all groups with two generators $h_2$ and $h_3$ and two relations which are between $14846$ cases of these $16462$ cases are solvable or finite. So there are $1616$ cases for the relations of the existence of $C_6--C_6$ in the graph $K(\alpha,\beta)$. Similar to the previous mentioned subgraphs it can be seen that $620$ cases of these relations lead to contradictions and  $996$ cases of them may lead to the  existence of a subgraph isomorphic to the graph $C_6--C_6$ in the graph $K(\alpha,\beta)$.

\subsection{$\mathbf{C_6--C_6(C_6--C_6)}$}
\begin{figure}[ht]
\psscalebox{0.9 0.9} 
{
\begin{pspicture}(0,-2.0985577)(5.194231,2.0985577)
\rput[bl](0.2,-2.0985577){$\mathbf{24) \ C_6--C_6(C_6--C_6)}$}
\psdots[linecolor=black, dotsize=0.4](2.5971153,1.9014423)
\psdots[linecolor=black, dotsize=0.4](2.5971153,0.70144236)
\psdots[linecolor=black, dotsize=0.4](2.5971153,-0.49855766)
\psdots[linecolor=black, dotsize=0.4](3.7971153,-0.098557666)
\psdots[linecolor=black, dotsize=0.4](3.7971153,0.70144236)
\psdots[linecolor=black, dotsize=0.4](3.7971153,1.5014423)
\psdots[linecolor=black, dotsize=0.4](1.3971153,1.5014423)
\psdots[linecolor=black, dotsize=0.4](1.3971153,0.70144236)
\psdots[linecolor=black, dotsize=0.4](1.3971153,-0.098557666)
\psdots[linecolor=black, dotsize=0.4](0.19711533,0.70144236)
\psdots[linecolor=black, dotsize=0.4](0.19711533,-0.49855766)
\psdots[linecolor=black, dotsize=0.4](1.3971153,-1.2985576)
\psdots[linecolor=black, dotsize=0.4](3.7971153,-1.2985576)
\psdots[linecolor=black, dotsize=0.4](4.997115,-0.49855766)
\psline[linecolor=black, linewidth=0.04](2.5971153,1.9014423)(2.5971153,0.70144236)(2.5971153,-0.49855766)(1.3971153,-0.098557666)(1.3971153,0.70144236)(1.3971153,1.5014423)(2.5971153,1.9014423)(3.7971153,1.5014423)(3.7971153,0.70144236)(3.7971153,-0.098557666)(2.5971153,-0.49855766)
\psline[linecolor=black, linewidth=0.04](2.5971153,0.70144236)(1.3971153,-1.2985576)(0.19711533,-0.49855766)(0.19711533,0.70144236)(1.3971153,1.5014423)
\psline[linecolor=black, linewidth=0.04](1.3971153,-1.2985576)(3.7971153,-1.2985576)(4.997115,-0.49855766)(3.7971153,1.5014423)
\end{pspicture}
}
\end{figure}
By considering the $996$ cases related to the existence of $C_6--C_6$ in the graph $K(\alpha,\beta)$, it can be seen that there are $5119$ cases for the relations of four $C_6$ cycles in this structure. Using GAP \cite{gap}, we see that all groups with two generators $h_2$ and $h_3$ and four relations which are between $4983$ cases of these $5119$ cases are finite and solvable, or just finite, that is a contradiction with the assumptions. So, there are just $136$ cases for the relations of these cycles which may lead to the existence of a subgraph isomorphic to the graph $C_6--C_6(C_6--C_6)$ in $K(\alpha,\beta)$. In the following, we show that these $136$ cases lead to contradictions and so, the graph $K(\alpha,\beta)$ contains no subgraph isomorphic to the graph $C_6--C_6(C_6--C_6)$.
\begin{enumerate}
\item[(1)]$ R_1:h_2^3h_3^3=1$, $R_2:h_2h_3h_2^{-1}h_3^{-1}h_2h_3^{-1}h_2^{-1}h_3=1$, $R_3:h_2^2h_3^{-2}(h_2^{-1}h_3)^2=1$,\\$R_4:h_2h_3h_2^{-1}h_3^{-1}h_2h_3^{-1}h_2^{-1}h_3=1$:\\
$\Rightarrow \langle h_2,h_3\rangle=\langle h_2\rangle$ is abelian, a contradiction.

\item[(2)]$ R_1:h_2^3h_3^3=1$, $R_2:h_2h_3^{-1}h_2^{-1}h_3h_2^{-1}h_3^{-1}h_2h_3=1$, $R_3:h_2^2h_3^{-2}(h_2^{-1}h_3)^2=1$,\\$R_4:h_2h_3^{-1}h_2^{-1}h_3h_2^{-1}h_3^{-1}h_2h_3=1$:\\
By interchanging $h_2$ and $h_3$ in (1) and with the same discussion, there is a contradiction.

\item[(3)]$ R_1:h_2^3h_3h_2^{-2}h_3=1$, $R_2:h_2h_3h_2h_3^{-1}h_2^{-1}h_3^2=1$, $R_3:h_2^2h_3^{-2}(h_2^{-1}h_3)^2=1$,\\$R_4:h_2h_3h_2h_3^{-1}h_2^{-1}h_3^2=1$:\\
$\Rightarrow\langle h_2,h_3\rangle\cong BS(2,1)$ is solvable, a contradiction.

\item[(4)]$ R_1:h_2^3h_3h_2^{-2}h_3=1$, $R_2:h_2h_3h_2h_3^{-1}h_2^{-1}h_3^2=1$, $R_3:h_2^2h_3^{-2}(h_2^{-1}h_3)^2=1$,\\$R_4:h_2h_3^{-1}h_2^{-1}h_3h_2^{-1}h_3^{-1}h_2h_3=1$:\\
$\Rightarrow\langle h_2,h_3\rangle\cong BS(2,1)$ is solvable, a contradiction.

\item[(5)]$ R_1:h_2^3h_3h_2^{-2}h_3=1$, $R_2:h_2h_3^2h_2^{-1}h_3^{-1}h_2h_3=1$, $R_3:h_2^2h_3^{-1}h_2(h_3h_2^{-1})^2h_3=1$,\\$R_4:h_2(h_3h_2^{-1})^2h_3^{-2}h_2=1$:\\
$\Rightarrow\langle h_2,h_3\rangle\cong BS(1,1)$ is solvable, a contradiction.

\item[(6)]$ R_1:h_2^3h_3^{-1}h_2^2h_3=1$, $R_2:h_2h_3h_2^{-1}h_3^{-2}h_2^{-1}h_3=1$, $R_3:h_2h_3h_2^{-1}h_3(h_2h_3^{-1})^2h_2=1$,\\$R_4:h_2h_3^{-2}h_2^{-1}h_3h_2^{-2}h_3=1$:\\
$\Rightarrow\langle h_2,h_3\rangle\cong BS(3,1)$ is solvable, a contradiction.

\item[(7)]$ R_1:h_2^3h_3^{-1}h_2^2h_3=1$, $R_2:h_2h_3h_2^{-1}h_3^{-1}h_2h_3^{-1}h_2^{-1}h_3=1$, $R_3:h_2h_3h_2^{-1}h_3(h_2h_3^{-1})^2h_2=1$,\\$R_4:h_2h_3^{-2}h_2^{-1}h_3h_2^{-2}h_3=1$:\\
$\Rightarrow\langle h_2,h_3\rangle\cong BS(3,1)$ is solvable, a contradiction.

\item[(8)]$ R_1:h_2^3h_3^{-1}h_2h_3^2=1$, $R_2:h_2^2h_3h_2h_3^{-1}h_2h_3=1$, $R_3:h_2^3h_3^{-1}h_2h_3^2=1$,\\$R_4:h_2(h_3h_2^{-1})^2h_3^{-2}h_2=1$:\\
$\Rightarrow \langle h_2,h_3\rangle=\langle h_2\rangle$ is abelian, a contradiction.

\item[(9)]$ R_1:h_2^2(h_2h_3^{-1})^3h_2=1$, $R_2:h_2^2h_3^{-2}(h_2^{-1}h_3)^2=1$, $R_3:h_2h_3^{-1}h_2^{-1}h_3^2h_2^{-1}h_3^{-1}h_2=1$,\\$R_4:h_2h_3^{-1}h_2^{-1}h_3^2h_2^{-1}h_3^{-1}h_2=1$:\\
$\Rightarrow \langle h_2,h_3\rangle=\langle h_2\rangle$ is abelian, a contradiction.

\item[(10)]$ R_1:h_2^2h_3h_2^2h_3^{-1}h_2=1$, $R_2:h_2^2h_3h_2^2h_3^{-1}h_2=1$, $R_3:h_2h_3^2(h_2h_3^{-1})^2h_2=1$,\\$R_4:h_2h_3^{-2}h_2^{-1}h_3h_2^{-2}h_3=1$:\\
$\Rightarrow\langle h_2,h_3\rangle\cong BS(1,3)$ is solvable, a contradiction.

\item[(11)]$ R_1:h_2^2h_3h_2^2h_3^{-1}h_2=1$, $R_2:h_2^2h_3h_2^2h_3^{-1}h_2=1$, $R_3:h_2h_3^2(h_2h_3^{-1})^2h_2=1$,\\$R_4:h_2h_3^{-1}h_2(h_3h_2^{-1})^2h_3^{-1}h_2=1$:\\
$\Rightarrow\langle h_2,h_3\rangle\cong BS(1,3)$ is solvable, a contradiction.

\item[(12)]$ R_1:h_2^2h_3h_2^2h_3^{-1}h_2=1$, $R_2:h_2^2h_3h_2^2h_3^{-1}h_2=1$, $R_3:h_2h_3^2h_2^{-1}h_3^{-2}h_2=1$,\\$R_4:h_2h_3^{-1}(h_2^{-1}h_3)^2h_2^{-1}h_3^{-1}h_2=1$:\\
$\Rightarrow\langle h_2,h_3\rangle\cong BS(1,-1)$ is solvable, a contradiction.

\item[(13)]$ R_1:h_2^2h_3(h_2h_3^{-1})^2h_2=1$, $R_2:h_2^2h_3^{-2}(h_2^{-1}h_3)^2=1$, $R_3:h_2h_3^{-1}h_2^{-1}h_3^2h_2^{-1}h_3^{-1}h_2=1$,\\$R_4:h_2h_3^{-1}h_2^{-1}h_3^2h_2^{-1}h_3^{-1}h_2=1$:\\
$\Rightarrow\langle h_2,h_3\rangle\cong BS(1,2)$ is solvable, a contradiction.

\item[(14)]$ R_1:h_2^2h_3(h_2h_3^{-1})^2h_2=1$, $R_2:h_2h_3h_2^{-1}h_3^{-1}h_2h_3^{-1}h_2^{-1}h_3=1$, $R_3:h_2^2h_3^{-2}(h_2^{-1}h_3)^2=1$,\\$R_4:h_2h_3h_2^{-1}h_3^{-1}h_2h_3^{-1}h_2^{-1}h_3=1$:\\
$\Rightarrow\langle h_2,h_3\rangle\cong BS(1,2)$ is solvable, a contradiction.

\item[(15)]$ R_1:h_2^2h_3^3h_2^{-1}h_3=1$, $R_2:h_2h_3^2h_2h_3h_2^{-1}h_3=1$, $R_3:h_2^2h_3^3h_2^{-1}h_3=1$,\\$R_4:h_2(h_3h_2^{-1})^2h_3^{-2}h_2=1$:\\
By interchanging $h_2$ and $h_3$ in (8) and with the same discussion, there is a contradiction.

\item[(16)]$ R_1:h_2^2h_3(h_3h_2^{-1})^2h_3=1$, $R_2:h_2h_3h_2^{-1}h_3^2h_2^{-1}h_3^{-1}h_2=1$, $R_3:(h_2h_3)^2h_2^{-2}h_3=1$,\\$R_4:h_2h_3^{-2}h_2^{-1}h_3^3=1$:\\
$\Rightarrow\langle h_2,h_3\rangle\cong BS(1,1)$ is solvable, a contradiction.

\item[(17)]$ R_1:h_2^2h_3h_2^{-2}h_3^{-1}h_2=1$, $R_2:h_2^2h_3h_2^{-2}h_3^{-1}h_2=1$, $R_3:h_2h_3^2h_2^{-1}(h_3^{-1}h_2)^2=1$,\\$R_4:h_2h_3^{-1}h_2^{-1}h_3^2h_2h_3^{-1}h_2=1$:\\
$\Rightarrow\langle h_2,h_3\rangle\cong BS(1,-1)$ is solvable, a contradiction.

\item[(18)]$ R_1:h_2^2h_3h_2^{-2}h_3^{-1}h_2=1$, $R_2:h_2^2h_3h_2^{-2}h_3^{-1}h_2=1$, $R_3:h_2h_3h_2^{-1}h_3h_2h_3^{-2}h_2=1$,\\$R_4:h_2h_3^{-1}h_2^{-1}h_3^2h_2h_3^{-1}h_2=1$:\\
$\Rightarrow\langle h_2,h_3\rangle\cong BS(1,-5)$ is solvable, a contradiction.

\item[(19)]$ R_1:h_2^2h_3h_2^{-2}h_3^{-1}h_2=1$, $R_2:h_2^2h_3h_2^{-2}h_3^{-1}h_2=1$, $R_3:h_2h_3h_2^{-1}h_3h_2h_3^{-2}h_2=1$,\\$R_4:h_2h_3^{-1}(h_2^{-1}h_3)^2h_2h_3^{-1}h_2=1$:\\
$\Rightarrow\langle h_2,h_3\rangle\cong BS(1,-3)$ is solvable, a contradiction.

\item[(20)]$ R_1:h_2^2h_3h_2^{-2}h_3^{-1}h_2=1$, $R_2:h_2^2h_3h_2^{-2}h_3^{-1}h_2=1$, $R_3:h_2h_3h_2^{-1}h_3(h_2h_3^{-1})^2h_2=1$,\\$R_4:(h_2h_3^{-1}h_2)^2h_3^2=1$:\\
$\Rightarrow\langle h_2,h_3\rangle\cong BS(1,-5)$ is solvable, a contradiction.

\item[(21)]$ R_1:h_2^2h_3h_2^{-2}h_3^{-1}h_2=1$, $R_2:h_2^2h_3h_2^{-2}h_3^{-1}h_2=1$, $R_3:h_2(h_3h_2^{-1})^2(h_3^{-1}h_2)^2=1$,\\$R_4:(h_2h_3^{-1}h_2)^2h_3h_2^{-1}h_3=1$:\\
$\Rightarrow\langle h_2,h_3\rangle\cong BS(1,1)$ is solvable, a contradiction.

\item[(22)]$ R_1:h_2^2h_3h_2^{-2}h_3^{-1}h_2=1$, $R_2:h_2h_3h_2h_3^{-2}h_2h_3=1$, $R_3:h_2^2h_3^{-2}(h_2^{-1}h_3)^2=1$,\\$R_4:h_2h_3h_2^{-1}h_3^{-2}h_2^{-1}h_3=1$:\\
$\Rightarrow\langle h_2,h_3\rangle\cong BS(1,2)$ is solvable, a contradiction.

\item[(23)]$ R_1:h_2^2h_3h_2^{-2}h_3^{-1}h_2=1$, $R_2:h_2h_3h_2h_3^{-2}h_2h_3=1$, $R_3:h_2^2h_3^{-2}h_2h_3h_2^{-1}h_3=1$,\\$R_4:h_2h_3h_2^{-1}h_3^{-2}h_2h_3=1$:\\
$\Rightarrow \langle h_2,h_3\rangle=\langle h_2\rangle$ is abelian, a contradiction.

\item[(24)]$ R_1:h_2^2h_3h_2^{-2}h_3^{-1}h_2=1$, $R_2:h_2h_3h_2h_3^{-2}h_2h_3=1$, $R_3:h_2(h_2h_3^{-1})^2h_2^{-1}h_3^2=1$,\\$R_4:(h_2h_3^{-1}h_2)^2h_3^2=1$:\\
$\Rightarrow\langle h_2,h_3\rangle\cong BS(1,4)$ is solvable, a contradiction.

\item[(25)]$ R_1:h_2^2h_3h_2^{-2}h_3^{-1}h_2=1$, $R_2:h_2h_3h_2h_3^{-2}h_2h_3=1$, $R_3:h_2(h_2h_3^{-1})^2h_2h_3^2=1$,\\$R_4:h_2h_3^{-1}h_2^{-1}h_3^2h_2h_3^{-1}h_2=1$:\\
By interchanging $h_2$ and $h_3$ in (16) and with the same discussion, there is a contradiction.

\item[(26)]$ R_1:h_2^2h_3h_2^{-2}h_3^{-1}h_2=1$, $R_2:h_2h_3h_2h_3^{-2}h_2h_3=1$, $R_3:h_2(h_2h_3^{-1})^2h_2h_3h_2^{-1}h_3=1$,\\$R_4:h_2h_3^{-1}h_2^{-1}h_3^2h_2h_3^{-1}h_2=1$:\\
$\Rightarrow\langle h_2,h_3\rangle\cong BS(1,-3)$ is solvable, a contradiction.

\item[(27)]$ R_1:h_2^2h_3h_2^{-2}h_3^{-1}h_2=1$, $R_2:h_2h_3h_2h_3^{-2}h_2h_3=1$, $R_3:h_2(h_2h_3^{-1})^2h_2h_3h_2^{-1}h_3=1$,\\$R_4:(h_2h_3^{-1}h_2)^2h_3^2=1$:\\
$\Rightarrow\langle h_2,h_3\rangle\cong BS(1,-3)$ is solvable, a contradiction.

\item[(28)]$ R_1:h_2^2h_3h_2^{-2}h_3^2=1$, $R_2:h_2h_3^{-2}h_2h_3^{-1}(h_2^{-1}h_3)^2=1$, $R_3:(h_2h_3)^2h_2^{-2}h_3=1$,\\$R_4:h_2(h_3h_2^{-1})^2h_3^{-2}h_2=1$:\\
$\Rightarrow\langle h_2,h_3\rangle\cong BS(1,1)$ is solvable, a contradiction.

\item[(29)]$ R_1:h_2^2h_3h_2^{-2}h_3^2=1$, $R_2:h_2h_3^{-2}h_2h_3^{-1}(h_2^{-1}h_3)^2=1$, $R_3:(h_2h_3)^2(h_2^{-1}h_3)^2=1$,\\$R_4:h_2h_3^{-1}h_2^{-1}h_3h_2^{-1}h_3^{-2}h_2=1$:\\
$\Rightarrow$  $G$ has a torsion element, a contradiction.

\item[(30)]$ R_1:h_2^2h_3h_2^{-1}h_3^{-1}h_2^{-1}h_3=1$, $R_2:h_2^2h_3^{-1}(h_2^{-1}h_3)^2h_3=1$, $R_3:h_2^2h_3^{-1}h_2^{-1}h_3^2h_2^{-1}h_3=1$,\\$R_4:h_2h_3^{-1}h_2(h_3h_2^{-1})^2h_3^2=1$:\\
$\Rightarrow\langle h_2,h_3\rangle\cong BS(1,1)$ is solvable, a contradiction.

\item[(31)]$ R_1:h_2^2h_3h_2^{-1}h_3^{-1}h_2^{-1}h_3=1$, $R_2:h_2^2h_3^{-2}(h_2^{-1}h_3)^2=1$, $R_3:(h_2h_3)^2h_2^{-2}h_3=1$,\\$R_4:h_2h_3^{-2}h_2^{-1}h_3^3=1$:\\
By interchanging $h_2$ and $h_3$ in (22) and with the same discussion, there is a contradiction.

\item[(32)]$ R_1:h_2^2h_3h_2^{-1}h_3^{-1}h_2^{-1}h_3=1$, $R_2:h_2h_3^3h_2^{-1}h_3^2=1$, $R_3:h_2^2h_3^{-1}h_2^{-1}h_3^2h_2^{-1}h_3=1$,\\$R_4:h_2h_3^{-1}h_2(h_3h_2^{-1})^2h_3^2=1$:\\
By interchanging $h_2$ and $h_3$ in (6) and with the same discussion, there is a contradiction.

\item[(33)]$ R_1:h_2^2h_3h_2^{-1}h_3^{-2}h_2=1$, $R_2:h_2^2h_3^{-1}h_2^{-1}h_3^{-1}h_2h_3=1$, $R_3:h_2^2h_3h_2^{-1}h_3^{-2}h_2=1$,\\$R_4:h_2h_3^{-2}h_2h_3h_2h_3^{-1}h_2=1$:\\
$\Rightarrow\langle h_2,h_3\rangle\cong BS(2,1)$ is solvable, a contradiction.

\item[(34)]$ R_1:h_2^2h_3h_2^{-1}h_3^{-2}h_2=1$, $R_2:h_2^2h_3^{-1}(h_3^{-1}h_2)^2h_3=1$, $R_3:h_2^2h_3h_2^{-1}h_3^{-2}h_2=1$,\\$R_4:h_2h_3^{-2}h_2h_3h_2h_3^{-1}h_2=1$:\\
$\Rightarrow$ $G$ is solvable, a contradiction.

\item[(35)]$ R_1:h_2^2h_3h_2^{-1}h_3^{-2}h_2=1$, $R_2:h_2(h_2h_3^{-1})^2h_3^{-1}h_2h_3=1$, $R_3:h_2^2h_3h_2^{-1}h_3^{-2}h_2=1$,\\$R_4:h_2h_3^{-2}h_2h_3h_2h_3^{-1}h_2=1$:\\
$\Rightarrow\langle h_2,h_3\rangle\cong BS(1,1)$ is solvable, a contradiction.

\item[(36)]$ R_1:h_2^2h_3h_2^{-1}(h_3^{-1}h_2)^2=1$, $R_2:(h_2^2h_3^{-1})^2h_2^{-1}h_3=1$, $R_3:h_2h_3h_2^{-1}h_3^{-1}h_2^2h_3^{-1}h_2=1$,\\$R_4:(h_2h_3^{-1})^3h_2h_3h_2^{-1}h_3=1$:\\
$\Rightarrow\langle h_2,h_3\rangle\cong BS(1,1)$ is solvable, a contradiction.

\item[(37)]$ R_1:h_2^2h_3h_2^{-1}h_3h_2h_3^{-1}h_2=1$, $R_2:h_2^2h_3^{-2}(h_2^{-1}h_3)^2=1$, $R_3:h_2h_3^{-1}h_2^{-1}h_3^2h_2^{-1}h_3^{-1}h_2=1$,\\$R_4:h_2h_3^{-1}h_2^{-1}h_3^2h_2^{-1}h_3^{-1}h_2=1$:\\
$\Rightarrow\langle h_2,h_3\rangle\cong BS(2,1)$ is solvable, a contradiction.

\item[(38)]$ R_1:h_2^2(h_3h_2^{-1}h_3)^2=1$, $R_2:h_2(h_2h_3^{-2})^2h_2=1$, $R_3:h_2h_3h_2^{-1}h_3^2h_2^{-1}h_3^{-1}h_2=1$,\\$R_4:h_2h_3^{-1}h_2^{-1}h_3^2h_2h_3^{-1}h_2=1$:\\
$\Rightarrow$  $G$ has a torsion element, a contradiction.

\item[(39)]$ R_1:h_2^2(h_3h_2^{-1}h_3)^2=1$, $R_2:h_2(h_2h_3^{-2})^2h_2=1$, $R_3:h_2h_3^{-1}h_2^{-1}h_3^2h_2h_3^{-1}h_2=1$,\\$R_4:h_2h_3^{-2}h_2^{-1}h_3h_2^{-2}h_3=1$:\\
$\Rightarrow$  $G$ has a torsion element, a contradiction.

\item[(40)]$ R_1:h_2^2(h_3h_2^{-1}h_3)^2=1$, $R_2:h_2h_3(h_3h_2^{-1})^2h_3^{-1}h_2=1$, $R_3:(h_2h_3)^2h_2^{-2}h_3=1$,\\$R_4:h_2h_3^{-2}h_2^{-1}h_3^3=1$:\\
By interchanging $h_2$ and $h_3$ in (24) and with the same discussion, there is a contradiction.

\item[(41)]$ R_1:h_2^2(h_3h_2^{-1}h_3)^2=1$, $R_2:h_2h_3^{-1}h_2h_3(h_3h_2^{-1})^2h_3=1$, $R_3:(h_2h_3)^2h_2^{-2}h_3=1$,\\$R_4:h_2h_3^{-2}h_2^{-1}h_3^3=1$:\\
By interchanging $h_2$ and $h_3$ in (27) and with the same discussion, there is a contradiction.

\item[(42)]$ R_1:h_2^2(h_3h_2^{-1}h_3)^2=1$, $R_2:h_2h_3^{-1}h_2(h_3h_2^{-1})^2h_3^2=1$, $R_3:h_2(h_3h_2^{-1})^2h_3^{-1}h_2h_3=1$,\\$R_4:h_2h_3^{-1}h_2^{-1}h_3h_2^{-1}h_3^{-2}h_2=1$:\\
$\Rightarrow$  $h_3=1$, a contradiction.

\item[(43)]$ R_1:h_2^2(h_3h_2^{-1}h_3)^2=1$, $R_2:h_2h_3^{-1}h_2(h_3h_2^{-1})^2h_3^2=1$, $R_3:h_2h_3^{-2}h_2^{-1}h_3^3=1$,\\$R_4:h_2h_3^{-2}h_2^{-1}h_3^3=1$:\\
By interchanging $h_2$ and $h_3$ in (20) and with the same discussion, there is a contradiction.

\item[(44)]$ R_1:h_2^2(h_3h_2^{-1})^2h_3^2=1$, $R_2:h_2^2h_3^{-1}h_2^{-1}h_3^2h_2^{-1}h_3=1$, $R_3:h_2^2(h_3h_2^{-1})^2h_3^2=1$,\\$R_4:h_2h_3h_2h_3^{-1}h_2^{-2}h_3=1$:\\
$\Rightarrow\langle h_2,h_3\rangle\cong BS(-1,1)$ is solvable, a contradiction.

\item[(45)]$ R_1:h_2^2(h_3h_2^{-1})^2h_3^2=1$, $R_2:h_2^2h_3^{-1}h_2^{-1}h_3^2h_2^{-1}h_3=1$, $R_3:h_2^2(h_3h_2^{-1})^2h_3^2=1$,\\$R_4:h_2h_3h_2h_3^{-1}(h_2^{-1}h_3)^2=1$:\\
$\Rightarrow$  $G$ has a torsion element, a contradiction.

\item[(46)]$ R_1:h_2^2(h_3h_2^{-1})^2h_3^2=1$, $R_2:h_2^2h_3^{-1}h_2^{-1}h_3^2h_2^{-1}h_3=1$, $R_3:(h_2h_3)^2h_2^{-2}h_3=1$,\\$R_4:h_2h_3^{-1}h_2(h_3h_2^{-1})^2h_3^2=1$:\\
$\Rightarrow\langle h_2,h_3\rangle\cong BS(1,1)$ is solvable, a contradiction.

\item[(47)]$ R_1:h_2^2(h_3h_2^{-1})^2h_3^2=1$, $R_2:h_2^2h_3^{-1}h_2^{-1}h_3^2h_2^{-1}h_3=1$, $R_3:(h_2h_3)^2h_2^{-2}h_3=1$,\\$R_4:(h_2h_3^{-1})^2(h_2^{-1}h_3)^2h_3=1$:\\
$\Rightarrow\langle h_2,h_3\rangle\cong BS(1,1)$ is solvable, a contradiction.

\item[(48)]$ R_1:h_2^2(h_3h_2^{-1})^2h_3^2=1$, $R_2:h_2^2h_3^{-1}h_2^{-1}h_3^2h_2^{-1}h_3=1$, $R_3:(h_2h_3)^2(h_2^{-1}h_3)^2=1$,\\$R_4:(h_2h_3^{-1})^2(h_2^{-1}h_3)^2h_3=1$:\\
$\Rightarrow\langle h_2,h_3\rangle\cong BS(1,1)$ is solvable, a contradiction.

\item[(49)]$ R_1:h_2^2(h_3h_2^{-1})^2h_3^2=1$, $R_2:h_2^2h_3^{-1}h_2^{-1}h_3^2h_2^{-1}h_3=1$, $R_3:h_2h_3^2h_2^{-1}h_3^3=1$,\\$R_4:h_2h_3^2h_2^{-1}h_3^3=1$:\\
By interchanging $h_2$ and $h_3$ in (10) and with the same discussion, there is a contradiction.

\item[(50)]$ R_1:h_2^2(h_3h_2^{-1})^2h_3^2=1$, $R_2:h_2h_3h_2h_3^{-1}h_2^{-2}h_3=1$, $R_3:h_2^2(h_3h_2^{-1})^2h_3^2=1$,\\$R_4:(h_2h_3^{-1})^2h_2^{-1}h_3^2h_2^{-1}h_3=1$:\\
$\Rightarrow\langle h_2,h_3\rangle\cong BS(1,3)$ is solvable, a contradiction.

\item[(51)]$ R_1:h_2^2(h_3h_2^{-1})^2h_3^2=1$, $R_2:h_2h_3h_2h_3^{-1}h_2^{-2}h_3=1$, $R_3:h_2^2h_3^{-1}h_2^{-2}h_3^2=1$,\\$R_4:h_2h_3^{-2}h_2(h_3h_2^{-1})^2h_3=1$:\\
$\Rightarrow\langle h_2,h_3\rangle\cong BS(1,1)$ is solvable, a contradiction.

\item[(52)]$ R_1:h_2^2(h_3h_2^{-1})^2h_3^2=1$, $R_2:h_2h_3h_2h_3^{-1}(h_2^{-1}h_3)^2=1$, $R_3:h_2^2(h_3h_2^{-1})^2h_3^2=1$,\\$R_4:(h_2h_3^{-1})^2h_2^{-1}h_3^2h_2^{-1}h_3=1$:\\
$\Rightarrow$  $h_3=1$, a contradiction.

\item[(53)]$ R_1:h_2^2(h_3h_2^{-1})^2h_3^2=1$, $R_2:h_2h_3h_2h_3^{-1}(h_2^{-1}h_3)^2=1$, $R_3:h_2^2h_3^{-1}h_2^{-2}h_3^2=1$,\\$R_4:h_2h_3^{-2}h_2(h_3h_2^{-1})^2h_3=1$:\\
$\Rightarrow\langle h_2,h_3\rangle\cong BS(1,2)$ is solvable, a contradiction.

\item[(54)]$ R_1:h_2^2(h_3h_2^{-1})^2h_3^2=1$, $R_2:h_2h_3^{-1}h_2^{-1}h_3h_2^{-1}h_3^{-1}h_2h_3=1$, $R_3:h_2^2h_3^{-1}h_2^{-1}h_3^2h_2^{-1}h_3=1$,\\$R_4:h_2h_3^{-1}h_2(h_3h_2^{-1})^2h_3^2=1$:\\
$\Rightarrow$  $h_3=1$, a contradiction.

\item[(55)]$ R_1:h_2^2(h_3h_2^{-1})^2h_3^2=1$, $R_2:(h_2h_3^{-1})^2h_2^{-1}h_3^2h_2^{-1}h_3=1$, $R_3:(h_2h_3)^2h_2^{-2}h_3=1$,\\$R_4:h_2h_3^{-1}h_2(h_3h_2^{-1})^2h_3^2=1$:\\
$\Rightarrow\langle h_2,h_3\rangle\cong BS(1,1)$ is solvable, a contradiction.

\item[(56)]$ R_1:h_2^2(h_3h_2^{-1})^2h_3^2=1$, $R_2:(h_2h_3^{-1})^2h_2^{-1}h_3^2h_2^{-1}h_3=1$, $R_3:(h_2h_3)^2h_2^{-2}h_3=1$,\\$R_4:(h_2h_3^{-1})^2(h_2^{-1}h_3)^2h_3=1$:\\
$\Rightarrow\langle h_2,h_3\rangle\cong BS(1,1)$ is solvable, a contradiction.

\item[(57)]$ R_1:h_2^2(h_3h_2^{-1})^2h_3^2=1$, $R_2:(h_2h_3^{-1})^2h_2^{-1}h_3^2h_2^{-1}h_3=1$, $R_3:(h_2h_3)^2(h_2^{-1}h_3)^2=1$,\\$R_4:(h_2h_3^{-1})^2(h_2^{-1}h_3)^2h_3=1$:\\
$\Rightarrow\langle h_2,h_3\rangle\cong BS(1,1)$ is solvable, a contradiction.

\item[(58)]$ R_1:h_2^2(h_3h_2^{-1})^2h_3^2=1$, $R_2:(h_2h_3^{-1})^2h_2^{-1}h_3^2h_2^{-1}h_3=1$, $R_3:h_2h_3^2h_2^{-1}h_3^3=1$,\\$R_4:h_2h_3^2h_2^{-1}h_3^3=1$:\\
By interchanging $h_2$ and $h_3$ in (11) and with the same discussion, there is a contradiction.

\item[(59)]$ R_1:h_2^2h_3^{-1}h_2^{-2}h_3^2=1$, $R_2:h_2h_3^{-2}h_2(h_3h_2^{-1})^2h_3=1$, $R_3:(h_2h_3)^2h_2^{-2}h_3=1$,\\$R_4:h_2h_3^{-1}h_2(h_3h_2^{-1})^2h_3^2=1$:\\
$\Rightarrow\langle h_2,h_3\rangle\cong BS(-3,1)$ is solvable, a contradiction.

\item[(60)]$ R_1:h_2^2h_3^{-1}h_2^{-2}h_3^2=1$, $R_2:h_2h_3^{-2}h_2(h_3h_2^{-1})^2h_3=1$, $R_3:(h_2h_3)^2h_2^{-2}h_3=1$,\\$R_4:(h_2h_3^{-1})^2(h_2^{-1}h_3)^2h_3=1$:\\
$\Rightarrow\langle h_2,h_3\rangle\cong BS(-3,1)$ is solvable, a contradiction.

\item[(61)]$ R_1:h_2^2h_3^{-1}h_2^{-2}h_3^2=1$, $R_2:h_2h_3^{-2}h_2(h_3h_2^{-1})^2h_3=1$, $R_3:(h_2h_3)^2(h_2^{-1}h_3)^2=1$,\\$R_4:(h_2h_3^{-1})^2(h_2^{-1}h_3)^2h_3=1$:\\
$\Rightarrow\langle h_2,h_3\rangle\cong BS(1,-1)$ is solvable, a contradiction.

\item[(62)]$ R_1:h_2^2h_3^{-1}h_2^{-2}h_3^2=1$, $R_2:h_2h_3^{-2}h_2(h_3h_2^{-1})^2h_3=1$, $R_3:h_2h_3^2h_2^{-1}h_3^3=1$,\\$R_4:h_2h_3^2h_2^{-1}h_3^3=1$:\\
By interchanging $h_2$ and $h_3$ in (12) and with the same discussion, there is a contradiction.

\item[(63)]$ R_1:h_2^2h_3^{-1}h_2^{-2}h_3^2=1$, $R_2:h_2h_3^{-2}h_2h_3^{-1}(h_2^{-1}h_3)^2=1$, $R_3:(h_2h_3)^2h_2^{-2}h_3=1$,\\$R_4:h_2(h_3h_2^{-1})^2h_3^{-2}h_2=1$:\\
$\Rightarrow\langle h_2,h_3\rangle\cong BS(-3,1)$ is solvable, a contradiction.

\item[(64)]$ R_1:h_2^2h_3^{-1}h_2^{-2}h_3^2=1$, $R_2:h_2h_3^{-2}h_2h_3^{-1}(h_2^{-1}h_3)^2=1$, $R_3:(h_2h_3)^2(h_2^{-1}h_3)^2=1$,\\$R_4:h_2h_3^{-1}h_2^{-1}h_3h_2^{-1}h_3^{-2}h_2=1$:\\
$\Rightarrow\langle h_2,h_3\rangle\cong BS(1,-1)$ is solvable, a contradiction.

\item[(65)]$ R_1:h_2^2(h_3^{-1}h_2^{-1})^2h_3=1$, $R_2:h_2h_3^{-1}h_2(h_3h_2^{-1})^2h_3^2=1$, $R_3:h_2^2h_3^{-1}h_2^{-1}h_3^2h_2^{-1}h_3=1$,\\$R_4:h_2h_3^{-1}h_2(h_3h_2^{-1})^2h_3^2=1$:\\
$\Rightarrow$  $G$ has a torsion element, a contradiction.

\item[(66)]$ R_1:h_2^2(h_3^{-1}h_2^{-1})^2h_3=1$, $R_2:(h_2h_3^{-1})^2(h_2^{-1}h_3)^2h_3=1$, $R_3:h_2^2h_3^{-1}h_2^{-1}h_3^2h_2^{-1}h_3=1$,\\$R_4:h_2h_3^{-1}h_2(h_3h_2^{-1})^2h_3^2=1$:\\
$\Rightarrow\langle h_2,h_3\rangle\cong BS(1,2)$ is solvable, a contradiction.

\item[(67)]$ R_1:h_2^2h_3^{-1}h_2^{-1}h_3^{-1}h_2h_3=1$, $R_2:h_2h_3(h_2h_3^{-1})^2h_3^{-1}h_2=1$, $R_3:h_2h_3h_2^{-1}h_3^{-1}h_2h_3^{-1}h_2^{-1}h_3=1$,\\$R_4:h_2h_3^{-2}h_2h_3h_2h_3^{-1}h_2=1$:\\
$\Rightarrow \langle h_2,h_3\rangle=\langle h_2\rangle$ is abelian, a contradiction.

\item[(68)]$ R_1:h_2^2h_3^{-1}h_2^{-1}h_3h_2h_3=1$, $R_2:h_2h_3^{-2}h_2h_3^3=1$, $R_3:h_2(h_3h_2^{-1})^2h_3^{-2}h_2=1$,\\$R_4:(h_2h_3^{-1})^2h_2h_3^2h_2^{-1}h_3=1$:\\
By interchanging $h_2$ and $h_3$ in (5) and with the same discussion, there is a contradiction.

\item[(69)]$ R_1:h_2^2h_3^{-1}h_2^{-1}h_3^2h_2^{-1}h_3=1$, $R_2:h_2h_3h_2^{-1}h_3^2h_2^{-1}h_3^{-1}h_2=1$, $R_3:h_2h_3^{-3}h_2^2h_3^{-1}h_2=1$,\\$R_4:(h_2h_3^{-1}h_2)^2h_3^2=1$:\\
By interchanging $h_2$ and $h_3$ in (39) and with the same discussion, there is a contradiction.

\item[(70)]$ R_1:h_2^2h_3^{-1}h_2^{-1}h_3^2h_2^{-1}h_3=1$, $R_2:h_2h_3^{-1}h_2(h_3h_2^{-1})^2h_3^2=1$, $R_3:h_2h_3^3h_2^{-1}h_3^2=1$,\\$R_4:h_2h_3^{-1}h_2^{-1}h_3h_2^{-1}h_3^{-1}h_2h_3=1$:\\
By interchanging $h_2$ and $h_3$ in (7) and with the same discussion, there is a contradiction.

\item[(71)]$ R_1:h_2^2h_3^{-1}h_2^{-1}h_3^2h_2^{-1}h_3=1$, $R_2:h_2h_3^{-1}h_2(h_3h_2^{-1})^2h_3^2=1$, $R_3:h_2(h_3h_2^{-1})^2h_3^{-1}h_2h_3=1$,\\$R_4:h_2h_3^{-1}h_2(h_3h_2^{-1})^2h_3^2=1$:\\
$\Rightarrow\langle h_2,h_3\rangle\cong BS(1,-2)$ is solvable, a contradiction.

\item[(72)]$ R_1:h_2^2h_3^{-1}(h_2^{-1}h_3)^2h_3=1$, $R_2:h_2h_3(h_2h_3^{-1})^2h_3^{-1}h_2=1$, $R_3:h_2(h_2h_3^{-2})^2h_2=1$,\\$R_4:(h_2h_3)^2(h_2^{-1}h_3)^2=1$:\\
$\Rightarrow\langle h_2,h_3\rangle\cong BS(1,-1)$ is solvable, a contradiction.

\item[(73)]$ R_1:h_2^2h_3^{-1}(h_2^{-1}h_3)^2h_3=1$, $R_2:h_2h_3h_2^{-1}h_3^2h_2^{-1}h_3^{-1}h_2=1$, $R_3:h_2(h_3h_2^{-1})^2h_3^{-1}h_2h_3=1$,\\$R_4:h_2h_3^{-1}h_2^{-1}h_3h_2^{-1}h_3^{-2}h_2=1$:\\
$\Rightarrow$  $G$ has a torsion element, a contradiction.

\item[(74)]$ R_1:h_2^2h_3^{-1}(h_2^{-1}h_3)^2h_3=1$, $R_2:h_2h_3h_2^{-1}h_3^2h_2^{-1}h_3^{-1}h_2=1$, $R_3:h_2h_3^{-2}h_2^{-1}h_3^3=1$,\\$R_4:h_2h_3^{-2}h_2^{-1}h_3^3=1$:\\
By interchanging $h_2$ and $h_3$ in (17) and with the same discussion, there is a contradiction.

\item[(75)]$ R_1:h_2^2h_3^{-2}(h_2^{-1}h_3)^2=1$, $R_2:(h_2h_3)^2h_2^{-1}h_3^{-1}h_2=1$, $R_3:(h_2h_3)^2h_2^{-1}h_3^{-1}h_2=1$,\\$R_4:h_2h_3^{-2}h_2h_3^3=1$:\\
By interchanging $h_2$ and $h_3$ in (3) and with the same discussion, there is a contradiction.

\item[(76)]$ R_1:h_2^2h_3^{-2}(h_2^{-1}h_3)^2=1$, $R_2:h_2h_3h_2^{-1}h_3^{-1}h_2h_3^{-1}h_2^{-1}h_3=1$, $R_3:(h_2h_3)^2h_2^{-1}h_3^{-1}h_2=1$,\\$R_4:h_2h_3^{-2}h_2h_3^3=1$:\\
By interchanging $h_2$ and $h_3$ in (4) and with the same discussion, there is a contradiction.

\item[(77)]$ R_1:h_2^2h_3^{-2}(h_2^{-1}h_3)^2=1$, $R_2:h_2(h_3h_2^{-1})^2h_3^3=1$, $R_3:h_2h_3^{-1}h_2^{-1}h_3^2h_2^{-1}h_3^{-1}h_2=1$,\\$R_4:h_2h_3^{-1}h_2^{-1}h_3^2h_2^{-1}h_3^{-1}h_2=1$:\\
By interchanging $h_2$ and $h_3$ in (13) and with the same discussion, there is a contradiction.

\item[(78)]$ R_1:h_2^2h_3^{-2}(h_2^{-1}h_3)^2=1$, $R_2:h_2h_3^{-1}h_2^{-1}h_3h_2^{-1}h_3^{-1}h_2h_3=1$, $R_3:h_2(h_3h_2^{-1})^2h_3^3=1$,\\$R_4:h_2h_3^{-1}h_2^{-1}h_3h_2^{-1}h_3^{-1}h_2h_3=1$:\\
By interchanging $h_2$ and $h_3$ in (14) and with the same discussion, there is a contradiction.

\item[(79)]$ R_1:h_2^2h_3^{-2}(h_2^{-1}h_3)^2=1$, $R_2:h_2h_3^{-1}h_2h_3h_2^{-1}h_3^3=1$, $R_3:h_2h_3^{-1}h_2^{-1}h_3^2h_2^{-1}h_3^{-1}h_2=1$,\\$R_4:h_2h_3^{-1}h_2^{-1}h_3^2h_2^{-1}h_3^{-1}h_2=1$:\\
By interchanging $h_2$ and $h_3$ in (37) and with the same discussion, there is a contradiction.

\item[(80)]$ R_1:h_2^2h_3^{-2}(h_2^{-1}h_3)^2=1$, $R_2:h_3^2(h_3h_2^{-1})^3h_3=1$, $R_3:h_2h_3^{-1}h_2^{-1}h_3^2h_2^{-1}h_3^{-1}h_2=1$,\\$R_4:h_2h_3^{-1}h_2^{-1}h_3^2h_2^{-1}h_3^{-1}h_2=1$:\\
By interchanging $h_2$ and $h_3$ in (9) and with the same discussion, there is a contradiction.

\item[(81)]$ R_1:h_2(h_2h_3^{-2})^2h_2=1$, $R_2:(h_2h_3)^2(h_2^{-1}h_3)^2=1$, $R_3:h_2^2h_3^{-1}h_2h_3h_2^{-1}h_3^2=1$,\\$R_4:h_2h_3(h_2h_3^{-1})^2h_3^{-1}h_2=1$:\\
$\Rightarrow\langle h_2,h_3\rangle\cong BS(1,-1)$ is solvable, a contradiction.

\item[(82)]$ R_1:h_2(h_2h_3^{-2})^2h_2=1$, $R_2:(h_2h_3)^2(h_2^{-1}h_3)^2=1$, $R_3:h_2h_3(h_2h_3^{-1})^2h_3^{-1}h_2=1$,\\$R_4:h_2h_3h_2^{-1}(h_2^{-1}h_3)^2h_3=1$:\\
$\Rightarrow\langle h_2,h_3\rangle\cong BS(1,-1)$ is solvable, a contradiction.

\item[(83)]$ R_1:(h_2h_3)^2h_2^{-2}h_3=1$, $R_2:h_2(h_3h_2^{-1})^2h_3^{-2}h_2=1$, $R_3:h_2(h_3h_2^{-1})^2h_3^{-2}h_2=1$,\\$R_4:h_2h_3^{-1}h_2(h_3h_2^{-1}h_3)^2=1$:\\
$\Rightarrow$  $G$ has a torsion element, a contradiction.

\item[(84)]$ R_1:(h_2h_3)^2h_2^{-2}h_3=1$, $R_2:h_2h_3^{-2}h_2^{-1}h_3^3=1$, $R_3:h_2h_3h_2h_3^{-1}h_2^{-2}h_3=1$,\\$R_4:h_2h_3^{-2}h_2^{-1}h_3h_2^{-1}h_3^{-1}h_2=1$:\\
By interchanging $h_2$ and $h_3$ in (23) and with the same discussion, there is a contradiction.

\item[(85)]$ R_1:(h_2h_3)^2h_2^{-2}h_3=1$, $R_2:h_2h_3^{-2}h_2^{-1}h_3^3=1$, $R_3:h_2h_3h_2^{-1}h_3^2h_2^{-1}h_3^{-1}h_2=1$,\\$R_4:h_2h_3^{-1}h_2h_3(h_3h_2^{-1})^2h_3=1$:\\
By interchanging $h_2$ and $h_3$ in (26) and with the same discussion, there is a contradiction.

\item[(86)]$ R_1:(h_2h_3)^2(h_2^{-1}h_3)^2=1$, $R_2:h_2h_3^{-1}h_2^{-1}h_3h_2^{-1}h_3^{-2}h_2=1$, $R_3:h_2(h_3h_2^{-1})^2h_3^{-2}h_2=1$,\\$R_4:h_2h_3^{-1}h_2(h_3h_2^{-1}h_3)^2=1$:\\
$\Rightarrow\langle h_2,h_3\rangle\cong BS(-3,1)$ is solvable, a contradiction.

\item[(87)]$ R_1:h_2h_3h_2h_3^{-1}h_2^{-1}h_3^{-1}h_2=1$, $R_2:h_2h_3h_2h_3^{-1}h_2^{-1}h_3^{-1}h_2=1$, $R_3:h_2h_3h_2^{-1}h_3^{-1}h_2h_3^{-1}h_2^{-1}h_3=1$,\\$R_4:h_2h_3^{-2}h_2h_3h_2h_3^{-1}h_2=1$:\\
$\Rightarrow \langle h_2,h_3\rangle=\langle h_2\rangle$ is abelian, a contradiction.

\item[(88)]$ R_1:h_2h_3h_2h_3^{-2}h_2^{-1}h_3=1$, $R_2:h_2h_3h_2^{-1}h_3(h_2h_3^{-1})^2h_2=1$, $R_3:h_2h_3h_2^{-1}h_3(h_2h_3^{-1})^2h_2=1$,\\$R_4:h_2h_3^{-2}h_2^{-1}h_3h_2^{-2}h_3=1$:\\
By interchanging $h_2$ and $h_3$ in (65) and with the same discussion, there is a contradiction.

\item[(89)]$ R_1:h_2h_3h_2h_3^{-2}h_2^{-1}h_3=1$, $R_2:h_2(h_3h_2^{-1})^2(h_3^{-1}h_2)^2=1$, $R_3:h_2h_3h_2^{-1}h_3(h_2h_3^{-1})^2h_2=1$,\\$R_4:h_2h_3^{-2}h_2^{-1}h_3h_2^{-2}h_3=1$:\\
By interchanging $h_2$ and $h_3$ in (66) and with the same discussion, there is a contradiction.

\item[(90)]$ R_1:h_2h_3h_2h_3^{-2}h_2h_3=1$, $R_2:h_2h_3h_2^{-1}h_3(h_2h_3^{-1})^2h_2=1$, $R_3:h_2h_3^2(h_2h_3^{-1})^2h_2=1$,\\$R_4:h_2h_3^{-2}h_2^{-1}h_3h_2^{-2}h_3=1$:\\
By interchanging $h_2$ and $h_3$ in (46) and with the same discussion, there is a contradiction.

\item[(91)]$ R_1:h_2h_3h_2h_3^{-2}h_2h_3=1$, $R_2:h_2h_3h_2^{-1}h_3(h_2h_3^{-1})^2h_2=1$, $R_3:h_2h_3^2(h_2h_3^{-1})^2h_2=1$,\\$R_4:h_2h_3^{-1}h_2(h_3h_2^{-1})^2h_3^{-1}h_2=1$:\\
By interchanging $h_2$ and $h_3$ in (55) and with the same discussion, there is a contradiction.

\item[(92)]$ R_1:h_2h_3h_2h_3^{-2}h_2h_3=1$, $R_2:h_2h_3h_2^{-1}h_3(h_2h_3^{-1})^2h_2=1$, $R_3:h_2h_3^2h_2^{-1}h_3^{-2}h_2=1$,\\$R_4:h_2h_3^{-1}(h_2^{-1}h_3)^2h_2^{-1}h_3^{-1}h_2=1$:\\
By interchanging $h_2$ and $h_3$ in (59) and with the same discussion, there is a contradiction.

\item[(93)]$ R_1:h_2h_3h_2h_3^{-2}h_2h_3=1$, $R_2:h_2(h_3h_2^{-1})^2h_3^{-2}h_2=1$, $R_3:h_2h_3^2h_2h_3^{-2}h_2=1$,\\$R_4:h_2h_3^{-1}(h_2^{-1}h_3)^2h_2h_3^{-1}h_2=1$:\\
By interchanging $h_2$ and $h_3$ in (28) and with the same discussion, there is a contradiction.

\item[(94)]$ R_1:h_2h_3h_2h_3^{-2}h_2h_3=1$, $R_2:h_2(h_3h_2^{-1})^2h_3^{-2}h_2=1$, $R_3:h_2h_3^2h_2^{-1}h_3^{-2}h_2=1$,\\$R_4:h_2h_3^{-1}(h_2^{-1}h_3)^2h_2h_3^{-1}h_2=1$:\\
By interchanging $h_2$ and $h_3$ in (63) and with the same discussion, there is a contradiction.

\item[(95)]$ R_1:h_2h_3h_2h_3^{-2}h_2h_3=1$, $R_2:h_2(h_3h_2^{-1})^2h_3^{-2}h_2=1$, $R_3:h_2(h_3h_2^{-1})^2h_3^{-2}h_2=1$,\\$R_4:(h_2h_3^{-1}h_2)^2h_3h_2^{-1}h_3=1$:\\
By interchanging $h_2$ and $h_3$ in (83) and with the same discussion, there is a contradiction.

\item[(96)]$ R_1:h_2h_3h_2h_3^{-2}h_2h_3=1$, $R_2:h_2(h_3h_2^{-1})^2(h_3^{-1}h_2)^2=1$, $R_3:h_2h_3^2(h_2h_3^{-1})^2h_2=1$,\\$R_4:h_2h_3^{-2}h_2^{-1}h_3h_2^{-2}h_3=1$:\\
By interchanging $h_2$ and $h_3$ in (47) and with the same discussion, there is a contradiction.

\item[(97)]$ R_1:h_2h_3h_2h_3^{-2}h_2h_3=1$, $R_2:h_2(h_3h_2^{-1})^2(h_3^{-1}h_2)^2=1$, $R_3:h_2h_3^2(h_2h_3^{-1})^2h_2=1$,\\$R_4:h_2h_3^{-1}h_2(h_3h_2^{-1})^2h_3^{-1}h_2=1$:\\
By interchanging $h_2$ and $h_3$ in (56) and with the same discussion, there is a contradiction.

\item[(98)]$ R_1:h_2h_3h_2h_3^{-2}h_2h_3=1$, $R_2:h_2(h_3h_2^{-1})^2(h_3^{-1}h_2)^2=1$, $R_3:h_2h_3^2h_2^{-1}h_3^{-2}h_2=1$,\\$R_4:h_2h_3^{-1}(h_2^{-1}h_3)^2h_2^{-1}h_3^{-1}h_2=1$:\\
By interchanging $h_2$ and $h_3$ in (60) and with the same discussion, there is a contradiction.

\item[(99)]$ R_1:h_2h_3h_2h_3^{-1}(h_3^{-1}h_2)^2=1$, $R_2:h_2(h_3h_2^{-1})^2h_2^{-1}h_3^2=1$, $R_3:h_2h_3(h_2h_3^{-1})^2h_2h_3=1$,\\$R_4:h_2h_3^{-3}h_2^2h_3^{-1}h_2=1$:\\
By interchanging $h_2$ and $h_3$ in (82) and with the same discussion, there is a contradiction.

\item[(100)]$ R_1:h_2h_3(h_2h_3^{-1})^2h_2^{-1}h_3=1$, $R_2:h_2h_3h_2^{-1}h_3h_2h_3^{-2}h_2=1$, $R_3:h_2h_3^2h_2^{-1}(h_3^{-1}h_2)^2=1$,\\$R_4:h_2h_3^{-1}h_2^{-1}h_3^2h_2h_3^{-1}h_2=1$:\\
By interchanging $h_2$ and $h_3$ in (73) and with the same discussion, there is a contradiction.

\item[(101)]$ R_1:h_2h_3(h_2h_3^{-1})^2h_2^{-1}h_3=1$, $R_2:h_2h_3h_2^{-1}h_3h_2h_3^{-2}h_2=1$, $R_3:h_2h_3h_2^{-1}h_3h_2h_3^{-2}h_2=1$,\\$R_4:h_2h_3^{-1}h_2^{-1}h_3^2h_2h_3^{-1}h_2=1$:\\
$\Rightarrow$  $G$ has a torsion element, a contradiction.

\item[(102)]$ R_1:h_2h_3(h_2h_3^{-1})^2h_2^{-1}h_3=1$, $R_2:h_2h_3h_2^{-1}h_3h_2h_3^{-2}h_2=1$, $R_3:h_2h_3h_2^{-1}h_3h_2h_3^{-2}h_2=1$,\\$R_4:h_2h_3^{-1}(h_2^{-1}h_3)^2h_2h_3^{-1}h_2=1$:\\
$\Rightarrow\langle h_2,h_3\rangle\cong BS(1,1)$ is solvable, a contradiction.

\item[(103)]$ R_1:h_2h_3(h_2h_3^{-1})^2h_2^{-1}h_3=1$, $R_2:h_2h_3h_2^{-1}h_3h_2h_3^{-2}h_2=1$, $R_3:h_2h_3h_2^{-1}h_3(h_2h_3^{-1})^2h_2=1$,\\$R_4:(h_2h_3^{-1}h_2)^2h_3^2=1$:\\
By interchanging $h_2$ and $h_3$ in (42) and with the same discussion, there is a contradiction.

\item[(104)]$ R_1:h_2h_3(h_2h_3^{-1})^2h_2^{-1}h_3=1$, $R_2:h_2h_3h_2^{-1}h_3h_2h_3^{-2}h_2=1$, $R_3:h_2(h_3h_2^{-1})^2(h_3^{-1}h_2)^2=1$,\\$R_4:(h_2h_3^{-1}h_2)^2h_3h_2^{-1}h_3=1$:\\
$\Rightarrow\langle h_2,h_3\rangle\cong BS(1,-2)$ is solvable, a contradiction.

\item[(105)]$ R_1:h_2h_3(h_2h_3^{-1})^2h_2^{-1}h_3=1$, $R_2:h_2h_3h_2^{-1}h_3(h_2h_3^{-1})^2h_2=1$, $R_3:h_2h_3h_2^{-1}h_3(h_2h_3^{-1})^2h_2=1$,\\$R_4:h_2h_3^{-2}h_2^{-1}h_3h_2^{-2}h_3=1$:\\
By interchanging $h_2$ and $h_3$ in (71) and with the same discussion, there is a contradiction.

\item[(106)]$ R_1:h_2h_3(h_2h_3^{-1})^2h_3^{-1}h_2=1$, $R_2:h_2h_3h_2^{-1}h_3^{-1}h_2h_3^{-1}h_2^{-1}h_3=1$, $R_3:h_2h_3h_2^{-1}h_3^{-1}h_2h_3^{-1}h_2^{-1}h_3=1$, $R_4:h_2h_3^{-2}h_2h_3h_2h_3^{-1}h_2=1$:\\
$\Rightarrow$ $G$ is solvable, a contradiction.

\item[(107)]$ R_1:h_2h_3(h_2h_3^{-1})^2h_2h_3=1$, $R_2:h_2h_3h_2^{-1}h_3h_2h_3^{-2}h_2=1$, $R_3:h_2h_3^2h_2h_3^{-2}h_2=1$,\\$R_4:h_2h_3^{-1}(h_2^{-1}h_3)^2h_2h_3^{-1}h_2=1$:\\
By interchanging $h_2$ and $h_3$ in (29) and with the same discussion, there is a contradiction.

\item[(108)]$ R_1:h_2h_3(h_2h_3^{-1})^2h_2h_3=1$, $R_2:h_2h_3h_2^{-1}h_3h_2h_3^{-2}h_2=1$, $R_3:h_2h_3^2h_2^{-1}h_3^{-2}h_2=1$,\\$R_4:h_2h_3^{-1}(h_2^{-1}h_3)^2h_2h_3^{-1}h_2=1$:\\
By interchanging $h_2$ and $h_3$ in (64) and with the same discussion, there is a contradiction.

\item[(109)]$ R_1:h_2h_3(h_2h_3^{-1})^2h_2h_3=1$, $R_2:h_2h_3h_2^{-1}h_3h_2h_3^{-2}h_2=1$, $R_3:h_2(h_3h_2^{-1})^2h_3^{-2}h_2=1$,\\$R_4:(h_2h_3^{-1}h_2)^2h_3h_2^{-1}h_3=1$:\\
By interchanging $h_2$ and $h_3$ in (86) and with the same discussion, there is a contradiction.

\item[(110)]$ R_1:h_2h_3(h_2h_3^{-1})^2h_2h_3=1$, $R_2:h_2(h_3h_2^{-1})^2(h_3^{-1}h_2)^2=1$, $R_3:h_2h_3^2(h_2h_3^{-1})^2h_2=1$,\\$R_4:h_2h_3^{-2}h_2^{-1}h_3h_2^{-2}h_3=1$:\\
By interchanging $h_2$ and $h_3$ in (48) and with the same discussion, there is a contradiction.

\item[(111)]$ R_1:h_2h_3(h_2h_3^{-1})^2h_2h_3=1$, $R_2:h_2(h_3h_2^{-1})^2(h_3^{-1}h_2)^2=1$, $R_3:h_2h_3^2(h_2h_3^{-1})^2h_2=1$,\\$R_4:h_2h_3^{-1}h_2(h_3h_2^{-1})^2h_3^{-1}h_2=1$:\\
By interchanging $h_2$ and $h_3$ in (57) and with the same discussion, there is a contradiction.

\item[(112)]$ R_1:h_2h_3(h_2h_3^{-1})^2h_2h_3=1$, $R_2:h_2(h_3h_2^{-1})^2(h_3^{-1}h_2)^2=1$, $R_3:h_2h_3^2h_2^{-1}h_3^{-2}h_2=1$,\\$R_4:h_2h_3^{-1}(h_2^{-1}h_3)^2h_2^{-1}h_3^{-1}h_2=1$:\\
By interchanging $h_2$ and $h_3$ in (61) and with the same discussion, there is a contradiction.

\item[(113)]$ R_1:h_2h_3(h_2h_3^{-1})^2h_2h_3=1$, $R_2:h_2h_3^{-3}h_2^2h_3^{-1}h_2=1$, $R_3:h_2h_3^2h_2^{-1}(h_3^{-1}h_2)^2=1$,\\$R_4:h_2(h_3h_2^{-1})^2h_2^{-1}h_3^2=1$:\\
By interchanging $h_2$ and $h_3$ in (72) and with the same discussion, there is a contradiction.

\item[(114)]$ R_1:h_2h_3(h_2h_3^{-1})^2h_2h_3=1$, $R_2:h_2h_3^{-3}h_2^2h_3^{-1}h_2=1$, $R_3:h_2h_3^2h_2^{-1}h_3h_2h_3^{-1}h_2=1$,\\$R_4:h_2(h_3h_2^{-1})^2h_2^{-1}h_3^2=1$:\\
By interchanging $h_2$ and $h_3$ in (81) and with the same discussion, there is a contradiction.

\item[(115)]$ R_1:h_2h_3^2(h_2h_3^{-1})^2h_2=1$, $R_2:h_2h_3h_2^{-1}h_3^{-2}h_2h_3=1$, $R_3:h_2h_3^2(h_2h_3^{-1})^2h_2=1$,\\$R_4:h_2h_3^{-2}h_2^{-1}h_3h_2^{-2}h_3=1$:\\
By interchanging $h_2$ and $h_3$ in (44) and with the same discussion, there is a contradiction.

\item[(116)]$ R_1:h_2h_3^2(h_2h_3^{-1})^2h_2=1$, $R_2:h_2h_3h_2^{-1}h_3^{-2}h_2h_3=1$, $R_3:h_2h_3^2(h_2h_3^{-1})^2h_2=1$,\\$R_4:h_2h_3^{-1}h_2(h_3h_2^{-1})^2h_3^{-1}h_2=1$:\\
By interchanging $h_2$ and $h_3$ in (50) and with the same discussion, there is a contradiction.

\item[(117)]$ R_1:h_2h_3^2(h_2h_3^{-1})^2h_2=1$, $R_2:h_2h_3h_2^{-1}h_3^{-2}h_2h_3=1$, $R_3:h_2h_3^2h_2^{-1}h_3^{-2}h_2=1$,\\$R_4:h_2h_3^{-1}(h_2^{-1}h_3)^2h_2^{-1}h_3^{-1}h_2=1$:\\
By interchanging $h_2$ and $h_3$ in (51) and with the same discussion, there is a contradiction.

\item[(118)]$ R_1:h_2h_3^2(h_2h_3^{-1})^2h_2=1$, $R_2:h_2h_3h_2^{-1}h_3^{-1}h_2h_3^{-1}h_2^{-1}h_3=1$, $R_3:h_2h_3h_2^{-1}h_3(h_2h_3^{-1})^2h_2=1$,\\$R_4:h_2h_3^{-2}h_2^{-1}h_3h_2^{-2}h_3=1$:\\
By interchanging $h_2$ and $h_3$ in (54) and with the same discussion, there is a contradiction.

\item[(119)]$ R_1:h_2h_3^2(h_2h_3^{-1})^2h_2=1$, $R_2:h_2h_3h_2^{-1}(h_3^{-1}h_2)^2h_3=1$, $R_3:h_2h_3^2(h_2h_3^{-1})^2h_2=1$,\\$R_4:h_2h_3^{-2}h_2^{-1}h_3h_2^{-2}h_3=1$:\\
By interchanging $h_2$ and $h_3$ in (45) and with the same discussion, there is a contradiction.

\item[(120)]$ R_1:h_2h_3^2(h_2h_3^{-1})^2h_2=1$, $R_2:h_2h_3h_2^{-1}(h_3^{-1}h_2)^2h_3=1$, $R_3:h_2h_3^2(h_2h_3^{-1})^2h_2=1$,\\$R_4:h_2h_3^{-1}h_2(h_3h_2^{-1})^2h_3^{-1}h_2=1$:\\
By interchanging $h_2$ and $h_3$ in (52) and with the same discussion, there is a contradiction.

\item[(121)]$ R_1:h_2h_3^2(h_2h_3^{-1})^2h_2=1$, $R_2:h_2h_3h_2^{-1}(h_3^{-1}h_2)^2h_3=1$, $R_3:h_2h_3^2h_2^{-1}h_3^{-2}h_2=1$,\\$R_4:h_2h_3^{-1}(h_2^{-1}h_3)^2h_2^{-1}h_3^{-1}h_2=1$:\\
By interchanging $h_2$ and $h_3$ in (53) and with the same discussion, there is a contradiction.

\item[(122)]$ R_1:h_2h_3^2h_2^{-1}(h_2^{-1}h_3)^2=1$, $R_2:h_2h_3h_2^{-1}h_3^{-3}h_2=1$, $R_3:h_2h_3h_2^{-1}h_3^{-3}h_2=1$,\\$R_4:h_2h_3h_2^{-1}h_3^2h_2^{-2}h_3=1$:\\
By interchanging $h_2$ and $h_3$ in (34) and with the same discussion, there is a contradiction.

\item[(123)]$ R_1:h_2h_3^2h_2^{-1}h_3^{-1}h_2^{-1}h_3=1$, $R_2:h_2h_3h_2^{-1}h_3^{-3}h_2=1$, $R_3:h_2h_3h_2^{-1}h_3^{-3}h_2=1$,\\$R_4:h_2h_3h_2^{-1}h_3^2h_2^{-2}h_3=1$:\\
By interchanging $h_2$ and $h_3$ in (33) and with the same discussion, there is a contradiction.

\item[(124)]$ R_1:h_2h_3^2h_2^{-1}h_3^{-1}h_2^{-1}h_3=1$, $R_2:h_2(h_3h_2^{-1})^2h_2^{-1}h_3^2=1$, $R_3:h_2h_3h_2^{-1}h_3^2h_2^{-2}h_3=1$,\\$R_4:h_2h_3^{-1}h_2^{-1}h_3h_2^{-1}h_3^{-1}h_2h_3=1$:\\
By interchanging $h_2$ and $h_3$ in (67) and with the same discussion, there is a contradiction.

\item[(125)]$ R_1:h_2h_3^2h_2^{-1}(h_3^{-1}h_2)^2=1$, $R_2:h_2h_3h_2^{-1}h_3^{-2}h_2^{-1}h_3=1$, $R_3:h_2h_3h_2^{-1}h_3(h_2h_3^{-1})^2h_2=1$,\\$R_4:h_2h_3^{-2}h_2^{-1}h_3h_2^{-2}h_3=1$:\\
By interchanging $h_2$ and $h_3$ in (30) and with the same discussion, there is a contradiction.

\item[(126)]$ R_1:h_2h_3(h_3h_2^{-1})^2h_2^{-1}h_3=1$, $R_2:h_2h_3h_2^{-1}h_3^{-3}h_2=1$, $R_3:h_2h_3h_2^{-1}h_3^{-3}h_2=1$,\\$R_4:h_2h_3h_2^{-1}h_3^2h_2^{-2}h_3=1$:\\
By interchanging $h_2$ and $h_3$ in (35) and with the same discussion, there is a contradiction.

\item[(127)]$ R_1:h_2h_3h_2^{-1}h_3^{-1}h_2^{-1}h_3^2=1$, $R_2:h_2h_3h_2^{-1}h_3^{-1}h_2^{-1}h_3^2=1$, $R_3:h_2h_3h_2^{-1}h_3^2h_2^{-2}h_3=1$,\\$R_4:h_2h_3^{-1}h_2^{-1}h_3h_2^{-1}h_3^{-1}h_2h_3=1$:\\
By interchanging $h_2$ and $h_3$ in (87) and with the same discussion, there is a contradiction.

\item[(128)]$ R_1:h_2h_3h_2^{-1}h_3^2h_2^{-2}h_3=1$, $R_2:h_2h_3^{-1}h_2^{-1}h_3h_2^{-1}h_3^{-1}h_2h_3=1$, $R_3:h_2(h_3h_2^{-1})^2h_2^{-1}h_3^2=1$,\\$R_4:h_2h_3^{-1}h_2^{-1}h_3h_2^{-1}h_3^{-1}h_2h_3=1$:\\
By interchanging $h_2$ and $h_3$ in (106) and with the same discussion, there is a contradiction.

\item[(129)]$ R_1:h_2h_3h_2^{-1}h_3^2h_2^{-1}h_3^{-1}h_2=1$, $R_2:h_2h_3^{-1}h_2^{-1}h_3^2h_2h_3^{-1}h_2=1$, $R_3:h_2h_3^{-3}h_2^2h_3^{-1}h_2=1$,\\$R_4:(h_2h_3^{-1}h_2)^2h_3^2=1$:\\
By interchanging $h_2$ and $h_3$ in (38) and with the same discussion, there is a contradiction.

\item[(130)]$ R_1:h_2h_3h_2^{-1}h_3^2h_2^{-1}h_3^{-1}h_2=1$, $R_2:h_2h_3^{-1}h_2^{-1}h_3h_2^{-1}h_3^{-2}h_2=1$, $R_3:h_2(h_3h_2^{-1})^2h_3^{-1}h_2h_3=1$,\\$R_4:h_2h_3^{-1}h_2^{-1}h_3h_2^{-1}h_3^{-2}h_2=1$:\\
By interchanging $h_2$ and $h_3$ in (101) and with the same discussion, there is a contradiction.

\item[(131)]$ R_1:h_2h_3h_2^{-1}h_3^2h_2^{-1}h_3^{-1}h_2=1$, $R_2:h_2h_3^{-1}h_2^{-1}h_3h_2^{-1}h_3^{-2}h_2=1$, $R_3:h_2h_3^{-2}h_2^{-1}h_3^3=1$,\\$R_4:h_2h_3^{-2}h_2^{-1}h_3^3=1$:\\
By interchanging $h_2$ and $h_3$ in (18) and with the same discussion, there is a contradiction.

\item[(132)]$ R_1:h_2(h_3h_2^{-1})^2h_3^{-1}h_2h_3=1$, $R_2:h_2h_3^{-1}h_2^{-1}h_3h_2^{-1}h_3^{-2}h_2=1$, $R_3:h_2h_3^{-1}h_2^{-1}h_3h_2^{-1}h_3^{-2}h_2=1$,\\$R_4:h_2h_3^{-2}h_2h_3^{-1}(h_2^{-1}h_3)^2=1$:\\
By interchanging $h_2$ and $h_3$ in (102) and with the same discussion, there is a contradiction.

\item[(133)]$ R_1:h_2(h_3h_2^{-1})^2h_3^{-1}h_2h_3=1$, $R_2:h_2h_3^{-1}h_2^{-1}h_3h_2^{-1}h_3^{-2}h_2=1$, $R_3:h_2h_3^{-1}h_2(h_3h_2^{-1}h_3)^2=1$,\\$R_4:(h_2h_3^{-1})^2(h_2^{-1}h_3)^2h_3=1$:\\
By interchanging $h_2$ and $h_3$ in (104) and with the same discussion, there is a contradiction.

\item[(134)]$ R_1:h_2h_3^{-1}(h_2^{-1}h_3^2)^2=1$, $R_2:h_2h_3^{-1}h_2(h_3h_2^{-1})^3h_3=1$, $R_3:h_2h_3^{-1}(h_2^{-1}h_3)^2h_3^2=1$,\\$R_4:(h_2h_3^{-2})^2h_2^{-1}h_3=1$:\\
By interchanging $h_2$ and $h_3$ in (36) and with the same discussion, there is a contradiction.

\item[(135)]$ R_1:h_2h_3^{-1}h_2^{-1}h_3h_2^{-1}h_3^{-2}h_2=1$, $R_2:h_2h_3^{-2}h_2h_3^{-1}(h_2^{-1}h_3)^2=1$, $R_3:h_2h_3^{-2}h_2^{-1}h_3^3=1$,\\$R_4:h_2h_3^{-2}h_2^{-1}h_3^3=1$:\\
By interchanging $h_2$ and $h_3$ in (19) and with the same discussion, there is a contradiction.

\item[(136)]$ R_1:h_2h_3^{-2}h_2^{-1}h_3^3=1$, $R_2:h_2h_3^{-2}h_2^{-1}h_3^3=1$, $R_3:h_2h_3^{-1}h_2(h_3h_2^{-1}h_3)^2=1$,\\$R_4:(h_2h_3^{-1})^2(h_2^{-1}h_3)^2h_3=1$:\\
By interchanging $h_2$ and $h_3$ in (21) and with the same discussion, there is a contradiction.

\end{enumerate}

\subsection{$\mathbf{C_6---C_6(C_6--C_6)}$}
\begin{figure}[ht]
\psscalebox{0.9 0.9} 
{
\begin{pspicture}(0,-2.0985577)(5.194231,2.0985577)
\psdots[linecolor=black, dotsize=0.4](2.5971153,1.9014423)
\psdots[linecolor=black, dotsize=0.4](2.5971153,-0.49855766)
\psdots[linecolor=black, dotsize=0.4](1.3971153,-1.2985576)
\psdots[linecolor=black, dotsize=0.4](3.7971153,-1.2985576)
\psdots[linecolor=black, dotsize=0.4](4.997115,-0.49855766)
\psdots[linecolor=black, dotsize=0.4](2.5971153,1.1014423)
\psdots[linecolor=black, dotsize=0.4](2.5971153,0.30144233)
\psdots[linecolor=black, dotsize=0.4](3.7971153,1.1014423)
\psdots[linecolor=black, dotsize=0.4](3.7971153,0.30144233)
\psdots[linecolor=black, dotsize=0.4](1.3971153,0.30144233)
\psdots[linecolor=black, dotsize=0.4](1.3971153,1.1014423)
\psline[linecolor=black, linewidth=0.04](2.5971153,1.9014423)(3.7971153,1.1014423)(3.7971153,0.30144233)(2.5971153,-0.49855766)(2.5971153,0.30144233)(2.5971153,1.1014423)(2.5971153,1.9014423)(1.3971153,1.1014423)(1.3971153,0.30144233)(2.5971153,-0.49855766)
\psdots[linecolor=black, dotsize=0.4](0.19711533,0.30144233)
\psdots[linecolor=black, dotsize=0.4](0.19711533,-0.89855766)
\psline[linecolor=black, linewidth=0.04](1.3971153,1.1014423)(0.19711533,0.30144233)(0.19711533,-0.89855766)(1.3971153,-1.2985576)(2.5971153,1.1014423)
\psline[linecolor=black, linewidth=0.04](1.3971153,-1.2985576)(3.7971153,-1.2985576)(4.997115,-0.49855766)(3.7971153,1.1014423)
\rput[bl](0.1,-2.0985577){$\mathbf{25) \ C_6---C_6(C_6--C_6)}$}
\end{pspicture}
}
\end{figure}
By considering the $996$ cases related to the existence of $C_6--C_6$ and $99$ cases related to the existence of $C_6---C_6$ in the graph $K(\alpha,\beta)$, it can be seen that there are $1594$ cases for the relations of four $C_6$ cycles in this structure. Using GAP \cite{gap}, we see that all groups with two generators $h_2$ and $h_3$ and four relations which are between $1446$ cases of these $1594$ cases are finite and solvable, or just finite, that is a contradiction with the assumptions. So, there are just $148$ cases for the relations of these cycles which may lead to the existence of a subgraph isomorphic to the graph $C_6---C_6(C_6--C_6)$ in $K(\alpha,\beta)$. In the following, we show that these $148$ cases lead to contradictions and so, the graph $K(\alpha,\beta)$ contains no subgraph isomorphic to the graph $C_6---C_6(C_6--C_6)$.
\begin{enumerate}
\item[(1)]$ R_1:h_2^3h_3^{-1}h_2^{-1}h_3^2=1$, $R_2:h_2^2h_3^{-2}(h_2^{-1}h_3)^2=1$, $R_3:h_2^3h_3^{-1}h_2^{-1}h_3^2=1$,\\$R_4:h_2^2h_3^2h_2^{-1}h_3^{-1}h_2=1$:\\
$\Rightarrow\langle h_2,h_3\rangle\cong BS(1,1)$ is solvable, a contradiction.

\item[(2)]$ R_1:h_2^2h_3h_2^{-1}h_3^{-2}h_2=1$, $R_2:h_2^2h_3^{-1}h_2^{-1}h_3^{-1}h_2h_3=1$, $R_3:h_2^2h_3^{-1}h_2^{-1}h_3^{-1}h_2h_3=1$,\\$R_4:h_2h_3^{-2}(h_2^{-1}h_3)^3=1$:\\
$\Rightarrow\langle h_2,h_3\rangle\cong BS(1,-2)$ is solvable, a contradiction.

\item[(3)]$ R_1:h_2^2h_3h_2^{-1}h_3^{-2}h_2=1$, $R_2:h_2^2h_3^{-1}h_2^{-1}h_3^{-1}h_2h_3=1$, $R_3:h_2^2h_3^{-1}h_2^{-1}h_3^{-1}h_2h_3=1$,\\$R_4:h_2h_3^{-1}(h_3^{-1}h_2)^2h_3h_2^{-1}h_3=1$:\\
$\Rightarrow \langle h_2,h_3\rangle=\langle h_2\rangle$ is abelian, a contradiction.

\item[(4)]$ R_1:h_2^2h_3h_2^{-1}h_3^{-2}h_2=1$, $R_2:h_2^2h_3^{-1}h_2^{-1}h_3^{-1}h_2h_3=1$, $R_3:h_2^2h_3^{-2}(h_2^{-1}h_3)^2=1$,\\$R_4:h_2h_3h_2^{-1}h_3^{-1}h_2h_3^{-1}h_2^{-1}h_3=1$:\\
$\Rightarrow \langle h_2,h_3\rangle=\langle h_2\rangle$ is abelian, a contradiction.

\item[(5)]$ R_1:h_2^2h_3h_2^{-1}h_3^{-2}h_2=1$, $R_2:h_2^2h_3^{-1}h_2^{-1}h_3^{-1}h_2h_3=1$, $R_3:h_2^2h_3^{-1}(h_3^{-1}h_2)^2h_3=1$,\\$R_4:h_2h_3h_2^{-1}h_3^{-1}h_2h_3^{-1}h_2^{-1}h_3=1$:\\
$\Rightarrow\langle h_2,h_3\rangle\cong BS(2,1)$ is solvable, a contradiction.

\item[(6)]$ R_1:h_2^2h_3h_2^{-1}h_3^{-2}h_2=1$, $R_2:h_2h_3h_2h_3^{-1}(h_3^{-1}h_2)^2=1$, $R_3:h_2^2h_3^{-1}h_2^{-1}h_3^{-1}h_2h_3=1$,\\$R_4:h_2h_3h_2h_3^{-1}(h_3^{-1}h_2)^2=1$:\\
$\Rightarrow\langle h_2,h_3\rangle\cong BS(1,1)$ is solvable, a contradiction.

\item[(7)]$ R_1:h_2^2h_3h_2^{-1}h_3^{-2}h_2=1$, $R_2:h_2h_3h_2h_3^{-1}(h_3^{-1}h_2)^2=1$, $R_3:h_2^2h_3^{-1}h_2^{-1}h_3^{-1}h_2h_3=1$,\\$R_4:h_2h_3^{-1}(h_3^{-1}h_2)^2h_3h_2^{-1}h_3=1$:\\
$\Rightarrow\langle h_2,h_3\rangle\cong BS(1,1)$ is solvable, a contradiction.

\item[(8)]$ R_1:h_2^2h_3h_2^{-1}h_3^{-2}h_2=1$, $R_2:h_2h_3h_2h_3^{-1}(h_3^{-1}h_2)^2=1$, $R_3:h_2^2h_3^{-1}(h_3^{-1}h_2)^2h_3=1$,\\$R_4:h_2h_3h_2^{-1}h_3^{-1}h_2h_3^{-1}h_2^{-1}h_3=1$:\\
$\Rightarrow\langle h_2,h_3\rangle\cong BS(2,1)$ is solvable, a contradiction.

\item[(9)]$ R_1:h_2^2h_3^{-1}h_2^{-1}(h_3^{-1}h_2)^2=1$, $R_2:h_2^2h_3^{-2}h_2h_3^{-1}h_2^{-1}h_3=1$, $R_3:h_2^2h_3^{-1}h_2^{-1}h_3h_2^{-2}h_3=1$,\\$R_4:(h_2h_3^{-1})^2(h_2^{-1}h_3)^3=1$:\\
$\Rightarrow\langle h_2,h_3\rangle\cong BS(1,-1)$ is solvable, a contradiction.

\item[(10)]$ R_1:h_2^2h_3^{-1}h_2^{-1}(h_3^{-1}h_2)^2=1$, $R_2:h_2^2h_3^{-2}h_2h_3^{-1}h_2^{-1}h_3=1$, $R_3:h_2(h_3h_2^{-2})^2h_3=1$,\\$R_4:(h_2h_3^{-1})^2(h_2^{-1}h_3)^3=1$:\\
$\Rightarrow\langle h_2,h_3\rangle\cong BS(1,-5)$ is solvable, a contradiction.

\item[(11)]$ R_1:h_2^2h_3^{-1}h_2^{-1}(h_3^{-1}h_2)^2=1$, $R_2:h_2^2h_3^{-2}h_2h_3^{-1}h_2^{-1}h_3=1$, $R_3:h_2(h_3h_2^{-2})^2h_3=1$,\\$R_4:(h_2h_3^{-1})^3(h_2^{-1}h_3)^2=1$:\\
$\Rightarrow\langle h_2,h_3\rangle\cong BS(1,-5)$ is solvable, a contradiction.

\item[(12)]$ R_1:h_2^2h_3^{-1}h_2^{-1}(h_3^{-1}h_2)^2=1$, $R_2:h_2^2h_3^{-2}h_2h_3^{-1}h_2^{-1}h_3=1$, $R_3:h_2h_3h_2^{-1}h_3^{-1}h_2h_3^{-2}h_2=1$,\\$R_4:h_2h_3^{-1}h_2^{-1}h_3^2h_2^{-2}h_3=1$:\\
$\Rightarrow\langle h_2,h_3\rangle\cong BS(1,2)$ is solvable, a contradiction.

\item[(13)]$ R_1:h_2^2h_3^{-1}h_2^{-1}(h_3^{-1}h_2)^2=1$, $R_2:h_2^2h_3^{-2}h_2h_3^{-1}h_2^{-1}h_3=1$, $R_3:h_2h_3^{-1}h_2^{-1}h_3h_2h_3^{-2}h_2=1$,\\$R_4:h_2h_3^{-2}h_2^{-1}h_3^2h_2^{-1}h_3=1$:\\
$\Rightarrow \langle h_2,h_3\rangle=\langle h_2\rangle$ is abelian, a contradiction.

\item[(14)]$ R_1:h_2^2h_3^{-1}h_2^{-1}(h_3^{-1}h_2)^2=1$, $R_2:h_2^2h_3^{-2}h_2h_3^{-1}h_2^{-1}h_3=1$, $R_3:h_2h_3^{-1}h_2^{-1}h_3^2h_2^{-1}h_3^{-1}h_2=1$,\\$R_4:h_2h_3^{-2}h_2^{-1}h_3^2h_2^{-1}h_3=1$:\\
$\Rightarrow$  $G$ has a torsion element, a contradiction.

\item[(15)]$ R_1:h_2^2h_3^{-1}h_2^{-1}(h_3^{-1}h_2)^2=1$, $R_2:h_2^2h_3^{-2}h_2h_3^{-1}h_2^{-1}h_3=1$, $R_3:h_2h_3^{-1}h_2^{-1}h_3^2h_2^{-1}h_3^{-1}h_2=1$,\\$R_4:h_2h_3^{-3}h_2h_3h_2^{-1}h_3=1$:\\
$\Rightarrow$  $G$ has a torsion element, a contradiction.

\item[(16)]$ R_1:h_2^2h_3^{-1}h_2^{-1}(h_3^{-1}h_2)^2=1$, $R_2:h_2^2h_3^{-2}h_2h_3^{-1}h_2^{-1}h_3=1$, $R_3:h_2h_3^{-1}h_2^{-1}h_3h_2^{-2}h_3^2=1$,\\$R_4:(h_2h_3^{-1})^2(h_2^{-1}h_3)^3=1$:\\
$\Rightarrow\langle h_2,h_3\rangle\cong BS(1,-5)$ is solvable, a contradiction.

\item[(17)]$ R_1:h_2^2h_3^{-1}h_2^{-1}(h_3^{-1}h_2)^2=1$, $R_2:h_2^2h_3^{-2}h_2h_3^{-1}h_2^{-1}h_3=1$, $R_3:h_2h_3^{-1}h_2^{-1}h_3h_2^{-2}h_3^2=1$,\\$R_4:(h_2h_3^{-1})^3(h_2^{-1}h_3)^2=1$:\\
$\Rightarrow\langle h_2,h_3\rangle\cong BS(1,-5)$ is solvable, a contradiction.

\item[(18)]$ R_1:h_2^2h_3^{-1}h_2^{-1}(h_3^{-1}h_2)^2=1$, $R_2:h_2^2h_3^{-2}h_2h_3^{-1}h_2^{-1}h_3=1$, $R_3:h_2h_3^{-3}h_2^2h_3^{-1}h_2=1$,\\$R_4:(h_2h_3^{-1})^2(h_2^{-1}h_3)^3=1$:\\
$\Rightarrow\langle h_2,h_3\rangle\cong BS(1,-5)$ is solvable, a contradiction.

\item[(19)]$ R_1:h_2^2h_3^{-1}h_2^{-1}(h_3^{-1}h_2)^2=1$, $R_2:h_2^2h_3^{-2}h_2h_3^{-1}h_2^{-1}h_3=1$, $R_3:h_2h_3^{-3}h_2^2h_3^{-1}h_2=1$,\\$R_4:(h_2h_3^{-1})^3(h_2^{-1}h_3)^2=1$:\\
$\Rightarrow\langle h_2,h_3\rangle\cong BS(1,-5)$ is solvable, a contradiction.

\item[(20)]$ R_1:h_2^2h_3^{-1}h_2^{-1}(h_3^{-1}h_2)^2=1$, $R_2:h_2h_3h_2^{-1}h_3^{-1}h_2h_3^{-2}h_2=1$, $R_3:h_2^2h_3^{-1}h_2^{-1}h_3h_2^{-1}h_3^{-1}h_2=1$,\\$R_4:h_2(h_3^{-1}h_2h_3^{-1})^2h_2h_3^{-1}h_2=1$:\\
$\Rightarrow\langle h_2,h_3\rangle\cong BS(1,5)$ is solvable, a contradiction.

\item[(21)]$ R_1:h_2^2h_3^{-1}h_2^{-1}(h_3^{-1}h_2)^2=1$, $R_2:h_2h_3h_2^{-1}h_3^{-1}h_2h_3^{-2}h_2=1$, $R_3:h_2^2h_3^{-1}h_2^{-1}h_3h_2^{-1}h_3^{-1}h_2=1$,\\$R_4:h_2h_3^{-1}h_2(h_3^{-1}h_2h_3^{-1})^2h_2=1$:\\
$\Rightarrow\langle h_2,h_3\rangle\cong BS(3,1)$ is solvable, a contradiction.

\item[(22)]$ R_1:h_2^2h_3^{-1}h_2^{-1}(h_3^{-1}h_2)^2=1$, $R_2:h_2h_3h_2^{-1}h_3^{-1}h_2h_3^{-2}h_2=1$, $R_3:h_2(h_3h_2^{-2})^2h_3=1$,\\$R_4:h_2h_3^{-1}h_2^{-1}h_3^2h_2^{-1}h_3^{-1}h_2=1$:\\
$\Rightarrow\langle h_2,h_3\rangle\cong BS(1,-5)$ is solvable, a contradiction.

\item[(23)]$ R_1:h_2^2h_3^{-1}h_2^{-1}(h_3^{-1}h_2)^2=1$, $R_2:h_2h_3h_2^{-1}h_3^{-1}h_2h_3^{-2}h_2=1$, $R_3:h_2h_3h_2^{-2}h_3h_2^{-1}h_3^{-1}h_2=1$,\\$R_4:h_2h_3^{-1}h_2^{-1}h_3^2h_2^{-1}h_3^{-1}h_2=1$:\\
$\Rightarrow\langle h_2,h_3\rangle\cong BS(1,5)$ is solvable, a contradiction.

\item[(24)]$ R_1:h_2^2h_3^{-1}h_2^{-1}(h_3^{-1}h_2)^2=1$, $R_2:h_2h_3h_2^{-1}h_3^{-1}h_2h_3^{-2}h_2=1$, $R_3:h_2h_3h_2^{-2}h_3h_2^{-1}h_3^{-1}h_2=1$,\\$R_4:h_2h_3^{-1}h_2(h_3^{-1}h_2h_3^{-1})^2h_2=1$:\\
$\Rightarrow\langle h_2,h_3\rangle\cong BS(3,1)$ is solvable, a contradiction.

\item[(25)]$ R_1:h_2^2h_3^{-1}h_2^{-1}(h_3^{-1}h_2)^2=1$, $R_2:h_2h_3h_2^{-1}h_3^{-1}h_2h_3^{-2}h_2=1$, $R_3:h_2h_3^{-1}h_2^{-1}h_3^2h_2^{-2}h_3=1$,\\$R_4:h_2h_3^{-3}h_2^2h_3^{-1}h_2=1$:\\
$\Rightarrow\langle h_2,h_3\rangle\cong BS(1,5)$ is solvable, a contradiction.

\item[(26)]$ R_1:h_2^2h_3^{-1}h_2^{-1}(h_3^{-1}h_2)^2=1$, $R_2:h_2h_3h_2^{-1}h_3^{-1}h_2h_3^{-2}h_2=1$, $R_3:h_2h_3^{-1}h_2^{-1}h_3h_2^{-1}h_3^{-1}h_2h_3=1$,\\$R_4:h_2(h_3^{-1}h_2h_3^{-1})^2h_2h_3^{-1}h_2=1$:\\
$\Rightarrow\langle h_2,h_3\rangle\cong BS(1,5)$ is solvable, a contradiction.

\item[(27)]$ R_1:h_2^2h_3^{-1}h_2^{-1}(h_3^{-1}h_2)^2=1$, $R_2:h_2h_3h_2^{-1}h_3^{-1}h_2h_3^{-2}h_2=1$, $R_3:h_2h_3^{-1}h_2^{-1}h_3h_2^{-1}h_3^{-1}h_2h_3=1$,\\$R_4:h_2h_3^{-1}h_2(h_3^{-1}h_2h_3^{-1})^2h_2=1$:\\
$\Rightarrow\langle h_2,h_3\rangle\cong BS(3,1)$ is solvable, a contradiction.

\item[(28)]$ R_1:h_2^2h_3^{-1}h_2^{-1}(h_3^{-1}h_2)^2=1$, $R_2:h_2h_3h_2^{-1}h_3^{-1}h_2h_3^{-2}h_2=1$, $R_3:h_2h_3^{-2}h_2^{-1}h_3h_2h_3^{-1}h_2=1$,\\$R_4:h_2h_3^{-1}h_2(h_3^{-1}h_2h_3^{-1})^2h_2=1$:\\
$\Rightarrow\langle h_2,h_3\rangle\cong BS(3,1)$ is solvable, a contradiction.

\item[(29)]$ R_1:h_2^2h_3^{-1}h_2^{-1}h_3h_2^{-2}h_3=1$, $R_2:h_2^2h_3^{-2}(h_3^{-1}h_2)^2=1$, $R_3:h_2^2h_3^{-1}h_2^{-1}h_3h_2^{-2}h_3=1$,\\$R_4:(h_2h_3^{-1})^2(h_2^{-1}h_3)^3=1$:\\
$\Rightarrow\langle h_2,h_3\rangle\cong BS(1,1)$ is solvable, a contradiction.

\item[(30)]$ R_1:h_2^2h_3^{-1}h_2^{-1}h_3h_2^{-2}h_3=1$, $R_2:h_2^2h_3^{-2}(h_3^{-1}h_2)^2=1$, $R_3:h_2(h_3h_2^{-2})^2h_3=1$,\\$R_4:(h_2h_3^{-1})^2(h_2^{-1}h_3)^3=1$:\\
$\Rightarrow\langle h_2,h_3\rangle\cong BS(1,1)$ is solvable, a contradiction.

\item[(31)]$ R_1:h_2^2h_3^{-1}h_2^{-1}h_3h_2^{-2}h_3=1$, $R_2:h_2^2h_3^{-2}(h_3^{-1}h_2)^2=1$, $R_3:h_2(h_3h_2^{-2})^2h_3=1$,\\$R_4:(h_2h_3^{-1})^3(h_2^{-1}h_3)^2=1$:\\
$\Rightarrow\langle h_2,h_3\rangle\cong BS(-3,1)$ is solvable, a contradiction.

\item[(32)]$ R_1:h_2^2h_3^{-1}h_2^{-1}h_3h_2^{-2}h_3=1$, $R_2:h_2^2h_3^{-2}(h_3^{-1}h_2)^2=1$, $R_3:h_2h_3h_2^{-1}h_3^{-1}h_2h_3^{-2}h_2=1$,\\$R_4:h_2h_3^{-1}h_2^{-1}h_3^2h_2^{-2}h_3=1$:\\
$\Rightarrow\langle h_2,h_3\rangle\cong BS(1,2)$ is solvable, a contradiction.

\item[(33)]$ R_1:h_2^2h_3^{-1}h_2^{-1}h_3h_2^{-2}h_3=1$, $R_2:h_2^2h_3^{-2}(h_3^{-1}h_2)^2=1$, $R_3:h_2h_3^{-1}h_2^{-1}h_3h_2h_3^{-2}h_2=1$,\\$R_4:h_2h_3^{-2}h_2^{-1}h_3^2h_2^{-1}h_3=1$:\\
$\Rightarrow\langle h_2,h_3\rangle\cong BS(4,1)$ is solvable, a contradiction.

\item[(34)]$ R_1:h_2^2h_3^{-1}h_2^{-1}h_3h_2^{-2}h_3=1$, $R_2:h_2^2h_3^{-2}(h_3^{-1}h_2)^2=1$, $R_3:h_2h_3^{-1}h_2^{-1}h_3^2h_2^{-1}h_3^{-1}h_2=1$,\\$R_4:h_2h_3^{-2}h_2^{-1}h_3^2h_2^{-1}h_3=1$:\\
$\Rightarrow\langle h_2,h_3\rangle\cong BS(4,1)$ is solvable, a contradiction.

\item[(35)]$ R_1:h_2^2h_3^{-1}h_2^{-1}h_3h_2^{-2}h_3=1$, $R_2:h_2^2h_3^{-2}(h_3^{-1}h_2)^2=1$, $R_3:h_2h_3^{-1}h_2^{-1}h_3^2h_2^{-1}h_3^{-1}h_2=1$,\\$R_4:h_2h_3^{-3}h_2h_3h_2^{-1}h_3=1$:\\
$\Rightarrow\langle h_2,h_3\rangle\cong BS(4,1)$ is solvable, a contradiction.

\item[(36)]$ R_1:h_2^2h_3^{-1}h_2^{-1}h_3h_2^{-2}h_3=1$, $R_2:h_2^2h_3^{-2}(h_3^{-1}h_2)^2=1$, $R_3:h_2h_3^{-1}h_2^{-1}h_3h_2^{-2}h_3^2=1$,\\$R_4:(h_2h_3^{-1})^2(h_2^{-1}h_3)^3=1$:\\
$\Rightarrow\langle h_2,h_3\rangle\cong BS(1,1)$ is solvable, a contradiction.

\item[(37)]$ R_1:h_2^2h_3^{-1}h_2^{-1}h_3h_2^{-2}h_3=1$, $R_2:h_2^2h_3^{-2}(h_3^{-1}h_2)^2=1$, $R_3:h_2h_3^{-1}h_2^{-1}h_3h_2^{-2}h_3^2=1$,\\$R_4:(h_2h_3^{-1})^3(h_2^{-1}h_3)^2=1$:\\
$\Rightarrow\langle h_2,h_3\rangle\cong BS(-3,1)$ is solvable, a contradiction.

\item[(38)]$ R_1:h_2^2h_3^{-1}h_2^{-1}h_3h_2^{-2}h_3=1$, $R_2:h_2^2h_3^{-2}(h_3^{-1}h_2)^2=1$, $R_3:h_2h_3^{-3}h_2^2h_3^{-1}h_2=1$,\\$R_4:(h_2h_3^{-1})^2(h_2^{-1}h_3)^3=1$:\\
$\Rightarrow\langle h_2,h_3\rangle\cong BS(1,1)$ is solvable, a contradiction.

\item[(39)]$ R_1:h_2^2h_3^{-1}h_2^{-1}h_3h_2^{-2}h_3=1$, $R_2:h_2^2h_3^{-2}(h_3^{-1}h_2)^2=1$, $R_3:h_2h_3^{-3}h_2^2h_3^{-1}h_2=1$,\\$R_4:(h_2h_3^{-1})^3(h_2^{-1}h_3)^2=1$:\\
$\Rightarrow\langle h_2,h_3\rangle\cong BS(-3,1)$ is solvable, a contradiction.

\item[(40)]$ R_1:h_2^2h_3^{-1}h_2^{-1}h_3h_2^{-1}h_3^{-1}h_2=1$, $R_2:h_2^2h_3^{-2}(h_2^{-1}h_3)^2=1$, $R_3:h_2^2h_3^{-1}h_2^{-1}h_3h_2^{-1}h_3^{-1}h_2=1$,\\$R_4:h_2(h_3^{-1}h_2h_3^{-1})^2h_2h_3^{-1}h_2=1$:\\
$\Rightarrow\langle h_2,h_3\rangle\cong BS(5,1)$ is solvable, a contradiction.

\item[(41)]$ R_1:h_2^2h_3^{-1}h_2^{-1}h_3h_2^{-1}h_3^{-1}h_2=1$, $R_2:h_2^2h_3^{-2}(h_2^{-1}h_3)^2=1$, $R_3:h_2^2h_3^{-1}h_2^{-1}h_3h_2^{-1}h_3^{-1}h_2=1$,\\$R_4:h_2h_3^{-1}h_2(h_3^{-1}h_2h_3^{-1})^2h_2=1$:\\
$\Rightarrow\langle h_2,h_3\rangle\cong BS(5,1)$ is solvable, a contradiction.

\item[(42)]$ R_1:h_2^2h_3^{-1}h_2^{-1}h_3h_2^{-1}h_3^{-1}h_2=1$, $R_2:h_2^2h_3^{-2}(h_2^{-1}h_3)^2=1$, $R_3:h_2(h_3h_2^{-2})^2h_3=1$,\\$R_4:h_2h_3^{-1}h_2^{-1}h_3^2h_2^{-1}h_3^{-1}h_2=1$:\\
$\Rightarrow\langle h_2,h_3\rangle\cong BS(5,1)$ is solvable, a contradiction.

\item[(43)]$ R_1:h_2^2h_3^{-1}h_2^{-1}h_3h_2^{-1}h_3^{-1}h_2=1$, $R_2:h_2^2h_3^{-2}(h_2^{-1}h_3)^2=1$, $R_3:h_2h_3h_2^{-2}h_3h_2^{-1}h_3^{-1}h_2=1$,\\$R_4:h_2h_3^{-1}h_2^{-1}h_3^2h_2^{-1}h_3^{-1}h_2=1$:\\
$\Rightarrow\langle h_2,h_3\rangle\cong BS(5,1)$ is solvable, a contradiction.

\item[(44)]$ R_1:h_2^2h_3^{-1}h_2^{-1}h_3h_2^{-1}h_3^{-1}h_2=1$, $R_2:h_2^2h_3^{-2}(h_2^{-1}h_3)^2=1$, $R_3:h_2h_3h_2^{-2}h_3h_2^{-1}h_3^{-1}h_2=1$,\\$R_4:h_2h_3^{-1}h_2(h_3^{-1}h_2h_3^{-1})^2h_2=1$:\\
$\Rightarrow\langle h_2,h_3\rangle\cong BS(5,1)$ is solvable, a contradiction.

\item[(45)]$ R_1:h_2^2h_3^{-1}h_2^{-1}h_3h_2^{-1}h_3^{-1}h_2=1$, $R_2:h_2^2h_3^{-2}(h_2^{-1}h_3)^2=1$, $R_3:h_2h_3^{-1}h_2^{-1}h_3^2h_2^{-2}h_3=1$,\\$R_4:h_2h_3^{-3}h_2^2h_3^{-1}h_2=1$:\\
$\Rightarrow\langle h_2,h_3\rangle\cong BS(5,1)$ is solvable, a contradiction.

\item[(46)]$ R_1:h_2^2h_3^{-1}h_2^{-1}h_3h_2^{-1}h_3^{-1}h_2=1$, $R_2:h_2^2h_3^{-2}(h_2^{-1}h_3)^2=1$, $R_3:h_2h_3^{-1}h_2^{-1}h_3h_2^{-1}h_3^{-1}h_2h_3=1$,\\$R_4:h_2(h_3^{-1}h_2h_3^{-1})^2h_2h_3^{-1}h_2=1$:\\
$\Rightarrow\langle h_2,h_3\rangle\cong BS(5,1)$ is solvable, a contradiction.

\item[(47)]$ R_1:h_2^2h_3^{-1}h_2^{-1}h_3h_2^{-1}h_3^{-1}h_2=1$, $R_2:h_2^2h_3^{-2}(h_2^{-1}h_3)^2=1$, $R_3:h_2h_3^{-1}h_2^{-1}h_3h_2^{-1}h_3^{-1}h_2h_3=1$,\\$R_4:h_2h_3^{-1}h_2(h_3^{-1}h_2h_3^{-1})^2h_2=1$:\\
$\Rightarrow\langle h_2,h_3\rangle\cong BS(5,1)$ is solvable, a contradiction.

\item[(48)]$ R_1:h_2^2h_3^{-1}h_2^{-1}h_3h_2^{-1}h_3^{-1}h_2=1$, $R_2:h_2^2h_3^{-2}(h_2^{-1}h_3)^2=1$, $R_3:h_2h_3^{-2}h_2^{-1}h_3h_2h_3^{-1}h_2=1$,\\$R_4:h_2h_3^{-1}h_2(h_3^{-1}h_2h_3^{-1})^2h_2=1$:\\
$\Rightarrow\langle h_2,h_3\rangle\cong BS(5,1)$ is solvable, a contradiction.

\item[(49)]$ R_1:h_2^2h_3^{-1}h_2h_3h_2^{-1}h_3^2=1$, $R_2:h_2h_3h_2^{-1}h_3h_2h_3^{-1}h_2h_3=1$, $R_3:h_2^3h_3^3=1$,\\$R_4:h_2h_3h_2^{-1}h_3^{-1}h_2h_3^{-1}h_2^{-1}h_3=1$:\\
$\Rightarrow\langle h_2,h_3\rangle\cong BS(1,2)$ is solvable, a contradiction.

\item[(50)]$ R_1:(h_2h_3)^2h_2h_3^{-1}h_2=1$, $R_2:(h_2h_3)^2h_2h_3^{-1}h_2=1$, $R_3:h_2(h_3h_2^{-1})^2h_3^{-2}h_2=1$,\\$R_4:(h_2h_3^{-1}h_2)^2h_3^{-2}h_2=1$:\\
$\Rightarrow$  $G$ has a torsion element, a contradiction.

\item[(51)]$ R_1:(h_2h_3)^2h_2h_3^{-1}h_2=1$, $R_2:h_2(h_3h_2^{-1})^2h_3^{-2}h_2=1$, $R_3:h_2(h_3h_2^{-1}h_3)^2h_3=1$,\\$R_4:h_2(h_3h_2^{-1})^2h_3^{-2}h_2=1$:\\
$\Rightarrow$  $G$ has a torsion element, a contradiction.

\item[(52)]$ R_1:h_2h_3(h_2h_3^{-1}h_2)^2=1$, $R_2:h_2h_3(h_2h_3^{-1}h_2)^2=1$, $R_3:(h_2h_3)^2h_2^{-1}h_3^2=1$,\\$R_4:h_2(h_3h_2^{-1})^2h_3^{-2}h_2=1$:\\
$\Rightarrow$  $G$ has a torsion element, a contradiction.

\item[(53)]$ R_1:h_2(h_3h_2^{-2})^2h_3=1$, $R_2:h_2h_3h_2^{-1}h_3^{-1}(h_3^{-1}h_2)^2=1$, $R_3:h_2^2h_3^{-1}h_2^{-1}h_3h_2^{-2}h_3=1$,\\$R_4:(h_2h_3^{-1})^2(h_2^{-1}h_3)^3=1$:\\
$\Rightarrow\langle h_2,h_3\rangle\cong BS(1,1)$ is solvable, a contradiction.

\item[(54)]$ R_1:h_2(h_3h_2^{-2})^2h_3=1$, $R_2:h_2h_3h_2^{-1}h_3^{-1}(h_3^{-1}h_2)^2=1$, $R_3:h_2(h_3h_2^{-2})^2h_3=1$,\\$R_4:(h_2h_3^{-1})^2(h_2^{-1}h_3)^3=1$:\\
$\Rightarrow\langle h_2,h_3\rangle\cong BS(-5,1)$ is solvable, a contradiction.

\item[(55)]$ R_1:h_2(h_3h_2^{-2})^2h_3=1$, $R_2:h_2h_3h_2^{-1}h_3^{-1}(h_3^{-1}h_2)^2=1$, $R_3:h_2(h_3h_2^{-2})^2h_3=1$,\\$R_4:(h_2h_3^{-1})^3(h_2^{-1}h_3)^2=1$:\\
$\Rightarrow\langle h_2,h_3\rangle\cong BS(-5,1)$ is solvable, a contradiction.

\item[(56)]$ R_1:h_2(h_3h_2^{-2})^2h_3=1$, $R_2:h_2h_3h_2^{-1}h_3^{-1}(h_3^{-1}h_2)^2=1$, $R_3:h_2h_3h_2^{-1}h_3^{-1}h_2h_3^{-2}h_2=1$,\\$R_4:h_2h_3^{-1}h_2^{-1}h_3^2h_2^{-2}h_3=1$:\\
$\Rightarrow\langle h_2,h_3\rangle\cong BS(-5,1)$ is solvable, a contradiction.

\item[(57)]$ R_1:h_2(h_3h_2^{-2})^2h_3=1$, $R_2:h_2h_3h_2^{-1}h_3^{-1}(h_3^{-1}h_2)^2=1$, $R_3:h_2h_3^{-1}h_2^{-1}h_3h_2h_3^{-2}h_2=1$,\\$R_4:h_2h_3^{-2}h_2^{-1}h_3^2h_2^{-1}h_3=1$:\\
$\Rightarrow\langle h_2,h_3\rangle\cong BS(-5,1)$ is solvable, a contradiction.

\item[(58)]$ R_1:h_2(h_3h_2^{-2})^2h_3=1$, $R_2:h_2h_3h_2^{-1}h_3^{-1}(h_3^{-1}h_2)^2=1$, $R_3:h_2h_3^{-1}h_2^{-1}h_3^2h_2^{-1}h_3^{-1}h_2=1$,\\$R_4:h_2h_3^{-2}h_2^{-1}h_3^2h_2^{-1}h_3=1$:\\
$\Rightarrow\langle h_2,h_3\rangle\cong BS(-5,1)$ is solvable, a contradiction.

\item[(59)]$ R_1:h_2(h_3h_2^{-2})^2h_3=1$, $R_2:h_2h_3h_2^{-1}h_3^{-1}(h_3^{-1}h_2)^2=1$, $R_3:h_2h_3^{-1}h_2^{-1}h_3^2h_2^{-1}h_3^{-1}h_2=1$,\\$R_4:h_2h_3^{-3}h_2h_3h_2^{-1}h_3=1$:\\
$\Rightarrow\langle h_2,h_3\rangle\cong BS(-5,1)$ is solvable, a contradiction.

\item[(60)]$ R_1:h_2(h_3h_2^{-2})^2h_3=1$, $R_2:h_2h_3h_2^{-1}h_3^{-1}(h_3^{-1}h_2)^2=1$, $R_3:h_2h_3^{-1}h_2^{-1}h_3h_2^{-2}h_3^2=1$,\\$R_4:(h_2h_3^{-1})^2(h_2^{-1}h_3)^3=1$:\\
$\Rightarrow\langle h_2,h_3\rangle\cong BS(-5,1)$ is solvable, a contradiction.

\item[(61)]$ R_1:h_2(h_3h_2^{-2})^2h_3=1$, $R_2:h_2h_3h_2^{-1}h_3^{-1}(h_3^{-1}h_2)^2=1$, $R_3:h_2h_3^{-1}h_2^{-1}h_3h_2^{-2}h_3^2=1$,\\$R_4:(h_2h_3^{-1})^3(h_2^{-1}h_3)^2=1$:\\
$\Rightarrow\langle h_2,h_3\rangle\cong BS(-5,1)$ is solvable, a contradiction.

\item[(62)]$ R_1:h_2(h_3h_2^{-2})^2h_3=1$, $R_2:h_2h_3h_2^{-1}h_3^{-1}(h_3^{-1}h_2)^2=1$, $R_3:h_2h_3^{-3}h_2^2h_3^{-1}h_2=1$,\\$R_4:(h_2h_3^{-1})^2(h_2^{-1}h_3)^3=1$:\\
$\Rightarrow\langle h_2,h_3\rangle\cong BS(-5,1)$ is solvable, a contradiction.

\item[(63)]$ R_1:h_2(h_3h_2^{-2})^2h_3=1$, $R_2:h_2h_3h_2^{-1}h_3^{-1}(h_3^{-1}h_2)^2=1$, $R_3:h_2h_3^{-3}h_2^2h_3^{-1}h_2=1$,\\$R_4:(h_2h_3^{-1})^3(h_2^{-1}h_3)^2=1$:\\
$\Rightarrow\langle h_2,h_3\rangle\cong BS(-5,1)$ is solvable, a contradiction.

\item[(64)]$ R_1:h_2h_3h_2^{-2}h_3h_2^{-1}h_3^{-1}h_2=1$, $R_2:h_2h_3h_2^{-1}h_3^{-1}(h_2^{-1}h_3)^2=1$, $R_3:h_2^2h_3^{-1}h_2^{-1}h_3h_2^{-1}h_3^{-1}h_2=1$,\\$R_4:h_2(h_3^{-1}h_2h_3^{-1})^2h_2h_3^{-1}h_2=1$:\\
$\Rightarrow\langle h_2,h_3\rangle\cong BS(-4,1)$ is solvable, a contradiction.

\item[(65)]$ R_1:h_2h_3h_2^{-2}h_3h_2^{-1}h_3^{-1}h_2=1$, $R_2:h_2h_3h_2^{-1}h_3^{-1}(h_2^{-1}h_3)^2=1$, $R_3:h_2^2h_3^{-1}h_2^{-1}h_3h_2^{-1}h_3^{-1}h_2=1$,\\$R_4:h_2h_3^{-1}h_2(h_3^{-1}h_2h_3^{-1})^2h_2=1$:\\
$\Rightarrow\langle h_2,h_3\rangle\cong BS(-4,1)$ is solvable, a contradiction.

\item[(66)]$ R_1:h_2h_3h_2^{-2}h_3h_2^{-1}h_3^{-1}h_2=1$, $R_2:h_2h_3h_2^{-1}h_3^{-1}(h_2^{-1}h_3)^2=1$, $R_3:h_2(h_3h_2^{-2})^2h_3=1$,\\$R_4:h_2h_3^{-1}h_2^{-1}h_3^2h_2^{-1}h_3^{-1}h_2=1$:\\
$\Rightarrow\langle h_2,h_3\rangle\cong BS(-4,1)$ is solvable, a contradiction.

\item[(67)]$ R_1:h_2h_3h_2^{-2}h_3h_2^{-1}h_3^{-1}h_2=1$, $R_2:h_2h_3h_2^{-1}h_3^{-1}(h_2^{-1}h_3)^2=1$, $R_3:h_2h_3h_2^{-2}h_3h_2^{-1}h_3^{-1}h_2=1$,\\$R_4:h_2h_3^{-1}h_2^{-1}h_3^2h_2^{-1}h_3^{-1}h_2=1$:\\
$\Rightarrow\langle h_2,h_3\rangle\cong BS(-4,1)$ is solvable, a contradiction.

\item[(68)]$ R_1:h_2h_3h_2^{-2}h_3h_2^{-1}h_3^{-1}h_2=1$, $R_2:h_2h_3h_2^{-1}h_3^{-1}(h_2^{-1}h_3)^2=1$, $R_3:h_2h_3h_2^{-2}h_3h_2^{-1}h_3^{-1}h_2=1$,\\$R_4:h_2h_3^{-1}h_2(h_3^{-1}h_2h_3^{-1})^2h_2=1$:\\
$\Rightarrow\langle h_2,h_3\rangle\cong BS(-4,1)$ is solvable, a contradiction.

\item[(69)]$ R_1:h_2h_3h_2^{-2}h_3h_2^{-1}h_3^{-1}h_2=1$, $R_2:h_2h_3h_2^{-1}h_3^{-1}(h_2^{-1}h_3)^2=1$, $R_3:h_2h_3^{-1}h_2^{-1}h_3^2h_2^{-2}h_3=1$,\\$R_4:h_2h_3^{-3}h_2^2h_3^{-1}h_2=1$:\\
$\Rightarrow\langle h_2,h_3\rangle\cong BS(-4,1)$ is solvable, a contradiction.

\item[(70)]$ R_1:h_2h_3h_2^{-2}h_3h_2^{-1}h_3^{-1}h_2=1$, $R_2:h_2h_3h_2^{-1}h_3^{-1}(h_2^{-1}h_3)^2=1$, $R_3:h_2h_3^{-1}h_2^{-1}h_3h_2^{-1}h_3^{-1}h_2h_3=1$,\\$R_4:h_2(h_3^{-1}h_2h_3^{-1})^2h_2h_3^{-1}h_2=1$:\\
$\Rightarrow\langle h_2,h_3\rangle\cong BS(-4,1)$ is solvable, a contradiction.

\item[(71)]$ R_1:h_2h_3h_2^{-2}h_3h_2^{-1}h_3^{-1}h_2=1$, $R_2:h_2h_3h_2^{-1}h_3^{-1}(h_2^{-1}h_3)^2=1$, $R_3:h_2h_3^{-1}h_2^{-1}h_3h_2^{-1}h_3^{-1}h_2h_3=1$,\\$R_4:h_2h_3^{-1}h_2(h_3^{-1}h_2h_3^{-1})^2h_2=1$:\\
$\Rightarrow\langle h_2,h_3\rangle\cong BS(-4,1)$ is solvable, a contradiction.

\item[(72)]$ R_1:h_2h_3h_2^{-2}h_3h_2^{-1}h_3^{-1}h_2=1$, $R_2:h_2h_3h_2^{-1}h_3^{-1}(h_2^{-1}h_3)^2=1$, $R_3:h_2h_3^{-2}h_2^{-1}h_3h_2h_3^{-1}h_2=1$,\\$R_4:h_2h_3^{-1}h_2(h_3^{-1}h_2h_3^{-1})^2h_2=1$:\\
$\Rightarrow\langle h_2,h_3\rangle\cong BS(-4,1)$ is solvable, a contradiction.

\item[(73)]$ R_1:h_2h_3h_2^{-1}h_3^{-1}h_2h_3^{-1}h_2^{-1}h_3=1$, $R_2:h_2h_3^{-1}h_2^{-1}h_3^2h_2^{-1}h_3^{-1}h_2=1$, $R_3:(h_2h_3)^2h_2^{-1}h_3^{-1}h_2=1$,\\$R_4:(h_2h_3^{-1})^2h_2h_3h_2^{-1}h_3^2=1$:\\
$\Rightarrow \langle h_2,h_3\rangle=\langle h_2\rangle$ is abelian, a contradiction.

\item[(74)]$ R_1:h_2(h_3h_2^{-1})^2h_3^{-2}h_2=1$, $R_2:(h_2h_3^{-1}h_2)^2h_3^{-2}h_2=1$, $R_3:(h_2h_3)^2h_2h_3^{-1}h_2=1$,\\$R_4:h_2(h_3h_2^{-1})^2h_3^{-2}h_2=1$:\\
$\Rightarrow$  $G$ has a torsion element, a contradiction.

\item[(75)]$ R_1:h_2^2h_3^{-2}(h_2^{-1}h_3)^2=1$, $R_2:h_2h_3^3h_2^{-1}h_3^{-1}h_2=1$, $R_3:h_2^2h_3^{-1}h_2^{-1}h_3^3=1$,\\$R_4:h_2h_3^3h_2^{-1}h_3^{-1}h_2=1$:\\
By interchanging $h_2$ and $h_3$ in (1) and with the same discussion, there is a contradiction.

\item[(76)]$ R_1:h_2h_3^2h_2^{-1}h_3^{-1}h_2^{-1}h_3=1$, $R_2:h_2h_3h_2^{-1}h_3^{-3}h_2=1$, $R_3:h_2^2h_3^{-1}(h_2^{-1}h_3)^3=1$,\\$R_4:h_2h_3^2h_2^{-1}h_3^{-1}h_2^{-1}h_3=1$:\\
By interchanging $h_2$ and $h_3$ in (2) and with the same discussion, there is a contradiction.

\item[(77)]$ R_1:h_2h_3^2h_2^{-1}h_3^{-1}h_2^{-1}h_3=1$, $R_2:h_2h_3h_2^{-1}h_3^{-3}h_2=1$, $R_3:h_2^2h_3^{-2}(h_2^{-1}h_3)^2=1$,\\$R_4:h_2h_3^{-1}h_2^{-1}h_3h_2^{-1}h_3^{-1}h_2h_3=1$:\\
By interchanging $h_2$ and $h_3$ in (4) and with the same discussion, there is a contradiction.

\item[(78)]$ R_1:h_2h_3^2h_2^{-1}h_3^{-1}h_2^{-1}h_3=1$, $R_2:h_2h_3h_2^{-1}h_3^{-3}h_2=1$, $R_3:h_2h_3^2h_2^{-1}(h_2^{-1}h_3)^2=1$,\\$R_4:h_2h_3^{-1}h_2^{-1}h_3h_2^{-1}h_3^{-1}h_2h_3=1$:\\
By interchanging $h_2$ and $h_3$ in (5) and with the same discussion, there is a contradiction.

\item[(79)]$ R_1:h_2h_3^2h_2^{-1}h_3^{-1}h_2^{-1}h_3=1$, $R_2:h_2h_3h_2^{-1}h_3^{-3}h_2=1$, $R_3:h_2h_3^2h_2^{-1}h_3^{-1}h_2^{-1}h_3=1$,\\$R_4:h_2h_3^{-1}h_2^{-1}h_3h_2^{-1}(h_3^{-1}h_2)^2=1$:\\
By interchanging $h_2$ and $h_3$ in (3) and with the same discussion, there is a contradiction.

\item[(80)]$ R_1:h_2h_3h_2^{-1}(h_2^{-1}h_3)^2h_3=1$, $R_2:h_2h_3h_2^{-1}h_3^{-3}h_2=1$, $R_3:h_2h_3^2h_2^{-1}(h_2^{-1}h_3)^2=1$,\\$R_4:h_2h_3^{-1}h_2^{-1}h_3h_2^{-1}h_3^{-1}h_2h_3=1$:\\
By interchanging $h_2$ and $h_3$ in (8) and with the same discussion, there is a contradiction.

\item[(81)]$ R_1:h_2h_3h_2^{-1}(h_2^{-1}h_3)^2h_3=1$, $R_2:h_2h_3h_2^{-1}h_3^{-3}h_2=1$, $R_3:h_2h_3^2h_2^{-1}h_3^{-1}h_2^{-1}h_3=1$,\\$R_4:h_2h_3h_2^{-1}(h_2^{-1}h_3)^2h_3=1$:\\
By interchanging $h_2$ and $h_3$ in (6) and with the same discussion, there is a contradiction.

\item[(82)]$ R_1:h_2h_3h_2^{-1}(h_2^{-1}h_3)^2h_3=1$, $R_2:h_2h_3h_2^{-1}h_3^{-3}h_2=1$, $R_3:h_2h_3^2h_2^{-1}h_3^{-1}h_2^{-1}h_3=1$,\\$R_4:h_2h_3^{-1}h_2^{-1}h_3h_2^{-1}(h_3^{-1}h_2)^2=1$:\\
By interchanging $h_2$ and $h_3$ in (7) and with the same discussion, there is a contradiction.

\item[(83)]$ R_1:h_2h_3^{-2}h_2^{-1}h_3h_2h_3^{-1}h_2=1$, $R_2:h_2h_3^{-2}(h_3^{-1}h_2)^2h_3=1$, $R_3:h_2^2h_3^{-1}h_2^{-1}h_3h_2^{-2}h_3=1$,\\$R_4:h_2h_3^{-1}h_2^{-1}h_3h_2h_3^{-2}h_2=1$:\\
By interchanging $h_2$ and $h_3$ in (13) and with the same discussion, there is a contradiction.

\item[(84)]$ R_1:h_2h_3^{-2}h_2^{-1}h_3h_2h_3^{-1}h_2=1$, $R_2:h_2h_3^{-2}(h_3^{-1}h_2)^2h_3=1$, $R_3:h_2^2h_3^{-1}h_2^{-1}h_3h_2^{-2}h_3=1$,\\$R_4:h_2h_3^{-1}h_2^{-1}h_3^2h_2^{-1}h_3^{-1}h_2=1$:\\
By interchanging $h_2$ and $h_3$ in (14) and with the same discussion, there is a contradiction.

\item[(85)]$ R_1:h_2h_3^{-2}h_2^{-1}h_3h_2h_3^{-1}h_2=1$, $R_2:h_2h_3^{-2}(h_3^{-1}h_2)^2h_3=1$, $R_3:h_2^2h_3^{-1}h_2^{-1}h_3h_2^{-1}h_3^{-1}h_2=1$,\\$R_4:h_2h_3^{-1}h_2^{-1}h_3^2h_2^{-1}h_3^{-1}h_2=1$:\\
By interchanging $h_2$ and $h_3$ in (15) and with the same discussion, there is a contradiction.

\item[(86)]$ R_1:h_2h_3^{-2}h_2^{-1}h_3h_2h_3^{-1}h_2=1$, $R_2:h_2h_3^{-2}(h_3^{-1}h_2)^2h_3=1$, $R_3:h_2(h_2h_3^{-2})^2h_2=1$,\\$R_4:(h_2h_3^{-1})^2(h_2^{-1}h_3)^3=1$:\\
By interchanging $h_2$ and $h_3$ in (18) and with the same discussion, there is a contradiction.

\item[(87)]$ R_1:h_2h_3^{-2}h_2^{-1}h_3h_2h_3^{-1}h_2=1$, $R_2:h_2h_3^{-2}(h_3^{-1}h_2)^2h_3=1$, $R_3:h_2(h_2h_3^{-2})^2h_2=1$,\\$R_4:(h_2h_3^{-1})^3(h_2^{-1}h_3)^2=1$:\\
By interchanging $h_2$ and $h_3$ in (19) and with the same discussion, there is a contradiction.

\item[(88)]$ R_1:h_2h_3^{-2}h_2^{-1}h_3h_2h_3^{-1}h_2=1$, $R_2:h_2h_3^{-2}(h_3^{-1}h_2)^2h_3=1$, $R_3:h_2h_3h_2^{-1}h_3^{-1}h_2h_3^{-2}h_2=1$,\\$R_4:(h_2h_3^{-1})^2(h_2^{-1}h_3)^3=1$:\\
By interchanging $h_2$ and $h_3$ in (16) and with the same discussion, there is a contradiction.

\item[(89)]$ R_1:h_2h_3^{-2}h_2^{-1}h_3h_2h_3^{-1}h_2=1$, $R_2:h_2h_3^{-2}(h_3^{-1}h_2)^2h_3=1$, $R_3:h_2h_3h_2^{-1}h_3^{-1}h_2h_3^{-2}h_2=1$,\\$R_4:(h_2h_3^{-1})^3(h_2^{-1}h_3)^2=1$:\\
By interchanging $h_2$ and $h_3$ in (17) and with the same discussion, there is a contradiction.

\item[(90)]$ R_1:h_2h_3^{-2}h_2^{-1}h_3h_2h_3^{-1}h_2=1$, $R_2:h_2h_3^{-2}(h_3^{-1}h_2)^2h_3=1$, $R_3:h_2h_3^{-1}h_2^{-1}h_3^2h_2^{-2}h_3=1$,\\$R_4:h_2h_3^{-1}h_2^{-1}h_3h_2^{-2}h_3^2=1$:\\
By interchanging $h_2$ and $h_3$ in (12) and with the same discussion, there is a contradiction.

\item[(91)]$ R_1:h_2h_3^{-2}h_2^{-1}h_3h_2h_3^{-1}h_2=1$, $R_2:h_2h_3^{-2}(h_3^{-1}h_2)^2h_3=1$, $R_3:h_2h_3^{-2}h_2^{-1}h_3^2h_2^{-1}h_3=1$,\\$R_4:(h_2h_3^{-1})^2(h_2^{-1}h_3)^3=1$:\\
By interchanging $h_2$ and $h_3$ in (9) and with the same discussion, there is a contradiction.

\item[(92)]$ R_1:h_2h_3^{-2}h_2^{-1}h_3h_2h_3^{-1}h_2=1$, $R_2:h_2h_3^{-2}(h_3^{-1}h_2)^2h_3=1$, $R_3:(h_2h_3^{-2})^2h_2h_3=1$,\\$R_4:(h_2h_3^{-1})^2(h_2^{-1}h_3)^3=1$:\\
By interchanging $h_2$ and $h_3$ in (10) and with the same discussion, there is a contradiction.

\item[(93)]$ R_1:h_2h_3^{-2}h_2^{-1}h_3h_2h_3^{-1}h_2=1$, $R_2:h_2h_3^{-2}(h_3^{-1}h_2)^2h_3=1$, $R_3:(h_2h_3^{-2})^2h_2h_3=1$,\\$R_4:(h_2h_3^{-1})^3(h_2^{-1}h_3)^2=1$:\\
By interchanging $h_2$ and $h_3$ in (11) and with the same discussion, there is a contradiction.

\item[(94)]$ R_1:h_2h_3^{-1}h_2^{-1}h_3h_2^{-2}h_3^2=1$, $R_2:h_2h_3^{-2}(h_3^{-1}h_2)^2h_3=1$, $R_3:h_2^2h_3^{-2}h_2h_3^{-1}h_2^{-1}h_3=1$,\\$R_4:h_2h_3^{-1}h_2(h_3^{-1}h_2h_3^{-1})^2h_2=1$:\\
By interchanging $h_2$ and $h_3$ in (28) and with the same discussion, there is a contradiction.

\item[(95)]$ R_1:h_2h_3^{-1}h_2^{-1}h_3h_2^{-2}h_3^2=1$, $R_2:h_2h_3^{-2}(h_3^{-1}h_2)^2h_3=1$, $R_3:h_2(h_2h_3^{-2})^2h_2=1$,\\$R_4:h_2h_3^{-1}h_2^{-1}h_3^2h_2^{-2}h_3=1$:\\
By interchanging $h_2$ and $h_3$ in (25) and with the same discussion, there is a contradiction.

\item[(96)]$ R_1:h_2h_3^{-1}h_2^{-1}h_3h_2^{-2}h_3^2=1$, $R_2:h_2h_3^{-2}(h_3^{-1}h_2)^2h_3=1$, $R_3:h_2h_3h_2^{-1}h_3^{-1}h_2h_3^{-1}h_2^{-1}h_3=1$,\\$R_4:h_2(h_3^{-1}h_2h_3^{-1})^2h_2h_3^{-1}h_2=1$:\\
By interchanging $h_2$ and $h_3$ in (26) and with the same discussion, there is a contradiction.

\item[(97)]$ R_1:h_2h_3^{-1}h_2^{-1}h_3h_2^{-2}h_3^2=1$, $R_2:h_2h_3^{-2}(h_3^{-1}h_2)^2h_3=1$, $R_3:h_2h_3h_2^{-1}h_3^{-1}h_2h_3^{-1}h_2^{-1}h_3=1$,\\$R_4:h_2h_3^{-1}h_2(h_3^{-1}h_2h_3^{-1})^2h_2=1$:\\
By interchanging $h_2$ and $h_3$ in (27) and with the same discussion, there is a contradiction.

\item[(98)]$ R_1:h_2h_3^{-1}h_2^{-1}h_3h_2^{-2}h_3^2=1$, $R_2:h_2h_3^{-2}(h_3^{-1}h_2)^2h_3=1$, $R_3:h_2h_3^{-1}h_2^{-1}h_3^2h_2^{-1}h_3^{-1}h_2=1$,\\$R_4:h_2h_3^{-2}h_2h_3^{-1}h_2^{-1}h_3^2=1$:\\
By interchanging $h_2$ and $h_3$ in (23) and with the same discussion, there is a contradiction.

\item[(99)]$ R_1:h_2h_3^{-1}h_2^{-1}h_3h_2^{-2}h_3^2=1$, $R_2:h_2h_3^{-2}(h_3^{-1}h_2)^2h_3=1$, $R_3:h_2h_3^{-1}h_2^{-1}h_3^2h_2^{-1}h_3^{-1}h_2=1$,\\$R_4:(h_2h_3^{-2})^2h_2h_3=1$:\\
By interchanging $h_2$ and $h_3$ in (22) and with the same discussion, there is a contradiction.

\item[(100)]$ R_1:h_2h_3^{-1}h_2^{-1}h_3h_2^{-2}h_3^2=1$, $R_2:h_2h_3^{-2}(h_3^{-1}h_2)^2h_3=1$, $R_3:h_2h_3^{-3}h_2h_3h_2^{-1}h_3=1$,\\$R_4:h_2(h_3^{-1}h_2h_3^{-1})^2h_2h_3^{-1}h_2=1$:\\
By interchanging $h_2$ and $h_3$ in (20) and with the same discussion, there is a contradiction.

\item[(101)]$ R_1:h_2h_3^{-1}h_2^{-1}h_3h_2^{-2}h_3^2=1$, $R_2:h_2h_3^{-2}(h_3^{-1}h_2)^2h_3=1$, $R_3:h_2h_3^{-3}h_2h_3h_2^{-1}h_3=1$,\\$R_4:h_2h_3^{-1}h_2(h_3^{-1}h_2h_3^{-1})^2h_2=1$:\\
By interchanging $h_2$ and $h_3$ in (21) and with the same discussion, there is a contradiction.

\item[(102)]$ R_1:h_2h_3^{-1}h_2^{-1}h_3h_2^{-2}h_3^2=1$, $R_2:h_2h_3^{-2}(h_3^{-1}h_2)^2h_3=1$, $R_3:h_2h_3^{-2}h_2h_3^{-1}h_2^{-1}h_3^2=1$,\\$R_4:h_2h_3^{-1}h_2(h_3^{-1}h_2h_3^{-1})^2h_2=1$:\\
By interchanging $h_2$ and $h_3$ in (24) and with the same discussion, there is a contradiction.

\item[(103)]$ R_1:h_2^2h_3^{-2}(h_3^{-1}h_2)^2=1$, $R_2:h_2h_3^{-2}h_2^{-1}h_3^2h_2^{-1}h_3=1$, $R_3:h_2^2h_3^{-1}h_2^{-1}h_3h_2^{-2}h_3=1$,\\$R_4:h_2h_3^{-1}h_2^{-1}h_3h_2h_3^{-2}h_2=1$:\\
By interchanging $h_2$ and $h_3$ in (33) and with the same discussion, there is a contradiction.

\item[(104)]$ R_1:h_2^2h_3^{-2}(h_3^{-1}h_2)^2=1$, $R_2:h_2h_3^{-2}h_2^{-1}h_3^2h_2^{-1}h_3=1$, $R_3:h_2^2h_3^{-1}h_2^{-1}h_3h_2^{-2}h_3=1$,\\$R_4:h_2h_3^{-1}h_2^{-1}h_3^2h_2^{-1}h_3^{-1}h_2=1$:\\
By interchanging $h_2$ and $h_3$ in (34) and with the same discussion, there is a contradiction.

\item[(105)]$ R_1:h_2^2h_3^{-2}(h_3^{-1}h_2)^2=1$, $R_2:h_2h_3^{-2}h_2^{-1}h_3^2h_2^{-1}h_3=1$, $R_3:h_2^2h_3^{-1}h_2^{-1}h_3h_2^{-1}h_3^{-1}h_2=1$,\\$R_4:h_2h_3^{-1}h_2^{-1}h_3^2h_2^{-1}h_3^{-1}h_2=1$:\\
By interchanging $h_2$ and $h_3$ in (35) and with the same discussion, there is a contradiction.

\item[(106)]$ R_1:h_2^2h_3^{-2}(h_3^{-1}h_2)^2=1$, $R_2:h_2h_3^{-2}h_2^{-1}h_3^2h_2^{-1}h_3=1$, $R_3:h_2(h_2h_3^{-2})^2h_2=1$,\\$R_4:(h_2h_3^{-1})^2(h_2^{-1}h_3)^3=1$:\\
By interchanging $h_2$ and $h_3$ in (38) and with the same discussion, there is a contradiction.

\item[(107)]$ R_1:h_2^2h_3^{-2}(h_3^{-1}h_2)^2=1$, $R_2:h_2h_3^{-2}h_2^{-1}h_3^2h_2^{-1}h_3=1$, $R_3:h_2(h_2h_3^{-2})^2h_2=1$,\\$R_4:(h_2h_3^{-1})^3(h_2^{-1}h_3)^2=1$:\\
By interchanging $h_2$ and $h_3$ in (39) and with the same discussion, there is a contradiction.

\item[(108)]$ R_1:h_2^2h_3^{-2}(h_3^{-1}h_2)^2=1$, $R_2:h_2h_3^{-2}h_2^{-1}h_3^2h_2^{-1}h_3=1$, $R_3:h_2h_3h_2^{-1}h_3^{-1}h_2h_3^{-2}h_2=1$,\\$R_4:(h_2h_3^{-1})^2(h_2^{-1}h_3)^3=1$:\\
By interchanging $h_2$ and $h_3$ in (36) and with the same discussion, there is a contradiction.

\item[(109)]$ R_1:h_2^2h_3^{-2}(h_3^{-1}h_2)^2=1$, $R_2:h_2h_3^{-2}h_2^{-1}h_3^2h_2^{-1}h_3=1$, $R_3:h_2h_3h_2^{-1}h_3^{-1}h_2h_3^{-2}h_2=1$,\\$R_4:(h_2h_3^{-1})^3(h_2^{-1}h_3)^2=1$:\\
By interchanging $h_2$ and $h_3$ in (37) and with the same discussion, there is a contradiction.

\item[(110)]$ R_1:h_2^2h_3^{-2}(h_3^{-1}h_2)^2=1$, $R_2:h_2h_3^{-2}h_2^{-1}h_3^2h_2^{-1}h_3=1$, $R_3:h_2h_3^{-1}h_2^{-1}h_3^2h_2^{-2}h_3=1$,\\$R_4:h_2h_3^{-1}h_2^{-1}h_3h_2^{-2}h_3^2=1$:\\
By interchanging $h_2$ and $h_3$ in (32) and with the same discussion, there is a contradiction.

\item[(111)]$ R_1:h_2^2h_3^{-2}(h_3^{-1}h_2)^2=1$, $R_2:h_2h_3^{-2}h_2^{-1}h_3^2h_2^{-1}h_3=1$, $R_3:h_2h_3^{-2}h_2^{-1}h_3^2h_2^{-1}h_3=1$,\\$R_4:(h_2h_3^{-1})^2(h_2^{-1}h_3)^3=1$:\\
By interchanging $h_2$ and $h_3$ in (29) and with the same discussion, there is a contradiction.

\item[(112)]$ R_1:h_2^2h_3^{-2}(h_3^{-1}h_2)^2=1$, $R_2:h_2h_3^{-2}h_2^{-1}h_3^2h_2^{-1}h_3=1$, $R_3:(h_2h_3^{-2})^2h_2h_3=1$,\\$R_4:(h_2h_3^{-1})^2(h_2^{-1}h_3)^3=1$:\\
By interchanging $h_2$ and $h_3$ in (30) and with the same discussion, there is a contradiction.

\item[(113)]$ R_1:h_2^2h_3^{-2}(h_3^{-1}h_2)^2=1$, $R_2:h_2h_3^{-2}h_2^{-1}h_3^2h_2^{-1}h_3=1$, $R_3:(h_2h_3^{-2})^2h_2h_3=1$,\\$R_4:(h_2h_3^{-1})^3(h_2^{-1}h_3)^2=1$:\\
By interchanging $h_2$ and $h_3$ in (31) and with the same discussion, there is a contradiction.

\item[(114)]$ R_1:h_2^2h_3^{-2}(h_2^{-1}h_3)^2=1$, $R_2:h_2h_3^{-3}h_2h_3h_2^{-1}h_3=1$, $R_3:h_2^2h_3^{-2}h_2h_3^{-1}h_2^{-1}h_3=1$,\\$R_4:h_2h_3^{-1}h_2(h_3^{-1}h_2h_3^{-1})^2h_2=1$:\\
By interchanging $h_2$ and $h_3$ in (48) and with the same discussion, there is a contradiction.

\item[(115)]$ R_1:h_2^2h_3^{-2}(h_2^{-1}h_3)^2=1$, $R_2:h_2h_3^{-3}h_2h_3h_2^{-1}h_3=1$, $R_3:h_2(h_2h_3^{-2})^2h_2=1$,\\$R_4:h_2h_3^{-1}h_2^{-1}h_3^2h_2^{-2}h_3=1$:\\
By interchanging $h_2$ and $h_3$ in (45) and with the same discussion, there is a contradiction.

\item[(116)]$ R_1:h_2^2h_3^{-2}(h_2^{-1}h_3)^2=1$, $R_2:h_2h_3^{-3}h_2h_3h_2^{-1}h_3=1$, $R_3:h_2h_3h_2^{-1}h_3^{-1}h_2h_3^{-1}h_2^{-1}h_3=1$,\\$R_4:h_2(h_3^{-1}h_2h_3^{-1})^2h_2h_3^{-1}h_2=1$:\\
By interchanging $h_2$ and $h_3$ in (46) and with the same discussion, there is a contradiction.

\item[(117)]$ R_1:h_2^2h_3^{-2}(h_2^{-1}h_3)^2=1$, $R_2:h_2h_3^{-3}h_2h_3h_2^{-1}h_3=1$, $R_3:h_2h_3h_2^{-1}h_3^{-1}h_2h_3^{-1}h_2^{-1}h_3=1$,\\$R_4:h_2h_3^{-1}h_2(h_3^{-1}h_2h_3^{-1})^2h_2=1$:\\
By interchanging $h_2$ and $h_3$ in (47) and with the same discussion, there is a contradiction.

\item[(118)]$ R_1:h_2^2h_3^{-2}(h_2^{-1}h_3)^2=1$, $R_2:h_2h_3^{-3}h_2h_3h_2^{-1}h_3=1$, $R_3:h_2h_3^{-1}h_2^{-1}h_3^2h_2^{-1}h_3^{-1}h_2=1$,\\$R_4:h_2h_3^{-2}h_2h_3^{-1}h_2^{-1}h_3^2=1$:\\
By interchanging $h_2$ and $h_3$ in (43) and with the same discussion, there is a contradiction.

\item[(119)]$ R_1:h_2^2h_3^{-2}(h_2^{-1}h_3)^2=1$, $R_2:h_2h_3^{-3}h_2h_3h_2^{-1}h_3=1$, $R_3:h_2h_3^{-1}h_2^{-1}h_3^2h_2^{-1}h_3^{-1}h_2=1$,\\$R_4:(h_2h_3^{-2})^2h_2h_3=1$:\\
By interchanging $h_2$ and $h_3$ in (42) and with the same discussion, there is a contradiction.

\item[(120)]$ R_1:h_2^2h_3^{-2}(h_2^{-1}h_3)^2=1$, $R_2:h_2h_3^{-3}h_2h_3h_2^{-1}h_3=1$, $R_3:h_2h_3^{-3}h_2h_3h_2^{-1}h_3=1$,\\$R_4:h_2(h_3^{-1}h_2h_3^{-1})^2h_2h_3^{-1}h_2=1$:\\
By interchanging $h_2$ and $h_3$ in (40) and with the same discussion, there is a contradiction.

\item[(121)]$ R_1:h_2^2h_3^{-2}(h_2^{-1}h_3)^2=1$, $R_2:h_2h_3^{-3}h_2h_3h_2^{-1}h_3=1$, $R_3:h_2h_3^{-3}h_2h_3h_2^{-1}h_3=1$,\\$R_4:h_2h_3^{-1}h_2(h_3^{-1}h_2h_3^{-1})^2h_2=1$:\\
By interchanging $h_2$ and $h_3$ in (41) and with the same discussion, there is a contradiction.

\item[(122)]$ R_1:h_2^2h_3^{-2}(h_2^{-1}h_3)^2=1$, $R_2:h_2h_3^{-3}h_2h_3h_2^{-1}h_3=1$, $R_3:h_2h_3^{-2}h_2h_3^{-1}h_2^{-1}h_3^2=1$,\\$R_4:h_2h_3^{-1}h_2(h_3^{-1}h_2h_3^{-1})^2h_2=1$:\\
By interchanging $h_2$ and $h_3$ in (44) and with the same discussion, there is a contradiction.

\item[(123)]$ R_1:h_2h_3h_2h_3^{-1}h_2h_3h_2^{-1}h_3=1$, $R_2:h_2h_3^2h_2^{-1}h_3h_2h_3^{-1}h_2=1$, $R_3:h_2^3h_3^3=1$,\\$R_4:h_2h_3^{-1}h_2^{-1}h_3h_2^{-1}h_3^{-1}h_2h_3=1$:\\
By interchanging $h_2$ and $h_3$ in (49) and with the same discussion, there is a contradiction.

\item[(124)]$ R_1:(h_2h_3)^2h_2^{-1}h_3^2=1$, $R_2:(h_2h_3)^2h_2^{-1}h_3^2=1$, $R_3:h_2(h_3h_2^{-1})^2h_3^{-2}h_2=1$,\\$R_4:h_2(h_3^{-1}h_2h_3^{-1})^2h_3^{-1}h_2=1$:\\
By interchanging $h_2$ and $h_3$ in (50) and with the same discussion, there is a contradiction.

\item[(125)]$ R_1:(h_2h_3)^2h_2^{-1}h_3^2=1$, $R_2:h_2(h_3h_2^{-1})^2h_3^{-2}h_2=1$, $R_3:h_2h_3(h_2h_3^{-1}h_2)^2=1$,\\$R_4:h_2(h_3h_2^{-1})^2h_3^{-2}h_2=1$:\\
By interchanging $h_2$ and $h_3$ in (51) and with the same discussion, there is a contradiction.

\item[(126)]$ R_1:h_2(h_3h_2^{-1}h_3)^2h_3=1$, $R_2:h_2(h_3h_2^{-1}h_3)^2h_3=1$, $R_3:(h_2h_3)^2h_2h_3^{-1}h_2=1$,\\$R_4:h_2(h_3h_2^{-1})^2h_3^{-2}h_2=1$:\\
By interchanging $h_2$ and $h_3$ in (52) and with the same discussion, there is a contradiction.

\item[(127)]$ R_1:h_2h_3h_2^{-1}h_3^{-1}(h_3^{-1}h_2)^2=1$, $R_2:(h_2h_3^{-2})^2h_2h_3=1$, $R_3:h_2^2h_3^{-1}h_2^{-1}h_3h_2^{-2}h_3=1$,\\$R_4:h_2h_3^{-1}h_2^{-1}h_3h_2h_3^{-2}h_2=1$:\\
By interchanging $h_2$ and $h_3$ in (57) and with the same discussion, there is a contradiction.

\item[(128)]$ R_1:h_2h_3h_2^{-1}h_3^{-1}(h_3^{-1}h_2)^2=1$, $R_2:(h_2h_3^{-2})^2h_2h_3=1$, $R_3:h_2^2h_3^{-1}h_2^{-1}h_3h_2^{-2}h_3=1$,\\$R_4:h_2h_3^{-1}h_2^{-1}h_3^2h_2^{-1}h_3^{-1}h_2=1$:\\
By interchanging $h_2$ and $h_3$ in (58) and with the same discussion, there is a contradiction.

\item[(129)]$ R_1:h_2h_3h_2^{-1}h_3^{-1}(h_3^{-1}h_2)^2=1$, $R_2:(h_2h_3^{-2})^2h_2h_3=1$, $R_3:h_2^2h_3^{-1}h_2^{-1}h_3h_2^{-1}h_3^{-1}h_2=1$,\\$R_4:h_2h_3^{-1}h_2^{-1}h_3^2h_2^{-1}h_3^{-1}h_2=1$:\\
By interchanging $h_2$ and $h_3$ in (59) and with the same discussion, there is a contradiction.

\item[(130)]$ R_1:h_2h_3h_2^{-1}h_3^{-1}(h_3^{-1}h_2)^2=1$, $R_2:(h_2h_3^{-2})^2h_2h_3=1$, $R_3:h_2(h_2h_3^{-2})^2h_2=1$,\\$R_4:(h_2h_3^{-1})^2(h_2^{-1}h_3)^3=1$:\\
By interchanging $h_2$ and $h_3$ in (62) and with the same discussion, there is a contradiction.

\item[(131)]$ R_1:h_2h_3h_2^{-1}h_3^{-1}(h_3^{-1}h_2)^2=1$, $R_2:(h_2h_3^{-2})^2h_2h_3=1$, $R_3:h_2(h_2h_3^{-2})^2h_2=1$,\\$R_4:(h_2h_3^{-1})^3(h_2^{-1}h_3)^2=1$:\\
By interchanging $h_2$ and $h_3$ in (63) and with the same discussion, there is a contradiction.

\item[(132)]$ R_1:h_2h_3h_2^{-1}h_3^{-1}(h_3^{-1}h_2)^2=1$, $R_2:(h_2h_3^{-2})^2h_2h_3=1$, $R_3:h_2h_3h_2^{-1}h_3^{-1}h_2h_3^{-2}h_2=1$,\\$R_4:(h_2h_3^{-1})^2(h_2^{-1}h_3)^3=1$:\\
By interchanging $h_2$ and $h_3$ in (60) and with the same discussion, there is a contradiction.

\item[(133)]$ R_1:h_2h_3h_2^{-1}h_3^{-1}(h_3^{-1}h_2)^2=1$, $R_2:(h_2h_3^{-2})^2h_2h_3=1$, $R_3:h_2h_3h_2^{-1}h_3^{-1}h_2h_3^{-2}h_2=1$,\\$R_4:(h_2h_3^{-1})^3(h_2^{-1}h_3)^2=1$:\\
By interchanging $h_2$ and $h_3$ in (61) and with the same discussion, there is a contradiction.

\item[(134)]$ R_1:h_2h_3h_2^{-1}h_3^{-1}(h_3^{-1}h_2)^2=1$, $R_2:(h_2h_3^{-2})^2h_2h_3=1$, $R_3:h_2h_3^{-1}h_2^{-1}h_3^2h_2^{-2}h_3=1$,\\$R_4:h_2h_3^{-1}h_2^{-1}h_3h_2^{-2}h_3^2=1$:\\
By interchanging $h_2$ and $h_3$ in (56) and with the same discussion, there is a contradiction.

\item[(135)]$ R_1:h_2h_3h_2^{-1}h_3^{-1}(h_3^{-1}h_2)^2=1$, $R_2:(h_2h_3^{-2})^2h_2h_3=1$, $R_3:h_2h_3^{-2}h_2^{-1}h_3^2h_2^{-1}h_3=1$,\\$R_4:(h_2h_3^{-1})^2(h_2^{-1}h_3)^3=1$:\\
By interchanging $h_2$ and $h_3$ in (53) and with the same discussion, there is a contradiction.

\item[(136)]$ R_1:h_2h_3h_2^{-1}h_3^{-1}(h_3^{-1}h_2)^2=1$, $R_2:(h_2h_3^{-2})^2h_2h_3=1$, $R_3:(h_2h_3^{-2})^2h_2h_3=1$,\\$R_4:(h_2h_3^{-1})^2(h_2^{-1}h_3)^3=1$:\\
By interchanging $h_2$ and $h_3$ in (54) and with the same discussion, there is a contradiction.

\item[(137)]$ R_1:h_2h_3h_2^{-1}h_3^{-1}(h_3^{-1}h_2)^2=1$, $R_2:(h_2h_3^{-2})^2h_2h_3=1$, $R_3:(h_2h_3^{-2})^2h_2h_3=1$,\\$R_4:(h_2h_3^{-1})^3(h_2^{-1}h_3)^2=1$:\\
By interchanging $h_2$ and $h_3$ in (55) and with the same discussion, there is a contradiction.

\item[(138)]$ R_1:h_2h_3h_2^{-1}h_3^{-1}(h_2^{-1}h_3)^2=1$, $R_2:h_2h_3^{-2}h_2h_3^{-1}h_2^{-1}h_3^2=1$, $R_3:h_2^2h_3^{-2}h_2h_3^{-1}h_2^{-1}h_3=1$,\\$R_4:h_2h_3^{-1}h_2(h_3^{-1}h_2h_3^{-1})^2h_2=1$:\\
By interchanging $h_2$ and $h_3$ in (72) and with the same discussion, there is a contradiction.

\item[(139)]$ R_1:h_2h_3h_2^{-1}h_3^{-1}(h_2^{-1}h_3)^2=1$, $R_2:h_2h_3^{-2}h_2h_3^{-1}h_2^{-1}h_3^2=1$, $R_3:h_2(h_2h_3^{-2})^2h_2=1$,\\$R_4:h_2h_3^{-1}h_2^{-1}h_3^2h_2^{-2}h_3=1$:\\
By interchanging $h_2$ and $h_3$ in (69) and with the same discussion, there is a contradiction.

\item[(140)]$ R_1:h_2h_3h_2^{-1}h_3^{-1}(h_2^{-1}h_3)^2=1$, $R_2:h_2h_3^{-2}h_2h_3^{-1}h_2^{-1}h_3^2=1$, $R_3:h_2h_3h_2^{-1}h_3^{-1}h_2h_3^{-1}h_2^{-1}h_3=1$,\\$R_4:h_2(h_3^{-1}h_2h_3^{-1})^2h_2h_3^{-1}h_2=1$:\\
By interchanging $h_2$ and $h_3$ in (70) and with the same discussion, there is a contradiction.

\item[(141)]$ R_1:h_2h_3h_2^{-1}h_3^{-1}(h_2^{-1}h_3)^2=1$, $R_2:h_2h_3^{-2}h_2h_3^{-1}h_2^{-1}h_3^2=1$, $R_3:h_2h_3h_2^{-1}h_3^{-1}h_2h_3^{-1}h_2^{-1}h_3=1$,\\$R_4:h_2h_3^{-1}h_2(h_3^{-1}h_2h_3^{-1})^2h_2=1$:\\
By interchanging $h_2$ and $h_3$ in (71) and with the same discussion, there is a contradiction.

\item[(142)]$ R_1:h_2h_3h_2^{-1}h_3^{-1}(h_2^{-1}h_3)^2=1$, $R_2:h_2h_3^{-2}h_2h_3^{-1}h_2^{-1}h_3^2=1$, $R_3:h_2h_3^{-1}h_2^{-1}h_3^2h_2^{-1}h_3^{-1}h_2=1$,\\$R_4:h_2h_3^{-2}h_2h_3^{-1}h_2^{-1}h_3^2=1$:\\
By interchanging $h_2$ and $h_3$ in (67) and with the same discussion, there is a contradiction.

\item[(143)]$ R_1:h_2h_3h_2^{-1}h_3^{-1}(h_2^{-1}h_3)^2=1$, $R_2:h_2h_3^{-2}h_2h_3^{-1}h_2^{-1}h_3^2=1$, $R_3:h_2h_3^{-1}h_2^{-1}h_3^2h_2^{-1}h_3^{-1}h_2=1$,\\$R_4:(h_2h_3^{-2})^2h_2h_3=1$:\\
By interchanging $h_2$ and $h_3$ in (66) and with the same discussion, there is a contradiction.

\item[(144)]$ R_1:h_2h_3h_2^{-1}h_3^{-1}(h_2^{-1}h_3)^2=1$, $R_2:h_2h_3^{-2}h_2h_3^{-1}h_2^{-1}h_3^2=1$, $R_3:h_2h_3^{-3}h_2h_3h_2^{-1}h_3=1$,\\$R_4:h_2(h_3^{-1}h_2h_3^{-1})^2h_2h_3^{-1}h_2=1$:\\
By interchanging $h_2$ and $h_3$ in (64) and with the same discussion, there is a contradiction.

\item[(145)]$ R_1:h_2h_3h_2^{-1}h_3^{-1}(h_2^{-1}h_3)^2=1$, $R_2:h_2h_3^{-2}h_2h_3^{-1}h_2^{-1}h_3^2=1$, $R_3:h_2h_3^{-3}h_2h_3h_2^{-1}h_3=1$,\\$R_4:h_2h_3^{-1}h_2(h_3^{-1}h_2h_3^{-1})^2h_2=1$:\\
By interchanging $h_2$ and $h_3$ in (65) and with the same discussion, there is a contradiction.

\item[(146)]$ R_1:h_2h_3h_2^{-1}h_3^{-1}(h_2^{-1}h_3)^2=1$, $R_2:h_2h_3^{-2}h_2h_3^{-1}h_2^{-1}h_3^2=1$, $R_3:h_2h_3^{-2}h_2h_3^{-1}h_2^{-1}h_3^2=1$,\\$R_4:h_2h_3^{-1}h_2(h_3^{-1}h_2h_3^{-1})^2h_2=1$:\\
By interchanging $h_2$ and $h_3$ in (68) and with the same discussion, there is a contradiction.

\item[(147)]$ R_1:h_2h_3^{-1}h_2^{-1}h_3^2h_2^{-1}h_3^{-1}h_2=1$, $R_2:h_2h_3^{-1}h_2^{-1}h_3h_2^{-1}h_3^{-1}h_2h_3=1$, $R_3:h_2h_3h_2h_3^{-1}h_2^{-1}h_3^2=1$,\\$R_4:h_2(h_3h_2^{-1})^2h_3h_2h_3^{-1}h_2=1$:\\
By interchanging $h_2$ and $h_3$ in (73) and with the same discussion, there is a contradiction.

\item[(148)]$ R_1:h_2(h_3h_2^{-1})^2h_3^{-2}h_2=1$, $R_2:h_2(h_3^{-1}h_2h_3^{-1})^2h_3^{-1}h_2=1$, $R_3:(h_2h_3)^2h_2^{-1}h_3^2=1$,\\$R_4:h_2(h_3h_2^{-1})^2h_3^{-2}h_2=1$:\\
By interchanging $h_2$ and $h_3$ in (74) and with the same discussion, there is a contradiction.

\end{enumerate}

\subsection{$\mathbf{C_6---C_6(-C_5-)}$}
\begin{figure}[ht]
\psscalebox{0.9 0.9} 
{
\begin{pspicture}(0,-2.0985577)(4.3771152,2.0985577)
\psdots[linecolor=black, dotsize=0.4](2.5971153,1.9014423)
\psdots[linecolor=black, dotsize=0.4](2.5971153,-0.49855766)
\psdots[linecolor=black, dotsize=0.4](2.5971153,1.1014423)
\psdots[linecolor=black, dotsize=0.4](2.5971153,0.30144233)
\psdots[linecolor=black, dotsize=0.4](3.7971153,1.1014423)
\psdots[linecolor=black, dotsize=0.4](3.7971153,0.30144233)
\psdots[linecolor=black, dotsize=0.4](1.3971153,0.30144233)
\psdots[linecolor=black, dotsize=0.4](1.3971153,1.1014423)
\psline[linecolor=black, linewidth=0.04](2.5971153,1.9014423)(3.7971153,1.1014423)(3.7971153,0.30144233)(2.5971153,-0.49855766)(2.5971153,0.30144233)(2.5971153,1.1014423)(2.5971153,1.9014423)(1.3971153,1.1014423)(1.3971153,0.30144233)(2.5971153,-0.49855766)
\psdots[linecolor=black, dotsize=0.4](0.19711533,0.30144233)
\psdots[linecolor=black, dotsize=0.4](0.9971153,-0.49855766)
\psline[linecolor=black, linewidth=0.04](1.3971153,1.1014423)(0.19711533,0.30144233)(0.9971153,-0.49855766)(2.5971153,1.1014423)
\rput[bl](0.1,-2.0985577){$\mathbf{26) \ C_6---C_6(-C_5-)}$}
\end{pspicture}
}
\end{figure}
By considering the $99$ cases related to the existence of $C_6---C_6$ in the graph $K(\alpha,\beta)$ and the relations from Table \ref{tab-C5} which are not disproved, it can be seen that there are $1482$ cases for the relations of two cycles $C_6$ and a cycle $C_5$ in the graph $C_6---C_6(-C_5-)$. Using GAP \cite{gap}, we see that all groups with two generators $h_2$ and $h_3$ and three relations which are between $1358$ cases of these $1482$ cases are finite and solvable, that is a contradiction with the assumptions. So, there are just $124$ cases for the relations of these cycles which may lead to the existence of a subgraph isomorphic to the graph $C_6---C_6(-C_5-)$ in $K(\alpha,\beta)$. In the following, we show that these $124$ cases lead to contradictions and so, the graph $K(\alpha,\beta)$ contains no subgraph isomorphic to the graph $C_6---C_6(-C_5-)$.
\begin{enumerate}
\item[(1)]$ R_1:h_2^3h_3^3=1$, $R_2:h_2^2(h_3h_2^{-1})^3h_3=1$, $R_3:h_2^2h_3^{-3}h_2=1$:\\
$\Rightarrow$  $G$ has a torsion element, a contradiction.

\item[(2)]$ R_1:h_2^3h_3^3=1$, $R_2:(h_2h_3^{-1})^3h_2h_3^2=1$, $R_3:h_2^2h_3^{-3}h_2=1$:\\
By interchanging $h_2$ and $h_3$ in (1) and with the same discussion, there is a contradiction.

\item[(3)]$ R_1:h_2^3(h_3h_2^{-1})^2h_3=1$, $R_2:h_3^2(h_3h_2^{-1})^3h_3=1$, $R_3:(h_2h_3^{-1})^2h_2h_3^2=1$:\\
$\Rightarrow$  $G$ has a torsion element, a contradiction.

\item[(4)]$ R_1:h_2^2(h_2h_3^{-1})^2h_2^{-1}h_3=1$, $R_2:(h_2h_3^{-1})^2h_3^{-1}h_2h_3^2=1$, $R_3:h_2^2h_3^{-2}h_2h_3=1$:\\
$\Rightarrow \langle h_2,h_3\rangle=\langle h_2\rangle$ is abelian, a contradiction.

\item[(5)]$ R_1:h_2^2(h_2h_3^{-1})^3h_2=1$, $R_2:(h_2h_3^{-1})^2h_2h_3^3=1$, $R_3:h_2^2(h_3h_2^{-1})^2h_3=1$:\\
By interchanging $h_2$ and $h_3$ in (3) and with the same discussion, there is a contradiction.

\item[(6)]$ R_1:h_2^2h_3h_2^{-2}h_3^2=1$, $R_2:h_2h_3h_2^{-1}h_3^{-1}h_2h_3^{-1}h_2^{-1}h_3=1$, $R_3:h_2h_3^{-1}h_2(h_3h_2^{-1})^2h_3=1$:\\
$\Rightarrow$  $h_3=1$, a contradiction.

\item[(7)]$ R_1:h_2^2h_3h_2^{-2}h_3^2=1$, $R_2:h_2h_3h_2^{-1}h_3^{-1}h_2h_3^{-1}h_2^{-1}h_3=1$, $R_3:(h_2h_3^{-1})^2(h_2^{-1}h_3)^2=1$:\\
$\Rightarrow$  $G$ has a torsion element, a contradiction.

\item[(8)]$ R_1:h_2^2h_3h_2^{-2}h_3^2=1$, $R_2:h_2h_3h_2^{-1}h_3^{-1}h_2h_3^{-1}h_2^{-1}h_3=1$, $R_3:h_2^2h_3h_2^{-2}h_3=1$:\\
$\Rightarrow$  $h_3=1$, a contradiction.

\item[(9)]$ R_1:h_2^2h_3h_2^{-2}h_3^2=1$, $R_2:h_2h_3h_2^{-1}h_3^{-1}h_2h_3^{-1}h_2^{-1}h_3=1$, $R_3:h_2^2h_3^{-1}h_2^{-2}h_3=1$:\\
$\Rightarrow$  $G$ has a torsion element, a contradiction.

\item[(10)]$ R_1:h_2^2h_3h_2^{-1}(h_2^{-1}h_3)^2=1$, $R_2:h_2h_3h_2^{-1}h_3^{-1}h_2h_3^{-1}h_2^{-1}h_3=1$, $R_3:(h_2h_3^{-1})^2(h_2^{-1}h_3)^2=1$:\\
$\Rightarrow\langle h_2,h_3\rangle\cong BS(1,1)$ is solvable, a contradiction.

\item[(11)]$ R_1:h_2^2h_3h_2^{-1}(h_2^{-1}h_3)^2=1$, $R_2:h_2h_3h_2^{-1}h_3^{-1}h_2h_3^{-1}h_2^{-1}h_3=1$, $R_3:h_2^2h_3^{-1}h_2^{-2}h_3=1$:\\
$\Rightarrow \langle h_2,h_3\rangle=\langle h_3\rangle$ is abelian, a contradiction.

\item[(12)]$ R_1:h_2^2h_3h_2^{-1}(h_2^{-1}h_3)^2=1$, $R_2:h_2h_3h_2^{-1}h_3^{-1}h_2h_3^{-1}h_2^{-1}h_3=1$, $R_3:h_2h_3^2h_2^{-2}h_3=1$:\\
$\Rightarrow \langle h_2,h_3\rangle=\langle h_3\rangle$ is abelian, a contradiction.

\item[(13)]$ R_1:h_2^2h_3h_2^{-1}h_3^{-2}h_2=1$, $R_2:h_2^2h_3^{-1}h_2^{-1}h_3^{-1}h_2h_3=1$, $R_3:h_2^2h_3^{-1}h_2^{-2}h_3=1$:\\
$\Rightarrow\langle h_2,h_3\rangle\cong BS(1,1)$ is solvable, a contradiction.

\item[(14)]$ R_1:h_2^2h_3h_2^{-1}h_3^{-2}h_2=1$, $R_2:h_2^2h_3^{-1}h_2^{-1}h_3^{-1}h_2h_3=1$, $R_3:h_2(h_2h_3^{-1})^2h_2^{-1}h_3=1$:\\
$\Rightarrow\langle h_2,h_3\rangle\cong BS(1,1)$ is solvable, a contradiction.

\item[(15)]$ R_1:h_2^2h_3h_2^{-1}h_3^{-2}h_2=1$, $R_2:h_2^2h_3^{-1}h_2^{-1}h_3^{-1}h_2h_3=1$, $R_3:h_2h_3h_2^{-1}(h_3^{-1}h_2)^2=1$:\\
$\Rightarrow\langle h_2,h_3\rangle\cong BS(1,1)$ is solvable, a contradiction.

\item[(16)]$ R_1:h_2^2h_3h_2^{-1}h_3^{-2}h_2=1$, $R_2:h_2^2h_3^{-1}h_2^{-1}h_3^{-1}h_2h_3=1$, $R_3:(h_2h_3^{-1})^2(h_2^{-1}h_3)^2=1$:\\
$\Rightarrow \langle h_2,h_3\rangle=\langle h_2\rangle$ is abelian, a contradiction.

\item[(17)]$ R_1:h_2^2h_3h_2^{-1}h_3^{-2}h_2=1$, $R_2:h_2^2h_3^{-1}h_2^{-1}h_3^{-1}h_2h_3=1$, $R_3:h_2^3h_3^{-2}h_2=1$:\\
$\Rightarrow \langle h_2,h_3\rangle=\langle h_2\rangle$ is abelian, a contradiction.

\item[(18)]$ R_1:h_2^2h_3h_2^{-1}h_3^{-2}h_2=1$, $R_2:h_2h_3h_2h_3^{-1}(h_3^{-1}h_2)^2=1$, $R_3:h_2h_3^{-2}h_2^{-1}h_3^2=1$:\\
$\Rightarrow \langle h_2,h_3\rangle=\langle h_2\rangle$ is abelian, a contradiction.

\item[(19)]$ R_1:h_2^2h_3h_2^{-1}h_3^{-2}h_2=1$, $R_2:h_2h_3^{-2}h_2^2h_3^{-1}h_2h_3=1$, $R_3:h_2(h_2h_3^{-1})^2h_2^{-1}h_3=1$:\\
$\Rightarrow \langle h_2,h_3\rangle=\langle h_2\rangle$ is abelian, a contradiction.

\item[(20)]$ R_1:h_2^2h_3h_2^{-1}h_3^{-2}h_2=1$, $R_2:h_2h_3^{-2}h_2^2h_3^{-1}h_2h_3=1$, $R_3:h_2h_3^{-1}h_2^{-1}h_3h_2h_3^{-1}h_2=1$:\\
$\Rightarrow \langle h_2,h_3\rangle=\langle h_2\rangle$ is abelian, a contradiction.

\item[(21)]$ R_1:h_2^2(h_3h_2^{-1}h_3)^2=1$, $R_2:(h_2h_3)^2(h_2^{-1}h_3)^2=1$, $R_3:h_2^2h_3^{-1}h_2^{-1}h_3^{-1}h_2=1$:\\
$\Rightarrow \langle h_2,h_3\rangle=\langle h_2\rangle$ is abelian, a contradiction.

\item[(22)]$ R_1:h_2^2(h_3h_2^{-1}h_3)^2=1$, $R_2:(h_2h_3)^2(h_2^{-1}h_3)^2=1$, $R_3:h_2h_3^{-2}h_2^{-1}h_3^2=1$:\\
$\Rightarrow\langle h_2,h_3\rangle\cong BS(1,1)$ is solvable, a contradiction.

\item[(23)]$ R_1:h_2^2(h_3h_2^{-1}h_3)^2=1$, $R_2:h_2h_3h_2^{-1}h_3^{-1}h_2h_3^{-1}h_2^{-1}h_3=1$, $R_3:h_2h_3^{-2}h_2^{-1}h_3^2=1$:\\
$\Rightarrow\langle h_2,h_3\rangle\cong BS(1,1)$ is solvable, a contradiction.

\item[(24)]$ R_1:h_2^2(h_3h_2^{-1}h_3)^2=1$, $R_2:h_2(h_3h_2^{-1})^2h_3^{-1}h_2h_3=1$, $R_3:h_2h_3^{-2}h_2^{-1}h_3^2=1$:\\
$\Rightarrow$ $G$ is solvable, a contradiction.

\item[(25)]$ R_1:h_2^2(h_3h_2^{-1})^2h_2^{-1}h_3=1$, $R_2:(h_2h_3^{-1})^2h_3^{-2}h_2^{-1}h_3=1$, $R_3:h_2h_3^2h_2^{-2}h_3=1$:\\
By interchanging $h_2$ and $h_3$ in (4) and with the same discussion, there is a contradiction.

\item[(26)]$ R_1:h_2^2h_3^{-1}h_2^{-2}h_3^2=1$, $R_2:h_2h_3h_2^{-1}h_3^{-1}h_2h_3^{-1}h_2^{-1}h_3=1$, $R_3:h_2^2h_3h_2^{-2}h_3=1$:\\
$\Rightarrow$  $G$ has a torsion element, a contradiction.

\item[(27)]$ R_1:h_2^2h_3^{-1}h_2^{-2}h_3^2=1$, $R_2:h_2h_3h_2^{-1}h_3^{-1}h_2h_3^{-1}h_2^{-1}h_3=1$, $R_3:h_2^2h_3^{-1}h_2^{-2}h_3=1$:\\
$\Rightarrow$  $h_3=1$, a contradiction.

\item[(28)]$ R_1:h_2^2h_3^{-1}h_2^{-2}h_3^2=1$, $R_2:h_2h_3h_2^{-1}h_3^{-1}h_2h_3^{-1}h_2^{-1}h_3=1$, $R_3:h_2h_3^{-1}h_2(h_3h_2^{-1})^2h_3=1$:\\
$\Rightarrow$  $G$ has a torsion element, a contradiction.

\item[(29)]$ R_1:h_2^2h_3^{-1}h_2^{-2}h_3^2=1$, $R_2:h_2h_3h_2^{-1}h_3^{-1}h_2h_3^{-1}h_2^{-1}h_3=1$, $R_3:(h_2h_3^{-1})^2(h_2^{-1}h_3)^2=1$:\\
$\Rightarrow$  $h_3=1$, a contradiction.

\item[(30)]$ R_1:h_2^2h_3^{-1}h_2^{-1}h_3^{-2}h_2=1$, $R_2:h_2^2h_3^{-2}(h_2^{-1}h_3)^2=1$, $R_3:h_2h_3h_2^{-1}(h_2^{-1}h_3)^2=1$:\\
$\Rightarrow\langle h_2,h_3\rangle\cong BS(1,1)$ is solvable, a contradiction.

\item[(31)]$ R_1:h_2^2h_3^{-1}h_2^{-1}h_3^2h_2^{-1}h_3=1$, $R_2:h_2h_3h_2h_3^{-1}(h_2^{-1}h_3)^2=1$, $R_3:h_2h_3^2h_2^{-1}h_3^2=1$:\\
$\Rightarrow$  $G$ has a torsion element, a contradiction.

\item[(32)]$ R_1:h_2^2h_3^{-1}h_2^{-1}h_3^2h_2^{-1}h_3=1$, $R_2:h_2h_3h_2^{-1}h_3^{-1}h_2h_3^{-1}h_2^{-1}h_3=1$, $R_3:h_2h_3^2h_2^{-1}h_3^2=1$:\\
$\Rightarrow$  $G$ has a torsion element, a contradiction.

\item[(33)]$ R_1:h_2^2h_3^{-2}(h_2^{-1}h_3)^2=1$, $R_2:h_2h_3h_2h_3^{-3}h_2=1$, $R_3:h_2h_3^{-1}(h_3^{-1}h_2)^2h_3=1$:\\
By interchanging $h_2$ and $h_3$ in (30) and with the same discussion, there is a contradiction.

\item[(34)]$ R_1:h_2(h_2h_3^{-2})^2h_2=1$, $R_2:h_2h_3h_2^{-1}h_3^2h_2h_3^{-1}h_2=1$, $R_3:h_2h_3^2h_2^{-1}h_3^2=1$:\\
$\Rightarrow\langle h_2,h_3\rangle\cong BS(1,-1)$ is solvable, a contradiction.

\item[(35)]$ R_1:h_2(h_2h_3^{-2})^2h_2=1$, $R_2:h_2h_3^{-1}h_2^{-1}h_3^2h_2^{-1}h_3^{-1}h_2=1$, $R_3:h_2^2h_3h_2^{-2}h_3=1$:\\
$\Rightarrow \langle h_2,h_3\rangle=\langle h_3\rangle$ is abelian, a contradiction.

\item[(36)]$ R_1:h_2^2h_3^{-1}h_2h_3h_2h_3^{-1}h_2=1$, $R_2:h_2(h_2h_3^{-1})^2h_2^2h_3=1$, $R_3:h_2h_3^{-2}h_2^{-1}h_3^2=1$:\\
$\Rightarrow\langle h_2,h_3\rangle\cong BS(1,1)$ is solvable, a contradiction.

\item[(37)]$ R_1:h_2^2h_3^{-1}h_2h_3h_2^{-1}h_3^2=1$, $R_2:h_2h_3h_2^{-1}h_3h_2h_3^{-1}h_2h_3=1$, $R_3:h_2^3h_3h_2^{-1}h_3=1$:\\
$\Rightarrow \langle h_2,h_3\rangle=\langle h_2\rangle$ is abelian, a contradiction.

\item[(38)]$ R_1:h_2^2h_3^{-1}h_2h_3h_2^{-1}h_3^2=1$, $R_2:h_2h_3h_2^{-1}h_3h_2h_3^{-1}h_2h_3=1$, $R_3:h_2^2h_3^{-1}h_2^{-2}h_3=1$:\\
$\Rightarrow \langle h_2,h_3\rangle=\langle h_2\rangle$ is abelian, a contradiction.

\item[(39)]$ R_1:h_2(h_2h_3^{-1})^2h_3^{-1}h_2h_3=1$, $R_2:h_2(h_3h_2^{-1})^2h_3^{-2}h_2=1$, $R_3:h_2^3h_3^{-2}h_2=1$:\\
$\Rightarrow\langle h_2,h_3\rangle\cong BS(1,1)$ is solvable, a contradiction.

\item[(40)]$ R_1:h_2^2(h_3^{-1}h_2h_3^{-1})^2h_2=1$, $R_2:h_2(h_3h_2^{-1})^2h_3^{-2}h_2=1$, $R_3:h_2h_3^{-2}h_2^{-1}h_3^2=1$:\\
$\Rightarrow\langle h_2,h_3\rangle\cong BS(1,1)$ is solvable, a contradiction.

\item[(41)]$ R_1:(h_2h_3)^2h_2h_3^{-1}h_2=1$, $R_2:(h_2h_3)^2h_2h_3^{-1}h_2=1$, $R_3:h_2h_3^{-2}h_2^{-1}h_3^2=1$:\\
$\Rightarrow\langle h_2,h_3\rangle\cong BS(1,1)$ is solvable, a contradiction.

\item[(42)]$ R_1:(h_2h_3)^2h_2h_3^{-1}h_2=1$, $R_2:h_2(h_3h_2^{-1})^2h_3^{-2}h_2=1$, $R_3:h_2^2h_3^{-1}h_2^{-2}h_3=1$:\\
$\Rightarrow \langle h_2,h_3\rangle=\langle h_2\rangle$ is abelian, a contradiction.

\item[(43)]$ R_1:(h_2h_3)^2h_2^{-1}h_3^2=1$, $R_2:(h_2h_3)^2h_2^{-1}h_3^2=1$, $R_3:h_2^2h_3^{-1}h_2^{-2}h_3=1$:\\
By interchanging $h_2$ and $h_3$ in (41) and with the same discussion, there is a contradiction.

\item[(44)]$ R_1:(h_2h_3)^2h_2^{-1}h_3^2=1$, $R_2:h_2(h_3h_2^{-1})^2h_3^{-2}h_2=1$, $R_3:h_2h_3^{-2}h_2^{-1}h_3^2=1$:\\
By interchanging $h_2$ and $h_3$ in (42) and with the same discussion, there is a contradiction.

\item[(45)]$ R_1:h_2h_3h_2h_3^{-1}(h_3^{-1}h_2)^2=1$, $R_2:h_2h_3^{-2}h_2h_3h_2h_3^{-1}h_2=1$, $R_3:h_2^3h_3^{-2}h_2=1$:\\
$\Rightarrow\langle h_2,h_3\rangle\cong BS(1,1)$ is solvable, a contradiction.

\item[(46)]$ R_1:h_2h_3h_2h_3^{-1}(h_3^{-1}h_2)^2=1$, $R_2:h_2h_3^{-2}h_2h_3h_2h_3^{-1}h_2=1$, $R_3:h_2^2h_3^{-1}h_2^{-2}h_3=1$:\\
$\Rightarrow\langle h_2,h_3\rangle\cong BS(1,1)$ is solvable, a contradiction.

\item[(47)]$ R_1:h_2h_3(h_2h_3^{-1}h_2)^2=1$, $R_2:h_2h_3(h_2h_3^{-1}h_2)^2=1$, $R_3:h_2^2h_3^{-1}h_2^{-2}h_3=1$:\\
$\Rightarrow \langle h_2,h_3\rangle=\langle h_2\rangle$ is abelian, a contradiction.

\item[(48)]$ R_1:h_2h_3h_2h_3^{-1}h_2h_3h_2^{-1}h_3=1$, $R_2:h_2h_3^2h_2^{-1}h_3h_2h_3^{-1}h_2=1$, $R_3:h_2h_3^{-2}h_2^{-1}h_3^2=1$:\\
By interchanging $h_2$ and $h_3$ in (38) and with the same discussion, there is a contradiction.

\item[(49)]$ R_1:h_2h_3h_2h_3^{-1}h_2h_3h_2^{-1}h_3=1$, $R_2:h_2h_3^2h_2^{-1}h_3h_2h_3^{-1}h_2=1$, $R_3:h_2h_3^{-1}h_2h_3^3=1$:\\
By interchanging $h_2$ and $h_3$ in (37) and with the same discussion, there is a contradiction.

\item[(50)]$ R_1:h_2h_3h_2h_3^{-1}h_2h_3h_2^{-1}h_3=1$, $R_2:h_2(h_3h_2^{-1})^2h_3^{-2}h_2=1$, $R_3:h_2^3h_3h_2^{-1}h_3=1$:\\
$\Rightarrow\langle h_2,h_3\rangle\cong BS(1,-1)$ is solvable, a contradiction.

\item[(51)]$ R_1:h_2h_3h_2h_3^{-1}h_2h_3h_2^{-1}h_3=1$, $R_2:h_2(h_3h_2^{-1})^2h_3^{-2}h_2=1$, $R_3:h_2^2h_3^{-1}h_2^{-2}h_3=1$:\\
$\Rightarrow\langle h_2,h_3\rangle\cong BS(1,1)$ is solvable, a contradiction.

\item[(52)]$ R_1:h_2h_3(h_2h_3^{-1})^2h_2^{-1}h_3=1$, $R_2:h_2h_3^{-1}h_2^{-1}h_3^2h_2h_3^{-1}h_2=1$, $R_3:h_2^2h_3^{-1}h_2^{-2}h_3=1$:\\
$\Rightarrow$  $G$ has a torsion element, a contradiction.

\item[(53)]$ R_1:h_2h_3(h_2h_3^{-1})^2h_2^{-1}h_3=1$, $R_2:(h_2h_3^{-1}h_2)^2h_3^2=1$, $R_3:h_2^2h_3^{-1}h_2^{-2}h_3=1$:\\
By interchanging $h_2$ and $h_3$ in (24) and with the same discussion, there is a contradiction.

\item[(54)]$ R_1:h_2h_3(h_2h_3^{-1})^2h_2h_3=1$, $R_2:(h_2h_3^{-1}h_2)^2h_3^2=1$, $R_3:h_2^2h_3^{-1}h_2^{-2}h_3=1$:\\
By interchanging $h_2$ and $h_3$ in (22) and with the same discussion, there is a contradiction.

\item[(55)]$ R_1:h_2h_3(h_2h_3^{-1})^2h_2h_3=1$, $R_2:(h_2h_3^{-1}h_2)^2h_3^2=1$, $R_3:h_2h_3^{-3}h_2h_3=1$:\\
By interchanging $h_2$ and $h_3$ in (21) and with the same discussion, there is a contradiction.

\item[(56)]$ R_1:h_2h_3^2h_2h_3^{-2}h_2=1$, $R_2:h_2h_3^{-1}h_2^{-1}h_3h_2^{-1}h_3^{-1}h_2h_3=1$, $R_3:(h_2h_3^{-1})^2(h_2^{-1}h_3)^2=1$:\\
By interchanging $h_2$ and $h_3$ in (7) and with the same discussion, there is a contradiction.

\item[(57)]$ R_1:h_2h_3^2h_2h_3^{-2}h_2=1$, $R_2:h_2h_3^{-1}h_2^{-1}h_3h_2^{-1}h_3^{-1}h_2h_3=1$, $R_3:(h_2h_3^{-1})^2h_2h_3h_2^{-1}h_3=1$:\\
By interchanging $h_2$ and $h_3$ in (6) and with the same discussion, there is a contradiction.

\item[(58)]$ R_1:h_2h_3^2h_2h_3^{-2}h_2=1$, $R_2:h_2h_3^{-1}h_2^{-1}h_3h_2^{-1}h_3^{-1}h_2h_3=1$, $R_3:h_2h_3^{-2}h_2^{-1}h_3^2=1$:\\
By interchanging $h_2$ and $h_3$ in (9) and with the same discussion, there is a contradiction.

\item[(59)]$ R_1:h_2h_3^2h_2h_3^{-2}h_2=1$, $R_2:h_2h_3^{-1}h_2^{-1}h_3h_2^{-1}h_3^{-1}h_2h_3=1$, $R_3:h_2h_3^{-2}h_2h_3^2=1$:\\
By interchanging $h_2$ and $h_3$ in (8) and with the same discussion, there is a contradiction.

\item[(60)]$ R_1:h_2h_3^2h_2^{-1}h_3^{-1}h_2^{-1}h_3=1$, $R_2:h_2h_3h_2^{-1}h_3^{-3}h_2=1$, $R_3:h_2h_3^{-4}h_2=1$:\\
By interchanging $h_2$ and $h_3$ in (17) and with the same discussion, there is a contradiction.

\item[(61)]$ R_1:h_2h_3^2h_2^{-1}h_3^{-1}h_2^{-1}h_3=1$, $R_2:h_2h_3h_2^{-1}h_3^{-3}h_2=1$, $R_3:h_2h_3^{-1}(h_2^{-1}h_3)^2h_3=1$:\\
By interchanging $h_2$ and $h_3$ in (15) and with the same discussion, there is a contradiction.

\item[(62)]$ R_1:h_2h_3^2h_2^{-1}h_3^{-1}h_2^{-1}h_3=1$, $R_2:h_2h_3h_2^{-1}h_3^{-3}h_2=1$, $R_3:h_2h_3^{-2}h_2^{-1}h_3^2=1$:\\
By interchanging $h_2$ and $h_3$ in (13) and with the same discussion, there is a contradiction.

\item[(63)]$ R_1:h_2h_3^2h_2^{-1}h_3^{-1}h_2^{-1}h_3=1$, $R_2:h_2h_3h_2^{-1}h_3^{-3}h_2=1$, $R_3:(h_2h_3^{-1})^2(h_2^{-1}h_3)^2=1$:\\
By interchanging $h_2$ and $h_3$ in (16) and with the same discussion, there is a contradiction.

\item[(64)]$ R_1:h_2h_3^2h_2^{-1}h_3^{-1}h_2^{-1}h_3=1$, $R_2:h_2h_3h_2^{-1}h_3^{-3}h_2=1$, $R_3:(h_2h_3^{-1})^2h_3^{-1}h_2^{-1}h_3=1$:\\
By interchanging $h_2$ and $h_3$ in (14) and with the same discussion, there is a contradiction.

\item[(65)]$ R_1:h_2h_3^2h_2^{-1}h_3^{-2}h_2=1$, $R_2:h_2h_3^{-1}h_2^{-1}h_3h_2^{-1}h_3^{-1}h_2h_3=1$, $R_3:h_2h_3^{-2}h_2^{-1}h_3^2=1$:\\
By interchanging $h_2$ and $h_3$ in (27) and with the same discussion, there is a contradiction.

\item[(66)]$ R_1:h_2h_3^2h_2^{-1}h_3^{-2}h_2=1$, $R_2:h_2h_3^{-1}h_2^{-1}h_3h_2^{-1}h_3^{-1}h_2h_3=1$, $R_3:h_2h_3^{-2}h_2h_3^2=1$:\\
By interchanging $h_2$ and $h_3$ in (26) and with the same discussion, there is a contradiction.

\item[(67)]$ R_1:h_2h_3^2h_2^{-1}h_3^{-2}h_2=1$, $R_2:h_2h_3^{-1}h_2^{-1}h_3h_2^{-1}h_3^{-1}h_2h_3=1$, $R_3:(h_2h_3^{-1})^2(h_2^{-1}h_3)^2=1$:\\
By interchanging $h_2$ and $h_3$ in (29) and with the same discussion, there is a contradiction.

\item[(68)]$ R_1:h_2h_3^2h_2^{-1}h_3^{-2}h_2=1$, $R_2:h_2h_3^{-1}h_2^{-1}h_3h_2^{-1}h_3^{-1}h_2h_3=1$, $R_3:(h_2h_3^{-1})^2h_2h_3h_2^{-1}h_3=1$:\\
By interchanging $h_2$ and $h_3$ in (28) and with the same discussion, there is a contradiction.

\item[(69)]$ R_1:h_2h_3(h_3h_2^{-1})^2h_2^{-1}h_3=1$, $R_2:h_2(h_3h_2^{-1})^2h_3^{-2}h_2=1$, $R_3:h_2h_3^{-4}h_2=1$:\\
By interchanging $h_2$ and $h_3$ in (39) and with the same discussion, there is a contradiction.

\item[(70)]$ R_1:h_2h_3(h_3h_2^{-1})^2h_3^2=1$, $R_2:h_2h_3h_2^{-1}h_3^3h_2^{-1}h_3=1$, $R_3:h_2^2h_3^{-1}h_2^{-2}h_3=1$:\\
By interchanging $h_2$ and $h_3$ in (36) and with the same discussion, there is a contradiction.

\item[(71)]$ R_1:h_2h_3h_2^{-2}h_3^2h_2^{-1}h_3=1$, $R_2:h_2h_3h_2^{-1}h_3^{-3}h_2=1$, $R_3:h_2h_3^{-2}h_2h_3^{-1}h_2^{-1}h_3=1$:\\
By interchanging $h_2$ and $h_3$ in (20) and with the same discussion, there is a contradiction.

\item[(72)]$ R_1:h_2h_3h_2^{-2}h_3^2h_2^{-1}h_3=1$, $R_2:h_2h_3h_2^{-1}h_3^{-3}h_2=1$, $R_3:(h_2h_3^{-1})^2h_3^{-1}h_2^{-1}h_3=1$:\\
By interchanging $h_2$ and $h_3$ in (19) and with the same discussion, there is a contradiction.

\item[(73)]$ R_1:h_2h_3h_2^{-1}(h_2^{-1}h_3)^2h_3=1$, $R_2:h_2h_3h_2^{-1}h_3^{-3}h_2=1$, $R_3:h_2^2h_3^{-1}h_2^{-2}h_3=1$:\\
By interchanging $h_2$ and $h_3$ in (18) and with the same discussion, there is a contradiction.

\item[(74)]$ R_1:h_2h_3h_2^{-1}(h_2^{-1}h_3)^2h_3=1$, $R_2:h_2h_3h_2^{-1}h_3^2h_2^{-2}h_3=1$, $R_3:h_2h_3^{-2}h_2^{-1}h_3^2=1$:\\
By interchanging $h_2$ and $h_3$ in (46) and with the same discussion, there is a contradiction.

\item[(75)]$ R_1:h_2h_3h_2^{-1}(h_2^{-1}h_3)^2h_3=1$, $R_2:h_2h_3h_2^{-1}h_3^2h_2^{-2}h_3=1$, $R_3:h_2h_3^{-4}h_2=1$:\\
By interchanging $h_2$ and $h_3$ in (45) and with the same discussion, there is a contradiction.

\item[(76)]$ R_1:h_2h_3h_2^{-1}h_3^{-1}h_2h_3^{-1}h_2^{-1}h_3=1$, $R_2:h_2h_3h_2^{-1}h_3^2h_2^{-1}h_3^{-1}h_2=1$, $R_3:h_2h_3^{-2}h_2^{-1}h_3^2=1$:\\
$\Rightarrow$  $G$ has a torsion element, a contradiction.

\item[(77)]$ R_1:h_2h_3h_2^{-1}h_3^{-1}h_2h_3^{-1}h_2^{-1}h_3=1$, $R_2:h_2h_3^{-1}h_2^{-1}h_3^2h_2^{-1}h_3^{-1}h_2=1$, $R_3:h_2^3h_3^2=1$:\\
$\Rightarrow$ $G$ is solvable, a contradiction.

\item[(78)]$ R_1:h_2h_3h_2^{-1}h_3^{-1}h_2h_3^{-1}h_2^{-1}h_3=1$, $R_2:h_2h_3^{-1}h_2^{-1}h_3^2h_2^{-1}h_3^{-1}h_2=1$, $R_3:h_2^2h_3^3=1$:\\
$\Rightarrow$ $G$ is solvable, a contradiction.

\item[(79)]$ R_1:h_2h_3h_2^{-1}h_3^{-1}h_2h_3^{-1}h_2^{-1}h_3=1$, $R_2:h_2h_3^{-1}h_2^{-1}h_3^2h_2^{-1}h_3^{-1}h_2=1$, $R_3:h_2^2h_3^2h_2^{-1}h_3=1$:\\
$\Rightarrow \langle h_2,h_3\rangle=\langle h_3\rangle$ is abelian, a contradiction.

\item[(80)]$ R_1:h_2h_3h_2^{-1}h_3^{-1}h_2h_3^{-1}h_2^{-1}h_3=1$, $R_2:h_2h_3^{-1}h_2^{-1}h_3^2h_2^{-1}h_3^{-1}h_2=1$, $R_3:h_2^2h_3h_2^{-1}h_3^2=1$:\\
$\Rightarrow \langle h_2,h_3\rangle=\langle h_3\rangle$ is abelian, a contradiction.

\item[(81)]$ R_1:h_2h_3h_2^{-1}h_3^{-1}h_2h_3^{-1}h_2^{-1}h_3=1$, $R_2:h_2h_3^{-1}h_2^{-1}h_3^2h_2^{-1}h_3^{-1}h_2=1$, $R_3:h_2^2h_3^{-1}h_2^2h_3=1$:\\
$\Rightarrow\langle h_2,h_3\rangle\cong BS(-1,1)$ is solvable, a contradiction.

\item[(82)]$ R_1:h_2h_3h_2^{-1}h_3^{-1}h_2h_3^{-1}h_2^{-1}h_3=1$, $R_2:h_2h_3^{-1}h_2^{-1}h_3^2h_2^{-1}h_3^{-1}h_2=1$, $R_3:h_2h_3h_2h_3^{-2}h_2=1$:\\
$\Rightarrow \langle h_2,h_3\rangle=\langle h_2\rangle$ is abelian, a contradiction.

\item[(83)]$ R_1:h_2h_3h_2^{-1}h_3^{-1}h_2h_3^{-1}h_2^{-1}h_3=1$, $R_2:h_2h_3^{-1}h_2^{-1}h_3^2h_2^{-1}h_3^{-1}h_2=1$, $R_3:h_2h_3(h_2h_3^{-1})^2h_2=1$:\\
$\Rightarrow$ $G$ is solvable, a contradiction.

\item[(84)]$ R_1:h_2h_3h_2^{-1}h_3^{-1}h_2h_3^{-1}h_2^{-1}h_3=1$, $R_2:h_2h_3^{-1}h_2^{-1}h_3^2h_2^{-1}h_3^{-1}h_2=1$, $R_3:h_2h_3^2h_2h_3^{-1}h_2=1$:\\
$\Rightarrow\langle h_2,h_3\rangle\cong BS(2,1)$ is solvable, a contradiction.

\item[(85)]$ R_1:h_2h_3h_2^{-1}h_3^{-1}h_2h_3^{-1}h_2^{-1}h_3=1$, $R_2:h_2h_3^{-1}h_2^{-1}h_3^2h_2^{-1}h_3^{-1}h_2=1$, $R_3:h_2h_3^2h_2^{-1}h_3^{-1}h_2=1$:\\
$\Rightarrow\langle h_2,h_3\rangle\cong BS(1,2)$ is solvable, a contradiction.

\item[(86)]$ R_1:h_2h_3h_2^{-1}h_3^{-1}h_2h_3^{-1}h_2^{-1}h_3=1$, $R_2:h_2h_3^{-1}h_2^{-1}h_3^2h_2^{-1}h_3^{-1}h_2=1$, $R_3:h_2h_3^2h_2^{-1}h_3^2=1$:\\
$\Rightarrow\langle h_2,h_3\rangle\cong BS(1,-1)$ is solvable, a contradiction.

\item[(87)]$ R_1:h_2h_3h_2^{-1}h_3^{-1}h_2h_3^{-1}h_2^{-1}h_3=1$, $R_2:h_2h_3^{-1}h_2^{-1}h_3^2h_2^{-1}h_3^{-1}h_2=1$, $R_3:h_2(h_3h_2^{-1})^2h_3^2=1$:\\
$\Rightarrow\langle h_2,h_3\rangle\cong BS(1,1)$ is solvable, a contradiction.

\item[(88)]$ R_1:h_2h_3h_2^{-1}h_3^{-1}h_2h_3^{-1}h_2^{-1}h_3=1$, $R_2:h_2h_3^{-1}h_2^{-1}h_3^2h_2^{-1}h_3^{-1}h_2=1$, $R_3:h_2h_3^{-1}(h_2^{-1}h_3)^2h_3=1$:\\
$\Rightarrow \langle h_2,h_3\rangle=\langle h_3\rangle$ is abelian, a contradiction.

\item[(89)]$ R_1:h_2h_3h_2^{-1}h_3^{-1}h_2h_3^{-1}h_2^{-1}h_3=1$, $R_2:h_2h_3^{-1}h_2^{-1}h_3^2h_2^{-1}h_3^{-1}h_2=1$, $R_3:h_2h_3^{-2}(h_3^{-1}h_2)^2=1$:\\
$\Rightarrow\langle h_2,h_3\rangle\cong BS(1,1)$ is solvable, a contradiction.

\item[(90)]$ R_1:h_2h_3h_2^{-1}h_3^{-1}h_2h_3^{-1}h_2^{-1}h_3=1$, $R_2:h_2h_3^{-1}h_2^{-1}h_3^2h_2^{-1}h_3^{-1}h_2=1$, $R_3:h_2h_3^{-1}h_2h_3h_2^{-1}h_3^2=1$:\\
$\Rightarrow\langle h_2,h_3\rangle\cong BS(1,1)$ is solvable, a contradiction.

\item[(91)]$ R_1:h_2h_3h_2^{-1}h_3^{-1}h_2h_3^{-1}h_2^{-1}h_3=1$, $R_2:h_2h_3^{-1}h_2^{-1}h_3^2h_2^{-1}h_3^{-1}h_2=1$, $R_3:h_3^2(h_3h_2^{-1})^2h_3=1$:\\
$\Rightarrow\langle h_2,h_3\rangle\cong BS(1,1)$ is solvable, a contradiction.

\item[(92)]$ R_1:h_2h_3h_2^{-1}h_3^{-1}h_2h_3^{-1}h_2^{-1}h_3=1$, $R_2:h_2h_3^{-1}h_2^{-1}h_3^2h_2^{-1}h_3^{-1}h_2=1$, $R_3:h_3(h_3h_2^{-1})^3h_3=1$:\\
$\Rightarrow\langle h_2,h_3\rangle\cong BS(1,1)$ is solvable, a contradiction.

\item[(93)]$ R_1:h_2h_3h_2^{-1}(h_3^{-1}h_2)^2h_3=1$, $R_2:h_2h_3^{-2}h_2^{-1}h_3h_2^{-2}h_3=1$, $R_3:h_2^2h_3^{-1}h_2^2h_3=1$:\\
By interchanging $h_2$ and $h_3$ in (31) and with the same discussion, there is a contradiction.

\item[(94)]$ R_1:h_2h_3h_2^{-1}h_3h_2h_3^{-1}h_2h_3=1$, $R_2:h_2(h_3h_2^{-1})^2h_3^{-2}h_2=1$, $R_3:h_2h_3^{-2}h_2^{-1}h_3^2=1$:\\
By interchanging $h_2$ and $h_3$ in (51) and with the same discussion, there is a contradiction.

\item[(95)]$ R_1:h_2h_3h_2^{-1}h_3h_2h_3^{-1}h_2h_3=1$, $R_2:h_2(h_3h_2^{-1})^2h_3^{-2}h_2=1$, $R_3:h_2h_3^{-1}h_2h_3^3=1$:\\
By interchanging $h_2$ and $h_3$ in (50) and with the same discussion, there is a contradiction.

\item[(96)]$ R_1:h_2h_3h_2^{-1}h_3^2h_2h_3^{-1}h_2=1$, $R_2:h_2h_3^{-3}h_2^2h_3^{-1}h_2=1$, $R_3:h_2^2h_3^{-1}h_2^2h_3=1$:\\
By interchanging $h_2$ and $h_3$ in (34) and with the same discussion, there is a contradiction.

\item[(97)]$ R_1:h_2h_3h_2^{-1}h_3^2h_2^{-1}h_3^{-1}h_2=1$, $R_2:h_2(h_3h_2^{-1})^2h_3^{-1}h_2h_3=1$, $R_3:h_2h_3^{-2}h_2^{-1}h_3^2=1$:\\
By interchanging $h_2$ and $h_3$ in (52) and with the same discussion, there is a contradiction.

\item[(98)]$ R_1:h_2(h_3h_2^{-1}h_3)^2h_3=1$, $R_2:h_2(h_3h_2^{-1}h_3)^2h_3=1$, $R_3:h_2h_3^{-2}h_2^{-1}h_3^2=1$:\\
By interchanging $h_2$ and $h_3$ in (47) and with the same discussion, there is a contradiction.

\item[(99)]$ R_1:h_2(h_3h_2^{-1})^2h_3^{-2}h_2=1$, $R_2:(h_2h_3^{-1}h_2)^2h_3^{-2}h_2=1$, $R_3:h_2^2h_3^{-1}h_2^{-2}h_3=1$:\\
$\Rightarrow\langle h_2,h_3\rangle\cong BS(1,1)$ is solvable, a contradiction.

\item[(100)]$ R_1:h_2(h_3h_2^{-1})^2h_3^{-2}h_2=1$, $R_2:(h_2h_3^{-1})^2h_3^{-1}(h_3^{-1}h_2)^2=1$, $R_3:h_2^2h_3^{-1}h_2^{-2}h_3=1$:\\
By interchanging $h_2$ and $h_3$ in (40) and with the same discussion, there is a contradiction.

\item[(101)]$ R_1:h_2(h_3h_2^{-1})^2h_3^{-2}h_2=1$, $R_2:h_2(h_3^{-1}h_2h_3^{-1})^2h_3^{-1}h_2=1$, $R_3:h_2h_3^{-2}h_2^{-1}h_3^2=1$:\\
By interchanging $h_2$ and $h_3$ in (99) and with the same discussion, there is a contradiction.

\item[(102)]$ R_1:h_2h_3^{-1}h_2^{-1}h_3^2h_2h_3^{-1}h_2=1$, $R_2:h_2h_3^{-1}h_2^{-1}h_3h_2^{-1}h_3^{-1}h_2h_3=1$, $R_3:h_2^2h_3^{-1}h_2^{-2}h_3=1$:\\
By interchanging $h_2$ and $h_3$ in (76) and with the same discussion, there is a contradiction.

\item[(103)]$ R_1:h_2h_3^{-1}h_2^{-1}h_3^2h_2^{-1}h_3^{-1}h_2=1$, $R_2:h_2h_3^{-1}h_2^{-1}h_3h_2^{-1}h_3^{-1}h_2h_3=1$, $R_3:h_2^3h_3^2=1$:\\
By interchanging $h_2$ and $h_3$ in (78) and with the same discussion, there is a contradiction.

\item[(104)]$ R_1:h_2h_3^{-1}h_2^{-1}h_3^2h_2^{-1}h_3^{-1}h_2=1$, $R_2:h_2h_3^{-1}h_2^{-1}h_3h_2^{-1}h_3^{-1}h_2h_3=1$, $R_3:h_2^2(h_2h_3^{-1})^2h_2=1$:\\
By interchanging $h_2$ and $h_3$ in (91) and with the same discussion, there is a contradiction.

\item[(105)]$ R_1:h_2h_3^{-1}h_2^{-1}h_3^2h_2^{-1}h_3^{-1}h_2=1$, $R_2:h_2h_3^{-1}h_2^{-1}h_3h_2^{-1}h_3^{-1}h_2h_3=1$, $R_3:h_2^2h_3^3=1$:\\
By interchanging $h_2$ and $h_3$ in (77) and with the same discussion, there is a contradiction.

\item[(106)]$ R_1:h_2h_3^{-1}h_2^{-1}h_3^2h_2^{-1}h_3^{-1}h_2=1$, $R_2:h_2h_3^{-1}h_2^{-1}h_3h_2^{-1}h_3^{-1}h_2h_3=1$, $R_3:h_2^2h_3h_2^{-1}h_3^2=1$:\\
By interchanging $h_2$ and $h_3$ in (84) and with the same discussion, there is a contradiction.

\item[(107)]$ R_1:h_2h_3^{-1}h_2^{-1}h_3^2h_2^{-1}h_3^{-1}h_2=1$, $R_2:h_2h_3^{-1}h_2^{-1}h_3h_2^{-1}h_3^{-1}h_2h_3=1$, $R_3:h_2^2h_3^{-1}h_2^{-1}h_3^2=1$:\\
By interchanging $h_2$ and $h_3$ in (85) and with the same discussion, there is a contradiction.

\item[(108)]$ R_1:h_2h_3^{-1}h_2^{-1}h_3^2h_2^{-1}h_3^{-1}h_2=1$, $R_2:h_2h_3^{-1}h_2^{-1}h_3h_2^{-1}h_3^{-1}h_2h_3=1$, $R_3:h_2^2h_3^{-1}(h_3^{-1}h_2)^2=1$:\\
By interchanging $h_2$ and $h_3$ in (89) and with the same discussion, there is a contradiction.

\item[(109)]$ R_1:h_2h_3^{-1}h_2^{-1}h_3^2h_2^{-1}h_3^{-1}h_2=1$, $R_2:h_2h_3^{-1}h_2^{-1}h_3h_2^{-1}h_3^{-1}h_2h_3=1$, $R_3:h_2^2h_3^{-1}h_2^2h_3=1$:\\
By interchanging $h_2$ and $h_3$ in (86) and with the same discussion, there is a contradiction.

\item[(110)]$ R_1:h_2h_3^{-1}h_2^{-1}h_3^2h_2^{-1}h_3^{-1}h_2=1$, $R_2:h_2h_3^{-1}h_2^{-1}h_3h_2^{-1}h_3^{-1}h_2h_3=1$, $R_3:h_2^2h_3^{-1}h_2h_3^2=1$:\\
By interchanging $h_2$ and $h_3$ in (79) and with the same discussion, there is a contradiction.

\item[(111)]$ R_1:h_2h_3^{-1}h_2^{-1}h_3^2h_2^{-1}h_3^{-1}h_2=1$, $R_2:h_2h_3^{-1}h_2^{-1}h_3h_2^{-1}h_3^{-1}h_2h_3=1$, $R_3:h_2(h_2h_3^{-1})^3h_2=1$:\\
By interchanging $h_2$ and $h_3$ in (92) and with the same discussion, there is a contradiction.

\item[(112)]$ R_1:h_2h_3^{-1}h_2^{-1}h_3^2h_2^{-1}h_3^{-1}h_2=1$, $R_2:h_2h_3^{-1}h_2^{-1}h_3h_2^{-1}h_3^{-1}h_2h_3=1$, $R_3:h_2h_3(h_2h_3^{-1})^2h_2=1$:\\
By interchanging $h_2$ and $h_3$ in (87) and with the same discussion, there is a contradiction.

\item[(113)]$ R_1:h_2h_3^{-1}h_2^{-1}h_3^2h_2^{-1}h_3^{-1}h_2=1$, $R_2:h_2h_3^{-1}h_2^{-1}h_3h_2^{-1}h_3^{-1}h_2h_3=1$, $R_3:h_2h_3^2h_2h_3^{-1}h_2=1$:\\
By interchanging $h_2$ and $h_3$ in (80) and with the same discussion, there is a contradiction.

\item[(114)]$ R_1:h_2h_3^{-1}h_2^{-1}h_3^2h_2^{-1}h_3^{-1}h_2=1$, $R_2:h_2h_3^{-1}h_2^{-1}h_3h_2^{-1}h_3^{-1}h_2h_3=1$, $R_3:h_2h_3^2h_2^{-1}h_3^2=1$:\\
By interchanging $h_2$ and $h_3$ in (81) and with the same discussion, there is a contradiction.

\item[(115)]$ R_1:h_2h_3^{-1}h_2^{-1}h_3^2h_2^{-1}h_3^{-1}h_2=1$, $R_2:h_2h_3^{-1}h_2^{-1}h_3h_2^{-1}h_3^{-1}h_2h_3=1$, $R_3:h_2h_3h_2^{-2}h_3^2=1$:\\
By interchanging $h_2$ and $h_3$ in (82) and with the same discussion, there is a contradiction.

\item[(116)]$ R_1:h_2h_3^{-1}h_2^{-1}h_3^2h_2^{-1}h_3^{-1}h_2=1$, $R_2:h_2h_3^{-1}h_2^{-1}h_3h_2^{-1}h_3^{-1}h_2h_3=1$, $R_3:h_2h_3h_2^{-1}(h_3^{-1}h_2)^2=1$:\\
By interchanging $h_2$ and $h_3$ in (88) and with the same discussion, there is a contradiction.

\item[(117)]$ R_1:h_2h_3^{-1}h_2^{-1}h_3^2h_2^{-1}h_3^{-1}h_2=1$, $R_2:h_2h_3^{-1}h_2^{-1}h_3h_2^{-1}h_3^{-1}h_2h_3=1$, $R_3:h_2h_3h_2^{-1}h_3h_2h_3^{-1}h_2=1$:\\
By interchanging $h_2$ and $h_3$ in (90) and with the same discussion, there is a contradiction.

\item[(118)]$ R_1:h_2h_3^{-1}h_2^{-1}h_3^2h_2^{-1}h_3^{-1}h_2=1$, $R_2:h_2h_3^{-1}h_2^{-1}h_3h_2^{-1}h_3^{-1}h_2h_3=1$, $R_3:h_2(h_3h_2^{-1})^2h_3^2=1$:\\
By interchanging $h_2$ and $h_3$ in (83) and with the same discussion, there is a contradiction.

\item[(119)]$ R_1:h_2h_3^{-1}h_2^{-1}h_3^2h_2^{-1}h_3^{-1}h_2=1$, $R_2:h_2h_3^{-3}h_2^2h_3^{-1}h_2=1$, $R_3:h_2h_3^{-2}h_2h_3^2=1$:\\
By interchanging $h_2$ and $h_3$ in (35) and with the same discussion, there is a contradiction.

\item[(120)]$ R_1:h_2h_3^{-1}h_2^{-1}h_3h_2^{-1}h_3^{-1}h_2h_3=1$, $R_2:h_2h_3^{-2}h_2^{-1}h_3h_2^{-2}h_3=1$, $R_3:h_2^2h_3^{-1}h_2^2h_3=1$:\\
By interchanging $h_2$ and $h_3$ in (32) and with the same discussion, there is a contradiction.

\item[(121)]$ R_1:h_2h_3^{-1}h_2^{-1}h_3h_2^{-1}h_3^{-1}h_2h_3=1$, $R_2:h_2h_3^{-1}(h_3^{-1}h_2)^2h_3^2=1$, $R_3:(h_2h_3^{-1})^2(h_2^{-1}h_3)^2=1$:\\
By interchanging $h_2$ and $h_3$ in (10) and with the same discussion, there is a contradiction.

\item[(122)]$ R_1:h_2h_3^{-1}h_2^{-1}h_3h_2^{-1}h_3^{-1}h_2h_3=1$, $R_2:h_2h_3^{-1}(h_3^{-1}h_2)^2h_3^2=1$, $R_3:h_2^2h_3^{-2}h_2h_3=1$:\\
By interchanging $h_2$ and $h_3$ in (12) and with the same discussion, there is a contradiction.

\item[(123)]$ R_1:h_2h_3^{-1}h_2^{-1}h_3h_2^{-1}h_3^{-1}h_2h_3=1$, $R_2:h_2h_3^{-1}(h_3^{-1}h_2)^2h_3^2=1$, $R_3:h_2h_3^{-2}h_2^{-1}h_3^2=1$:\\
By interchanging $h_2$ and $h_3$ in (11) and with the same discussion, there is a contradiction.

\item[(124)]$ R_1:h_2h_3^{-1}h_2^{-1}h_3h_2^{-1}h_3^{-1}h_2h_3=1$, $R_2:(h_2h_3^{-1}h_2)^2h_3^2=1$, $R_3:h_2^2h_3^{-1}h_2^{-2}h_3=1$:\\
By interchanging $h_2$ and $h_3$ in (23) and with the same discussion, there is a contradiction.

\end{enumerate}

\subsection{$\mathbf{C_4-C_6(--C_7--)(---C_6)}$}
\begin{figure}[ht]
\psscalebox{0.9 0.9} 
{
\begin{pspicture}(0,-2.0985577)(5.59,2.0985577)
\psdots[linecolor=black, dotsize=0.4](4.0,1.9014423)
\psdots[linecolor=black, dotsize=0.4](3.2,1.5014423)
\psdots[linecolor=black, dotsize=0.4](3.2,0.70144236)
\psdots[linecolor=black, dotsize=0.4](4.0,0.30144233)
\psdots[linecolor=black, dotsize=0.4](4.8,0.70144236)
\psdots[linecolor=black, dotsize=0.4](4.8,1.5014423)
\psdots[linecolor=black, dotsize=0.4](2.4,1.5014423)
\psdots[linecolor=black, dotsize=0.4](2.4,0.70144236)
\psdots[linecolor=black, dotsize=0.4](1.6,-0.49855766)
\psline[linecolor=black, linewidth=0.04](2.4,1.5014423)(3.2,1.5014423)(4.0,1.9014423)(4.8,1.5014423)(4.8,0.70144236)(4.0,0.30144233)(3.2,0.70144236)(2.4,0.70144236)(2.4,1.5014423)(1.6,-0.49855766)
\psline[linecolor=black, linewidth=0.04](3.2,1.5014423)(3.2,0.70144236)
\psline[linecolor=black, linewidth=0.04](1.6,-0.49855766)(4.4,-0.49855766)
\psline[linecolor=black, linewidth=0.04](4.4,-0.49855766)(4.8,0.70144236)
\psdots[linecolor=black, dotsize=0.4](4.4,-0.49855766)
\rput[bl](-0.4,-2.0985577){$\mathbf{27) \ C_4-C_6(--C_7--)(---C_6)}$}
\psdots[linecolor=black, dotsize=0.4](1.2,0.30144233)
\psdots[linecolor=black, dotsize=0.4](1.2,1.9014423)
\psline[linecolor=black, linewidth=0.04](4.0,1.9014423)(1.2,1.9014423)(1.2,0.30144233)(1.6,-0.49855766)
\end{pspicture}
}
\end{figure}
By considering the $16$ cases related to the existence of $C_4-C_6(--C_7--)$ in the graph $K(\alpha,\beta)$ from Table \ref{tab-C4-C6(--C7--)} and the relations from Table \ref{tab-C6} which are not disproved, it can be seen that there are $124$ cases for the relations of a cycle $C_4$, two cycles $C_6$ and a cycle $C_7$ in the graph $C_4-C_6(--C_7--)(---C_6)$. Using GAP \cite{gap}, we see that all groups with two generators $h_2$ and $h_3$ and four relations which are between $112$ cases of these $124$ cases are finite and solvable, that is a contradiction with the assumptions. So, there are just $12$ cases for the relations of these cycles which may lead to the existence of a subgraph isomorphic to the graph $C_4-C_6(--C_7--)(---C_6)$ in $K(\alpha,\beta)$. In the following, we show that these $12$ cases lead to contradictions and so, the graph $K(\alpha,\beta)$ contains no subgraph isomorphic to the graph $C_4-C_6(--C_7--)(---C_6)$.
\begin{enumerate}
\item[(1)]$ R_1:h_2h_3h_2^{-2}h_3=1$, $R_2:h_2h_3^2h_2^{-1}(h_2^{-1}h_3)^2=1$, $R_3:h_2h_3^{-5}h_2h_3=1$,\\ $R_4:h_3^2h_2^{-1}(h_3h_2^{-1}h_3)^2=1$:\\
$\Rightarrow\langle h_2,h_3\rangle\cong BS(-2,1)$ is solvable, a contradiction.

\item[(2)]$ R_1:h_2h_3h_2^{-2}h_3=1$, $R_2:h_2h_3^2h_2^{-1}(h_2^{-1}h_3)^2=1$, $R_3:h_2h_3^{-4}h_2h_3h_2^{-1}h_3=1$,\\ $R_4:h_2^2h_3^{-1}(h_2^{-1}h_3)^3=1$:\\
$\Rightarrow\langle h_2,h_3\rangle\cong BS(1,-1)$ is solvable, a contradiction.

\item[(3)]$ R_1:h_2h_3h_2^{-2}h_3=1$, $R_2:h_2h_3^2h_2^{-1}(h_2^{-1}h_3)^2=1$, $R_3:h_2h_3^{-2}h_2(h_3h_2^{-1})^3h_3=1$,\\ $R_4:h_2^2h_3^{-3}h_2^{-1}h_3=1$:\\
$\Rightarrow\langle h_2,h_3\rangle\cong BS(1,-1)$ is solvable, a contradiction.

\item[(4)]$ R_1:h_2h_3h_2^{-2}h_3=1$, $R_2:h_2h_3^2h_2^{-1}(h_2^{-1}h_3)^2=1$, $R_3:h_2h_3^{-1}h_2(h_3h_2^{-1})^4h_3=1$,\\ $R_4:h_3^2h_2^{-1}(h_3h_2^{-1}h_3)^2=1$:\\
$\Rightarrow\langle h_2,h_3\rangle\cong BS(-2,1)$ is solvable, a contradiction.

\item[(5)]$ R_1:h_2h_3h_2^{-2}h_3=1$, $R_2:h_2(h_3h_2^{-1})^2h_2^{-1}h_3^2=1$, $R_3:(h_2h_3^{-2})^2h_2^{-1}h_3^2=1$,\\ $R_4:h_2(h_3h_2^{-1})^3h_3^{-1}h_2=1$:\\
$\Rightarrow\langle h_2,h_3\rangle\cong BS(1,-1)$ is solvable, a contradiction.

\item[(6)]$ R_1:h_2h_3h_2^{-2}h_3=1$, $R_2:h_2(h_3h_2^{-1})^2h_2^{-1}h_3^2=1$, $R_3:h_2h_3^{-1}h_2h_3^{-1}h_2^{-1}(h_3h_2^{-1}h_3)^2=1$,\\ $R_4:h_2h_3h_2^{-1}h_3^{-3}h_2=1$:\\
$\Rightarrow\langle h_2,h_3\rangle\cong BS(1,-1)$ is solvable, a contradiction.

\item[(7)]$ R_1:h_2h_3^{-2}h_2h_3=1$, $R_2:h_2h_3(h_2h_3^{-1})^2h_3^{-1}h_2=1$, $R_3:(h_2^2h_3^{-1})^2h_2^{-2}h_3=1$,\\ $R_4:(h_2h_3^{-1})^3h_2^{-1}h_3^2=1$:\\
By interchanging $h_2$ and $h_3$ in (5) and with the same discussion, there is a contradiction.

\item[(8)]$ R_1:h_2h_3^{-2}h_2h_3=1$, $R_2:h_2h_3(h_2h_3^{-1})^2h_3^{-1}h_2=1$, $R_3:h_2h_3^{-1}h_2(h_3h_2^{-1})^2(h_3^{-1}h_2)^2=1$,\\ $R_4:h_2^2h_3h_2^{-1}h_3^{-2}h_2=1$:\\
By interchanging $h_2$ and $h_3$ in (6) and with the same discussion, there is a contradiction.

\item[(9)]$ R_1:h_2h_3^{-2}h_2h_3=1$, $R_2:h_2^2h_3^{-1}(h_3^{-1}h_2)^2h_3=1$, $R_3:h_2^4h_3^{-1}h_2^{-1}h_3^{-1}h_2=1$,\\ $R_4:h_2^2h_3^{-1}(h_2h_3^{-1}h_2)^2=1$:\\
By interchanging $h_2$ and $h_3$ in (1) and with the same discussion, there is a contradiction.

\item[(10)]$ R_1:h_2h_3^{-2}h_2h_3=1$, $R_2:h_2^2h_3^{-1}(h_3^{-1}h_2)^2h_3=1$, $R_3:h_2^3h_3^{-1}h_2^{-1}h_3h_2^{-1}h_3^{-1}h_2=1$,\\ $R_4:h_2h_3^{-2}(h_2^{-1}h_3)^3=1$:\\
By interchanging $h_2$ and $h_3$ in (2) and with the same discussion, there is a contradiction.

\item[(11)]$ R_1:h_2h_3^{-2}h_2h_3=1$, $R_2:h_2^2h_3^{-1}(h_3^{-1}h_2)^2h_3=1$, $R_3:h_2h_3^{-1}(h_2^{-1}h_3)^3h_2^{-1}h_3^{-1}h_2=1$,\\ $R_4:h_2^3h_3^{-2}h_2^{-1}h_3=1$:\\
By interchanging $h_2$ and $h_3$ in (3) and with the same discussion, there is a contradiction.

\item[(12)]$ R_1:h_2h_3^{-2}h_2h_3=1$, $R_2:h_2^2h_3^{-1}(h_3^{-1}h_2)^2h_3=1$, $R_3:(h_2h_3^{-1})^4h_2h_3h_2^{-1}h_3=1$,\\ $R_4:h_2^2h_3^{-1}(h_2h_3^{-1}h_2)^2=1$:\\
By interchanging $h_2$ and $h_3$ in (4) and with the same discussion, there is a contradiction.

\end{enumerate}

\subsection{$\mathbf{C_4-C_6(--C_7--)(C_4)(C_4)}$}
\begin{figure}[ht]
\psscalebox{0.9 0.9} 
{
\begin{pspicture}(0,-2.0985577)(6.71,2.0985577)
\psdots[linecolor=black, dotsize=0.4](4.0,1.9014423)
\psdots[linecolor=black, dotsize=0.4](3.2,1.5014423)
\psdots[linecolor=black, dotsize=0.4](3.2,0.70144236)
\psdots[linecolor=black, dotsize=0.4](4.0,0.30144233)
\psdots[linecolor=black, dotsize=0.4](4.8,0.70144236)
\psdots[linecolor=black, dotsize=0.4](4.8,1.5014423)
\psdots[linecolor=black, dotsize=0.4](2.4,1.5014423)
\psdots[linecolor=black, dotsize=0.4](2.4,0.70144236)
\psdots[linecolor=black, dotsize=0.4](1.6,-0.49855766)
\psline[linecolor=black, linewidth=0.04](2.4,1.5014423)(3.2,1.5014423)(4.0,1.9014423)(4.8,1.5014423)(4.8,0.70144236)(4.0,0.30144233)(3.2,0.70144236)(2.4,0.70144236)(2.4,1.5014423)(1.6,-0.49855766)
\psline[linecolor=black, linewidth=0.04](3.2,1.5014423)(3.2,0.70144236)
\psline[linecolor=black, linewidth=0.04](1.6,-0.49855766)(4.4,-0.49855766)
\psline[linecolor=black, linewidth=0.04](4.4,-0.49855766)(4.8,0.70144236)
\psdots[linecolor=black, dotsize=0.4](4.4,-0.49855766)
\psdots[linecolor=black, dotsize=0.4](2.4,-0.098557666)
\psline[linecolor=black, linewidth=0.04](2.4,0.70144236)(2.4,-0.098557666)(1.6,-0.49855766)(1.6,-0.49855766)
\psline[linecolor=black, linewidth=0.04](4.0,0.30144233)(3.6,-0.098557666)(4.4,-0.49855766)
\psdots[linecolor=black, dotsize=0.4](3.6,-0.098557666)
\rput[bl](0.4,-2.0985577){$\mathbf{28) \ C_4-C_6(--C_7--)(C_4)(C_4)}$}
\end{pspicture}
}
\end{figure}
By considering the $16$ cases related to the existence of $C_4-C_6(--C_7--)$ in the graph $K(\alpha,\beta)$ from Table \ref{tab-C4-C6(--C7--)} and the relations from Table \ref{tab-C4} which are not disproved, it can be seen that there are $8$ cases for the relations of three cycles $C_4$, a cycle $C_7$ and a cycle $C_6$ in this structure. By considering all groups with two generators $h_2$ and $h_3$ and five relations which are between these cases and by using GAP \cite{gap}, we see that all of these groups are finite and solvable. So, the graph $K(\alpha,\beta)$ contains no subgraph isomorphic to the graph $C_4-C_6(--C_7--)(C_4)(C_4)$.

\subsection{$\mathbf{C_6---C_6(-C_5--)}$}
\begin{figure}[ht]
\psscalebox{0.9 0.9} 
{
\begin{pspicture}(0,-1.8985577)(4.21,1.8985577)
\psdots[linecolor=black, dotsize=0.4](2.4,1.7014424)
\psdots[linecolor=black, dotsize=0.4](2.4,0.90144235)
\psdots[linecolor=black, dotsize=0.4](2.4,0.10144234)
\psdots[linecolor=black, dotsize=0.4](2.4,-0.6985577)
\psdots[linecolor=black, dotsize=0.4](3.6,0.10144234)
\psdots[linecolor=black, dotsize=0.4](3.6,0.90144235)
\psdots[linecolor=black, dotsize=0.4](1.2,0.90144235)
\psdots[linecolor=black, dotsize=0.4](1.2,0.10144234)
\psdots[linecolor=black, dotsize=0.4](0.4,-0.6985577)
\psline[linecolor=black, linewidth=0.04](1.2,0.90144235)(2.4,1.7014424)(2.4,0.90144235)(2.4,0.10144234)(2.4,-0.6985577)(1.2,0.10144234)(1.2,0.90144235)(0.4,-0.6985577)(2.4,0.10144234)
\psline[linecolor=black, linewidth=0.04](2.4,-0.6985577)(3.6,0.10144234)(3.6,0.90144235)(2.4,1.7014424)
\rput[bl](-0.2,-1.8985577){$\mathbf{29) \ C_6---C_6(-C_5--)}$}
\end{pspicture}
}
\end{figure}
By considering the $99$ cases related to the existence of $C_6---C_6$ in the graph $K(\alpha,\beta)$ and the relations from Table \ref{tab-C5} which are not disproved, it can be seen that there are $418$ cases for the relations of two cycles $C_6$ and a cycle $C_5$ in the graph $C_6---C_6(-C_5--)$. Using GAP \cite{gap}, we see that all groups with two generators $h_2$ and $h_3$ and three relations which are between $358$ cases of these $418$ cases are finite and solvable, that is a contradiction with the assumptions. So, there are just $60$ cases for the relations of these cycles which may lead to the existence of a subgraph isomorphic to the graph $C_6---C_6(-C_5--)$ in $K(\alpha,\beta)$. In the following, we show that these $60$ cases lead to contradictions and so, the graph $K(\alpha,\beta)$ contains no subgraph isomorphic to the graph $C_6---C_6(-C_5--)$.
\begin{enumerate}
\item[(1)]$ R_1:h_2^3h_3^3=1$, $R_2:h_2^2(h_2h_3^{-1})^3h_2=1$, $R_3:h_2^2h_3^{-3}h_2=1$:\\
$\Rightarrow$  $G$ has a torsion element, a contradiction.

\item[(2)]$ R_1:h_2^3h_3^3=1$, $R_2:h_3^2(h_3h_2^{-1})^3h_3=1$, $R_3:h_2^2h_3^{-3}h_2=1$:\\
By interchanging $h_2$ and $h_3$ in (1) and with the same discussion, there is a contradiction.

\item[(3)]$ R_1:h_2^2h_3h_2^{-1}(h_2^{-1}h_3)^2=1$, $R_2:h_2h_3h_2^{-1}h_3^{-1}h_2h_3^{-1}h_2^{-1}h_3=1$, $R_3:h_2h_3^2h_2^{-2}h_3=1$:\\
$\Rightarrow \langle h_2,h_3\rangle=\langle h_3\rangle$ is abelian, a contradiction.

\item[(4)]$ R_1:h_2^2h_3h_2^{-1}h_3^{-2}h_2=1$, $R_2:h_2^2h_3^{-1}h_2^{-1}h_3^{-1}h_2h_3=1$, $R_3:h_2^3h_3^{-2}h_2=1$:\\
$\Rightarrow \langle h_2,h_3\rangle=\langle h_2\rangle$ is abelian, a contradiction.

\item[(5)]$ R_1:h_2^2h_3h_2^{-1}h_3^{-2}h_2=1$, $R_2:h_2h_3h_2h_3^{-1}(h_3^{-1}h_2)^2=1$, $R_3:h_2^3h_3^{-2}h_2=1$:\\
$\Rightarrow \langle h_2,h_3\rangle=\langle h_2\rangle$ is abelian, a contradiction.

\item[(6)]$ R_1:h_2^2(h_3h_2^{-1}h_3)^2=1$, $R_2:(h_2h_3)^2(h_2^{-1}h_3)^2=1$, $R_3:h_2h_3^2h_2^{-1}h_3^2=1$:\\
$\Rightarrow$ $G$ is solvable, a contradiction.

\item[(7)]$ R_1:h_2^2(h_3h_2^{-1}h_3)^2=1$, $R_2:h_2h_3h_2^{-1}h_3^{-1}h_2h_3^{-1}h_2^{-1}h_3=1$, $R_3:h_2h_3^{-1}h_2(h_3h_2^{-1})^2h_3=1$:\\
$\Rightarrow$  $G$ has a torsion element, a contradiction.

\item[(8)]$ R_1:h_2^2(h_3h_2^{-1}h_3)^2=1$, $R_2:h_2h_3h_2^{-1}h_3^{-1}h_2h_3^{-1}h_2^{-1}h_3=1$, $R_3:(h_2h_3^{-1})^2(h_2^{-1}h_3)^2=1$:\\
$\Rightarrow$ $G$ is solvable, a contradiction.

\item[(9)]$ R_1:h_2^2(h_3h_2^{-1}h_3)^2=1$, $R_2:h_2(h_3h_2^{-1})^2h_3^{-1}h_2h_3=1$, $R_3:h_2h_3^2h_2^{-1}h_3^2=1$:\\
$\Rightarrow$ $G$ is solvable, a contradiction.

\item[(10)]$ R_1:h_2^2h_3^{-1}h_2^{-1}h_3^2h_2^{-1}h_3=1$, $R_2:h_2h_3h_2h_3^{-1}(h_2^{-1}h_3)^2=1$, $R_3:h_2h_3^{-2}h_2^{-1}h_3^2=1$:\\
$\Rightarrow$  $G$ has a torsion element, a contradiction.

\item[(11)]$ R_1:h_2^2h_3^{-1}h_2^{-1}h_3^2h_2^{-1}h_3=1$, $R_2:h_2h_3h_2^{-1}h_3^{-1}h_2h_3^{-1}h_2^{-1}h_3=1$, $R_3:h_2h_3^{-1}h_2(h_3h_2^{-1})^2h_3=1$:\\
$\Rightarrow$  $G$ has a torsion element, a contradiction.

\item[(12)]$ R_1:h_2^2h_3^{-1}h_2^{-1}h_3^2h_2^{-1}h_3=1$, $R_2:h_2h_3h_2^{-1}h_3^{-1}h_2h_3^{-1}h_2^{-1}h_3=1$, $R_3:(h_2h_3^{-1})^2(h_2^{-1}h_3)^2=1$:\\
$\Rightarrow$  $G$ has a torsion element, a contradiction.

\item[(13)]$ R_1:h_2^2h_3^{-2}h_2h_3^{-1}h_2^{-1}h_3=1$, $R_2:h_2h_3^{-1}h_2^{-1}h_3h_2h_3^{-2}h_2=1$, $R_3:h_2^2h_3^{-1}h_2^{-1}h_3^{-1}h_2=1$:\\
$\Rightarrow$ $G$ is solvable, a contradiction.

\item[(14)]$ R_1:h_2^2h_3^{-2}h_2h_3^{-1}h_2^{-1}h_3=1$, $R_2:h_2h_3^{-1}h_2^{-1}h_3^2h_2^{-1}h_3^{-1}h_2=1$, $R_3:h_2^2h_3^{-1}h_2^{-1}h_3^{-1}h_2=1$:\\
$\Rightarrow\langle h_2,h_3\rangle\cong BS(1,1)$ is solvable, a contradiction.

\item[(15)]$ R_1:h_2(h_2h_3^{-2})^2h_2=1$, $R_2:h_2h_3h_2^{-1}h_3^2h_2h_3^{-1}h_2=1$, $R_3:h_2(h_3^{-1}h_2h_3^{-1})^2h_2=1$:\\
$\Rightarrow\langle h_2,h_3\rangle\cong BS(1,-1)$ is solvable, a contradiction.

\item[(16)]$ R_1:h_2(h_2h_3^{-2})^2h_2=1$, $R_2:h_2h_3^{-1}h_2^{-1}h_3h_2h_3^{-2}h_2=1$, $R_3:h_2^2h_3^{-1}h_2^{-2}h_3=1$:\\
$\Rightarrow$ $G$ is solvable, a contradiction.

\item[(17)]$ R_1:h_2(h_2h_3^{-2})^2h_2=1$, $R_2:h_2h_3^{-1}h_2^{-1}h_3^2h_2^{-1}h_3^{-1}h_2=1$, $R_3:h_2^2h_3^{-1}h_2^{-2}h_3=1$:\\
$\Rightarrow$ $G$ is solvable, a contradiction.

\item[(18)]$ R_1:(h_2h_3)^2h_2h_3^{-1}h_2=1$, $R_2:(h_2h_3)^2h_2h_3^{-1}h_2=1$, $R_3:h_2^2h_3^{-1}h_2^{-2}h_3=1$:\\
$\Rightarrow \langle h_2,h_3\rangle=\langle h_2\rangle$ is abelian, a contradiction.

\item[(19)]$ R_1:(h_2h_3)^2h_2^{-1}h_3^2=1$, $R_2:(h_2h_3)^2h_2^{-1}h_3^2=1$, $R_3:h_2h_3^{-2}h_2^{-1}h_3^2=1$:\\
By interchanging $h_2$ and $h_3$ in (18) and with the same discussion, there is a contradiction.

\item[(20)]$ R_1:h_2h_3(h_2h_3^{-1})^2h_2^{-1}h_3=1$, $R_2:h_2h_3^{-1}h_2^{-1}h_3^2h_2h_3^{-1}h_2=1$, $R_3:h_2^2h_3^{-1}h_2^2h_3=1$:\\
$\Rightarrow$  $G$ has a torsion element, a contradiction.

\item[(21)]$ R_1:h_2h_3(h_2h_3^{-1})^2h_2^{-1}h_3=1$, $R_2:(h_2h_3^{-1}h_2)^2h_3^2=1$, $R_3:h_2^2h_3^{-1}h_2^2h_3=1$:\\
By interchanging $h_2$ and $h_3$ in (9) and with the same discussion, there is a contradiction.

\item[(22)]$ R_1:h_2h_3(h_2h_3^{-1})^2h_2h_3=1$, $R_2:(h_2h_3^{-1}h_2)^2h_3^2=1$, $R_3:h_2^2h_3^{-1}h_2^2h_3=1$:\\
By interchanging $h_2$ and $h_3$ in (6) and with the same discussion, there is a contradiction.

\item[(23)]$ R_1:h_2h_3^2h_2^{-1}h_3^{-1}h_2^{-1}h_3=1$, $R_2:h_2h_3h_2^{-1}h_3^{-3}h_2=1$, $R_3:h_2h_3^{-4}h_2=1$:\\
By interchanging $h_2$ and $h_3$ in (4) and with the same discussion, there is a contradiction.

\item[(24)]$ R_1:h_2(h_3h_2^{-2})^2h_3=1$, $R_2:h_2h_3^{-1}h_2^{-1}h_3^2h_2^{-2}h_3=1$, $R_3:h_2h_3h_2^{-1}h_3^{-2}h_2=1$:\\
$\Rightarrow$  $G$ has a torsion element, a contradiction.

\item[(25)]$ R_1:h_2h_3h_2^{-2}h_3h_2^{-1}h_3^{-1}h_2=1$, $R_2:h_2h_3^{-1}h_2^{-1}h_3^2h_2^{-2}h_3=1$, $R_3:h_2h_3h_2^{-1}h_3^{-1}h_2^{-1}h_3=1$:\\
$\Rightarrow$  $G$ has a torsion element, a contradiction.

\item[(26)]$ R_1:h_2h_3h_2^{-1}(h_2^{-1}h_3)^2h_3=1$, $R_2:h_2h_3h_2^{-1}h_3^{-3}h_2=1$, $R_3:h_2h_3^{-4}h_2=1$:\\
By interchanging $h_2$ and $h_3$ in (5) and with the same discussion, there is a contradiction.

\item[(27)]$ R_1:h_2h_3h_2^{-1}h_3^{-1}h_2h_3^{-1}h_2^{-1}h_3=1$, $R_2:h_2h_3h_2^{-1}h_3^2h_2^{-1}h_3^{-1}h_2=1$, $R_3:h_2^2h_3h_2^{-2}h_3=1$:\\
$\Rightarrow$  $G$ has a torsion element, a contradiction.

\item[(28)]$ R_1:h_2h_3h_2^{-1}h_3^{-1}h_2h_3^{-1}h_2^{-1}h_3=1$, $R_2:h_2h_3h_2^{-1}h_3^2h_2^{-1}h_3^{-1}h_2=1$, $R_3:h_2^2h_3^{-1}h_2^{-2}h_3=1$:\\
$\Rightarrow$  $G$ has a torsion element, a contradiction.

\item[(29)]$ R_1:h_2h_3h_2^{-1}h_3^{-1}h_2h_3^{-1}h_2^{-1}h_3=1$, $R_2:h_2h_3^{-1}h_2^{-1}h_3h_2h_3^{-2}h_2=1$, $R_3:h_2^2h_3^{-1}h_2^{-2}h_3=1$:\\
$\Rightarrow$ $G$ is solvable, a contradiction.

\item[(30)]$ R_1:h_2h_3h_2^{-1}h_3^{-1}h_2h_3^{-1}h_2^{-1}h_3=1$, $R_2:h_2h_3^{-1}h_2^{-1}h_3h_2h_3^{-2}h_2=1$, $R_3:h_2h_3^{-1}(h_3^{-1}h_2)^3=1$:\\
$\Rightarrow\langle h_2,h_3\rangle\cong BS(-3,1)$ is solvable, a contradiction.

\item[(31)]$ R_1:h_2h_3h_2^{-1}h_3^{-1}h_2h_3^{-1}h_2^{-1}h_3=1$, $R_2:h_2h_3^{-1}h_2^{-1}h_3^2h_2^{-1}h_3^{-1}h_2=1$, $R_3:h_2^2h_3h_2^{-2}h_3=1$:\\
$\Rightarrow$  $G$ has a torsion element, a contradiction.

\item[(32)]$ R_1:h_2h_3h_2^{-1}h_3^{-1}h_2h_3^{-1}h_2^{-1}h_3=1$, $R_2:h_2h_3^{-1}h_2^{-1}h_3^2h_2^{-1}h_3^{-1}h_2=1$, $R_3:h_2^2h_3^{-1}h_2^{-2}h_3=1$:\\
$\Rightarrow$ $G$ is solvable, a contradiction.

\item[(33)]$ R_1:h_2h_3h_2^{-1}h_3^{-1}h_2h_3^{-1}h_2^{-1}h_3=1$, $R_2:h_2h_3^{-1}h_2^{-1}h_3^2h_2^{-1}h_3^{-1}h_2=1$, $R_3:h_2^2h_3^{-1}(h_3^{-1}h_2)^2=1$:\\
$\Rightarrow\langle h_2,h_3\rangle\cong BS(1,1)$ is solvable, a contradiction.

\item[(34)]$ R_1:h_2h_3h_2^{-1}h_3^{-1}h_2h_3^{-1}h_2^{-1}h_3=1$, $R_2:h_2h_3^{-1}h_2^{-1}h_3^2h_2^{-1}h_3^{-1}h_2=1$, $R_3:h_2h_3^{-1}(h_3^{-1}h_2)^3=1$:\\
$\Rightarrow\langle h_2,h_3\rangle\cong BS(-3,1)$ is solvable, a contradiction.

\item[(35)]$ R_1:h_2h_3h_2^{-1}h_3^{-1}h_2h_3^{-2}h_2=1$, $R_2:h_2h_3^{-1}h_2^{-1}h_3^2h_2^{-1}h_3^{-1}h_2=1$, $R_3:h_2^2h_3^{-1}h_2^{-1}h_3^{-1}h_2=1$:\\
$\Rightarrow\langle h_2,h_3\rangle\cong BS(1,1)$ is solvable, a contradiction.

\item[(36)]$ R_1:h_2h_3h_2^{-1}h_3^{-1}h_2h_3^{-2}h_2=1$, $R_2:h_2h_3^{-1}h_2^{-1}h_3^2h_2^{-1}h_3^{-1}h_2=1$, $R_3:h_2h_3^{-1}(h_3^{-1}h_2)^3=1$:\\
$\Rightarrow\langle h_2,h_3\rangle\cong BS(-3,1)$ is solvable, a contradiction.

\item[(37)]$ R_1:h_2h_3h_2^{-1}(h_3^{-1}h_2)^2h_3=1$, $R_2:h_2h_3^{-2}h_2^{-1}h_3h_2^{-2}h_3=1$, $R_3:h_2^2h_3^{-1}h_2^{-2}h_3=1$:\\
By interchanging $h_2$ and $h_3$ in (10) and with the same discussion, there is a contradiction.

\item[(38)]$ R_1:h_2h_3h_2^{-1}h_3^2h_2h_3^{-1}h_2=1$, $R_2:h_2h_3^{-3}h_2^2h_3^{-1}h_2=1$, $R_3:h_2(h_3^{-1}h_2h_3^{-1})^2h_2=1$:\\
By interchanging $h_2$ and $h_3$ in (15) and with the same discussion, there is a contradiction.

\item[(39)]$ R_1:h_2h_3h_2^{-1}h_3^2h_2^{-1}h_3^{-1}h_2=1$, $R_2:h_2(h_3h_2^{-1})^2h_3^{-1}h_2h_3=1$, $R_3:h_2h_3^2h_2^{-1}h_3^2=1$:\\
By interchanging $h_2$ and $h_3$ in (20) and with the same discussion, there is a contradiction.

\item[(40)]$ R_1:h_2h_3^{-1}h_2^{-1}h_3h_2h_3^{-2}h_2=1$, $R_2:h_2h_3^{-1}h_2^{-1}h_3h_2^{-1}h_3^{-1}h_2h_3=1$, $R_3:h_2h_3^{-2}h_2^{-1}h_3^2=1$:\\
By interchanging $h_2$ and $h_3$ in (29) and with the same discussion, there is a contradiction.

\item[(41)]$ R_1:h_2h_3^{-1}h_2^{-1}h_3h_2h_3^{-2}h_2=1$, $R_2:h_2h_3^{-1}h_2^{-1}h_3h_2^{-1}h_3^{-1}h_2h_3=1$, $R_3:h_2h_3^{-1}(h_3^{-1}h_2)^3=1$:\\
By interchanging $h_2$ and $h_3$ in (30) and with the same discussion, there is a contradiction.

\item[(42)]$ R_1:h_2h_3^{-1}h_2^{-1}h_3h_2h_3^{-2}h_2=1$, $R_2:h_2h_3^{-2}h_2^{-1}h_3h_2h_3^{-1}h_2=1$, $R_3:h_2h_3^{-3}h_2h_3=1$:\\
By interchanging $h_2$ and $h_3$ in (13) and with the same discussion, there is a contradiction.

\item[(43)]$ R_1:h_2h_3^{-1}h_2^{-1}h_3h_2h_3^{-2}h_2=1$, $R_2:h_2h_3^{-3}h_2^2h_3^{-1}h_2=1$, $R_3:h_2h_3^{-2}h_2^{-1}h_3^2=1$:\\
By interchanging $h_2$ and $h_3$ in (16) and with the same discussion, there is a contradiction.

\item[(44)]$ R_1:h_2h_3^{-1}h_2^{-1}h_3^2h_2h_3^{-1}h_2=1$, $R_2:h_2h_3^{-1}h_2^{-1}h_3h_2^{-1}h_3^{-1}h_2h_3=1$, $R_3:h_2h_3^{-2}h_2^{-1}h_3^2=1$:\\
By interchanging $h_2$ and $h_3$ in (28) and with the same discussion, there is a contradiction.

\item[(45)]$ R_1:h_2h_3^{-1}h_2^{-1}h_3^2h_2h_3^{-1}h_2=1$, $R_2:h_2h_3^{-1}h_2^{-1}h_3h_2^{-1}h_3^{-1}h_2h_3=1$, $R_3:h_2h_3^{-2}h_2h_3^2=1$:\\
By interchanging $h_2$ and $h_3$ in (27) and with the same discussion, there is a contradiction.

\item[(46)]$ R_1:h_2h_3^{-1}h_2^{-1}h_3^2h_2^{-2}h_3=1$, $R_2:h_2h_3^{-2}h_2h_3^{-1}h_2^{-1}h_3^2=1$, $R_3:h_2h_3h_2^{-1}h_3^{-1}h_2^{-1}h_3=1$:\\
By interchanging $h_2$ and $h_3$ in (25) and with the same discussion, there is a contradiction.

\item[(47)]$ R_1:h_2h_3^{-1}h_2^{-1}h_3^2h_2^{-2}h_3=1$, $R_2:(h_2h_3^{-2})^2h_2h_3=1$, $R_3:h_2h_3h_2^{-1}h_3^{-2}h_2=1$:\\
By interchanging $h_2$ and $h_3$ in (24) and with the same discussion, there is a contradiction.

\item[(48)]$ R_1:h_2h_3^{-1}h_2^{-1}h_3^2h_2^{-1}h_3^{-1}h_2=1$, $R_2:h_2h_3^{-1}h_2^{-1}h_3h_2^{-2}h_3^2=1$, $R_3:h_2h_3^{-3}h_2h_3=1$:\\
By interchanging $h_2$ and $h_3$ in (35) and with the same discussion, there is a contradiction.

\item[(49)]$ R_1:h_2h_3^{-1}h_2^{-1}h_3^2h_2^{-1}h_3^{-1}h_2=1$, $R_2:h_2h_3^{-1}h_2^{-1}h_3h_2^{-2}h_3^2=1$, $R_3:h_2h_3^{-1}(h_3^{-1}h_2)^3=1$:\\
By interchanging $h_2$ and $h_3$ in (36) and with the same discussion, there is a contradiction.

\item[(50)]$ R_1:h_2h_3^{-1}h_2^{-1}h_3^2h_2^{-1}h_3^{-1}h_2=1$, $R_2:h_2h_3^{-1}h_2^{-1}h_3h_2^{-1}h_3^{-1}h_2h_3=1$, $R_3:h_2h_3^{-2}h_2^{-1}h_3^2=1$:\\
By interchanging $h_2$ and $h_3$ in (32) and with the same discussion, there is a contradiction.

\item[(51)]$ R_1:h_2h_3^{-1}h_2^{-1}h_3^2h_2^{-1}h_3^{-1}h_2=1$, $R_2:h_2h_3^{-1}h_2^{-1}h_3h_2^{-1}h_3^{-1}h_2h_3=1$, $R_3:h_2h_3^{-2}(h_3^{-1}h_2)^2=1$:\\
By interchanging $h_2$ and $h_3$ in (33) and with the same discussion, there is a contradiction.

\item[(52)]$ R_1:h_2h_3^{-1}h_2^{-1}h_3^2h_2^{-1}h_3^{-1}h_2=1$, $R_2:h_2h_3^{-1}h_2^{-1}h_3h_2^{-1}h_3^{-1}h_2h_3=1$, $R_3:h_2h_3^{-2}h_2h_3^2=1$:\\
By interchanging $h_2$ and $h_3$ in (31) and with the same discussion, there is a contradiction.

\item[(53)]$ R_1:h_2h_3^{-1}h_2^{-1}h_3^2h_2^{-1}h_3^{-1}h_2=1$, $R_2:h_2h_3^{-1}h_2^{-1}h_3h_2^{-1}h_3^{-1}h_2h_3=1$, $R_3:h_2h_3^{-1}(h_3^{-1}h_2)^3=1$:\\
By interchanging $h_2$ and $h_3$ in (34) and with the same discussion, there is a contradiction.

\item[(54)]$ R_1:h_2h_3^{-1}h_2^{-1}h_3^2h_2^{-1}h_3^{-1}h_2=1$, $R_2:h_2h_3^{-2}h_2^{-1}h_3h_2h_3^{-1}h_2=1$, $R_3:h_2h_3^{-3}h_2h_3=1$:\\
By interchanging $h_2$ and $h_3$ in (14) and with the same discussion, there is a contradiction.

\item[(55)]$ R_1:h_2h_3^{-1}h_2^{-1}h_3^2h_2^{-1}h_3^{-1}h_2=1$, $R_2:h_2h_3^{-3}h_2^2h_3^{-1}h_2=1$, $R_3:h_2h_3^{-2}h_2^{-1}h_3^2=1$:\\
By interchanging $h_2$ and $h_3$ in (17) and with the same discussion, there is a contradiction.

\item[(56)]$ R_1:h_2h_3^{-1}h_2^{-1}h_3h_2^{-1}h_3^{-1}h_2h_3=1$, $R_2:h_2h_3^{-2}h_2^{-1}h_3h_2^{-2}h_3=1$, $R_3:(h_2h_3^{-1})^2(h_2^{-1}h_3)^2=1$:\\
By interchanging $h_2$ and $h_3$ in (12) and with the same discussion, there is a contradiction.

\item[(57)]$ R_1:h_2h_3^{-1}h_2^{-1}h_3h_2^{-1}h_3^{-1}h_2h_3=1$, $R_2:h_2h_3^{-2}h_2^{-1}h_3h_2^{-2}h_3=1$, $R_3:(h_2h_3^{-1})^2h_2h_3h_2^{-1}h_3=1$:\\
By interchanging $h_2$ and $h_3$ in (11) and with the same discussion, there is a contradiction.

\item[(58)]$ R_1:h_2h_3^{-1}h_2^{-1}h_3h_2^{-1}h_3^{-1}h_2h_3=1$, $R_2:h_2h_3^{-1}(h_3^{-1}h_2)^2h_3^2=1$, $R_3:h_2^2h_3^{-2}h_2h_3=1$:\\
By interchanging $h_2$ and $h_3$ in (3) and with the same discussion, there is a contradiction.

\item[(59)]$ R_1:h_2h_3^{-1}h_2^{-1}h_3h_2^{-1}h_3^{-1}h_2h_3=1$, $R_2:(h_2h_3^{-1}h_2)^2h_3^2=1$, $R_3:(h_2h_3^{-1})^2(h_2^{-1}h_3)^2=1$:\\
By interchanging $h_2$ and $h_3$ in (8) and with the same discussion, there is a contradiction.

\item[(60)]$ R_1:h_2h_3^{-1}h_2^{-1}h_3h_2^{-1}h_3^{-1}h_2h_3=1$, $R_2:(h_2h_3^{-1}h_2)^2h_3^2=1$, $R_3:(h_2h_3^{-1})^2h_2h_3h_2^{-1}h_3=1$:\\
By interchanging $h_2$ and $h_3$ in (7) and with the same discussion, there is a contradiction.

\end{enumerate}

$\mathbf{C_6--C_6(--C_5-)}$ \textbf{subgraph:} 
By considering the $996$ cases related to the existence of $C_6--C_6$ in the graph $K(\alpha,\beta)$ and the relations from Table \ref{tab-C5} which are not disproved, it can be seen that there are $3267$ cases for the relations of two cycles $C_6$ and a cycle $C_5$ in this structure. Using GAP \cite{gap}, we see that all groups with two generators $h_2$ and $h_3$ and three relations which are between $2910$ cases of these $3267$ cases are finite or solvable, that is a contradiction with the assumptions. So, there are just $357$ cases for the relations of these cycles which may lead to the existence of a subgraph isomorphic to the graph $C_6--C_6(--C_5-)$ in $K(\alpha,\beta)$. Similar to the previous mentioned subgraphs it can be seen that $349$ cases of these relations lead to contradictions and  $8$ cases of them may lead to the  existence of a subgraph isomorphic to the graph $C_6--C_6(--C_5-)$ in the graph $K(\alpha,\beta)$.

\subsection{$\mathbf{C_6--C_6(--C_5-)(-C_5-)}$}
\begin{figure}[ht]
\psscalebox{0.9 0.9} 
{
\begin{pspicture}(0,-2.0985577)(5.01,2.0985577)
\psdots[linecolor=black, dotsize=0.4](2.4,1.9014423)
\psdots[linecolor=black, dotsize=0.4](2.4,0.70144236)
\psdots[linecolor=black, dotsize=0.4](2.4,-0.49855766)
\psdots[linecolor=black, dotsize=0.4](3.6,-0.098557666)
\psdots[linecolor=black, dotsize=0.4](3.6,0.70144236)
\psdots[linecolor=black, dotsize=0.4](3.6,1.5014423)
\psdots[linecolor=black, dotsize=0.4](1.2,1.5014423)
\psdots[linecolor=black, dotsize=0.4](1.2,0.70144236)
\psdots[linecolor=black, dotsize=0.4](1.2,-0.098557666)
\psline[linecolor=black, linewidth=0.04](2.4,1.9014423)(2.4,0.70144236)(2.4,-0.49855766)(3.6,-0.098557666)(3.6,0.70144236)(3.6,1.5014423)(2.4,1.9014423)(1.2,1.5014423)(1.2,0.70144236)(1.2,-0.098557666)(2.4,-0.49855766)
\psbezier[linecolor=black, linewidth=0.04](3.6,-0.098557666)(3.6,-0.89855766)(0.4,-1.6985576)(0.4,-0.89855766)
\psdots[linecolor=black, dotsize=0.4](0.4,-0.89855766)
\psline[linecolor=black, linewidth=0.04](1.2,0.70144236)(0.4,-0.89855766)
\psdots[linecolor=black, dotsize=0.4](3.6,-0.89855766)
\rput[bl](-0.5,-2.0985577){$\mathbf{30) \ C_6--C_6(--C_5-)(-C_5-)}$}
\psdots[linecolor=black, dotsize=0.4](4.8,-0.49855766)
\psline[linecolor=black, linewidth=0.04](2.4,0.70144236)(3.6,-0.89855766)(4.8,-0.49855766)(3.6,1.5014423)
\end{pspicture}
}
\end{figure}
By considering the $8$ cases related to the existence of $C_6--C_6(--C_5-)$ in the graph $K(\alpha,\beta)$ and the relations from Table \ref{tab-C5} which are not disproved, it can be seen that there are $62$ cases for the relations of two cycles $C_5$ and two cycles $C_6$ in this structure. Using GAP \cite{gap}, we see that all groups with two generators $h_2$ and $h_3$ and four relations which are between $56$ cases of these $62$ cases are finite or solvable, that is a contradiction with the assumptions. So, there are just $6$ cases for the relations of these cycles which may lead to the existence of a subgraph isomorphic to the graph $C_6--C_6(--C_5-)(-C_5-)$ in $K(\alpha,\beta)$. In the following, we show that these $6$ cases lead to contradictions and so, the graph $K(\alpha,\beta)$ contains no subgraph isomorphic to the graph $C_6--C_6(--C_5-)(-C_5-)$.
\begin{enumerate}
\item[(1)]$ R_1:h_2h_3^{-1}h_2^{-1}h_3^2h_2^{-1}h_3^{-1}h_2=1$, $R_2:h_2h_3^{-1}h_2^{-1}h_3^2h_2^{-1}h_3^{-1}h_2=1$, $R_3:(h_2h_3^{-1})^2h_2h_3^2=1$,\\$R_4:h_2^2h_3^{-1}h_2^2h_3=1$:\\
$\Rightarrow$  $G$ has a torsion element, a contradiction.

\item[(2)]$ R_1:h_2h_3^{-1}h_2^{-1}h_3^2h_2^{-1}h_3^{-1}h_2=1$, $R_2:h_2h_3^{-1}h_2^{-1}h_3^2h_2^{-1}h_3^{-1}h_2=1$, $R_3:(h_2h_3^{-1})^3h_2h_3=1$,\\$R_4:h_2h_3h_2^{-1}(h_3^{-1}h_2)^2=1$:\\
$\Rightarrow \langle h_2,h_3\rangle=\langle h_2\rangle$ is abelian, a contradiction.

\item[(3)]$ R_1:h_2h_3^{-1}h_2^{-1}h_3^2h_2^{-1}h_3^{-1}h_2=1$, $R_2:h_2h_3^{-1}h_2^{-1}h_3^2h_2^{-1}h_3^{-1}h_2=1$, $R_3:h_2^2(h_3h_2^{-1})^2h_3=1$,\\$R_4:h_2h_3^2h_2^{-1}h_3^2=1$:\\
By interchanging $h_2$ and $h_3$ in (1) and with the same discussion, there is a contradiction.

\item[(4)]$ R_1:h_2h_3^{-1}h_2^{-1}h_3^2h_2^{-1}h_3^{-1}h_2=1$, $R_2:h_2h_3^{-1}h_2^{-1}h_3^2h_2^{-1}h_3^{-1}h_2=1$, $R_3:h_2(h_3h_2^{-1})^3h_3=1$,\\$R_4:h_2h_3^{-1}(h_2^{-1}h_3)^2h_3=1$:\\
By interchanging $h_2$ and $h_3$ in (2) and with the same discussion, there is a contradiction.

\item[(5)]$ R_1:h_2h_3^{-1}h_2(h_3h_2^{-1})^2h_3^{-1}h_2=1$, $R_2:h_2h_3^{-1}h_2(h_3h_2^{-1})^2h_3^{-1}h_2=1$, $R_3:h_2h_3^{-2}h_2h_3^2=1$,\\$R_4:h_2^2h_3^{-1}h_2^2h_3=1$:\\
$\Rightarrow$  $G$ has a torsion element, a contradiction.

\item[(6)]$ R_1:(h_2h_3^{-1})^2h_2^{-1}h_3^2h_2^{-1}h_3=1$, $R_2:(h_2h_3^{-1})^2h_2^{-1}h_3^2h_2^{-1}h_3=1$, $R_3:h_2^2h_3h_2^{-2}h_3=1$,\\$R_4:h_2h_3^2h_2^{-1}h_3^2=1$:\\
By interchanging $h_2$ and $h_3$ in (5) and with the same discussion, there is a contradiction.

\end{enumerate}

\subsection{$\mathbf{C_6--C_6(--C_5-)(C_6---)}$}
\begin{figure}[ht]
\psscalebox{0.9 0.9} 
{
\begin{pspicture}(0,-2.0985577)(5.44,2.0985577)
\psdots[linecolor=black, dotsize=0.4](2.8,1.9014423)
\psdots[linecolor=black, dotsize=0.4](2.8,0.70144236)
\psdots[linecolor=black, dotsize=0.4](2.8,-0.49855766)
\psdots[linecolor=black, dotsize=0.4](4.0,-0.098557666)
\psdots[linecolor=black, dotsize=0.4](4.0,0.70144236)
\psdots[linecolor=black, dotsize=0.4](4.0,1.5014423)
\psdots[linecolor=black, dotsize=0.4](1.6,1.5014423)
\psdots[linecolor=black, dotsize=0.4](1.6,0.70144236)
\psdots[linecolor=black, dotsize=0.4](1.6,-0.098557666)
\psline[linecolor=black, linewidth=0.04](2.8,1.9014423)(2.8,0.70144236)(2.8,-0.49855766)(4.0,-0.098557666)(4.0,0.70144236)(4.0,1.5014423)(2.8,1.9014423)(1.6,1.5014423)(1.6,0.70144236)(1.6,-0.098557666)(2.8,-0.49855766)
\psbezier[linecolor=black, linewidth=0.04](4.0,-0.098557666)(4.0,-0.89855766)(0.8,-1.6985576)(0.8,-0.89855766)
\psdots[linecolor=black, dotsize=0.4](0.8,-0.89855766)
\psline[linecolor=black, linewidth=0.04](1.6,0.70144236)(0.8,-0.89855766)
\psdots[linecolor=black, dotsize=0.4](2.0,-0.49855766)
\psdots[linecolor=black, dotsize=0.4](1.6,-0.89855766)
\psline[linecolor=black, linewidth=0.04](2.8,0.70144236)(2.0,-0.49855766)(1.6,-0.89855766)(0.8,-0.89855766)
\rput[bl](-0.4,-2.0985577){$\mathbf{31) \ C_6--C_6(--C_5-)(C_6---)}$}
\end{pspicture}
}
\end{figure}
By considering the $8$ cases related to the existence of $C_6--C_6(--C_5-)$ in the graph $K(\alpha,\beta)$ and the relations from Table \ref{tab-C6} which are not disproved, it can be seen that there are $76$ cases for the relations of a cycle $C_5$ and three cycles $C_6$ in the graph $C_6--C_6(--C_5-)(C_6---)$. Using GAP \cite{gap}, we see that all groups with two generators $h_2$ and $h_3$ and four relations which are between $64$ cases of these $76$ cases are finite or solvable, that is a contradiction with the assumptions. So, there are just $12$ cases for the relations of these cycles which may lead to the existence of a subgraph isomorphic to the graph $C_6--C_6(--C_5-)(C_6---)$ in $K(\alpha,\beta)$. In the following, we show that these $12$ cases lead to contradictions and so, the graph $K(\alpha,\beta)$ contains no subgraph isomorphic to the graph $C_6--C_6(--C_5-)(C_6---)$.
\begin{enumerate}
\item[(1)]$ R_1:h_2h_3^{-1}h_2^{-1}h_3^2h_2^{-1}h_3^{-1}h_2=1$, $R_2:h_2h_3^{-1}h_2^{-1}h_3^2h_2^{-1}h_3^{-1}h_2=1$, $R_3:(h_2h_3^{-1})^2h_2h_3^2=1$,\\$R_4:h_2h_3^{-2}h_2^{-1}h_3h_2^{-2}h_3=1$:\\
$\Rightarrow$  $G$ has a torsion element, a contradiction.

\item[(2)]$ R_1:h_2h_3^{-1}h_2^{-1}h_3^2h_2^{-1}h_3^{-1}h_2=1$, $R_2:h_2h_3^{-1}h_2^{-1}h_3^2h_2^{-1}h_3^{-1}h_2=1$, $R_3:(h_2h_3^{-1})^2h_2h_3^2=1$,\\$R_4:(h_2h_3^{-1}h_2)^2h_3^2=1$:\\
$\Rightarrow$  $h_2=1$, a contradiction.

\item[(3)]$ R_1:h_2h_3^{-1}h_2^{-1}h_3^2h_2^{-1}h_3^{-1}h_2=1$, $R_2:h_2h_3^{-1}h_2^{-1}h_3^2h_2^{-1}h_3^{-1}h_2=1$, $R_3:(h_2h_3^{-1})^3h_2h_3=1$,\\$R_4:h_2^2h_3^{-2}(h_2^{-1}h_3)^2=1$:\\
$\Rightarrow\langle h_2,h_3\rangle\cong BS(1,1)$ is solvable, a contradiction.

\item[(4)]$ R_1:h_2h_3^{-1}h_2^{-1}h_3^2h_2^{-1}h_3^{-1}h_2=1$, $R_2:h_2h_3^{-1}h_2^{-1}h_3^2h_2^{-1}h_3^{-1}h_2=1$, $R_3:(h_2h_3^{-1})^3h_2h_3=1$,\\$R_4:h_2^2h_3^{-1}(h_3^{-1}h_2)^2h_3=1$:\\
$\Rightarrow\langle h_2,h_3\rangle\cong BS(1,1)$ is solvable, a contradiction.

\item[(5)]$ R_1:h_2h_3^{-1}h_2^{-1}h_3^2h_2^{-1}h_3^{-1}h_2=1$, $R_2:h_2h_3^{-1}h_2^{-1}h_3^2h_2^{-1}h_3^{-1}h_2=1$, $R_3:h_2^2(h_3h_2^{-1})^2h_3=1$,\\$R_4:h_2^2(h_3h_2^{-1}h_3)^2=1$:\\
By interchanging $h_2$ and $h_3$ in (2) and with the same discussion, there is a contradiction.

\item[(6)]$ R_1:h_2h_3^{-1}h_2^{-1}h_3^2h_2^{-1}h_3^{-1}h_2=1$, $R_2:h_2h_3^{-1}h_2^{-1}h_3^2h_2^{-1}h_3^{-1}h_2=1$, $R_3:h_2^2(h_3h_2^{-1})^2h_3=1$,\\$R_4:h_2^2h_3^{-1}h_2^{-1}h_3^2h_2^{-1}h_3=1$:\\
By interchanging $h_2$ and $h_3$ in (1) and with the same discussion, there is a contradiction.

\item[(7)]$ R_1:h_2h_3^{-1}h_2^{-1}h_3^2h_2^{-1}h_3^{-1}h_2=1$, $R_2:h_2h_3^{-1}h_2^{-1}h_3^2h_2^{-1}h_3^{-1}h_2=1$, $R_3:h_2(h_3h_2^{-1})^3h_3=1$,\\$R_4:h_2^2h_3^{-2}(h_2^{-1}h_3)^2=1$:\\
By interchanging $h_2$ and $h_3$ in (3) and with the same discussion, there is a contradiction.

\item[(8)]$ R_1:h_2h_3^{-1}h_2^{-1}h_3^2h_2^{-1}h_3^{-1}h_2=1$, $R_2:h_2h_3^{-1}h_2^{-1}h_3^2h_2^{-1}h_3^{-1}h_2=1$, $R_3:h_2(h_3h_2^{-1})^3h_3=1$,\\$R_4:h_2h_3^2h_2^{-1}(h_2^{-1}h_3)^2=1$:\\
By interchanging $h_2$ and $h_3$ in (4) and with the same discussion, there is a contradiction.

\item[(9)]$ R_1:h_2h_3^{-1}h_2(h_3h_2^{-1})^2h_3^{-1}h_2=1$, $R_2:h_2h_3^{-1}h_2(h_3h_2^{-1})^2h_3^{-1}h_2=1$, $R_3:h_2h_3^{-2}h_2h_3^2=1$,\\$R_4:h_2h_3h_2h_3^{-2}h_2h_3=1$:\\
$\Rightarrow$  $h_2=1$, a contradiction.

\item[(10)]$ R_1:h_2h_3^{-1}h_2(h_3h_2^{-1})^2h_3^{-1}h_2=1$, $R_2:h_2h_3^{-1}h_2(h_3h_2^{-1})^2h_3^{-1}h_2=1$, $R_3:h_2h_3^{-2}h_2h_3^2=1$,\\$R_4:h_2h_3h_2^{-1}h_3^{-2}h_2h_3=1$:\\
$\Rightarrow\langle h_2,h_3\rangle\cong BS(1,2)$ is solvable, a contradiction.

\item[(11)]$ R_1:(h_2h_3^{-1})^2h_2^{-1}h_3^2h_2^{-1}h_3=1$, $R_2:(h_2h_3^{-1})^2h_2^{-1}h_3^2h_2^{-1}h_3=1$, $R_3:h_2^2h_3h_2^{-2}h_3=1$,\\$R_4:(h_2h_3)^2h_2^{-2}h_3=1$:\\
By interchanging $h_2$ and $h_3$ in (9) and with the same discussion, there is a contradiction.

\item[(12)]$ R_1:(h_2h_3^{-1})^2h_2^{-1}h_3^2h_2^{-1}h_3=1$, $R_2:(h_2h_3^{-1})^2h_2^{-1}h_3^2h_2^{-1}h_3=1$, $R_3:h_2^2h_3h_2^{-2}h_3=1$,\\$R_4:h_2h_3h_2h_3^{-1}h_2^{-2}h_3=1$:\\
By interchanging $h_2$ and $h_3$ in (10) and with the same discussion, there is a contradiction.

\end{enumerate}

\subsection{$\mathbf{C_5(--C_6--)C_5(C_6)}$}
\begin{figure}[ht]
\psscalebox{0.9 0.9} 
{
\begin{pspicture}(0,-1.8985577)(4.1971154,1.8985577)
\psdots[linecolor=black, dotsize=0.4](1.6,1.7014424)
\psdots[linecolor=black, dotsize=0.4](2.8,1.7014424)
\psdots[linecolor=black, dotsize=0.4](3.2,0.90144235)
\psdots[linecolor=black, dotsize=0.4](2.8,0.10144234)
\psdots[linecolor=black, dotsize=0.4](1.6,0.10144234)
\psdots[linecolor=black, dotsize=0.4](1.2,0.90144235)
\psdots[linecolor=black, dotsize=0.4](0.4,1.3014424)
\psdots[linecolor=black, dotsize=0.4](0.4,0.5014423)
\psdots[linecolor=black, dotsize=0.4](4.0,0.5014423)
\psdots[linecolor=black, dotsize=0.4](4.0,1.3014424)
\psdots[linecolor=black, dotsize=0.4](1.2,-0.6985577)
\psdots[linecolor=black, dotsize=0.4](3.2,-0.6985577)
\psline[linecolor=black, linewidth=0.04](1.6,1.7014424)(1.2,0.90144235)(1.6,0.10144234)(2.8,0.10144234)(3.2,0.90144235)(2.8,1.7014424)(1.6,1.7014424)(0.4,1.3014424)(0.4,0.5014423)(1.6,0.10144234)
\psline[linecolor=black, linewidth=0.04](2.8,0.10144234)(4.0,0.5014423)(4.0,1.3014424)(2.8,1.7014424)
\psline[linecolor=black, linewidth=0.04](4.0,0.5014423)(3.2,-0.6985577)(1.2,-0.6985577)(0.4,0.5014423)
\rput[bl](-0.4,-1.8985577){$\mathbf{32) \ C_5(--C_6--)C_5(C_6)}$}
\end{pspicture}
}
\end{figure}
By considering the $6$ cases related to the existence of $C_5(--C_6--)C_5$ in the graph $K(\alpha,\beta)$ and the relations from Table \ref{tab-C6} which are not disproved, it can be seen that there are $120$ cases for the relations of two cycles $C_5$ and two cycles $C_6$ in this structure. Using GAP \cite{gap}, we see that all groups with two generators $h_2$ and $h_3$ and four relations which are between $104$ cases of these $120$ cases are finite and solvable, or just finite, that is a contradiction with the assumptions. So, there are just $16$ cases for the relations of these cycles which may lead to the existence of a subgraph isomorphic to the graph $C_5(--C_6--)C_5(C_6)$ in $K(\alpha,\beta)$.  In the following, we show that these $16$ cases lead to contradictions and so, the graph $K(\alpha,\beta)$ contains no subgraph isomorphic to the graph $C_5(--C_6--)C_5(C_6)$.
\begin{enumerate}
\item[(1)]$ R_1:h_2h_3^{-3}h_2h_3=1$, $R_2:h_2^2h_3^{-1}h_2^{-1}h_3^{-1}h_2=1$, $R_3:h_2h_3^{-1}h_2^{-1}h_3^2h_2^{-1}h_3^{-1}h_2=1$,\\$R_4:h_2(h_2h_3^{-1})^2h_2^{-2}h_3=1$:\\
$\Rightarrow$  $G$ has a torsion element, a contradiction.

\item[(2)]$ R_1:h_2h_3^{-3}h_2h_3=1$, $R_2:h_2^2h_3^{-1}h_2^{-1}h_3^{-1}h_2=1$, $R_3:h_2h_3^{-1}h_2^{-1}h_3^2h_2^{-1}h_3^{-1}h_2=1$,\\$R_4:h_2(h_2h_3^{-1})^2h_2^{-1}h_3^{-1}h_2=1$:\\
$\Rightarrow$  $h_2=1$, a contradiction.

\item[(3)]$ R_1:h_2^2h_3^{-1}h_2^{-1}h_3^{-1}h_2=1$, $R_2:h_2h_3^{-3}h_2h_3=1$, $R_3:h_2h_3^{-1}h_2^{-1}h_3^2h_2^{-1}h_3^{-1}h_2=1$,\\$R_4:(h_2h_3^{-1})^2h_3^{-1}h_2^{-1}h_3^2=1$:\\
By interchanging $h_2$ and $h_3$ in (1) and with the same discussion, there is a contradiction.

\item[(4)]$ R_1:h_2^2h_3^{-1}h_2^{-1}h_3^{-1}h_2=1$, $R_2:h_2h_3^{-3}h_2h_3=1$, $R_3:h_2h_3^{-1}h_2^{-1}h_3^2h_2^{-1}h_3^{-1}h_2=1$,\\$R_4:(h_2h_3^{-1})^2h_3^{-2}h_2h_3=1$:\\
By interchanging $h_2$ and $h_3$ in (2) and with the same discussion, there is a contradiction.

\item[(5)]$ R_1:h_2^2h_3^{-2}h_2^{-1}h_3=1$, $R_2:h_2h_3h_2^{-1}h_3^{-2}h_2=1$, $R_3:h_2h_3^{-1}h_2^{-1}h_3^2h_2^{-2}h_3=1$,\\$R_4:(h_2h_3^{-1})^2h_3^{-1}h_2^{-1}h_3^2=1$:\\
$\Rightarrow\langle h_2,h_3\rangle\cong BS(1,1)$ is solvable, a contradiction.

\item[(6)]$ R_1:h_2^2h_3^{-2}h_2^{-1}h_3=1$, $R_2:h_2h_3h_2^{-1}h_3^{-2}h_2=1$, $R_3:h_2h_3^{-1}h_2^{-1}h_3^2h_2^{-2}h_3=1$,\\$R_4:(h_2h_3^{-1})^2h_3^{-2}h_2h_3=1$:\\
$\Rightarrow\langle h_2,h_3\rangle\cong BS(1,1)$ is solvable, a contradiction.

\item[(7)]$ R_1:h_2h_3h_2^{-1}h_3^{-2}h_2=1$, $R_2:h_2^2h_3^{-2}h_2^{-1}h_3=1$, $R_3:h_2h_3^{-1}h_2^{-1}h_3^2h_2^{-2}h_3=1$,\\$R_4:h_2(h_2h_3^{-1})^2h_2^{-2}h_3=1$:\\
By interchanging $h_2$ and $h_3$ in (5) and with the same discussion, there is a contradiction.

\item[(8)]$ R_1:h_2h_3h_2^{-1}h_3^{-2}h_2=1$, $R_2:h_2^2h_3^{-2}h_2^{-1}h_3=1$, $R_3:h_2h_3^{-1}h_2^{-1}h_3^2h_2^{-2}h_3=1$,\\$R_4:h_2(h_2h_3^{-1})^2h_2^{-1}h_3^{-1}h_2=1$:\\
By interchanging $h_2$ and $h_3$ in (6) and with the same discussion, there is a contradiction.

\item[(9)]$ R_1:h_2h_3h_2^{-1}(h_3^{-1}h_2)^2=1$, $R_2:h_2(h_2h_3^{-1})^2h_2^{-1}h_3=1$, $R_3:h_2^2h_3^{-1}(h_3^{-1}h_2)^2h_3=1$,\\$R_4:h_2h_3^{-2}h_2h_3h_2h_3^{-1}h_2=1$:\\
$\Rightarrow\langle h_2,h_3\rangle\cong BS(1,2)$ is solvable, a contradiction.

\item[(10)]$ R_1:h_2h_3h_2^{-1}(h_3^{-1}h_2)^2=1$, $R_2:h_2(h_2h_3^{-1})^2h_2^{-1}h_3=1$, $R_3:h_2^2h_3^{-1}(h_3^{-1}h_2)^2h_3=1$,\\$R_4:h_2h_3^{-2}h_2^2h_3^{-1}h_2h_3=1$:\\
$\Rightarrow\langle h_2,h_3\rangle\cong BS(1,2)$ is solvable, a contradiction.

\item[(11)]$ R_1:h_2h_3h_2^{-1}(h_3^{-1}h_2)^2=1$, $R_2:h_2h_3h_2^{-1}(h_3^{-1}h_2)^2=1$, $R_3:h_2^3h_3^{-2}h_2^{-1}h_3=1$,\\$R_4:h_2h_3h_2h_3^{-1}(h_3^{-1}h_2)^2=1$:\\
$\Rightarrow \langle h_2,h_3\rangle=\langle h_2\rangle$ is abelian, a contradiction.

\item[(12)]$ R_1:h_2h_3h_2^{-1}(h_3^{-1}h_2)^2=1$, $R_2:h_2h_3h_2^{-1}(h_3^{-1}h_2)^2=1$, $R_3:h_2^3h_3^{-2}h_2^{-1}h_3=1$,\\$R_4:h_2h_3^{-1}(h_3^{-1}h_2)^2h_3h_2^{-1}h_3=1$:\\
$\Rightarrow \langle h_2,h_3\rangle=\langle h_2\rangle$ is abelian, a contradiction.

\item[(13)]$ R_1:h_2h_3^{-1}(h_2^{-1}h_3)^2h_3=1$, $R_2:h_2h_3^{-1}(h_2^{-1}h_3)^2h_3=1$, $R_3:h_2^2h_3^{-3}h_2^{-1}h_3=1$,\\$R_4:h_2h_3h_2^{-1}(h_2^{-1}h_3)^2h_3=1$:\\
By interchanging $h_2$ and $h_3$ in (11) and with the same discussion, there is a contradiction.

\item[(14)]$ R_1:h_2h_3^{-1}(h_2^{-1}h_3)^2h_3=1$, $R_2:h_2h_3^{-1}(h_2^{-1}h_3)^2h_3=1$, $R_3:h_2^2h_3^{-3}h_2^{-1}h_3=1$,\\$R_4:h_2h_3^{-1}h_2^{-1}h_3h_2^{-1}(h_3^{-1}h_2)^2=1$:\\
By interchanging $h_2$ and $h_3$ in (12) and with the same discussion, there is a contradiction.

\item[(15)]$ R_1:h_2h_3^{-1}(h_2^{-1}h_3)^2h_3=1$, $R_2:(h_2h_3^{-1})^2h_3^{-1}h_2^{-1}h_3=1$, $R_3:h_2h_3^2h_2^{-1}(h_2^{-1}h_3)^2=1$,\\$R_4:h_2h_3h_2^{-2}h_3^2h_2^{-1}h_3=1$:\\
By interchanging $h_2$ and $h_3$ in (10) and with the same discussion, there is a contradiction.

\item[(16)]$ R_1:h_2h_3^{-1}(h_2^{-1}h_3)^2h_3=1$, $R_2:(h_2h_3^{-1})^2h_3^{-1}h_2^{-1}h_3=1$, $R_3:h_2h_3^2h_2^{-1}(h_2^{-1}h_3)^2=1$,\\$R_4:h_2h_3h_2^{-1}h_3^2h_2^{-2}h_3=1$:\\
By interchanging $h_2$ and $h_3$ in (9) and with the same discussion, there is a contradiction.

\end{enumerate}

\subsection{$\mathbf{C_5(--C_6--)C_5(C_7)}$}
\begin{figure}[ht]
\psscalebox{0.9 0.9} 
{
\begin{pspicture}(0,-1.8985577)(4.794231,1.8985577)
\psdots[linecolor=black, dotsize=0.4](1.3971155,1.7014424)
\psdots[linecolor=black, dotsize=0.4](2.5971155,1.7014424)
\psdots[linecolor=black, dotsize=0.4](2.9971154,0.90144235)
\psdots[linecolor=black, dotsize=0.4](2.5971155,0.10144234)
\psdots[linecolor=black, dotsize=0.4](1.3971155,0.10144234)
\psdots[linecolor=black, dotsize=0.4](0.9971155,0.90144235)
\psdots[linecolor=black, dotsize=0.4](0.19711548,1.3014424)
\psdots[linecolor=black, dotsize=0.4](0.19711548,0.5014423)
\psdots[linecolor=black, dotsize=0.4](3.7971156,0.5014423)
\psdots[linecolor=black, dotsize=0.4](3.7971156,1.3014424)
\psline[linecolor=black, linewidth=0.04](1.3971155,1.7014424)(0.9971155,0.90144235)(1.3971155,0.10144234)(2.5971155,0.10144234)(2.9971154,0.90144235)(2.5971155,1.7014424)(1.3971155,1.7014424)(0.19711548,1.3014424)(0.19711548,0.5014423)(1.3971155,0.10144234)
\psline[linecolor=black, linewidth=0.04](2.5971155,0.10144234)(3.7971156,0.5014423)(3.7971156,1.3014424)(2.5971155,1.7014424)
\psdots[linecolor=black, dotsize=0.4](0.59711546,-0.6985577)
\psdots[linecolor=black, dotsize=0.4](4.5971155,-0.6985577)
\psline[linecolor=black, linewidth=0.04](0.19711548,0.5014423)(0.59711546,-0.6985577)(4.5971155,-0.6985577)(3.7971156,1.3014424)
\rput[bl](-0.1,-1.8985577){$\mathbf{33) \ C_5(--C_6--)C_5(C_7)}$}
\end{pspicture}
}
\end{figure}
By considering the $6$ cases related to the existence of $C_5(--C_6--)C_5$ in the graph $K(\alpha,\beta)$ and the relations of $C_7$ cycles, it can be seen that there are $248$ cases for the relations of two cycles $C_5$, a cycle $C_6$ and a cycle $C_7$ in this structure. Using GAP \cite{gap}, we see that all groups with two generators $h_2$ and $h_3$ and four relations which are between $220$ cases of these $248$ cases are finite and solvable, or just finite, that is a contradiction with the assumptions. So, there are just $28$ cases for the relations of these cycles which may lead to the existence of a subgraph isomorphic to the graph $C_5(--C_6--)C_5(C_7)$ in $K(\alpha,\beta)$. In the following, we show that these $28$ cases lead to contradictions and so, the graph $K(\alpha,\beta)$ contains no subgraph isomorphic to the graph $C_5(--C_6--)C_5(C_7)$.
\begin{enumerate}
\item[(1)]$ R_1:h_2h_3^{-3}h_2h_3=1$, $R_2:h_2^2h_3^{-1}h_2^{-1}h_3^{-1}h_2=1$, $R_3:h_2h_3^{-1}h_2^{-1}h_3^2h_2^{-1}h_3^{-1}h_2=1$,\\$R_4:h_2^2(h_3^{-1}h_2h_3^{-1})^2h_2^{-1}h_3=1$:\\
$\Rightarrow\langle h_2,h_3\rangle\cong BS(-5,1)$ is solvable, a contradiction.

\item[(2)]$ R_1:h_2h_3^{-3}h_2h_3=1$, $R_2:h_2^2h_3^{-1}h_2^{-1}h_3^{-1}h_2=1$, $R_3:h_2h_3^{-1}h_2^{-1}h_3^2h_2^{-1}h_3^{-1}h_2=1$,\\$R_4:h_2^2(h_3^{-1}h_2h_3^{-1})^2h_3^{-1}h_2=1$:\\
$\Rightarrow\langle h_2,h_3\rangle\cong BS(-2,1)$ is solvable, a contradiction.

\item[(3)]$ R_1:h_2h_3^{-3}h_2h_3=1$, $R_2:h_2^2h_3^{-1}h_2^{-1}h_3^{-1}h_2=1$, $R_3:h_2h_3^{-1}h_2^{-1}h_3^2h_2^{-1}h_3^{-1}h_2=1$,\\$R_4:h_2(h_3h_2^{-1})^2h_3^{-1}h_2h_3^{-1}h_2^{-1}h_3=1$:\\
$\Rightarrow\langle h_2,h_3\rangle\cong BS(2,1)$ is solvable, a contradiction.

\item[(4)]$ R_1:h_2h_3^{-3}h_2h_3=1$, $R_2:h_2^2h_3^{-1}h_2^{-1}h_3^{-1}h_2=1$, $R_3:h_2h_3^{-1}h_2^{-1}h_3^2h_2^{-1}h_3^{-1}h_2=1$,\\$R_4:h_2(h_3h_2^{-1})^2h_3^{-1}h_2h_3^{-2}h_2=1$:\\
$\Rightarrow\langle h_2,h_3\rangle\cong BS(5,1)$ is solvable, a contradiction.

\item[(5)]$ R_1:h_2h_3^{-3}h_2h_3=1$, $R_2:h_2^2h_3^{-1}h_2^{-1}h_3^{-1}h_2=1$, $R_3:h_2h_3^{-1}h_2^{-1}h_3^2h_2^{-1}h_3^{-1}h_2=1$,\\$R_4:h_2h_3^{-1}h_2h_3^{-1}(h_2^{-1}h_3)^4=1$:\\
$\Rightarrow$  $G$ has a torsion element, a contradiction.

\item[(6)]$ R_1:h_2h_3^{-3}h_2h_3=1$, $R_2:h_2^2h_3^{-1}h_2^{-1}h_3^{-1}h_2=1$, $R_3:h_2h_3^{-1}h_2^{-1}h_3^2h_2^{-1}h_3^{-1}h_2=1$,\\$R_4:h_2h_3^{-1}h_2h_3^{-1}(h_3^{-1}h_2)^4=1$:\\
$\Rightarrow$  $G$ has a torsion element, a contradiction.

\item[(7)]$ R_1:h_2^2h_3^{-1}h_2^{-1}h_3^{-1}h_2=1$, $R_2:h_2h_3^{-3}h_2h_3=1$, $R_3:h_2h_3^{-1}h_2^{-1}h_3^2h_2^{-1}h_3^{-1}h_2=1$,\\$R_4:(h_2h_3^{-1}h_2)^2h_3^{-2}h_2^{-1}h_3=1$:\\
By interchanging $h_2$ and $h_3$ in (1) and with the same discussion, there is a contradiction.

\item[(8)]$ R_1:h_2^2h_3^{-1}h_2^{-1}h_3^{-1}h_2=1$, $R_2:h_2h_3^{-3}h_2h_3=1$, $R_3:h_2h_3^{-1}h_2^{-1}h_3^2h_2^{-1}h_3^{-1}h_2=1$,\\$R_4:(h_2h_3^{-1}h_2)^2h_3^{-3}h_2=1$:\\
By interchanging $h_2$ and $h_3$ in (2) and with the same discussion, there is a contradiction.

\item[(9)]$ R_1:h_2^2h_3^{-1}h_2^{-1}h_3^{-1}h_2=1$, $R_2:h_2h_3^{-3}h_2h_3=1$, $R_3:h_2h_3^{-1}h_2^{-1}h_3^2h_2^{-1}h_3^{-1}h_2=1$,\\$R_4:h_2h_3^{-1}h_2(h_3h_2^{-1})^2h_3^{-1}h_2^{-1}h_3=1$:\\
By interchanging $h_2$ and $h_3$ in (3) and with the same discussion, there is a contradiction.

\item[(10)]$ R_1:h_2^2h_3^{-1}h_2^{-1}h_3^{-1}h_2=1$, $R_2:h_2h_3^{-3}h_2h_3=1$, $R_3:h_2h_3^{-1}h_2^{-1}h_3^2h_2^{-1}h_3^{-1}h_2=1$,\\$R_4:h_2h_3^{-1}h_2(h_3h_2^{-1})^2h_3^{-2}h_2=1$:\\
By interchanging $h_2$ and $h_3$ in (4) and with the same discussion, there is a contradiction.

\item[(11)]$ R_1:h_2^2h_3^{-2}h_2^{-1}h_3=1$, $R_2:h_2h_3h_2^{-1}h_3^{-2}h_2=1$, $R_3:h_2h_3^{-1}h_2^{-1}h_3^2h_2^{-2}h_3=1$,\\$R_4:(h_2h_3^{-1}h_2)^2h_3^{-2}h_2^{-1}h_3=1$:\\
$\Rightarrow\langle h_2,h_3\rangle\cong BS(1,-1)$ is solvable, a contradiction.

\item[(12)]$ R_1:h_2^2h_3^{-2}h_2^{-1}h_3=1$, $R_2:h_2h_3h_2^{-1}h_3^{-2}h_2=1$, $R_3:h_2h_3^{-1}h_2^{-1}h_3^2h_2^{-2}h_3=1$,\\$R_4:(h_2h_3^{-1}h_2)^2h_3^{-3}h_2=1$:\\
$\Rightarrow\langle h_2,h_3\rangle\cong BS(1,-1)$ is solvable, a contradiction.

\item[(13)]$ R_1:h_2^2h_3^{-2}h_2^{-1}h_3=1$, $R_2:h_2h_3h_2^{-1}h_3^{-2}h_2=1$, $R_3:h_2h_3^{-1}h_2^{-1}h_3^2h_2^{-2}h_3=1$,\\$R_4:h_2h_3^{-1}h_2(h_3h_2^{-1})^2h_3^{-1}h_2^{-1}h_3=1$:\\
$\Rightarrow\langle h_2,h_3\rangle\cong BS(2,1)$ is solvable, a contradiction.

\item[(14)]$ R_1:h_2^2h_3^{-2}h_2^{-1}h_3=1$, $R_2:h_2h_3h_2^{-1}h_3^{-2}h_2=1$, $R_3:h_2h_3^{-1}h_2^{-1}h_3^2h_2^{-2}h_3=1$,\\$R_4:h_2h_3^{-1}h_2(h_3h_2^{-1})^2h_3^{-2}h_2=1$:\\
$\Rightarrow\langle h_2,h_3\rangle\cong BS(1,1)$ is solvable, a contradiction.

\item[(15)]$ R_1:h_2^2h_3^{-2}h_2^{-1}h_3=1$, $R_2:h_2h_3h_2^{-1}h_3^{-2}h_2=1$, $R_3:h_2h_3^{-1}h_2^{-1}h_3^2h_2^{-2}h_3=1$,\\$R_4:h_2h_3^{-1}h_2h_3^{-1}(h_2^{-1}h_3)^4=1$:\\
$\Rightarrow$  $G$ has a torsion element, a contradiction.

\item[(16)]$ R_1:h_2^2h_3^{-2}h_2^{-1}h_3=1$, $R_2:h_2h_3h_2^{-1}h_3^{-2}h_2=1$, $R_3:h_2h_3^{-1}h_2^{-1}h_3^2h_2^{-2}h_3=1$,\\$R_4:h_2h_3^{-1}h_2h_3^{-1}(h_3^{-1}h_2)^4=1$:\\
$\Rightarrow$  $G$ has a torsion element, a contradiction.

\item[(17)]$ R_1:h_2h_3h_2^{-1}h_3^{-2}h_2=1$, $R_2:h_2^2h_3^{-2}h_2^{-1}h_3=1$, $R_3:h_2h_3^{-1}h_2^{-1}h_3^2h_2^{-2}h_3=1$,\\$R_4:h_2^2(h_3^{-1}h_2h_3^{-1})^2h_2^{-1}h_3=1$:\\
By interchanging $h_2$ and $h_3$ in (11) and with the same discussion, there is a contradiction.

\item[(18)]$ R_1:h_2h_3h_2^{-1}h_3^{-2}h_2=1$, $R_2:h_2^2h_3^{-2}h_2^{-1}h_3=1$, $R_3:h_2h_3^{-1}h_2^{-1}h_3^2h_2^{-2}h_3=1$,\\$R_4:h_2^2(h_3^{-1}h_2h_3^{-1})^2h_3^{-1}h_2=1$:\\
By interchanging $h_2$ and $h_3$ in (12) and with the same discussion, there is a contradiction.

\item[(19)]$ R_1:h_2h_3h_2^{-1}h_3^{-2}h_2=1$, $R_2:h_2^2h_3^{-2}h_2^{-1}h_3=1$, $R_3:h_2h_3^{-1}h_2^{-1}h_3^2h_2^{-2}h_3=1$,\\$R_4:h_2(h_3h_2^{-1})^2h_3^{-1}h_2h_3^{-1}h_2^{-1}h_3=1$:\\
By interchanging $h_2$ and $h_3$ in (13) and with the same discussion, there is a contradiction.

\item[(20)]$ R_1:h_2h_3h_2^{-1}h_3^{-2}h_2=1$, $R_2:h_2^2h_3^{-2}h_2^{-1}h_3=1$, $R_3:h_2h_3^{-1}h_2^{-1}h_3^2h_2^{-2}h_3=1$,\\$R_4:h_2(h_3h_2^{-1})^2h_3^{-1}h_2h_3^{-2}h_2=1$:\\
By interchanging $h_2$ and $h_3$ in (14) and with the same discussion, there is a contradiction.

\item[(21)]$ R_1:h_2h_3h_2^{-1}(h_3^{-1}h_2)^2=1$, $R_2:h_2(h_2h_3^{-1})^2h_2^{-1}h_3=1$, $R_3:h_2^2h_3^{-1}(h_3^{-1}h_2)^2h_3=1$,\\$R_4:h_2h_3^{-1}(h_2^{-1}h_3h_2^{-1})^2h_2^{-1}h_3=1$:\\
$\Rightarrow\langle h_2,h_3\rangle\cong BS(1,1)$ is solvable, a contradiction.

\item[(22)]$ R_1:h_2h_3h_2^{-1}(h_3^{-1}h_2)^2=1$, $R_2:h_2(h_2h_3^{-1})^2h_2^{-1}h_3=1$, $R_3:h_2^2h_3^{-1}(h_3^{-1}h_2)^2h_3=1$,\\$R_4:h_2h_3^{-1}h_2^{-1}h_3(h_2h_3^{-1}h_2)^2=1$:\\
$\Rightarrow\langle h_2,h_3\rangle\cong BS(2,1)$ is solvable, a contradiction.

\item[(23)]$ R_1:h_2h_3h_2^{-1}(h_3^{-1}h_2)^2=1$, $R_2:h_2h_3h_2^{-1}(h_3^{-1}h_2)^2=1$, $R_3:h_2^3h_3^{-2}h_2^{-1}h_3=1$,\\$R_4:h_2h_3^{-1}h_2^{-1}h_3^2h_2h_3^{-2}h_2=1$:\\
$\Rightarrow\langle h_2,h_3\rangle\cong BS(1,-1)$ is solvable, a contradiction.

\item[(24)]$ R_1:h_2h_3h_2^{-1}(h_3^{-1}h_2)^2=1$, $R_2:h_2h_3h_2^{-1}(h_3^{-1}h_2)^2=1$, $R_3:h_2^3h_3^{-2}h_2^{-1}h_3=1$,\\$R_4:(h_2h_3^{-1})^3h_3^{-1}h_2h_3^2=1$:\\
$\Rightarrow \langle h_2,h_3\rangle=\langle h_2\rangle$ is abelian, a contradiction.

\item[(25)]$ R_1:h_2h_3^{-1}(h_2^{-1}h_3)^2h_3=1$, $R_2:h_2h_3^{-1}(h_2^{-1}h_3)^2h_3=1$, $R_3:h_2^2h_3^{-3}h_2^{-1}h_3=1$,\\$R_4:h_2^2(h_3h_2^{-1})^3h_2^{-1}h_3=1$:\\
By interchanging $h_2$ and $h_3$ in (24) and with the same discussion, there is a contradiction.

\item[(26)]$ R_1:h_2h_3^{-1}(h_2^{-1}h_3)^2h_3=1$, $R_2:h_2h_3^{-1}(h_2^{-1}h_3)^2h_3=1$, $R_3:h_2^2h_3^{-3}h_2^{-1}h_3=1$,\\$R_4:h_2h_3h_2^{-2}h_3^2h_2^{-1}h_3^{-1}h_2=1$:\\
By interchanging $h_2$ and $h_3$ in (23) and with the same discussion, there is a contradiction.

\item[(27)]$ R_1:h_2h_3^{-1}(h_2^{-1}h_3)^2h_3=1$, $R_2:(h_2h_3^{-1})^2h_3^{-1}h_2^{-1}h_3=1$, $R_3:h_2h_3^2h_2^{-1}(h_2^{-1}h_3)^2=1$,\\$R_4:(h_2h_3^{-2})^2h_2h_3^{-1}h_2^{-1}h_3=1$:\\
By interchanging $h_2$ and $h_3$ in (22) and with the same discussion, there is a contradiction.

\item[(28)]$ R_1:h_2h_3^{-1}(h_2^{-1}h_3)^2h_3=1$, $R_2:(h_2h_3^{-1})^2h_3^{-1}h_2^{-1}h_3=1$, $R_3:h_2h_3^2h_2^{-1}(h_2^{-1}h_3)^2=1$,\\$R_4:h_2h_3^{-1}(h_2^{-1}h_3^2)^2h_2^{-1}h_3=1$:\\
By interchanging $h_2$ and $h_3$ in (21) and with the same discussion, there is a contradiction.

\end{enumerate}

\begin{figure}[ht]
\psscalebox{0.9 0.9} 
{
\begin{pspicture}(0,-2.0985577)(3.7971156,2.0985577)
\psdots[linecolor=black, dotsize=0.4](0.4,1.1014423)
\psdots[linecolor=black, dotsize=0.4](1.2,1.9014423)
\psdots[linecolor=black, dotsize=0.4](1.2,0.30144233)
\psdots[linecolor=black, dotsize=0.4](2.0,0.70144236)
\psdots[linecolor=black, dotsize=0.4](2.0,1.5014423)
\psdots[linecolor=black, dotsize=0.4](2.8,1.9014423)
\psdots[linecolor=black, dotsize=0.4](3.6,1.1014423)
\psdots[linecolor=black, dotsize=0.4](2.8,0.30144233)
\psdots[linecolor=black, dotsize=0.4](1.2,-0.89855766)
\psdots[linecolor=black, dotsize=0.4](2.8,-0.89855766)
\psline[linecolor=black, linewidth=0.04](2.0,1.5014423)(1.2,1.9014423)(0.4,1.1014423)(1.2,0.30144233)(2.0,0.70144236)(2.0,1.5014423)(2.8,1.9014423)(3.6,1.1014423)(2.8,0.30144233)(2.0,0.70144236)
\psline[linecolor=black, linewidth=0.04](0.4,1.1014423)(1.2,-0.89855766)(2.8,-0.89855766)(3.6,1.1014423)
\rput[bl](0.0,-2.0985577){$\mathbf{C_5-C_5(--C_7--)}$}
\end{pspicture}
}
\end{figure}
$\mathbf{C_5-C_5(--C_7--)}$ \textbf{subgraph:} 
By considering the relations from Table \ref{tab-C5-C5} which are not disproved and the relations of $C_7$ cycles, it can be seen that there are $1981$ cases for the relations of two cycles $C_5$ and a cycle $C_7$ in this structure. Using GAP \cite{gap}, we see that all groups with two generators $h_2$ and $h_3$ and three relations which are between $1785$ cases of these $1981$ cases are finite or solvable, that is a contradiction with the assumptions. So, there are just $196$ cases for the relations of these cycles which may lead to the existence of a subgraph isomorphic to the graph $C_5-C_5(--C_7--)$ in $K(\alpha,\beta)$. Similar to the previous mentioned subgraphs it can be seen that $177$ cases of these relations lead to contradictions and  $19$ cases of them may lead to the  existence of a subgraph isomorphic to the graph $C_5-C_5(--C_7--)$ in the graph $K(\alpha,\beta)$.

\subsection{$\mathbf{C_5-C_5(--C_7--)(--C_5)}$}
\begin{figure}[ht]
\psscalebox{0.9 0.9} 
{
\begin{pspicture}(0,-1.8985577)(5.32,1.8985577)
\psdots[linecolor=black, dotsize=0.4](2.4,1.3014424)
\psdots[linecolor=black, dotsize=0.4](2.4,0.5014423)
\psdots[linecolor=black, dotsize=0.4](3.6,1.7014424)
\psdots[linecolor=black, dotsize=0.4](4.4,0.90144235)
\psdots[linecolor=black, dotsize=0.4](3.6,0.10144234)
\psdots[linecolor=black, dotsize=0.4](1.2,0.10144234)
\psdots[linecolor=black, dotsize=0.4](0.4,0.90144235)
\psdots[linecolor=black, dotsize=0.4](1.2,1.7014424)
\psline[linecolor=black, linewidth=0.04](1.2,1.7014424)(2.4,1.3014424)(3.6,1.7014424)(4.4,0.90144235)(3.6,0.10144234)(2.4,0.5014423)(2.4,1.3014424)
\psline[linecolor=black, linewidth=0.04](2.4,0.5014423)(1.2,0.10144234)(0.4,0.90144235)(1.2,1.7014424)
\psdots[linecolor=black, dotsize=0.4](1.2,-0.6985577)
\psdots[linecolor=black, dotsize=0.4](3.6,-0.6985577)
\psdots[linecolor=black, dotsize=0.4](4.8,-0.29855767)
\psdots[linecolor=black, dotsize=0.4](4.8,1.3014424)
\psline[linecolor=black, linewidth=0.04](0.4,0.90144235)(1.2,-0.6985577)(3.6,-0.6985577)(4.4,0.90144235)
\psline[linecolor=black, linewidth=0.04](3.6,1.7014424)(4.8,1.3014424)(4.8,-0.29855767)(3.6,-0.6985577)
\rput[bl](-0.3,-1.8985577){$\mathbf{34) \ C_5-C_5(--C_7--)(--C_5)}$}
\end{pspicture}
}
\end{figure}
By considering the $19$ cases related to the existence of $C_5-C_5(--C_7--)$ in the graph $K(\alpha,\beta)$ and the relations from Table \ref{tab-C5} which are not disproved, it can be seen that there are $256$ cases for the relations of a cycle $C_7$ and three cycles $C_5$ in the graph $C_5-C_5(--C_7--)(--C_5)$. Using GAP \cite{gap}, we see that all groups with two generators $h_2$ and $h_3$ and four relations which are between $232$ cases of these $256$ cases are finite or solvable, that is a contradiction with the assumptions. So, there are just $24$ cases for the relations of these cycles which may lead to the existence of a subgraph isomorphic to the graph $C_5-C_5(--C_7--)(--C_5)$ in $K(\alpha,\beta)$. In the following, we show that these $24$ cases lead to contradictions and so, the graph $K(\alpha,\beta)$ contains no subgraph isomorphic to the graph $C_5-C_5(--C_7--)(--C_5)$.
\begin{enumerate}
\item[(1)]$ R_1:h_2^2h_3h_2^{-2}h_3=1$, $R_2:h_2^2h_3h_2^{-2}h_3=1$, $R_3:(h_2h_3)^2(h_3h_2^{-1})^2h_3=1$,\\$R_4:h_2h_3^2h_2^{-1}h_3^2=1$:\\
$\Rightarrow$  $G$ has a torsion element, a contradiction.

\item[(2)]$ R_1:h_2^2h_3^{-1}h_2^{-2}h_3=1$, $R_2:h_2^2h_3^{-1}h_2^{-2}h_3=1$, $R_3:h_2h_3h_2h_3^{-2}(h_2^{-1}h_3)^2=1$,\\$R_4:h_2h_3^{-2}h_2^{-1}h_3^2=1$:\\
$\Rightarrow$ $G$ is solvable, a contradiction.

\item[(3)]$ R_1:h_2^2h_3^{-1}h_2^{-2}h_3=1$, $R_2:h_2^2h_3^{-1}h_2^{-2}h_3=1$, $R_3:h_2h_3h_2h_3^{-3}h_2^{-1}h_3=1$,\\$R_4:h_2h_3h_2^{-1}h_3^{-1}h_2^{-1}h_3=1$:\\
$\Rightarrow\langle h_2,h_3\rangle\cong BS(1,2)$ is solvable, a contradiction.

\item[(4)]$ R_1:h_2^2h_3^{-1}h_2^{-2}h_3=1$, $R_2:h_2^2h_3^{-1}h_2^{-2}h_3=1$, $R_3:h_2h_3h_2h_3^{-3}h_2^{-1}h_3=1$,\\$R_4:h_2h_3h_2^{-1}h_3^{-2}h_2=1$:\\
$\Rightarrow \langle h_2,h_3\rangle=\langle h_2\rangle$ is abelian, a contradiction.

\item[(5)]$ R_1:h_2^2h_3^{-1}h_2^{-2}h_3=1$, $R_2:h_2^2h_3^{-1}h_2^{-2}h_3=1$, $R_3:h_2h_3h_2h_3^{-2}h_2h_3^2=1$,\\$R_4:(h_2h_3^{-1})^2(h_2^{-1}h_3)^2=1$:\\
$\Rightarrow$ $G$ is solvable, a contradiction.

\item[(6)]$ R_1:h_2^2h_3^{-1}h_2^{-2}h_3=1$, $R_2:h_2^2h_3^{-1}h_2^{-2}h_3=1$, $R_3:h_2h_3h_2h_3^{-2}h_2h_3h_2^{-1}h_3=1$,\\$R_4:(h_2h_3^{-1})^2(h_2^{-1}h_3)^2=1$:\\
$\Rightarrow$ $G$ is solvable, a contradiction.

\item[(7)]$ R_1:h_2^2h_3^{-1}h_2^{-2}h_3=1$, $R_2:h_2^2h_3^{-1}h_2^{-2}h_3=1$, $R_3:h_2h_3h_2h_3^{-1}(h_3^{-1}h_2)^2h_3=1$,\\$R_4:(h_2h_3^{-1}h_2)^2h_3=1$:\\
$\Rightarrow\langle h_2,h_3\rangle\cong BS(1,-1)$ is solvable, a contradiction.

\item[(8)]$ R_1:h_2^2h_3^{-1}h_2^{-1}h_3^{-1}h_2=1$, $R_2:h_2^2h_3^{-1}h_2^{-1}h_3^{-1}h_2=1$, $R_3:h_2h_3^2h_2^{-1}h_3^{-3}h_2=1$,\\$R_4:h_2h_3h_2^{-1}h_3^{-1}h_2^{-1}h_3=1$:\\
$\Rightarrow \langle h_2,h_3\rangle=\langle h_2\rangle$ is abelian, a contradiction.

\item[(9)]$ R_1:h_2^2h_3^{-1}h_2^{-1}h_3^{-1}h_2=1$, $R_2:h_2^2h_3^{-1}h_2^{-1}h_3^{-1}h_2=1$, $R_3:h_2h_3^2h_2^{-1}h_3^{-3}h_2=1$,\\$R_4:h_2h_3h_2^{-1}h_3^{-2}h_2=1$:\\
$\Rightarrow\langle h_2,h_3\rangle\cong BS(-2,1)$ is solvable, a contradiction.

\item[(10)]$ R_1:h_2^2h_3^{-2}h_2^{-1}h_3=1$, $R_2:h_2^2h_3^{-2}h_2^{-1}h_3=1$, $R_3:h_2^2h_3h_2^{-1}h_3^{-3}h_2=1$,\\$R_4:h_2h_3^{-2}h_2^{-1}h_3^2=1$:\\
$\Rightarrow \langle h_2,h_3\rangle=\langle h_3\rangle$ is abelian, a contradiction.

\item[(11)]$ R_1:h_2^2h_3^{-2}h_2^{-1}h_3=1$, $R_2:h_2^2h_3^{-2}h_2^{-1}h_3=1$, $R_3:h_2^2h_3h_2^{-1}h_3^{-3}h_2=1$,\\$R_4:h_2h_3^{-3}h_2h_3=1$:\\
$\Rightarrow\langle h_2,h_3\rangle\cong BS(2,1)$ is solvable, a contradiction.

\item[(12)]$ R_1:h_2^2h_3^{-2}h_2^{-1}h_3=1$, $R_2:h_2^2h_3^{-2}h_2^{-1}h_3=1$, $R_3:h_2^2h_3h_2^{-1}h_3^{-3}h_2=1$,\\$R_4:h_2^2h_3^{-1}h_2^{-2}h_3=1$:\\
$\Rightarrow \langle h_2,h_3\rangle=\langle h_2\rangle$ is abelian, a contradiction.

\item[(13)]$ R_1:h_2^2h_3^{-2}h_2^{-1}h_3=1$, $R_2:h_2^2h_3^{-2}h_2^{-1}h_3=1$, $R_3:h_2^2h_3h_2^{-1}h_3^{-3}h_2=1$,\\$R_4:h_2^2h_3^{-1}h_2^{-1}h_3^{-1}h_2=1$:\\
$\Rightarrow\langle h_2,h_3\rangle\cong BS(2,1)$ is solvable, a contradiction.

\item[(14)]$ R_1:h_2h_3h_2h_3^{-2}h_2=1$, $R_2:h_2h_3h_2h_3^{-2}h_2=1$, $R_3:h_2^2(h_3^{-1}h_2^{-1})^2h_3^2=1$,\\$R_4:(h_2h_3^{-1})^2(h_2^{-1}h_3)^2=1$:\\
$\Rightarrow\langle h_2,h_3\rangle\cong BS(2,1)$ is solvable, a contradiction.

\item[(15)]$ R_1:h_2h_3h_2^{-2}h_3^2=1$, $R_2:h_2h_3h_2^{-2}h_3^2=1$, $R_3:h_2h_3h_2h_3^{-2}h_2^{-2}h_3=1$,\\$R_4:(h_2h_3^{-1})^2(h_2^{-1}h_3)^2=1$:\\
By interchanging $h_2$ and $h_3$ in (14) and with the same discussion, there is a contradiction.

\item[(16)]$ R_1:h_2h_3^{-2}h_2^{-1}h_3^2=1$, $R_2:h_2h_3^{-2}h_2^{-1}h_3^2=1$, $R_3:h_2^3(h_3^{-1}h_2^{-1})^2h_3=1$,\\$R_4:h_2h_3h_2^{-1}h_3^{-1}h_2^{-1}h_3=1$:\\
By interchanging $h_2$ and $h_3$ in (3) and with the same discussion, there is a contradiction.

\item[(17)]$ R_1:h_2h_3^{-2}h_2^{-1}h_3^2=1$, $R_2:h_2h_3^{-2}h_2^{-1}h_3^2=1$, $R_3:h_2^3(h_3^{-1}h_2^{-1})^2h_3=1$,\\$R_4:h_2h_3h_2^{-1}h_3^{-2}h_2=1$:\\
By interchanging $h_2$ and $h_3$ in (4) and with the same discussion, there is a contradiction.

\item[(18)]$ R_1:h_2h_3^{-2}h_2^{-1}h_3^2=1$, $R_2:h_2h_3^{-2}h_2^{-1}h_3^2=1$, $R_3:h_2(h_2h_3)^2h_2^{-2}h_3=1$,\\$R_4:(h_2h_3^{-1})^2(h_2^{-1}h_3)^2=1$:\\
By interchanging $h_2$ and $h_3$ in (5) and with the same discussion, there is a contradiction.

\item[(19)]$ R_1:h_2h_3^{-2}h_2^{-1}h_3^2=1$, $R_2:h_2h_3^{-2}h_2^{-1}h_3^2=1$, $R_3:h_2^2(h_3^{-1}h_2^{-1})^2h_3h_2^{-1}h_3=1$,\\$R_4:h_2^2h_3^{-1}h_2^{-2}h_3=1$:\\
By interchanging $h_2$ and $h_3$ in (2) and with the same discussion, there is a contradiction.

\item[(20)]$ R_1:h_2h_3^{-2}h_2^{-1}h_3^2=1$, $R_2:h_2h_3^{-2}h_2^{-1}h_3^2=1$, $R_3:(h_2h_3)^2h_2^{-1}(h_2^{-1}h_3)^2=1$,\\$R_4:h_2(h_3h_2^{-1}h_3)^2=1$:\\
By interchanging $h_2$ and $h_3$ in (7) and with the same discussion, there is a contradiction.

\item[(21)]$ R_1:h_2h_3^{-2}h_2^{-1}h_3^2=1$, $R_2:h_2h_3^{-2}h_2^{-1}h_3^2=1$, $R_3:h_2h_3h_2^{-2}h_3h_2h_3^{-1}h_2h_3=1$,\\$R_4:(h_2h_3^{-1})^2(h_2^{-1}h_3)^2=1$:\\
By interchanging $h_2$ and $h_3$ in (6) and with the same discussion, there is a contradiction.

\item[(22)]$ R_1:h_2h_3^{-3}h_2h_3=1$, $R_2:h_2h_3^{-3}h_2h_3=1$, $R_3:h_2^2h_3h_2^{-2}h_3^{-2}h_2=1$,\\$R_4:h_2h_3h_2^{-1}h_3^{-1}h_2^{-1}h_3=1$:\\
By interchanging $h_2$ and $h_3$ in (8) and with the same discussion, there is a contradiction.

\item[(23)]$ R_1:h_2h_3^{-3}h_2h_3=1$, $R_2:h_2h_3^{-3}h_2h_3=1$, $R_3:h_2^2h_3h_2^{-2}h_3^{-2}h_2=1$,\\$R_4:h_2h_3h_2^{-1}h_3^{-2}h_2=1$:\\
By interchanging $h_2$ and $h_3$ in (9) and with the same discussion, there is a contradiction.

\item[(24)]$ R_1:h_2h_3^{-2}h_2h_3^2=1$, $R_2:h_2h_3^{-2}h_2h_3^2=1$, $R_3:h_2(h_2h_3^{-1})^2(h_2h_3)^2=1$,\\$R_4:h_2^2h_3^{-1}h_2^2h_3=1$:\\
By interchanging $h_2$ and $h_3$ in (1) and with the same discussion, there is a contradiction.

\end{enumerate}

\subsection{$\mathbf{C_5-C_5(--C_7--)(-C_5-)}$}
\begin{figure}[ht]
\psscalebox{0.9 0.9} 
{
\begin{pspicture}(0,-1.8985577)(5.16,1.8985577)
\psdots[linecolor=black, dotsize=0.4](2.4,1.3014424)
\psdots[linecolor=black, dotsize=0.4](2.4,0.5014423)
\psdots[linecolor=black, dotsize=0.4](3.6,1.7014424)
\psdots[linecolor=black, dotsize=0.4](4.4,0.90144235)
\psdots[linecolor=black, dotsize=0.4](3.6,0.10144234)
\psdots[linecolor=black, dotsize=0.4](1.2,0.10144234)
\psdots[linecolor=black, dotsize=0.4](0.4,0.90144235)
\psdots[linecolor=black, dotsize=0.4](1.2,1.7014424)
\psline[linecolor=black, linewidth=0.04](1.2,1.7014424)(2.4,1.3014424)(3.6,1.7014424)(4.4,0.90144235)(3.6,0.10144234)(2.4,0.5014423)(2.4,1.3014424)
\psline[linecolor=black, linewidth=0.04](2.4,0.5014423)(1.2,0.10144234)(0.4,0.90144235)(1.2,1.7014424)
\psdots[linecolor=black, dotsize=0.4](1.2,-0.6985577)
\psdots[linecolor=black, dotsize=0.4](3.6,-0.6985577)
\psline[linecolor=black, linewidth=0.04](0.4,0.90144235)(1.2,-0.6985577)(3.6,-0.6985577)(4.4,0.90144235)
\psdots[linecolor=black, dotsize=0.4](1.6,-0.29855767)
\psdots[linecolor=black, dotsize=0.4](3.2,-0.29855767)
\psline[linecolor=black, linewidth=0.04](1.2,0.10144234)(1.6,-0.29855767)(3.2,-0.29855767)(3.6,0.10144234)
\rput[bl](-0.3,-1.8985577){$\mathbf{35) \ C_5-C_5(--C_7--)(-C_5-)}$}
\end{pspicture}
}
\end{figure}
By considering the $19$ cases related to the existence of $C_5-C_5(--C_7--)$ in the graph $K(\alpha,\beta)$ and the relations from Table \ref{tab-C5} which are not disproved, it can be seen that there are $138$ cases for the relations of a cycle $C_7$ and three cycles $C_5$ in the graph $C_5-C_5(--C_7--)(-C_5-)$. Using GAP \cite{gap}, we see that all groups with two generators $h_2$ and $h_3$ and four relations which are between $132$ cases of these $138$ cases are finite or solvable, that is a contradiction with the assumptions. So, there are just $6$ cases for the relations of these cycles which may lead to the existence of a subgraph isomorphic to the graph $C_5-C_5(--C_7--)(-C_5-)$ in $K(\alpha,\beta)$. In the following, we show that these $6$ cases lead to contradictions and so, the graph $K(\alpha,\beta)$ contains no subgraph isomorphic to the graph $C_5-C_5(--C_7--)(-C_5-)$.
\begin{enumerate}
\item[(1)]$ R_1:h_2^2h_3h_2^{-2}h_3=1$, $R_2:h_2^2h_3h_2^{-2}h_3=1$, $R_3:(h_2h_3)^2(h_3h_2^{-1})^2h_3=1$,\\$R_4:h_2h_3^{-2}h_2^{-1}h_3^2=1$:\\
$\Rightarrow$  $G$ has a torsion element, a contradiction.

\item[(2)]$ R_1:h_2^2h_3^{-1}h_2^{-2}h_3=1$, $R_2:h_2^2h_3^{-1}h_2^{-2}h_3=1$, $R_3:h_2h_3h_2h_3^{-2}(h_2^{-1}h_3)^2=1$,\\$R_4:h_2h_3^2h_2^{-1}h_3^2=1$:\\
$\Rightarrow$ $G$ is solvable, a contradiction.

\item[(3)]$ R_1:h_2^2h_3^{-1}h_2^{-2}h_3=1$, $R_2:h_2^2h_3^{-1}h_2^{-2}h_3=1$, $R_3:h_2h_3h_2h_3^{-2}h_2h_3^2=1$,\\$R_4:h_2^2h_3^{-1}h_2h_3^2=1$:\\
$\Rightarrow \langle h_2,h_3\rangle=\langle h_2\rangle$ is abelian, a contradiction.

\item[(4)]$ R_1:h_2h_3^{-2}h_2^{-1}h_3^2=1$, $R_2:h_2h_3^{-2}h_2^{-1}h_3^2=1$, $R_3:h_2(h_2h_3)^2h_2^{-2}h_3=1$,\\$R_4:h_2^2h_3^2h_2^{-1}h_3=1$:\\
By interchanging $h_2$ and $h_3$ in (3) and with the same discussion, there is a contradiction.

\item[(5)]$ R_1:h_2h_3^{-2}h_2^{-1}h_3^2=1$, $R_2:h_2h_3^{-2}h_2^{-1}h_3^2=1$, $R_3:h_2^2(h_3^{-1}h_2^{-1})^2h_3h_2^{-1}h_3=1$,\\$R_4:h_2^2h_3^{-1}h_2^2h_3=1$:\\
By interchanging $h_2$ and $h_3$ in (2) and with the same discussion, there is a contradiction.

\item[(6)]$ R_1:h_2h_3^{-2}h_2h_3^2=1$, $R_2:h_2h_3^{-2}h_2h_3^2=1$, $R_3:h_2(h_2h_3^{-1})^2(h_2h_3)^2=1$,\\$R_4:h_2^2h_3^{-1}h_2^{-2}h_3=1$:\\
By interchanging $h_2$ and $h_3$ in (1) and with the same discussion, there is a contradiction.

\end{enumerate}

\subsection{$\mathbf{C_5-C_5(-C_6--)(--C_6-2)}$}
\begin{figure}[ht]
\psscalebox{0.9 0.9} 
{
\begin{pspicture}(0,-1.8985577)(5.5,1.8985577)
\psdots[linecolor=black, dotsize=0.4](2.8,1.3014424)
\psdots[linecolor=black, dotsize=0.4](2.0,1.7014424)
\psdots[linecolor=black, dotsize=0.4](1.2,0.90144235)
\psdots[linecolor=black, dotsize=0.4](2.0,0.10144234)
\psdots[linecolor=black, dotsize=0.4](2.8,0.5014423)
\psdots[linecolor=black, dotsize=0.4](3.6,1.7014424)
\psdots[linecolor=black, dotsize=0.4](4.4,0.90144235)
\psdots[linecolor=black, dotsize=0.4](3.6,0.10144234)
\psdots[linecolor=black, dotsize=0.4](2.4,-0.6985577)
\psdots[linecolor=black, dotsize=0.4](4.4,-0.6985577)
\psline[linecolor=black, linewidth=0.04](2.8,1.3014424)(2.8,0.5014423)(2.0,0.10144234)(1.2,0.90144235)(2.0,1.7014424)(2.8,1.3014424)(3.6,1.7014424)(4.4,0.90144235)(3.6,0.10144234)(2.8,0.5014423)
\psline[linecolor=black, linewidth=0.04](2.0,0.10144234)(2.4,-0.6985577)(4.4,-0.6985577)(4.4,0.90144235)
\psdots[linecolor=black, dotsize=0.4](1.2,-1.0985577)
\psdots[linecolor=black, dotsize=0.4](3.2,-1.0985577)
\psline[linecolor=black, linewidth=0.04](3.6,0.10144234)(3.2,-1.0985577)(1.2,-1.0985577)(1.2,0.90144235)
\rput[bl](-0.5,-1.8985577){$\mathbf{36) \ C_5-C_5(-C_6--)(--C_6-2)}$}
\end{pspicture}
}
\end{figure}
By considering the relations from Tables \ref{tab-C5-C5(-C6--)} and \ref{tab-C6} which are not disproved, it can be seen that there are $22$ cases for the relations of two cycles $C_5$ and two cycles $C_6$ in the graph $C_5-C_5(-C_6--)(--C_6-2)$. By considering all groups with two generators $h_2$ and $h_3$ and four relations which are between these cases and by using GAP \cite{gap}, we see that all of these groups are finite and solvable. So, the graph $K(\alpha,\beta)$ contains no subgraph isomorphic to the graph $C_5-C_5(-C_6--)(--C_6-2)$.

\begin{figure}[ht]
\psscalebox{0.9 0.9} 
{
\begin{pspicture}(0,-1.8985577)(2.8471155,1.8985577)
\psdots[linecolor=black, dotsize=0.4](0.59711546,1.7014424)
\psdots[linecolor=black, dotsize=0.4](0.59711546,0.90144235)
\psdots[linecolor=black, dotsize=0.4](1.3971155,0.90144235)
\psdots[linecolor=black, dotsize=0.4](1.3971155,1.7014424)
\psdots[linecolor=black, dotsize=0.4](2.1971154,1.7014424)
\psdots[linecolor=black, dotsize=0.4](2.1971154,0.90144235)
\psdots[linecolor=black, dotsize=0.4](0.19711548,0.10144234)
\psdots[linecolor=black, dotsize=0.4](0.59711546,-0.6985577)
\psdots[linecolor=black, dotsize=0.4](2.1971154,-0.6985577)
\psdots[linecolor=black, dotsize=0.4](2.5971155,0.10144234)
\psline[linecolor=black, linewidth=0.04](1.3971155,1.7014424)(0.59711546,1.7014424)(0.59711546,0.90144235)(0.19711548,0.10144234)(0.59711546,-0.6985577)(2.1971154,-0.6985577)(2.5971155,0.10144234)(2.1971154,0.90144235)(2.1971154,1.7014424)(1.3971155,1.7014424)(1.3971155,0.90144235)(0.59711546,0.90144235)
\psline[linecolor=black, linewidth=0.04](1.3971155,0.90144235)(2.1971154,0.90144235)
\rput[bl](0.19711548,-1.8985577){$\mathbf{C_4-C_4(-C_7-)}$}
\end{pspicture}
}
\end{figure}
$\mathbf{C_4-C_4(-C_7-)}$ \textbf{subgraph:} 
By considering the possible cases of $C_4-C_4$ from Remark \ref{r-1} and the relations of $C_7$ cycles, it can be seen that there are $258$ cases for the relations of two cycles $C_4$ and a cycle $C_7$ in this structure. Using GAP \cite{gap}, we see that all groups with two generators $h_2$ and $h_3$ and three relations which are between $236$ cases of these $258$ cases are finite or solvable, that is a contradiction with the assumptions. So, there are just $22$ cases for the relations of these cycles which may lead to the existence of a subgraph isomorphic to the graph $C_4-C_4(-C_7-)$ in $K(\alpha,\beta)$. Similar to the previous mentioned subgraphs it can be seen that $16$ cases of these relations lead to contradictions and  $6$ cases of them may lead to the  existence of a subgraph isomorphic to the graph $C_4-C_4(-C_7-)$ in the graph $K(\alpha,\beta)$.

\subsection{$\mathbf{C_4-C_4(-C_7-)(C_4)}$}
\begin{figure}[ht]
\psscalebox{0.9 0.9} 
{
\begin{pspicture}(0,-1.8985577)(3.7971156,1.8985577)
\psdots[linecolor=black, dotsize=0.4](0.8,0.90144235)
\psdots[linecolor=black, dotsize=0.4](1.6,0.90144235)
\psdots[linecolor=black, dotsize=0.4](2.4,0.90144235)
\psdots[linecolor=black, dotsize=0.4](2.8,0.10144234)
\psdots[linecolor=black, dotsize=0.4](2.4,-0.6985577)
\psdots[linecolor=black, dotsize=0.4](0.8,-0.6985577)
\psdots[linecolor=black, dotsize=0.4](0.4,0.10144234)
\psdots[linecolor=black, dotsize=0.4](3.2,-1.0985577)
\psdots[linecolor=black, dotsize=0.4](3.6,-0.29855767)
\psdots[linecolor=black, dotsize=0.4](0.8,1.7014424)
\psdots[linecolor=black, dotsize=0.4](1.6,1.7014424)
\psdots[linecolor=black, dotsize=0.4](2.4,1.7014424)
\psline[linecolor=black, linewidth=0.04](0.8,1.7014424)(0.8,0.90144235)(1.6,0.90144235)(1.6,1.7014424)(0.8,1.7014424)
\psline[linecolor=black, linewidth=0.04](1.6,1.7014424)(2.4,1.7014424)(2.4,0.90144235)(1.6,0.90144235)
\psline[linecolor=black, linewidth=0.04](0.8,0.90144235)(0.4,0.10144234)(0.8,-0.6985577)(2.4,-0.6985577)(3.2,-1.0985577)(3.6,-0.29855767)(2.8,0.10144234)(2.4,0.90144235)
\psline[linecolor=black, linewidth=0.04](2.8,0.10144234)(2.4,-0.6985577)
\rput[bl](-0.3,-1.8985577){$\mathbf{37) \ C_4-C_4(-C_7-)(C_4)}$}
\end{pspicture}
}
\end{figure}
By considering the $6$ cases related to the existence of $C_4-C_4(-C_7-)$ in the graph $K(\alpha,\beta)$ and the relations from Table \ref{tab-C4} which are not disproved, it can be seen that there are $32$ cases for the relations of a cycle $C_7$ and three cycles $C_4$ in the graph $C_4-C_4(-C_7-)(C_4)$. By considering all groups with two generators $h_2$ and $h_3$ and four relations which are between these cases and by using GAP \cite{gap}, we see that all of these groups are finite and solvable. So, the graph $K(\alpha,\beta)$ contains no subgraph isomorphic to the graph $C_4-C_4(-C_7-)(C_4)$.

\subsection{$\mathbf{C_5--C_5(-C_5--)}$}
\begin{figure}[ht]
\psscalebox{0.9 0.9} 
{
\begin{pspicture}(0,-1.8985577)(3.78,1.8985577)
\psdots[linecolor=black, dotsize=0.4](1.6,1.7014424)
\psdots[linecolor=black, dotsize=0.4](1.6,0.5014423)
\psdots[linecolor=black, dotsize=0.4](1.6,-0.6985577)
\psdots[linecolor=black, dotsize=0.4](2.8,0.10144234)
\psdots[linecolor=black, dotsize=0.4](2.8,0.90144235)
\psdots[linecolor=black, dotsize=0.4](0.4,0.90144235)
\psdots[linecolor=black, dotsize=0.4](0.4,0.10144234)
\psline[linecolor=black, linewidth=0.04](1.6,1.7014424)(0.4,0.90144235)(0.4,0.10144234)(1.6,-0.6985577)(1.6,0.5014423)(1.6,1.7014424)(2.8,0.90144235)(2.8,0.10144234)(1.6,-0.6985577)
\psdots[linecolor=black, dotsize=0.4](2.8,-1.0985577)
\psbezier[linecolor=black, linewidth=0.04](2.8,0.90144235)(3.6944273,0.45422873)(3.7486832,-0.7823299)(2.8,-1.0985577)(1.8513167,-1.4147854)(0.4,-0.89855766)(0.4,0.10144234)
\rput[bl](-0.3,-1.8985577){$\mathbf{38) \ C_5--C_5(-C_5--)}$}
\end{pspicture}
}
\end{figure}
By considering the relations from Tables \ref{tab-C5--C5} and \ref{tab-C5} which are not disproved, it can be seen that there are $64$ cases for the relations of three cycles $C_5$ in the graph $C_5--C_5(-C_5--)$. Using GAP \cite{gap}, we see that all groups with two generators $h_2$ and $h_3$ and three relations which are between $58$ cases of these $64$ cases are finite and solvable, that is a contradiction with the assumptions. So, there are just $6$ cases for the relations of these cycles which may lead to the existence of a subgraph isomorphic to the graph $C_5--C_5(-C_5--)$ in $K(\alpha,\beta)$. In the following, we show that these $6$ cases lead to contradictions and so, the graph $K(\alpha,\beta)$ contains no subgraph isomorphic to the graph $C_5--C_5(-C_5--)$.
\begin{enumerate}
\item[(1)]$ R_1:h_2^2h_3^{-1}h_2^{-2}h_3=1$, $R_2:(h_2h_3)^2h_2^{-1}h_3=1$, $R_3:h_2^2h_3^2h_2^{-1}h_3=1$:\\
$\Rightarrow \langle h_2,h_3\rangle=\langle h_3\rangle$ is abelian, a contradiction.

\item[(2)]$ R_1:h_2^2h_3^{-1}h_2^{-1}h_3^{-1}h_2=1$, $R_2:h_2h_3^{-3}h_2h_3=1$, $R_3:h_2h_3^{-1}(h_3^{-1}h_2)^3=1$:\\
$\Rightarrow\langle h_2,h_3\rangle\cong BS(-3,1)$ is solvable, a contradiction.

\item[(3)]$ R_1:h_2^2h_3^{-2}h_2^{-1}h_3=1$, $R_2:h_2h_3h_2^{-1}h_3^{-2}h_2=1$, $R_3:h_2h_3^{-1}(h_3^{-1}h_2)^3=1$:\\
$\Rightarrow\langle h_2,h_3\rangle\cong BS(-3,1)$ is solvable, a contradiction.

\item[(4)]$ R_1:h_2h_3h_2h_3^{-1}h_2h_3=1$, $R_2:h_2h_3^{-2}h_2^{-1}h_3^2=1$, $R_3:h_2^2h_3^{-1}h_2h_3^2=1$:\\
By interchanging $h_2$ and $h_3$ in (1) and with the same discussion, there is a contradiction.

\item[(5)]$ R_1:h_2(h_3h_2^{-1}h_3)^2=1$, $R_2:(h_2h_3^{-1})^2(h_2^{-1}h_3)^2=1$, $R_3:h_2^2h_3^{-1}h_2^{-2}h_3=1$:\\
$\Rightarrow \langle h_2,h_3\rangle=\langle h_3\rangle$ is abelian, a contradiction.

\item[(6)]$ R_1:(h_2h_3^{-1}h_2)^2h_3=1$, $R_2:(h_2h_3^{-1})^2(h_2^{-1}h_3)^2=1$, $R_3:h_2h_3^{-2}h_2^{-1}h_3^2=1$:\\
By interchanging $h_2$ and $h_3$ in (5) and with the same discussion, there is a contradiction.

\end{enumerate}

\subsection{$\mathbf{C_6---C_6(-C_4)}$}
\begin{figure}[ht]
\psscalebox{0.9 0.9} 
{
\begin{pspicture}(0,-1.8985577)(3.7071154,1.8985577)
\psdots[linecolor=black, dotsize=0.4](2.1971154,1.7014424)
\psdots[linecolor=black, dotsize=0.4](2.1971154,0.90144235)
\psdots[linecolor=black, dotsize=0.4](2.1971154,0.10144234)
\psdots[linecolor=black, dotsize=0.4](2.1971154,-0.6985577)
\psdots[linecolor=black, dotsize=0.4](3.3971152,0.10144234)
\psdots[linecolor=black, dotsize=0.4](3.3971152,0.90144235)
\psdots[linecolor=black, dotsize=0.4](0.9971153,0.90144235)
\psdots[linecolor=black, dotsize=0.4](0.9971153,0.10144234)
\psdots[linecolor=black, dotsize=0.4](0.19711533,0.10144234)
\psdots[linecolor=black, dotsize=0.4](0.19711533,0.90144235)
\psline[linecolor=black, linewidth=0.04](2.1971154,1.7014424)(2.1971154,0.90144235)(2.1971154,0.10144234)(2.1971154,-0.6985577)(3.3971152,0.10144234)(3.3971152,0.90144235)(2.1971154,1.7014424)(0.9971153,0.90144235)(0.9971153,0.10144234)(2.1971154,-0.6985577)
\psline[linecolor=black, linewidth=0.04](0.9971153,0.10144234)(0.19711533,0.10144234)(0.19711533,0.90144235)(0.9971153,0.90144235)
\rput[bl](-0.1,-1.8985577){$\mathbf{39) \ C_6---C_6(-C_4)}$}
\end{pspicture}
}
\end{figure}
By considering the $99$ cases related to the existence of $C_6---C_6$ in the graph $K(\alpha,\beta)$ and the relations from Table \ref{tab-C4} which are not disproved, it can be seen that there are $420$ cases for the relations of two cycles $C_6$ and a cycle $C_4$ in the graph $C_6---C_6(-C_4)$. Using GAP \cite{gap}, we see that all groups with two generators $h_2$ and $h_3$ and three relations which are between $398$ cases of these $420$ cases are finite or solvable, that is a contradiction with the assumptions. So, there are just $22$ cases for the relations of these cycles which may lead to the existence of a subgraph isomorphic to the graph $C_6---C_6(-C_4)$ in $K(\alpha,\beta)$. In the following, we show that these $22$ cases lead to contradictions and so, the graph $K(\alpha,\beta)$ contains no subgraph isomorphic to the graph $C_6---C_6(-C_4)$.
\begin{enumerate}
\item[(1)]$ R_1:h_2^2h_3h_2^{-1}h_3^{-2}h_2=1$, $R_2:h_2^2h_3^{-1}h_2^{-1}h_3^{-1}h_2h_3=1$, $R_3:h_2h_3^{-2}h_2h_3=1$:\\
$\Rightarrow \langle h_2,h_3\rangle=\langle h_2\rangle$ is abelian, a contradiction.

\item[(2)]$ R_1:h_2^2h_3h_2^{-1}h_3^{-2}h_2=1$, $R_2:h_2h_3h_2h_3^{-1}(h_3^{-1}h_2)^2=1$, $R_3:h_2h_3^{-2}h_2h_3=1$:\\
$\Rightarrow \langle h_2,h_3\rangle=\langle h_2\rangle$ is abelian, a contradiction.

\item[(3)]$ R_1:h_2^2h_3h_2^{-1}h_3^{-2}h_2=1$, $R_2:h_2h_3^{-2}h_2^2h_3^{-1}h_2h_3=1$, $R_3:h_2h_3^{-2}h_2h_3=1$:\\
$\Rightarrow \langle h_2,h_3\rangle=\langle h_2\rangle$ is abelian, a contradiction.

\item[(4)]$ R_1:h_2^2h_3^{-1}h_2h_3h_2^{-1}h_3^2=1$, $R_2:h_2h_3h_2^{-1}h_3h_2h_3^{-1}h_2h_3=1$, $R_3:h_2h_3^{-1}h_2h_3h_2^{-1}h_3=1$:\\
$\Rightarrow \langle h_2,h_3\rangle=\langle h_2\rangle$ is abelian, a contradiction.

\item[(5)]$ R_1:h_2(h_2h_3^{-1})^2h_3^{-1}h_2h_3=1$, $R_2:h_2(h_3h_2^{-1})^2h_3^{-2}h_2=1$, $R_3:h_2h_3^{-2}h_2h_3=1$:\\
$\Rightarrow \langle h_2,h_3\rangle=\langle h_2\rangle$ is abelian, a contradiction.

\item[(6)]$ R_1:h_2h_3h_2h_3^{-1}(h_3^{-1}h_2)^2=1$, $R_2:h_2h_3^{-2}h_2h_3h_2h_3^{-1}h_2=1$, $R_3:h_2h_3^{-2}h_2h_3=1$:\\
$\Rightarrow \langle h_2,h_3\rangle=\langle h_2\rangle$ is abelian, a contradiction.

\item[(7)]$ R_1:h_2h_3h_2h_3^{-1}h_2h_3h_2^{-1}h_3=1$, $R_2:h_2h_3^2h_2^{-1}h_3h_2h_3^{-1}h_2=1$, $R_3:h_2h_3^{-1}h_2h_3h_2^{-1}h_3=1$:\\
By interchanging $h_2$ and $h_3$ in (4) and with the same discussion, there is a contradiction.

\item[(8)]$ R_1:h_2h_3h_2h_3^{-1}h_2h_3h_2^{-1}h_3=1$, $R_2:h_2(h_3h_2^{-1})^2h_3^{-2}h_2=1$, $R_3:h_2h_3^{-1}h_2h_3h_2^{-1}h_3=1$:\\
$\Rightarrow \langle h_2,h_3\rangle=\langle h_2\rangle$ is abelian, a contradiction.

\item[(9)]$ R_1:h_2h_3^2h_2^{-1}h_3^{-1}h_2^{-1}h_3=1$, $R_2:h_2h_3h_2^{-1}h_3^{-3}h_2=1$, $R_3:h_2h_3h_2^{-2}h_3=1$:\\
By interchanging $h_2$ and $h_3$ in (1) and with the same discussion, there is a contradiction.

\item[(10)]$ R_1:h_2h_3(h_3h_2^{-1})^2h_2^{-1}h_3=1$, $R_2:h_2(h_3h_2^{-1})^2h_3^{-2}h_2=1$, $R_3:h_2h_3h_2^{-2}h_3=1$:\\
By interchanging $h_2$ and $h_3$ in (5) and with the same discussion, there is a contradiction.

\item[(11)]$ R_1:h_2h_3h_2^{-2}h_3^2h_2^{-1}h_3=1$, $R_2:h_2h_3h_2^{-1}h_3^{-3}h_2=1$, $R_3:h_2h_3h_2^{-2}h_3=1$:\\
By interchanging $h_2$ and $h_3$ in (3) and with the same discussion, there is a contradiction.

\item[(12)]$ R_1:h_2h_3h_2^{-1}(h_2^{-1}h_3)^2h_3=1$, $R_2:h_2h_3h_2^{-1}h_3^{-3}h_2=1$, $R_3:h_2h_3h_2^{-2}h_3=1$:\\
By interchanging $h_2$ and $h_3$ in (2) and with the same discussion, there is a contradiction.

\item[(13)]$ R_1:h_2h_3h_2^{-1}(h_2^{-1}h_3)^2h_3=1$, $R_2:h_2h_3h_2^{-1}h_3^2h_2^{-2}h_3=1$, $R_3:h_2h_3h_2^{-2}h_3=1$:\\
By interchanging $h_2$ and $h_3$ in (6) and with the same discussion, there is a contradiction.

\item[(14)]$ R_1:h_2h_3h_2^{-1}h_3^{-1}h_2h_3^{-1}h_2^{-1}h_3=1$, $R_2:h_2h_3^{-1}h_2^{-1}h_3^2h_2^{-1}h_3^{-1}h_2=1$, $R_3:h_2(h_3h_2^{-1})^2h_3=1$:\\
$\Rightarrow \langle h_2,h_3\rangle=\langle h_3\rangle$ is abelian, a contradiction.

\item[(15)]$ R_1:h_2h_3h_2^{-1}h_3^{-1}h_2h_3^{-1}h_2^{-1}h_3=1$, $R_2:h_2h_3^{-1}h_2^{-1}h_3^2h_2^{-1}h_3^{-1}h_2=1$, $R_3:h_2h_3^{-3}h_2=1$:\\
$\Rightarrow \langle h_2,h_3\rangle=\langle h_2\rangle$ is abelian, a contradiction.

\item[(16)]$ R_1:h_2h_3h_2^{-1}h_3^{-1}h_2h_3^{-1}h_2^{-1}h_3=1$, $R_2:h_2h_3^{-1}h_2^{-1}h_3^2h_2^{-1}h_3^{-1}h_2=1$, $R_3:h_2h_3^{-1}h_2h_3^2=1$:\\
$\Rightarrow \langle h_2,h_3\rangle=\langle h_2\rangle$ is abelian, a contradiction.

\item[(17)]$ R_1:h_2h_3h_2^{-1}h_3^{-1}h_2h_3^{-1}h_2^{-1}h_3=1$, $R_2:h_2h_3^{-1}h_2^{-1}h_3^2h_2^{-1}h_3^{-1}h_2=1$, $R_3:(h_2h_3^{-1})^2h_2h_3=1$:\\
$\Rightarrow \langle h_2,h_3\rangle=\langle h_2\rangle$ is abelian, a contradiction.

\item[(18)]$ R_1:h_2h_3h_2^{-1}h_3h_2h_3^{-1}h_2h_3=1$, $R_2:h_2(h_3h_2^{-1})^2h_3^{-2}h_2=1$, $R_3:h_2h_3^{-1}h_2h_3h_2^{-1}h_3=1$:\\
By interchanging $h_2$ and $h_3$ in (8) and with the same discussion, there is a contradiction.

\item[(19)]$ R_1:h_2h_3^{-1}h_2^{-1}h_3^2h_2^{-1}h_3^{-1}h_2=1$, $R_2:h_2h_3^{-1}h_2^{-1}h_3h_2^{-1}h_3^{-1}h_2h_3=1$, $R_3:h_2^2h_3h_2^{-1}h_3=1$:\\
By interchanging $h_2$ and $h_3$ in (16) and with the same discussion, there is a contradiction.

\item[(20)]$ R_1:h_2h_3^{-1}h_2^{-1}h_3^2h_2^{-1}h_3^{-1}h_2=1$, $R_2:h_2h_3^{-1}h_2^{-1}h_3h_2^{-1}h_3^{-1}h_2h_3=1$, $R_3:h_2^2h_3^{-2}h_2=1$:\\
By interchanging $h_2$ and $h_3$ in (15) and with the same discussion, there is a contradiction.

\item[(21)]$ R_1:h_2h_3^{-1}h_2^{-1}h_3^2h_2^{-1}h_3^{-1}h_2=1$, $R_2:h_2h_3^{-1}h_2^{-1}h_3h_2^{-1}h_3^{-1}h_2h_3=1$, $R_3:h_2(h_3h_2^{-1})^2h_3=1$:\\
By interchanging $h_2$ and $h_3$ in (17) and with the same discussion, there is a contradiction.

\item[(22)]$ R_1:h_2h_3^{-1}h_2^{-1}h_3^2h_2^{-1}h_3^{-1}h_2=1$, $R_2:h_2h_3^{-1}h_2^{-1}h_3h_2^{-1}h_3^{-1}h_2h_3=1$, $R_3:(h_2h_3^{-1})^2h_2h_3=1$:\\
By interchanging $h_2$ and $h_3$ in (14) and with the same discussion, there is a contradiction.

\end{enumerate}

\subsection{$\mathbf{C_6--C_6(C_4)}$}
\begin{figure}[ht]
\psscalebox{0.9 0.9} 
{
\begin{pspicture}(0,-1.8985577)(2.9971154,1.8985577)
\psdots[linecolor=black, dotsize=0.4](1.3971153,1.7014424)
\psdots[linecolor=black, dotsize=0.4](1.3971153,0.5014423)
\psdots[linecolor=black, dotsize=0.4](1.3971153,-0.6985577)
\psdots[linecolor=black, dotsize=0.4](2.5971153,-0.29855767)
\psdots[linecolor=black, dotsize=0.4](2.5971153,0.5014423)
\psdots[linecolor=black, dotsize=0.4](2.5971153,1.3014424)
\psdots[linecolor=black, dotsize=0.4](0.19711533,1.3014424)
\psdots[linecolor=black, dotsize=0.4](0.19711533,0.5014423)
\psdots[linecolor=black, dotsize=0.4](0.19711533,-0.29855767)
\psline[linecolor=black, linewidth=0.04](1.3971153,1.7014424)(1.3971153,0.5014423)(1.3971153,-0.6985577)(0.19711533,-0.29855767)(0.19711533,0.5014423)(0.19711533,1.3014424)(1.3971153,1.7014424)(2.5971153,1.3014424)(2.5971153,0.5014423)(2.5971153,-0.29855767)(1.3971153,-0.6985577)
\psline[linecolor=black, linewidth=0.04](1.3971153,0.5014423)(0.19711533,0.5014423)
\rput[bl](-0.3,-1.8985577){$\mathbf{40) \ C_6--C_6(C_4)}$}
\end{pspicture}
}
\end{figure}
By considering the $996$ cases related to the existence of $C_6--C_6$ in the graph $K(\alpha,\beta)$ and the relations from Table \ref{tab-C4} which are not disproved, it can be seen that there are $279$ cases for the relations of two cycles $C_6$ and a cycle $C_4$ in the graph $C_6--C_6(C_4)$. Using GAP \cite{gap}, we see that all groups with two generators $h_2$ and $h_3$ and three relations which are between $268$ cases of these $279$ cases are finite or solvable, that is a contradiction with the assumptions. So, there are just $11$ cases for the relations of these cycles which may lead to the existence of a subgraph isomorphic to the graph $C_6--C_6(C_4)$ in $K(\alpha,\beta)$. In the following, we show that these $11$ cases lead to contradictions and so, the graph $K(\alpha,\beta)$ contains no subgraph isomorphic to the graph $C_6--C_6(C_4)$.
\begin{enumerate}
\item[(1)]$ R_1:h_2^3h_3^3=1$, $R_2:h_2^2h_3^{-1}h_2h_3^2h_2^{-1}h_3=1$, $R_3:h_2h_3^{-1}h_2h_3h_2^{-1}h_3=1$:\\
$\Rightarrow \langle h_2,h_3\rangle=\langle h_2\rangle$ is abelian, a contradiction.

\item[(2)]$ R_1:h_2^2(h_2h_3^{-1})^3h_2=1$, $R_2:h_2^2h_3^{-1}(h_3^{-1}h_2)^2h_3=1$, $R_3:h_2h_3^{-2}h_2h_3=1$:\\
$\Rightarrow \langle h_2,h_3\rangle=\langle h_2\rangle$ is abelian, a contradiction.

\item[(3)]$ R_1:h_2^2h_3h_2^{-1}h_3^{-2}h_2=1$, $R_2:h_2^2h_3^{-1}(h_3^{-1}h_2)^2h_3=1$, $R_3:h_2h_3^{-2}h_2h_3=1$:\\
$\Rightarrow\langle h_2,h_3\rangle\cong BS(1,-1)$ is solvable, a contradiction.

\item[(4)]$ R_1:h_2^2h_3^{-1}h_2^{-1}h_3^{-1}h_2h_3=1$, $R_2:(h_2h_3^{-1})^3h_2^{-1}h_3^2=1$, $R_3:h_2h_3^{-2}h_2h_3=1$:\\
$\Rightarrow\langle h_2,h_3\rangle\cong BS(1,-1)$ is solvable, a contradiction.

\item[(5)]$ R_1:h_2^2h_3^{-1}(h_3^{-1}h_2)^2h_3=1$, $R_2:(h_2h_3^{-1})^3h_2^{-1}h_3^2=1$, $R_3:h_2h_3^{-2}h_2h_3=1$:\\
$\Rightarrow\langle h_2,h_3\rangle\cong BS(1,-1)$ is solvable, a contradiction.

\item[(6)]$ R_1:h_2(h_2h_3^{-1})^2h_3^{-1}h_2h_3=1$, $R_2:(h_2h_3^{-1})^3h_2^{-1}h_3^2=1$, $R_3:h_2h_3^{-2}h_2h_3=1$:\\
$\Rightarrow\langle h_2,h_3\rangle\cong BS(1,-1)$ is solvable, a contradiction.

\item[(7)]$ R_1:h_2h_3^2h_2^{-1}(h_2^{-1}h_3)^2=1$, $R_2:h_2h_3h_2^{-1}h_3^{-3}h_2=1$, $R_3:h_2h_3h_2^{-2}h_3=1$:\\
By interchanging $h_2$ and $h_3$ in (3) and with the same discussion, there is a contradiction.

\item[(8)]$ R_1:h_2h_3^2h_2^{-1}(h_2^{-1}h_3)^2=1$, $R_2:h_2(h_3h_2^{-1})^3h_3^{-1}h_2=1$, $R_3:h_2h_3h_2^{-2}h_3=1$:\\
By interchanging $h_2$ and $h_3$ in (5) and with the same discussion, there is a contradiction.

\item[(9)]$ R_1:h_2h_3^2h_2^{-1}(h_2^{-1}h_3)^2=1$, $R_2:h_3^2(h_3h_2^{-1})^3h_3=1$, $R_3:h_2h_3h_2^{-2}h_3=1$:\\
By interchanging $h_2$ and $h_3$ in (2) and with the same discussion, there is a contradiction.

\item[(10)]$ R_1:h_2h_3^2h_2^{-1}h_3^{-1}h_2^{-1}h_3=1$, $R_2:h_2(h_3h_2^{-1})^3h_3^{-1}h_2=1$, $R_3:h_2h_3h_2^{-2}h_3=1$:\\
By interchanging $h_2$ and $h_3$ in (4) and with the same discussion, there is a contradiction.

\item[(11)]$ R_1:h_2h_3(h_3h_2^{-1})^2h_2^{-1}h_3=1$, $R_2:h_2(h_3h_2^{-1})^3h_3^{-1}h_2=1$, $R_3:h_2h_3h_2^{-2}h_3=1$:\\
By interchanging $h_2$ and $h_3$ in (6) and with the same discussion, there is a contradiction.

\end{enumerate}

\begin{figure}[ht]
\psscalebox{0.9 0.9} 
{
\begin{pspicture}(0,-1.8985577)(3.194231,1.8985577)
\psdots[linecolor=black, dotsize=0.4](1.3971155,1.7014424)
\psdots[linecolor=black, dotsize=0.4](0.59711546,1.3014424)
\psdots[linecolor=black, dotsize=0.4](0.59711546,0.5014423)
\psdots[linecolor=black, dotsize=0.4](1.3971155,0.10144234)
\psdots[linecolor=black, dotsize=0.4](2.1971154,0.5014423)
\psdots[linecolor=black, dotsize=0.4](2.1971154,1.3014424)
\psdots[linecolor=black, dotsize=0.4](2.9971154,1.3014424)
\psdots[linecolor=black, dotsize=0.4](2.9971154,0.5014423)
\psdots[linecolor=black, dotsize=0.4](0.19711548,-0.29855767)
\psdots[linecolor=black, dotsize=0.4](0.9971155,-0.6985577)
\psline[linecolor=black, linewidth=0.04](1.3971155,1.7014424)(0.59711546,1.3014424)(0.59711546,0.5014423)(0.19711548,-0.29855767)(0.9971155,-0.6985577)(1.3971155,0.10144234)(2.1971154,0.5014423)(2.9971154,0.5014423)(2.9971154,1.3014424)(2.1971154,1.3014424)(1.3971155,1.7014424)(1.3971155,1.7014424)
\psline[linecolor=black, linewidth=0.04](2.1971154,1.3014424)(2.1971154,0.5014423)
\psline[linecolor=black, linewidth=0.04](0.59711546,0.5014423)(1.3971155,0.10144234)
\rput[bl](0.59711546,-1.8985577){$\mathbf{C_4-C_6(-C_4)}$}
\end{pspicture}
}
\end{figure}
$\mathbf{C_4-C_6(-C_4)}$ \textbf{subgraph:} 
By considering the relations from Tables \ref{tab-C4-C6} and \ref{tab-C4} which are not disproved, it can be seen that there are $94$ cases for the relations of two cycles $C_4$ and a cycle $C_6$ in this structure. Using GAP \cite{gap}, we see that all groups with two generators $h_2$ and $h_3$ and three relations which are between $84$ cases of these $94$ cases are finite or solvable, that is a contradiction with the assumptions. So, there are just $10$ cases for the relations of these cycles which may lead to the existence of a subgraph isomorphic to the graph $C_4-C_6(-C_4)$ in $K(\alpha,\beta)$.

\subsection{$\mathbf{C_4-C_6(-C_4)(-C_4)}$}
\begin{figure}[ht]
\psscalebox{0.9 0.9}  
{
\begin{pspicture}(0,-2.0985577)(3.76,2.0985577)
\psdots[linecolor=black, dotsize=0.4](1.6,1.1014423)
\psdots[linecolor=black, dotsize=0.4](1.6,-0.4985577)
\psdots[linecolor=black, dotsize=0.4](2.4,-0.0985577)
\psdots[linecolor=black, dotsize=0.4](2.4,0.7014423)
\psdots[linecolor=black, dotsize=0.4](0.8,0.7014423)
\psdots[linecolor=black, dotsize=0.4](0.8,-0.0985577)
\psdots[linecolor=black, dotsize=0.4](0.4,1.5014423)
\psdots[linecolor=black, dotsize=0.4](1.2,1.9014423)
\psdots[linecolor=black, dotsize=0.4](3.2,0.7014423)
\psdots[linecolor=black, dotsize=0.4](3.2,-0.0985577)
\psdots[linecolor=black, dotsize=0.4](1.2,-1.2985578)
\psdots[linecolor=black, dotsize=0.4](0.4,-0.8985577)
\psline[linecolor=black, linewidth=0.04](0.8,0.7014423)(1.6,1.1014423)(2.4,0.7014423)(2.4,-0.0985577)(1.6,-0.4985577)(0.8,-0.0985577)(0.8,0.7014423)(0.4,1.5014423)(1.2,1.9014423)(1.6,1.1014423)
\psline[linecolor=black, linewidth=0.04](2.4,0.7014423)(3.2,0.7014423)(3.2,-0.0985577)(2.4,-0.0985577)
\psline[linecolor=black, linewidth=0.04](1.6,-0.4985577)(1.2,-1.2985578)(0.4,-0.8985577)(0.8,-0.0985577)
\rput[bl](-0.3,-2.0985577){$\mathbf{41) \ C_4-C_6(-C_4)(-C_4)}$}
\end{pspicture}
}
\end{figure}
By considering the $10$ cases related to the existence of $C_4-C_6(-C_4)$ in the graph $K(\alpha,\beta)$ and the relations from Table \ref{tab-C4} which are not disproved, it can be seen that there are $36$ cases for the relations of a cycle $C_6$ and three cycles $C_4$ in the graph $C_4-C_6(-C_4)(-C_4)$. By considering all groups with two generators $h_2$ and $h_3$ and four relations which are between these cases and by using GAP \cite{gap}, we see that all of these groups are solvable. So, the graph $K(\alpha,\beta)$ contains no subgraph isomorphic to the graph $C_4-C_6(-C_4)(-C_4)$.

\subsection{$\mathbf{C_4-C_6(--C_7--)(-C_5-)}$}
\begin{figure}[ht]
\psscalebox{0.9 0.9} 
{
\begin{pspicture}(0,-1.8985577)(5.16,1.8985577)
\psdots[linecolor=black, dotsize=0.4](2.0,1.3014424)
\psdots[linecolor=black, dotsize=0.4](2.0,0.5014423)
\psdots[linecolor=black, dotsize=0.4](2.8,0.5014423)
\psdots[linecolor=black, dotsize=0.4](2.8,1.3014424)
\psdots[linecolor=black, dotsize=0.4](3.6,1.7014424)
\psdots[linecolor=black, dotsize=0.4](4.4,1.3014424)
\psdots[linecolor=black, dotsize=0.4](4.4,0.5014423)
\psdots[linecolor=black, dotsize=0.4](3.6,0.10144234)
\psline[linecolor=black, linewidth=0.04](2.8,0.5014423)(3.6,0.10144234)(4.4,0.5014423)
\psline[linecolor=black, linewidth=0.04](2.0,1.3014424)(2.0,0.5014423)(2.8,0.5014423)(2.8,1.3014424)(3.6,1.7014424)(4.4,1.3014424)(4.4,0.5014423)
\psline[linecolor=black, linewidth=0.04](2.8,1.3014424)(2.0,1.3014424)
\psdots[linecolor=black, dotsize=0.4](1.6,0.10144234)
\psdots[linecolor=black, dotsize=0.4](1.2,-0.29855767)
\psdots[linecolor=black, dotsize=0.4](3.6,-0.6985577)
\psline[linecolor=black, linewidth=0.04](2.0,1.3014424)(0.4,-0.6985577)(3.6,-0.6985577)(4.4,0.5014423)
\psdots[linecolor=black, dotsize=0.4](0.4,-0.6985577)
\psline[linecolor=black, linewidth=0.04](0.4,-0.6985577)(1.2,-0.29855767)(1.6,0.10144234)(2.0,0.5014423)
\rput[bl](-0.3,-1.8985577){$\mathbf{42) \ C_4-C_6(--C_7--)(-C_5-)}$}
\end{pspicture}
}
\end{figure}
By considering the $16$ cases related to the existence of $C_4-C_6(--C_7--)$ in the graph $K(\alpha,\beta)$ from Table \ref{tab-C4-C6(--C7--)} and the relations from Table \ref{tab-C5} which are not disproved, it can be seen that there are $62$ cases for the relations of a cycle $C_4$, a cycle $C_6$, a cycle $C_7$ and a cycle $C_5$ in the graph $C_4-C_6(--C_7--)(-C_5-)$. Using GAP \cite{gap}, we see that all groups with two generators $h_2$ and $h_3$ and four relations which are between $58$ cases of these $62$ cases are solvable, that is a contradiction with the assumptions. So, there are just $4$ cases for the relations of these cycles which may lead to the existence of a subgraph isomorphic to the graph $C_4-C_6(--C_7--)(-C_5-)$ in $K(\alpha,\beta)$. In the following, we show that these $4$ cases lead to contradictions and so, the graph $K(\alpha,\beta)$ contains no subgraph isomorphic to the graph $C_4-C_6(--C_7--)(-C_5-)$.
\begin{enumerate}
\item[(1)]$ R_1:h_2h_3h_2^{-2}h_3=1$, $R_2:h_2h_3^2h_2^{-1}(h_2^{-1}h_3)^2=1$, $R_3:h_2h_3^{-4}h_2h_3h_2^{-1}h_3=1$,\\$R_4:h_2^2h_3^{-1}h_2^{-2}h_3=1$:\\
$\Rightarrow \langle h_2,h_3\rangle=\langle h_3\rangle$ is abelian, a contradiction.

\item[(2)]$ R_1:h_2h_3h_2^{-2}h_3=1$, $R_2:h_2(h_3h_2^{-1})^2h_2^{-1}h_3^2=1$, $R_3:(h_2h_3^{-2})^2h_2h_3^{-1}h_2^{-1}h_3=1$,\\$R_4:h_2(h_3h_2^{-1})^3h_3=1$:\\
$\Rightarrow \langle h_2,h_3\rangle=\langle h_3\rangle$ is abelian, a contradiction.

\item[(3)]$ R_1:h_2h_3^{-2}h_2h_3=1$, $R_2:h_2h_3(h_2h_3^{-1})^2h_3^{-1}h_2=1$, $R_3:h_2h_3^{-1}h_2^{-1}h_3(h_2h_3^{-1}h_2)^2=1$,\\$R_4:(h_2h_3^{-1})^3h_2h_3=1$:\\
By interchanging $h_2$ and $h_3$ in (2) and with the same discussion, there is a contradiction.

\item[(4)]$ R_1:h_2h_3^{-2}h_2h_3=1$, $R_2:h_2^2h_3^{-1}(h_3^{-1}h_2)^2h_3=1$, $R_3:h_2^3h_3^{-1}h_2^{-1}h_3h_2^{-1}h_3^{-1}h_2=1$,\\$R_4:h_2h_3^{-2}h_2^{-1}h_3^2=1$:\\
By interchanging $h_2$ and $h_3$ in (1) and with the same discussion, there is a contradiction.

\end{enumerate}

\subsection{$\mathbf{C_5-C_5(-C_6--)(--C_5-)}$}
\begin{figure}[ht]
\psscalebox{0.9 0.9} 
{
\begin{pspicture}(0,-2.0985577)(5.16,2.0985577)
\psdots[linecolor=black, dotsize=0.4](2.4,0.7014423)
\psdots[linecolor=black, dotsize=0.4](1.6,1.1014423)
\psdots[linecolor=black, dotsize=0.4](0.8,0.3014423)
\psdots[linecolor=black, dotsize=0.4](1.6,-0.4985577)
\psdots[linecolor=black, dotsize=0.4](2.4,-0.0985577)
\psdots[linecolor=black, dotsize=0.4](3.2,1.1014423)
\psdots[linecolor=black, dotsize=0.4](4.0,0.3014423)
\psdots[linecolor=black, dotsize=0.4](3.2,-0.4985577)
\psdots[linecolor=black, dotsize=0.4](1.6,1.9014423)
\psline[linecolor=black, linewidth=0.04](2.4,0.7014423)(1.6,1.1014423)(0.8,0.3014423)(1.6,-0.4985577)(2.4,-0.0985577)(2.4,0.7014423)(3.2,1.1014423)(4.0,0.3014423)(3.2,-0.4985577)(2.4,-0.0985577)
\psdots[linecolor=black, dotsize=0.4](2.0,-1.2985578)
\psdots[linecolor=black, dotsize=0.4](4.0,-1.2985578)
\psline[linecolor=black, linewidth=0.04](1.6,-0.4985577)(2.0,-1.2985578)(4.0,-1.2985578)(4.0,0.3014423)
\psline[linecolor=black, linewidth=0.04](0.8,0.3014423)(1.6,1.9014423)(3.2,1.1014423)
\rput[bl](-0.4,-2.0985577){$\mathbf{43) \ C_5-C_5(-C_6--)(--C_5-)}$}
\end{pspicture}
}
\end{figure}
By considering the relations from Tables \ref{tab-C5-C5(-C6--)} and \ref{tab-C5} which are not disproved, it can be seen that there are $14$ cases for the relations of three cycles $C_5$ and a cycle $C_6$ in the graph $C_5-C_5(-C_6--)(--C_5-)$. By considering all groups with two generators $h_2$ and $h_3$ and four relations which are between these cases and by using GAP \cite{gap}, we see that all of these groups are solvable. So, the graph $K(\alpha,\beta)$ contains no subgraph isomorphic to the graph $C_5-C_5(-C_6--)(--C_5-)$.

\subsection{$\mathbf{C_4-C_6(--C_7--)(C_7-2)}$}
\begin{figure}[ht]
\psscalebox{0.9 0.9}  
{
\begin{pspicture}(0,-1.8985577)(4.62,1.8985577)
\psdots[linecolor=black, dotsize=0.4](1.6,1.3014424)
\psdots[linecolor=black, dotsize=0.4](1.6,0.5014423)
\psdots[linecolor=black, dotsize=0.4](2.4,0.5014423)
\psdots[linecolor=black, dotsize=0.4](2.4,1.3014424)
\psdots[linecolor=black, dotsize=0.4](3.2,1.7014424)
\psdots[linecolor=black, dotsize=0.4](4.0,1.3014424)
\psdots[linecolor=black, dotsize=0.4](4.0,0.5014423)
\psdots[linecolor=black, dotsize=0.4](3.2,0.10144234)
\psline[linecolor=black, linewidth=0.04](2.4,0.5014423)(3.2,0.10144234)(4.0,0.5014423)
\psline[linecolor=black, linewidth=0.04](1.6,1.3014424)(1.6,0.5014423)(2.4,0.5014423)(2.4,1.3014424)(3.2,1.7014424)(4.0,1.3014424)(4.0,0.5014423)
\psline[linecolor=black, linewidth=0.04](2.4,1.3014424)(1.6,1.3014424)
\psdots[linecolor=black, dotsize=0.4](2.0,0.10144234)
\psdots[linecolor=black, dotsize=0.4](2.4,-0.29855767)
\psdots[linecolor=black, dotsize=0.4](3.6,-0.6985577)
\psdots[linecolor=black, dotsize=0.4](0.8,-0.6985577)
\psline[linecolor=black, linewidth=0.04](1.6,1.3014424)(0.8,-0.6985577)(3.6,-0.6985577)(4.0,0.5014423)
\psline[linecolor=black, linewidth=0.04](1.6,0.5014423)(2.0,0.10144234)(2.4,-0.29855767)(3.6,-0.6985577)
\rput[bl](-0.5,-1.8985577){$\mathbf{44) \ C_4-C_6(--C_7--)(C_7-2)}$}
\end{pspicture}
}
\end{figure}
By considering the $16$ cases related to the existence of $C_4-C_6(--C_7--)$ in the graph $K(\alpha,\beta)$ from Table \ref{tab-C4-C6(--C7--)} and the relations of $C_7$ cycles, it can be seen that there are $168$ cases for the relations of a cycle $C_4$, a cycle $C_6$ and two cycles $C_7$ in the graph $C_4-C_6(--C_7--)(C_7-2)$. Using GAP \cite{gap}, we see that all groups with two generators $h_2$ and $h_3$ and four relations which are between $152$ cases of these $168$ cases are solvable, that is a contradiction with the assumptions. So, there are just $16$ cases for the relations of these cycles which may lead to the existence of a subgraph isomorphic to the graph $C_4-C_6(--C_7--)(C_7-2)$ in $K(\alpha,\beta)$. In the following, we show that these $16$ cases lead to contradictions and so, the graph $K(\alpha,\beta)$ contains no subgraph isomorphic to the graph $C_4-C_6(--C_7--)(C_7-2)$.
\begin{enumerate}
\item[(1)]$ R_1:h_2h_3h_2^{-2}h_3=1$, $R_2:h_2h_3^2h_2^{-1}(h_2^{-1}h_3)^2=1$, $R_3:h_2h_3^{-4}h_2h_3h_2^{-1}h_3=1$,\\$R_4:h_2h_3^{-4}(h_3^{-1}h_2)^2=1$:\\
$\Rightarrow \langle h_2,h_3\rangle=\langle h_3\rangle$ is abelian, a contradiction.

\item[(2)]$ R_1:h_2h_3h_2^{-2}h_3=1$, $R_2:h_2h_3^2h_2^{-1}(h_2^{-1}h_3)^2=1$, $R_3:h_2h_3^{-2}h_2(h_3h_2^{-1})^3h_3=1$,\\$R_4:h_2h_3(h_3h_2^{-1})^4h_3=1$:\\
$\Rightarrow \langle h_2,h_3\rangle=\langle h_3\rangle$ is abelian, a contradiction.

\item[(3)]$ R_1:h_2h_3h_2^{-2}h_3=1$, $R_2:h_2(h_3h_2^{-1})^2h_2^{-1}h_3^2=1$, $R_3:(h_2h_3^{-2})^2h_2^{-1}h_3^2=1$,\\$R_4:(h_2h_3^{-3})^2h_2=1$:\\
$\Rightarrow \langle h_2,h_3\rangle=\langle h_3\rangle$ is abelian, a contradiction.

\item[(4)]$ R_1:h_2h_3h_2^{-2}h_3=1$, $R_2:h_2(h_3h_2^{-1})^2h_2^{-1}h_3^2=1$, $R_3:(h_2h_3^{-2})^2h_2h_3^{-1}h_2^{-1}h_3=1$,\\$R_4:h_2h_3^{-2}(h_2^{-1}h_3^2)^2=1$:\\
$\Rightarrow \langle h_2,h_3\rangle=\langle h_3\rangle$ is abelian, a contradiction.

\item[(5)]$ R_1:h_2h_3h_2^{-2}h_3=1$, $R_2:h_2(h_3h_2^{-1})^2h_2^{-1}h_3^2=1$, $R_3:(h_2h_3^{-2})^2h_2h_3^{-1}h_2^{-1}h_3=1$,\\$R_4:h_2h_3^{-2}(h_3^{-2}h_2)^2=1$:\\
$\Rightarrow\langle h_2,h_3\rangle\cong BS(1,-1)$ is solvable, a contradiction.

\item[(6)]$ R_1:h_2h_3h_2^{-2}h_3=1$, $R_2:h_2(h_3h_2^{-1})^2h_2^{-1}h_3^2=1$, $R_3:h_2h_3^{-1}(h_2^{-1}h_3^2)^2h_2^{-1}h_3=1$,\\$R_4:h_2(h_3h_2^{-1}h_3)^2h_2^{-1}h_3h_2^{-1}h_3=1$:\\
$\Rightarrow\langle h_2,h_3\rangle\cong BS(1,-1)$ is solvable, a contradiction.

\item[(7)]$ R_1:h_2h_3h_2^{-2}h_3=1$, $R_2:h_2(h_3h_2^{-1})^2h_2^{-1}h_3^2=1$, $R_3:h_2h_3^{-1}(h_2^{-1}h_3^2)^2h_2^{-1}h_3=1$,\\$R_4:h_2(h_3^{-1}h_2h_3^{-1})^2h_2^{-1}h_3h_2^{-1}h_3=1$:\\
$\Rightarrow\langle h_2,h_3\rangle\cong BS(1,-1)$ is solvable, a contradiction.

\item[(8)]$ R_1:h_2h_3h_2^{-2}h_3=1$, $R_2:h_2(h_3h_2^{-1})^2h_2^{-1}h_3^2=1$, $R_3:h_2h_3^{-1}h_2h_3^{-1}h_2^{-1}(h_3h_2^{-1}h_3)^2=1$,\\$R_4:h_2((h_3h_2^{-1})^2h_3)^2=1$:\\
$\Rightarrow \langle h_2,h_3\rangle=\langle h_3\rangle$ is abelian, a contradiction.

\item[(9)]$ R_1:h_2h_3^{-2}h_2h_3=1$, $R_2:h_2h_3(h_2h_3^{-1})^2h_3^{-1}h_2=1$, $R_3:(h_2^2h_3^{-1})^2h_2^{-2}h_3=1$,\\$R_4:h_2^2h_3^{-2}h_2^3h_3^{-1}h_2=1$:\\
By interchanging $h_2$ and $h_3$ in (3) and with the same discussion, there is a contradiction.

\item[(10)]$ R_1:h_2h_3^{-2}h_2h_3=1$, $R_2:h_2h_3(h_2h_3^{-1})^2h_3^{-1}h_2=1$, $R_3:h_2h_3^{-1}h_2^{-1}h_3(h_2h_3^{-1}h_2)^2=1$,\\$R_4:h_2^3h_3^{-2}h_2^2h_3^{-1}h_2=1$:\\
By interchanging $h_2$ and $h_3$ in (5) and with the same discussion, there is a contradiction.

\item[(11)]$ R_1:h_2h_3^{-2}h_2h_3=1$, $R_2:h_2h_3(h_2h_3^{-1})^2h_3^{-1}h_2=1$, $R_3:h_2h_3^{-1}h_2^{-1}h_3(h_2h_3^{-1}h_2)^2=1$,\\$R_4:h_2^2h_3^{-1}(h_2^{-2}h_3)^2=1$:\\
By interchanging $h_2$ and $h_3$ in (4) and with the same discussion, there is a contradiction.

\item[(12)]$ R_1:h_2h_3^{-2}h_2h_3=1$, $R_2:h_2h_3(h_2h_3^{-1})^2h_3^{-1}h_2=1$, $R_3:h_2h_3^{-1}(h_2^{-1}h_3h_2^{-1})^2h_2^{-1}h_3=1$,\\$R_4:(h_2h_3^{-1}h_2)^2h_3^{-1}h_2^{-1}h_3h_2^{-1}h_3=1$:\\
By interchanging $h_2$ and $h_3$ in (7) and with the same discussion, there is a contradiction.

\item[(13)]$ R_1:h_2h_3^{-2}h_2h_3=1$, $R_2:h_2h_3(h_2h_3^{-1})^2h_3^{-1}h_2=1$, $R_3:h_2h_3^{-1}(h_2^{-1}h_3h_2^{-1})^2h_2^{-1}h_3=1$,\\$R_4:(h_2h_3^{-1}h_2)^2h_3^{-1}h_2h_3^{-1}h_2h_3=1$:\\
By interchanging $h_2$ and $h_3$ in (6) and with the same discussion, there is a contradiction.

\item[(14)]$ R_1:h_2h_3^{-2}h_2h_3=1$, $R_2:h_2h_3(h_2h_3^{-1})^2h_3^{-1}h_2=1$, $R_3:h_2h_3^{-1}h_2(h_3h_2^{-1})^2(h_3^{-1}h_2)^2=1$,\\$R_4:((h_2h_3^{-1})^2h_2)^2h_3=1$:\\
By interchanging $h_2$ and $h_3$ in (8) and with the same discussion, there is a contradiction.

\item[(15)]$ R_1:h_2h_3^{-2}h_2h_3=1$, $R_2:h_2^2h_3^{-1}(h_3^{-1}h_2)^2h_3=1$, $R_3:h_2^3h_3^{-1}h_2^{-1}h_3h_2^{-1}h_3^{-1}h_2=1$,\\$R_4:h_2^4h_3^{-1}(h_3^{-1}h_2)^2=1$:\\
By interchanging $h_2$ and $h_3$ in (1) and with the same discussion, there is a contradiction.

\item[(16)]$ R_1:h_2h_3^{-2}h_2h_3=1$, $R_2:h_2^2h_3^{-1}(h_3^{-1}h_2)^2h_3=1$, $R_3:h_2h_3^{-1}(h_2^{-1}h_3)^3h_2^{-1}h_3^{-1}h_2=1$,\\$R_4:h_2(h_2h_3^{-1})^4h_2h_3=1$:\\
By interchanging $h_2$ and $h_3$ in (2) and with the same discussion, there is a contradiction.

\end{enumerate}


\end{document}